\documentclass[a4paper,reqno,12pt,openany]{amsbook}
\pdfoutput=1
\usepackage[all, 2cell, knot,]{xy} \UseAllTwocells \SilentMatrices
\usepackage{xspace}
\usepackage{amsmath}
\usepackage{amstext}
\usepackage{amsfonts}
\usepackage[mathscr]{euscript}
\usepackage{amscd}
\usepackage{latexsym}
\usepackage{amssymb}
\usepackage{layout}
\usepackage{amsthm}
\usepackage{amsfonts}
\usepackage{amsxtra}
\usepackage{color}
\usepackage[usenames,dvipsnames,svgnames,table]{xcolor}
\usepackage{graphics}
\usepackage{bm, amsbsy}
\usepackage{tabulary}
\usepackage{etoolbox}
\usepackage{fancyhdr}
\usepackage[bookmarks=true,hyperfootnotes=false]{hyperref}
\usepackage{tikz-cd}
\usepackage{tikz}
\usepackage{xargs}
\usepackage{standalone}
\usepackage{relsize}
\usepackage{chngcntr}
\usepackage{longtable}
\usepackage[strict]{changepage}
\usepackage{bbm}
\usepackage[T1]{fontenc}
\usepackage[utf8]{inputenc}
\usepackage[english]{babel}
\usepackage[babel]{csquotes}
\usepackage[style=numeric]{biblatex}
\renewbibmacro{in:}{}

\hypersetup{
			linktoc=all,
            colorlinks=true,
			linkcolor=blue,
			anchorcolor=black,
			citecolor=green,
			urlcolor=black,
}

\usetikzlibrary{arrows}
\usetikzlibrary{arrows.meta}
\usepgflibrary{plotmarks}
\pgfsetplotmarksize{0.035cm}
\usetikzlibrary{positioning}
\usetikzlibrary{calc} % for manipulation of coordinates

\tikzset{
  commutative diagrams/.cd,
  arrow style=tikz,
  diagrams={->}
}

\pgfsetplotmarksize{0.045cm}

\theoremstyle{plain}% default, italic fonts for text of the following items
\newtheorem{theorem}{Theorem}[section]
\newtheorem*{theorem*}{Theorem}  %  Unnumbered theorem
\newtheorem{lemma}[theorem]{Lemma}
\newtheorem{proposition}[theorem]{Proposition}
\newtheorem{corollary}[theorem]{Corollary}
\theoremstyle{definition}% roman font for text of the following items
\newtheorem{definition}[theorem]{Definition}
\newtheorem{example}[theorem]{Example}

\newtheorem{notation}[theorem]{Notation}
\newtheorem{remark}[theorem]{Remark}

\numberwithin{section}{chapter}
\numberwithin{equation}{chapter}
\counterwithin{figure}{chapter}

\newcommand{\clA}{\mathcal{A}}
\newcommand{\clB}{\mathcal{B}}
\newcommand{\clC}{\mathcal{C}}
\newcommand{\clD}{\mathcal{D}}
\newcommand{\clE}{\mathcal{E}}

\newcommand{\clG}{\mathcal{G}}
\newcommand{\clH}{\mathcal{H}}
\newcommand{\clI}{\mathcal{I}}

\newcommand{\clL}{\mathcal{L}}

\newcommand{\clP}{\mathcal{P}}
\newcommand{\clR}{\mathcal{R}}
\newcommand{\clS}{\mathcal{S}}
\newcommand{\clT}{\mathcal{T}}
\newcommand{\clU}{\mathcal{U}}
\newcommand{\clV}{\mathcal{V}}

\newcommand{\za}{\alpha}
\newcommand{\zb}{\beta}
\newcommand{\zg}{\gamma}

\newcommand{\zd}{\delta}
\newcommand{\zD}{\Delta}

\newcommand{\zve}{\varepsilon}

\newcommand{\zs}{\sigma}

\newcommand{\pt}{\partial}
\newcommand{\zgu}[1]{\zg^{(#1)}}

\newcommand{\Dec}{\operatorname{Dec}}

\newcommand{\Ho}{\operatorname{Ho}}
\newcommand{\Id}{\operatorname{Id}}

\newcommand{\pr}{\operatorname{pr}}

\newcommand{\nequ}{\mbox{$n$-equivalence}}
\newcommand{\equ}[1]{\mbox{$(#1)$-equivalence}}
\newcommand{\bsim}{/\!\!\sim}

\newcommand{\nid}{\noindent}
\newcommand{\ds}{\displaystyle}
\newcommand{\bk}{\bigskip}
\newcommand{\mk}{\medskip}
\newcommand{\sk}{\smallskip}
\newcommand{\ovl}[1]{\overline{#1}}
\newcommand{\ovll}[1]{\overset{=}{#1}}
\newcommand{\up}[1]{^{(#1)}}
\newcommand{\lo}[1]{_{(#1)}}
\newcommand{\rw}{\rightarrow}
\newcommand{\Rw}{\Rightarrow}
\newcommand{\lw}{\leftarrow}

\newcommand{\xrw}{\xrightarrow} % use as follows: \xrw{} brackets possibly %empty
\newcommand{\xlw}{\xleftarrow} % use as follows: \xlw{} brackets possibly %empty
\newcommand{\hxrw}[1]{\xymatrix{\ \ar@{^{(}->}^{#1}[r] & \ }}% Extensible right hookarrow
\newcommand{\tiund}[1]{{\times}_{#1}\:}
\newcommand{\pro}[3]{#1\tiund{#2}\overset{#3}{\cdots}\tiund{#2}#1}
\newcommand{\tens}[2]{#1\,\tiund{#2}\,#1}
\newcommand{\uset}[2]{\underset{#1}{#2}}
\newcommand{\oset}[2]{\overset{#1}{#2}}
\newcommand{\ob}{ob\,}
\newcommand{\ata}{A\tiund{B}A}
\newcommand{\mi}{\text{-}}
\newcommand{\nm}{(n-1)}
\newcommand{\bl}{\bullet}
\newcommand{\cop}{\textstyle{\,\coprod\,}}
\newcommand{\seq}[3]{{#1}_{#2}...{#1}_{#3}}
\newcommand{\seqc}[3]{{#1}_{#2},...,{#1}_{#3}}
\newcommand{\ssr}{\!\!\!\!\!\!\!\!\!\!\!\!\!\!}
\newcommand{\cirsm}{\scriptstyle{\circ}\textstyle}
\newcommand{\mlg}[1]{\mathlarger{\mathlarger{\bm #1}}}
\newlength{\myline}% line thickness
\newcommandx*{\doublearrow}[4][1=0, 2=1]{% #1 = shorten left (optional), #2 = shorten right (optionsl),
  \draw[black,line width=\myline,double distance=3pt,#3] #4;
}

\newcommandx*{\triplearrow}[4][1=0, 2=3]{% #1 = shorten left (optional), #2 = shorten right (optionsl),
  \draw[black,line width=\myline,double distance=4pt,#3] #4;
  \draw[black,line width=0.7pt] #4;
}

\newcommand{\dop}[1]{\Delta^{{#1}^{op}}}
\newcommand{\Dop}{\Delta^{op}}
\newcommand{\Dnop}{\Delta^{{n}^{op}}}
\newcommand{\Dmenop}{\Delta^{{n-1}^{op}}}
\newcommand{\cat}[1]{\mbox{$\mathsf{Cat^{#1}}$}}
\newcommand{\Cat}{\mbox{$\mathsf{Cat}\,$}}
\newcommand{\Catn}{\mbox{$\mathsf{Cat^n}$}}
\newcommand{\Gp}{\mbox{$\mathsf{Gp}$}}
\newcommand{\Gpd}{\mbox{$\mathsf{Gpd}$}}
\newcommand{\Af}{A[f]}
\newcommand{\eqr}[1]{\mbox{$\mathsf{EqRel^{#1}}$}}
\newcommand{\cathd}[1]{\mbox{$\mathsf{Cat_{hd}^{#1}}$}}

\newcommand{\catwg}[1]{\mbox{$\mathsf{Cat_{wg}^{#1}}$}}
\newcommand{\gcatwg}[1]{\mbox{$\mathsf{GCat_{wg}^{#1}}$}}
\newcommand{\tawg}[1]{\mbox{$\mathsf{Ta_{wg}^{#1}}$}}
\newcommand{\lta}[1]{\mbox{$\mathsf{LTa_{wg}^{#1}}$}}
\newcommand{\ftawg}[1]{\mbox{$\mathsf{FCat_{wg}^{#1}}$}}
\newcommand{\gtawg}[1]{\mbox{$\mathsf{GTa_{wg}^{#1}}$}}
\newcommand{\gseg}[1]{\mbox{$\mathsf{GSeg_{#1}}$}}

\newcommand{\lnta}[2]{\mbox{$\mathsf{L\lo{#1}Ta_{wg}^{#2}}$}}

\newcommand{\seg}[1]{\mbox{$\mathsf{Seg_{#1}}$}}
\newcommand{\segpsc}[2]{\mbox{$\mathsf{SegPs}$}\funcat{#1}{#2}}
\newcommand{\Ps}{\mbox{$\mathsf{Ps}$}}
\newcommand{\psc}[2]{\mbox{$\mathsf{Ps}$}\funcat{#1}{#2}}
\newcommand{\PsTalg}{\mbox{\sf Ps-}T\mbox{\sf -alg}}
\newcommand{\Top}{\mbox{$\mathsf{Top}$}}
\newcommand{\muk}{\mu_k}
\newcommand{\hmu}[1]{\hat\mu_{#1}}
\newcommand{\hmuk}{\hat{\mu}_k}

\newcommand{\Nn}{N_{(n)}}
\newcommand{\N}[1]{N_{(#1)}}
\newcommand{\Nb}[1]{N_{(#1)}}
\newcommand{\Nu}[1]{N^{(#1)}}
\newcommand{\di}[1]{d^{(#1)}}
\newcommand{\dn}{d^{(n)}}
\newcommand{\bd}{\bar d}
\newcommand{\Dn}{D_{n}}
\newcommand{\Dnm}{D_{n-1}}
\newcommand{\tld}{\tilde{d}}
\newcommand{\D}[1]{D_{#1}}

\newcommand{\p}[1]{p^{(#1)}}
\newcommand{\bp}{\bar p}

\newcommand{\pn}{p^{(n)}}
\newcommand{\op}[1]{\bar{p}^{(#1)}}
\newcommand{\q}[1]{q^{(#1)}}

\newcommand{\qn}{q^{(n)}}

\newcommand{\Tan}{\mbox{$\mathsf{Ta^{n}}$}}
\newcommand{\ta}[1]{\mbox{$\mathsf{Ta^{#1}}$}}
\newcommand{\gta}[1]{\mbox{$\mathsf{GTa^{#1}}$}}
\newcommand{\Set}{\mbox{$\mathsf{Set}$}}
\newcommand{\St}{St\,}

\newcommand{\Tr}{Tr\,}
\newcommand{\tr}[1]{Tr_{#1}}
\newcommand{\uh}{\underline{h}}
\newcommand{\uk}{\underline{k}}
\newcommand{\uv}{\underline{v}}
\newcommand{\ur}{\underline{r}}
\newcommand{\us}{\underline{s}}
\newcommand{\rz}{R_0}
\newcommand{\Qn}{Q_{n}}
\newcommand{\Qnm}{Q_{n-1}}
\newcommand{\Discn}{Disc_{n}}

\newcommand{\funcat}[2]{[\Delta^{{#1}^{op}},#2]}
\newcommand{\ps}{\sf{ps}}

\newcommand{\Lb}[1]{\mbox{$\mathsf{L_{(#1)}}$}}
\newcommand{\nfol}{$n$-fold }

\newcommand{\tms}[2]{\overset{#1}{\times}_{#2}}
\newcommand{\btil}{\widetilde{B}}
\newcommand{\ord}{\operatorname{Or}}
\newcommand{\orn}[1]{\ord_{(#1)}}
\newcommand{\nty}{\mbox{$n$-types}}
\newcommand{\gpd}[1]{\mbox{$\mathsf{Gpd}^{#1}$}}
\newcommand{\gpdwg}[1]{\mbox{$\mathsf{Gpd_{wg}^{#1}}$}}
\newcommand{\Rbt}[1]{\clR_{#1}}
\newcommand{\Lbt}[1]{\clL_{#1}}
\newcommand{\Cube}[1]{\mathrm{Cube}(#1)}
\newcommand{\Sc}[1]{\scriptstyle{#1}}
\newcommand{\Scc}[1]{\scriptscriptstyle{#1}}

\newcommand{\pagelabel}[1]{\phantomsection\label{#1}}
\makeindex

\bibliography{C:/Users/SIMONA/Documents/BOOK ALL FILES/VERSION 5 Feb/bibliography/BIB_LIB}{}

\makeatletter
\renewcommand{\l@figure}{\@tocline{0}{3pt plus2pt}{0pt}{2.5pc}{}}
\makeatother

\let\LaTeXStandardTableOfContents\tableofcontents
\renewcommand{\tableofcontents}{%
\begingroup%
\renewcommand{\pmb}{\relax}%
\LaTeXStandardTableOfContents%
\endgroup%
}%

\makeatletter
\def\@tocline#1#2#3#4#5#6#7{\relax
   \ifnum #1>\c@tocdepth % then omit
   \else
     \par \addpenalty\@secpenalty\addvspace{#2}%
     \begingroup \hyphenpenalty\@M
     \@ifempty{#4}{%
       \@tempdima\csname r@tocindent\number#1\endcsname\relax
     }{%
       \@tempdima#4\relax
     }%
     \parindent\z@ \leftskip#3\relax \advance\leftskip\@tempdima\relax
     \rightskip\@pnumwidth plus4em \parfillskip-\@pnumwidth
     #5\leavevmode\hskip-\@tempdima #6\nobreak\relax
     \ifnum#1<0\hfill\else\dotfill\fi\hbox to\@pnumwidth{\@tocpagenum{#7}}\par
     \nobreak
     \endgroup
   \fi}
\makeatother

\begin{document}
\frontmatter

\title {Segal-type models of higher categories}

\author{Simona Paoli}
\address{\small{Department of Mathematics, University of Leicester, UK}}
% \email{sp424@le.ac.uk}

\subjclass[2010]{18D05, 55U10 }
\keywords{weak $n$-categories, $n$-fold categories, higher groupoids, multi-simplicial objects, homotopy types }
\date{21 June 2017}

\begin{abstract}

Higher category theory is an exceedingly active area of research, whose rapid growth has been driven by
its penetration into a diverse range of scientific fields. Its influence extends through key mathematical
disciplines, notably homotopy theory, algebraic geometry and algebra, mathematical physics, to encompass important
applications in logic, computer science and beyond. Higher categories
provide a unifying language whose greatest strength lies in its ability to bridge between diverse areas and
uncover novel applications.

In this foundational work we introduce a new approach to higher categories. It builds upon the theory of
iterated internal categories, one of the simplest possible higher categorical structures available, by
adopting a novel and remarkably simple "weak globularity" postulate and demonstrating that the resulting
model provides a fully general theory of weak $n$-categories. The latter are among the most complex of the higher structures, and are crucial for applications. We show that this new model of
"weakly globular $n$-fold categories" is suitably equivalent to the well studied model of weak $n$-categories
due to Tamsamani and Simpson.

\end{abstract}

\maketitle

\tableofcontents

%%%%%%%%%%%%%%%%%%%%%%%%%%%%%%%%%%%%%%%%%%%%%%%%%%%%%%%%%%%%%%%%%%%%%%%%%%

\newpage

% 2nd page, thanks message
%-------------------------------------------------------------------------------
%\thispagestyle{empty}

%%%%%%%%%%%%%%%%%%%%%%%%%%%%%%%%%%%%%%%%%%%%%%%%%%%%%%%%%%%%%%%%%%%%%%%%%%%%
% General definitions for all Chapters
%-------------------------------------------------------------------------------

%% Define Page style for all chapters
\pagestyle{fancy}
%% Delete the current section for header and footer
\fancyhf{}
%% Set custom header
\lhead[Segal-type models of higher categories]{\thepage}
\rhead[\thepage]{Simona Paoli}

\setlength{\baselineskip}{1.1\baselineskip}
%%%%%%%%%%%%%%%%%%%%%%%%%%%%%%%%%%%%%%%%%%%%%%%%%%%%%%%%%%%%%%%%%%%%%%%%%%%%%%

\chapter*{Preface}
The theory of higher categories is a very active area of research and is penetrating diverse fields of science.

Historically the subject was motivated by questions in algebraic topology and mathematical physics, two areas where the most important applications are currently found. Algebraic geometry also makes use of higher categorical notions. More recently higher categories are penetrating logic and computer science, and also start to appear in algebra and representation theory. Higher categories can sometimes be used as a common language to describe complex phenomena occurring in these areas.

A plethora of different approaches to higher categories have been developed over the years. Each one
represents certain relevant aspects of the abstract notion being modelled, often with a view to supporting
a particular ecosystem of applications. At this stage no one approach suits all such contexts, and indeed
one might doubt the viability of a universally applicable model. It appears instead that the most prudent
approach is to continue the development of all of these important strands, relating them where necessary
by explicit comparisons.

The purpose of this monograph is to introduce a new approach to working with higher categories: this is based on a simple higher categorical structure consisting of iterated internal categories (also called $n$-fold categories) as well as on a new paradigm  to weaker higher categorical structures which is the idea of weak globularity.

We show that our new model, called weakly globular \nfol categories, is suitably equivalent to one of the models of higher categories that was studied in greatest depth, the one introduced by Tamsamani \cite{Ta} and further studied by Simpson \cite{Simp}.

  We achieve this comparison by developing a larger context of 'Segal-type models of weak $n$-categories', based on multi-simplicial structures, of which both the Tamsamani model and weakly globular \nfol categories are special cases.

The use of simplicial structures to capture higher coherence phenomena has a long history in algebraic topology, starting with the work of Graeme Segal \cite{Segal1974}, and then the study of categories enriched in simplicial sets by Vogt, Dwyer, Kan, Smith \cite{DwyKan1980-1}, \cite{DwyKanSmi1989} and by others. More recently, simplicial techniques underpin the development of so called $(\infty, n)$-categories, where several models have been developed and studied by Bergner \cite{Be2} \cite{BeRe}, Barwick and Kan \cite{Bk1}, Lurie \cite{Lu1}, Joyal \cite{Jo}, Rezk \cite{Re1} \cite{Re2} and others.

Simplicial models of $(\infty,\infty)$-categories have been developed by Verity \cite{Vecom}, building upon insights
from the study of simplicial nerves of strict n-categories initiated by Street \cite{Street1987}.

In this work we concentrate on higher structures in the 'truncated' case, where there are higher morphisms only in dimensions $0$ up to $n$. This is intimately connected to the Postnikov tower in algebraic topology. In fact the algebraic modelling of the building blocks of spaces, the $n$-types, which are Postnikov sections of spaces, is related to models of weak $n$-categories via the so called 'homotopy hypothesis': a good model of weak $n$-categories should give an algebraic model of $n$-types in the weak $n$-groupoid case. We show that our model does satisfy the homotopy hypothesis.

There are long standing open questions about weak $n$-categories, both within category theory and in its applications to homotopy theory:  for instance the comparison between the simplicial and higher operadic models of higher categories and the algebraic description of the $k$-invariants of spaces.

The present work provides a platform where these and other open questions can be studied. This however goes beyond the scope of this work, whose goal is to lay the foundations of this theory.

The potential of our model to tackle these open questions comes from one of the main novelties of our approach here: the use of an entirely rigid structure, namely a subcategory of $n$-fold categories, to model weak $n$-categories.

 The terminology 'rigid structure' refers to the fact that $n$-fold categories, being iterated internal categories, have associative and unital compositions in $n$ different simplicial directions. In this sense, \nfol categories are a strict higher categorical structure, though not the same as strict $n$-categories, since the higher morphisms in dimensions $0$ to up $n$ do not form just a set. In our model, the higher morphisms in dimension $k$ have themselves a $(n-1+k)$-fold categorical structure of a special type which is suitably equivalent to a discrete structure (that is, a set): this is the 'weak globularity condition'.

$n$-Fold structures were used in homotopy theory by Loday \cite{Loday1} for the modelling of connected $(n+1)$-types via cat$^n$-groups. The idea of weak globularity was first introduced by the author in \cite{Pa} in an internal setting for the category of cat$^n$-groups: weakly globular cat$^n$-groups were shown in \cite{Pa} to be algebraic models of connected $(n+1)$-types; weak globularity  was extended and further studied by Blanc and the author in \cite{BP} in the context of general $n$-types, for which an analogue of Loday's model was not available.

 None of these works however captured the general categorical case. This necessitates many novel ideas and techniques, such as the use of pseudo-functors to model higher structures and the construction of a rigidification functor from Tamsamani model to weakly globular $n$-fold categories.

  This work uses a blend of techniques from category theory and simplicial homotopy theory, reviewed in the background Part \ref{part-back}. We therefore hope this work will be accessible both to category theorists and to algebraic topologists.

\bk

There are five parts in the organization of this work:
\bk

 \begin{itemize}
 \item [] \textbf{Part \ref{part-1} Higher categories}.

  \nid This part aims to provide the reader with a guide for the rest of this work. It contains a broad introduction to higher categories, some historical development of the notion of weak globularity and a non-technical overview of the main ideas and results of this work. \bk

   \item [] \textbf{Part \ref{part-back} Background}.

    \nid This part covers the main techniques from category theory and simplicial homotopy theory used in this work.

    \bk

   \item []  \textbf{Part \ref{part-3} Weakly globular $\pmb{n}$-fold categories and Segalic pseudo-functors}.

    \nid In this part we introduce the main new structure of this work, the category $\catwg{n}$ of weakly globular \nfol categories, and we study its relation to a class of pseudo-functors which we call Segalic pseudo-functors.

    \bk

   \item [] \textbf{Part \ref{part-4} Weakly globular Tamsamani $\pmb{n}$-categories and their rigidification}.

   \nid This part is devoted to another new structure, the category $\tawg{n}$ of weakly globular Tamsamani $n$-categories.  The main goal of this part is the construction of the rigidification functor from weakly globular Tamsamani $n$-categories to weakly globular \nfol categories.

   \bk

   \item []  \textbf{Part \ref{part-5} Weakly globular $\pmb{n}$-fold categories as a model of weak $n$-categories}.

\nid This part contains the construction of the discretization functor from weakly globular $n$-fold categories to Tamsamani $n$-categories, and the final results: the comparison between $\catwg{n}$ and $\ta{n}$, exhibiting $\catwg{n}$ as a model of weak $n$-categories, and the proof of the homotopy hypothesis.

 \end{itemize}

\chapter*{Acknowledgements}
This work was supported by the EU International Reintegration Grant HOMALGHIGH No 256341 of which I was principal investigator. It also received financial support from the University of Leicester, where I worked since 2011 and which supported my study leave in 2015; from the Centre of Australian Category Theory at Macquarie University which hosted me during August-December 2015 and from the University of Chicago where I visited in 2016.

I am grateful for the opportunity I had to present this work in several talks, where I could gather invaluable feedback from colleagues. In particular I thank the Centre of Australian Category Theory for the opportunity to give a long series of talks about this work at the Australian Category Seminar, and I am particularly indebted for useful feedback from its members, especially Michael Batanin, Steve Lack, Richard Garner, Ross Street, Mark Weber, Dominic Verity.

 I also would like to thank Peter May for giving me the chance to give a series of talks about this work at the University of Chicago and for many useful comments.

 I am grateful for the opportunity to present this work at international conferences and at seminars and for the discussions with experts in the field that followed. I thank in particular Joachim Kock, Martin Hyland, Dorette Pronk for interesting discussions and feedback. I also thank David Blanc, Frank Neumann, Alexander Kurz and Dominic Verity for reading the introductory part of this work and for helpful comments, and Nicola Gambino for some bibliographical suggestions.

I finally thank all my colleagues in the Department of Mathematics of the University of Leicester for ongoing support and encouragement during this work.

\clearpage

%%%%%%%%%%%%%%%%%%%%%%%%%%%%%%%%%%%%%%%%%%%%%%%%%%%%%%%%%%%%%%%%%%%%%%%%%%
% Define block styles for flowchart
%
% Rounded rectangular blocks
%
% Blocks definition

\tikzstyle{nullblock} = [rectangle, fill=white,
    text width=5mm, text centered, minimum height=5mm]

% block = size w=45mm, h=35mm
\tikzstyle{block} = [rectangle, draw, fill=gray!10,
    text width=45mm, text centered, rounded corners, minimum height=30mm]

% block1 = size w=65mm, h=25mm
\tikzstyle{block1} = [rectangle, draw, fill=gray!10,
    text width=65mm, text centered, rounded corners, minimum height=25mm]

%% block2 = size w=65mm, h=20mm
\tikzstyle{block2} = [rectangle, draw, fill=gray!10,
    text width=65mm, text centered, rounded corners, minimum height=20mm]

% block3 = size w=45mm, h=20mm
\tikzstyle{block3} = [rectangle, draw, fill=gray!10,
    text width=45mm, text centered, rounded corners, minimum height=20mm]
%
% block4 = size w=70mm, h=20mm
\tikzstyle{block4} = [rectangle, draw, fill=gray!10,
    text width=100mm, text centered, rounded corners, minimum height=20mm]
%
% block5 = size w=70mm, h=40mm
\tikzstyle{block5} = [rectangle, draw, fill=gray!10,
    text width=100mm, text centered, rounded corners, minimum height=40mm]
\tikzstyle{line} = [draw, -latex']

%%%%%%%%%%%%%%%%%%%%%%%%%%%%%%%%%%%%%%%%%%%%%%%%%%%%%%%%%%%%%%%%%%%%%%%%%%

%%%%%%%%%%%%%%%%%%%%%%%%%%%%%%%%%%%%%%%%%%%%%%%%%%%%%%%%%%%%%%%%%%%%%%%%%%%%%%%

\mainmatter

% Set arabic (1,2,3...) page numbering
\pagenumbering{arabic}

%%%%%%%%%%%%%%%%%%%%%%%%%%%%%%%%%%%%%%%%%%%%%%%%%%%%%%%%%%%%%%%%%%%%%%%%%%%%%%%
\part{Higher categories}\label{part-1}

Part \ref{part-1} aims to provide the reader with a guide for the rest of this work. Chapter \ref{chap1} contains an introduction to higher categories, while Chapter \ref{chap1-2} is an introduction to the three Segal-type models studied in this work: the Tamsamani model and the two new models we introduce here, called weakly globular \nfol categories and weakly globular Tamsamani $n$-categories.

Our overview of higher categories in Chapter \ref{chap1} aims to highlight the open questions that led to our approach to weak higher categories alongside with providing some context of their development. We do not aim to give a comprehensive detailed survey of different models of weak higher categories, and we have provided several bibliographical references where further information can be found.

After giving some motivation and context, we describe in Chapter \ref{chap1} the main classes of higher structures and highlight their differences: strict versus weak higher structures, truncated versus non-truncated case, and the class of $n$-fold categories, which is central to this work.

We then concentrate in Section \ref{homotyp} on one of the most important connections between higher category theory and homotopy theory, which is the algebraic modelling of the building blocks of spaces, the $n$-types. Once again, rather than giving a detailed comprehensive survey of all the different models, we point the reader towards bibliographical references and then concentrate on a model of path-connected $(n+1)$-types via $n$-fold structures due to Loday, called cat$^n$-groups, which developed independently on models of weak higher categories. The question of how cat$^n$-groups compare to models of weak higher groupoids is what first led the author to introduce the notion of weak globularity, in an internal context inside the category of groups.

In Chapter \ref{chap1-2} we explain the idea of weak globularity in the general categorical context. Our aim in this chapter is to convey some of the ideas behind our constructions and to give a summary of the main results. In particular, in Section \ref{sec-comm-fea} we give an account of the main common features of the three Segal-type models of this work.

Although Chapter \ref{chap1-2} is written in a way that avoids technical details, at times we need to refer to some basic notions and notations which are fully explained in Part \ref{part-back} on background techniques, and we refer to the appropriate sections there when needed.

We end in Section \ref{sub-org} with a description of the overall organization of the rest of this work, and with a  diagrammatic summary in Figure \ref{FigIntro-2}.

\clearpage

\chapter{An introduction to higher categories}\label{chap1}

In this chapter we give a non-technical introduction to higher categories. In Section \ref{sec-motiv-cx} we describe some of the contexts that inspired and motivated their development. In Section \ref{sec-diff} we explain the idea of higher categories, and the different classes of higher structures.

In Section \ref{homotyp} we discuss one of the most important occurrences of higher categories in algebraic topology, which is the algebraic modelling of homotopy types. We give an account of the use of $n$-fold structures in modelling $n$-types, including the work by the author. This provides an historical development of the notion of weak globularity, which is central to this work and is discussed in more details in the next chapter.

\section{Motivation and context}\label{sec-motiv-cx}

The language of categories and functors permeates modern mathematics. In a category we have objects, morphisms, compositions of morphisms and identity morphisms for each object, such that compositions are associative and unital . When each morphism is invertible we obtain a groupoid. A one-object groupoid is the familiar notion of a group, while a one-object category is a monoid.
\index{Groupoid}
\index{Monoid}

The maps between categories are the functors: These associate objects to objects and arrows to arrows in a way that is compatible with the composition and with the identities.

Many familiar mathematical structures form a category: for instance vector spaces and linear maps, topological spaces and continuous maps and so on.

The idea of a higher category was prompted by several inputs. First of all, there are many natural examples of higher structures. An important one is the $2$-category $\Cat$ of small categories. This comprises $0$-dimensional data  which are the categories (the objects), $1$-dimensional data which are the functors ($1$-morphisms between the objects), and $2$-dimensional data which are the natural transformations between functors ($2$-morphisms between $1$-morphisms).

 Functors can be composed: given functors $F:\mathcal{A}\rw \mathcal{B}$ and $G:\mathcal{B}\rw \mathcal{C}$ the composite functor $G\circ F:\mathcal{A}\rw \mathcal{C}$ associates to each object $a \in \mathcal{A}$ the object $G(F(a))\in \mathcal{C}$ and to each morphism $f \in \mathcal{A}$ the morphism $G(F(f))\in \mathcal{C}$. Natural transformations can also be composed, but in two different ways. Given functors $F,G,H: \mathcal{A} \rw \mathcal{B}$ between categories $\mathcal{A}$ and $\mathcal{B}$, and natural transformations $\alpha: F \Rightarrow G$ and $\beta:F \Rightarrow G $ we can form a 'vertical' composite natural transformation $\beta \circ_v \alpha: G \rw H$ with components given by the composites
  \begin{equation*}
  (\beta \circ_v \alpha)_a= \beta_a \alpha_a : F(a) \rw H(a)
  \end{equation*}
   for each object $a\in \mathcal{A}$. Given functors $F,F': \mathcal{A} \rw \mathcal{B}$ between categories $\mathcal{A}$ and $\mathcal{B}$, functors $G,G': \mathcal{B} \rw \mathcal{C}$ between categories $\mathcal{B}$ and $\mathcal{C}$ and natural transformations $\alpha: F \Rightarrow F'$ and $\beta: G\Rightarrow G'$ we can form the 'horizontal' composite natural transformation $\beta \circ_h \alpha: GF \rw G'F'$ with components given by the composites
    \begin{equation*}
    (\beta \circ_h \alpha)_a= (\beta_{F'(a)})\: (G\alpha_a): GF(a) \rw G'F'(a)
    \end{equation*}
     for each object $a\in \mathcal{A}$.
   We can associate a geometric picture to these data by associating points to categories (the objects), arrows to functors (the $1$-morphisms) and globes to natural transformations (the $2$-morphisms), with the two different compositions of natural transformations pictured as vertical and horizontal compositions of globes, as illustrated below:

 \bk
\pagelabel{2cat-begin}
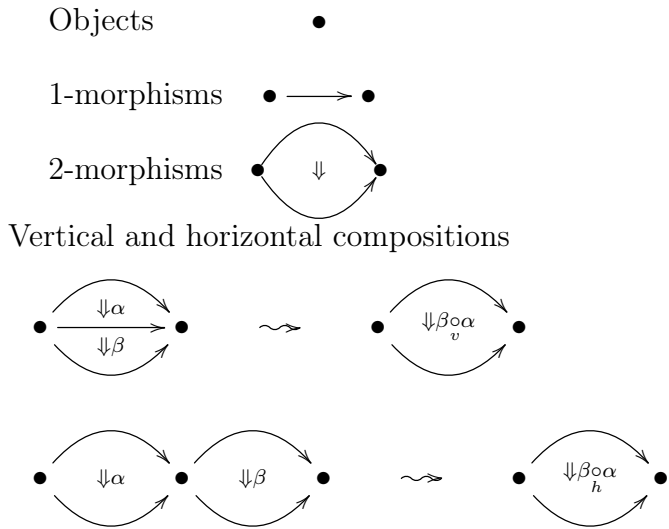
\begin{tabular}{ l l c}
    & Objects  & $\bullet$ \\
    & & \\
    & 1-morphisms &
   $
   \xymatrix{
   \bullet  \ar[r] & \bullet
   }
   $  \\
   &&\\
    & 2-morphisms &
    $
    \entrymodifiers={[o]}
   \xymatrix@C=40pt{{\bullet} \ar@/^1.5pc/[r]\ar@{}[r]|-{\Downarrow} \ar@/_1.5pc/[r] & {\bullet}}
   $
\end{tabular}
\bk

Vertical and horizontal compositions
  \begin{flalign*}
  & \qquad
    \xymatrix@C=40pt{
    {\bullet}  \ar@/^1.5pc/[r]_{} \ar@/_1.5pc/[r]^{} \ar@{->}^{\Downarrow\alpha}_{\Downarrow\beta}[r] & {\bullet}}
    \quad
    \xymatrix@C=40pt@R=15pt @C=15pt{
    \text{ } \ar@{~>}[r]& \text{ }}
    \quad
    \xymatrix@C=40pt{
    {\bullet}  \ar@/^1.5pc/[r]_{} \ar@/_1.5pc/[r]^{} \ar@{}|{\Downarrow\beta\underset{v}{\circ}\alpha}[r] & {\bullet}}
  &
  \end{flalign*}

  \begin{flalign*}
  & \qquad
    \xymatrix@C=40pt{{\bullet} \ar@/^1.5pc/[r]\ar@{}[r]|-{\Downarrow\alpha} \ar@/_1.5pc/[r] & {\bullet}\ar@/^1.5pc/[r]\ar@{}[r]|-{\Downarrow\beta} \ar@/_1.5pc/[r] & {\bullet}}
    \quad
    \xymatrix@C=40pt@R=15pt @C=15pt{
    \text{ } \ar@{~>}[r]& \text{ }}
    \quad
    \xymatrix@C=40pt{
    {\bullet}  \ar@/^1.5pc/[r]_{} \ar@/_1.5pc/[r]^{} \ar@{}|{\Downarrow\beta\underset{h}{\circ}\alpha}[r] & {\bullet}}
  &
  \end{flalign*}

Every $2$-category comprises data as above, and we also call the objects '$0$-cells', the $1$-morphisms '$1$-cells', the $2$-morphisms '$2$-cells'. So we see that every $k$-cell (for $k=1,2$) has a $(k-1)$-cells as its source and target. The highest dimension of cells in this example is $2$, and this is called the 'dimension' of the higher structure.\index{Cells}

The idea of a higher category in dimension greater than 2 is to have higher cells of globular shape, with each $k$-dimensional cell having $(k-1)$-cells as its source and target. When cells are present only in dimensions $0$ up to $n$ we say that the higher structure is $n$-dimensional. So for instance when $n=3$ we can think of a $3$-cell as a 'globe between globe', that is a sphere.

Many higher categorical ideas have their root in the notion of homotopy coherence \index{Homotopy coherence} in algebraic topology. The latter developed along several directions, leading to diverse applications. One was the study of loop spaces,\index{Loop space} leading in an algebraic setting to the notions of $H_\infty$ and $E_\infty$ spaces: the works of Boardman and Vogt \cite{Boardman1973}, May \cite{May1977}, Stasheff \cite{Stasheff}, Segal \cite{Segal1968}, Sugawara \cite{Sugawara} are relevant here. Operads \index{Operads}also developed as a way to encode higher homotopy coherences, see for instance the works of May \cite{May1972}, Loday \cite{Loday1996}, Markl \cite{MarklShnider}.

The abstraction and 'categorification' of these models of homotopy coherence led to several combinatorial approaches to higher categorical structures, in particular the Segal-type model of Tamsamani and Simpson \cite{Simp}, \cite{Ta} and the higher operadic models of Batanin \cite{B}, also studied by Leinster \cite{Le2}, Weber \cite{BatWeb2011}, \cite{BatCisWeb2013} and others.

Another way to encode higher order homotopical information is via the notion of model category pioneered by Quillen \cite{QuillenLectureNotes1967}, and the one of simplicial categories and their localizations, which was studied by Dwyer, Kan and others \cite{DwyKan1980-1}, \cite{DwyKan1980-2}, \cite{DwyKan1980-3}, \cite{DwyKanSmi1989}, \cite{DwyerKanSmith1986}.

 Quillen model categories \index{Quillen model structure} remain a key tool in algebraic topology and are increasingly tackled from a categorical perspective, as shown for instance in the works of Garner \cite{Garner2012} and Riehl \cite{Riehl2014}.

 Simplicial categories\index{Simplicial!- category} are one of the models of a class of higher structures called $(\infty,1)$-categories, which have become a central object of study in modern homotopy theory, encoding the idea of 'homotopy theory of homotopy theories'. Several other models have also been developed, leading to a variety of applications, as further explained in the next section.

 Algebraic structures in $(\infty,1)$-categories led to the notion of $\infty$-operads\index{Models of $\infty$-operads}: See the work of Lurie \cite{LuHiAlg}, as well as the dendroidal sets model\index{Dendroidal sets} of Moerdijk and Weiss \cite{MoerdijkWeiss2009}, further studied by Heuts, Hinich, Cisinski and others \cite{HeutsHinich}, \cite{CisinskiMoerdijk2013}, \cite{CisinskiMoerdijk2011}. A further extension of these approaches is developed by Hackney and Robertson \cite{Hackney2015}.

Another motivating force for the development of higher category theory coming from algebraic topology was the algebraic modelling of the Postnikov systems of spaces. This is quite central to this work, and we explain this context in further detail in section \ref{homotyp}.

Mathematical physics has also been inspiring many developments in higher categories, in the pursuit of models for TQFT \index{TQFT} and higher cobordism categories. Several conjectures in this direction were formulated by Baez and Dolan \cite{BD}, \index{Baez-Dolan conjectures}the cobordism hypothesis was recently proved by Lurie \cite{Lu3}, higher categories and low-dimensional TQFTs were investigated by Schommer-Pries \cite{Schommer2009}, \cite{SchommerChristopher}.

Algebraic geometry also has seen the use of higher categorical ideas in the pursuit of the notion of higher and derived stacks \index{Higher stacks} as well as higher non-abelian cohomology \index{Non-abelian cohomology}: see for instance the works of Hirschowitz, Simpson \cite{Simp}, \cite{HirschowitzSimpson}, Pridham \cite{Pridham}, To\"{e}n \cite{MoerToen2010}, \cite{Toen2009}, \cite{Toen2014}, To\"{e}n and Vezzosi \cite{ToenVezzosi2005}, \cite{ToenVezzosi2008}.

More recently higher categories entered logic and computer science in the area of homotopy type theory,\index{Homotopy type theory} see for instance the book of the Univalent Foundations Project \cite{HTT}, and the works of Voevodsky \cite{Voe}, LeFanu Lumsdaine and Kapulkin \cite{Kapulkin}\cite{LeFanu2010},  Awodey and Warren\cite{AwodeyWarren2009}, van den Berg and Garner \cite{BergGarner2011}. Other recent applications of categorical structures are in the areas of quantum computing,\index{Quantum computing} see for instance the works of Coecke, Kissinger and Vicary \cite{CoeckeKissinger}, \cite{JamieVicary}. Higher categories have also given rise to interesting software implementations, see the works of Bar, Kissinger and Vicary \cite{BarKissingerVicary}, \cite{BarVicary}.

\section{Different types of higher structures}\label{sec-diff}
The behaviour of compositions of cells in a higher category determine two main classes: strict and weak higher categories. As further explained in Section \ref{subs-svw} below, in the strict case, compositions are associative and unital; in the weak case, they are associative and unital only up to coherent isomorphims.

For each of these classes, there are higher categories which admit cells in every dimension ($\omega$-categories), and those that have cells in dimensions only $0$ up to $n$ (truncated $n$-categories). A further class of higher structures central to this work is the one of $n$-fold structures. Below we give a description of these different types of higher structures and some of their relationships.

\subsection{$\pmb{\omega}$-Categories}
These have been studied extensively in relation to applications to homotopy theory, mathematical physics and algebraic geometry, giving rise to several models of $(\infty,n)$-category; intuitively, the latter are weak higher categories admitting cells in all dimensions and with weakly invertible arrows in dimension higher than $n$.

  There are several models of $(\infty,1)$-category,\index{Models of $(\infty,1)$-category} all of which are Quillen equivalent: quasicategories,\index{Quasicategories} introduced by Boardman and Vogt \cite{Boardman1973} under the name of 'weak Kan complexes' and much developed by Joyal \cite{Jo} and Lurie \cite{Lu1}; simplicial categories\index{Simplicial!- category} introduced by Dwyer and Kan \cite{DwyKan1980-1}, \cite{DwyKan1980-2}, more recently studied with a model category approach by Bergner \cite{Be2}; complete Segal spaces \index{Complete Segal spaces} studied by Rezk \cite{Re2}; relative categories,\index{Relative categories} studied by Barwick and Kan \cite{Bk1}. The survey paper of Bergner \cite{Bergner2010} gives a description of these different models and their Quillen equivalences.
  More recently quasi-categories have been studied using techniques of 2-category theory and monad theory by Riehl and Verity \cite{RiehlVerity2015}, \cite{RiehlVerity2015b}, and a model of $(\infty,1)$-categories in terms of internal categories in simplicial sets was studied by Horel \cite{Horel2015}.

   Models of $(\infty,n)$-categories \index{Models of $(\infty,n)$-categories} for $n>1$ have been studied by Ara, \cite{Ara2014}, Barwick and Kan \cite{BarwickKan2013} Bergner and  Rezk \cite{Be3} \cite{BeRe} \cite{Re1}, Lurie \cite{Lu3} and played an important role in Lurie's proof of the cobordism hypothesis \index{Cobordism hypothesis}. An axiomatic approach to $(\infty,n)$-categories was developed by Barwick and Schommer-Pries \cite{BarwickSchommer}.

  The most general kind of weak $\omega$-category possible would admit cells at all dimensions without
stipulating that all cells should be weakly invertible above some finite dimension. Verity developed the theory of complicial sets \cite{Vecom} \cite{Ve} to model these $(\infty,\infty)$-categories as an adaptation of the theory of strict complicial sets \cite{Verity2008}, which he developed to prove the the Street-Roberts conjecture \cite{Street1987}\index{Street-Roberts conjecture} on the characterization of nerves of strict $n$-categories.
\index{Infinity categories}
\index{Weak $\omega$-category}
\index{Complicial sets}

\subsection{Truncated higher categories} \index{Truncated higher categories}In this work we concentrate on 'truncated' higher categories, with cells only in dimensions $0$ up to $n$. This relates to one of the original motivations for the development of higher categories, namely the algebraic modelling of the Postnikov systems of spaces, whose sections are the $n$-types, that is spaces with trivial homotopy groups in dimension higher than $n$.

The largely open problem of understanding algebraic invariants such as the higher homotopy and
cohomology operations leads us inexorably to the question of unravelling the combinatorics of Postnikov
systems of spaces and simplicial categories. This is connected with achieving a useful
combinatorial description of the $k$-invariants of spaces, another open problem of some significant merit.
The work of Baues \cite{Baues2006}  provides a low dimensional and stable exemplar in this direction, by
demonstrating the utility of this approach to computations of some differentials in the Adams spectral
sequence.\index{Higher homotopy operations} \index{Higher cohomology operations}\index{Adams spectral sequence}

Applications of the truncated case to $E_{n}$-structures and to Hochschild cohomology \index{Hochschild cohomology}were developed in the context of the higher operadic model of Batanin \cite{Batanin2007}, see also the work of Tamarkin on the Deligne conjecture \index{Deligne conjecture} \cite{Tamarkin2007}. Simpson in \cite{Simp} envisages the use of Tamsamani $n$-categories for applications to algebraic geometry, in the theory of higher stacks and of non-abelian cohomology.\bk

\subsection{Strict versus weak $\pmb{n}$-categories}\label{subs-svw} In a strict higher category, compositions of cells are associative and unital, and there is a simple way to describe strict $n$-categories via iterated enrichment.\index{Strict $n$-category} Although simple to define, strict $n$-categories are insufficient for many applications and the wider class of weak $n$-categories is needed. For instance, strict $n$-groupoids \index{Strict $n$-groupoid}do not model $n$-types in dimension $n>2$ (see \cite{S2} for a counterexample in dimension $n=3$).

In a weak $n$-category,\index{Weak $n$-category} higher cells compose in a way that is associative and unital only up to an invertible cell in the next dimension, and these associativity and unit isomorphisms are suitably compatible or coherent.

In dimensions $n=2$ and $n=3$ the idea of a weak $n$-category is embodied in the classical notions of bicategory due to B\'{e}nabou \cite{Ben} and tricategory due to Gordon, Power and Street \cite{GPS}, and more recently studied by Gurski \cite{Gurski2013}, Garner \cite{GarnerGurski2009} and others. In these structures, explicit diagrams encode the coherence axioms for the associativity and unit isomorphisms.

Capturing the coherence axioms explicitly in dimension $n>3$ seems intractable. Instead, various combinatorial machines have emerged to automate the process of defining weak $n$-categories \cite{L1}: in these approaches the coherence data for the higher associativity are not given explicitly but they are automatically encoded in the combinatorics defining the models.
\index{Bicategory} \index{Tricategory} \index{Coherence axioms}

Different types of combinatorics have been used, including multi-simplicial structures as in Tamsamani and Simpson \cite{S2}, \cite{Ta}, higher operads as in Batanin  \cite{B}, Leinster \cite{L1}, Weber \cite{BatWeb2011}, \cite{BatCisWeb2013} and Trimble  \cite{Cheng2}, opetopes\index{Opetopes} as in Baez, Dolan\cite{BD2} and as in Cheng \cite{Cheng1} and several others. \index{Models of higher categories} The comparison between these different approaches is still largely an open problem.

\subsection{$\pmb{n}$-Fold categories} There is a third class of higher structures besides strict $n$-categories and weak $n$-categories which is central to this work: the class of $n$-fold categories.
 $n$-Fold categories were introduced by Ehresmann \cite{ACEhresII}, \cite{ACEhresIII}, \cite{ACEhresIV}. There is an extensive literature for the case $n=2$ (when they are called double categories), developed among others by Dawson, Grandis, Pare, Pronk, see for instance \cite{DawParPro2010}, \cite{GrandisPare1999}, \cite{GrandisPare2004}. Model structures on double categories \index{Quillen model structure}were developed in joint work by the author \cite{FiorePaoliPronk2008} while Fiore and the author built a  model structure on $n$-fold categories \cite{FiorePaoli2010} generalizing Thomason's model structure on categories \cite{Thomason1980}. A recent application of $n$-fold categories to algebraic geometry is found in \cite{Weizhe2017}.

 The definition of an $n$-fold category is elementary, and is based on the notion internal category (see Section \ref{sbs-nint-cat} for more details).
 To understand the latter, remember that the data for a small category can be presented by a diagram
\begin{equation*}
\xymatrix{
\tens{X_1}{X_0} \ar^(0.65){m}[r] & X_1 \ar^{d_0}[r]<2.5ex> \ar^{d_1}[r] & X_0
\ar^{s}[l]<2ex>
}
\end{equation*}
where $X_0$ is the set of objects, $X_1$ the set of arrows, the maps $d_0, d_1$ are the source and target maps, $s$ is the identity map and $m$ is the composition. These maps satisfy the axioms of a category, giving associativity of composition and identity laws.
Such a diagram and axioms make sense in any category $\mathcal{C}$ with pullbacks and this defines the notion of an internal category in $\mathcal{C}$.

 Let $\zD$ be the category of non-empty finite ordinals and morphisms the non-decreasing maps between them. This category is the basis for the combinatorial models of topological spaces called simplicial sets, which are functors from $\zD^{op}$ to the category of sets. These are ubiquitous in algebraic topology.

Functors from the product of $n$ copies of $\dop{}$ (which we denote by $\Dnop$) to $\Set$ are called multi-simplicial sets and are also prominent in algebraic topology, more specifically in simplicial homotopy theory.

The relation between categories and simplicial sets comes from the nerve functor
\begin{equation*}
  N:\Cat \rw \funcat{}{\Set}
\end{equation*}
Given a category $X$, $(NX)[0]$ is the set of objects of $X$ while for $k\geq 1$ $(NX)[k]$ consists of the set of sequences of arrows in $X$ of length $k$. Conversely, a characterization of those simplicial sets that are nerves of categories can be given in terms of the so called Segal condition on a simplicial set.

Similarly, for any category $\mathcal{C}$ with pullbacks there is a nerve functor
\begin{equation*}
  N:\Cat(\mathcal{C}) \rw \funcat{}{\mathcal{C}}
\end{equation*}
whose image can be characterized in terms of Segal maps.

The definition of internal categories can be iterated: if $\mathcal{C}$ is the category $\Cat$, we can consider internal categories in $\Cat$, then repeating this until the $n$-th iteration affords the category $\cat{n}$ of \nfol categories. By iterating the nerve construction one obtains a full and faithful nerve functor
\begin{equation*}
    \Nn :\cat{n}\rw \funcat{n}{\Cat}\;.
\end{equation*}
and a characterization of the image of this functor can be given in terms of (iterated) Segal conditions (see Lemma \ref{lem-multin-iff} for more details).

Thus we can think of \nfol categories as structures whose elements carry an intrinsic $n$-cube geometry
and which may be composed along any one of the $n$ axes of that geometry.

  For instance, when $n=2$, we can visualize a double category as having objects, $1$-morphisms in two different directions (horizontal and vertical), both composable in the respective directions, and squares which can be composed both horizontally and vertically, as in the picture below:

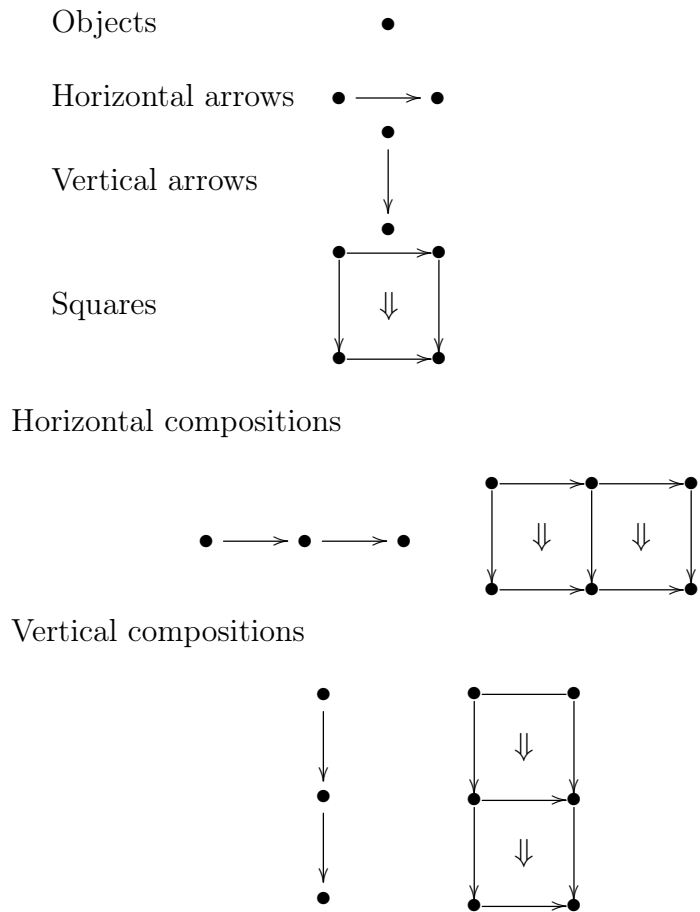
\bk

\pagelabel{double-cat-begin}
\begin{tabular}{ c l c}
    & Objects  & $\bullet$ \\
    & \\
    & Horizontal arrows &
   $
   \xymatrix{
   \bullet  \ar[r] & \bullet
   }
   $  \\
    & Vertical arrows &
   \parbox{10mm}
   {\centering $
   \xymatrix{
   \bullet  \ar[d] \\
   \bullet
   }$}
   \\
    & Squares &
   \parbox{15mm}{\centering
   $
   \entrymodifiers={[]}
   \xymatrix@R12pt @C12pt{
   \bl  \ar[rr] \ar[dd] &  & \bl \ar[dd]\\
   \  & \Downarrow &\ \\
   \bl \ar[rr] & & \bl
   }$} \\
\end{tabular}

\bk

Horizontal compositions
  \begin{equation*}
  \xymatrix@R12pt{
  & &\\
  \bullet \ar[r] & \bullet \ar[r] &\bullet
  }
  \qquad
  \entrymodifiers={[]}
\xymatrix@R12pt @C12pt{
   \bl  \ar@{->}[rr] \ar@{->}[dd] & \ & \bl \ar@{->}[dd]\\
   \  & \Downarrow &\ \\
   \bl \ar@{->}[rr] & & \bl
   }\!\!
\xymatrix@R12pt @C12pt{
   \  \ar@{->}[rr]  & \ & \bl \ar@{->}[dd]\\
   \  & \Downarrow &\ \\
   \ \ar@{->}[rr] & & \bl
   }
\end{equation*}

Vertical compositions
  \begin{equation*}
  \xymatrix@R26pt{
   \bullet  \ar[d] \\
   \bullet \ar[d] \\ \bullet}
   \qquad\qquad
 \entrymodifiers={[]}
\xymatrix@R12pt @C12pt{
   \bl  \ar@{-}[rr] \ar@{->}[dd] & \ & \bl \ar@{->}[dd]\\
   \  & \Downarrow &\ \\
   \bl  \ar@{->}[rr] \ar@{->}[dd] & \ & \bl \ar@{->}[dd]\\
   \  & \Downarrow &\ \\
   \bl \ar@{->}[rr]  & & \bl
   }
\end{equation*}

 These compositions are all associative and unital, hence these structures are, in this sense, 'strict'. However, they are much wider than strict $n$-categories because, unlike in strict $n$-categories, these $n$ directions are completely symmetric and we cannot identify in a \nfol category any 'sets of $k$-cells' for $k=0,\cdots, n$.

\subsection{$\pmb{n}$-Fold structures versus strict and weak $\pmb{n}$-categories}
An important application of $n$-fold structures appeared in the context of homotopy theory, in the modelling of connected $(n+1)$-types using $n$-fold categories internal to the category of groups \cite{Loday1}, as we are going to explain in Section \ref{subs-mod-nfol} below.

This leads to the question of how \nfol categories relate to strict and weak $n$-categories. The first question has an easy answer, and as we will illustrate in detail in Section \ref{sbs-multi-strict} there is a full and faithful embedding
\begin{equation*}
  n\mi\Cat \hookrightarrow \cat{n}
\end{equation*}
For instance, a strict $2$-category is a double category in which the category of objects and vertical arrows is discrete.
Thus the vertical sides of the squares in  Figure \ref{corner2} are identities and thus can be represented as globes under the identification

Consequently the vertical sides of the squares in the picture on page \pageref{double-cat-begin} are identities which we may contract in our
diagrams, thereby depicting those squares as globes:
%
%% Create 1-edge
\tikzset{edge/.pic={
\filldraw (0,0)  circle[radius=0.035cm] -- (1,0)  circle[radius=0.035cm]; %% centered
}}
%% Create  1-cell
\tikzset{Cell1/.pic={
\draw(0,0) ellipse [x radius=0.5cm, y radius=0.30cm];
\doublearrow{arrows={-Implies}}{(0.0,0.1) -- (0.0,-0.1)};
\draw plot [mark=*] coordinates {(-0.5,0)};
\draw plot [mark=*] coordinates {(0.5,0)};
}}
\begin{center}
\begin{tikzpicture}[thick,scale=2]]
\pic [black][scale=2] at (0,0) {edge};
\pic [black][scale=2] at (0,1) {edge};
\draw[thick,double distance=2pt] (0,0.2) -- (0,0.8);
\draw[thick,double distance=2pt] (1,0.2) -- (1,0.8);
\doublearrow{arrows={-Implies}}{(0.5,0.6) -- (0.5,0.4)};
\pic [black][scale=2] at (3,0.5) {Cell1};
\node (1) at (0,1.2) {a};
\node (2) at (1,1.2) {b};
\node (3) at (0,-0.2) {a};
\node (4) at (1,-0.2) {b};
\node (5) at (0.5,1.1) {$g$};
\node (6) at (0.5,-0.12) {$f$};
\node (7) at (1.7,0.5) {$\equiv$};
\node (8) at (2.3,0.5) {a};
\node (9) at (3.7,0.5) {b};
\node (10) at (3,0.9) {$g$};
\node (10) at (3,0.05) {$f$};
\end{tikzpicture}
\end{center}
Thus the picture on page \pageref{double-cat-begin} for a double category reduces to the picture on page \pageref{2cat-begin} for a strict $2$-category.

The second question, of the relationship between $n$-fold categories and weak $n$-categories, is in general
much harder to answer. Indeed much of the work presented here addresses itself directly to answering
the following motivating question of this kind:
\medskip

\emph{Can we identity a suitable subcategory of \nfol categories which gives a model of weak $n$-categories?}
\medskip

We positively answer this question with the introduction of the category of weakly globular $n$-fold categories and the proof of a suitable equivalence to the Tamsamani model of weak $n$-category. Our model is based on a new paradigm to weaken higher categorical structures which is the notion of weak globularity. We explain the idea of weak globularity in Chapter \ref{chap1-2}, while in what follows we give an account of its first appearance in the context of the algebraic modelling of homotopy types.

\section{The homotopy hypothesis}\label{homotyp}
\index{Homotopy hypothesis} As already pointed out, one of the most important connections between homotopy theory and higher category theory is the 'homotopy hypothesis': any good model of weak $n$-category should give an algebraic model of $n$-types in the weak $n$-groupoid case.

In this section we first give a summary of the notion of algebraic models of $n$-types, and we then concentrate on the use of $n$-fold structures to give such models in the path-connected case. This is the first context which saw the development of the notion of weak globularity in the work by the author \cite{Pa}. In the next chapter we will introduce the idea if weak globularity in a more general categorical context.

\subsection{Homotopy types and their algebraic models}\index{Algebraic model of $n$-types}
A topological space whose homotopy groups vanish in dimension higher than $n$ is called an $n$-type.\index{n-type} The $n$-types are building blocks of spaces thanks to a classical construction in algebraic topology, which is the Postnikov decomposition \index{Postnikov decomposition} \cite{Jard}, \cite{MayBook1967}, \cite{Postnikov}, \cite{Spanier1966}. More precisely, the Postnikov system of a space consists of its $n$-type constituents together with some cohomological invariants called $k$-invariants. One of the fundamental theorems in algebraic topology is that the Postnikov system of a space determines its homotopy type \cite{Jard}.

Postnikov systems exist for more complex structures than spaces, for instance for categories enriched in simplicial sets, also called simplicial categories, which are models of $(\infty,1)$-categories. Their Postnikov sections consist of categories enriched in (simplicial) $n$-types, while the $k$-invariants are cohomology classes in the Dwyer-Kan-Smith cohomology \cite{DwyerKanSmith1986}.

A fundamental question in algebraic topology is the search for algebraic models of $n$-types. By this we mean a category $\mathcal{G}_n$, built only from combinatorial and categorical data, together with a pair of functors
\begin{equation*}
  \Pi_n: \text{$n$-types} \rw \mathcal{G}_n   \qquad  B: \mathcal{G}_n\rw \text{$n$-types}
\end{equation*}
inducing an equivalence of categories
  \begin{equation}\label{eq-model-nty}
   \mathcal{G}_n  \bsim^n\;\simeq\;\Ho\text{($n$-types)}\;.
\end{equation}
 Here $\mathcal{G}_n  \bsim^n$ is the localization\index{Localization} of $\mathcal{G}_n$ with respect to some algebraically defined weak equivalences, and $\Ho\text{($n$-types)}$ is the homotopy category of $n$-types. The equivalence of categories \eqref{eq-model-nty} means that the right hand side, which has a purely topological input, is described by the left hand side in an entirely algebraic or categorical way. Hence we call $\mathcal{G}_n$ an 'algebraic model of $n$-types'. \index{Algebraic model of $n$-types}

Since the Postnikov sections of simplicial categories are categories enriched in (simplicial) $n$-types, an algebraic model of the latter via a product-preserving functor $\Pi_n$ gives rise to an algebraic model for the Postnikov sections of simplicial categories, namely consisting of categories enriched in $\mathcal{G}_n$. It is crucial for this application to work with a model of general $n$-types, nor merely path-connected ones, since the mapping spaces of simplicial categories are general simplicial sets, not merely path-connected ones.

The simplest case is $n=1$. In this case, $\mathcal{G}_1$ is the category of groupoids. The functor  $B$ is the classifying space functor, which can be obtained by taking the geometric realization of the nerve of the groupoid; equivalence of groupoids are equivalence of categories. The functor $\Pi_1$ is the classical fundamental groupoid functor \cite{Jard}.\index{Fundamental groupoid} The latter associates to a space the groupoid whose objects are the points of the space and whose morphisms are the homotopy classes of paths.

The extension to the case $n>1$ amounts to finding higher order analogues of the fundamental groupoid. It was noted by Grothendieck in 'Pursuing Stacks'  \cite{Grothen} \index{Pursuing Stacks}that this would naturally lead to some type of higher groupoidal structure. For a given space, there are notions not only of 'homotopy between paths' but also of 'homotopies between homotopies' and then 'homotopies between homotopies between homotopies' and so on. Intuitively we expect the 'fundamental $n$-groupoid' $\Pi_n(X)$ of a space $X$ to have these higher homotopies as higher cells, with only the $n$-cells requiring dividing out by the homotopies one level higher.

Converting this intuition into rigorously defined models is far from trivial and has been the object of much study both by algebraic topologists and category theorists.

The category $n\mi\Gpd$ of strict $n$-groupoids \index{Strict $n$-groupoid} is insufficient for the purpose of modelling $n$-types. For instance, in \cite{S2} there is a counterexample showing that one cannot use strict 3-groupoids to model 3-types of spaces with non-trivial Whitehead products.\index{Whitehead products} Thus a more complex notion of higher groupoid is needed.

 A good notion of weak $n$-category needs to satisfy the homotopy hypothesis: namely, it should provide an algebraic model of  $n$-types when the cells in the structure are weakly invertible, that is in the weak $n$-groupoid case. For several models of higher categories this hypothesis has been shown to hold.

  Some algebraic models of $n$-types arose as the higher groupoidal case of models of weak $n$-categories, such as the Tamsamani \cite{Ta} and the Batanin \cite{B} models. Other algebraic models of $n$-types developed independently within algebraic topology; these models were build algebraically and combinatorially and were shown to satisfy the equivalence of categories \eqref{eq-model-nty} but they were not built as the higher groupoid version of a model of weak $n$-category.  Examples of these include the cat$^{n}$-groups of Loday \cite{Loday1}, the hypercrossed complexes \index{Hypercrossed complexes} of Carrasco and Cegarra \cite{CarrascoCegarra1991}, the crossed $n$-cubes \index{Crossed N-cubes} of Porter \cite{Porter1993} and of Ellis and Steiner \cite{EllStein1987}. In low dimension other examples are the crossed modules in groups introduced by MacLane and Whitehead \cite{MLaneWhite1950}, the double groupoids with connections of Brown and Spencer \cite{BrowSpenc1976}, the quadratic modules of Baues \cite{Baues1991}, and the crossed modules of length 2 of Conduch\'{e} \cite{Conduche1984}.

 Of particular relevance for this work is Loday's model of connected $(n+1)$-types,\index{Connected $n$-types} which is based on \nfol categories internal to the category of groups. This was developed independently in an algebraic topological context and gave rise to interesting topological applications, such as a higher order version of the Van Kampen theorem \cite{BrownLoday1987}. \index{Van Kampen theorem}

 A crucial question, which lies at the very origin of this work, is the following:

 \smallskip
  \emph{Is Loday's model of connected $(n+1)$-types, based on $n$-fold categories internal to groups, the higher groupoid version of some general model of weak higher category?}

  \smallskip

  In order to answer this question one must first understand how the \nfol model of Loday relates to some
model of higher categories that is known to satisfy the homotopy hypothesis, such as the Tamsamani
model.

  This is far from trivial. As we explain in the next section, the key to this  comparison is the notion of weak globularity, introduced in the context of cat$^n$-groups by the author in \cite{Pa}.

\subsection{Modelling homotopy types with $\pmb{n}$-fold structures}\label{subs-mod-nfol}

One of the earliest appearances of higher categorical structures in algebraic topology is due to MacLane and Whitehead \cite{MLaneWhite1950} who used crossed modules in groups \index{Crossed modules} to model connected $2$-types. It is well known that the category of crossed modules in groups is equivalent to the category of internal categories in groups. The latter is the same as the category of strict $2$-groupoids with one object, also called strict $2$-groups. Thus the Maclane-Whitehead model amounts to modelling of connected $2$-types with a strict $2$-dimensional categorical structure.

The modelling of homotopy types becomes easier to handle in the path connected case,
 since it allows us to work within the entirely algebraic context of the category of groups.
  The Kan loop group functor \index{Kan loop group functor} from the category of based topological spaces to simplicial groups \cite{Jard} is the reason that internalisation in groups works in the path-connected case.

As mentioned in the previous section, several models of connected $(n+1)$-types were developed independently of the pursuit of higher categorical models, even though they all exhibit some type of higher structures.

Among these, cat$^{n}$-groups are very appealing from a categorical perspective as their underlying higher categorical structure is particularly simple: it is easy to prove (see for instance \cite{Pa} for details) that they are equivalent to the category $\cat{n}(\Gp)$ of \nfol categories internal to groups.
\index{Cat$^{n}$-groups }

At the same time, models of higher categories were being developed, and in particular the Tamsamani model for which a proof of the homotopy hypothesis is given \cite{Ta}: thus Tamsamani $(n+1)$-groupoids model $(n+1)$-types. A suitable subcategory of Tamsamani $(n+1)$-groupoids can be identified to model connected $(n+1)$-types.

This naturally leads to the question of finding an explicit comparison between the Tamsamani model and the cat$^{n}$-groups model for the path connected case. Such a comparison is highly non-trivial because the two higher categorical structures have a fundamental difference: the Tamsamani model has sets of cells in dimensions $0$ up to $n$, together with compositions coming from a multi-simplicial structure. The higher cells have a globular shape as in the case of strict $n$-groupoids, but the compositions are no longer strictly associative and unital.

 In the cat$^{n}$-groups model, however, there is no immediate way to identify sets of higher cells. The $n$-fold structure is symmetric in all $n$ different simplicial directions, and the only higher cells which we can identify in the structure have a (hyper-)cubical shape. Hence the question:

\medskip
\emph{How can we connect the cubical model of cat$^{n}$-groups to a globular model like the Tamsamani model of connected $(n+1)$-types?}
\medskip

The notion of weak globularity, introduced by the author \cite{Pa} for the category of cat$^n$-groups, is the key to answering this question.

We showed in \cite{Pa} that the category $\cat{n}(\Gp)$ of \nfol categories in groups can be replaced by the category $\cat{n}(\Gp)_{wg}$ of weakly globular cat$^{n}$-groups (which is also called in \cite{Pa} category of 'special cat$^{n}$-groups') without loss of homotopical information. That is, weakly globular cat$^{n}$-groups model connected $(n+1)$-types.

The idea of weak globularity is to impose additional conditions on cat$^{n}$-groups to break the symmetry of the $n$ simplicial directions and identify 'substructures' in a cat$^{n}$-group from which one can recover the sets of higher cells. To identity which substructures in a cat$^n$-group should play this role we consider the full and faithful embedding
  \begin{equation*}
  n\mi\Cat(\Gp) \hookrightarrow \cat{n}(\Gp)
\end{equation*}
of strict $n$-categories internal to groups to \nfol categories internal to groups, which is formally analogous to the embedding of strict $n$-categories in $n$-fold categories mentioned in Section \ref{sec-diff}.

It is not difficult to see (see Section \ref{sbs-multi-strict} for more details) that a strict $n$-category in groups amounts to an $n$-fold category in groups in which certain substructures (which are cat$^{k}$-groups, for $0<k<n-2$) are discrete. By discrete here we mean they are groups, seen as discrete cat$^{k}$-groups in which all structure maps are identities. The corresponding underlying sets are the sets of higher cells in a strict $n$-category internal to groups.
 We call this the \emph{globularity condition} since it is the condition that gives rise to the globular shape of higher cells in the structure.

 In a weakly globular cat$^{n}$-group these substructures are no longer required to be  discrete, but 'homotopically discrete' in a precise sense that allows iteration (these are called in \cite{Pa} 'strongly contractible cat$^{n}$-groups').  Each of these substructures has a higher groupoidal structure of its own which is however not a general $k$-fold structure (for the respective dimension $k$) like in a general internal cat$^n$-group, but is equivalent to a discrete one. The sets underlying the latter correspond to the 'sets of higher cells'.

 Further, we show in \cite{Pa} that there is a comparison functor
\begin{equation*}
 { \cat{n}(\Gp)}_{wg}\rw \gta{n+1}
\end{equation*}
 which preserves the homotopy type. Here $\gta{n+1}$ denotes groupoidal Tamsamani $(n+1)$-groupoids. This functor is obtained by discretizing the homotopically discrete substructures in a weakly globular cat$^n$-group, so as to recover the globularity condition. Full details can be found in \cite{Pa}.

 Besides the comparison with the Tamsamani model, another advantage of weakly globular cat$^{n}$-groups over general cat$^{n}$-groups to model $(n+1)$-types is that the latter do not have an algebraic version of the Postonikov tower functor from cat$^{n}$-groups to cat$^{n-1}$-groups: in other words this functor can be produced only by passing to classifying spaces and then applying the fundamental  cat$^{n-1}$-group functor. Instead, weakly globular cat$^{n}$-groups come equipped with a very simple and entirely algebraically defined functor
 \begin{equation*}
 \p{n}:{ \cat{n}(\Gp)}_{wg}\rw { \cat{n-1}(\Gp)}_{wg}
\end{equation*}
 exactly corresponding to the Postnikov truncation, that is such that for each $X\in { \cat{n}(\Gp)}_{wg}$ (whose classifying space $BX$ is therefore a connected $(n+1)$-type) there is a map $BX \rw B\p{n}X$ inducing isomorphisms of homotopy groups in dimensions $0$ up to $n$.

 In \cite{BP} Blanc and the author developed the notion of weakly globular \nfol groupoids with the goal of using \nfol structures for modelling general $n$-types. This case is considerably more complex than the path-connected case. In particular, several features of weakly globular cat$^n$-groups that could be deduced from their definition now need to become part of the definition of weakly globular $n$-fold groupoid.

 In \cite{BP} we build a functor from spaces to weakly globular \nfol groupoids as well as from  weakly globular \nfol groupoids to $n$-types, but we do not exhibit a proof that this gives an equivalence of categories after localization, that is that weakly globular \nfol groupoids are an algebraic model of $n$-types.

  More precisely, in \cite{BP} we show that the functor from the localization of weakly globular \nfol groupoids to the homotopy category of $n$-types is essentially surjective on objects, but we do not give a proof of the fully faithfulness needed to realize an equivalence of categories. We do realize an equivalence of categories in \cite{BP} by enlarging the category of weakly globular \nfol groupoids, but this comes at the price that the larger category is no longer an \nfol structure.

  The realization of the modelling of general $n$-types via weakly globular $n$-fold structures remained open until the present work, where it is one of our main results. This needs the category of \emph{groupoidal weakly globular \nfol categories}, which strictly contains the weakly globular \nfol groupoids of \cite{BP}. Further, we show that the functor from spaces to the latter given in \cite{BP} can still be used as a convenient model: such a model is very explicit and it is independent of other models of $n$-types.

%%%%%%%%%%%%%%%%%%%%%%%%%%%%%%%%%%%%%%%%%%%%%%%%%%%%%%%%%%%%%%%%%%%%%%%%%%%%%%%%%%%%%

\chapter{An introduction to the three Segal-type models}\label{chap1-2}

In this chapter we give an overview of the approach to higher categories developed in this book, which is based on \nfol structures and on the idea of weak globularity. We develop this in the context of three models based on multi-simplicial structures, which we call Segal-type models since the compositions of higher cells is related to the notion of Segal maps. One of these models is due to Tamsamani \cite{Ta}, the other two are new. Our aim in this chapter is to convey the main intuitions and ideas, referring the reader to the later chapters for the precise definitions and results.

In Section \ref{sec-geo-vs} we highlight one of the issues that we face when dealing with the general higher categorical case, instead of the higher groupoidal one, which is the notion of higher categorical equivalence. In Section \ref{sbs-multisimpl-env} we explain why multi-simplicial structures are a natural environment for the development of models of higher categories, referring to the crucial notion of Segal maps.

In Section \ref{sbs-idea-of} we explain the intuition behind the idea of weak globularity, which is central to this work in developing the new paradigm to weaker higher categorical structures. Finally in Section \ref{three-seg} we give a broad overview of the main features of the three Segal-type models of higher categories treated in this work and of the main result. We also give an overview of the organization of this work.

In this chapter we sometimes needs notation and definitions recalled in the background part, and we refer the reader to Part \ref{part-back} for further details.

\section{Geometric versus higher categorical equivalences}\label{sec-geo-vs}\index{Geometric weak equivalences}

An important aspect of the passage from the higher groupoidal case to the general higher categorical case in modelling weak higher structures concerns the notion of equivalence. It is well known that a functor between groupoids is an equivalence of categories (that is, fully faithful and essentially surjective on objects) if and only if it induces a weak homotopy equivalence of classifying spaces; here the classifying space of a groupoid is the geometric realization of its simplicial nerve, and a weak homotopy equivalence is a map inducing isomorphisms of homotopy groups in all dimensions (see \cite{Jard}). So in the category of groupoids, categorical and geometric weak equivalences of coincide.

In fact, let $F: X \rw Y$ be a functor between groupoids which is a weak homotopy equivalence. Then $\pi_0(BX)= q(X)$, where $q$ denotes the connected components functor, and similarly for $Y$. Thus $\pi_0(BF)$ being an isomorphism implies that $F$ is essentially surjective on objects. Given $a,b \in X$, since $X$ and $Y$ are groupoids there are bijections
\begin{equation*}
  X(a,b)\cong X(a,a), \qquad Y(Fa,Fb)\cong Y(Fa,Fa)
\end{equation*}
On the other hand, $\pi_1(BX, a)=X(a,a)$ and $\pi_1(BY, Fa)=Y(Fa,Fa)$. This $\pi_1 (BF, a)$ being an isomorphism implies from above that there is a bijection
\begin{equation*}
  X(a,a)\cong Y(Fa, Fa)
\end{equation*}
that is, $F$ is fully faithful. In conclusion, $F$ is an equivalence of categories.

Conversely, if $F:X \rw Y$ is an equivalence of categories, then $q(F)$ is an isomorphism, hence from above $\pi_0 B(F)$ is an isomorphism. Also, since $F$ is fully faithful, for each $a\in X$, $\pi_1(BX, a)=X(a,a)\cong Y(Fa, Fa)=\pi_1(BY, Fa)$. That is, $\pi_1(B(F),a)$ is an isomorphism. In conclusion, $B(F)$ is a weak homotopy equivalence, since it induces isomorphisms of all homotopy groups (recall that $\pi_i(BX,a)=0$ for each $i>1$ since $X$ is a groupoid).

This relation between categorical and weak homotopy equivalences for groupoidal structures extends to higher dimensions: for instance for weakly globular cat$^n$-groups in \cite{Pa} and for weakly globular \nfol groupoids in \cite{BP}, weak equivalences are defined using classifying spaces, the latter being the geometric realization of a multi-diagonal, and it is shown that these can also be described in terms of higher categorical equivalences.

The general categorical case is more complex: an equivalence of categories is also a weak homotopy equivalence of its simplicial nerves, but not conversely. Thus a notion of higher categorical equivalence needs to be defined alongside with notions of weak $n$-categories. In the case of our Segal-type models, a notion of $n$-equivalence is given as part of the inductive definition of the structure.

\section{Multi-simplicial structures as an environment for higher categories}\label{sbs-multisimpl-env}

Our Segal-type models of weak $n$-categories are based on multi-simplicial objects, more precisely functors from $\Dnop$ to $\Cat$, where $\Dnop$ denotes the product of $n$ copies of $\Dop$. The reason why multi-simplicial objects are a good environment for building models of higher categories is the fact that their combinatorics exhibits certain maps, called Segal maps, as natural candidates for the compositions of higher cells.

Let first illustrate the case $n=1$. As already recalled in Section \ref{sec-diff}, there is a nerve functor
\begin{equation*}
  N:\Cat\rw\funcat{}{\Set}
\end{equation*}
where
\begin{equation*}
    (N X)_k=
    \left\{
      \begin{array}{ll}
        X_0, & \hbox{$k=0$;} \\
        X_1, & \hbox{$k=1$;} \\
        \pro{X_1}{X_0}{k}, & \hbox{$k>1$.}
      \end{array}
    \right.
\end{equation*}
It is possible to give a characterization of the simplicial sets that are nerves of small categories by using the notion of Segal maps. For each $k\geq 2$ these are maps
\begin{equation*}
  \mu_k :X_k\rw\pro{X_1}{X_0}{k}
\end{equation*}
which arise from the commuting diagram

 
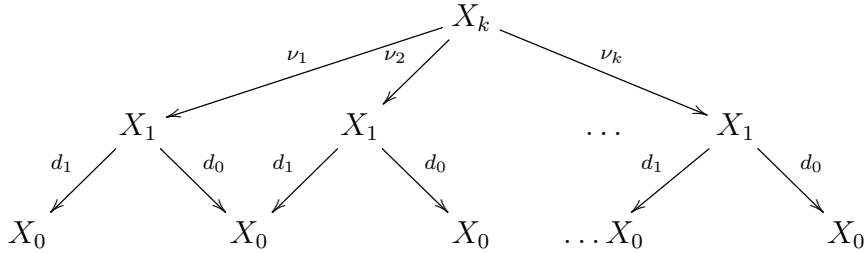
\begin{equation*}
\xymatrix@C=20pt{
&&&& X\sb{k} \ar[llld]_{\nu\sb{1}} \ar[ld]_{\nu\sb{2}} \ar[rrd]^{\nu\sb{k}} &&&& \\
& X\sb{1} \ar[ld]_{d\sb{1}} \ar[rd]^{d\sb{0}} &&
X\sb{1} \ar[ld]_{d\sb{1}} \ar[rd]^{d\sb{0}} && \dotsc &
X\sb{1} \ar[ld]_{d\sb{1}} \ar[rd]^{d\sb{0}} & \\
X\sb{0} && X\sb{0} && X\sb{0} &\dotsc X\sb{0} && X\sb{0}
}
\end{equation*}
where for each ${1\leq j\leq k}$ and $k\geq 2$, let ${\nu_j:X_k\rw X_1}$ be induced by the map  $[1]\rw[k]$ in $\Delta$ sending $0$ to ${j-1}$ and $1$ to $j$.

 Then a simplicial set is the nerve of a small category if and only if it satisfies the \emph{Segal condition} \index{Segal condition} that all the Segal maps are isomorphisms. Under this isomorphisms the composition map
\begin{equation*}
  \tens{X_1}{X_0}\xrw{c} X_1
\end{equation*}
corresponds to the face map $X_2\rw X_1$ induced by the map $[1]\rw [2]$ in $\Delta$ sending 0 to 1 and 1 to 2.

 Thus given a simplicial set, its simplicial structure together with the Segal maps contain a natural candidate for the composition, and the Segal condition ensures that such composition satisfies the axioms of a categorical composition (that is, it is associative and unital).
This paradigm carries on in higher dimensions and it can be used to give notions of strict as well as weak higher categories.

Namely, Segal maps can be defined for any simplicial object $\funcat{}{\clC}$ in a category $\clC$ with finite limits, and they can be used to characterize nerves of internal categories in $\clC$ using the Segal condition that all Segal maps are isomorphisms, in a way formally analogous to the case of $\Cat$ (see Chapter \ref{prelim} for more details). For instance internal categories in $\Cat$ (also called double categories) can be described as simplicial objects in $\Cat$ such that all the Segal maps are isomorphisms.

Strict 2-categories can be described as simplicial objects $X$ in $\Cat$ satisfying the Segal condition as well as the condition that $X_0$ is a discrete category. Once again, the simplicial structure and the Segal maps give the candidates for the compositions of higher cells while the Segal condition implies the associativity and unitality of the compositions. We refer to Examples \ref{ex-mult-ner-dim2} and \ref{ex-1-multi-strict} in the next chapter for a pictorial representation of double categories and strict 2-categories.

Another approach to the use of Segal maps is possible: given a simplicial object $X\in\funcat{}{\Cat}$ such that $X_0$ is a discrete category, we can relax the Segal condition and require that Segal maps for each $k\geq 2$
\begin{equation*}
  \mu_k :X_k\rw\pro{X_1}{X_0}{k}
\end{equation*}
are not isomorphisms but merely equivalences of categories. We obtain a structure with sets of 0-cells $X_{00}$, sets of 1-cells $X_{10}$, sets of 2-cells $X_{11}$. We can also define a composition
\begin{equation*}
  \tens{X_1}{X_0}\xrw{\mu'_{2}} X_2 \rw X_1
\end{equation*}
where $\mu'_{2}$ is a pseudo-inverse to $\mu_2$ and the map $X_2\rw X_1$ is induced by the map $[1]\rw[2]$ in $\Delta$ sending 0 to 1 and 1 to 2. Since $\mu'_{2}$ is no longer an isomorphism, this composition is no longer associative and unital. What we obtain is a Tamsamani 2-category \cite{Ta}. Lack and the author proved in \cite{LackPaoli2008} that Tamsamani 2-categories are suitably equivalent to bicategories.

In dimensions higher than 2, we use the category $\funcat{n-1}{\Cat}$ of multi-simplicial objects in $\Cat$ to obtain suitable notions of $n$-categories. Segal maps can be defined for multi-simplicial objects, and it is easy to see that there are fully faithful multi-nerve functors
\begin{align*}
    & \cat{n} \rw \funcat{n-1}{\Cat} \\
    & n\mi\Cat \rw \funcat{n-1}{\Cat}
\end{align*}
(see Sections \ref{sus-ner-funct} and \ref{sbs-multi-strict} for more details). A characterization of the image of these functors is also easy to give in terms of Segal maps (see Lemma \ref{lem-multin-iff} and Corollary \ref{cor-multi-strict}). However, as in the case $n=2$, we can impose different kind of Segal conditions to obtain a different behavior of the composition of higher cells and therefore define models of weak higher categories.

In summary, multi-simplicial objects are a good environment for building models of higher categories since the multi-simplicial maps and the corresponding Segal maps provide suitable candidates for the compositions of higher cells.

To obtain models of weak $n$-categories using multi-simplicial structures we need to impose extra conditions to encode:
\begin{itemize}
  \item [a)] The sets of cells in dimensions $0,\ldots,n$.

  \item [b)] The behavior of the compositions, giving weak associative and unit laws.

  \item [c)] A notion of higher categorical equivalence.
\end{itemize}
As outlined in Section \ref{three-seg}, we adopt three different approaches that encode this extra structure, giving rise to three different Segal-type models of weak $n$-categories.

\section{The idea of weak globularity}\label{sbs-idea-of}
Weakly globular \nfol categories form a full subcategory of \nfol categories. They are based on a new paradigm to weaken higher categorical structures: the idea of \emph{weak globularity}.

 In a strict $n$-category, the $k$-cells, for each for each $0\leq k\leq n$, form a set. Equivalently, considering strict $n$-categories as embedded in  $\funcat{n-1}{\Cat}$ via the multinerve functor, the $k$-cells for each $0\leq k\leq (n-2)$ form a discrete $(n-k-1)$-fold category, that is one in which all structure maps are identities.

  In the weakly globular approach the $k$ cells (for each $0\leq k\leq (n-2)$) no longer form a set but have a higher categorical structure on their own. More precisely, they form a 'homotopically discrete $(n-k-1)$-fold category': this is a  $(n-k-1)$-fold category which is suitably equivalent to a discrete one (see Chapter \ref{chap3} for more details). We call this the \emph{weak globularity condition}.
\index{Weak globularity condition}

 The weakness in a weakly globular \nfol category is encoded by the weak globularity condition. Further, the weak globularity condition allows to recover the notion of 'sets of higher cells' in a $n$-fold category. More precisely, in a weakly globular \nfol category, for each $0\leq k\leq (n-2)$ there are substructures which are homotopically discrete $(n-k-1)$-fold categories, equivalent to discrete structures: the underlying sets of these discrete structures correspond to the sets of $k$-cells. The definition of weakly globular \nfol category also requires several additional conditions to obtain well behaved compositions of higher cells.

 The main result of this work is that there is a suitable equivalence between weakly globular $n$-fold categories and Tamsamani $n$-categories. To establish this result, we work in a broader context of three Segal-types models: weakly globular $n$-fold categories, Tamsamani $n$-categories, and a further new model which we call weakly globular Tamsamani $n$-categories, containing the previous two as special cases.

 Below we give a summary account of the main features of these three models.

 \section{The three Segal-type models}\label{three-seg}
We identify three multi-simplicial models based on the notion of Segal maps, which we therefore call Segal-type models. The first is the category $\Tan$ of Tamsamani $n$-categories introduced by Tamsamani \cite{Ta} and further studied by Simpson \cite{Simp}. The second is the category $\catwg{n}$ of weakly globular \nfol categories, introduced first in this work for any $n\geq 1$, which is a full subcategory of the category $\cat{n}$ of \nfol categories. The third is another new model, the category $\tawg{n}$ of weakly globular Tamsamani $n$-categories.\index{Strict $n$-category}

There is a third higher categorical structure which embeds in all three of the above, which is the category $n\mi\Cat$ of strict $n$-categories. There are full and faithful inclusions
\begin{equation*}
\xymatrix@R15pt @C30pt{
& \tawg{n} & \\
\Tan \ar@{^{(}->}[ur] & & \catwg{n} \ar@{_{(}->}[ul]\\
& n\mi\Cat \ar@{_{(}->}[ul]\ar@{^{(}->}[ur] &
}
\end{equation*}
The category $n\mi\Cat$ admits a multi-simplicial description (see Section \ref{sbs-multi-strict} for more details) as the full subcategory of $\nm$-fold simplicial objects $X\in \funcat{n-1}{\Cat}$ satisfying the following
\begin{itemize}
  \item [(i)] $X_0 \in \funcat{n-2}{\Cat}$ and $X_{\uset{r}{1...1}0} \in \funcat{n-r-2}{\Cat}$ are constant multi-simplicial objects taking value in a discrete category, for all $1\leq r < n-2$. Here we use  Notation \ref{not-simp}.
      \smallskip
  \item [(ii)] The Segal maps (see Definition \ref{def-seg-map}) in all directions are isomorphisms.
\end{itemize}
The underlying sets of the discrete structures $X_0$ (resp. $X_{\uset{r}{1...1}0}$) in (i) correspond to the sets of 0-cells (resp. $r$-cells) for $1\leq r \leq n-2$; the sets of $\nm$ and of $n$-cells are given by $ob(X_{\uset{n-1}{1...1}})$ and $mor(X_{\uset{n-1}{1...1}})$ respectively.

The isomorphisms of the Segal maps (condition (ii)) ensures that the composition of cells is associative and unital.

The discreteness condition (i) is also called the \emph{globularity condition}.\index{Globularity condition}  The name comes from the fact that it determines the globular shape of the cells in a strict $n$-category. For instance, when $n=2$, we can picture 2-cells as globes
\begin{equation*}
\xymatrix{
\bullet \ar@/^2pc/^f[rr] \ar@/_2pc/_g[rr] &  \Downarrow \xi & \bullet
}
\end{equation*}
See Examples \ref{ex-1-multi-strict} and \ref{ex-2-multi-strict} for more detailed descriptions of the cases $n=2$ and $n=3$.

Strict $n$-categories found several applications, for example in the groupoidal case where they are equivalent to crossed $n$-complexes (see \cite{BHS}). \index{Crossed $n$-complexes} However, they do not satisfy the homotopy hypothesis (see \cite{S2} for a counterexample showing that strict 3-groupoids do not model 3-types).

Therefore we must relax the structure to obtain a model of weak $n$-category. Using the multi-simplicial framework, we consider three approaches to this:
\begin{itemize}
  \item [a)] \emph{Category $\Tan$of Tamsamani $n$-categories}. In the first approach, we preserve the globularity condition (i) and we relax the Segal map condition (ii) by allowing the Segal maps to be suitably defined higher categorical equivalences. This makes the composition of cells no longer strictly associative and unital.
\bigskip
  \item [b)] \emph{Category $\catwg{n}$ of weakly globular \nfol categories}. In the second approach condition (ii) is preserved so we obtain a subcategory of \nfol categories. However, the globularity condition (i) is replaced by weak globularity: the sub-structures $X_0$, $X_{\uset{r}{1...1}0}$ ($1\leq r\leq n-2$) are no longer discrete but 'homotopically discrete' in a higher categorical sense that allows iterations (see Chapter \ref{chap3} for more details).

      The notion of homotopically discrete \nfol category is a higher order version of equivalence relations. In particular, if an \nfol category $X$ is homotopically discrete, it is suitably equivalent to a discrete \nfol category $X^d$ via a 'discretization map' $\zg : X\rw X^d$.
\bigskip
  \item [c)] \emph{Category $\tawg{n}$ of weakly globular Tamsamani $n$-categories}. In the third approach, both conditions (i) and (ii) are relaxed.
\end{itemize}
\subsection{Common features of the three Segal-type models}\label{sec-comm-fea}\index{Segal-type model} In this section we describe the main common features of the three models, which we denote collectively by $\seg{n}$. Our purpose is to convey some of the main ideas underpinning the three Segal-type models. We point out, however, that to prove our results, the main structures need to be developed in the order presented in this work: namely we first need to define the categories $\cathd{n}$ and $\catwg{n}$ and only later $\tawg{n}$ and $\ta{n}$.
\medskip
\begin{itemize}
  \item [(1)] \emph{Inductive definition} \index{Inductive definition}

   \nid $\seg{n}$ is defined inductively on dimension, starting with $\seg{1}=\Cat$. For each $n>1$ we set
  \begin{equation*}
    \seg{n}\subset \funcat{}{\seg{n-1}}.
  \end{equation*}
    Unravelling this definition gives an embedding \index{Functor!- $J_n$}
      \begin{equation*}
        J_{n}:\seg{n}\rw\funcat{n-1}{\Cat}\;
      \end{equation*}
      Thus our models have a multi-simplicial structure. The simplicial maps are the candidates for encoding the compositions of higher cells.
An object $X$ of $\seg{n}$ is called \emph{discrete} if $J_{n}X$ is a constant functor taking value in a discrete category.
The functor $J_n$ is defined in Notation \ref{not-ner-func-dirk} for $n$-fold categories (this also for $\catwg{n}$), and in  Definition \ref{def-wg-ps-cat} for the category $\tawg{n}$ (and thus also for $\ta{n}$).
\medskip
  \item [(2)] \emph{Weak globularity condition}\index{Weak globularity condition}

   \nid This condition encodes the sets of higher cells and is based on our new paradigm to weaken higher categorical structures. Namely, if $X\in\seg{n}$, then $X_0$ is a homotopically discrete $\nm$-fold category, and it is discrete if $X\in\Tan$.

   We refer to Chapter \ref{chap3} for more details on the notion of homotopically discrete $\nm$-fold category. If $X$ is homotopically discrete, there is a discretization map
    \begin{equation}\label{eq-intro-dis-ma}
      \zg:X\rw X^d
    \end{equation}
    where $X^d$ is discrete (as defined above). Thus, if $X\in\seg{n}$, $X_0$ and $X_{\uset{r}{1\ldots 1}0}$ (for $1 \leq r\leq n-2$) are homotopically discrete. The sets underlying the discrete structures $X_0^d$, $X^d_{\uset{r}{1\ldots 1}0}$ correspond to the sets of $r$-cells for $0\leq r\leq n-2$.

    The weak globularity condition is used in Definition \ref{def-n-equiv} (for the category $\catwg{n}$) and in Definition \ref{def-wg-ps-cat} (for the category \tawg{n}).
\medskip
  \item [(3)] \emph{Truncation functor $\p{n}$} \index{Truncation functor}

   \nid Given $X\in\seg{n}$, we can apply levelwise the isomorphism classes of objects functor $p:\Cat\rw\Set$ to $J_n X\in \funcat{n-1}{\Cat}$ to obtain $\ovl{p}J_n X\in\funcat{n-1}{\Set}$. We require that this is the multinerve of an object of $\seg{n-1}$; that is, there is a truncation functor
      \begin{equation*}
        \p{n}:\seg{n}\rw\seg{n-1}
      \end{equation*}
      making the following diagram commute:
      \begin{equation*}
      \xymatrix{
      \seg{n} \ar@{^{(}->}^{J_n}[rr]  \ar_{\p{n}}[d] && \funcat{n-1}{\Cat} \ar^{\ovl{p}}[d] \\
      \seg{n-1} \ar@{^{(}->}_{\Nb{n-1}}[rr]  && \funcat{n-2}{\Set}
      }
      \end{equation*}
    The effect of the functor $\p{n}$ is of dividing out by the highest dimensional invertible cells.

    \nid The definition of truncation functor is found in Definition \ref{def-hom-dis-ncat} (for the category $\cathd{n}$), in Definition \ref{def-n-equiv} (for the category $\catwg{n}$), in Definition \ref{def-wg-ps-cat}(for the category $\tawg{n}$) and in Example \ref{ex-tam} (for the category $\ta{n}$).
   \medskip
  \item [(4)] \emph{Higher categorical equivalences} \index{n-equivalences}

      \nid Given $X\in\seg{n}$ and $(a,b)\in X_0^d$, let $X(a,b)\subset X_1$ be the fiber of the map
      \begin{equation*}
        X_1 \xrw{(d_0,d_1)} X_0\times X_0 \xrw{\zg\times\zg} X_0^d\times X_0^d\;.
      \end{equation*}
      where $\zg$ is the discretization map as in \eqref{eq-intro-dis-ma}. Each $X(a,b)\in \seg{n-1}$ should be thought of as a hom-$(n-1)$-category. \index{Hom-$(n-1)$-category}

      The 1-equivalences in $\seg{1}$ are equivalences of categories. Inductively, if we defined $\nm$-equivalences in $\seg{n-1}$, we say that a map $f:X\rw Y$ in $\seg{n}$ is a $n$-equivalence if the following conditions hold

      \begin{itemize}
        \item [i)] For all $a,b \in X_0^d$,
      \begin{equation*}
        f(a,b):X(a,b)\rw Y(fa,fb)
      \end{equation*}
      are $\nm$-equivalences

        \item [ii)]$\p{n}f$ is an $\nm$-equivalence.
      \end{itemize}

\smallskip
        \nid This definition is a higher dimensional generalization of a functor which is an equivalence of categories. Indeed the latter can be formulated by saying that a functor $F:X \rw Y$ is such that $X(a,b)\cong Y(Fa, Fb)$, for all $a,b \in X$ and $p(F)$ is an isomorphism.

         \nid Condition i) is a higher dimensional generalization of fully faithfullness while condition ii) generalizes essential surjectivity on objects. As in the case of an equivalence of categories, condition ii) can be weakened.

        \nid The notion of $n$-equivalence is given in Definition \ref{def-hom-dis-ncat-3} (for the category $\cathd{n}$) and it is part of the inductive Definitions  \ref{def-n-equiv}, \ref{def-wg-ps-cat} of $\catwg{n}$ and $\tawg{n}$ respectively.

\medskip
  \item [(5)] \emph{Induced Segal maps condition} \index{Induced Segal maps condition}

   \nid This condition regulates the behaviour of the compositions. Given $X\in \seg{n}$, since $X\in\funcat{}{\seg{n-1}}$ and there is a map $\zg:X_0\rw X_0^d$, we can consider the induced Segal maps for $k\geq 2$ (see Definition \ref{def-ind-seg-map} for more details)
       \begin{equation*}
         \hmu{k}:X_k \rw \pro{X_1}{X_0^d}{k}\;.
       \end{equation*}
       In defining $\seg{n}$ we require these maps to be $\nm$-equivalences. Note that when $X\in\Tan$, $\zg=\Id$, so $\hmu{k}$ are the same as the Segal maps.

       \nid The induced Segal maps condition is shown to hold for the category $\cathd{n}$ in Proposition \ref{pro-ind-map-hom-disc}, while it is part of the inductive Definitions \ref{def-n-equiv}, \ref{def-wg-ps-cat} of $\catwg{n}$ and $\tawg{n}$ respectively.

\medskip
  \item [(6)] \emph{The functor $\q{n}$}

  \nid Given $X\in\seg{n}$, we can apply levelwise the isomorphism classes of objects functor $q:\Cat\rw\Set$ to $J_n X\in \funcat{n-1}{\Cat}$ to obtain $\ovl{q}J_n X\in\funcat{n-1}{\Set}$. We can prove that this is the multinerve of an object of $\seg{n-1}$; that is, there is a functor
      \begin{equation*}
        \q{n}:\seg{n}\rw\seg{n-1}
      \end{equation*}
      making the following diagram commute;
      \begin{equation*}
      \xymatrix{
      \seg{n} \ar@{^{(}->}^{J_n}[rr]  \ar_{\q{n}}[d] && \funcat{n-1}{\Cat} \ar^{\ovl{q}}[d] \\
      \seg{n-1} \ar@{^{(}->}_{\Nb{n-1}}[rr]  && \funcat{n-2}{\Set}
      }
      \end{equation*}

  \nid There is also a functor

  \begin{equation*}
    \di{n}:\seg{n-1}\rw \seg{n}
  \end{equation*}
making the following diagram commute

\begin{equation*}
      \xymatrix{
      \seg{n-1} \ar@{^{(}->}^{J_n}[rr]  \ar_{\di{n}}[d] && \funcat{n-2}{\Cat} \ar^{\ovl{d}}[d] \\
      \seg{n} \ar@{^{(}->}_{\Nb{n}}[rr]  && \funcat{n-1}{\Set}
      }
      \end{equation*}
  where $d:\Set \rw \Cat$ is the discrete category functor. The effect of the functor $\di{n}$ is to include objects of $\seg{n-1}$ in $\seg{n}$ as discrete in the $n$th direction. Since $q$ is left adjoint to $d$ there is a map, natural in $X\in \seg{n-1}$

  \begin{equation*}
    X \rw \di{n} \q{n} X.
  \end{equation*}

\nid This map plays an important role in the construction of the rigidification functor $Q_n$.

\nid The existence of the functor $\q{n}$ is proved in Proposition \ref{pro-post-trunc-fun} (for the category $\tawg{n}$) and in Corollary \ref{pro-post-wg-ncat} (for the categories $\cathd{n}$, $\catwg{n}$, $\ta{n}$). The functor $\di{n}$ is discussed in Definition \ref{dn} for internal $n$-fold categories, in Lemma \ref{lem-wg-ps-cat-b} for $\tawg{n}$ and in Remark \ref{rem-wg-ps-cat-b} for $\catwg{n}$ and $\ta{n}$.
\medskip
\item [(7)] \emph{The functors $\p{j,n}$, $\q{j,n}$ and $\di{n,j}$} \index{Functor!- $\p{j,n}$}
 \index{Functor!- $\q{j,n}$}\index{Functor!- $\di{n,j}$}

\nid By repeatedly applying the functors $\p{n}$ and $\q{n}$ we obtain functors, for each $1 \leq j < n$
\begin{align*}
        \p{j,n} & = \p{j}\p{j+1}\cdots \p{n-1}\p{n}:\seg{n}\rw \seg{j-1}\\
        \q{j,n}& = \q{j}\q{j+1}\cdots\q{n-1}\q{n}:\seg{n}\rw \seg{j-1}\;.
    \end{align*}
These are often used in proving the results of this work. The definitions for the different models are found in Definition \ref{def-pn} and Notation \ref{not-ex-tam} (for $\p{j,n}$) and Notation \ref{not-ex-tam-q}, Remark \ref{rem-qrn} (for $\q{j,n}$).

\nid We can also repeatedly apply the functor $\di{n}$ and obtain functors
\begin{equation*}
  \di{n,j}=\di{n}...\di{j+1}\di{j}: \seg{j-1}\rw \seg{n}
\end{equation*}
see Notation \ref{dnot-ex-tam}.

\medskip
\item [(8)] \emph{Groupoidal Segal-type models}

\nid The definition of the groupoidal version $\gseg{n}$ of the three Segal-type models is obtained inductively starting from groupoids as follows. We set $\gseg{1}=\Gpd$. Suppose inductively we defined $\gseg{n-1}\subset\seg{n-1}$, then $X\in\gseg{n}\subset \seg{n}$ if
\begin{itemize}
  \item [i)] $X_k\in\gseg{n-1}$ for all $k\geq 0$.

  \item [ii)] $\p{n}X\in\gseg{n-1}$.
\end{itemize}
We show in Corollary \ref{cor-three-models} that $\gseg{n}$ is a model of $n$-types, that is there is an equivalence of categories
\begin{equation*}
  \gseg{n}\bsim^n\;\simeq\;\Ho(\nty)\;.
\end{equation*}
Thus $\seg{n}$ is a model of weak $n$-categories, satisfying the homotopy hypothesis.

\end{itemize}
\subsection{Main Results}\label{subs-main-res}
The central result of this work is a model comparison between weakly globular \nfol categories and Tamsamani $n$-categories showing that the two models are suitably equivalent after localization. Our main theorems are as follows
(see Theorem \ref{the-funct-Qn}, Theorem \ref{the-disc-func}, Theorem \ref{cor-the-disc-func} and Theorem \ref{cor-gta-2}).

\begin{theorem*}\rm{A}.
There is a functor \emph{rigidification}
\begin{equation*}
\Qn: \tawg{n}\rw\catwg{n}
\end{equation*}
and for each $X\in\tawg{n}$ an $n$-equivalence natural in $X$
\begin{equation*}
  s_n(X):\Qn X\rw X.
\end{equation*}
\end{theorem*}

\begin{theorem*}\rm{B}.
There is a functor \emph{discretization}
\begin{equation*}
  \Discn:\catwg{n}\rw \ta{n}
\end{equation*}
and, for each $X\in\catwg{n}$, a zig-zag of $n$-equivalences in $\tawg{n}$ between $X$ and $\Discn X$.
\end{theorem*}

\begin{theorem*}\rm{C}.
The functors
    \begin{equation*}
        \Qn:\ta{n}\leftrightarrows \catwg{n}:\Discn
    \end{equation*}
    induce an equivalence of categories after localization with respect to the $n$-equivalences
    \begin{equation*}
        \ta{n}\bsim^n\;\simeq \; \catwg{n}\bsim^n
    \end{equation*}
\end{theorem*}

We also identify a subcategory $\gcatwg{n}\subset \catwg{n}$ of groupoidal weakly globular $n$-fold categories and we show that it gives an algebraic model of $n$-types. This means that the category $\catwg{n}$ satisfies the homotopy hypothesis:

\begin{theorem*}\rm{D}.
  There is an equivalence of categories
  \begin{equation*}
    \gcatwg{n}\bsim^n\;\simeq \; \Ho\mbox{($n$-types)}
  \end{equation*}
\end{theorem*}
In Corollary \ref{cor1-mod-fund-wg-group} we also show that the equivalence of categories of Theorem D can be realized using the fundamental weakly globular \nfol groupoid functor of Blanc and the author \cite{BP}.
\bk

We call the functor $Q_n$ \emph{rigidification functor} because it replaces a globular and weak structure with an equivalent more rigid structure (an \nfol category), which is no longer globular but it is weakly globular. For $n=2$ the functor $Q_2$ was constructed by Pronk and the author in \cite{PP}. The construction of $Q_n$ and $\Discn$ for $n>2$ is much more complex and it requires several novel ideas and techniques, as presented in this work.

One of the key features in the construction of $Q_n$ is the use of pseudo-functors. Pseudo-functors feature prominently in homotopy theory, for instance in iterated loop space theory \cite{Thomas}. They are also ubiquitous in category theory \cite{Borc}, and can be described with the language of $2$-monad and their pseudo-algebras \cite{PW}.

  Weakly globular $n$-fold categories are a full subcategory of $(n-1)$-fold simplicial objects in $\Cat$, that is functors $\funcat{n-1}{\Cat}$. We consider the pseudo-version of these, that is pseudo-functors $\psc{n-1}{\Cat}$.

  Crucial to this work is the use of the strictification of pseudo-functors into strict functors

  \begin{equation*}
  \St : \psc{n-1}{\Cat} \rw \funcat{n-1}{\Cat}
\end{equation*}
left adjoint to the inclusion. This topic had many contributions in category theory, including \cite{Str}. We use in this work the elegant formulation of Power \cite{PW}, further refined by Lack in \cite{Lack}.

We introduce a subcategory
\begin{equation*}
  \segpsc{n-1}{\Cat}\subset \psc{n-1}{\Cat}
\end{equation*}
of \emph{Segalic pseudo-functors}. We show in Theorem \ref{the-strict-funct} is that the  strictification functor $\St$ restricts to a functor
\begin{equation*}
    \St:\segpsc{n-1}{\Cat} \rw\catwg{n} \subset \funcat{n-1}{\Cat}\;.
\end{equation*}

The rigidification functor factors through the subcategory of Segalic pseudo-functors. That is, $Q_n$ is a composite
\begin{equation*}
    Q_n:\tawg{n}\xrw{\ \;\;} \segpsc{n-1}{\Cat}\xrw{\St} \catwg{n}\;.
\end{equation*}

In the case $n=2$, it is easy to build pseudo-functors from $\tawg{2}$. More precisely, given $X\in\tawg{2}$, define $\tr{2}X\in[ob(\Dop),\Cat]$ by
\begin{equation}\label{eq1-the-XXXX}
(\tr{2}X)_k=
\left\{
  \begin{array}{ll}
    X_0^d & k=0 \\
    X_1 & k=1 \\
    \pro{X_1}{X_0}{k} & k>1\;.
  \end{array}
\right.
\end{equation}
Since $X\in\tawg{2}$, $X_0\in\cathd{}$ so there are equivalences of categories
\begin{equation*}
\begin{split}
    & X_0\simeq X_0^d \\
    & X_k\simeq \pro{X_1}{X_0^d}{k}\quad \text{for }\; k>1.
\end{split}
\end{equation*}
Thus, for all $k\geq 0$ there is an equivalence of categories
\begin{equation*}
 (\tr{2}X)_k\simeq X_k\;.
\end{equation*}
By using transport of structure (more precisely Lemma \ref{lem-PP} with $\clC=\Dop$) we can lift $\tr{2}X$ to a pseudo-functor
\begin{equation*}
    \tr{2}X \in \psc{}{\Cat}
\end{equation*}
and by construction $\tr{2}X \in\segpsc{}{\Cat}$.

Building pseudo-functors from $\tawg{n}$ when $n>2$ is much more complex. The above approach cannot be applied directly because the induced Segal maps, when $n>2$ are $(n-1)$-equivalences but not in general levelwise equivalence of categories. For this reason we need to introduce an intermediate category $\lta{n}$, from which it is possible to build pseudo-functors using transport of structure. The functor from $\tawg{n}$ to Segalic pseudo-functors factorizes as
\begin{equation*}
    \tawg{n}\xrw{P_n} \lta{n} \xrw{\Tr_{n}} \segpsc{n-1}{\Cat}
\end{equation*}
 The functor $P_n$ produces a functorial approximation (up to an $n$-equivalence) of an object of $\tawg{n}$ with an object of $\lta{n}$, while the functor $\Tr_{n}$ is built using transport of structure.

 We call the functor $\Discn$ \emph{discretization functor} because it replaces a weakly globular structure with a globular one.

The idea of the functor $\Discn$ is to replace the homotopically discrete substructures in a weakly globular \nfol category by their discretization. This goes at the expenses of the Segal maps, which from being isomorphisms become higher categorical equivalences.

 In the higher groupoidal case, this idea had already appeared in the work by the author in \cite{Pa} and \cite{BP}, but further work is needed in the general categorical case to deal with the functorial choice of the discretization maps for the homotopically discrete substructures.

 We illustrate this in the case $n=2$. Given $X\in\catwg{2}$, by definition $X_0\in\cathd{}$, so there is a discretization map $\zg:X_0\rw X_0^d$ which is an equivalence of categories. Given a choice $\zg'$ of pseudo-inverse, we have $\zg\zg'=\Id$ since $X_0^d$ is discrete.

We can therefore construct $D_0 X\in\funcat{}{\Cat}$ as follows
\begin{equation*}
    (D_0 X)_k = \left\{
                  \begin{array}{ll}
                    X_0^d, & k=0 \\
                    X_k, & k>0\;.
                  \end{array}
                \right.
\end{equation*}
The face maps
 \begin{equation*}
 (D_0 X)_1\rightrightarrows (D_0 X)_0
 \end{equation*}
  are given by $\zg \pt_i$\; $i=0,1$ (where $\pt_i:X_1 \rightrightarrows X_0$ are face maps of $X$) while the degeneracy map
   \begin{equation*}
   (D_0 X)_0 \rw (D_0 X)_1
   \end{equation*}
    is $\zs_0\zg'$ (where $\zs_0:X_0 \rw X_1$ if the degeneracy map of $X$). All other face and degeneracy maps in $D_0 X$ are as in $X$. Since $\zg\zg'=\Id$, all simplicial identities are satisfied for $D_0X$. By construction, $(D_0 X)_0$ is discrete while the Segal maps are given, for each $k\geq 2$, by
\begin{equation*}
    \pro{X_1}{X_0}{k}\rw \pro{X_1}{X_0^d}{k}
\end{equation*}
and these are equivalences of categories since $X\in\catwg{2}$. Thus, by definition, $D_0 X\in\ta{2}$. This construction however does not afford a functor
 \begin{equation*}
 D_0:\catwg{2}\rw\ta{2}
 \end{equation*}
  but only a functor
\begin{equation*}
    D_0:\catwg{2}\rw(\ta{2})_{\ps}
\end{equation*}
where $(\ta{2})_{\ps}$ is the full subcategory of $\psc{}{\Cat}$ whose objects are in $\ta{2}$. This is because, for any morphism $F:X\rw Y$ in $\ta{2}$, the diagram in $\Cat$
\begin{equation*}
\xymatrix@C=35pt{
X_0^d \ar^{f^d}[r] \ar_{\zg'(X_0)}[d] & Y_0^d \ar^{\zg'(Y_0)}[d]\\
X_0 \ar_{f}[r] & Y_0
}
\end{equation*}
in general only pseudo-commutes.

   To remedy this problem we introduce the category $\ftawg{n}$ which exhibits functorial choices of the homotopically discrete  substructures in $\catwg{n}$ and from which the discretization process can be done functorially, using an iteration of the above idea, via a functor
   \begin{equation*}
     D_n:\ftawg{n} \rw \ta{n}.
   \end{equation*}
  We then show that we can approximate any object of $\catwg{n}$ with an $n$-equivalent object of $\ftawg{n}$. Namely we prove in Theorem \ref{pro-fta-1} that there is a functor
\begin{equation*}
  G_n:\catwg{n}\rw \ftawg{n}
\end{equation*}
and an $n$-equivalence $G_n X\rw X$. The discretization functor $\Discn$ is then defined as the composite
 \begin{equation*}
    \catwg{n}\xrw{G_n}\ftawg{n}\xrw{D_n}\ta{n}\;.
\end{equation*}

It is straightforward (Corollary \ref{cor-gta-1}) that $\Qn$ and $\Discn$ restrict to the groupoidal versions of the models
\begin{equation*}
  \Qn: \gta{n} \leftrightarrows \gcatwg{n}:\Discn
\end{equation*}
and that these induce equivalence of categories after localization with respect to the $n$-equivalences (Proposition \ref{pro-gta-1}).
 Using the result of Tamsamani \cite{Ta} it then follows (Theorem \ref{cor-gta-2}) that
\begin{equation}\label{eq2-the-cor-gta-2}
  \gcatwg{n}\bsim^n\;\simeq\;\Ho\mbox{($n$-types)}\;.
\end{equation}
In Corollary \ref{cor1-mod-fund-wg-group} the equivalence of categories \eqref{eq2-the-cor-gta-2} is obtained with a different fundamental higher groupoid functor
\begin{equation*}
  \Top \xrw{\clH_n}\gpdwg{n}\overset{j}{\hookrightarrow}\gcatwg{n}
\end{equation*}
where $\clH_n$ is the functor from spaces to weakly globular $n$-fold groupoids of Blanc and the author in \cite{BP} and $j$ is the inclusion. This has the advantage that $\clH_n$ is independent of \cite{Ta} and has a very explicit form, as we illustrate in some low-dimensional examples at the end of Section \ref{mod-fund-wg-group}.

%%%%%%%%%%%%%%%%%%%%%%%%%%%%%%%%%%%%%%%%%%%%%%%%%%%%%%%%%%%%%%%%%%%%%%%%%%%%%%%%%%%%%%%%

\subsection{Organization of this work}\label{sub-org}

We conclude this part with an account of the overall organization of this work. A more detailed synopsis of the content of each part can be found at the beginning of each, where we also provide some diagrammatic summaries.

Part \ref{part-back}  contains an account of the main background techniques from simplicial homotopy theory and from category theory used in this work and also fixes the relative notation. The content of this part is essentially expository, and aims to facilitate the reader with following the technical aspects of this work.

The core of this work is developed in Parts \ref{part-3}, \ref{part-4}, \ref{part-5}.

In part \ref{part-3} we introduce some of the new structures in this work and study their properties: the category $\cathd{n}$ of homotopically discrete \nfol categories (Chapter \ref{chap3}), the category $\catwg{n}$ of weakly globular \nfol categories (Chapter \ref{chap4}) and the category $\segpsc{n-1}{\Cat}$ of Segalic pseudo-functors (Chapter \ref{chap5}).

The main result of this part is Theorem \ref{the-strict-funct}, establishing that the classical strictification of pseudo-functors, when restricted to Segalic pseudo-functors, yields weakly globular \nfol categories; that is, there is a functor
 \begin{equation*}
        \St: \segpsc{n-1}{\Cat}\rw \catwg{n}
    \end{equation*}

 Part \ref{part-4} introduces the more general of the three Segal-type models of this work, the category $\tawg{n}$ of weakly globular Tamsamani $n$-categories (Chapter \ref{chap6}). The main result of this part is  Theorem \ref{the-funct-Qn} which constructs the rigidification functor
    \begin{equation*}
        Q_n:\tawg{n} \rw \catwg{n}
    \end{equation*}
This uses crucially the functor $\St$ of Theorem \ref{the-strict-funct} as well as several other intermediate steps, such as the subcategory $\lta{n}\subseteq \tawg{n}$. These are developed in Chapter \ref{chap7}.

In Part \ref{part-5} we establish the main comparison result between Tamsamani $n$-categories and weakly globular \nfol categories, exhibiting the latter as a model of weak $n$-categories. One of the main constructions is the discretization functor
 \begin{equation*}
 \Discn : \catwg{n}\rw\ta{n}
 \end{equation*}
 and its properties established in Theorem \ref{the-disc-func}. This needs several intermediate steps, in particular the new category $\ftawg{n}$, developed in Chapter \ref{chap8}.

  The rigidification and discretization functors then lead to the main comparison result, Theorem \ref{cor-the-disc-func}, on the equivalence of categories
  \begin{equation*}
   \ta{n}\bsim^n\;\simeq \; \catwg{n}\bsim^n
  \end{equation*}

  Chapter \ref{chap9} concludes with a proof of the homotopy hypothesis, after introducing the higher groupoidal version of the three Segal-type models.

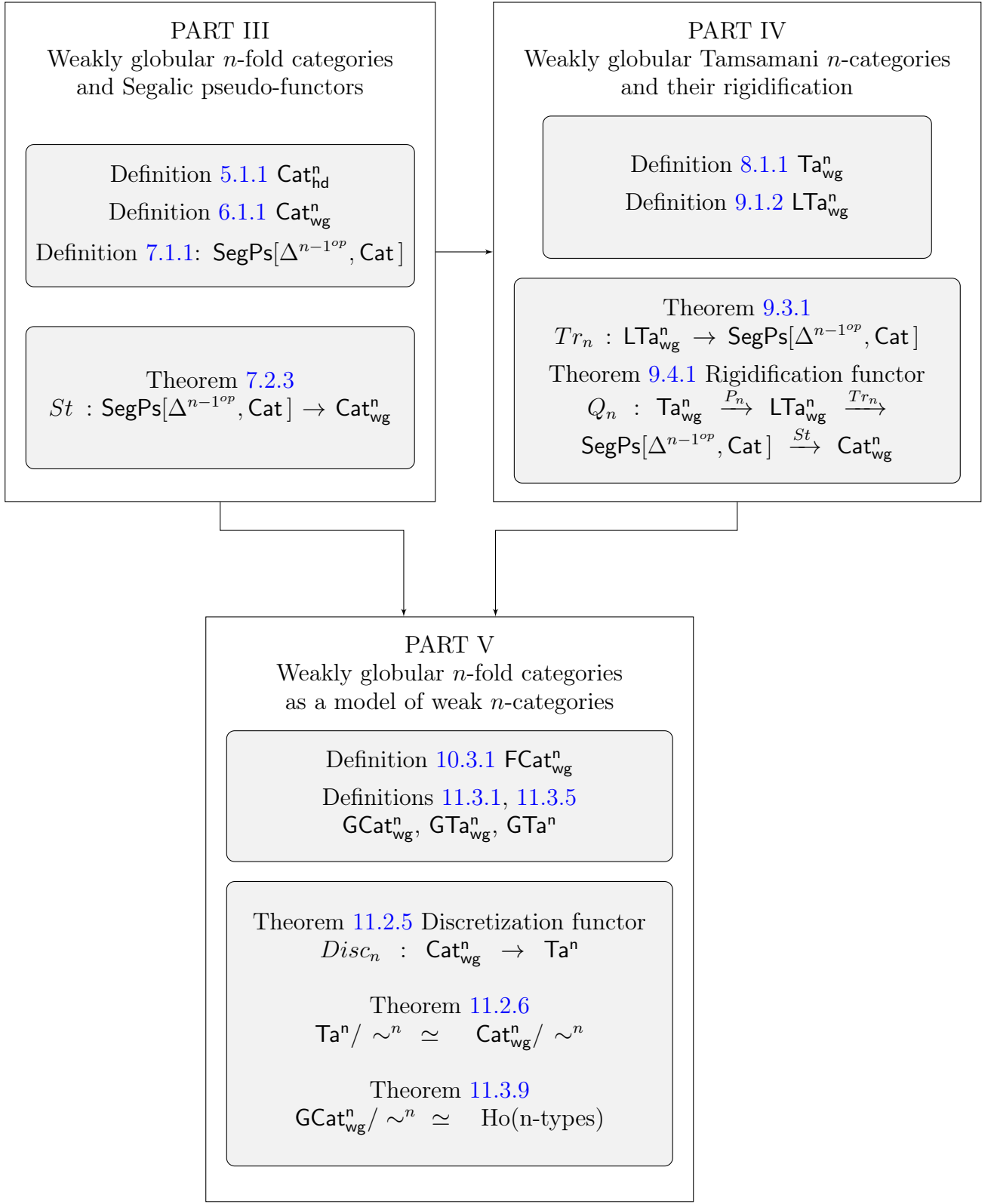
\begin{figure}[ht]
\begin{adjustwidth}{-1.5cm}{-1cm}
  \centering
\vspace{10mm}
\begin{tikzpicture}[node distance = 15mm and 10mm][width=16cm]% First vertical distance and second horizontal from border
\node[draw, text centered, fill=white,text width = 70mm,text depth = 85mm, minimum width=75mm] (maina){\\PART III\\Weakly globular $n$-fold categories\\and Segalic pseudo-functors};
\node [block1] (BLK1) at (maina.north) [below=25mm] {Definition \ref{def-hom-dis-ncat} $\cathd{n}$\\ \mk Definition \ref{def-n-equiv} $\catwg{n}$\\ \mk Definition \ref{def-seg-ps-fun}: $\segpsc{n-1}{\Cat}$};
\node [block1, below of = BLK1, yshift=-17mm] (BLK2){Theorem \ref{the-strict-funct}\\$\St:\segpsc{n-1}{\Cat}\rw\catwg{n}$};
\node[draw, text centered, fill=white,text width = 75mm,text depth = 85mm, minimum width=85mm, right = of maina] (mainb){\\PART IV\\Weakly globular Tamsamani $n$-categories\\and their rigidification};
\node [block1] (BLK2) at (mainb.north) [below=20mm] {Definition \ref{def-wg-ps-cat} $\tawg{n}$\\ \mk Definition \ref{def-ind-sub-ltawg} $\lta{n}$};
\node [rectangle, draw, fill=gray!10,
    text width=75mm, text centered, rounded corners, minimum height=36mm] [below of = BLK2, yshift=-19mm] (BLK3){Theorem \ref{the-XXXX}\\$\tr{n}:\lta{n}\rw\segpsc{n-1}{\Cat}$\\ \mk Theorem \ref{the-funct-Qn} Rigidification functor\\ $\Qn:\tawg{n}\xrw{P_n}\lta{n}\xrw{\tr{n}}\segpsc{n-1}{\Cat}\xrw{\St}\catwg{n}$};
\node [nullblock, below= of maina , xshift=40mm, yshift=15mm] (null) {};
\node[draw, text centered, fill=white,text width = 75mm,text depth = 100mm, minimum width=85mm, below = of null] (mainc){\\PART V\\Weakly globular  $n$-fold categories\\as a model of weak $n$-categories};
\node [rectangle, draw, fill=gray!10,
    text width=75mm, text centered, rounded corners, minimum height=23mm] (BLK4) at (mainc.north) [below=20mm]  {Definition \ref{def-fta-1} $\ftawg{n}$\\ \mk Definitions \ref{def-gta-1}, \ref{def-gta-2}\\$\gcatwg{n}$, $\gtawg{n}$, $\gta{n}$};
\node [rectangle, draw, fill=gray!10,
    text width=75mm, text centered, rounded corners, minimum height=50mm] [below of = BLK4, yshift=-25mm] (BLK5)  {Theorem \ref{the-disc-func} Discretization functor\\$\Discn:\catwg{n}\rw\ta{n}$\\ \ \\ Theorem \ref{cor-the-disc-func}\\$\ta{n}\bsim^n\;\simeq \; \catwg{n}\bsim^n$\\ \ \\ Theorem \ref{cor-gta-2}\\$\gcatwg{n}\bsim^n\;\simeq\;\mbox{Ho(n-types)}$};
%

%% Draw edges
% \draw [line] (pro8_1_5.east)++(0,-0.5) -- ++(1,0) --++(0,-3) |- (pro8_2_5.west);
\path [line] (maina) -- (mainb);
%\draw [line] (maina.south)++(0,0) -- ++(0,1) --++(3,1) |- (mainc.north);
\path [line] (maina.south) -- ++(0,-.5cm) -| node [near start] {} ($(mainc.north) + (-0.8cm, 0em)$);
\path [line] (mainb.south) -- ++(0,-.5cm) -| node [near start] {} ($(mainc.north) + (0.8cm, 0em)$);

\end{tikzpicture}
  \caption{Summary of overall organization and main results}
  \label{FigIntro-2}
\end{adjustwidth}
\end{figure}

  In Figure \ref{FigIntro-2} below we give a schematic account of the main results of this work.
%%%%%%%%%%%%%%%%%%%%%%%%%%%%%%%%%%%%%%%%%%%%%%%%%%%%%%%%%%%%%%%%%%%%%%%%%%
\subsection{List of informal discussions}\index{List of informal discussions} We have included a number of informal discussions throughout the text to convey the ideas and intuitions behind the main definitions and constructions. A list of the main ones is as follows:

\bk

Section \ref{subs-idea-cathd}: The idea of homotopically discrete \nfol category.\bk

Section \ref{subs-idea-catwg}: The idea of weakly globular \nfol categories.\bk

Section \ref{subs-idea-segps}: The idea of Segalic pseudo-functor.\bk

Section \ref{sub-idea-tawgn}: The idea of weakly globular Tamsamani $n$-categories.\bk

Section \ref{subs-idea-ltan}: The idea of the category $\lta{n}$.\bk

Section \ref{subs-approx-idea}: Main steps in approximating $\tawg{n}$ with $\lta{n}$.\bk

Section \ref{subs-idea-trn}: The idea of the functor $\tr{n}$.\bk

Section \ref{sub-idea-qn} The rigidification functor $Q_n$: main steps.\bk

Section \ref{subs-idea-constrxf}: The idea of the construction $X(f_0)$. \bk

Section \ref{subs-idea-fn}: The idea of the functors $V_n$ and $F_n$.\bk

Section \ref{sub-idea-ftawgn}: The idea of the category $\ftawg{n}$.\bk

Section \ref{subs-idea-gn} The idea of the functor $G_n$.\bk

Section \ref{subs-fn-dn-idea}: The idea of the functor $\D{n}$. \bk

Section \ref{subs-idea-discn}: The idea of the functor $\Discn$.

\label{organization}

%

%%%%%%%%%%%%%%%%%%%%%%%%%%%%%%%%%%%%%%%%%%%%%%%%%%%%%%%%%%%%%%%%%%%%%%%%%%%%%

\part{Background techniques} \label{part-back}

The techniques used throughout this work draw from two sources: simplicial homotopy theory and category theory. Part \ref{part-back} reviews these techniques, more precisely  we cover multi-simplicial techniques in Chapter \ref{prelim}, and categorical background in Chapter \ref{chap-ps}.

  We recall in Section \ref{sbs-simp-tech} the central notion of Segal maps for (multi-simplicial) objects and fix our notation for multi-simplicial objects. We then review in Section \ref{sbs-nint-cat} internal $n$-fold categories and in Section \ref{sus-ner-funct} their multi-nerves. We use these multi-simplicial notions to give in Section \ref{sbs-multi-strict} the multi-simplicial description of strict $n$-categories. This is very important to build a geometric intuition around the Segal-type models of this work. For this reason, we illustrate in detail some low-dimensional cases and draw corresponding pictures. In Section \ref{decalage} we recall the functor d\'{e}calage, which is used in Chapters \ref{chap8} and \ref{chap9}.

 We recall in Section \ref{sbs-funct-cat} two important functors from $\Cat$ to $\Set$, the isomorphism classes of objects and the connected components functors, which play a crucial role in the Segal-type models of this work.

 The last part of Chapter \ref{chap-ps} reviews the notion of pseudo-functor and their strictification, as well as a standard technique to produce pseudo-functors, which is an instance of 'transport of structure along an adjunction' in the sense of \cite{lk}. These techniques play a crucial role in this work, in relation to the new notion of Segalic pseudo-functor studied in Chapters \ref{chap5} and \ref{chap7}.

 The material in Part \ref{part-back} is essentially known, but we have presented it in a way that is best suited for the rest of this work, with emphasis on the multi-simplicial description of strict $n$-categories and $n$-fold categories: the latter is not always spelled out at this level of detail in the literature, and it is crucial to build the intuition around our Segal-type models.

 Suitable references for Part \ref{part-back} include \cite{Borc}, \cite{Jard}, \cite{lk}, \cite{Lack}, \cite{MacLane1998}, \cite{MayBook1967}, \cite{PW}.

%%%%%%%%%%%%%%%%%%%%%%%%%%%%%%%%%%%%%%%%%%%%%%%%%%%%%%%%%%%%%%%%%%%%%%%%%%%%%%%

%%%%%%%%%%%%%%%%%%%%%%%%%%%%%%%%%%%%%%%%%%%%%%%%%%%%%%%%%%%%%%%%%%%%%%%%%%%%%

\chapter{Multi-simplicial techniques}\label{prelim}

 This chapter reviews (multi)-simplicial techniques. Multi-simplicial objects are functors from $\Dnop$ to a category $\mathcal{C}$ with finite limits. When $\mathcal{C}$ is the category $\Set$, multi-simplicial sets are combinatorial models of spaces and are used in homotopy theory, see for instance \cite{Jard}.

As outlined in Chapter \ref{chap1} the connection between higher categorical structures and multi-simplicial structures arises via the notions of nerves (and their iterated versions, called multi-nerves) and of Segal maps. We review these in this chapter, where we also spell out in detail the multi-simplicial description of strict $n$-categories and $n$-fold categories. In the cases $n=2$ and $n=3$ these descriptions admit a geometric interpretation which is easy to visualize (see Examples \ref{ex-mult-ner-dim3} and \ref{ex-2-multi-strict}). Such a geometric visualization is very helpful to build the intuition around the modifications needed to build a weak model of higher categories using multi-simplicial structures.

\section{Multi-simplicial objects and Segal maps}\label{sbs-simp-tech}

We start by reviewing the notion of multi-simplicial objects and their associated Segal maps. These play a crucial role throughout this work since, as explained in the previous chapter, multi-simplicial objects form a natural environment for the definition of higher categorical structures.

\subsection{Simplicial objects and their Segal maps}\label{subs-simplob}

\index{Segal map}

 \index{Simplicial!- object} \index{Simplicial!- map}

Let $\Delta$ be the simplicial category. Its objects are finite ordered sets \index{Simplicial! - category}
 \begin{equation*}
   [n]=\{0 < 1 < \cdots < n\}
 \end{equation*}
 for integers $n\geq 0$ and its morphisms are non decreasing monotone functions. If $\clC$ is any category, a simplicial object $X$ in $\clC$ is a contravariant functor from $\Delta$ to $\clC$, that is $X:\dop{}\rw\clC$. We write $X_n$ for $X([n])$. A simplicial map $F:X \rw Y$ between simplicial objects in $\clC$ is a natural transformation. We denote by $\funcat{}{\clC}$ the category of simplicial objects and simplicial maps.

It is well known (see for instance \cite{MayBook1967}) that to give a simplicial object $X$ in $\clC$  is the same as to give a sequence of objects $X_0, X_1, X_2, \ldots$ together with face operators $\pt_i:X_n\rw X_{n-1}$ and degeneracy operators $\zs_i: X_n\rw X_{n+1}$ ($i=0,\ldots,n$) which satisfy the following simplicial identities:
\begin{equation}\label{eq-simp-ob-segmap}
\begin{split}
  \pt_i\pt_j = &  \pt_{j-1}\pt_{i} \;\;\qquad \mbox{if $i<j$}\\
  \zs_i\zs_j = &  \zs_{j-1}\zs_{i} \;\;\qquad \mbox{if $i<j$}\\
  \pt_i\zs_j= &
  \left\{
    \begin{array}{ll}
      \zs_{j-1}\pt_i, & \mbox{if $i<j$} \\
      \Id, & \mbox{if $i=j$ or $i=j+1$} \\
      \zs_{j}\pt_{i-1}, & \mbox{if $i>j+1$}
    \end{array}
  \right.
\end{split}
\end{equation}
Under this correspondence $\pt_i=A(\zve_i)$ and $\zs_i=A(\eta_i)$ where
\begin{equation*}
\begin{split}
  \zve_i: [n-1]\rw [n] \qquad\qquad & \qquad\qquad\eta_i: [n+1]\rw [n] \\
   \zve_i(j)=
   \left\{
     \begin{array}{ll}
       j, & \mbox{if $j<i$} \\
       j+1, & \mbox{if $j\geq i$}
     \end{array}
   \right.
    &
   \qquad\qquad\eta_i(j)=
   \left\{
     \begin{array}{ll}
       j, & \mbox{if $j\leq i$} \\
       j-1, & \mbox{if $j>i$.}
     \end{array}
   \right.
\end{split}
\end{equation*}

\begin{definition}\label{def-seg-map}
    Let ${X\in\funcat{}{\clC}}$ be a simplicial object in any category $\clC$ with pullbacks. For each ${1\leq j\leq k}$ and $k\geq 2$, let ${\nu_j:X_k\rw X_1}$ be induced by the map  $[1]\rw[k]$ in $\Delta$ sending $0$ to ${j-1}$ and $1$ to $j$. Then the following diagram commutes:
\begin{equation}\label{eq-seg-map}
\xymatrix@C=20pt{
&&&& X\sb{k} \ar[llld]_{\nu\sb{1}} \ar[ld]_{\nu\sb{2}} \ar[rrd]^{\nu\sb{k}} &&&& \\
& X\sb{1} \ar[ld]_{d\sb{1}} \ar[rd]^{d\sb{0}} &&
X\sb{1} \ar[ld]_{d\sb{1}} \ar[rd]^{d\sb{0}} && \dotsc &
X\sb{1} \ar[ld]_{d\sb{1}} \ar[rd]^{d\sb{0}} & \\
X\sb{0} && X\sb{0} && X\sb{0} &\dotsc X\sb{0} && X\sb{0}
}
\end{equation}

If  ${\pro{X_1}{X_0}{k}}$ denotes the limit of the lower part of the
diagram \eqref{eq-seg-map}, the \emph{$k$-th Segal map for $X$} is the unique map
$$
\muk:X\sb{k}~\rw~\pro{X\sb{1}}{X\sb{0}}{k}
$$
\noindent such that ${\pr_j\,\muk=\nu\sb{j}}$ where
${\pr\sb{j}}$ is the $j^{th}$ projection.
\end{definition}
\begin{definition}\label{def-ind-seg-map}

    Let ${X\in\funcat{}{\clC}}$ and suppose that there is a map
     \begin{equation*}
       \zg: X_0 \rw Y
     \end{equation*}
      in $\clC$    such that the limit of the diagram
\begin{equation*}
\xymatrix@R25pt@C16pt{
& X\sb{1} \ar[ld]_{\zg d\sb{1}} \ar[rd]^{\zg d\sb{0}} &&
X\sb{1} \ar[ld]_{\zg d\sb{1}} \ar[rd]^{\zg d\sb{0}} &\cdots& k &\cdots&
X\sb{1} \ar[ld]_{\zg d\sb{1}} \ar[rd]^{\zg d\sb{0}} & \\
Y && Y && Y\cdots &&\cdots Y && Y
    }
\end{equation*}
exists; denote the latter by $\pro{X_1}{Y}{k}$. Then the following diagram commutes, where $\nu_j$ is as in Definition \ref{def-seg-map}, and $k\geq 2$
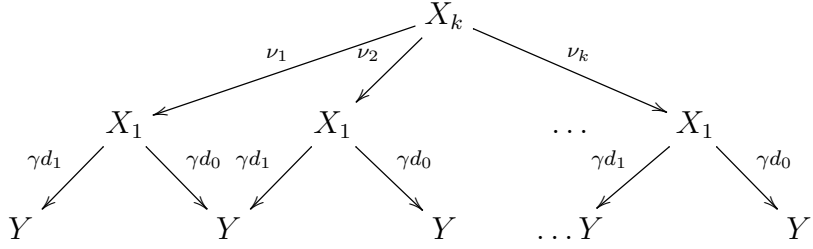
\begin{equation*}
\xymatrix@C=20pt{
&&&& X\sb{k} \ar[llld]_{\nu\sb{1}} \ar[ld]_{\nu\sb{2}} \ar[rrd]^{\nu\sb{k}} &&&& \\
& X\sb{1} \ar[ld]_{\zg d\sb{1}} \ar[rd]^{\zg d\sb{0}} &&
X\sb{1} \ar[ld]_{\zg d\sb{1}} \ar[rd]^{\zg d\sb{0}} && \dotsc &
X\sb{1} \ar[ld]_{\zg d\sb{1}} \ar[rd]^{\zg d\sb{0}} & \\
Y && Y && Y &\dotsc Y && Y
}
\end{equation*}
The \emph{$k$-th induced Segal map for $X$} is the unique map \index{Induced Segal map}
\begin{equation*}
\hmuk:X\sb{k}~\rw~\pro{X\sb{1}}{Y}{k}
\end{equation*}
such that ${\pr_j\,\hmuk=\nu\sb{j}}$ where ${\pr\sb{j}}$ is the $j^{th}$ projection. Clearly if $Y=X_0$ and $\gamma$ is the identity, the induced Segal map coincides with the Segal map of Definition \ref{def-seg-map}.
\end{definition}
\subsection{Multi-simplicial objects}\label{subs-multisimplob} \index{Multi-simplicial object}
Let $\Dnop$ denote the product of $n$ copies of $\Dop$, that is
 \begin{equation*}
   \Dnop=(\Delta^{\text{op}})^n= \dop{}\times \oset{n}{\cdots} \times \dop{}
 \end{equation*}
 Given a category $\clC$, $\funcat{n}{\clC}$ is called the category of $n$-fold simplicial objects in $\clC$ (or multi-simplicial objects for short).

\index{n-fold!- simplicial object}

\begin{notation}\label{not-simp}
    If $X\in \funcat{n}{\clC}$ and $\uk=([k_1],\ldots,[k_n])\in \Dnop$, we shall use several alternative notations as follows

     \begin{equation*}
       X ([k_1],\ldots,[k_n])=X(k_1,\ldots,k_n)=X_{k_1,\ldots,k_n}=X_{\uk}
     \end{equation*}

     We shall also denote for $1\leq i\leq n$

      \begin{equation*}
      \uk(1,i)=([k_1],\ldots,[k_{i-1}],1,[k_{i+1}],\ldots,[k_n]) \in \Dnop.
      \end{equation*}

\end{notation}
The following is an elementary fact:

\begin{lemma}\label{lem-multi-simpl-as}
  Every $n$-simplicial object in $\clC$ can be regarded as a simplicial object in $\funcat{n-1}{\clC}$ in $n$ possible ways.
\end{lemma}

     \begin{proof} For each $1\leq i\leq n$ there is an isomorphism
    \begin{equation*}
        \xi_i:\funcat{n}{\clC}\rw\funcat{}{\funcat{n-1}{\clC}}
    \end{equation*}
    given by
    \begin{equation*}
        (\xi_i X)_{r}(k_1,\ldots,k_{n-1})=X(k_1,\ldots,k_{i-1},r,k_{i+1},\ldots,k_{n-1})
    \end{equation*}
    for $X\in\funcat{n}\clC$ and $r\in\Dop$.
\end{proof}

From Lemma \ref{lem-multi-simpl-as}, given an $n$-fold simplicial object $X\in\funcat{n}{\clC}$, for each simplicial direction $1\leq i \leq n$ we have Segal maps for $\xi_i X\in[\dop{},\funcat{n-1}{\clC}]$. For each $r\geq 2$ these are maps in $\funcat{n-1}{\clC}$
\begin{equation*}
  (\xi_i X)_r \rw \pro{(\xi_i X)_1}{(\xi_i X)_0}{r}
\end{equation*}
where $(\xi_i X)_r$ is as in the proof of Lemma \ref{lem-multi-simpl-as}.
\begin{example}\label{corner-n2}
Consider the case $n=2$. Then $X\in\funcat{2}{\clC}$ is called a bisimplicial object in $\clC$. \index{Bisimplicial object} It is equivalent to a bigraded sequence of objects $X_{st}$ of $\clC$ ($s,t \geq 0$) together with horizontal face and degeneracy maps
\begin{equation*}
  \pt^h_i : X_{st}\rw X_{s-1,t} \qquad \zs^h_i : X_{st}\rw X_{s+1,t}
\end{equation*}
as well as vertical face and degeneracy maps
\begin{equation*}
  \pt^v_i : X_{st}\rw X_{s,t-1} \qquad \zs^v_i : X_{st}\rw X_{s,t+1}\;.
\end{equation*}
These face and degeneracy maps must satisfy the simplicial identities (horizontally and vertically) and every horizontal map must commute with every vertical map. We call these equalities bisimplicial identities. In Figure \ref{fig-bisimpl-set} we draw a pictorial representation of the corner of a bisimplicial object $X$, where direction 1 is horizontal and direction 2 is vertical:
\begin{figure}[ht]
\centering
\begin{equation*}
\xymatrix@R=5pt @C=75pt{ &\!\!\!\!\!\! \vdots & \vdots & \vdots}
\end{equation*}
\begin{equation*}
\entrymodifiers={+++[o]}
 \xymatrix@R=40pt@C=50pt{
    \cdots \ssr & X_{22} \ar@<2ex>[r] \ar[r] \ar@<1ex>[r] \ar@<-2ex>[d] \ar[d]
    \ar@<-1ex>[d]& X_{12} \ar@<1ex>[r] \ar@<0ex>[r]  \ar@<-2ex>[d] \ar[d] \ar@<-1ex>[d] \ar@<1ex>[l] \ar@<2ex>[l] &
    X_{02}  \ar@<-2ex>[d] \ar[d] \ar@<-1ex>[d] \ar@<1ex>[l]\\
    \cdots \ssr  & X_{21}\ar@<2ex>[r] \ar[r] \ar@<1ex>[r] \ar@<-1ex>[d] \ar@<0.ex>[d] \ar@<-1ex>[u] \ar@<-2.ex>[u] &
    X_{11} \ar@<1ex>[r] \ar@<0ex>[r]  \ar@<-1ex>[d] \ar@<0.ex>[d] \ar@<1ex>[l] \ar@<2ex>[l] \ar@<-1ex>[u] \ar@<-2.ex>[u] &
    X_{01}  \ar@<0.ex>[d] \ar@<-1ex>[d] \ar@<1ex>[l] \ar@<-1ex>[u] \ar@<-2.ex>[u] \\
    \cdots \ssr & X_{20} \ar@<2ex>[r] \ar[r] \ar@<1ex>[r] \ar@<-1ex>[u] &
     X_{10} \ar@<1ex>[r] \ar@<0ex>[r] \ar@<1ex>[l] \ar@<2ex>[l] \ar@<-1ex>[u] &
    X_{00}\ar@<1ex>[l] \ar@<-1ex>[u]\\
}
\end{equation*}
\caption{Corner of a bisimplicial set $X$}
\label{fig-bisimpl-set}
\end{figure}
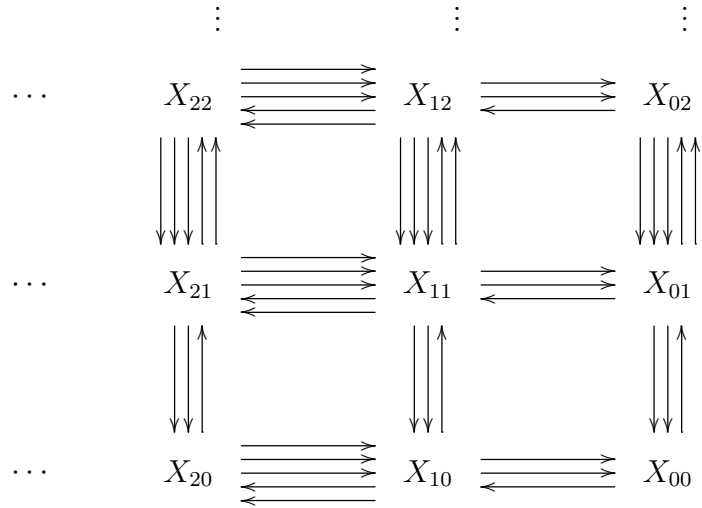

\bk
We see from the picture that a bisimplicial object is a simplicial object in simplicial objects in two directions, vertical and horizontal.

For each $r\geq 0$ denote by $X_{r*}\in\funcat{}{\clC}$ the vertical simplicial object with
\begin{equation*}
  (X_{r*})_k=X_{rk}
\end{equation*}
and by $X_{*r}\in\funcat{}{\clC}$ the horizontal simplicial object with
\begin{equation*}
  (X_{*r})_k=X_{kr}\;.
\end{equation*}
Then for each $r\geq 2$ there are horizontal Segal maps in $\funcat{}{\clC}$
\begin{equation*}
  X_{r*}=\pro{X_{1*}}{X_{0*}}{r}
\end{equation*}
and vertical Segal maps in $\funcat{}{\clC}$
\begin{equation*}
  X_{*r}=\pro{X_{*1}}{X_{*0}}{r}
\end{equation*}
\end{example}
\section{$\pmb{n}$-Fold internal categories}\label{sbs-nint-cat}

\index{n-fold!- category} \index{Internal!- category} \index{Internal!- groupoid} \index{Internal!- functor}

\begin{definition}\label{def-intercat}
Let  $\clC$ be a category with finite limits. An internal category $X$ in $\clC$ is a diagram in $\clC$
\begin{equation}\label{eq-nint-cat}
\xymatrix{
\tens{X_1}{X_0} \ar^(0.65){m}[r] & X_1 \ar^{d_0}[r]<2.5ex> \ar^{d_1}[r] & X_0
\ar^{s}[l]<2ex>
}
\end{equation}
where the maps $d_0, d_1, s, m$ satisfy the axioms:
\begin{itemize}
  \item [(1)] $d_0 s=d_1 s= \Id_{X_0}$.\bk

  \item [(2)] $d_1 p_2 = d_1 m$, $d_0 p_1 = d_0 m$, where $p_i: \tens{X_1}{X_0}\rw X_1$, $i=1,2$ are the two projections.\bk

  \item [(3)] $\ds m\cirsm \begin{pmatrix}\Id_{X_1} \\ s\cirsm d_0 \\\end{pmatrix}=
  \Id_{X_1}=\begin{pmatrix}s\cirsm d_1\\ \Id_{X_1} \\\end{pmatrix}$.\bk

  \item [(4)] $m\cirsm(\Id_{X_1}\tiund{X_0} m)= m\cirsm(m \tiund{X_0} \Id_{X_1})$
\end{itemize}
\end{definition}
\bk

\nid $X_0$ is called 'object of objects', $X_1$ 'object of arrows', $d_0$ and $d_1$ are called respectively 'source' and 'target' maps, $s$ is the 'identity map' and $m$ is the 'composition map'. When $\clC = \Set$ this gives the axioms of a small category and $\tens{X_1}{X_0}$ is the set of composable arrows.

\begin{definition}\label{def-interfunc}
If $X$ and $Y$ are internal categories in $\clC$ , an internal functor $F:X \rw Y$ is a diagram in $\clC$
\begin{equation*}
\xymatrix@R=45pt@C=35pt{
\tens{X_1}{X_0} \ar^(0.6){m}[r] \ar_{\tens{F_1}{F_0}}[d] & X_1 \ar@<2ex>^{d_0}[r] \ar@<0ex>^{d_1}[r] \ar_{F_1}[d] & X_0 \ar@<1.5ex>^{s}[l] \ar^{F_0}[d]\\
\tens{Y_1}{Y_0} \ar^(0.6){m'}[r]  & Y_1 \ar@<2.5ex>^{d'_0}[r] \ar@<-0.5ex>^{d'_1}[r]  & Y_0 \ar@<2ex>^{s'}[l]
}
\end{equation*}
where the corresponding maps commute, that is they satisfy
\begin{itemize}
  \item [(1)] $d'_0\cirsm F'_1 = F_0\cirsm d_0$, \qquad $d'_1\cirsm F'_1 = F_0\cirsm d_1$.\bk

  \item [(2)] $F_1\cirsm s = s'\cirsm F_0$.\bk

  \item [(3)] $F_1\cirsm m = m' \cirsm (\tens{F_1}{F_0})$.
\end{itemize}
\end{definition}

\bk

\nid We denote by $\Cat \clC$ the category of internal categories and internal functors. When $\clC=\Set$, this is the category $\Cat$ of small categories.

\begin{definition}\label{def-intergroup}
The category $\Gpd\,\clC$ of internal groupoids in $\clC$ is the full subcategory of $\Cat\clC$ whose objects $X$ are such that there is a map
\begin{equation*}
  i : X_1 \rw X_1
\end{equation*}
such that
\begin{equation*}
  m(\Id_{X_1},i)=s d_0,\qquad m(i,\Id_{X_1})=s d_1\;.
\end{equation*}
\end{definition}
The map $i$ gives the inverses to each (internal) arrow. When $\clC=\Set$, $\Gpd \clC$ is the category $\Gpd$ of small groupoids.
\begin{definition}\label{nfold internal}
The category $\cat{n}(\clC)$ of \nfol categories in $\clC$ is defined inductively by iterating $n$ times the internal category construction. That is, $\cat{1}(\clC)=\Cat \clC$ and, for $n>1$,

\index{n-fold!- internal category}

\begin{equation*}
  \cat{n}(\clC)= \Cat(\cat{n-1}(\clC)).
\end{equation*}
\end{definition}

When $\clC=\Set$, $\cat{n}(\Set)$ is simply denoted by $\cat{n}$ and called the category of \nfol categories (double categories when $n=2$).

\begin{definition}\label{nfoldgpd internal}
The category $\gpd{n}(\clC)$ of \nfol groupoids in $\clC$ is defined inductively by iterating $n$ times the internal groupoid construction. That is, $\gpd{1}(\clC)=\Gpd \clC$ and, for $n>1$,

\index{n-fold!- internal groupoid}

\begin{equation*}
  \gpd{n}(\clC)= \Gpd(\gpd{n-1}(\clC)).
\end{equation*}
\end{definition}

When $\clC=\Set$, $\gpd{n}(\Set)$ is simply denoted by $\gpd{n}$ and called the category of \nfol groupoids (double groupoids when $n=2$).

%%%%%%%%%%%%%%%%%%%%%%%%%%%%%%%%%%%%%%%%%%%%%%%%%%%%%%%%%%%%%%%%%%%%%%%%%%%%
\section{Multi-nerve functors}\label{sus-ner-funct}

\index{Multi-nerve functor}

\index{Nerve functor}

There is a nerve functor
\begin{equation*}
    N:\Cat\clC \rw \funcat{}{\clC}
\end{equation*}
such that, for $X\in\Cat\clC$
\begin{equation*}
    (N X)_k=
    \left\{
      \begin{array}{ll}
        X_0, & \hbox{$k=0$;} \\
        X_1, & \hbox{$k=1$;} \\
        \pro{X_1}{X_0}{k}, & \hbox{$k>1$.}
      \end{array}
    \right.
\end{equation*}
When no ambiguity arises, we shall sometimes denote $(NX)_k$ by $X_k$ for all $k\geq 0$.

The following fact is well known:
\begin{proposition}\label{pro-ner-int-cat}
    A simplicial object in $\clC$ is the nerve of an internal category in $\clC$ if and only if all the Segal maps are isomorphisms.
\end{proposition}

\begin{definition}\label{def-fun-smacat}
    Let $F:\clC \rw \clD$ be a functor, $\clI$ a small category. Denote
    \begin{equation*}
        \ovl{F}:[\clI,\clC]\rw [\clI,\clD]
    \end{equation*}
    the functor given, for all $i\in\clI$,  by
    \begin{equation*}
        (\ovl{F} X)_i = F(X(i)).
    \end{equation*}

\end{definition}

\begin{definition}\label{def-multinerve}
  By iterating the nerve construction, we obtain the multinerve functor
\begin{equation*}
    \Nn :\cat{n}(\clC)\rw \funcat{n}{\clC}\;.
\end{equation*}
More precisely, $\Nn$ is defined recursively as
\begin{equation*}
    \Nb{1}= N: \Cat \clC \rw \funcat{}{\clC}\;.
\end{equation*}
\end{definition}

Given $\Nb{n-1}: \cat{n-1}\clC \rw \funcat{n-1}{\clC}$ let
\begin{equation*}
    \Nn = \ovl{\Nb{n-1}}\cirsm N_1
\end{equation*}
where $\ovl{\Nb{n-1}}$ is as in Definition \ref{def-fun-smacat}.

\begin{definition}\label{def-discr-nfold}
  An internal \nfol category $X\in \cat{n}(\clC)$ is discrete when $\Nb{n-1} X$ is a constant functor.

\index{Discrete!- internal $n$-fold category}

\end{definition}

Thus in a discrete internal \nfol categories all structure maps are identities.

\begin{remark}\label{rem-mult-ner}
  Let $X\in\funcat{n}{\clC}$. For each $\uk= (\seqc{k}{1}{n})\in\dop{n}$ and $1\leq i\leq n$ we have
\begin{equation*}
    X(k_1,...,k_{i-1},-,k_{i+1},...,k_n)\in \funcat{ }{\clC}
\end{equation*}
and there are corresponding Segal maps for each $k_i\geq 2$
\begin{equation}\label{eq-rem-mult-ner}
  X_{\uk}\rw \pro{X_{\uk(1,i)}}{X_{\uk(0,i)}}{k_i}
\end{equation}
where $\uk(1,i)$ and $\uk(0,i)$ are as in Notation \ref{not-simp}.
\end{remark}

We claim that the condition that the maps \eqref{eq-rem-mult-ner} are isomorphisms (for all $\uk\in\dop{n}$) is equivalent to the condition that for each $1\leq i\le n$ the Segal maps for $\xi_i X\in\funcat{}{\funcat{n-1}{\clC}}$  are isomorphisms. This follows immediately from the fact that
\begin{align*}
    & (\xi_i X)_{k_i}(k_1,...,k_{i-1},k_{i+1},...,k_n)=X_{\uk} \\
    & (\xi_i X)_{0}(k_1,...,k_{i-1},k_{i+1},...,k_n)=X_{\uk(i,0)}\\
    & (\xi_i X)_{1}(k_1,...,k_{i-1},k_{i+1},...,k_n)=X_{\uk(i,1)}\;.
\end{align*}

\bk

Using this observation we can now give a characterization of those multi-simplicial objects $X\in\funcat{n}{\clC}$ which are in the image of the multinerve functor $\Nn :\cat{n}(\clC)\rw\funcat{n}{\clC}$. Namely, they are those multi-simplicial objects whose Segal maps in all directions are isomorphisms, as illustrated in the following Lemma.
\begin{lemma}\label{lem-multin-iff}\
\begin{itemize}
      \item [a)] $X\in\funcat{n}{\clC}$ is the multinerve of an object of $\cat{n}(\clC)$ if and only if, for each $\uk\in\dop{n}$, $k_i\geq 2$ and $1\leq i\leq n$
\begin{equation}\label{eq1-lem-multin-iff}
X_{\uk}\cong \pro{X_{\uk(1,i)}}{X_{\uk(0,i)}}{k_i}.
\end{equation}

      \item [b)] $X\in \funcat{n}{\clC}$ is the multinerve of an object of $\cat{n}(\clC)$ if and only if, for each $1\leq i\leq n$, all the Segal maps of $\xi_i X\in\funcat{}{\cat{n-1}(\clC)}$ are isomorphisms.
    \end{itemize}
\end{lemma}
\begin{proof}\

a) By induction on $n$. For $n=1$ this is Proposition \ref{pro-ner-int-cat}. Suppose the lemma holds for $n-1$ and let $\uk=(\seq{k}{1}{n})\in\dop{n}$. Denote $\ur=(\seq{k}{2}{n})\in\dop{n-1}$. Observe that
\begin{equation}\label{eq2-lem-multin-iff}
  X_{\uk}=(X_{k_1})_{\ur}, \qquad X_{\uk(1,i)}=(X_{k_1})_{\ur(1,i-1)}\qquad X_{\uk(1,0)}=(X_{k_1})_{\ur(0,i-1)}\;.
\end{equation}
Thus if \eqref{eq1-lem-multin-iff} holds, we have
\begin{equation*}
  (X_{k_1})_{\ur}=\pro{(X_{k_1})_{\ur(1,i)}}{(X_{k_1})_{\ur(0,i)}}{k_1}
\end{equation*}
for each $k_1\geq 2$. Hence $X_{k_1}$ satisfies the inductive hypothesis and we conclude that $X_{k_1}\in\cat{n-1}(\clC)$. Further, taking $i=1$ in \eqref{eq1-lem-multin-iff} we see that
\begin{equation*}
  (X_{k_1})_{\ur}=\pro{(X_1)_{\ur}}{(X_0)_{\ur}}{k_1}
\end{equation*}
that is, we have isomorphisms in $\cat{n-1}(\clC)$
\begin{equation*}
  X_{k_1}\cong \pro{X_1}{X_0}{k_1}
\end{equation*}
for each $k_1\geq 2$. By Proposition \ref{pro-ner-int-cat} it follows that $X\in\cat{n}(\clC)$.

Conversely, let $X=\Nn Y$ with $Y\in\cat{n}(\clC)$. Since $\Nn = \ovl{\Nb{n-1}}N_1$, it is
\begin{equation}\label{eq3-lem-multin-iff}
\begin{split}
    & X_{\uk}=(\Nb{n-1}(N_1 Y)_{k_1})_{\ur} \\
    & X_{\uk(1,i)}=(\Nb{n-1}(N_1 Y)_{k_1})_{\ur(1,i-1)}\\
    & X_{\uk(0,i)}=(\Nb{n-1}(N_1 Y)_{k_1})_{\ur(0,i-1)}\;.
\end{split}
\end{equation}
By induction hypothesis applied to $(N_1 Y)_{k_1}\in\cat{n-1}(\clC)$ and by \eqref{eq3-lem-multin-iff}, \eqref{eq1-lem-multin-iff} follows.

\bk

b) This follows immediately from a) and from Remark \ref{rem-mult-ner}.

\end{proof}
We now deduce that every object of $\cat{n}(\clC)$ can be considered as an internal category in $\cat{n-1}(\clC)$ in $n$ possible ways, corresponding to the $n$ simplicial directions of the multinerve.
\begin{proposition}\label{pro-mult-ner}
  For each $1\leq i\leq n$ there is an isomorphism
  \begin{equation*}
    \tilde{\xi_i}: \cat{n}(\clC) \rw \Cat(\cat{n-1}(\clC))
  \end{equation*}
making the following diagram commute:
\begin{equation}\label{eq1-pro-mult-ner}
\xymatrix{
\cat{n}(\clC) \ar^{\Nn}[r] \ar_{\tilde{\xi_i}}[d] & \funcat{n}{\clC} \ar^{\xi_i}[d] \\
\Cat(\cat{n-1}(\clC)) \ar_{N_1}[r] & \funcat{}{\funcat{n-1}{\clC}}
}
\end{equation}
where $\xi_i$ is an in the proof of Lemma \ref{lem-multi-simpl-as}.
\end{proposition}
\begin{proof}
We claim that, for each $Y\in\cat{n}(\clC)$, $(\xi_i \Nn Y)_k$ is the multinerve of an object $Z_{ik}\in \cat{n-1}(\clC)$. In fact, by Lemma \ref{lem-multin-iff}, $\Nn Y$ has all the Segal maps isomorphisms, hence the same holds for $(\xi_i \Nn Y)_k$ so that,  by Lemma \ref{lem-multin-iff} again, the claim follows. Consider the simplicial object
 \begin{equation*}
   (\xi_i \Nn Y)=Z_{i*}\in\funcat{}{\cat{n-1}(\clC)}.
 \end{equation*}
  By Lemma \ref{lem-multin-iff} b) its Segal maps are isomorphisms. Hence $Z_{ik}$ is the nerve of an internal category $\tilde{\xi_i}Y\in \Cat(\cat{n-1}(\clC))$, that is
\begin{equation*}
  N_1 \tilde{\xi_i}Y=\xi_i \Nb{n}Y\;.
\end{equation*}
This defines $\tilde{\xi_i}$ and shows that \eqref{eq1-pro-mult-ner} commutes.
\end{proof}
We next consider some low-dimensional examples in which the multinerves can be visualized geometrically.
\begin{example}\label{ex-mult-ner-dim2}
\emph{Double nerves of double categories.} \index{Double nerve}

\index{Double category}

Taking $\clC=\Set$ and $n=2$ above, we obtain the double nerve functor
\begin{equation*}
  \Nb{2}:\Cat(\Cat)\rw\funcat{2}{\Set}
\end{equation*}
from double categories to bisimplicial sets. Given $X\in\Cat(\Cat)$ we can visualize the corner of its double nerve as in Figure \ref{corner-double-nerve} on page \pageref{corners2}.

We  see that all the Segal maps in both the vertical and horizontal directions are isomorphisms. It is possible to give a geometric interpretation to this double nerve by thinking of $X_{00}$ as  sets of points, $X_{01}$ as sets of vertical arrows, $X_{10}$ as sets of horizontal arrows, $X_{11}$ as sets of squares. Then horizontal and vertical arrows compose (respectively horizontally and vertically) while squares compose horizontally via
\begin{equation*}
  \tens{X_{11}}{X_{01}}\rw X_{11}
\end{equation*}
and vertically via
\begin{equation*}
  \tens{X_{11}}{X_{10}}\rw X_{11}\;.
\end{equation*}
This geometric data, and the relative axioms, are an alternative way to define double categories: see for instance \cite{CEhres1963} for more details.

Although this presentation is useful for the study of some aspects of double categories, for this work the presentation in terms of internal categories in $\Cat$ and the one in terms of their double nerves is sufficient. Figure \ref{corner2} on page \pageref{corners2} is a geometric picture of the corner of the double nerve of a double category.

\end{example}

\begin{example}\label{ex-mult-ner-dim3}
\emph{3-Fold nerve of 3-fold categories.}

Taking $\clC=\Set$  and $n=3$ we obtain the 3-fold nerve of 3-fold categories
\begin{equation*}
  \Nb{3}:\cat{3}\rw\funcat{3}{\Set}\;.
\end{equation*}

Given $X\in\cat{3}$, a picture of the corner of $\Nb{3}X$ is given in Figure ~\ref{corner3X} on page \pageref{corners3}, where direction 1 is horizontal and direction 3 is vertical. We have omitted drawing the degeneracy operators for simplicity.

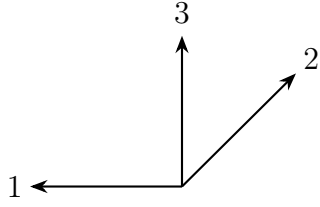
\begin{center}
\begin{tikzpicture}[>=Stealth,thick]
\node (A) at (-2.2,0) {1};
\node (B) at (0,2.3) {3};
\node (C) at (1.7,1.7) {2};
\node (OR) at (0,0) {};
\draw[->] (0,0) -- (-2,0);
\draw[->] (0,0) -- (0,2);
\draw[->] (0,0) -- (1.5,1.5);
\end{tikzpicture}
\end{center}

The isomorphisms above Figure ~\ref{corner3X} correspond to the Segal condition.

We can obtain a geometric visualisation of Figure ~\ref{corner3X} by setting

\begin{itemize}
  \item[]  $X_{000}=$ set of objects;
  \item[]  $X_{010}=$ set of arrows in direction 2;
  \item[]  $X_{100}=$ set of arrows in direction 1;
  \item[]  $X_{001}=$ set of arrows in direction 3;
  \item[]  $X_{110}=$ set of squares in directions 1,2;
  \item[]  $X_{011}=$ set of squares in directions 2,3;
  \item[]  $X_{101}=$ set of squares in directions 1,3;
  \item[]  $X_{111}=$ set of cubes.
\end{itemize}

Arrows in direction $i$ can be composed in direction $i$ (for $i=1,2,3$). Squares in direction $i,j$ ($i,j=1,2,3$) can be composed in direction $i$ and in direction $j$. Cubes can be composed in all three directions. A geometric picture is found in Figure \ref{corner3} on page \pageref{corners3}.

\end{example}

\bk
\begin{definition}\label{def-ner-func-dirk}
    Let $1 \leq k \leq n$. The nerve functor in the $k^{th}$ direction is defined as the composite
    \begin{equation*}
        \Nu{k}:\cat{n}(\clC)\xrw{\tilde{\xi_k}}\Cat(\cat{n-1}(\clC))\xrw{N}\funcat{}{\cat{n-1}(\clC)}
    \end{equation*}
    where $\tilde{\xi_k}$ is as in Proposition \ref{pro-mult-ner}.

    Note that
    \begin{equation*}
      \Nb{n}=\Nu{n}...\Nu{2}\Nu{1}.
    \end{equation*}
\end{definition}

\begin{remark}\label{rem-obj-fun}
  Let $ob : \Cat \clC \rw \clC$ be the object of objects functor. The left adjoint to $ob$ is the discrete internal category functor
  \begin{equation*}
    d:\clC \rw \Cat \clC
  \end{equation*}
  associating to an object $X\in \clC$ the discrete internal category on $X$. By Definition \ref{def-discr-nfold}, the nerve of the discrete internal category on $X$ is the constant simplicial object on $X$.

\index{Discrete!- internal category}

\end{remark}

 By Proposition \ref{pro-mult-ner} we then have
\begin{equation*}
\xymatrix{
    \cat{n}\clC \oset{\tilde{\xi_n}}{\cong}\Cat(\cat{n-1}\clC) \ar@<1ex>[r]^(0.65){ob} & \cat{n-1}\clC \ar@<1ex>[l]^(0.35){d}\;.
}
\end{equation*}

\begin{definition}\label{dn}
We define $\di{1}=d $ and for $n>1$,
\begin{equation*}
\di{n}=\tilde{\xi}^{-1}_{n}\circ d :\cat{n-1}\clC \rw \cat{n}\clC;.
\end{equation*}
\end{definition}
Thus $\di{n}$ is the discrete inclusion of $\cat{n-1}\clC$ into $\cat{n}\clC$ in the $n^{th}$ direction.

\begin{notation}\label{not-ner-func-dirk}\index{Functor!- $J_n$}
    When $\clC=\Set$ we shall denote
    \begin{equation*}
        J_n=\Nu{n-1}\ldots \Nu{1}:\cat{n}\rw\funcat{n-1}{\Cat}\;.
    \end{equation*}
\end{notation}
 Thus $J_n$ amounts to taking the nerve construction in all but the last simplicial direction. The functor $J_n$ is fully faithful, thus we can identify $\cat{n}$ with the image $J_n(\cat{n})$ of the functor $J_n$.

 Given $X\in\cat{n}$, when no ambiguity arises we shall denote, for each $(s_1,\ldots ,s_{n-1})\in\Dmenop$
\begin{equation*}
    X_{s_1,\ldots ,s_{n-1}}=(J_n X)_{s_1,\ldots ,s_{n-1}}\in\Cat
\end{equation*}
and more generally, if $1\leq j \leq n-1$,
\begin{equation*}
    X_{s_1,\ldots ,s_{j}}=(\Nu{j}\ldots \Nu{1} X)_{s_1,\ldots ,s_{j}}\in\cat{n-j}\;.
\end{equation*}

\begin{lemma}\label{lem-char-obj-III}
  There is a commuting diagram
  \begin{equation}\label{eq-char-obj-III}
  \xymatrix @C=50pt{
  \cat{n-1} \ar^{\Nb{n-1}}[r] \ar_{\di{n}}[d] & \funcat{n-1}{\Set} \ar^{\bar d}[d] \\
  \cat{n} \ar_{J_n}[r] & \funcat{n-1}{\Set}
  }
  \end{equation}
  where $\di{n}$ is as in Definition \ref{dn} and $d:\Set\rw\Cat$ is the discrete category functor.
\end{lemma}
\begin{proof}
Let $X\in\cat{n}$ and $\uk\in\dop{n}$. Then
\begin{equation*}
  (\bar d \Nb{n-1}X)_{\uk}=d(\Nb{n-1}X)_{\uk}=d X_{\uk}=(\xi_n^{-1}d_1 X)_{\uk}=(J_n \di{n}X)_{\uk}\;.
\end{equation*}
Since this holds for each $X\in\cat{n}$ and $\uk\in\dop{n}$, it follows that
\begin{equation*}
  \bar d \Nb{n-1}=J_n \di{n}
\end{equation*}
that is, \ref{eq-char-obj-III} commutes.
\end{proof}
\begin{lemma}\label{lem-char-obj-II}
    Let $P$ be the pullback in $\funcat{n-1}{\Cat}$ of the diagram in $A\rw C \lw B$. Suppose that for each $k\geq 2$ there are isomorphisms of Segal maps in $\funcat{n-2}{\Cat}$
    \begin{equation*}
        A_k\cong \pro{A_1}{A_0}{k},\;\;  C_k\cong \pro{C_1}{C_0}{k},\;\;  B_k\cong \pro{B_1}{B_0}{k}\;.
    \end{equation*}
    Then $P_k\cong \pro{P_1}{P_0}{k}$.
\end{lemma}
\begin{proof}
We show this for $k=2$, the case $k>2$ being similar. Since the nerve functor $N:\Cat\rw\funcat{}{\Set}$ commutes with pullbacks (as it is right adjoint) and pullbacks in $\funcat{n-1}{\Cat}$ are computed pointwise, for each $\us\in\dop{n-1}$ we have a pullback in $\Set$
\begin{equation*}
    \xymatrix{
    (NP)_{2\us} \ar[r] \ar[d] & (NA)_{2\us} \ar[d]\\
    (NC)_{2\us} \ar[r] & (NC)_{2\us}
    }
\end{equation*}
where
\begin{equation*}
    (NA)_{2\us}=\tens{(NA)_{1\us}}{(NA)_{0\us}}
\end{equation*}
and similarly for $NC$ and $NB$. We then calculate
\begin{align*}
    &(NP)_{2\us}=(NA)_{2\us} \tiund{(NC)_{2\us}} (NB)_{2\us} =\\
    & \resizebox{1.0\hsize}{!}{$
    =\{ \tens{(NA)_{1\us}}{(NA)_{0\us}} \}\tiund{\tens{(NC)_{1\us}}{(NC)_{0\us}}} \{ \tens{(NB)_{1\us}}{(NB)_{0\us}}\}\cong $}\\
    & \resizebox{1.0\hsize}{!}{$
    \cong\{ \tens{(NB)_{1\us}}{(NC)_{0\us}} \}\tiund{\tens{(NA)_{1\us}}{(NC)_{0\us}}} \{ \tens{(NB)_{1\us}}{(NC)_{0\us}}\}= $}\\
    &= \tens{(NP)_{1\us}}{(NP)_{0\us}}\;.
\end{align*}
In the above, the isomorphism before the last takes $(x_1,x_2,x_3,x_4)$ to $(x_1,x_3,x_2,x_4)$. Since this holds for all $\us$, it follows that
\begin{equation*}
    P_2\cong \tens{P_1}{P_0}\;.
\end{equation*}
The case $k>2$ is similar.

\end{proof}

%%%%%%%%%%%%%%%%%%
\section{Multisimplicial description of strict $\pmb{n}$-categories}\label{sbs-multi-strict}
%\begin{definition}\label{def-discrete-nfold}
%An internal $n$-fold category $X\in \cat{n}(\clC)$ is said to be discrete if $\Nn X$ is a constant %functor.
%\end{definition}

In this section we establish the analogue of Lemma \ref{lem-multin-iff} for strict $n$-categories. We recall the definition of strict $n$-category, which is given by iterated enrichment.\index{Enriched category} We refer to \cite{Kelly1982} for background on enriched categories.

\begin{definition}\label{def-strict}\index{Strict $n$-category}
  The category $n\mi\Cat$ of strict $n$-categories is defined by induction of $n$. For $n=1$,  $1\mi\Cat=\Cat$. Given $(n-1)\mi\Cat$, let
  \begin{equation*}
    n\mi\Cat=((n-1)\mi\Cat)\mi\Cat
  \end{equation*}
  where the enrichment is with respect to the cartesian monoidal structure.
\end{definition}

Unravelling this inductive definition we see that strict $n$-categories have sets of cells in dimensions $0$ up to $n$, and these compose in a way which is associative and unital. Below we make this more precise using the multi-simplicial language.

 A strict $n$-category in which all cells are invertible is called a strict $n$-groupoid. \index{Strict $n$-groupoid} We denote by $n\mi\Gpd$ the category of strict $n$-groupoids, where $1\mi\Gpd=\Gpd$.

 We prove a characterization of the image of the functors
\begin{equation*}
  J_n: n\mi\Cat \rw\funcat{n-1}{\Cat} \quad \mbox{and}\quad \Nn:n\mi\Cat \rw\funcat{n}{\Set}\;.
\end{equation*}
The resulting multi-simplicial description of strict $n$-categories is helpful to build the intuition about the weakening of the structure achieved with the Segal-type models of this work, which are also full subcategories of $\funcat{n-1}{\Cat}$.

\index{Enriched category}To establish this characterization, we use the fact, proved in the Appendix of \cite{ACEhresIII} that if $\clV$ is a category satisfying mild conditions, the categories enriched in $\clV$ are the internal categories in $\clV$ whose object of objects is discrete, that is it is the coproduct of copies of the terminal object. As observed in \cite{ACEhresIII}, these conditions are satisfied when $\clV=\cat{n}$ and $\clV=n\mi\Cat$.
\begin{remark}\label{rem-multi-strict}
Consider the composite
\begin{equation*}
  N:\clV\mi\Cat \rw \Cat \clV \rw \funcat{}{\clV}
\end{equation*}
where $\clV$ is as in \cite{ACEhresIII}. It follows by \cite{ACEhresIII} and by Proposition \ref{pro-ner-int-cat} that $X\in\funcat{}{\clV}$ is in the image of $N$ if and only if the Segal maps of $X$ are isomorphisms and $X_0$ is discrete.
\end{remark}

\begin{figure}[ht]
  \centering
\begin{equation*}
\entrymodifiers={+++[o]}
 \xymatrix@R=35pt@C=25pt{
   \quad \cdots \quad \ar@<2ex>[r] \ar[r] \ar@<1ex>[r] \ar@<-2ex>[d] \ar[d] \ar@<-1ex>[d]&
   \tens{X_{11}}{X_{10}} \ar@<1ex>[r] \ar@<0ex>[r]  \ar@<-2ex>[d] \ar[d] \ar@<-1ex>[d] \ar@<1ex>[l] \ar@<2ex>[l]   &
    X_{00}  \ar@<-2ex>[d] \ar[d] \ar@<-1ex>[d] \ar@<1ex>[l]\\
    \tens{X_{11}}{X_{00}}\ar@<2ex>[r] \ar[r] \ar@<1ex>[r] \ar@<-1ex>[d] \ar@<0.ex>[d] \ar@<-1ex>[u] \ar@<-2.ex>[u] &
    X_{11} \ar@<1ex>[r] \ar@<0ex>[r]  \ar@<-1ex>[d] \ar@<0.ex>[d] \ar@<1ex>[l] \ar@<2ex>[l] \ar@<-1ex>[u] \ar@<-2.ex>[u] &
    X_{00}  \ar@<0.ex>[d] \ar@<-1ex>[d] \ar@<1ex>[l] \ar@<-1ex>[u] \ar@<-2.ex>[u] \\
    \tens{X_{10}}{X_{00}} \ar@<2ex>[r] \ar[r] \ar@<1ex>[r] \ar@<-1ex>[u] &
     X_{10} \ar@<1ex>[r] \ar@<0ex>[r] \ar@<1ex>[l] \ar@<2ex>[l] \ar@<-1ex>[u] &
    X_{00}\ar@<1ex>[l] \ar@<-1ex>[u]\\
}
\end{equation*}
  \caption{Corner of the double nerve of a strict 2-category $X$}
  \label{corner-double-strict}
\end{figure}
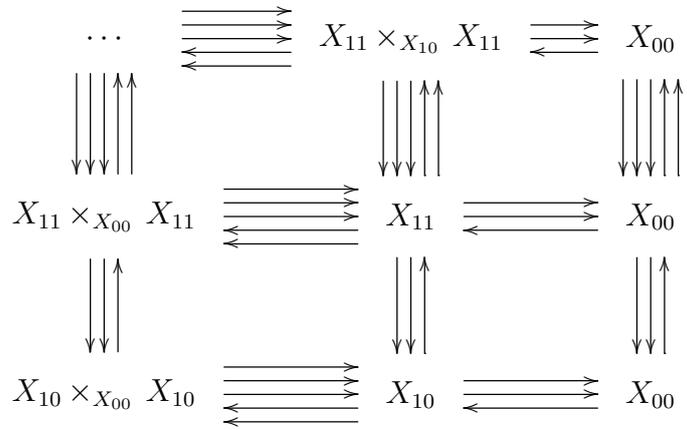

\begin{figure}[ht]
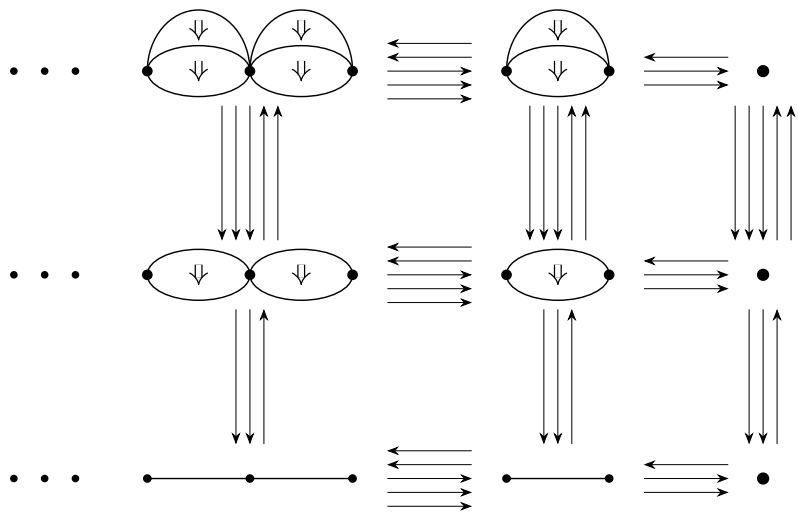

  \centering
  \includestandalone[width=1.0\textwidth]{Cells2}
  \caption{Geometric picture of the corner of the double nerve of a strict $2$-category.}
  \label{cells2}
\end{figure}

\begin{figure}[ht]
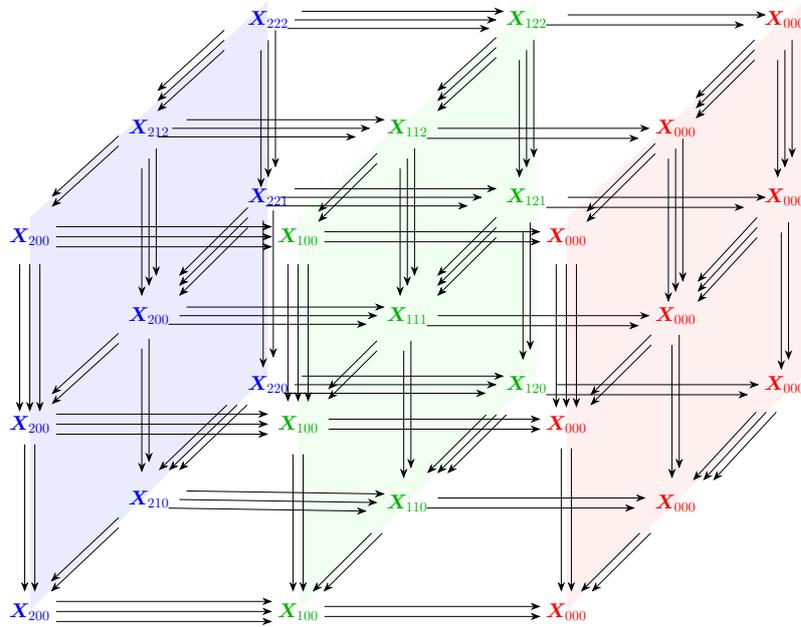

  \centering
  \includestandalone[width=0.9\textwidth]{Corner3XGLOB}
  \caption{Corner of the $3$-fold nerve of a strict $3$-category $X$.}
  \label{corner3XGLOB}
\end{figure}

\begin{figure}[ht]
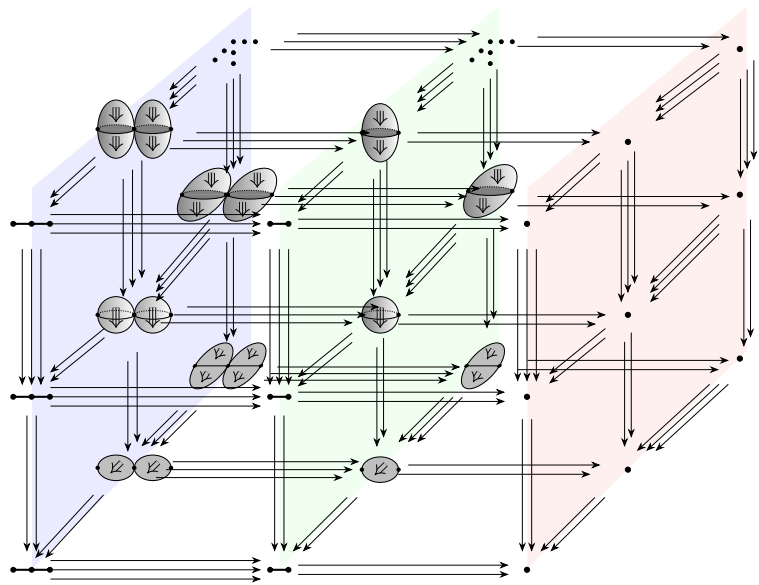

  \centering
  \includestandalone[width=0.9\textwidth]{Globular3}
  \caption{Geometric picture of the corner of the $3$-fold nerve of a strict $3$-category.}
  \label{corner-glob}
\end{figure}

When $\clV=(n-1)\mi\Cat$, we obtain the functor
\begin{equation*}
N:n\mi\Cat=((n-1)\mi\Cat)\mi\Cat \rw \funcat{}{(n-1)\mi\Cat}\;.
\end{equation*}
Iterating this construction we obtain the multinerve
\begin{equation}\label{jn-str}
  J_n: n\mi\Cat \rw \funcat{n-1}{\Cat}\;.
\end{equation}
We next give a characterization of the image of $J_n$.
\begin{lemma}\label{lem-1-multi-strict}
  Let $X\in\funcat{n-1}{\Cat}$. Then $X$ is in the image of the functor $J_n$ as in \eqref{jn-str} if and only if:
  \begin{itemize}
    \item [a)] The Segal maps in all directions are isomorphisms.\bk

    \item [b)] $X_0\in\funcat{n-2}{\Cat}$ and $X_{\oset{r}{1...1}0}\in \funcat{n-r-2}{\Cat }$ are constant functors taking value in a discrete category for all $1\leq r \leq n-2$.
  \end{itemize}
\end{lemma}
\begin{proof}
By induction on $n$. For $n=2$ it follows by Remark \ref{rem-multi-strict} taking $\clV=\Cat$. Suppose, inductively, that the lemma holds for $(n-1)$ and let $X=J_n Y$ with $Y\in n\mi\Cat$. By construction, $J_n$ is the composite:
\begin{equation*}
  n\mi\Cat \xrw{N}\funcat{}{(n-1)\mi\Cat}\xrw{\ovl{J}_{n-1}} \funcat{}{\funcat{n-2}{\Cat}}\cong \funcat{n-1}{\Cat}\;.
\end{equation*}
By Remark \ref{rem-multi-strict} when $\clV=(n-1)\mi\Cat$, $(N Y)_0$ is a constant functor taking value in a discrete and the Segal maps of $N Y$ are isomorphisms. Hence $X_0=(J_n Y)_0=J_{n-1}(NY)_0$ is constant with values in a discrete category and the Segal maps of $X$ in direction 1 are isomorphisms.

The Segal maps of $J_{n}X$ in direction $k>1$ are levelwise the Segal maps of $X_t=J_{n-1}(NY)_t$ for each $t\geq 0$ and these are isomorphisms by inductive hypothesis applied to $(NY)_t\in(n-1)\mi\Cat$. So in conclusion the Segal maps in all directions are isomorphisms, which is a).

Further $X_1=J_{n-1}(NY)_1$ and since $(NY)_1\in (n-1)\mi\Cat$, by inductive hypothesis $X_{10}$ and $X_{\oset{s}{1...1}0}$ are constant with values in a discrete category for all $1\leq s\leq n-3$. So in conclusion $X_{\oset{r}{1...1}0}$ is constant with values in a discrete category for all $1\leq r\leq n-2$, proving that condition b) holds.

Conversely, let $X\in\funcat{n-1}{\Cat}$ satisfy a) and b). By Lemma \ref{lem-multin-iff} it follows from condition a) that $X=J_n Y$ for $Y\in\cat{n}$. Further, for each $t\geq 0$, $(NY)_t\in \cat{n-1}$ and $J_{n-1}(NY)_t$ satisfies a) and b), so by induction hypothesis $(NY)_t\in (n-1)\mi\Cat$. Also, since by b) $X_0$ is constant with values in a discrete category, $(NY)_0$ is discrete.

In conclusion, $NY\in\funcat{}{(n-1)\mi\Cat}$ is such that $(NY)_0$ is discrete and the Segal maps are isomorphisms. From Remark  \ref{rem-multi-strict}  it follows that $Y\in n\mi\Cat$.
\end{proof}
By taking nerves of categories dimensionwise in $\funcat{n-1}{\Cat}$ we obtain the multinerve functor:
\begin{equation*}
  \Nn : n\mi\Cat \xrw{} \funcat{n-1}{\Cat}\xrw{\ovl{N}} \funcat{n-1}{\funcat{}{\Set}}\cong \funcat{n}{\Set}\;.
\end{equation*}
Using Lemma \ref{lem-1-multi-strict} we immediately deduce the following characterization of the image of $\Nn$, which affords a multi-simplicial description of strict $n$-categories.
\begin{corollary}\label{cor-multi-strict}
Let $X\in\funcat{n}{\Set}$. Then $X$ is in the image of the functor
\begin{equation*}
  \Nn : n\mi\Cat \rw \funcat{n}{\Set}
\end{equation*}
if and only if
\begin{itemize}
  \item [a)] The Segal maps of $X$ in all directions are isomorphisms.\bk

  \item [b)] $X_0\in\funcat{n-1}{\Set}$ and $X_{\oset{r}{1...1}0}\in\funcat{n-r-1}{\Set}$ are constant functors for all $1\leq r\leq n-2$.
\end{itemize}
\end{corollary}
\begin{proof}
It follows immediately from Lemma \ref{lem-1-multi-strict}, Proposition \ref{pro-ner-int-cat} and the fact that the nerve of a discrete category is a constant simplicial set.
\end{proof}
The condition b) in Lemma \ref{lem-1-multi-strict} and Corollary \ref{cor-multi-strict} is called \emph{globularity condition}, since it gives rise to the globular shape of the higher cells in a strict $n$-category. We illustrate this in low-dimensional examples.
\begin{example}\label{ex-1-multi-strict}
\emph{Double nerves of strict 2-categories.}
The double nerve functor
\begin{equation*}
  \Nb{2} : 2\mi\Cat\rw \funcat{2}{\Set}
\end{equation*}
associates to $X\in 2\mi\Cat$ the bisimplicial set $\Nb{2} X$ such that the Segal maps in both horizontal and vertical directions are isomorphisms and $(\Nb{2}X)_{0*}$ is a constant functor.

Denoting $(\Nb{2}X)_{ij}=X_{ij}$, we can visualize the corner of its double nerve as in Figure \ref{corner-double-strict} on page \pageref{corner-double-strict-p}.
%
%\begin{equation*}
%\xymatrix@R=5pt @C=90pt{\hspace{30pt}  \vdots & \hspace{10pt}\vdots & \!\!\vdots}
%\end{equation*}

%
A geometric picture can be obtained by thinking of $X_{00}$ as sets of objects or 0-cells, $X_{10}$ as sets of 1-cells and $X_{20}$ as sets of 2-cells, see Figure \ref{cells2} on page \pageref{corner-double-strict-p}.
\bk

This should be compared with Figure \ref{corner2} in Example \ref{ex-mult-ner-dim2}: the globularity condition that $X_{0*}$ is a constant simplicial set means that the vertical sides of the squares in  Figure \ref{corner2} are identities and thus can be represented as globes under the identification
%
%% Create 1-edge
\tikzset{edge/.pic={
\filldraw (0,0)  circle[radius=0.035cm] -- (1,0)  circle[radius=0.035cm]; %% centered
}}
%% Create  1-cell
\tikzset{Cell1/.pic={
\draw(0,0) ellipse [x radius=0.5cm, y radius=0.30cm];
\doublearrow{arrows={-Implies}}{(0.0,0.1) -- (0.0,-0.1)};
\draw plot [mark=*] coordinates {(-0.5,0)};
\draw plot [mark=*] coordinates {(0.5,0)};
}}
\begin{center}
\begin{tikzpicture}[thick,scale=2]]
\pic [black][scale=2] at (0,0) {edge};
\pic [black][scale=2] at (0,1) {edge};
\draw[thick,double distance=2pt] (0,0.2) -- (0,0.8);
\draw[thick,double distance=2pt] (1,0.2) -- (1,0.8);
\doublearrow{arrows={-Implies}}{(0.5,0.6) -- (0.5,0.4)};
\pic [black][scale=2] at (3,0.5) {Cell1};
\node (1) at (0,1.2) {a};
\node (2) at (1,1.2) {b};
\node (3) at (0,-0.2) {a};
\node (4) at (1,-0.2) {b};
\node (5) at (0.5,1.1) {$g$};
\node (6) at (0.5,-0.12) {$f$};
\node (7) at (1.7,0.5) {$\equiv$};
\node (8) at (2.3,0.5) {a};
\node (9) at (3.7,0.5) {b};
\node (10) at (3,0.9) {$g$};
\node (10) at (3,0.05) {$f$};
\end{tikzpicture}
\end{center}
%% End of TIKX picture
\end{example}
\begin{example}\label{ex-2-multi-strict}
\emph{3-Fold nerves of strict 3-categories.}\index{Strict $3$-category}
The 3-fold nerve functor
\begin{equation*}
  \Nb{3}: 3\mi\cat \;\rw \funcat{3}{\Set}
\end{equation*}
associates to $Y\in 3\mi\Cat$ the 3-fold simplicial set $\Nb{3}Y=X$ such that
\begin{itemize}
  \item [a)] All Segal maps in directions 1,2,3 are isomorphisms.\bk

  \item [b)] $X_{0**}$ and $X_{10*}$ are constant functors.\bk
\end{itemize}
The corner of $X$ can be visualized as in Figure \ref{corner3XGLOB} on page \pageref{globular3}, where for simplicity we omitted drawing the face operators. The isomorphisms above Figure \ref{corner3XGLOB} correspond to the Segal condition.
\bk

We can obtain a geometric visualization of Figure \ref{corner3XGLOB} by setting
\begin{itemize}
  \item[]  $X_{000}=$ set of points;
  \item[]  $X_{100}=$ set of arrows in direction 1;
  \item[]  $X_{010}=$ set of arrows in direction 2;
  \item[]  $X_{001}=$ set of arrows in direction 3;
  \item[]  $X_{110}=$ set of globes in directions 1,2;
  \item[]  $X_{011}=$ set of globes in directions 2,3;
  \item[]  $X_{101}=$ set of globes in directions 1,3;
  \item[]  $X_{111}=$ set of spheres.
\end{itemize}

The points are the 0-cells. The arrows (1-cells) in direction $i$ can be composed in direction $i$ ($i=1,2,3$). The globes (2-cells) in direction $i,j$ ($i,j=1,...,3$) can be composed in direction $i$ and in direction $j$. The spheres (3-cells) can be composed in all three directions. A geometric picture is found in Figure \ref{corner-glob} on page \pageref{globular3}.

This should be compared with Figure \ref{corner3} in Example \ref{ex-mult-ner-dim3}: the globularity condition identifies the squares with globes and the cubes with spheres.

\end{example}

\section{The functor d\'{e}calage}\label{decalage}

\index{D\'{e}calage functor}

Recall from Duskin \cite{Dusk} the d\'{e}calage comonad
\begin{equation*}
  \Dec: \funcat{}{\Set} \rw \funcat{}{\Set}
\end{equation*}
on simplicial sets. Given $X\in\funcat{}{\Set}$, $\Dec X\in\funcat{}{\Set}$ has
\begin{equation*}
  (\Dec X)_n = X_{n+1} \qquad n\geq 0
\end{equation*}
with face and degeneracy given by
\begin{equation*}
\begin{split}
    & d_i: (\Dec X)_n = X_{n+1}\rw (\Dec X)_{n-1} = X_n, \qquad 1\leq i\leq n \\
    & s_i: (\Dec X)_n = X_{n+1}\rw (\Dec X)_{n+1} = X_{n+2}, \qquad 0\leq i\leq n\;.
\end{split}
\end{equation*}
That is, the last face and degeneracy operators are omitted in each dimension. The omitted last face map $d_n: X_n \rw X_{n-1}$ defines a simplicial map
\begin{equation*}
  u_X: \Dec X\rw X
\end{equation*}
natural in $X$, levelwise surjective, which is the component of the counit of the comonad.

The composition of $d_0:X_1 \rw d(X_0)$ with the retained face maps gives a simplicial map
\begin{equation*}
  \Dec X\rw d(X_0)
\end{equation*}
natural in $X$, where $d(X_0)$ is the constant simplicial set at $X_0$. Composition of $s_0: X_0\rw X_1$ with the degeneracies gives a unique map $s\up{n}: X_0\rw X_n$ which gives a simplicial map
\begin{equation*}
  d(X_0) \rw \Dec X
\end{equation*}
natural in $X$ that is a section for $\Dec X\rw d(X_0)$; it is in fact a contracting homotopy, so $d(X_0)$ is a deformation retract of $\Dec X$.

Suppose that $X$ is the nerve of a groupoid. Then $\Dec X$ is the equivalence relation corresponding to the surjective map of sets $d_0: X_1\rw X_0$ (see \cite{Bourne1987}). If $\tens{X_1}{X_0}$ is the pullback of $X_1\xrw{\pt_0} X_0 \xlw{\pt_1} X_1$ and $\tens{X_1}{d_0}$ is the kernel pair of $d_0: X_1\rw X_1$, in this case there is an isomorphism
\begin{equation*}
  \tens{X_1}{X_0} \cong \tens{X_1}{d_0}
\end{equation*}
sending $(f,g)\in \tens{X_1}{X_0}$ to $(f,g\cirsm f)\in \tens{X_1}{d_0}$; $\pt_0,\pt_1$ correspond to the two projections $p_1, p_2: \tens{X_1}{d_0}\rw X_1$ while the identity map $X_1\rw \tens{X_1}{X_0}$ corresponds to the diagonal map $X_1\rw \tens{X_1}{d_0}$  sending $f$ to $(f,f)$.

There is also a version of the d\'{e}calage comonad forgetting the first face and degeneracy operators, which we denote

\begin{equation*}
  \Dec': \funcat{}{\Set} \rw \funcat{}{\Set}.
\end{equation*}

%%%%%%%%%%%%%%%%%%%%%%%%%%%%%%%%%%%%%%%%%%%%%%%%%%%%%%%%%%%%%%%%%%%%%%%%%%%%
\clearpage

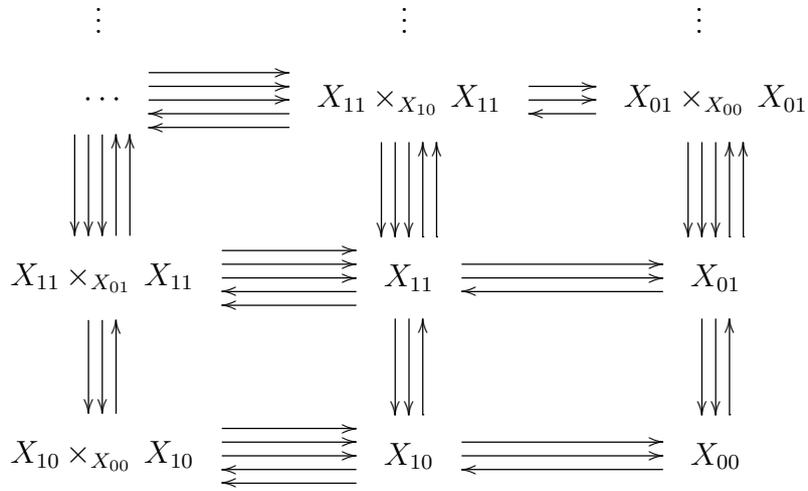
\begin{figure}[ht]
  \centering

\begin{equation*}
\xymatrix@R=5pt @C=100pt{ \!\!\!\!\!\! \vdots & \vdots & \vdots}
\end{equation*}
\begin{equation*}
\entrymodifiers={+++[o]}
 \xymatrix@R=35pt@C=25pt{
   \cdots \ar@<2ex>[r] \ar[r] \ar@<1ex>[r] \ar@<-2ex>[d] \ar[d]
    \ar@<-1ex>[d]& \tens{X_{11}}{X_{10}} \ar@<1ex>[r] \ar@<0ex>[r]  \ar@<-2ex>[d] \ar[d] \ar@<-1ex>[d] \ar@<1ex>[l] \ar@<2ex>[l] &
    \tens{X_{01}}{X_{00}}  \ar@<-2ex>[d] \ar[d] \ar@<-1ex>[d] \ar@<1ex>[l]\\
    \tens{X_{11}}{X_{01}}\ar@<2ex>[r] \ar[r] \ar@<1ex>[r] \ar@<-1ex>[d] \ar@<0.ex>[d] \ar@<-1ex>[u] \ar@<-2.ex>[u] &
    X_{11} \ar@<1ex>[r] \ar@<0ex>[r]  \ar@<-1ex>[d] \ar@<0.ex>[d] \ar@<1ex>[l] \ar@<2ex>[l] \ar@<-1ex>[u] \ar@<-2.ex>[u] &
    X_{01}  \ar@<0.ex>[d] \ar@<-1ex>[d] \ar@<1ex>[l] \ar@<-1ex>[u] \ar@<-2.ex>[u] \\
    \tens{X_{10}}{X_{00}} \ar@<2ex>[r] \ar[r] \ar@<1ex>[r] \ar@<-1ex>[u] &
     X_{10} \ar@<1ex>[r] \ar@<0ex>[r] \ar@<1ex>[l] \ar@<2ex>[l] \ar@<-1ex>[u] &
    X_{00}\ar@<1ex>[l] \ar@<-1ex>[u]\\
}
\end{equation*}
  \caption{Corner of the double nerve of a double category $X$.}
  \label{corner-double-nerve}
\end{figure}

\bk

\bk
%% Enter a TIKX picture
%%%%%%%%%%%%%%%%%%%%%%%%%%%%%%%%%%%%%%%%%%%%%%%%%%%%%%%%%%%%%%%
\begin{figure}[ht]
  \centering
  %% Begin a TIKX picture
%%%%%%%%%%%%%%%%%%%%%%%%%%%%%%%%%%%%%%%%%%%%%%%%%%%%%%%%%%%%%%%
%% Create objects \tickzset
%% Create single edge
\tikzset{Edge/.pic={
\draw plot [mark=*] coordinates {(0,0) (0.5,0)};
}}

%% Create a face
\tikzset{Face/.pic={
\draw plot [mark=*] coordinates {(0,0) (0.5,0)  (0.5,0.5)  (0,0.5) (0,0)};
\node[rotate=-90,scale=0.7] at (0.25,0.25) () {$\Rightarrow$};
}}
%%
%% Create 3 dots
\tikzset{Dots/.pic={
\filldraw [black]
(0,0) circle [radius=1.0pt]
(0.3,0) circle [radius=1.0pt]
(0.6,0) circle [radius=1.0pt];
}}

\begin{tikzpicture}[thick]
% Objects in row 1 (bottom)
\draw[fill=red] (0,0)  circle(0.05)[black]; %% (nodo 18)
\pic [black] at (-2,0) {Edge};
\pic [black] at (-4.0,0) {Edge};\pic [black] at (-4.5,0) {Edge};
% Object in row 2 (center)
\pic [black,rotate=-90] at (0,2) [xshift=0cm] {Edge};
\pic [black] at (-2,1.5) {Face};
\pic [black] at (-4.0,1.5) {Face};\pic [black] at (-4.5,1.5) {Face};

% Objects in row 3 (upper)
\pic [black,rotate=-90] at (0,4.5) [xshift=0cm] {Edge};
\pic [black,rotate=-90] at (0,4.0) [xshift=0cm] {Edge};
\pic [black] at (-2,3.5) {Face};\pic [black] at (-2,4.0) {Face};
\pic [black] at (-4.0,3.5) {Face};\pic [black] at (-4.0,4.0) {Face};
\pic [black] at (-4.5,3.5) {Face};\pic [black] at (-4.5,4.0) {Face};
\pic [black] at (-6,0) {Dots};
\pic [black] at (-6,1.75) {Dots};
\pic [black] at (-6,4) {Dots};
%% Positioning nodes
\node (1)[minimum size=10mm] at (0,0) {};
\node (2)[minimum size=10mm] at (-1.75,0) {};
\node (3)[minimum size=10mm] at (-4,0) {};
\node (4)[minimum size=10mm] at (0,1.75) {};
\node (5)[minimum size=10mm] at (-1.75,1.75) {};
\node (6)[minimum size=10mm] at (-4,1.75) {};
\node (7)[minimum size=10mm] at (0,4) {};
\node (8)[minimum size=10mm] at (-1.75,4) {};
\node (9)[minimum size=10mm] at (-4,4) {};

%% Arrows
%% 1-2
\draw[transform canvas={yshift=1.0ex,xshift=0ex},->][arrows = {-Stealth}][thin] (2) -- (1);
\draw[transform canvas={yshift=0.0ex,xshift=0ex},->][arrows = {-Stealth}][thin] (2) -- (1);
\draw[transform canvas={yshift=-1.0ex,xshift=0ex},->][arrows = {-Stealth}][thin] (1) -- (2);
%% 4-5
\draw[transform canvas={yshift=1.0ex,xshift=0ex},->][arrows = {-Stealth}][thin] (5) -- (4);
\draw[transform canvas={yshift=0.0ex,xshift=0ex},->][arrows = {-Stealth}][thin] (5) -- (4);
\draw[transform canvas={yshift=-1.0ex,xshift=0ex},->][arrows = {-Stealth}][thin] (4) -- (5);
%% 7-8
\draw[transform canvas={yshift=1.0ex,xshift=0ex},->][arrows = {-Stealth}][thin] (8) -- (7);
\draw[transform canvas={yshift=0.0ex,xshift=0ex},->][arrows = {-Stealth}][thin] (8) -- (7);
\draw[transform canvas={yshift=-1.0ex,xshift=0ex},->][arrows = {-Stealth}][thin] (7) -- (8);
%% 2-3
\draw[transform canvas={yshift=2.0ex,xshift=1ex},->][arrows = {-Stealth},shorten >=0.2cm, shorten <=0.2cm][thin] (3) -- (2);
\draw[transform canvas={yshift=1.0ex,xshift=1ex},->][arrows = {-Stealth},shorten >=0.2cm, shorten <=0.2cm][thin] (3) -- (2);
\draw[transform canvas={yshift=0.0ex,xshift=1ex},->][arrows = {-Stealth},shorten >=0.2cm, shorten <=0.2cm][thin] (3) -- (2);
\draw[transform canvas={yshift=-1.0ex,xshift=1ex},->][arrows = {-Stealth},shorten >=0.2cm, shorten <=0.2cm][thin] (2) -- (3);
\draw[transform canvas={yshift=-2.0ex,xshift=1ex},->][arrows = {-Stealth},shorten >=0.2cm, shorten <=0.2cm][thin] (2) -- (3);
%% 5-6
\draw[transform canvas={yshift=2.0ex,xshift=1ex},->][arrows = {-Stealth},shorten >=0.2cm, shorten <=0.2cm][thin] (6) -- (5);
\draw[transform canvas={yshift=1.0ex,xshift=1ex},->][arrows = {-Stealth},shorten >=0.2cm, shorten <=0.2cm][thin] (6) -- (5);
\draw[transform canvas={yshift=0.0ex,xshift=1ex},->][arrows = {-Stealth},shorten >=0.2cm, shorten <=0.2cm][thin] (6) -- (5);
\draw[transform canvas={yshift=-1.0ex,xshift=1ex},->][arrows = {-Stealth},shorten >=0.2cm, shorten <=0.2cm][thin] (5) -- (6);
\draw[transform canvas={yshift=-2.0ex,xshift=1ex},->][arrows = {-Stealth},shorten >=0.2cm, shorten <=0.2cm][thin] (5) -- (6);
%% 8-9
\draw[transform canvas={yshift=2.0ex,xshift=1ex},->][arrows = {-Stealth},shorten >=0.2cm, shorten <=0.2cm][thin] (9) -- (8);
\draw[transform canvas={yshift=1.0ex,xshift=1ex},->][arrows = {-Stealth},shorten >=0.2cm, shorten <=0.2cm][thin] (9) -- (8);
\draw[transform canvas={yshift=0.0ex,xshift=1ex},->][arrows = {-Stealth},shorten >=0.2cm, shorten <=0.2cm][thin] (9) -- (8);
\draw[transform canvas={yshift=-1.0ex,xshift=1ex},->][arrows = {-Stealth},shorten >=0.2cm, shorten <=0.2cm][thin] (8) -- (9);
\draw[transform canvas={yshift=-2.0ex,xshift=1ex},->][arrows = {-Stealth},shorten >=0.2cm, shorten <=0.2cm][thin] (8) -- (9);
%% 7-4
\draw[transform canvas={xshift=-2.0ex,yshift=0ex},->][arrows = {-Stealth},shorten >=0.2cm, shorten <=0.2cm][thin] (7) -- (4);
\draw[transform canvas={xshift=-1.0ex,yshift=0ex},->][arrows = {-Stealth},shorten >=0.2cm, shorten <=0.2cm][thin] (7) -- (4);
\draw[transform canvas={xshift=0.0ex,yshift=0ex},->][arrows = {-Stealth},shorten >=0.2cm, shorten <=0.2cm][thin] (7) -- (4);
\draw[transform canvas={xshift=1.0ex,yshift=0ex},->][arrows = {-Stealth},shorten >=0.2cm, shorten <=0.2cm][thin] (4) -- (7);
\draw[transform canvas={xshift=2.0ex,yshift=0ex},->][arrows = {-Stealth},shorten >=0.2cm, shorten <=0.2cm][thin] (4) -- (7);
%% 7-4
\draw[transform canvas={xshift=-2.0ex,yshift=0ex},->][arrows = {-Stealth},shorten >=0.2cm, shorten <=0.2cm][thin] (8) -- (5);
\draw[transform canvas={xshift=-1.0ex,yshift=0ex},->][arrows = {-Stealth},shorten >=0.2cm, shorten <=0.2cm][thin] (8) -- (5);
\draw[transform canvas={xshift=0.0ex,yshift=0ex},->][arrows = {-Stealth},shorten >=0.2cm, shorten <=0.2cm][thin] (8) -- (5);
\draw[transform canvas={xshift=1.0ex,yshift=0ex},->][arrows = {-Stealth},shorten >=0.2cm, shorten <=0.2cm][thin] (5) -- (8);
\draw[transform canvas={xshift=2.0ex,yshift=0ex},->][arrows = {-Stealth},shorten >=0.2cm, shorten <=0.2cm][thin] (5) -- (8);
%% 9-6
\draw[transform canvas={xshift=-2.0ex,yshift=0ex},->][arrows = {-Stealth},shorten >=0.2cm, shorten <=0.2cm][thin] (9) -- (6);
\draw[transform canvas={xshift=-1.0ex,yshift=0ex},->][arrows = {-Stealth},shorten >=0.2cm, shorten <=0.2cm][thin] (9) -- (6);
\draw[transform canvas={xshift=0.0ex,yshift=0ex},->][arrows = {-Stealth},shorten >=0.2cm, shorten <=0.2cm][thin] (9) -- (6);
\draw[transform canvas={xshift=1.0ex,yshift=0ex},->][arrows = {-Stealth},shorten >=0.2cm, shorten <=0.2cm][thin] (6) -- (9);
\draw[transform canvas={xshift=2.0ex,yshift=0ex},->][arrows = {-Stealth},shorten >=0.2cm, shorten <=0.2cm][thin] (6) -- (9);
%% 4-1
\draw[transform canvas={xshift=-1.0ex,yshift=0ex},->][arrows = {-Stealth},shorten >=0.cm, shorten <=0.cm][thin] (4) -- (1);
\draw[transform canvas={xshift=0.0ex,yshift=0ex},->][arrows = {-Stealth},shorten >=0.cm, shorten <=0.cm][thin] (4) -- (1);
\draw[transform canvas={xshift=1.0ex,yshift=0ex},->][arrows = {-Stealth},shorten >=0.cm, shorten <=0.cm][thin] (1) -- (4);
%% 5-2
\draw[transform canvas={xshift=-1.0ex,yshift=0ex},->][arrows = {-Stealth},shorten >=0.cm, shorten <=0.cm][thin] (5) -- (2);
\draw[transform canvas={xshift=0.0ex,yshift=0ex},->][arrows = {-Stealth},shorten >=0.cm, shorten <=0.cm][thin] (5) -- (2);
\draw[transform canvas={xshift=1.0ex,yshift=0ex},->][arrows = {-Stealth},shorten >=0.cm, shorten <=0.cm][thin] (2) -- (5);
%% 6-3
\draw[transform canvas={xshift=-1.0ex,yshift=0ex},->][arrows = {-Stealth},shorten >=0.cm, shorten <=0.cm][thin] (6) -- (3);
\draw[transform canvas={xshift=0.0ex,yshift=0ex},->][arrows = {-Stealth},shorten >=0.cm, shorten <=0.cm][thin] (6) -- (3);
\draw[transform canvas={xshift=1.0ex,yshift=0ex},->][arrows = {-Stealth},shorten >=0.cm, shorten <=0.cm][thin] (3) -- (6);
\end{tikzpicture}
  \caption{Geometric picture of the corner of the double nerve of a double category}
  \label{corner2}
\end{figure}

\label{corners2}

%%%%%%%%%%%%%%%%%%%%%%%%%%%%%%%%%%%%%%%%%%%%%%%%%%%%%%%%%%%%%%%%%%%%%
\clearpage
In the following picture, for all $i,j,k \in\Delta^{op}$

$X_{2jk}\cong \tens{X_{1jk}}{X_{0jk}}, \;X_{i2k}\cong \tens{X_{i1k}}{X_{i0k}},\;X_{ij2}\cong \tens{X_{ij1}}{X_{ij0}}$\:.

\begin{figure}[ht]
  \centering
  \includestandalone[width=0.9\textwidth]{Corner3X}
  \caption{Corner of the 3-fold nerve of a 3-fold category $X$}
  \label{corner3X}
\end{figure}
\vspace{-2mm}

\mk

\begin{figure}[ht]
  \centering
  \includestandalone[width=0.9\textwidth]{Corner3}
  \caption{Geometric picture of the corner of the 3-fold nerve of a 3-fold category}
  \label{corner3}
\end{figure}

\label{corners3}

%%%%%%%%%%%%%%%%%%%%%%%%%%%%%%%%%%%%%%%%%%%%%%%%%%%%%%%%%%%%%%%%%%
\clearpage

\bk

\bk

\label{corner-double-strict-p}

\clearpage
%%%%%%%%%%%%%%%%%%%%%%%%%%%%%%%%%%%%%%%%%%%%%%%%%%%%%%%%%%%%%%%%%%%%%%%%%
%%

In the following picture, for all $i,j,k \in\Delta^{op}$

$X_{2jk}\cong \tens{X_{1jk}}{X_{0jk}}, \;X_{i2k}\cong \tens{X_{i1k}}{X_{i0k}},\;X_{ij2}\cong \tens{X_{ij1}}{X_{ij0}}$\:.

\vspace{-2mm}

\bk

\label{globular3}

\clearpage

%%%%%%%%%%%%%%%%%%%%%%%%%%%%%%%%%%%%%%%%%%%%%%%%%%%%%%%%%%%%%%%%%%%%%%%%%%%%

\chapter{Categorical background}\label{chap-ps}

In this chapter we cover some categorical background which is needed in the rest of this work. We first recall two important functors from $\Cat$ to $\Set$, associating to a category $X$ the sets of isomorphism classes of objects $p(X)$ and the set of connected components $q(X)$. Their properties are recalled in Section \ref{sbs-funct-cat}. They will give rise to functors
\begin{equation*}
        \p{n}:\seg{n}\rw\seg{n-1} \qquad \q{n}:\seg{n}\rw\seg{n-1}
      \end{equation*}
for each of the three Segal-type models. These functors are crucial to this work. The functor $\p{n}$ is a 'truncation functor' that divides out by the highest dimensional invertible cells, and is used to give the notion of $n$-equivalence in $\seg{n}$, while the functor $\q{n}$ plays an important role in Chapter \ref{chap7} in the construction of the rigidification functor $Q_n$.

The second set of techniques we review in this chapter is about the theory of pseudo-functors. Although they are used in homotopy theory (see for instance \cite{Thomas}), the theory of pseudo-functors mostly developed within category theory, see for instance \cite{Borc}. A crucial technique in this work is the use of the strictification functor from pseudo-functors to strict functors. Several versions of this exist in the literature (see for instance \cite{Str}), and we review here the one due to Power \cite{PW} and refined by Lack \cite{Lack}, as these are the forms most suitable for our calculations.

As outlined in Chapter \ref{chap1} we will show that the rigidification functor from weakly globular Tamsamani $n$-categories to weakly globular $n$-fold categories factors through a subcategory of pseudo-functors (called Segalic pseudo-functors) and the Power-Lack form of the strictification identifies the image as multi-nerves of weakly globular $n$-fold categories.

Another crucial technique is a way to create pseudo-functors from a small category $\mathcal{C}$ to $\Cat$ out of simpler data of a functor from the objects of $\mathcal{C}$ (viewed as a discrete category) to $\Cat$. This is an instance of a more general categorical technique called 'transport of structure along an adjunction' \cite{lk} which is about $2$-dimensional monad theory, and it uses the description of pseudo-functors as pseudo-algebras for a $2$-monad.

%%%%%%%%%%%%%%%%%%%%%%%%%%%%%%%%%%%%%%%%%%%%%%%%%%%%%%%%%%%%%%%%%%%%%%%
\section{Some functors on $\pmb{\Cat}$}\label{sbs-funct-cat}

The connected component functor \index{Connected components functor}
\begin{equation*}
    q: \Cat\rw \Set
\end{equation*}
associates to a category its set of paths components. This is left adjoint to the discrete category functor
\begin{equation*}
    \di{1}:\Set \rw \Cat
\end{equation*}
associating to a set $X$ the discrete category on that set. We denote by
\begin{equation*}
    \zgu{1}:\Id\Rw \di{1}q
\end{equation*}
the unit of the adjunction $q\dashv \di{1}$.
\begin{lemma}\label{lem-q-pres-fib-pro}
    $q$ preserves fiber products over discrete objects and sends
    equivalences of categories to isomorphisms.
\end{lemma}

\begin{proof}
We claim that $q$ preserves products; that is, given categories
$\clC$ and $\clD$, there is a bijection
\begin{equation*}
    q(\clC\times \clD)=q(\clC)\times q(\clD)\;.
\end{equation*}
In fact, given $(c,d)\in q(\clC\times \clD)$ the map
$q(\clC\times\clD)\rw q(\clC)\times q(\clD)$ given by
$[(c,d)]=([c],[d])$ is well defined and is clearly surjective. On
the other hand, this map is also injective: given $[(c,d)]$ and
$[(c',d')]$ with $[c]=[c']$ and $[d]=[d']$, we have paths in $\clC$
\newcommand{\lin}{-\!\!\!-\!\!\!-}
\begin{equation*}
\xymatrix @R5pt{c \hspace{2mm} \lin \hspace{2mm}\cdots \hspace{2mm}
\lin
\hspace{2mm}c'\\
d \hspace{2mm} \lin \hspace{2mm}\cdots \hspace{2mm} \lin
\hspace{2mm}d' }
\end{equation*}
and hence a path in $\clC\times \clD$
\begin{equation*}
(c,d)\hspace{2mm}\lin\hspace{2mm}\cdots\hspace{2mm}\lin\hspace{2mm}(c',d)
\hspace{2mm}\lin\hspace{2mm}\cdots\hspace{2mm}\lin\hspace{2mm}(c',d')\;.
\end{equation*}
Thus $[(c,d)]=[(c',d')]$ and so the map is also injective, hence it
is a bijection, as claimed.

Given a diagram in $\Cat$ $\xymatrix{\clC\ar_{f}[r] & \clE & \clD
\ar^{g}[l]}$ with $\clE$ discrete, we have
\begin{equation}\label{eq-q-pres-fib-pro}
    \clC\tiund{\clE}\clD=\underset{x\in\clE}{\coprod}\clC_x\times
    \clD_x
\end{equation}
where $\clC_x,\;\clD_x$ are the full subcategories of $\clC$ and
$\clD$ with objects $c,\;d$ such that \;$f(c)=x=g(d)$. Since $q$ preserves
products and (being left adjoint) coproducts, we conclude by
\eqref{eq-q-pres-fib-pro} that
\begin{equation*}
    q(\clC\tiund{\clE}\clD)\cong q(\clC)\tiund{\clE}\,q(\clD)\;.
\end{equation*}
Finally, if $F:\clC\simeq \clD:G$ is an equivalence of categories,
$FG\,\clC\cong\clC$ and $FG\,\clD\cong \clD$ which implies that
$qF\,qG\,\clC\cong q\clC$ and $qF\,qG\,\clD\cong q\clD$, so $q\clC$
and $q\clD$ are isomorphic.
\end{proof}
The isomorphism classes of objects functor  \index{Isomorphism classes of objects functor}
\begin{equation*}
    p:\Cat\rw\Set
\end{equation*}
associates to a category the set of isomorphism classes of its objects. Note that if $\clC$ is a groupoid, $p\clC=q\clC$.
\begin{lemma}\label{lem-p-pres-fib-pro}
    $p$ preserves pullbacks over discrete objects and sends
    equivalences of categories to isomorphisms.
\end{lemma}
\begin{proof}
For a category $\clC$, let $m\clC$ be its maximal sub-groupoid. Then $p\clC=qm\clC$. Given a diagram in $\Cat$ $\xymatrix{\clC\ar_{f}[r] & \clE & \clD
\ar^{g}[l]}$ with $\clE$ discrete, we have
\begin{equation*}
    \clC\tiund{\clE}\clD=\underset{x\in\clE}{\coprod}\clC_x\times
    \clD_x\;.
\end{equation*}
Since, as easily seen, $m$ commutes with (co)products, and $m\clE=\clE$, we obtain $m(\clC\tiund{\clE}\clD)=m\clC\tiund{\clE}m\clD$; so by Lemma \ref{lem-q-pres-fib-pro},
\begin{align*}
    p(\clC\tiund{\clE}\clD)&=qm(\clC\tiund{\clE}\clD)=q(m\clC\tiund{\clE}m\clD)
    =qm\clC\tiund{q\clE}qm\clD=p\clC\tiund{\clE}p\clD\;.
\end{align*}
Finally, if $F:\clC\simeq \clD:G$ is an equivalence of categories, $FG\clC\cong \clC$ and $FG\clD\cong\clD$ which implies that $pF\,pG\,\clC\cong p\clC$ and $qF\,qG\,\clD\cong q\clD$, so $q\clC$ and $q\clD$ are isomorphic.
\end{proof}

\begin{lemma}\label{lem-iso-cla-obj-fib-pro}\

    \begin{itemize}
      \item [a)] Let $X\xrw{f} Z \xlw{g}Y$ be a diagram in $\Cat$. Then
      \begin{equation*}
        p(X\tiund{Z}Y)\subseteq pX\tiund{pZ}pY\;.
      \end{equation*}

      \item [b)] Suppose, further, that $g_0=\Id$. Then
      \begin{equation*}
        p(X\tiund{Z}Y)\cong pX\tiund{pZ}pY\;.
      \end{equation*}
    \end{itemize}
\end{lemma}
\begin{proof}\

\nid a) The map
\begin{equation*}
    j:p(X\tiund{Z}Y)\rw pX\tiund{pZ}pY
\end{equation*}
is determined by the maps
\begin{equation*}
    p(X\tiund{Z}Y)\rw pX\quad \text{and}\quad p(X\tiund{Z}Y)\rw pY
\end{equation*}
induced by the projections
\begin{equation*}
    X\tiund{Z}Y \rw X \quad \text{and}\quad X\tiund{Z}Y \rw Y\;.
\end{equation*}
Thus, for each $(a,b)\in X\tiund{Z}Y$,
\begin{equation}\label{eq1-lem-iso-cla}
    j\,p(a,b)=(p(a),p(b))\;.
\end{equation}
Let $(a,b), (a',b')\in X\tiund{Z}Y$ be such that $p(a,b)=p(a',b')$. It follows by \eqref{eq1-lem-iso-cla} that $p(a)=p(a')$ and $p(b)=p(b')$. Thus there are isomorphisms $\za:a\cong a'$ in $X$ and $\zb: b\cong b'$ in $Y$ and in $Z$ we have
\begin{equation*}
    f\za=g\zb: fa=gb\cong fa'=gb'\;.
\end{equation*}
Thus $(\za,\zb):(a,b)\cong (a',b')$ is an isomorphism in $X\tiund{Z} Y$ and so $p(a,b) = p(a',b')$. This shows that $j$ is injective, proving a).
\mk

\nid b) By a) the map $j$ is injective. We now show that, if $g_0=\Id$, then $j$ is also surjective. Let $(px,py)\in pX\tiund{pZ} pY$ be such that $g_0=\Id$. Then
\begin{equation*}
    f p(x)= p f(x)=gp(y)=pg(y)=p(y)\;.
\end{equation*}
Also, $f(x)=g(f(x))$, so that $(x,f(x))\in X\tiund{Z}Y$. It follows that
\begin{equation*}
    j(x,f(x))=(p(x),pf(x))=(p(x),p(y))
\end{equation*}
so that $j$ is surjective. Hence $j$ is a bijection.

\end{proof}

\begin{definition}\label{def-isof}\index{Isofibration}
  Recall that the functor $F:X\rw Y$ is an isofibration if for each $x\in X$ and isomorphism $\alpha:Fx\cong y$ in $Y$, there is an isomorphism $\beta:x\cong z$ in $X$ with $F\beta=\alpha$.
 \end{definition}

\begin{lemma}\label{lem-gen-const-1}
Let
\begin{equation*}
    \xymatrix{
    A \ar_{g'}[r] \ar_{f'}[d] & B \ar^{f}[d]\\
    C \ar_{g}[r] & D
    }
\end{equation*}
be a pullback in $\Cat$ with $f$ an isofibration. Then
\begin{itemize}
  \item [a)]\begin{equation*}
    \xymatrix{
    pA \ar[r] \ar[d] & pB \ar^{}[d]\\
    pC \ar_{}[r] & pD
    }
\end{equation*}
is a pullback in $Set$.
  \item [b)] Suppose that $D$ is a groupoid. Then

  \begin{equation*}
    \xymatrix{
    qA \ar[r] \ar[d] & qB \ar^{}[d]\\
    qC \ar_{}[r] & qD
    }
\end{equation*}
is a pullback in $Set$.
\end{itemize}

\end{lemma}

\begin{proof}
a) Since $f$ is an isofibration, by \cite{JoyStr1993} $A$ is equivalent to the pseudo-pullback $\ds A \simeq C \tms{ps}{D} B$. The functor $p$ sends pseudo-pullbacks to pullbacks. In fact, suppose we are given a commuting diagram in $\Set$
\begin{equation*}
    \xymatrix{
    X \ar@/_/[ddr]_r \ar@/^/[drr]^s \\
    & p(\ds{C\tms{ps}{D} B}) \ar[d] \ar[r]  & pB \ar^{p(f)}[d] \\
    & pC \ar_{p(g)}[r] & pD }
\end{equation*}
$p(g)r =p(f)s$. If we choose maps $b:dpB\rw B$ and $c:dpC\rw C$ (so that $p(b)=id$ and $p(c)=id$) we have
\begin{equation*}
  p(f)p(b)p(d(s))=p(f)s=p(g)r=p(g)p(c)p(d(r)).
\end{equation*}

It follows that, for each $x\in X$
\begin{equation*}
  (fbd(s))(x)\cong (gcd(r))(x)
\end{equation*}
 Therefore, there is $\ds v:X \rw C \tms{ps}{D} B$ such that $g'v=bd(s)$ and $f' v=cd(r)$. Hence
\begin{equation*}
    p(g')p(v)=p(b)p(d(s))=s,\qquad p(f')p(v)=p(c)p(d(r))=r\;.
\end{equation*}
This shows that
\begin{equation*}
    p(A) \cong p(\ds{C\tms{ps}{D} B}) \cong pC\tiund{pD}pB\;.
\end{equation*}
b) Since $f$ is an isofibration, by \cite{JoyStr1993} $A$ is equivalent to the pseudo-pullback $\ds A \simeq C \tms{ps}{D} B$. The functor $q$ sends pseudo-pullbacks to pullbacks. In fact, suppose we are given a commuting diagram in $\Set$
\begin{equation*}
    \xymatrix{
    X \ar@/_/[ddr]_r \ar@/^/[drr]^s \\
    & q(\ds{C\tms{ps}{D} B}) \ar[d] \ar[r]  & qB \ar^{q(f)}[d] \\
    & qC \ar_{q(g)}[r] & qD }
\end{equation*}
$q(g)r =q(f)s$. If we choose maps $b:dqB\rw B$ and $c:dqC\rw C$ (so that $q(b)=id$ and $q(c)=id$) we have
\begin{equation*}
  q(f)q(b)q(d(s))=q(f)s=q(g)r=q(g)q(c)q(d(r)).
\end{equation*}

It follows that, for each $x\in X$, $(fbd(s))(x)$ and $(gcd(r))(x)$ are in the same connected component. Since $D$ is a groupoid, this means that there is an isomorphism
\begin{equation*}
  (fbd(s))(x)\cong (gcd(r))(x)
\end{equation*}
 Therefore, there is $\ds v:X \rw C \tms{ps}{D} B$ such that $g'v=bd(s)$ and $f' v=cd(r)$. Hence
\begin{equation*}
    q(g')q(v)=q(b)q(d(s))=s,\qquad q(f')q(v)=q(c)q(d(r))=r\;.
\end{equation*}
This shows that
\begin{equation*}
    q(A)\cong q(\ds{C\tms{ps}{D} B})\cong qC\tiund{qD}qB\;.
\end{equation*}
\end{proof}
\begin{lemma}\label{lem-gen-const-2}
    Let
    \begin{equation*}
        \xymatrix{
         A \ar^{s}[r] \ar_{r}[d] & B \ar^{f}[d]\\
    C \ar_{g}[r] & D
        }
    \end{equation*}
    be a pullback in $\Cat$, and suppose that $f$ is fully faithful. Then so is $r$.
\end{lemma}
\begin{proof}
For all $x,y\in A_0$ there is a pullback in $\Cat$
\begin{equation*}
    \xymatrix{
    A(x,y) \ar[rr] \ar[d] && B(sx,sy) \ar[d]\\
    C(rx.ry) \ar[rr] && D(grx,gry)=D(fsx,fsy)
    }
\end{equation*}
Since the right vertical map is an isomorphism (as $f$ is fully faithful), so is the left vertical map, showing that $r$ is fully faithful.
\end{proof}

%%%%
\section{Pseudo-functors and their strictification}\label{sbs-pseudo-functors}
\index{Pseudo-functor} \index{Strictification functor}

We recall the notion of pseudo-functor and the classical theory of strictification of pseudo-functors, see \cite{Borc},\cite{PW}, \cite{Lack}.

\subsection{The notion of pseudo-functor}
\begin{definition}[\cite{Borc}]\label{def-pseudo-fun}
A pseudo-functor $F:\clA \rw \clB$ between 2-categories $\clA, \clB$ consists of the following data:
\begin{itemize}
  \item [(1)] For every $A\in \clA$, an object $FA\in\clB$.\bk

  \item [(2)] For every pair of objects $A,B \in\clA$, a functor
  \begin{equation*}
    F_{AB}:\clA(A,B) \rw \clB(FA,FB)\;.
  \end{equation*}

  \item [(3)] For every triple of objects $A, B, C\in\clA$, a natural isomorphism $\zg_{ABC}$:
  \begin{equation*}
  \xymatrix@C=80pt@R=40pt{
    \Sc{\clA(A,B)\times\clA(B,C)} \ar^{\Scc{C_{ABC}}}[r] \ar_{\Scc{F_{AB}\times \Scc{F_{BC}}}}[d]  &
    \Sc{\clA(A,C)} \ar^{\Scc{F_{AC}}}[d]\\
    %%%%%%%%%%%%%%%%%%%%%%%%%%%%%%%%%%%%%%%%%%%%%%%%%%%%%%%%
    \Sc{\clB(FA,FB)\times\clB(FB,FC)} \ar_{\Scc{C_{FA,FB,FC}}}[r] \ar@{}[ur]^(.35){}="a"^(.65){}="b" \ar@{=>}^{\zg_{\Scc{ABC}}}"a";"b"  &
    \Sc{\clB(FA,FC)}
  }
  \end{equation*}

  \item [(4)] For every object $A\in\clA$, a natural isomorphism $\zd_A$:
  \begin{equation*}
  \xymatrix@C=50pt@R=40pt{
    \mathbbm{1} \ar^(0.4){u_A}[r] \ar@{=}[d] & \clA(A,A) \ar^{F_{AA}}[d]\\
    \mathbbm{1} \ar_(0.4){u_{FA}}[r]\ar@{}[ur]^(.35){}="a"^(.65){}="b" \ar@{=>}^{\zd_{A}}"a";"b" & \clB(FA,FA)
  }
  \end{equation*}
  \end{itemize}

  such that the following coherence axioms are satisfied:

  \begin{itemize}
    \item [(i)] Composition axiom: for every triple of arrows
    \begin{equation*}
      A\xrw{f} B \xrw{g} C \xrw{h} D
    \end{equation*}
    in $\clA$, the following equality between 2-cells holds
  \begin{equation*}
  \xymatrix@C=80pt@R=40pt{
    Fh \cirsm Fg\cirsm Ff \ar^{i_{Fh}*\zg_{f,g}}[r] \ar_{\zg_{g,h}*i_{Ff}}[d]  &
    Fh\cirsm F(g\cirsm f) \ar^{\zg_{g{\cirsm} f, h}}[d]\\
    %%%%%%%%%%%%%%%%%%%%%%%%%%%%%%%%%%%%%%%%%%%%%%%%%%%%%%%%
    F(h\cirsm g)\cirsm Ff \ar_{\zg_{f,h{\cirsm} g}}[r] &
    F(h\cirsm g\cirsm f)
  }
  \end{equation*}

    \item [(ii)] Unit axiom: for every arrow $f:A\rw B$ in $\clA$, the following equalities between 2-cells hold
    \begin{equation*}
    \xymatrix@C=40pt@R=40pt{
    Ff\cirsm 1_{FA} \ar^{i_{Ff}*\zd_A}[r] \ar_{i_{Ff}}[d] & Ff\cirsm F 1_A \ar^{\zg_{1_A,f}}[d] \\
    Ff \ar_{i_{Ff}}[r] & F(f\cirsm 1_A)
    }
    \qquad
    \xymatrix@C=40pt@R=40pt{
    1_{FB}\cirsm Ff \ar^{\zd_A*i_{Ff}}[r] \ar_{i_{Ff}}[d] & F 1_B\cirsm Ff \ar^{\zg_{f, 1_B}}[d] \\
    Ff \ar_{}[r] & F(1_B\cirsm f)
    }
    \end{equation*}
\end{itemize}
\end{definition}

There is also a notion of pseudo-natural transformation between pseudo-functors (see \cite[\S 7.5]{Borc} for more details)  and thus a category $Ps[\clA,\clB]$ of pseudo-functors and pseudo-natural transformations.

\subsection{Strictification of pseudo-functors}
The functor 2-category $\funcat{n}{\Cat}$ is 2-monadic over $[ob(\Dnop),\Cat]$ where $ob(\Dnop)$ is the set of objects of $\Dnop$. Let
\begin{equation*}
    U:\funcat{n}{\Cat}\rw [ob(\Dnop),\Cat]
\end{equation*}
be the forgetful functor $(UX)_{\uk}=X_{\uk}$. Its left adjoint $F$ is given on objects by
\begin{equation*}
    (FH)_{\uk}=\underset{\ur\in ob(\Dnop)}{\coprod}\Dnop(\ur,\uk)\times H_{\ur}
\end{equation*}
for $H\in [ob(\Dmenop),\Cat]$, $\uk\in ob(\Dmenop)$. If $T$ is the monad corresponding to the adjunction $F\dashv U$, then
\begin{equation*}
    (TH)_{\uk}=\underset{\ur\in ob(\Dnop)}{\coprod}\Dnop(\ur,\uk)\times H_{\ur}
\end{equation*}
A pseudo $T$-algebra is given by $H\in [ob(\Dnop),\Cat]$,
\begin{equation*}
    h_{n}: \underset{\ur\in ob(\Dnop)}{\coprod}\Dnop(\ur,\uk)\times H_{\ur} \rw H_{\uk}
\end{equation*}
and additional data, as described in \cite{PW}. This amounts precisely to functors from $\Dnop$ to $\Cat$ and the 2-category $\sf{Ps\mi T\mi alg}$ of pseudo $T$-algebras corresponds to the 2-category $\Ps\funcat{n}{\Cat}$ of pseudo-functors, pseudo-natural transformations and modifications.

The strictification result proved in \cite{PW} yields that every pseudo-functor from $\Dnop$ to $\Cat$ is equivalent, in $\Ps\funcat{n}{\Cat}$, to a 2-functor.

Given a pseudo $T$-algebra as above, \cite{PW} consider the factorization  of $h:TH\rw H$ as
\begin{equation*}
    TH\xrw{v}L\xrw{g}H
\end{equation*}
with $v_{\uk}$ bijective on objects and $g_{\uk}$ fully faithful, for each $\uk\in\Dnop$. It is shown in \cite{PW} that it is possible to give a strict $T$-algebra structure $TL\rw L$ such that $(g,Tg)$ is an equivalence of pseudo $T$-algebras. It is immediate to see that, for each $\uk\in\Dnop$, $g_{\uk}$ is an equivalence of categories.

Further, it is shown in \cite{Lack} that $\St:\psc{n}{\Cat}\rw\funcat{n}{\Cat}$ as described above is left adjoint to the inclusion
\begin{equation*}
  J:\funcat{n}{\Cat}\rw\psc{n}{\Cat}
\end{equation*}
 and that the components of the units are equivalences in $\psc{n}{\Cat}$.
\section{Transport of structure}\label{transport-structure}
\index{Transport of structure}

We now recall a general categorical technique, known as transport of structure along an adjunction, with one of its applications. This will be used crucially in the proof of Theorem \ref{the-XXXX}.
\begin{theorem}\rm{\cite[Theorem 6.1]{lk}}\em\label{s2.the1}
    Given an equivalence $\;\eta,\;\zve : f \dashv f^* : A\rw B$ in the complete and locally small
    2-category $\clA$, and an algebra $(A,a)$ for the monad $T=(T,i,m)$ on $\clA$, the
    equivalence enriches to an equivalence
\begin{equation*}
  \eta,\zve:(f,\ovll{f})\vdash (f^*,\ovll{f^*}):(A,a)\rw(B,b,\hat{b},\ovl{b})
\end{equation*}
in $\PsTalg$, where $\hat{b}=\eta$, $\;\ovl{b}=f^* a\cdot T\zve \cdot
Ta\cdot T^2 f$, $\;\ovll{f}=\zve^{-1} a\cdot Tf$, $\;\ovll{f^*}=f^* a\cdot
T\zve$.
\end{theorem}
Let $\eta',\zve':f'\vdash f'^{*}:A'\rw B'$ be another equivalence in $\clA$ and
let $(B',b',\hat{b'},\ovl{b'})$ be the corresponding pseudo-$T$-algebra as in
Theorem \ref{s2.the1}. Suppose $g:(A,a)\rw(A',a')$ is a morphism in $\clA$ and
$\gamma$ is an invertible 2-cell in $\clA$
\begin{equation*}
  \xy
    0;/r.10pc/:
    (-20,20)*+{B}="1";
    (-20,-20)*+{B'}="2";
    (20,20)*+{A}="3";
    (20,-20)*+{A'}="4";
    {\ar_{f^*}"3";"1"};
    {\ar_{h}"1";"2"};
    {\ar^{f'^*}"4";"2"};
    {\ar^{g}"3";"4"};
    {\ar@{=>}^{\gamma}(0,3);(0,-3)};
\endxy
\end{equation*}
Let $\ovl{\gamma}$ be the invertible 2-cell given by the following pasting:
\begin{equation*}
    \xy
    0;/r.15pc/:
    (-40,40)*+{TB}="1";
    (40,40)*+{TB'}="2";
    (-40,-40)*+{B}="3";
    (40,-40)*+{B'}="4";
    (-20,20)*+{TA}="5";
    (20,20)*+{TA'}="6";
    (-20,-20)*+{A}="7";
    (20,-20)*+{A'}="8";
    {\ar^{Th}"1";"2"};
    {\ar_{b}"1";"3"};
    {\ar^{b'}"2";"4"};
    {\ar_{h}"3";"4"};
    {\ar_{Tg}"5";"6"};
    {\ar^{}"5";"7"};
    {\ar_{}"6";"8"};
    {\ar^{g}"7";"8"};
    {\ar_{Tf^*}"5";"1"};
    {\ar^{f^*}"7";"3"};
    {\ar^{Tf'^*}"6";"2"};
    {\ar^{f'^*}"8";"4"};
    {\ar@{=>}^{(T\gamma)^{-1}}(0,33);(0,27)};
    {\ar@{=>}^{\gamma}(0,-27);(0,-33)};
    {\ar@{=>}^{\ovll{f'^*}}(30,3);(30,-3)};
    {\ar@{=>}^{\ovll{f^*}}(-30,3);(-30,-3)};
\endxy
\end{equation*}
Then it is not difficult to show that
$(h,\ovl{\gamma}):(B,b,\hat{b},\ovl{b})\rw(B',b',\hat{b'},\ovl{b'})$ is a
pseudo-$T$-algebra morphism.

The following fact is essentially known and, as sketched in the proof below, it is an instance of Theorem \ref{s2.the1}
\begin{lemma}{\rm{\cite{PP}}}\label{lem-PP}
     Let $\clC$ be a small 2-category, $F,F':\clC\rw\Cat$ be 2-functors, $\alpha:F\rw F'$
    a 2-natural transformation. Suppose that, for all objects $C$ of $\clC$, the
    following conditions hold:
\begin{itemize}
  \item [i)] $G(C),\;G'(C)$ are objects of $\Cat$ and there are adjoint equivalences of
  categories $\mu_C\vdash\eta_C$, $\mu'_C\vdash\eta'_C$,
\begin{equation*}
  \mu_C:G(C)\;\rightleftarrows\;F(C):\eta_C\qquad\qquad
  \mu'_C:G'(C)\;\rightleftarrows\;F'(C):\eta'_C,
\end{equation*}
  \item [ii)] there are functors $\beta_C:G(C)\rw G'(C),$
  \item [iii)] there is an invertible 2-cell
\begin{equation*}
  \gamma_C:\beta_C\,\eta_C\Rightarrow\eta'_C\,\alpha_C.
\end{equation*}
\end{itemize}
Then
\begin{itemize}
  \item [a)] There exists a pseudo-functor $G:\clC\rw\Cat$ given on objects by $G(C)$,
  and pseudo-natural transformations $\eta:F\rw G$, $\mu:G\rw F$ with
  $\eta(C)=\eta_C$, $\mu(C)=\mu_C$; these are part of an adjoint equivalence
  $\mu\vdash\eta$ in the 2-category $\Ps[\clC,\Cat]$.
  \item [b)] There is a pseudo-natural transformation $\beta:G\rw G'$ with
  $\beta(C)=\beta_C$ and an invertible 2-cell in $\Ps[\clC,\Cat]$,
  $\gamma:\beta\eta\Rightarrow\eta\alpha$ with $\gamma(C)=\gamma_C$.
\end{itemize}
\end{lemma}
\begin{proof}
Recall \cite{PW} that the functor 2-category $[\clC,\Cat]$ is 2-monadic over
$[ob(\clC),\Cat]$, where $ob(\clC)$ is the set of objects in $\clC$. Let
\begin{equation*}
  \clU:[\clC,\Cat]\rw[ob(\clC),\Cat]
\end{equation*}
be the forgetful functor. Let $T$ be the 2-monad; then the pseudo-$T$-algebras are precisely the pseudo-functors from
$\clC$ to $\Cat$.

Then the adjoint equivalences $\mu_C\vdash\eta_C$ amount precisely to an
adjoint equivalence in $[ob(\clC),\Cat]$, $\;\mu_0\vdash\eta_0$,
$\;\mu_0:G_0\;\;\rightleftarrows\;\;\clU F:\eta_0$ where $\;G_0(C)=G(C)$ for
all $C\in ob(\clC)$. This equivalence enriches to an
adjoint equivalence $\mu\vdash\eta$ in $\Ps[\clC,\Cat]$
\begin{equation*}
  \mu:G\;\rightleftarrows\; F:\eta
\end{equation*}
between $F$ and a pseudo-functor $G$; it is $\clU G=G_0$, $\;\clU\eta=\eta_0$,
$\;\clU\mu=\mu_0$; hence on objects $G$ is given by $G(C)$, and
$\eta(C)=\clU\eta(C)=\eta_C$, $\;\mu(C)=\clU\mu(C)=\mu_C$.

Let $\nu_C:\Id_{G(C)}\Rw\eta_C\mu_C$ and $\zve_C:\mu_C\eta_C\Rw\Id_{F(C)}$ be
the unit and counit of the adjunction $\mu_C\vdash\eta_C$. Given a morphism $f:C\rw D$ in $\clC$, it is
\begin{equation*}
  G(f)=\eta_D F(f)\mu_C
\end{equation*}
and we have natural isomorphisms:
\begin{align*}
   & \eta_f : G(f)\eta_C=\eta_D F(f)\mu_C\eta_C\overset{\eta_D F(f)\zve_C}{=\!=\!=\!=\!\Rw} \eta_D F(f)\\
   & \mu_f : F(f)\mu_C\overset{\nu_{F(f)}\mu_C}{=\!=\!=\!\Rw}\mu_D\eta_D F(f)\mu_C=\mu_D
   G(f).
\end{align*}
Also, the natural isomorphism
\begin{equation*}
  \beta_f: G'(f)\beta_C\Rw\beta_D G(f)
\end{equation*}
is the result of the following pasting
\begin{equation*}
    \xy
    0;/r.13pc/:
    (-40,40)*+{G(C)}="1";
    (40,40)*+{G'(C)}="2";
    (-40,-40)*+{G(D)}="3";
    (40,-40)*+{G'(D)}="4";
    (-20,20)*+{F(C)}="5";
    (20,20)*+{F'(C)}="6";
    (-20,-20)*+{F(D)}="7";
    (20,-20)*+{F'(D)}="8";
    {\ar^{\beta_C}"1";"2"};
    {\ar_{G(f)}"1";"3"};
    {\ar^{G'(f)}"2";"4"};
    {\ar_{\beta_D}"3";"4"};
    {\ar^{\alpha_C}"5";"6"};
    {\ar^{F(f)}"5";"7"};
    {\ar_{F'(f)}"6";"8"};
    {\ar_{\alpha'_D}"7";"8"};
    {\ar^{}"5";"1"};
    {\ar^{}"7";"3"};
    {\ar^{}"6";"2"};
    {\ar^{}"8";"4"};
    {\ar@{=>}^{\gamma_C}(0,33);(0,27)};
    {\ar@{=>}^{\gamma_D^{-1}}(0,-27);(0,-33)};
    {\ar@{=>}^{\eta_f'}(30,3);(30,-3)};
    {\ar@{=>}^{\eta_f}(-30,3);(-30,-3)};
\endxy
\end{equation*}
\end{proof}
%%

%%%%%%%%%%%%%%%%%%%%%%%%%%%%%%%%%%%%%%%%%%%%%%%%%%%%%%%%%%%%%%%%%%%%%%%%%

\part{Weakly globular n-fold categories and Segalic pseudo-functors}\label{part-3}

\newpage

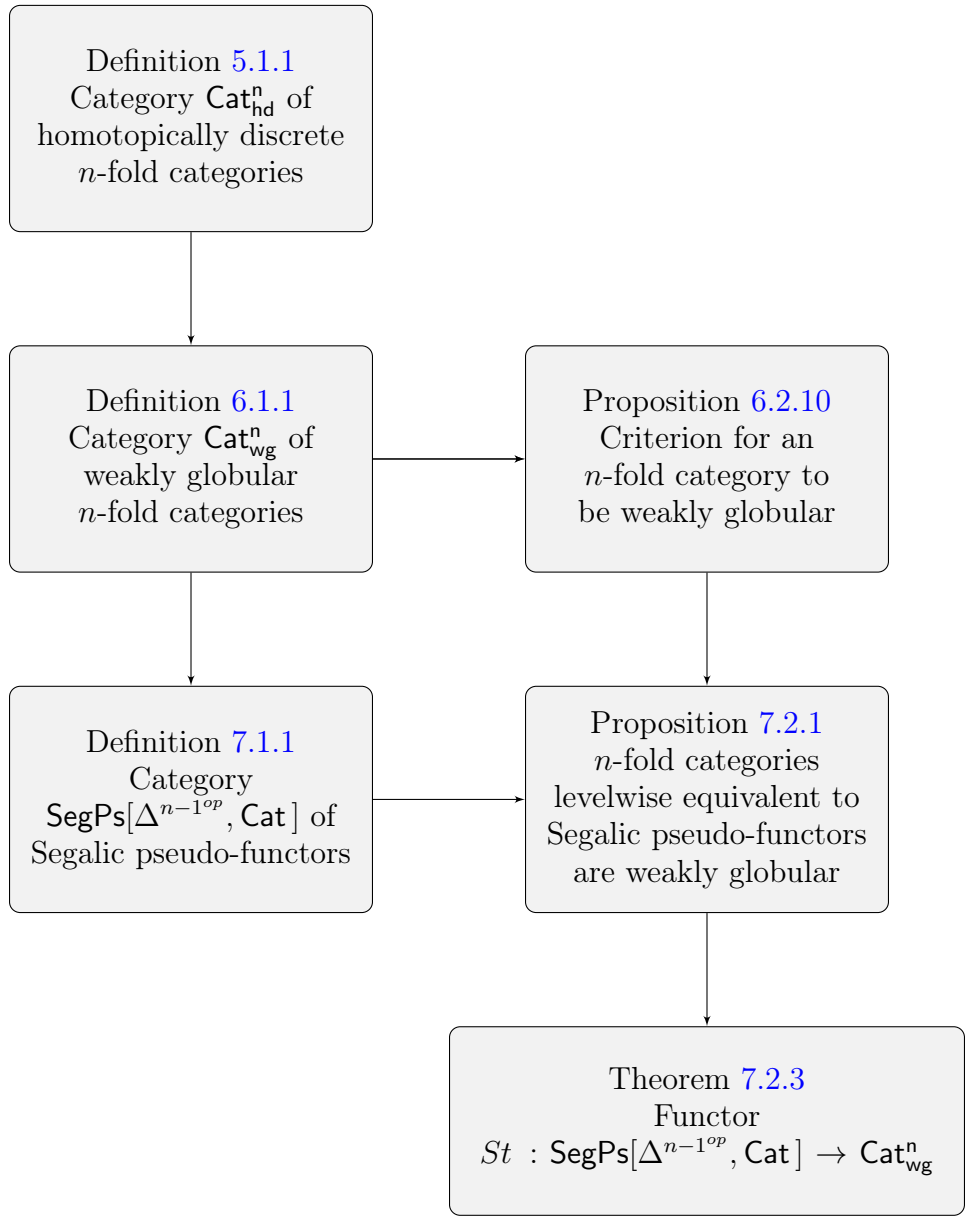
\begin{figure}[ht]
  \centering
\vspace{10mm}
\begin{tikzpicture}[node distance = 15mm and 20mm]% First vertical distance and second horizontal from border
\node [block] (def3_1_1) {Definition \ref{def-hom-dis-ncat}\\Category $\cathd{n}$ of homotopically discrete $n$-fold categories};
\node [block, below= of def3_1_1] (def4_1_1) {Definition \ref{def-n-equiv}\\Category $\catwg{n}$ of weakly globular $n$-fold categories};
\node [block, right= of def4_1_1] (pro4_2_10) {Proposition \ref{pro-crit-ncat-be-wg}\\Criterion for an\\ $n$-fold category to be weakly globular};
\node [block, below= of def4_1_1] (def5_1_1) {Definition \ref{def-seg-ps-fun}\\Category $\segpsc{n-1}{\Cat}$ of Segalic pseudo-functors};
\node [block, right= of def5_1_1] (pro5_2_1) {Proposition \ref{pro-transf-wg-struc}\\ $n$-fold categories levelwise equivalent to Segalic pseudo-functors are weakly globular};
\node [block1, below= of pro5_2_1] (the5_2_3) {Theorem \ref{the-strict-funct}\\Functor\\ $\St:\segpsc{n-1}{\Cat}\rw\catwg{n}$};

% Draw edges
\path [line] (def3_1_1) -- (def4_1_1);
\path [line] (def4_1_1) -- (pro4_2_10);
\path [line] (def4_1_1) -- (def5_1_1);
\path [line] (def4_1_1) -- (pro4_2_10);
\path [line] (pro4_2_10) -- (pro5_2_1);
\path [line] (def5_1_1) -- (pro5_2_1);
\path [line] (pro5_2_1) -- (the5_2_3);
\end{tikzpicture}
  \caption{Weakly globular $n$-fold categories and Segalic pseudo-functors.}
  \label{FigIntro-4}
\end{figure}

In Part \ref{part-3} we introduce the main new structure of this work, the category $\catwg{n}$ of weakly globular \nfol categories, and study its relation to the special class of pseudo-functors
 \begin{equation*}
   \segpsc{n-1}{\Cat} \subset \psc{n-1}{\Cat}
 \end{equation*}
 called Segalic pseudo-functors. The main result of this part is Theorem \ref{the-strict-funct}, establishing that the classical strictification of pseudo-functors, when restricted to Segalic pseudo-functors, yields weakly globular \nfol categories; that is, we have a functor
 \begin{equation*}
        \St: \segpsc{n-1}{\Cat}\rw \catwg{n}
    \end{equation*}
 A summary of the main steps in this part is given in Figure \ref{FigIntro-4}.

 In Chapter \ref{chap3} we introduce the category $\cathd{n}$ of homotopically discrete \nfol categories and study its properties. The idea of the category $\cathd{n}$ is introduced in Section \ref{subs-idea-cathd} before the formal definition. This category is needed for the precise formulation of the weak globularity condition in the definitions of the categories $\catwg{n}$ of weakly globular \nfol categories and $\tawg{n}$ of weakly globular Tamsamani $n$-categories. We also introduce $n$-equivalences of homotopically discrete \nfol categories, and show in Lemma \ref{lem-neq-hom-disc} that they are detected by isomorphisms of their discretizations.

 In Chapter \ref{chap4} we introduce the category $\catwg{n}$ of weakly globular \nfol categories and study its properties. The idea of this category is introduced in Section \ref{subs-idea-catwg} before the formal definition. The main result of this chapter is Proposition \ref{pro-crit-ncat-be-wg}. Part a) of this Proposition establishes that the nerve functor in direction $2$ on the category $\catwg{n}$ is levelwise a weakly globular $(n-1)$-fold category: this will be used at several occasions in the rest of the work, for instance in Corollary \ref{pro-post-wg-ncat} in showing the existence of the functor
  \begin{equation*}
     \qn:\catwg{n}\rw \catwg{n-1}
  \end{equation*}
  which plays an important role in the construction of the rigidification functor.

  Part b) of Proposition \ref{pro-crit-ncat-be-wg} gives a sufficient criterion for a \nfol category to be weakly globular: this requires  certain sub-structures to be homotopically discrete as well as the fact that applying levelwise in the $n$-th direction the functor isomorphism classes of objects $p$ yields a weakly globular $(n-1)$-fold  category.

   In Chapter \ref{chap5} we introduce the category $\segpsc{n}{\Cat}$ of Segalic pseudo-functors and study its main properties. The idea of this category is introduced in Section \ref{subs-idea-segps} before the formal definition. In Proposition \ref{pro-transf-wg-struc} we establish that an \nfol category whose multinerve is levelwise equivalent to a Segalic pseudo-functor is weakly globular: the proof of this result is based on the criterion of Proposition \ref{pro-crit-ncat-be-wg} b) together with Corollary \ref{cor-crit-nequiv-rel} which gives a sufficient condition for a weakly globular \nfol category to be homotopically discrete.

   In Lemma \ref{lem-strict-mon-wg-seg-ps-fun} we study the properties of the monad corresponding to Segalic pseudo-functors. Proposition \ref{pro-transf-wg-struc} and Lemma \ref{lem-strict-mon-wg-seg-ps-fun} lead to the main result Theorem \ref{the-strict-funct} on the strictification of Segalic pseudo-functors.\bk

% Figure

\clearpage

%%%%%%%%%%%%%%%%%%%%%%%%%%%%%%%%%%%%%%%%%%%%%%%%%%%%%%%%%%%%%%%%%%%%%%%%%
\chapter{Homotopically discrete $\pmb{n}$-fold categories}\label{chap3}
In this chapter we introduce the category $\cathd{n}$ of homotopically discrete \nfol categories, which will be needed in Chapter \ref{chap4} and \ref{chap6} to define the categories of weakly globular \nfol categories and of weakly globular Tamsamani $n$-categories.

Recall (see Definition \ref{def-discr-nfold}) that a discrete $n$-fold category is one in which the multinerve is a constant functor: thus it amounts to a set and identity structure maps in all the simplicial directions. Homotopically discrete $n$-fold categories are 'discrete up to homotopy' in a specified way. An object $X\in \cathd{n}$ comes equipped with a discretization map
 \begin{equation*}
 \zg\lo{n}:X\rw X^d
 \end{equation*}
  where $X^d$ is a discrete \nfol category and $\zg\lo{n}$ is a suitably defined higher categorical equivalence.

   In the case $n=1$, an object of $\cathd{}$ is a groupoid equivalent to a discrete category, that is an equivalence relation, or equivalently a groupoid with no non-trivial loops.

As outlined in Section \ref{sbs-idea-of}, homotopically discrete $k$-fold categories (for $1 \leq k \leq n-1$) are sub-structures of our two Segal-type models $\catwg{n}$ and $\tawg{n}$. The sets underlying the discretizations of the homotopically discrete sub-structures in objects of $\catwg{n}$ and $\tawg{n}$ play the role of sets of higher cells. In the category $\catwg{n}$ the weakness of the structure is encoded in these homotopically discrete objects: indeed an object of $\catwg{n}$ in which all these homotopically discrete sub-structures are actually discrete is a strict $n$-category.

Homotopically discrete \nfol categories are higher groupoidal structures to which we can associate a classifying space which is a $0$-type (see Proposition \ref{pro-sign-cathd}), that is a topological space whose homotopy groups are zero in dimension greater than 0.

We present two equivalent descriptions of homotopically discrete $n$-fold categories: one is a multi-simplicial description built inductively on dimension (see Definition \ref{def-hom-dis-ncat}). As explained in the next chapter (see Remark \ref{rem-n-equiv}) this has the advantage of making it evident that $\cathd{n}$ is a subcategory of $\catwg{n}$. The second description is more conceptual and uses an iteration of the notion of internal equivalence relation. We establish in Theorem \ref{the-hom-disc-neq-rel} that these two descriptions are equivalent.

This chapter is organized as follows. In Section \ref{sec-hom-dis-ncat} we give an inductive definition of the category $\cathd{n}$ and of $n$-equivalences, and we define the discretization map $\zg\lo{n}:X \rw X^d$. In Section \ref{sbs-proper-nfol} we establish the main properties of the category $\cathd{n}$. In Lemma \ref{lem-pos-grou-hom-disc} we show that $\cathd{n}$ can be viewed as a diagram of equivalence relations. In Lemma \ref{lem-neq-hom-disc} we show that $n$-equivalences in $\cathd{n}$ are detected by isomorphisms of their discretizations, and we deduce in Proposition \ref{pro-ind-map-hom-disc} that the induced Segal maps of objects of $\cathd{n}$ are $(n-1)$-equivalences. In Section \ref{highrel} we give a different description of weakly globular \nfol categories via a notion of iterated internal equivalence relations.

%%%%%%%%%%%%%%%%%%%%%%%%%%%%%%%%%%%%%%%%%%%%%%%%%%%%%%%%%%%%%%%%%%%%%%%%%%%%%%%%%%%%%%%%%%%%%

\section{The definition of homotopically discrete $\pmb{n}$-fold categories}\label{sec-hom-dis-ncat}
In this section we give an inductive definition of the category $\cathd{n}$ of homotopically discrete \nfol categories and of $n$-equivalences between them. The latter are a higher dimensional generalization of a functor which is fully faithful and essentially surjective on objects.

The definition of $\cathd{n}$ is in the spirit of the Segal-type models outlined in Section \ref{three-seg}, and we will see in Chapter \ref{chap4} that the category of homotopically discrete \nfol categories is a full subcategory of the category of weakly globular \nfol categories.

\subsection{The idea of homotopically discrete $\pmb{n}$-fold category}\label{subs-idea-cathd}
When $n=1$, a homotopically discrete category is simply a groupoid with no non-trivial loops, in other words an equivalence relation. The idea of a homotopically discrete \nfol category when $n>1$ is that it is an \nfol category $X$ which is suitably equivalent to a discrete one both 'globally' and in each simplicial dimension along $\Nu{1}X \in \funcat{}{\Cat^{n-1}}$ so that in fact
 \begin{equation*}
   \Nu{1}X\in \funcat{}{\cathd{n-1}}.
 \end{equation*}
 This also implies that $J_n X \in \funcat{n-1}{\Cat}$ is levelwise an equivalence relation, that is
 \begin{equation*}
   J_n X \in \funcat{n-1}{\cathd{}}.
 \end{equation*}
 Further, we impose in the definition the condition that taking isomorphisms classes of objects in each dimension in $J_n X$ gives (the multinerve of) a homotopically discrete $(n-1)$-fold category. That is, we have truncation functors
 \begin{equation*}
        \cathd{n}\xrw{\p{n}}\cathd{n-1} \xrw{\p{n-1}}  \cdots \cathd{} \xrw{p} \Set
    \end{equation*}
and corresponding maps for each $X\in\cathd{n}$
 \begin{equation}\label{maphd}
        X\xrw{\zgu{n}}\di{n}\p{n}X \xrw{\di{n}\zgu{n-1}} \di{n}\di{n-1}\p{n-1}\p{n}X \rw \cdots \rw X^d\;.
    \end{equation}
where
\begin{equation*}
        X^d =\di{n}\di{n-1}...\di{1}\p{1}\p{2}...\p{n}X
    \end{equation*}
is a discrete \nfol category.

The notion of $n$-equivalence in $\cathd{n}$ is given according to the general pattern of the Segal-type models, see Section \ref{sec-comm-fea}.
We show in Corollary \ref{gamma} that the maps \eqref{maphd} and their composite (called discretization map) are $n$-equivalences so $X$ is $n$-equivalent to a discrete category.

 \subsection{The formal definition of $\pmb{\cathd{n}}$}

\index{Homotopically discrete $n$-fold category}
\begin{definition}\label{def-hom-dis-ncat}
    Define inductively the full subcategory $\cathd{n}\subset\cat{n}$ of homotopically discrete \nfol categories.

    For $n=1$, $\cathd{1}=\cathd{}$ is the category of  equivalence relations that is, groupoids equivalent to discrete ones. Denote by $\p{1}=p:\Cat\rw\Set$ the isomorphism classes of object functor.

    Suppose, inductively, that for each $1\leq k\leq n-1$ we defined $\cathd{k}\subset\cat{k}$  such that the following holds:

    \bigskip
    \begin{itemize}
      \item [a)] \emph{Groupoidal directions}

      \nid The $k^{th}$ direction in $\cathd{k}$ is groupoidal; that is, if $X\in\cathd{k}$, $\tilde{\xi_{k}}X\in\Gpd(\cat{k-1})$ (where $\tilde{\xi_{k}}X$ is as in Proposition \ref{pro-mult-ner}).

\bigskip
      \item [b)] \emph{Truncation functor} \index{Truncation functor}

      \nid There is a functor $\p{k}:\cathd{k}\rw\cathd{k-1}$ making the following diagram commute:\index{Functor!- $\p{n}$}

          \begin{equation}\label{eq-p-fun-def}%
            \xymatrix{
            \cathd{k} \ar^{J_{k-1}}[rrr] \ar_{p^{(k)}}[d] &&& \funcat{k-1}{\Cat} \ar^{\bar p}[d]\\
            \cathd{k-1} \ar_{\N{k-1}}[rrr] &&& \funcat{k-1}{\Set}
            }
          \end{equation}
          Note that this implies that for all $(s_1 ... s_{k-1})\in\dop{k-1}$

          \begin{equation*}
            (\p{k}X)_{s_1 ... s_{k-1}}=p X_{s_1 ... s_{k-1}}.
          \end{equation*}

    \end{itemize}

    We define $\cathd{n}$ to be the full subcategory of $\funcat{}{\cathd{n-1}}$ whose objects $X$ are such that
    \bigskip
    \begin{itemize}
      \item [(i)] \emph{Segal condition}

          \begin{equation*}
            X_s\cong\pro{X_1}{X_0}{s} \quad \mbox{for all} \; s \geq 2.
          \end{equation*}
    In particular this implies that $X\in \Cat(\Gpd(\cat{n-2})) =\Gpd(\cat{n-1})$ and the $n^{th}$ direction in $X$ is groupoidal.
    \bigskip
      \item [(ii)] \emph{Truncation functor}

       \nid The functor
      \begin{equation*}
        \op{n-1}:\cathd{n}\subset \funcat{}{\cathd{n-1}}\rw\funcat{}{\cathd{n-2}}
      \end{equation*}
      restricts to a functor
      \begin{equation*}
        \p{n}:\cathd{n}\rw\cathd{n-1}
      \end{equation*}
     Note that this implies that
      \begin{equation*}
        (\p{n}X)_{s_1 ... s_{n-1}}=p X_{s_1 ... s_{n-1}}
      \end{equation*}
      for all $s_1, ..., s_{n-1}\in\dop{n-1}$ and that the following diagram commutes

          \begin{equation}\label{eq-p-fun-def-a}
            \xymatrix{
            \cathd{n} \ar^{J_{n-1}}[rrr] \ar_{p^{(n)}}[d] &&& \funcat{n-1}{\Cat} \ar^{\bar p}[d]\\
            \cathd{n-1} \ar_{\N{n-1}}[rrr] &&& \funcat{n-1}{\Set}
            }
          \end{equation}
     \end{itemize}
\end{definition}
\mk
\label{C2}
\begin{definition}\label{def-hom-dis-ncat-2}

\index{Discretization!- map}

    Let $\di{n}$ be as in Definition \ref{dn} and $X\in \cathd{n}$. Denote by $\zgu{n}_X:X\rw \di{n}\p{n}X$ the morphism given by
    \begin{equation*}
        (\zgu{n}_X)_{s_1...s_{n-1}} :X_{s_1...s_{n-1}} \rw d p X_{s_1...s_{n-1}}
    \end{equation*}
    for all  $(s_1,...,s_{n-1})\in \dop{n-1}$, where $p X_{s_1...s_{n-1}}= q X_{s_1...s_{n-1}}$ since $X_{s_1...s_{n-1}}$ is a groupoid. Denote by
    \begin{equation*}
        X^d =\di{n}\di{n-1}...\di{1}\p{1}\p{2}...\p{n}X
    \end{equation*}
    and by $\zg\lo{n}$ the composite
    \begin{equation*}
        X\xrw{\zgu{n}}\di{n}\p{n}X \xrw{\di{n}\zgu{n-1}} \di{n}\di{n-1}\p{n-1}\p{n}X \rw \cdots \rw X^d\;.
    \end{equation*}
    We call $\zg\lo{n}$ the discretization map of $X$.
\end{definition}
\begin{notation}\label{not-fiber}
    Given $X\in\cathd{n}$, for each $a,b\in X_0^d$ denote by $X(a,b)$ the fiber at $(a,b)$ of the map
    \begin{equation*}
        X_1 \xrw{(d_0,d_1)} X_0\times X_0 \xrw{\zg\lo{n}\times\zg\lo{n}} X_0^d\times X_0^d\;.
    \end{equation*}
    The object $X(a,b)\in \cathd{n-1}$ should be thought of as a hom-$(n-1)$-category. \index{Hom-$(n-1)$-category}
\end{notation}

\begin{definition}\label{def-hom-dis-ncat-3} \index{n-equivalences}
Define inductively $n$-equivalences in $\cathd{n}$. For $n=1$, a 1-equivalence is an equivalence of categories. Suppose we defined $\nm$-equivalences in $\cathd{n-1}$. Then a map $f:X\rw Y$ in $\cathd{n}$ is an $n$-equivalence if

\begin{itemize}
  \item [a)] for all $a,b \in X_0^d$,
   \begin{equation*}
   f(a,b):X(a,b) \rw Y(fa,fb)
   \end{equation*}
    is a $\nm$-equivalence.
  \item [b)]$\p{n}f$ is a  $\nm$-equivalence.

  \end{itemize}

\end{definition}
\begin{remark}\label{rem-hom-dis-ncat}
By definition, the functor $\p{n}$ sends $n$-equivalences to $\nm$-equivalences. We observe that $\p{n}$ commutes with pullbacks over discrete objects. In fact, if $X\rw Z \lw Y$ is a diagram in $\cathd{n}$ with $Z$ discrete and $X\tiund{Z}Y\in\cathd{n}$, by Definition \ref{def-hom-dis-ncat}, for all $(s_1...s_{n-1})\in \dop{n-1}$,
\begin{align*}
    & (\p{n}(X\tiund{Z}Y))_{s_1...s_{n-1}}=p(X_{s_1...s_{n-1}}\tiund{Z}Y_{s_1...s_{n-1}})=\\
   =\  & p X_{s_1...s_{n-1}}\tiund{p Z} p Y_{s_1...s_{n-1}}=(\p{n}X \tiund{\p{n}Z} \p{n}Y)_{s_1...s_{n-1}}
\end{align*}
where we used the fact (Lemma \ref{lem-p-pres-fib-pro}) that $p$ commutes with pullbacks over discrete objects. Since this holds for each $s_1...s_{n-1}$ we conclude that
\begin{equation*}
    \p{n}(X\tiund{Z}Y) \cong \p{n} X \tiund{\p{n} Z} \p{n} Y\;.
\end{equation*}
\end{remark}
\begin{example}\label{ex-hom-disc-neq-rel}
    Let $X\in\cathd{2}$; then $\p{2}X$ is the equivalence relation associated to the surjective map of sets (in the sense of Definition \ref{def-int-eq-rel})
     \begin{equation*}
       \zg:pX_{0*}=(\p{2}X)_0\rw p(\p{2}X)=X^d
     \end{equation*}
     and $X$ has the form
\tiny{
\begin{equation*}
\xymatrix@C10pt{
% First row
% first term empty
& X_{10}\tiund{(\tens{pX_{0*}}{X^d})} X_{10}\tiund{(\tens{pX_{0*}}{X^d})} X_{10} \ar^{}[rr]<1.5ex> \ar^{}[rr]<0ex> \ar^{}[d] &&
X_{00}\tiund{pX_{0*}}X_{00}\tiund{pX_{0*}}X_{00} \ar^{}[ll]<1.5ex>  \ar^{}[d]\\
%% Second row
\cdots  \ar^{}[r] \ar_{}[d]<-2ex>\ar^{}[d] &
X_{10}\tiund{(\tens{pX_{0*}}{X^d})}X_{10}\ar^{}[rr]<1.5ex> \ar^{}[rr]<0ex>  \ar_{}[d]<-2ex>\ar^{}[d] &&
\tens{X_{00}}{pX_{0*}}\ar_{}[d]<-2ex>\ar^{}[d] \ar^{}[ll]<1.5ex> \\
%% Third row
\tens{X_{10}}{X_{00}} \ar^{}[r] \ar_{}[u]<-2ex> & X_{10} \ar^{}[rr]<1.5ex> \ar^{}[rr]<0ex>  \ar_{}[u]<-2ex> && X_{00}\ar_{}[u]<-2ex> \ar^{}[ll]<1.5ex>
}
\end{equation*}}
\end{example}

\nid The vertical structure is groupoidal and the horizontal  nerve $\Nu{1}X\in\funcat{}{\Cat}$ has in each component an equivalence relation. The horizontal structure is not in general groupoidal; however $\p{2}X$ is an equivalence relation, so in particular a groupoid. This means that the horizontal arrows in the double category $X$ have inverses after dividing out by the double cells. This structure is a special case of what called in \cite{PP} a groupoidal weakly globular double category.

\begin{remark}
  Homotopically discrete $n$-fold categories are more general than the homotopically discrete $n$-fold groupoids of \cite{BP}. In particular, they are \nfol categories but not in general \nfol groupoids since only some but not all of the $n$ different simplicial directions in the structure are required to be groupoidal.

   This added generality makes them more suitable to construct  the notion of weakly globular \nfol category. In particular, the way weakly globular \nfol category arise as strictification of Segalic pseudo-functors (see Theorem \ref{the-strict-funct}) requires this added generality to the notion of homotopically discrete $n$-fold category compared to the homotopically discrete $n$-fold groupoids of \cite{BP}.
\end{remark}

\section{Properties of homotopically discrete $\pmb{n}$-fold categories.}\label{sbs-proper-nfol}
In this section we establish the main properties of homotopically discrete \nfol categories.
In Lemma \ref{lem-pos-grou-hom-disc} we show that $\cathd{n}$ can be viewed as a diagram of equivalence relations.

An equivalence relation is categorically equivalent to a discrete category, so in particular a categorical equivalence between two equivalence relations is detected by an isomorphism of their discretizations. Similarly one expects a homotopically discrete $n$-fold category to be equivalent to a discrete structure in a higher categorical sense, and their higher categorical equivalences to be detected by isomorphisms of their discretizations. We show in Lemma \ref{lem-neq-hom-disc} that this is indeed the case. This criterion will be used throughout this work.

Using this characterization we show that every  homotopically discrete \nfol category $X$ is $n$-equivalent to a discrete $n$-fold category $X^d$ via the discretization map $\zg\lo{n}$ of Definition \ref{def-hom-dis-ncat-2}.

 Together with the good behavior of homotopically discrete \nfol categories with respect to pullbacks over discrete objects (Lemma \ref{lem-copr-hom-disc}), this implies that the induced Segal maps in a homotopically discrete \nfol category are $(n-1)$-equivalences (Proposition \ref{pro-ind-map-hom-disc}). We will see in Chapter \ref{chap4} that this makes homotopically discrete \nfol categories a subcategory of weakly globular \nfol categories.

\begin{lemma}\label{lem-pos-grou-hom-disc}
   The functor $J_n:\cat{n}\rw \funcat{n-1}{\Cat}$ restricts to a functor
   \begin{equation*}
    J_n:\cathd{n}\rw \funcat{n-1}{\cathd{}}\;.
   \end{equation*}
\end{lemma}
\begin{proof}
By induction on $n$. For $n=2$ if $X\in\cathd{2}$ then by definition $X_s\in\cathd{}$ for all $s\geq 0$. Suppose the lemma holds for $n-1$ and let $X\in \cathd{n}$. Then for all $s_1\geq 0$, $X_{s_1}\in\cathd{n-1}$ so by induction hypothesis
\begin{equation*}
    (X_{s_1})_{s_2...s_{n-1}}=X_{s_1...s_{n-1}}\in \cathd{}\;.
\end{equation*}
\end{proof}
\begin{lemma}\label{lem-neq-hom-disc}
    A map $f:X\rw Y$ in $\cathd{n}$ is a $\nequ$ if and only if $X^d \cong Y^d$.
\end{lemma}
\begin{proof}
By induction on $n$. For $n=1$, $f$ is a map of equivalence relations, so the statement is true by Lemma \ref{lem-p-pres-fib-pro}. Suppose the lemma holds for $(n-1)$ and let $f:X\rw Y$ be a $\nequ$ in $\cathd{n}$. Then by definition $\p{n}f$ is a $\nm$-equivalence; therefore by induction hypothesis
\begin{equation*}
    X^d=(\p{n}X)^d\cong(\p{n}Y)^d=Y^d\;.
\end{equation*}
Conversely, suppose that $f:X\rw Y$ is such that $X^d\cong Y^d$. This is the same as $(\p{n}X)^d=(\p{n}Y)^d$, so, by induction, $\p{n}f$ is a $(n-1)$-equivalence. This implies that, for each $a,b\in X^d$, $(\p{n}f)(a,b)$ is a $(n-2)$-equivalence. But
\begin{equation*}
    (\p{n}f)(a,b)=(\p{n-1}f)(a,b)
\end{equation*}
so $(\p{n-1}f)(a,b)$ is a $(n-2)$-equivalence. This implies that
\begin{equation*}
    X(a,b)^d = (\p{n-1}X(a,b))^d \cong (\p{n-1}Y(fa,fb))^d = Y(fa,fb)^d\;.
\end{equation*}
By induction hypothesis, we deduce that
\begin{equation*}
    f_{(a,b)}:X(a,b) \rw Y(fa,fb)
\end{equation*}
is a $\nm$-equivalence. We conclude that $f$ is a $\nequ$.
\end{proof}
\begin{corollary}\label{gamma}
  Let $X\in\cathd{n}$. Then the maps $\zg\up{n}:X\rw \di{n}\p{n}X$ and $\zg\lo{n}:X\rw X^d$ are $n$-equivalences.
\end{corollary}
\begin{proof}
This follows from Lemma \ref{lem-neq-hom-disc} since
\begin{equation*}
  X^d\cong (\p{n}X)^d\cong (X^d)^d.
\end{equation*}
\end{proof}
\begin{remark}\label{rem-neq-hom-disc}
    It follows from Lemma \ref{lem-neq-hom-disc} that $\nequ$s in $\cathd{n}$ have the 2-out-of-3 property.
\end{remark}
\begin{lemma}\label{lem-copr-hom-disc}
    Let $X\xrw{f} Z \xlw{g}Y$ be a diagram in $\cathd{n}$ with $Z$ discrete. Then
    \begin{itemize}
      \item [a)] $X \coprod Y \in \cathd{n}$.\mk

      \item [b)] $X \times Y \in \cathd{n}$.\mk

      \item [c)] $X \tiund{Z} Y \in \cathd{n}$ \quad \text{and} \quad $(X\tiund{Z}Y)^d=X^d\tiund{Z^d}Y^d$.\mk
    \end{itemize}
\end{lemma}
\begin{proof}\

\mk
a) By induction on $n$. It is clear for $n=1$. Suppose it holds for $n-1$ and let $X,Y\in \cathd{n}$. Since $\cathd{n}\subset\funcat{}{\cathd{n-1}}$ and coproducts in functor categories are computed pointwise, for each $s\geq 0$ we have, by induction hypothesis
\begin{equation*}
    (X \cop Y)_s=X_s \cop Y_s\in \cathd{n-1}\;.
\end{equation*}
Since $p$ commutes with coproducts, the same holds for $\p{n}$, thus by induction hypothesis
\begin{equation*}
    \p{n}(X \cop Y)= \p{n} X \cop \p{n} Y\in \cathd{n-1}\;.
\end{equation*}
this proves that $X \cop Y\in \cathd{n}$.

\mk
b) By induction on $n$. It is clear for $n=1$; suppose it holds for $n-1$. Then for each $s\geq 0$, by induction hypothesis
\begin{equation*}
    (X\times Y)_s=X_s \times Y_s\in \cathd{n-1}\;.
\end{equation*}
Since $\p{n}$ commutes with pullbacks over discrete objects (see Remark \ref{rem-hom-dis-ncat}) and therefore with products, by the induction hypothesis
\begin{equation*}
    \p{n}(X\times Y)=\p{n}X\times \p{n}Y\in \cathd{n-1}\;.
\end{equation*}
This proves that $X \times Y\in\cathd{n}$.

\mk
c) Since $Z$ is discrete.
\begin{equation*}
    X\tiund{Z}Y=\underset{c\in Z}{\cop}X(c)\times Y(c)
\end{equation*}
where $X(c)$ (resp. $Y(c)$) is the pre-image of $c$ under $f$ (resp. $g$). Since $X(c), \; Y(c) \in \cathd{n}$,  from a) and b) it follows that $X\tiund{Z}Y\in\cathd{n}$. Since by Remark \ref{rem-hom-dis-ncat} $\p{n}$ commutes with pullbacks over discrete objects for all $n$, we have
\begin{align*}
    &(X\tiund{Z}Y)^d= p\cdots \p{n}(X\tiund{Z}Y)= \\
    & =p\cdots \p{n}X \tiund{p\cdots \p{n}Z} p\cdots \p{n}Y=X^d\tiund{Z^d}Y^d\;.
\end{align*}
\end{proof}
Given $X\in\cathd{n}$, since $X_0^d$ is discrete and $X_1\in\cathd{n-1}$, by Lemma \ref{lem-copr-hom-disc}, for all $s\geq 2$,
\begin{equation*}
    \pro{X_1}{X_0^d}{s}\in\cathd{n-1}\;.
\end{equation*}
We can therefore consider the induced Segal maps
\begin{equation*}
    \hmu{s}:X_s= \pro{X_1}{X_0}{s}\rw \pro{X_1}{X_0^d}{s}
\end{equation*}
(see Definition \ref{def-ind-seg-map}). Using Lemma \ref{lem-neq-hom-disc} we show next that the induced Segal maps in a homotopically discrete \nfol category are $\nm$-equivalences.
\begin{proposition}\label{pro-ind-map-hom-disc} \index{Induced Segal maps condition}
    Let $X\in\cathd{n}$. For each $s\geq 2$ the induced Segal maps
    \begin{equation*}
        \hat\mu_s:\pro{X_1}{X_0}{s}\rw \pro{X_1}{X^d_0}{s}
    \end{equation*}
    are $\nm$-equivalences.
\end{proposition}
\begin{proof}
We show this for $s=2$, the case $s>2$ being similar.
By Lemma \ref{lem-neq-hom-disc} it is enough to show that
\begin{equation}\label{eq-ind-map-hom-disc}
    (\tens{X_1}{X_0})^d\cong (\tens{X_1^d}{X^d_0})\;.
\end{equation}
Denote
\begin{equation*}
  \begin{split}
 & \p{j,n-1}=\p{j}...\p{n-1} \; \; \mathrm{for} \; \; 1\leq j\leq n-1, \\
 & \p{n-1,n-1}=\p{n-1}.
  \end{split}
\end{equation*}
 We claim that
\begin{equation}\label{eq2-ind-map-hom-disc}
    \p{j,n-1}(\tens{X_1}{X_0})=\tens{\p{j,n-1}X_1}{\p{j,n-1}X_0}\;.
\end{equation}
We prove this by induction on $n$. When $n=2$, $X\in\catwg{2}$ so that
 \begin{equation*}
   p(\tens{X_1}{X_0})\cong p(\tens{X_1}{X_0^d}).
 \end{equation*}
 Suppose, inductively, the claim holds for $n-1$. Since $\p{n}X\in\cathd{n-1}$,
\begin{equation*}
    \p{n-1}(\tens{X_1}{X_0})= \tens{\p{n-1}X_1}{\p{n-1}X_0}\;.
\end{equation*}
By induction hypothesis applied to $\p{n}X$ we therefore obtain
\begin{equation*}
\begin{split}
    & \p{j,n-1}(\tens{X_1}{X_0})=\p{j,n-2}\p{n-1}(\tens{X_1}{X_0})= \\
   = \,& \p{j,n-2}(\tens{\p{n-1}X_1}{\p{n-1}X_0})= \\
   =\, & \tens{\p{j,n-2}\p{n-1}X_1}{\p{j,n-2}\p{n-1}X_0} = \\
   =\, & \tens{\p{j,n-1}X_1}{\p{j,n-1}X_0}\;.
\end{split}
\end{equation*}
This proves \eqref{eq2-ind-map-hom-disc}. In the case $j=1$ we obtain
\begin{equation*}
    (\tens{X_1}{X_0})^d=\tens{X_1^d}{X_0^d}\;.
\end{equation*}
Since $\p{n}$ commutes with pullbacks over discrete objects (see Remark \ref{rem-hom-dis-ncat}) we have
\begin{equation*}
\tens{X_1^d}{X_0^d}=(\tens{X_1}{X_0^d})^d
\end{equation*}
so that, from above, we conclude
\begin{equation*}
    (\tens{X_1}{X_0})^d \cong (\tens{X_1}{X_0^d})^d
\end{equation*}
as required.

\end{proof}
\section{Homotopically discrete $\pmb{n}$-fold categories and $\pmb{0}$-types}\label{sec-ho-as} \index{Zero-types}
We now discuss the homotopical significance of the category $\cathd{n}$. We introduce a classifying space functor from $\cathd{n}$ to spaces and we show that the classifying space of a homotopically discrete \nfol category is a $0$-type.
\begin{definition}\label{def-class-sp-funct}
    The classifying space functor is the composite \index{Classifying space}
    \begin{equation*}
         B:\cathd{n}\xrw{\Nb{n}}\funcat{n}{\Set}\xrw{Diag_n}\funcat{}{\Set}
    \end{equation*}
    where  $Diag_n$ denotes the multi-diagonal defined by \index{Multi-diagonal}
    \begin{equation*}
    (Diag_n Y)_k=Y_{k\oset{n}{...}k}
   \end{equation*}
   for $Y\in\funcat{n}{\Set}$ and $k\geq 0$.
\end{definition}
\begin{proposition}\label{pro-sign-cathd}
    If $X\in\cathd{n}$, $B\zg_X:BX\rw BX^d$ is a weak homotopy equivalence. In particular, $BX$ is a 0-type with $ \pi_{i}(BX,x)=0$ for $i>0$ and $\pi_{0}BX= UX^d$ where $UX^d$ is the set underlying the discrete \nfol category $X^d$.
\end{proposition}
\begin{proof}
By induction on $n$. For $n=1$, $X$ is a groupoid with no non-trivial loops, hence, since groupoids are models of $1$-types,
 $\pi_{i}(BX,x)=0$ for $i>0$ while $\pi_0 BX=UX^d$; suppose the statement holds for $\nm$.

The functor $B$ is also the composite
\begin{equation*}
    B:\cathd{n}\xrw{N_{1}}\funcat{}{\cathd{n-1}}\xrw{\ovl{B}} \funcat{}{\funcat{}{\Set}}\xrw{Diag_2}\funcat{}{\Set}\;.
\end{equation*}
Thus $B\zg_{X}$ is obtained by applying $Diag_2$ to the map of bisimplicial sets $N_1\ovl{B}\zg_{X}$. For each $s\geq 0$ the latter is given by
\begin{equation*}
    (N_1\ovl{B}\zg_{X})_s=B\zg_s: B X_s \rw B X_s^d=B(\p{2}...\p{n}X)_s
\end{equation*}
and this is a weak homotopy equivalence by induction hypothesis.

A map of bisimplicial sets which is a levelwise weak homotopy equivalence induces a weak homotopy equivalence of diagonals (see \cite{Jard}). \index{Diagonal} Hence
\begin{equation*}
    Diag_2 N_1\ovl{B}\zg_X=B\zg_X
\end{equation*}
is a weak homotopy equivalence, as required. Thus $BX$ is weakly homotopy equivalent to $B(\p{2}...\p{n}X)$, which is a 0-type since $\p{2}...\p{n}X\in\cathd{}$. Further,
\begin{equation*}
  \pi_{0}BX\cong \pi_0 B(\p{2}...\p{n}X) \cong Up\p{2}...\p{n}X \cong UX^d.
\end{equation*}

\end{proof}

\section{Higher equivalence relations}\label{highrel}
\index{Higher equivalence relations}
In this section we give a different description of weakly globular \nfol categories via a notion of iterated internal equivalence relation. The notion of internal equivalence relation (Definition \ref{def-int-eq-rel}) associated to a morphism $f:A\rw B$ in a category $\clC$ with finite limits is known.\index{Internal! -equivalence relation} When $\clC=\Set$ and $f:A\rw B$ is surjective, this affords the category $\cathd{}$ and the category $A[f]\in\cathd{}$ corresponding to $f$ has set of connected components given by $q A[f]= p A[f]=B$.

We define $\eqr{n}$ by iterating this notion in $(n-1)$-fold categories in such a way that the target $Y$ of the morphism $f:X\rw Y$ in $\cat{n-1}$ belongs to $\eqr{n-1}$ and $\Nb{n-1}f$ is a levelwise surjection in $\Set$. This surjectivity condition ensures that there is a functor
\begin{equation*}
  \p{n}:\eqr{n}\rw\eqr{n-1}
\end{equation*}
 with
 \begin{equation*}
   \p{n}X[f]=Y.
 \end{equation*}

In Theorem \ref{the-hom-disc-neq-rel} we reconcile the definition of $\eqr{n}$ with the definition of $\cathd{n}$ of the previous section.

\begin{definition}\label{def-int-eq-rel}
Let $A\rw B$ be a morphism in a category $\clC$ with finite limits. The diagonal map defines a unique section
\begin{equation*}
 s:A\rw\ata
\end{equation*}
  (so that ${p_{1}s=\Id_{A}=p_{2}s}$ where ${\ata}$ is the pullback of  ${A\xrw{f}B\xlw{f}A}$ and ${p_{1},p_{2}:\ata\rw A}$ are the two projections). The commutative diagram
\begin{equation*}
    \xymatrix{
    \ata \ar[rr]^{p_{1}} \ar[d]_{p_{2}} && A \ar[d]_{f} &&
    \ata \ar[ll]_{p_{2}} \ar[d]^{p_{1}}\\
    A \ar[rr]_{f} && B && A \ar[ll]^{f}
    }
\end{equation*}
defines a unique morphism
 \begin{equation*}
 {m:(\ata)\tiund{A}(\ata)\rw\ata}
 \end{equation*}
  such that ${p_{2}m=p_{2}\pi_{2}}$ and $p_{1}m=p_{1}\pi_{1}$ where $\pi_{1}$ and ${\pi_{2}}$ are the two projections. We  denote by $\Af$ the following object of $\Cat(\clC)$
\begin{equation*}
\xymatrix{
(\ata)\tiund{A}(\ata)\ar[rr]^<<<<<<<<{m} && \ata \ar@<1.5ex>[rr]^{p_{1}}
\ar@<-.5ex>[rr]^{p_{2}} && A \ar@<1.5ex>[ll]^{s} }
\end{equation*}
\end{definition}
\begin{lemma}\label{inter-eq}
  Let ${\Af}$ be as in Definition \ref{def-int-eq-rel}. Then ${\Af}$ is an internal groupoid in $\clC$.
\end{lemma}
\begin{proof}
We need to prove the axioms of internal groupoid as in Section \ref{sbs-nint-cat}. Axioms (1) of (2) in Definition \ref{def-intercat} hold by construction. As for (3), denote by
\begin{equation*}
  \pi_i : \tens{(\tens{A}{B})}{A}\qquad i=1,2
\end{equation*}
the two projections. Then
\begin{align*}
    & f p_1 m \begin{pmatrix}\Id \\ s p_1 \\ \end{pmatrix} = f p_1 \pi_1 \begin{pmatrix}\Id \\ s p_1 \\ \end{pmatrix} = f p_1 \Id = f p_1\\
    & f p_2 m \begin{pmatrix}\Id \\ s p_1 \\ \end{pmatrix} = f p_2 \pi_2 \begin{pmatrix}\Id \\ s p_2 \\ \end{pmatrix} = f p_2 s p_2 \Id = f p_2\;.
\end{align*}
Therefore
\begin{equation*}
  m\cirsm \begin{pmatrix}\Id \\ s p_1 \\ \end{pmatrix}= \Id_{\tens{A}{B}}
\end{equation*}
which is the first half of axiom (3). The second half is proved similarly.

To show (4) denote by
\begin{equation*}
  r_i: \tens{(\tens{A}{B})}{A}\tiund{A}(\tens{A}{B})\rw (\tens{A}{B})\qquad i=1,2,3
\end{equation*}
the three projections and by
\begin{equation*}
\begin{split}
   & r_{12}:  \tens{(\tens{A}{B})}{A}\tiund{A}(\tens{A}{B})\rw \tens{(\tens{A}{B})}{A} \\
   & r_{23}:  \tens{(\tens{A}{B})}{A}\tiund{A}(\tens{A}{B})\rw \tens{(\tens{A}{B})}{A}
\end{split}
\end{equation*}
the projections
\begin{equation*}
  r_{12}=r_1\times r_2\times\Id \qquad r_{23}=\Id\times r_2\times r_3\;.
\end{equation*}
To show axim (3) in Definition \ref{def-intercat}, it is enough to prove that
\begin{equation}\label{eq1-axiom}
   f p_2 m(\Id\times m) = f p_2 m(m\times \Id)
\end{equation}
\begin{equation}\label{eq2-axiom}
   f p_1 m(\Id\times m) = f p_1 m(m\times \Id)\;.
\end{equation}
We calculate
\begin{equation*}
  f p_2 m r_{23}= f p_2 \pi_2 r_{23}= f p_2 r_3 = f p_2 \Id r_3
\end{equation*}
which is \eqref{eq1-axiom}. Also,
\begin{equation*}
  f p_1 m r_{12}= f p_1 \pi_1 r_{12}= f p_1 r_1 = f p_1 \Id r_1
\end{equation*}
which is \eqref{eq2-axiom}.

Thus $A[f]\in\Cat \clC$. To show that $A[f]$ is an internal groupoid, we identify the inverses as given by the map
\begin{equation*}
  i : \tens{A}{B} \rw \tens{A}{B}
\end{equation*}
determined by the diagram
\begin{equation*}
\xymatrix@R=30pt{
\tens{A}{B} \ar@/_1pc/[ddr]_{p_2} \ar@/^1pc/[drr]^{p_1}
\ar@{->}[dr]|-{\;i\;} \\
& \tens{A}{B} \ar[d]_{p_1} \ar[r]_(0.6){p_2}
& A \ar[d]^f \\
& A \ar[r]_f & B
}
\end{equation*}
To show that $m(\Id,i)=s p_1$ it is enough to show that
\begin{equation}\label{eq3-axiom}
  f p_i m(\Id,i) =f p_i s p_1
\end{equation}
for $i=1,2$. We calculate
\begin{align*}
    & f p_1 m(\Id,i)=f p_1 \pi_1(\Id,i)=f p_1 \Id=f p_1=f p_1 s p_1 \\
    & f p_2 m(\Id,i)=f p_2 \pi_2(\Id,i)=f p_2 i =f p_1=f p_2 s p_1\;.
\end{align*}
Thus \eqref{eq3-axiom} holds. The proof that $m(i,\Id)=s p_2$ is similar.
\end{proof}
\begin{remark}\label{rem-cathd}
If $X\in\cathd{}$ then
\begin{equation*}
  X=X_0 [\gamma]
\end{equation*}
 where $\gamma:X_0\rw pX$ is a surjective map of sets.
\end{remark}

In what follows $ \Nb{n-1}:\cat{n-1}\rw \funcat{n-1}{\Set}$ is the multinerve and we denote for all $\us\in\Dmenop$
\begin{equation*}
  (\Nb{n-1}X)_{\us}=X_{\us}.
\end{equation*}

\begin{definition}\label{def1-int-eq-rel}
    We define
     \begin{equation*}
       \eqr{n}\subset \cat{n}
     \end{equation*}
      by induction on $n$. For $n=1$, $\eqr{1}=\cathd{}$. Suppose, inductively, we defined $\eqr{n-1}\subset\cat{n-1}$ and let $f:X\rw Y$ be a morphism in $\cat{n-1}$ with $Y\in\eqr{n-1}$ such that, for all $\us\in\Dmenop$, $f_{\us}:X_{\us}\rw Y_{\us}$ is surjective. We define $\eqr{n}$ to be the full subcategory of $\cat{n}$ whose objects have the form $X[f]$, where $X[f]$ is as in Definition \ref{def-int-eq-rel}.
\end{definition}
\begin{remark}\label{rem-int-eq-rel}
    Let $\za:X[f]\rw X'[f']$ be a morphism in $\cathd{n}$, with $f:X\rw Y$ and $f':X'\rw Y'$.
    \begin{equation*}
        q X_{\us}[f_{\us}]= p X_{\us}[f_{\us}]= Y_{\us}
    \end{equation*}
    and there is a functor
    \begin{equation*}
        X_{\us}[f_{\us}]\rw d Y_{\us}\;.
    \end{equation*}
    We therefore have a commuting diagram in $\Cat$
    \begin{equation*}
        \xymatrix@C=40pt{
        X_{\us}[f_{\us}] \ar[r] \ar_{\za_{\us}}[d] & d Y_{\us} \ar^{\ovl{\za}_{\us}}[d]\\
        X'_{\us}[f'_{\us}] \ar[r] & d Y'_{\us}
        }
    \end{equation*}
    inducing a commuting diagram in $\Set$
    \begin{equation*}
        \xymatrix@C=40pt{
        X_{\us} \ar^{f_{\us}}[r] \ar_{\za_{\us}}[d] &  Y_{\us} \ar^{\ovl{\za}_{\us}}[d]\\
        X'_{\us} \ar^{f'_{\us}}[r] &  Y'_{\us}\;.
        }
    \end{equation*}
    Since this holds for all $\us$, we conclude that there is a commuting diagram in $\cat{n-1}$
    \begin{equation*}
        \xymatrix@C=40pt{
        X \ar^{f}[r] \ar_{\za}[d] &  Y \ar^{\ovl{\za}}[d]\\
        X' \ar^{f'}[r] &  Y'\;.
        }
    \end{equation*}

\end{remark}

\bigskip
We now show that the categories  $\eqr{n}$ and $\cathd{n}$ are isomorphic. The proof of this result relies on Lemma \ref{lem-pos-grou-hom-disc} saying that a homotopically discrete \nfol category is a diagram of equivalence relations in direction $n$. This allows to write an object of $\cathd{n}$ as an internal equivalence relation (in the sense of Definition \ref{def-int-eq-rel}) in $\cat{n-1}$ corresponding to a morphism in $\cat{n-1}$ whose target is homotopically discrete and whose multinerve is a levelwise surjection in $\Set$. An inductive arguments then reconciles this with the definition of  $\eqr{n}$.
\begin{theorem}\label{the-hom-disc-neq-rel}
   There is an isomorphism of categories
   \begin{equation*}
     \eqr{n} \cong \cathd{n}.
   \end{equation*}
\end{theorem}
\begin{proof}
By induction on $n$. For $n=1$, it holds by definition. Suppose this is true for all $k\leq n-1$. Let $X[f]\in\eqr{n}$ with $f: X\rw Y$ a morphism in $\cat{n-1}$ and $Y\in\eqr{n-1}$. So by inductive hypothesis, $Y\in \cathd{n}$.

To show that $X[f]\in\cathd{n}$ we need to show that, for all $s_1\geq 0$, $(X[f])_{s_1}\in\cathd{n-1}$ and $\p{n}X[f]\in\cathd{n-1}$ where, for all $\us\in\Dmenop$,
\begin{equation*}
    (\p{n}X[f])_{\us}=p(X[f])_{\us}\;.
\end{equation*}
%
%%, $$ where $f_{s_1}:X_{s_1}\rw Y_{s_1}$ is such that $Y_{s_1}\in\cathd{n-2}=\eqr{n-2}$ (since $Y\in\cathd{n-1}=\eqr{n-1}$) and $(f_{s_1})_{s_2,\ldots,s_{n-1}}=f_{s_1,\ldots,s_{n-1}}$ is surjective.

For all $s_1\geq 0$
\begin{equation*}
    (X[f])_{s_1}=X_{s_1}[f_{s_1}]\;.
\end{equation*}
 where $f_{s_1}:X_{s_1}\rw Y_{s_1}$ in a morphism in $\cat{n-2}$; since $Y\in\eqr{n-1}$, by induction hypothesis $Y\in\cathd{n-1}$, thus by definition $Y_{s_1}\in\cathd{n-2}$ and, by induction hypothesis again, $Y_{s_1}\in\eqr{n-2}$. Further, $(f_{s_1})_{s_2...s_{n-1}}=f_{s_1s_2...s_{n-1}}$ is surjective.
Thus, by definition,
\begin{equation*}
  (X[f])_{s_1}=X_{s_1}[f_{s_1}]\in\eqr{n-1}\cong\cathd{n-1}.
\end{equation*}
Since $f_{\us}$ is surjective, we have
\begin{equation*}
    p X_{\us}[f_{\us}]=Y_{\us}
\end{equation*}
which implies $\p{n}X[f]=Y\in\cathd{n-1}$, as required.

Conversely, let $X\in\cathd{n}$. Consider the morphism in $\cat{n}$
\begin{equation*}
    \tilde{\xi_n}\zg_X\up{n}:\tilde{\xi_n} X\rw \tilde{\xi_n} \di{n}\p{n}X
\end{equation*}
where $\tilde{\xi_n}$ is as in Proposition \ref{pro-mult-ner}, that is $\tilde{\xi_n} X$ in an internal category in $\cat{n-1}$ in direction $n$, and $\di{n}$ is as in Definition \ref{dn}.
 At the object of objects level this gives a morphism in $\cat{n-1}$
\begin{equation*}
    (\tilde{\xi_n}\zg_X\up{n})_0:(\tilde{\xi_n} X)_0\rw \p{n} X\;.
\end{equation*}
Denote
\begin{equation*}
  f_n X=(\tilde{\xi_n}\zg_X\up{n})_0.
\end{equation*}
We claim that
\begin{equation}\label{eq1-hom-disc-neq-rel}
    X=(\tilde{\xi_n} X)_0[f_n X]\;.
\end{equation}
To show this, let $\ur=(s_1,...,s_{n-1})\in\dop{n-1}$, then
\begin{equation*}
  (f_n X)_{\ur} : (\tilde{\xi}_n X)_{0\ur}= X_{(\ur,0)} \rw p X_{\ur}
\end{equation*}
is surjective and by Lemma \ref{lem-pos-grou-hom-disc} $X_{\ur *}\in\cathd{}$. Hence by Remark \ref{rem-cathd}
\begin{equation*}
  X_{\ur*}=X_{(\ur,0)}[(f_n X)_{\ur}]=(\tilde{\xi}_n X)_{0\ur}[(f_n X)_{\ur}]\;.
\end{equation*}
Since this holds for each $\ur\in\dop{n-1}$, \eqref{eq1-hom-disc-neq-rel} follows.

By Lemma \ref{lem-pos-grou-hom-disc} the map
\begin{equation*}
(f_n X)_{\ur}= X_{(\ur,0)}\rw p X_{\ur} = (\p{n}X)_{\ur}
\end{equation*}
is surjective. Also, $\p{n}X\in\cathd{n-1}$ thus by inductive hypothesis $\p{n}X\in\eqr{n-1}$. By \eqref{eq1-hom-disc-neq-rel} and by definition we conclude that $X\in\eqr{n}$.

\end{proof}
%%

%%
%%%%%%%%%%%%%%%%%%%%%%%%%%%%%%%%%%%%%%%%%%%%%%%%%%%%%%%%%%%%%%%%%%%%%%%%%%%%%%%%%
\chapter{Weakly globular $\pmb{n}$-fold categories}\label{chap4}
%%
%%%%%%%%%%%%%%%%%%%%%%%%%%%%%%%%%%%%%%%%%%%%%%%%%%%%%%%%%%%%%%%%%%%%%%%%%%%%%%%

%
In this chapter we introduce the central higher categorical structure of this work, the category $\catwg{n}$ of weakly globular \nfol categories and we establish its main properties.

Weakly globular \nfol categories form a full subcategory of the category $\cat{n}$ of \nfol categories; they are therefore rigid structures in which there are compositions in $n$ different directions and all these compositions are associative and unital. The weakness in a weakly globular \nfol category is encoded by the weak globularity condition. The latter is formulated using the category $\cathd{n}$ of homotopically discrete \nfol categories introduced in Chapter \ref{chap3}.
In our approach the cells in dimension $k$ (for $0\leq k\leq n-2$) no longer form a set but have the higher categorical structure of homotopically discrete $(n-1-k)$-fold categories.

The underlying set of the discretizations of the homotopically discrete substructures in a weakly globular $n$-fold category play the role of sets of cells in the respective dimensions. We also impose additional conditions in the definition of weakly globular \nfol category to obtain well behaved compositions of higher cells.

In the case $n=2$, weakly globular double categories were introduced in joint work by the author in \cite{BP} and shown to be biequivalent to bicategories. The generalization to the case $n>2$ is much more complex.

This chapter is organized as follows. In Section \ref{sec-wg-nfold-categ} we inductively define weakly globular \nfol categories and $n$-equivalences between them. In Section \ref{sbs-prop-wg-nf-cat} we establish the main properties of weakly globular \nfol categories. We show in Proposition \ref{pro-crit-ncat-be-wg} b) a criterion for a \nfol category to be weakly globular playing a crucial role in the proof of the main result of the next chapter, Theorem \ref{the-strict-funct}, showing that weakly globular \nfol categories arise as strictification of certain types of pseudo-functors.

%%%%%%%%%%%%%%%%%%%%%%%%%%%%%%%%%%%%%%%%%%%%%%%%%%%%%%%%%%%%%%%%%%%%%%%%%%%%%%%%%%%%%%%%%%%%%%%%%%%%%%%%%%%

%%
\section{The definition of weakly globular $\pmb{n}$-fold categories}\label{sec-wg-nfold-categ}
In this section we define the category $\catwg{n}$ of weakly globular \nfol categories and we also define $n$-equivalences in this category. The definition uses the category $\cathd{n}$ of homotopically discrete \nfol categories defined in Chapter \ref{chap3}.

 \subsection{The idea of weakly globular $\pmb{n}$-fold categories}\label{subs-idea-catwg}
 The idea of the definition of the category $\catwg{n}$ is to build the structure by induction on dimension starting with the category $\Cat$ with equivalences of categories.

  At dimension $n$, the structure is a full subcategory of simplicial objects in $\catwg{n-1}$. Unraveling this definition, this affords an embedding
\begin{equation*}
        \xymatrix{
        \catwg{n}\ar@{^(->}^(0.35){J_n}[r] & \funcat{n-1}{\Cat}\;.
        }
\end{equation*}
The first condition for $X\in\funcat{}{\catwg{n-1}}$ to be an object of $\catwg{n}$ is the weak globularity condition that $X_0$ is homotopically discrete.

 The set underlying the discrete $(n-1)$-fold category $X_0^d$ plays the role of set of cells in dimension $0$. When $1\leq r\leq n-2$, the set underlying $(X)^d_{1\oset{r}{\cdots}10}$ corresponds to the set of $r$-cells.

    The next condition in the definition of $\catwg{n}$ is that the Segal maps
    \begin{equation*}
        X_k\rw \pro{X_1}{X_0}{k}
    \end{equation*}
are isomorphisms for all $k\geq 2$. Since each $X_k\in\cat{n-1}$, by the characterization of internal categories via the Segal condition (Proposition \ref{pro-ner-int-cat}) it follows that $X$ is an \nfol category.

We further require the induced Segal map condition stating that, for each $k\geq 2$, the maps in $\catwg{n-1}$
\begin{equation*}
    X_k \rw \pro{X_1}{X_0^d}{k}
\end{equation*}
are $(n-1)$-equivalences. This condition controls the compositions of higher cells and is the analogue of the Segal condition in the Tamsamani-Simpson model \cite{Ta}, \cite{Simp}.

 We finally require the existence of a truncation functor $\p{n}$ from $\catwg{n}$ to $\catwg{n-1}$ obtained by applying dimensionwise the isomorphism classes of object functor to the corresponding diagram in $\funcat{n-1}{\Cat}$. In the case $n=2$, this last condition is redundant. The effect of the functor $\p{n}$ on  weakly globular \nfol category $X$ is to quotient by the highest dimensional invertible cells of $X$.

The functor $\p{n}$ is used to define $n$-equivalences, thus completing the inductive step in the definition of $\catwg{n}$. The definition of $n$-equivalences is given in terms of two conditions: the first is a higher dimensional generalization of the notion of fully faithfulness of a functor, the second is a generalization of 'essentially surjective on objects'.

\subsection{The formal definition of the category $\pmb{\catwg{n}}$}
\begin{definition}\label{def-n-equiv}\index{Functor!- $\p{n}$}
    For $n=1$, $\catwg{1}=\Cat$ and $1$-equivalences are equivalences of categories. \index{Weakly globular!- $n$-fold category}

    Suppose, inductively, that we defined $\catwg{n-1}$ and $(n-1)$-equivalences. Then $\catwg{n}$ is the full subcategory of $\funcat{}{\catwg{n-1}}$ whose objects $X$ are such that
    \begin{itemize}
      \item [a)] \textsl{Weak globularity condition} $X_0\in\cathd{n-1}$.\mk \index{{Weak globularity condition}}
      \item [b)] \textsl{Segal maps condition} For all $k\geq 2$ the Segal maps are isomorphisms:\index{Segal maps condition}      %
      \begin{equation*}
        X_k\cong\pro{X_1}{X_0}{k}\;.
      \end{equation*}

      \item [c)] \textsl{Induced Segal maps condition} For all $k\geq 2$ the induced Segal maps \index{Induced Segal maps condition}

      \begin{equation*}
        \hmuk :X_k\rw\pro{X_1}{X^d_0}{k}
      \end{equation*}
      (induced by the map $\zg:X_0\rw X_0^d$) are $(n-1)$-equivalences.\mk

      \item [d)] \textsl{Truncation functor} There is a functor \index{Truncation functor}
       \begin{equation*}
         \p{n}:\catwg{n}\rw\catwg{n-1}
       \end{equation*}
        making the following diagram commute
      \begin{equation*}
        \xymatrix{
        \catwg{n} \ar^{J_n}[rr] \ar_{\p{n}}[d] && \funcat{n-1}{\Cat} \ar^{\ovl p}[d]\\
        \catwg{n-1} \ar^{\Nb{n-1}}[rr]  && \funcat{n-1}{\Set}
        }
    \end{equation*}
    \end{itemize}
    Given $a,b\in X_0^d$, denote by $X(a,b)$ the fiber at $(a,b)$ of the map
    \begin{equation*}
         X_1\xrw{(\pt_0,\pt_1)} X_0\times X_0 \xrw{\zg\times \zg}  X^d_0\times X^d_0\;.
    \end{equation*}
    The object $X(a,b)\in \catwg{n-1}$ should be thought of as hom-$(n-1)$-category. \index{Hom-$(n-1)$-category}
    We say that a map $f:X\rw Y$ in $\catwg{n}$ is an $n$-equivalence if \index{n-equivalences}
    \begin{itemize}
      \item [i)] For all $a,b\in X_0^d$
      \begin{equation*}
        f(a,b): X(a,b) \rw Y(fa,fb)
      \end{equation*}
      is an $(n-1)$-equivalence.\mk

      \item [ii)] $\p{n}f$ is an $(n-1)$-equivalence.
    \end{itemize}
    This completes the inductive step in the definition of $\catwg{n}$.
\end{definition}

\begin{remark}\label{rem-n-equiv}
    It follows by Definition \ref{def-n-equiv}, Definition \ref{def-hom-dis-ncat} and Proposition \ref{pro-ind-map-hom-disc} that $\cathd{n}\subset \catwg{n}$.
\end{remark}

\begin{figure}[ht]
  \centering
\begin{equation*}
\entrymodifiers={+++[o]}
 \xymatrix@R=0.5pc@C=20pt{
  \vdots  & \vdots & \vdots\\
\cdots \ar@<2ex>[r] \ar[r] \ar@<1ex>[r] \ar@<-2ex>[dd] \ar[dd]
    \ar@<-1ex>[dd]
    & \bm{\tens{X_{11}}{X_{10}}} \ar@<1ex>[r] \ar@<0ex>[r]  \ar@<-2ex>[dd] \ar[dd] \ar@<-1ex>[dd] \ar@<1ex>[l] \ar@<2ex>[l]
    & \textcolor[rgb]{1.00,0.00,0.00}{\bm{\tens{X_{01}}{X_{00}}}}  \ar@<-2ex>@[red][dd] \ar@[red][dd] \ar@<-1ex>@[red][dd] \ar@<1ex>[l]\\
    &&\\
    \bm{\tens{X_{11}}{X_{01}}}\ar@<2ex>[r] \ar[r] \ar@<1ex>[r] \ar@<-1ex>[dd] \ar@<0.ex>[dd] \ar@<-1ex>[uu] \ar@<-2.ex>[uu]
    & \bm{X_{11}} \ar@<1ex>[r] \ar@<0ex>[r]  \ar@<-1ex>[dd] \ar@<0.ex>[dd] \ar@<1ex>[l] \ar@<2ex>[l] \ar@<-1ex>[uu] \ar@<-2.ex>[uu]
    & \textcolor[rgb]{1.00,0.00,0.00}{\bm{X_{01}}}  \ar@<0.ex>@[red][dd] \ar@<-1ex>@[red][dd] \ar@<1ex>[l] \ar@<-1ex>@[red][uu] \ar@<-2.ex>@[red][uu] \\
    &&\\
    \bm{\tens{X_{10}}{X_{00}}} \ar@<2ex>[r] \ar[r] \ar@<1ex>[r] \ar@<-1ex>[uu]
    & \bm{X_{10}} \ar@<1ex>[r] \ar@<0ex>[r] \ar@<1ex>[l] \ar@<2ex>[l] \ar@<-1ex>[uu]
    & \textcolor[rgb]{1.00,0.00,0.00}{\bm{X_{00}}}\ar@<1ex>[l] \ar@<-1ex>@[red][uu]\\
}
\end{equation*}
  \caption{Corner of the double nerve of a weakly globular double category $X$}
  \label{corner2ng}
\end{figure}
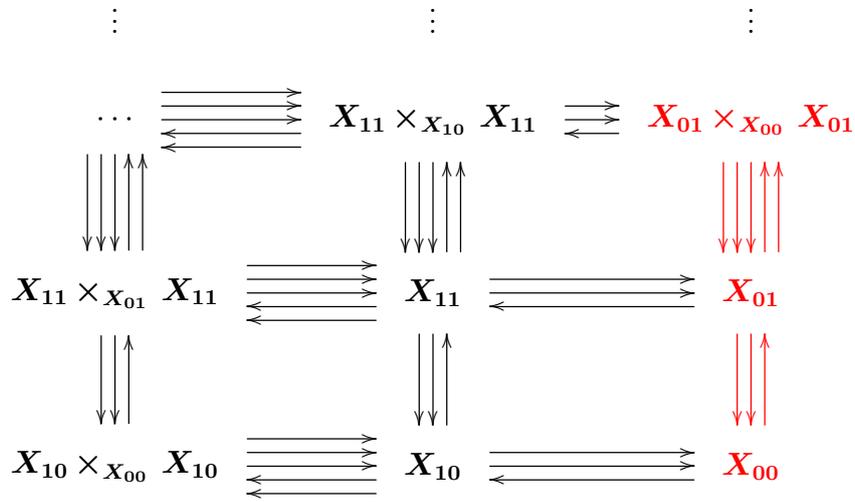

\begin{example}\label{ex1-n-equival}
\emph{Weakly globular double categories.}

\index{Weakly globular!- double category}

Let $X\in\catwg{2}$. Then, by definition, $X\in\funcat{}{\Cat}$ is such that
\begin{itemize}
  \item [a)] $X_0 \in \cathd{}$. \medskip

  \item [b)] For all $k\geq 0$ \; $X_k\cong \pro{X_1}{X_0}{k}$.\medskip

  \item [c)] For all $k\geq 2$ the induced Segal maps
  \begin{equation*}
   \hmu{k}:X_k\rw \pro{X_1}{X^d_0}{k}
  \end{equation*}
  are equivalences of categories.\medskip

  \item [d)] There is a functor $\p{2}:\catwg{2}\rw\Cat$ making the following diagram commute
  \begin{equation*}
    \xymatrix@C=40pt{
    \catwg{2} \ar^(0.4){J_2}[r] \ar_{\p{2}}[d]
    & \funcat{}{\Cat} \ar^{\bar p}[d] \\
    \Cat \ar_(0.4){N}[r]
    & \funcat{}{\Set}
    }
  \end{equation*}
\end{itemize}
Note that in the case $n=2$ condition d) is redundant. In fact, by b) and c) and the fact (Lemma \ref{lem-p-pres-fib-pro}) that $p$ sends equivalences to isomorphisms and commutes with pullbacks over discrete objects, given $X\in\funcat{}{\Cat}$ satisfying conditions  a),b) c), for each $k\geq 2$
\begin{equation*}
\begin{split}
    &  pX_k \cong p(\pro{X_1}{X^d_0}{k}) = \pro{p X_1}{X^d_0}{k} \cong\\
    & \cong \pro{p X_1}{p X^d_0}{k}\;.
\end{split}
\end{equation*}
Hence d) holds.

On page \pageref{page-corner2red}, Figure \ref{corner2ng} is a picture of the corner of $X\in\catwg{2}$, where the red structure is homotopically discrete. The corresponding geometric picture is on page \pageref{page-corner2red}, Figure \ref{corner2red}.
Condition c) also has a geometric interpretation. For this, note that, given $X\in\funcat{}{\Cat}$ satisfying a) and b), the induced Segal maps
\begin{equation*}
  \hmu{k}:\pro{X_1}{X_0}{k} \rw \pro{X_1}{X^d_0}{k}
\end{equation*}
are fully faithful. We show this for $k=2$, the case $k>2$ being similar. Given $(a,b),(c,d)\in \tens{X_{10}}{X_{00}}$ we have
\begin{equation*}
\begin{split}
    & (\tens{X_1}{X_0})\{(a,b),(c,d)\} \cong X_{11}(a,c)\tiund{X_{01}(\pt_0 a,\pt_0 c)}X_{11}(b,d) \cong \\
    & \cong X_{11}(a,c)\times X_{11}(b,d)
\end{split}
\end{equation*}
since $X_{01}(\pt_0 a,\pt_0 c)=\{\cdot\}$ as $X_{0*}\in\cathd{}$. Hence
\begin{equation*}
\begin{split}
    &  (\tens{X_1}{X_0})\{(a,b),(c,d)\} \cong \\
    &  \cong (\tens{X_1}{X^d_0}) \{\hmu{k}(a,b),\hmu{k}(c,d)\}\cong X_{11}(a,c) \times X_{11}(b,d)\;.
\end{split}
\end{equation*}
Since, when conditions a) and b) hold, $\hmu{k}$ is always fully and faithful, condition c) that $\hmu{k}$ it is an equivalence of categories is equivalent to the requirement that it is essentially surjective on objects.

An object of $\pro{X_1}{X_0^d}{k}$ is a staircase of horizontal arrows of length $k$ whose source and targets match up in vertical connected component, as in the following picture:
%%
%% Create 1-edge
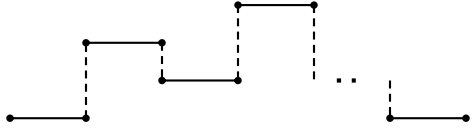
\tikzset{edge/.pic={
\filldraw (0,0)  circle[radius=0.035cm] -- (1,0)  circle[radius=0.035cm]; %% centered
}}
\begin{center}
\begin{tikzpicture}[thick,scale=1]]
\pic [black][scale=1] at (0,0) {edge};
\pic [black][scale=1] at (1,1) {edge};
\pic [black][scale=1] at (2,0.5) {edge};
\pic [black][scale=1] at (3,1.5) {edge};
\pic [black][scale=1] at (5,0) {edge};
\draw[thick, densely dashed] (1,0) -- (1,1);
\draw[thick, densely dashed] (2,1) -- (2,0.5);
\draw[thick, densely dashed] (3,1.5) -- (3,0.5);
\draw[thick, densely dashed] (4,1.5) -- (4,0.5);
\draw[thick, densely dashed] (5,0.5) -- (5,0);
\draw[ultra thick, loosely dotted] (4.3,0.5) -- (4.55,0.5);
\end{tikzpicture}
\end{center}

\bk
Essential surjectivity of $\hmu{k}$ means that this staircase can be lifted to a sequence of horizontally composable arrows through vertically invertible squares:
\bk

\tikzset{edge/.pic={
\filldraw (0,0)  circle[radius=0.04cm] -- (1,0)  circle[radius=0.04cm]; %% centered
}}
\begin{center}
\begin{tikzpicture}[thick,scale=1]]
\pic [black][scale=1] at (0,0) {edge};\pic [black][scale=1] at (0,2.5) {edge};
\node[rotate=-90,scale=1.2] at (0.5,1.5) () {$\bm\Rightarrow$};
\pic [black][scale=1] at (1,1) {edge};\pic [black][scale=1] at (1,2.5) {edge};
\node[rotate=-90,scale=1.2] at (1.5,1.8) () {$\bm\Rightarrow$};
\pic [black][scale=1] at (2,0.5) {edge};\pic [black][scale=1] at (2,2.5) {edge};
\node[rotate=-90,scale=1.2] at (2.5,1.5) () {$\bm\Rightarrow$};
\pic [black][scale=1] at (3,1.5) {edge};\pic [black][scale=1] at (3,2.5) {edge};
\node[rotate=-90,scale=1.2] at (3.5,2.0) () {$\bm\Rightarrow$};
\pic [black][scale=1] at (5,0) {edge};\pic [black][scale=1] at (5,2.5) {edge};
\node[rotate=-90,scale=1.2] at (5.5,1.5) () {$\bm\Rightarrow$};
\draw[ultra thick, loosely dotted] (4.3,2.5) -- (4.6,2.5);
\draw[thick, densely dashed] (0,0) -- (0,2.5);
\draw[thick, densely dashed] (1,0) -- (1,2.5);
\draw[thick, densely dashed] (2,0.5) -- (2,2.5);
\draw[thick, densely dashed] (3,0.5) -- (3,2.5);
\draw[thick, densely dashed] (4,0.5) -- (4,2.5);
\draw[thick, densely dashed] (5,0) -- (5,2.5);
\draw[thick, densely dashed] (6,0) -- (6,2.5);
\draw[ultra thick, loosely dotted] (4.3,0.5) -- (4.6,0.5);
\end{tikzpicture}
\end{center}
\end{example}

\begin{figure}[h]
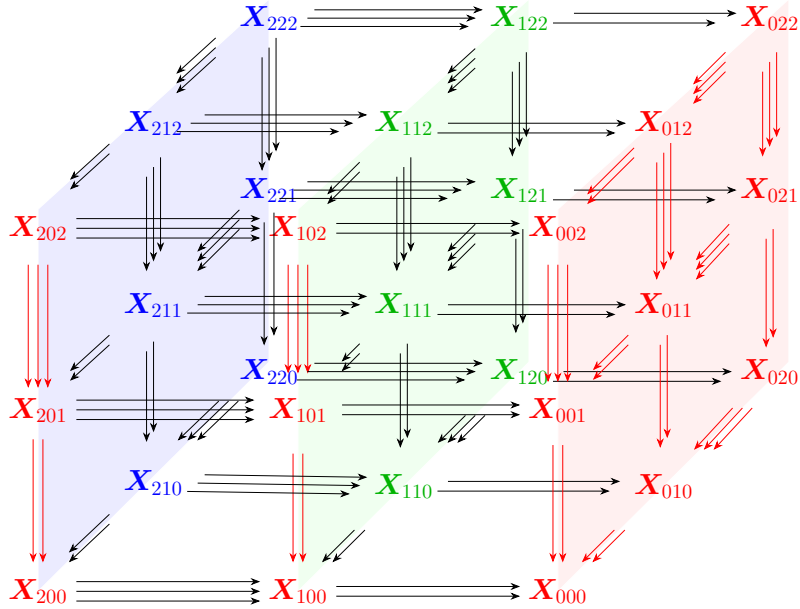

%  \centering
\hspace*{-2cm}
  \includestandalone[width=0.9\textwidth]{Corner3Xred}
  \caption{Corner of the multinerve of a weakly globular $3$-fold category $X$}
  \label{corner3Xred}
\end{figure}

\begin{figure}[h]
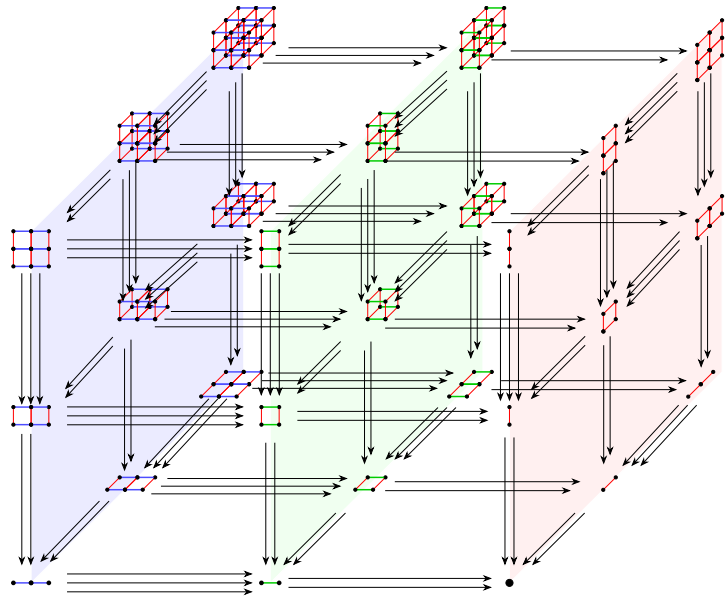

%  \centering
\hspace*{-2cm}
  \includestandalone[width=0.9\textwidth]{Corner3red}
  \caption{Geometric picture of the corner of the multinerve of a weakly globular $3$-fold category}
  \label{corner3red}
\end{figure}

\bk
\begin{example}\label{ex2-n-equival}
\emph{Weakly globular 3-fold categories.}

\index{Weakly globular!- 3-fold category}

A weakly globular 3-fold category $X\in\catwg{3}$ is given by $X\in\funcat{}{\catwg{2}}$ such that
\begin{itemize}
  \item [a)] $X_0 \in \cathd{2}$.\bk

  \item [b)] For each $k\geq 0$, $X_k\cong \pro{X_1}{X_0}{k}$.\bk

  \item [c)] For each $k\geq 0$ the induced Segal maps
  \begin{equation*}
    \hmu{k}: X_k\rw \pro{X_1}{X^d_0}{k}
  \end{equation*}
  are 2-equivalences in $\catwg{2}$.\bk

  \item [d)] There is a functor $\p{3}: \catwg{3}\rw \catwg{2}$ making the following diagram commute
      \begin{equation*}
      \xymatrix{
      \catwg{3} \ar@{^{(}->}^{J_{3}}[rr] \ar_{\p{3}}[d] & & \funcat{2}{\Cat} \ar^{\bar p}[d]\\
      \catwg{2} \ar_{\Nb{2}}[rr] & & \funcat{2}{\Set}
      }
      \end{equation*}
\end{itemize}
It follows from the definition that $X_{k0}\in\cathd{}$ for all $k\geq 0$. On page \pageref{page-corner3red}, Figure \ref{corner3Xred} is a picture of the corner of $\Nb{3}X\in\funcat{3}{\Set}$, where we omitted drawing the degeneracy operators for simplicity. The structures in red are homotopically discrete. A corresponding geometric picture (again with omitted degeneracy operators) in Figure \ref{corner3red} on page \pageref{page-corner3red}.

\end{example}

\section{Properties of weakly globular $\pmb{n}$-fold categories}\label{sbs-prop-wg-nf-cat}
In this section we discuss the main properties of weakly globular \nfol categories. In Proposition \ref{pro-nequiv-to-obj} we show that a weakly globular \nfol category $n$-equivalent to a homotopically discrete one is homotopically discrete. This generalizes to higher dimensions the fact that a category equivalent to an equivalence relation is an equivalence relation. We deduce in Corollary \ref{cor-crit-nequiv-rel} a criterion for a weakly globular \nfol category to be homotopically discrete.

The main result of this section, Proposition \ref{pro-crit-ncat-be-wg} b), gives a criterion for an \nfol category to be weakly globular. This criterion requires certain sub-structures in the \nfol category being homotopically discrete as well as $\bar p J_n X$ (obtained applying levelwise $p$ to $J_n X$ for $X \in \catwg{n}$) to be the multinerve of an object of $\catwg{n-1}$.

This criterion will be used crucially in the proof of Proposition \ref{pro-transf-wg-struc} to characterize \nfol categories levelwise equivalent to Segalic pseudo-functors. This leads to the main result Theorem \ref{the-strict-funct} on the strictification of Segalic pseudo-functors.

The proof of Proposition \ref{pro-crit-ncat-be-wg} b) uses an inductive argument in conjunction with the proof of a property of the category $\catwg{n}$ (Proposition \ref{pro-crit-ncat-be-wg} a)): the fact that the nerve functor in direction 2, when applied to $\catwg{n}$ takes values in $\funcat{}{\catwg{n-1}}$.

\begin{definition}\label{def-pn}\index{Functor!- $\p{j,n}$}
    For each $1\leq j \leq n$ denote
    \begin{align*}
        \p{j,n} & = \p{j}\p{j+1}\cdots \p{n-1}\p{n}:\catwg{n}\rw \catwg{j-1}\\
        \p{n,n}& = \p{n}\;.
    \end{align*}

    Note that this restricts to
    \begin{equation*}
      \p{j,n}= \p{j}\p{j+1}\cdots \p{n-1}\p{n}:\cathd{n}\rw \cathd{j-1}
    \end{equation*}
\end{definition}
\begin{lemma}\label{lem-prop-pn}
    For each $X\in\catwg{n}$, $1\leq j < n$ and $s\geq 2$ it is
    \begin{equation}\label{eq-lem-prop-pn}
    \begin{split}
        &\p{j,n-1}X_s \cong  \p{j,n-1}(\pro{X_1}{X_0}{s})=\\
         & =\pro{\p{j,n-1} X_1}{\p{j,n-1} X_0}{s}\;.
    \end{split}
    \end{equation}
\end{lemma}
\begin{proof}
Since $X\in\catwg{n}$ by definition $\p{n}X\in\catwg{n-1}$, hence
\begin{equation*}
    \p{n-1}(\pro{X_1}{X_0}{s})= \pro{\p{n-1} X_1}{\p{n-1} X_0}{s}
\end{equation*}
which is \eqref{eq-lem-prop-pn} for $j=n-1$. Since $\p{j+1,n}X \in\catwg{j}$ for $1\leq j\leq (n-1)$,  its Segal maps are isomorphisms. Further for all $s\geq 0$
\begin{equation*}
    (\p{j+1,n}X)_s=(\p{j+1}...\p{n}X)_s = \p{j}...\p{n-1}X_s = \p{j,n-1} X_s
\end{equation*}
with $X_s=\pro{X_1}{X_0}{s}$ for $s\geq 2$. This proves \eqref{eq-lem-prop-pn}.
\end{proof}
\begin{remark}\label{rem-eq-def-wg-ncat}
    It follows immediately from Lemma \ref{lem-prop-pn} that if $X\in\catwg{n}$, for all $s\geq 2$
    \begin{equation}\label{eq1-rem-eq-def-wg-ncat}
        X_{s0}^d=(\pro{X_{10}}{X_{00}}{s})^d=\pro{X^d_{10}}{X^d_{00}}{s}\;.
    \end{equation}
    In fact, by \eqref{eq-lem-prop-pn} in the case $j=2$, taking the 0-component, we obtain
    \begin{equation*}
    \begin{split}
        & \p{1,n-2}(\pro{X_{10}}{X_{00}}{s})=\p{1,n-2}X_{s0}=(\p{2,n-1}X_s)_0 \\
        =\ & (\pro{\p{2,n-1} X_1}{\p{2,n-1} X_0}{s})_0= \\
        =\ & \pro{\p{1,n-2}X_{10}}{\p{1,n-2}X_{00}}{s}
    \end{split}
    \end{equation*}
    which is the same as \eqref{eq1-rem-eq-def-wg-ncat}.
\end{remark}
The following proposition is a higher dimensional generalization of the fact that, if a category is equivalent to an equivalence relation, it is itself an equivalence relation.
\begin{proposition}\label{pro-nequiv-to-obj} %5.6
    Let $f:X\rw Y$ be a $\nequ$ in $\catwg{n}$ with $Y\in\cathd{n}$, then $X\in\cathd{n}$.
\end{proposition}
\begin{proof}
By induction on $n$. It is clear for $n=1$. Suppose it is true for $n-1$ and let $f$ be as in the hypothesis. Then $\p{n}f:\p{n} X\rw\p{n}Y$ is a $\equ{n-1}$ with $\p{n}Y\in\cathd{n-1}$ since $Y\in\cathd{n}$. It follows by induction hypothesis that $\p{n}X\in\cathd{n-1}$. We have
\begin{equation}\label{eq1-pro-nequiv}
    X_1=\uset{a,b\in X_0^d}{\cop}X(a,b)\;.
\end{equation}
Since $f$ is a $\nequ$, there are $\equ{n-1}$s
\begin{equation*}
    f(a,b):X(a,b)\rw Y(fa,fb)
\end{equation*}
where $Y(fa,fb)\in\cathd{n-1}$ since $Y\in\cathd{n}$. By induction hypothesis, it follows that $X(a,b)\in\cathd{n-1}$. From \eqref{eq1-pro-nequiv} and the fact that $\cathd{n-1}$ is closed under coproducts (see Lemma \ref{lem-copr-hom-disc}), we conclude that $X_1\in\cathd{n-1}$.

Since $X\in \catwg{n}$, the induced Segal map
\begin{equation*}
    \hmu{s}:X_s=\pro{X_1}{X_0}{s}\rw \pro{X_1}{X^d_0}{s}
\end{equation*}
is a $\equ{n-1}$. Since, from above, $X_1$ is homotopically discrete and $X_0^d$ is discrete, by Lemma \ref{lem-copr-hom-disc},
\begin{equation*}
  \pro{X_1}{X^d_0}{s}\in\cathd{n-1}
\end{equation*}

Thus by induction hypothesis applied to the induced Segal map $\hmu{s}$ we conclude that $X_s\in\cathd{n-1}$ for all $s\geq 0$.

In summary, we showed that $X\in\catwg{n}$ is such that $X_s\in\cathd{n-1}$ for all $s\geq 0$ and $\p{n}X\in\cathd{n-1}$. Therefore, by definition, $X\in\cathd{n}$.
\end{proof}
\begin{corollary}\label{cor-crit-nequiv-rel}
    Let $X\in\catwg{n}$ be such that $X_1$ and $\pn X$ are in $\cathd{n-1}$. Then $X\in\cathd{n}$.
\end{corollary}
\begin{proof}
Since $X\in\catwg{n}$, the induced Segal maps
\begin{equation*}
    \hmu{s}: X_s\rw \pro{X_1}{X^d_0}{s}
\end{equation*}
are $\equ{n-1}$s for all $s\geq 2$. Since by hypothesis $X_1\in\cathd{n-1}$ and $X_0^d$ is discrete, by Lemma \ref{lem-copr-hom-disc},
 \begin{equation*}
    \pro{X_1}{X^d_0}{s}\in\cathd{n-1}\;.
\end{equation*}
By Proposition \ref{pro-nequiv-to-obj} applied to $\hmu{s}$ we conclude that $X_s\in\cathd{n-1}$  for all $s\geq 2$. Therefore $X\in\catwg{n}$ is such that $X_s\in \cathd{n-1}$ and $p^{(n)}X\in\cathd{n-1}$. By definition then $X\in\cathd{n}$.
\end{proof}
\begin{corollary}\label{cor2-crit-nequiv-rel}
    Let $X\in\catwg{n}$, then $X\in\cathd{n}$ if and only if there is an $n$-equivalence $\zg:X\rw Y$ with $Y$ discrete.
\end{corollary}
\begin{proof}
If $X\in\cathd{n}$ then by Corollary \ref{gamma}, $\zg\lo{n}:X\rw X^d$ is an $n$-equivalence. Conversely, suppose that there is an $n$-equivalence $\zg:X\rw Y$ with $Y$ discrete, then in particular $Y\in\cathd{n}$ so, by Proposition \ref{pro-nequiv-to-obj}, $X\in\cathd{n}$.
\end{proof}
\begin{definition}\label{def-kdir-wg-ncat}
Given $X\in\cat{n}$ and $k\geq 0$, let $\Nu{2}X\in\funcat{}{\cat{n-1}}$ as in Definition \ref{def-ner-func-dirk}. Denote for each $k\geq 0$,
 \begin{equation*}
   (\Nu{2}X)_k=X_k\up{2}\in \funcat{}{\cat{n-2}}
 \end{equation*}
 so that
    \begin{equation*}
        (X_k\up{2})_s =
        \left\{
          \begin{array}{ll}
            X_{0k}, & \hbox{$s=0$;} \\
            X_{1k}, & \hbox{$s=1$;} \\
            X_{sk}=\pro{X_{1k}}{X_{0k}}{s}, & \hbox{$s\geq 2$.}
          \end{array}
        \right.
    \end{equation*}
\end{definition}
\bigskip
We denote by $\Nb{n}\catwg{n}$ the image of the multinerve functor $\Nb{n}:\catwg{n}\rw \funcat{n-1}{\Cat}$. Note that, since $\Nb{n}$ is fully faithful, we have an isomorphism $\catwg{n}\cong \Nb{n}\catwg{n}$.

\bigskip
The following lemmas are needed in the initial steps of the induction in the proof of Proposition \ref{pro-crit-ncat-be-wg}.

\begin{lemma}\label{lem-crit-doucat-wg}
    Let $X\in\cat{2}$ be such that
    \begin{itemize}
      \item [i)] $X_0\in\cathd{}$,\mk

      \item [ii)] $\bar p J_2 X\in N\Cat$.\mk
    \end{itemize}
    Then $X\in \catwg{2}$.
\end{lemma}
\begin{proof}
Since $X_0\in \cathd{}$, $pX_0 = X_0^d$. By hypothesis, $p X_2\cong \tens{p X_1}{p X_0}$ and $X_2\cong \tens{X_1}{X_0}$. Using the fact (Lemma \ref{lem-p-pres-fib-pro}) that $p$ commutes with pullbacks over discrete objects, we obtain
\begin{equation*}
    p(\tens{X_1}{X_0})\cong pX_2\cong \tens{p X_1}{p X_0}=\tens{p X_1}{p X^d_0}=p(\tens{X_1}{X^d_0})\;.
\end{equation*}
This shows that the map
\begin{equation*}
    \hmu2:\tens{X_1}{X_0} \rw \tens{X_1}{X^d_0}
\end{equation*}
is essentially surjective on objects. On the other hand, this map is also fully faithful. In fact, given $(a,b),(c,d)\in \tens{X_{10}}{X_{00}}$, we have
\begin{equation*}
\begin{split}
    & (\tens{X_1}{X_0})((a,b),(c,d))\cong X_1(a,c)\tiund{X_0(\pt_0 a,\pt_0 c)} X_1(b,d)\cong \\
    & \cong X_1(a,c)\times X_1(b,d) \cong (\tens{X_1}{X^d_0})(\hmu{2}(a,b),\hmu{2}(c,d))
\end{split}
\end{equation*}
where we used the fact that $X_0(\pt_0 a,\pt_0 c)$ is the one-element set, since $X_0\in\cathd{}$. We conclude that $\hmu{2}$ is an equivalence of categories.

Similarly one shows that for all $k\geq 2$
\begin{equation*}
    \hmuk : \pro{X_1}{X_0}{k} \rw \pro{X_1}{X^d_0}{k}
\end{equation*}
is an equivalence of categories. By definition (see also Example \ref{ex1-n-equival}), this means that $X\in\catwg{2}$.
\end{proof}
\begin{lemma}\label{lem-crit-3cat-wg}\
\begin{itemize}
  \item [a)] The functor $\Nu{2}:\cat{3}\rw \funcat{}{\cat{2}}$ restricts to
  \begin{equation*}
    \Nu{2}:\catwg{3}\rw \funcat{}{\catwg{2}}.
  \end{equation*}

  \item [b)] Let $X\in\cat{3}$ be such that
  \begin{itemize}
    \item [i)] $X_0\in\cathd{2}$, $X_{s0}\in\cathd{}$.

    \item [ii)] $\bar p J_3 X\in \Nb{2}\catwg{2}$.
  \end{itemize}
  \mk

 \nid Then $X\in\catwg{3}$.
\end{itemize}

\end{lemma}
\begin{proof}\

\nid a) We show that, if $X\in\catwg{3}$ and $k\geq 0$, $(\Nu{2}X)_k\in \cat{2}$ satisfies the hypotheses of Lemma \ref{lem-crit-doucat-wg}, so that $(\Nu{2}X)_k\in \catwg{2}$.

Since $X\in\catwg{3}$, $X_0\in\cathd{2}$ and thus for each $k\geq 0$
\begin{equation*}
    (\Nu{2}X)_{k0}=X_{0k}\in\cathd{}\;.
\end{equation*}
So hypothesis i) of Lemma \ref{lem-crit-doucat-wg} holds.
Further,
\begin{equation*}
    \ovl{p}J_2(\Nu{2}X)_k=(\p{3}X)_k\up{2}
\end{equation*}
is the nerve of a category since $\p{3}X\in\catwg{2}$, so hypothesis ii) of Lemma \ref{lem-crit-doucat-wg} also holds. We conclude that $(\Nu{2}X)_k\in\catwg{2}$.

\bk

\nid b) By hypothesis, $X_s\in\cat{2}$ is such that $X_{s0}\in\cathd{}$ and $\bar p J_2 X_s$ is the nerve of a category. Thus by  Lemma \ref{lem-crit-doucat-wg}, $X_s\in\catwg{2}$.

Also by hypothesis $\p{3}X\in\catwg{2}$. To show that $X\in\catwg{3}$ it remains to prove that for each $s\geq 2$ the induced Segal map
\begin{equation*}
    \hmu{s}:\pro{X_1}{X_0}{s}\rw \pro{X_1}{X^d_0}{s}
\end{equation*}
is a 2-equivalence. We show this for $s=2$, the case $s>2$ being similar. We first show that it is a local equivalence. By part a) $(\Nu{2}X)_1\in \catwg{2}$. Thus there is an equivalence of categories
\begin{equation}\label{eq1-lem-crit-3cat}
    \tens{X_{11}}{X_{01}}\rw \tens{X_{11}}{X^d_{01}}=\tens{X_{11}}{(\di{2}\p{2}X_0)_1}\;.
\end{equation}
From hypothesis ii) by Remark \ref{rem-eq-def-wg-ncat} using the fact that $p\p{2}X_{s0}=X_{s0}^d$ we have
\begin{equation*}
   X_{20}^{d}=(\tens{X_{10}}{X_{00}})^d \cong \tens{X^d_{10}}{X^d_{00}}\;.
\end{equation*}
Let $(a,b),(c,d)\in \tens{X^d_{10}}{X^d_{00}}$\;. By \eqref{eq1-lem-crit-3cat} there is an equivalence of categories
\begin{equation}\label{eq2-lem-crit-3cat}
\begin{split}
    & (\tens{X_1}{X_0})((a,b),(c,d)) = \\
    & = X_1(a,c) \tiund{X_0(\pt_0 a, \pt_0 c)} X_1 (b,d) \rw X_1(a,c) \tiund{(\di{2}\p{2}X_0)(\tilde\pt_0 a, \tilde\pt_0 c)} X_1 (b,d)
\end{split}
\end{equation}
On the other hand, since $\p{2}X_0 \in \cathd{}$, $\p{2}X_0(\tilde\pt_0 a, \tilde\pt_0 c)$ is the one-element set. Therefore
\begin{equation}\label{eq2A-lem-crit-3cat}
\begin{split}
    &  X_1(a,c) \tiund{(\di{2}\p{2}X_0)(\tilde\pt_0 a, \tilde\pt_0 c)} X_1 (b,d) \cong \\
    & \cong X_1(a,c)\times X_1(b,d)\cong(\tens{X_1}{X^d_0})((a,b),(c,d))\;.
\end{split}
\end{equation}
From \eqref{eq2-lem-crit-3cat} and \eqref{eq2A-lem-crit-3cat} we conclude that $\hmu{2}$ is a local equivalence. Further, by hypothesis ii), there is an equivalence of categories
\begin{equation*}
\begin{split}
    & \p{2}\hmu{2}:\p{2}(\tens{X_1}{X_0})= \tens{\p{2} X_1}{\p{2} X_0} \xrw{\sim} \\
    & \rw \tens{\p{2} X_1}{(\p{2} X_0)^d}=\p{2}(\tens{X_1}{X^d_0})\;.
\end{split}
\end{equation*}
In conclusion, $\hmu{2}$ is a 2-equivalence, as required.
\end{proof}
\begin{proposition}\label{pro-crit-ncat-be-wg}\

    \begin{itemize}
      \item [a)] The functor $\Nu{2}:\cat{n}\rw\funcat{}{\cat{n-1}}$ restricts to a functor
       \begin{equation*}
         \Nu{2}:\catwg{n}\rw \funcat{}{\catwg{n-1}}.
       \end{equation*}

      \item [b)] Let $X\in\cat{n}$ be such that\mk
  \begin{itemize}
    \item [i)] $X_0\in\cathd{n-1}$, $X_{s0}\in\cathd{n-2}$ for all $s \geq 0$.\mk

    \item [ii)] $\bar p J_n X\in \Nb{n-1}\catwg{n-1}$. \mk
  \end{itemize}
  \mk
   Then $X\in\catwg{n}$.
    \end{itemize}
\end{proposition}
\begin{proof}
By induction on $n$. For $n=2,3$ see Lemmas \ref{lem-crit-doucat-wg} and \ref{lem-crit-3cat-wg}. Suppose, inductively, that it holds for $(n-1)$. We shall denote $(\Nu{2}X)_k=X_k\up{2}$.
\bk

a) Clearly $X_k\up{2}\in\cat{n-1}$; we show that $X_k\up{2}$ satisfies the inductive hypothesis b) and thus conclude that $X_k\up{2}\in\catwg{n-1}$.

We have $(X_k\up{2})_0=X_{0k}\in\cathd{n-2}$ since $X_0\in\cathd{n-1}$ (as $X\in\catwg{n}$). Further,
\begin{equation*}
    (X_k\up{2})_{s0}=X_{sk0}\in \cathd{n-3}
\end{equation*}
since $X_{sk}\in\catwg{n-2}$ (as $X_s\in \catwg{n-1}$ because $X\in\catwg{n}$). Thus condition i) in the inductive hypothesis b)  holds for $X_{k}\up{2}$. To show that condition ii) holds, note that
\begin{equation}\label{eq1-lem-xk-catwg}
    \bar p J_{n-1}X_k\up{2}=\Nb{n-1}(\p{n}X)_k\up{2}
\end{equation}
In fact, for all $(r_1,...,r_{n-2})\in {{\zD}^{n-2}}^{op}$,
\begin{align*}
    &(\bar p J_{n-1}X_k\up{2})_{r_1...r_{n-2}}=p(X_k\up{2})_{r_1...r_{n-2}}=\\
    =\ & pX_{r_1,k,r_2...r_{n-2}}=(\bar p J_{n-2}X_{r_1,k})_{r_2...r_{n-2}}=\\
    =\ & (\Nb{n-2}\p{n-2}X_{r_1,k})_{r_2...r_{n-2}}=\\
    =\ & (\Nb{n-2}((\p{n}X)\up{2}_k)_{r_1})_{r_2...r_{n-2}}=\\
    =\ & (\Nb{n-1}(\p{n}X)\up{2}_k)_{r_1...r_{n-2}}\;.
\end{align*}
Since this holds for all $r_1,...,r_{n-2}$, \eqref{eq1-lem-xk-catwg} follows.

By induction hypothesis a) applied to $\p{n}X$, $(\p{n}X)_k\up{2}\in\catwg{n-2}$. Therefore \eqref{eq1-lem-xk-catwg} means that $X_k\up{2}\in\cat{n-1}$ satisfies condition ii) in the inductive hypothesis b). Thus we conclude that $X_k\up{2}\in\catwg{n-1}$ proving a).
\bk

b) Suppose, inductively, that the statement holds for $n-1$ and let $X$ be as in the hypothesis. For each $s\geq 0$ consider $X_s\in\cat{n-1}$. By hypothesis, $X_{s0}\in\cathd{n-2}$ and
\begin{equation*}
    \bar p J_{n-1}X_s =(\bar p J_n X)_s \in \Nb{n-2}\catwg{n-2}
\end{equation*}
since $\bar p J_n X\in \Nb{n-1}\catwg{n-1}$. Thus $X_s$ satisfies the induction hypothesis and we conclude that $X_s\in\catwg{n-1}$. Further, for each $(k_1,\ldots,k_{n-2}) \in \dop{n-2}$ we have
\begin{equation*}
\begin{split}
    & (\Nb{n-1}\ovl{\p{n-1}}X)_{k_1...k_{n-2}} = (\Nb{n-2}\p{n-1}X_{k_1})_{k_2...k_{n-2}}=\\
    =\;& (\bar p J_{n-1}X_{k_1})_{k_2...k_{n-2}}=p X_{k_1...k_{n-2}}=(\bar p J_n X)_{k_1...k_{n-2}}\;.
\end{split}
\end{equation*}
Since, by hypothesis, $\bar p J_n X\in \Nb{n-1}\catwg{n-1}$, we conclude that $\ovl{\p{n-1}}X\in\catwg{n-1}$. We can therefore define \begin{equation*}
\p{n}X = \ovl{\p{n-1}}X\in\catwg{n-1}.
\end{equation*}

To prove that $X\in\catwg{n}$ it remains to prove that the induced Segal maps
\begin{equation*}
    \hmu{s}:\pro{X_1}{X_0}{s}\rw \pro{X_1}{X_0^d}{s}
\end{equation*}
are $\equ{n-1}$s for all $s\geq 2$. We prove this for $s=2$ the case $s>2$ being similar. We claim that $X_k\up{2}\in\cat{n-1}$ satisfies the inductive hypothesis b). In fact, $(X_k\up{2})_0 = X_{0k}\in \cathd{n-2}$ since $X_0\in\cathd{n-1}$; for each $s\geq 0$, $(X_k\up{2})_{s0}=X_{sk0}\in\cathd{n-3}$ since, from above, $X_s\in\catwg{n-1}$.

Also, from a) and the fact that, by hypothesis, $\bar p J_n X\in\Nb{n-1}\catwg{n-1}$, we conclude that
\begin{equation*}
    \bar p J_{n-1} X_k\up{2} \in \Nb{n-2}\catwg{n-2} \;.
\end{equation*}
Thus $X_k\up{2}$ satisfies the inductive hypothesis b) and we conclude that $X_k\up{2}\in\catwg{n-1}$. It follows that the induced Segal map
\begin{equation}\label{eq3-crit-ncat-be-wg}
    \tens{X_{1k}}{X_{0k}}\rw \tens{X_{1k}}{X^d_{0k}}=\tens{X_{1k}}{(\p{2,n-1}X_0)_k}
\end{equation}
is a $\equ{n-2}$. Since $\p{n}X\in\catwg{n-1}$, using Remark \ref{rem-eq-def-wg-ncat} and the fact that $(\p{n}X)_{s0}^d = (\p{n-2}\p{n-1}X_{s0})^d = X_{s0}^d$ we obtain
\begin{equation*}
    (\tens{X_{10}}{X_{00}})^d=\tens{X^d_{10}}{X^d_{00}}\;.
\end{equation*}
Let $(a,b),(c,d)\in(\tens{X_{10}}{X_{00}})^d=\tens{X^d_{10}}{X^d_{00}}$. By \eqref{eq3-crit-ncat-be-wg} there is a $(n-2)$-equivalence
\begin{equation}\label{eq3A-crit-ncat-be-wg}
\begin{split}
    & (\tens{X_1}{X_0})((a,b),(c,d))=X_1(a,c)\tiund{X_0(\pt_0 a,\pt_0 c)}X_1(b,d)\rw \\
    & \rw X_1(a,c)\tiund{(\p{2,n-1}X_0)(\tilde\pt_0 a,\tilde\pt_0 c)}X_1(b,d)\;.
\end{split}
\end{equation}
On the other hand, $\p{2,n-1}X_0\in\cathd{}$ is an equivalence relation, therefore
\begin{equation*}
    \p{2,n-1}X_0(\tilde\pt_0 a,\tilde\pt_0 c)
\end{equation*}
is the one-element set. It follows that
\bk
\begin{equation}\label{eq3B-crit-ncat-be-wg}
\begin{split}
    & X_1(a,c)\tiund{(\p{2,n-1}X_0)(\tilde\pt_0 a,\tilde\pt_0 c)}X_1(b,d)\cong\\
    & \cong X_1(a,c)\times X_1(b,d)\cong (\tens{X_1}{X_0^d})((a,b),(c,d))\;.
\end{split}
\end{equation}
Thus \eqref{eq3A-crit-ncat-be-wg}  and \eqref{eq3B-crit-ncat-be-wg} imply that $\hmu{2}$ is a local $\equ{n-2}$.

To show that $\hmu{2}$ is a $\equ{n-1}$ it remains to prove that $\p{n-1}\hmu{2}$ is a $\equ{n-2}$. Since from above, $\p{n}X=\ovl{\p{n-1}}X\in\catwg{n-1}$, we have
\begin{equation*}
\begin{split}
    & \p{n-1}\hmu{2}:\p{n-1}(\tens{X_1}{X_0})\cong  \tens{\p{n-1}X_1}{\p{n-1}X_0}\rw\\
    & \rw \tens{\p{n-1}X_1}{(\p{n-1}X_0)^d}=\p{n-1}(\tens{X_1}{X^d_0})\;.
\end{split}
\end{equation*}
is a $\equ{n-2}$, as required.
\end{proof}

\begin{corollary}\label{rem-last-ins}
Let $f:X\rw Y$ be a morphism in $\cat{n}$ such that $(J_n f)_{\uk}$ is an equivalence of categories for all $\uk\in\dop{n-1}$. Then
\begin{itemize}
  \item [a)] If $X\in\catwg{n}$, then $Y\in\catwg{n}$.\mk

  \item [b)] If $X\in\cathd{n}$, then $Y\in\cathd{n}$.\mk
\end{itemize}
\end{corollary}
\begin{proof}
In this proof we shall denote $(J_n X)_{\uk}=X_{\uk}$, and similarly for $Y$. By proceed by induction on $n$.

 When $n=2$ let $X\in\catwg{2}$. Since $X_0\simeq Y_0$ and $X_0\in\cathd{}$, also $Y_0\in\cathd{}$. Also $\bp J_2 Y\cong \p{2}X\in N\Cat$. Thus by Lemma \ref{lem-crit-doucat-wg}, $Y\in\catwg{2}$. If $X\in\cathd{2}$, in particular $X\in\catwg{2}$, so $Y\in\catwg{2}$; but we also have $X_1\simeq Y_1$ so, since $X_1\in\cathd{}$, $Y_1\in\cathd{}$. Further, $\p{2}Y\cong\p{2}X\in\cathd{}$ (as $X\in\cathd{2}$).
Thus the hypotheses of Corollary \ref{cor-crit-nequiv-rel} are satisfied for $Y$ and we conclude that $Y\in\cathd{2}$.

Suppose, inductively, that the corollary holds for $1\leq k\leq n-1$.

a) We verify that $Y\in\cat{n}$ satisfies the hypotheses of Proposition \ref{pro-crit-ncat-be-wg} b). Since $X\in\catwg{n}$, $X_0\in\cathd{n-1}$. So by inductive hypothesis b) applied to the map $f_0: X_0\rw Y_0$ we conclude that $Y_0\in\cathd{n-1}$. Similarly, $X_{s0}\in\cathd{n-2}$, so the inductive hypothesis b) applied to the map $f_{s0}: X_{s0}\rw Y_{s0}$ affords $Y_{s0}\in\cathd{n-2}$. So hypotheses b) ii) of Proposition \ref{pro-crit-ncat-be-wg} are satisfied for $Y$.

 As for hypothesis b) i), this also holds because, since $(J_n f)_{\uk}$ is an equivalence of categories for each $\uk$, $\bp (J_n f)_{\uk}$ is an isomorphism, so that
\begin{equation*}
  \bp J_n Y \cong \bp J_n X = \p{n}X \in \Nb{n-1}\catwg{n-1}\;.
\end{equation*}
We conclude by Proposition \ref{pro-crit-ncat-be-wg} b) that $Y\in\catwg{n-1}$.

\mk

b) if $X\in\cathd{n}$, in particular $X\in\catwg{n}$, so by a) $Y\in\catwg{n}$. In addition, $X_1\in\cathd{n-1}$ (since $X\in\cathd{n}$), so by inductive hypothesis b) applied to $f_1:X_1\rw Y_1$, we conclude that $Y_1\in\cathd{n-1}$. Since $(J_n f)_{\uk}$ is an equivalence of categories for all $\uk\in\dop{n-1}$ we also have
\begin{equation*}
  \p{n}Y\cong \p{n}X\in\cathd{n-1}\;.
\end{equation*}
Thus $Y$ satisfies all the hypotheses of Corollary \ref{cor-crit-nequiv-rel} and we conclude that $Y\in\cathd{n}$.
\end{proof}

%%
%%%%%%%%%%%%%%%%%%%%%%%%%%%%%%%%%%%%%%%%%%%%%%%%%%%%%%%%%%%%%%%%%%%%%%%%%%%%%%%%%%
\clearpage

%%%%%%%%%%%%%%%%%%%%%%%%%%%%%%%%%%%%%%%%%%%%%%%%%%%%%%%%%%%%%%%
\bk

\bk

%A geometric picture is as follows, where the red structures are homotopically discrete

%% Enter a TIKX picture
%%%%%%%%%%%%%%%%%%%%%%%%%%%%%%%%%%%%%%%%%%%%%%%%%%%%%%%%%%%%%%%
\begin{figure}[ht]
%  \centering
  %% Begin a TIKX picture
%%%%%%%%%%%%%%%%%%%%%%%%%%%%%%%%%%%%%%%%%%%%%%%%%%%%%%%%%%%%%%%
%% Create objects \tickzset
%% Create single edge
\tikzset{Edge/.pic={
\draw plot [mark=*] coordinates {(0,0) (0.5,0)};
}}

%% Create a face
\tikzset{Face/.pic={
\draw plot [mark=*] coordinates {(0,0) (0.5,0)  (0.5,0.5)  (0,0.5) (0,0)};
}}
\tikzset{Facebr/.pic={
\draw[black] plot [mark=*] coordinates {(0,0) (0.5,0)};
\draw[red] plot [mark=*] coordinates {(0.5,0) (0.5,0.5)};
\draw[red] plot [mark=*] coordinates {(0,0.5) (0,0)};
\draw[black] plot [mark=*] coordinates {(0.5,0.5) (0,0.5)};
%\draw[black] plot $Rightarrow$ coordinates {(0.25,0.25)};
\node[rotate=-90,scale=0.7] at (0.25,0.25) () {$\Rightarrow$};
}}

%% Create 3 dots
\tikzset{Dots/.pic={
\filldraw [black]
(0,0) circle [radius=0.5pt]
(0.3,0) circle [radius=0.5pt]
(0.6,0) circle [radius=0.5pt];
}}

\begin{tikzpicture}[thick,mark options={color=black}]
% Objects in row 1 (bottom)
\draw[fill=red] (0,0)  circle(0.05)[black]; %% (nodo 18)
\pic [black] at (-2,0) {Edge};
\pic [black] at (-4.0,0) {Edge};\pic [black] at (-4.5,0) {Edge};
% Object in row 2 (center)
\pic [red,rotate=-90] at (0,2) [xshift=0cm] {Edge};
\pic [black] at (-2,1.5) {Facebr};
\pic [black] at (-4.0,1.5) {Facebr};\pic [black] at (-4.5,1.5) {Facebr};

% Objects in row 3 (upper)
\pic [red,rotate=-90] at (0,4.5) [xshift=0cm] {Edge};
\pic [red,rotate=-90] at (0,4.0) [xshift=0cm] {Edge};
\pic [black] at (-2,3.5) {Facebr};\pic [black] at (-2,4.0) {Facebr};
\pic [black] at (-4.0,3.5) {Facebr};\pic [black] at (-4.0,4.0) {Facebr};
\pic [black] at (-4.5,3.5) {Facebr};\pic [black] at (-4.5,4.0) {Facebr};
\pic [black] at (-6,0) {Dots};
\pic [black] at (-6,1.75) {Dots};
\pic [black] at (-6,4) {Dots};

%% Positioning nodes
\node (1)[minimum size=10mm] at (0,0) {};
\node (2)[minimum size=10mm] at (-1.75,0) {};
\node (3)[minimum size=10mm] at (-4,0) {};
\node (4)[minimum size=10mm] at (0,1.75) {};
\node (5)[minimum size=10mm] at (-1.75,1.75) {};
\node (6)[minimum size=10mm] at (-4,1.75) {};
\node (7)[minimum size=10mm] at (0,4) {};
\node (8)[minimum size=10mm] at (-1.75,4) {};
\node (9)[minimum size=10mm] at (-4,4) {};

%% Arrows
%% 1-2
\draw[transform canvas={yshift=1.0ex,xshift=0ex},->][arrows = {-Stealth}][thin] (2) -- (1);
\draw[transform canvas={yshift=0.0ex,xshift=0ex},->][arrows = {-Stealth}][thin] (2) -- (1);
\draw[transform canvas={yshift=-1.0ex,xshift=0ex},->][arrows = {-Stealth}][thin] (1) -- (2);
%% 4-5
\draw[transform canvas={yshift=1.0ex,xshift=0ex},->][arrows = {-Stealth}][thin] (5) -- (4);
\draw[transform canvas={yshift=0.0ex,xshift=0ex},->][arrows = {-Stealth}][thin] (5) -- (4);
\draw[transform canvas={yshift=-1.0ex,xshift=0ex},->][arrows = {-Stealth}][thin] (4) -- (5);
%% 7-8
\draw[transform canvas={yshift=1.0ex,xshift=0ex},->][arrows = {-Stealth}][thin] (8) -- (7);
\draw[transform canvas={yshift=0.0ex,xshift=0ex},->][arrows = {-Stealth}][thin] (8) -- (7);
\draw[transform canvas={yshift=-1.0ex,xshift=0ex},->][arrows = {-Stealth}][thin] (7) -- (8);
%% 2-3
\draw[transform canvas={yshift=2.0ex,xshift=1ex},->][arrows = {-Stealth},shorten >=0.2cm, shorten <=0.2cm][thin] (3) -- (2);
\draw[transform canvas={yshift=1.0ex,xshift=1ex},->][arrows = {-Stealth},shorten >=0.2cm, shorten <=0.2cm][thin] (3) -- (2);
\draw[transform canvas={yshift=0.0ex,xshift=1ex},->][arrows = {-Stealth},shorten >=0.2cm, shorten <=0.2cm][thin] (3) -- (2);
\draw[transform canvas={yshift=-1.0ex,xshift=1ex},->][arrows = {-Stealth},shorten >=0.2cm, shorten <=0.2cm][thin] (2) -- (3);
\draw[transform canvas={yshift=-2.0ex,xshift=1ex},->][arrows = {-Stealth},shorten >=0.2cm, shorten <=0.2cm][thin] (2) -- (3);
%% 5-6
\draw[transform canvas={yshift=2.0ex,xshift=1ex},->][arrows = {-Stealth},shorten >=0.2cm, shorten <=0.2cm][thin] (6) -- (5);
\draw[transform canvas={yshift=1.0ex,xshift=1ex},->][arrows = {-Stealth},shorten >=0.2cm, shorten <=0.2cm][thin] (6) -- (5);
\draw[transform canvas={yshift=0.0ex,xshift=1ex},->][arrows = {-Stealth},shorten >=0.2cm, shorten <=0.2cm][thin] (6) -- (5);
\draw[transform canvas={yshift=-1.0ex,xshift=1ex},->][arrows = {-Stealth},shorten >=0.2cm, shorten <=0.2cm][thin] (5) -- (6);
\draw[transform canvas={yshift=-2.0ex,xshift=1ex},->][arrows = {-Stealth},shorten >=0.2cm, shorten <=0.2cm][thin] (5) -- (6);
%% 8-9
\draw[transform canvas={yshift=2.0ex,xshift=1ex},->][arrows = {-Stealth},shorten >=0.2cm, shorten <=0.2cm][thin] (9) -- (8);
\draw[transform canvas={yshift=1.0ex,xshift=1ex},->][arrows = {-Stealth},shorten >=0.2cm, shorten <=0.2cm][thin] (9) -- (8);
\draw[transform canvas={yshift=0.0ex,xshift=1ex},->][arrows = {-Stealth},shorten >=0.2cm, shorten <=0.2cm][thin] (9) -- (8);
\draw[transform canvas={yshift=-1.0ex,xshift=1ex},->][arrows = {-Stealth},shorten >=0.2cm, shorten <=0.2cm][thin] (8) -- (9);
\draw[transform canvas={yshift=-2.0ex,xshift=1ex},->][arrows = {-Stealth},shorten >=0.2cm, shorten <=0.2cm][thin] (8) -- (9);
%% 7-4
\draw[red, transform canvas={xshift=-2.0ex,yshift=0ex},->][arrows = {-Stealth},shorten >=0.2cm, shorten <=0.2cm][thin] (7) -- (4);
\draw[red, transform canvas={xshift=-1.0ex,yshift=0ex},->][arrows = {-Stealth},shorten >=0.2cm, shorten <=0.2cm][thin] (7) -- (4);
\draw[red, transform canvas={xshift=0.0ex,yshift=0ex},->][arrows = {-Stealth},shorten >=0.2cm, shorten <=0.2cm][thin] (7) -- (4);
\draw[red, transform canvas={xshift=1.0ex,yshift=0ex},->][arrows = {-Stealth},shorten >=0.2cm, shorten <=0.2cm][thin] (4) -- (7);
\draw[red, transform canvas={xshift=2.0ex,yshift=0ex},->][arrows = {-Stealth},shorten >=0.2cm, shorten <=0.2cm][thin] (4) -- (7);
%% 8-5
\draw[transform canvas={xshift=-2.0ex,yshift=0ex},->][arrows = {-Stealth},shorten >=0.2cm, shorten <=0.2cm][thin] (8) -- (5);
\draw[transform canvas={xshift=-1.0ex,yshift=0ex},->][arrows = {-Stealth},shorten >=0.2cm, shorten <=0.2cm][thin] (8) -- (5);
\draw[transform canvas={xshift=0.0ex,yshift=0ex},->][arrows = {-Stealth},shorten >=0.2cm, shorten <=0.2cm][thin] (8) -- (5);
\draw[transform canvas={xshift=1.0ex,yshift=0ex},->][arrows = {-Stealth},shorten >=0.2cm, shorten <=0.2cm][thin] (5) -- (8);
\draw[transform canvas={xshift=2.0ex,yshift=0ex},->][arrows = {-Stealth},shorten >=0.2cm, shorten <=0.2cm][thin] (5) -- (8);
%% 9-6
\draw[transform canvas={xshift=-2.0ex,yshift=0ex},->][arrows = {-Stealth},shorten >=0.2cm, shorten <=0.2cm][thin] (9) -- (6);
\draw[transform canvas={xshift=-1.0ex,yshift=0ex},->][arrows = {-Stealth},shorten >=0.2cm, shorten <=0.2cm][thin] (9) -- (6);
\draw[transform canvas={xshift=0.0ex,yshift=0ex},->][arrows = {-Stealth},shorten >=0.2cm, shorten <=0.2cm][thin] (9) -- (6);
\draw[transform canvas={xshift=1.0ex,yshift=0ex},->][arrows = {-Stealth},shorten >=0.2cm, shorten <=0.2cm][thin] (6) -- (9);
\draw[transform canvas={xshift=2.0ex,yshift=0ex},->][arrows = {-Stealth},shorten >=0.2cm, shorten <=0.2cm][thin] (6) -- (9);
%% 4-1
\draw[red, transform canvas={xshift=-1.0ex,yshift=0ex},->][arrows = {-Stealth},shorten >=0.cm, shorten <=0.cm][thin] (4) -- (1);
\draw[red, transform canvas={xshift=0.0ex,yshift=0ex},->][arrows = {-Stealth},shorten >=0.cm, shorten <=0.cm][thin] (4) -- (1);
\draw[red, transform canvas={xshift=1.0ex,yshift=0ex},->][arrows = {-Stealth},shorten >=0.cm, shorten <=0.cm][thin] (1) -- (4);
%% 5-2
\draw[transform canvas={xshift=-1.0ex,yshift=0ex},->][arrows = {-Stealth},shorten >=0.cm, shorten <=0.cm][thin] (5) -- (2);
\draw[transform canvas={xshift=0.0ex,yshift=0ex},->][arrows = {-Stealth},shorten >=0.cm, shorten <=0.cm][thin] (5) -- (2);
\draw[transform canvas={xshift=1.0ex,yshift=0ex},->][arrows = {-Stealth},shorten >=0.cm, shorten <=0.cm][thin] (2) -- (5);
%% 6-3
\draw[transform canvas={xshift=-1.0ex,yshift=0ex},->][arrows = {-Stealth},shorten >=0.cm, shorten <=0.cm][thin] (6) -- (3);
\draw[transform canvas={xshift=0.0ex,yshift=0ex},->][arrows = {-Stealth},shorten >=0.cm, shorten <=0.cm][thin] (6) -- (3);
\draw[transform canvas={xshift=1.0ex,yshift=0ex},->][arrows = {-Stealth},shorten >=0.cm, shorten <=0.cm][thin] (3) -- (6);
\end{tikzpicture}
  \caption{Geometric picture of the corner of the double nerve of a weakly globular double category}
  \label{corner2red}
\end{figure}
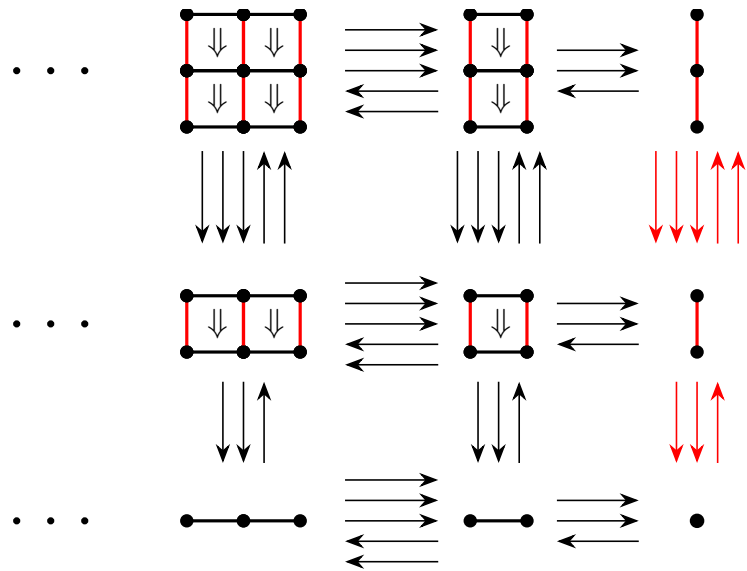

\label{page-corner2red}

%%
%%%%%%%%%%%%%%%%%%%%%%%%%%%%%%%%%%%%%%%%%%%%%%%%%%%%%%%%%%%%%%%%%%%%%%%%%%%%%%%%%%%

\clearpage

In the following picture, for all $i,j,k \in\Delta^{op}$

$X_{2jk}\cong \tens{X_{1jk}}{X_{0jk}}, \;X_{i2k}\cong \tens{X_{i1k}}{X_{i0k}},\;X_{ij2}\cong \tens{X_{ij1}}{X_{ij0}}$\:.

%The corresponding geometric picture is as follows, where the structures in red are homotopically discrete.

\bk

\label{page-corner3red}

\clearpage

%%%%%%%%%%%%%%%%%%%%%%%%%%%%%%%%%%%%%%%%%%%%%%%%%%%%%%%%%%%%%%%%%%%%%%%%%%%%%%%%%%

\chapter{Pseudo-functors modelling higher structures}\label{chap5}
In this chapter we connect the category $\catwg{n}$ of weakly globular $n$-fold categories introduced in Chapter \ref{chap4} to the notion of pseudo-functor. Our main result in this chapter is that weakly globular $n$-fold categories arise as strictification of a special class of pseudo-functors, which we call Segalic pseudo-functors. This result is used in Chapters \ref{chap7} and \ref{chap9} to build the comparison functors between $\catwg{n}$ and Tamsamani $n$-categories and show that they are suitably equivalent.

The classical theory of strictification of pseudo-algebras \cite{PW}, \cite{Lack} affords the strictification functor

 \begin{equation}\label{eq1-intro}
    \St:\psc{n-1}{\Cat} \rw \funcat{n-1}{\Cat}
\end{equation}
left adjoint to the inclusion.

The coherence axioms in a pseudo-functor are reminiscent of the coherence data for the compositions of higher cells in a weak higher category. So it is natural to ask if a subcategory of pseudo-functors can model, in a suitable sense, higher structures. In this chapter we positively answer this question by introducing a subcategory

\begin{equation*}
  \segpsc{n-1}{\Cat}\subset \psc{n-1}{\Cat}
\end{equation*}
of \emph{Segalic pseudo-functors}. Our main result, Theorem \ref{the-strict-funct} is that the strictification functor \ref{eq1-intro} restricts to a functor
\begin{equation*}
    \St:\segpsc{n-1}{\Cat} \rw\catwg{n}\;.
\end{equation*}
(where we identity \catwg{n} with the image $J_{n}\catwg{n}$ of the multinerve functor $J_n:\catwg{n}\rw\funcat{n-1}{\Cat}$\;).

In Chapter \ref{chap7} we will associate to a weakly globular Tamsamani $n$-category a Segalic pseudo-functor and build the rigidification functor $Q_n$ from weakly globular Tamsamani $n$-categories to weakly globular \nfol categories as a composite
\begin{equation*}
    Q_n:\tawg{n}\xrw{\ \;\;} \segpsc{n-1}{\Cat}\xrw{\St} \catwg{n}\;.
\end{equation*}

This chapter is organized as follows. In Section \ref{sec-pseud-strict} we define Segalic pseudo-functors, in Section \ref{sec-stric-pse-fun} we discuss their strictification. We show in Proposition \ref{pro-transf-wg-struc} that an \nfol category whose multinerve is levelwise equivalent to a Segalic pseudo-functor satisfies the hypothesis of Proposition \ref{pro-crit-ncat-be-wg} b) and is therefore weakly globular.

 In Theorem \ref{the-strict-funct}, using  the properties of the monad corresponding to Segalic pseudo-functors proved in Lemma \ref{lem-strict-mon-wg-seg-ps-fun}, we show that the strictification of a Segalic pseudo-functor is an \nfol category and that it satisfies the hypotheses of Proposition \ref{pro-transf-wg-struc}. We therefore conclude that the strictification functor restricts to the functor $\St$ from Segalic pseudo-functors to weakly globular \nfol categories.

\section{The definition of Segalic pseudo-functor}\label{sec-pseud-strict}

In this section we give the definition of the category of Segalic pseudo-functors as a full subcategory of the category
 \begin{equation*}
   \psc{n-1}{\Cat}
 \end{equation*}
 of pseudo-functors and pseudo-natural transformations \cite{Borc}.

\subsection{The idea of Segalic pseudo-functor}\label{subs-idea-segps}
A topological intuition about an object of $\psc{n-1}{\Cat}$ is that it consists of categories $X_{\uk}$ for each object $\uk$ of $\dop{n-1}$ together with multi-simplicial face and degeneracy maps satisfying the multi-simplicial identities not as equalities but as isomorphisms, and these isomorphisms satisfy coherence axioms. Guided by this intuition, we generalize to certain pseudo-functors the multi-simplicial notion of Segal map.

For this purpose, consider a  functor $H\in\funcat{n-1}{\Cat}$. For each $\uk=(k_1,...,k_{n-1})\in\dop{n-1}$ and $1\leq i\leq n-1$ we have
\begin{equation*}
    H(k_1,...,k_{i-1},\mi,k_{i+1},...,k_{n-1})\in\funcat{}{\Cat}
\end{equation*}
Denoting $\uk(1,i)$ and $\uk(0,i)$ as in Notation \ref{not-simp} there is a corresponding Segal map for each $k_i\geq 2$
\begin{equation*}
    H_{\uk}\rw \pro{H_{\uk(1,i)}}{H_{\uk(0,i)}}{k_i}
\end{equation*}
identified by the commuting diagram
\begin{equation}\label{eq2-sbs-seg-map}
    \xy
    0;/r.60pc/:
    (0,0)*+{H_{\uk}}="1";
    (-7,-5)*+{H_{\uk(1,i)}}="2";
    (-2,-5)*+{H_{\uk(1,i)}}="3";
    (9,-5)*+{H_{\uk(1,i)}}="4";
    (-10,-9)*+{H_{\uk(0,i)}}="5";
    (-5,-9)*+{H_{\uk(0,i)}}="6";
    (0,-9)*+{H_{\uk(0,i)}}="7";
    (6,-9)*+{H_{\uk(0,i)}}="8";
    (12,-9)*+{H_{\uk(0,i)}}="9";
    (3,-5)*+{\cdots}="10";
    (3,-9)*+{\cdots}="11";
    {\ar_{\nu_1}"1";"2"};
    {\ar^{\nu_2}"1";"3"};
    {\ar^{\nu_k}"1";"4"};
    {\ar_{d_1}"2";"5"};
    {\ar^{d_0}"2";"6"};
    {\ar^{d_1}"3";"6"};
    {\ar^{d_0}"3";"7"};
    {\ar_{d_1}"4";"8"};
    {\ar^{d_0}"4";"9"};
    \endxy
\end{equation}
If $H$ is not a  functor but a pseudo-functor $H\in\psc{n-1}{\Cat}$, diagram \eqref{eq2-sbs-seg-map} no longer commutes but pseudo-commutes and thus we can no longer define Segal maps.
However, if $H_{\uk(0,i)}$ is a discrete category, then diagram  \eqref{eq2-sbs-seg-map} commutes and therefore we can define Segal maps for $H$.

In the definition of Segalic pseudo-functor we require the above discreetness conditions to be satisfied to be able to define Segal maps and then we require all Segal maps to be isomorphisms.

The last condition in the definition of Segalic pseudo-functor is about the existence of a truncation functor. Applying the isomorphism classes of objects functor $p:\Cat\rw\Set$ levelwise to a pseudo-functor in $\psc{n}{\Cat}$ produces a strict functor; that is, there is a functor
\begin{equation*}
    \ovl{p}:\psc{n}{\Cat}\rw\funcat{n-1}{\Set}
\end{equation*}
such that, for all $\uk\in\dop{n-1}$
\begin{equation*}
    (\ovl{p}X)_{\uk}=pX_{\uk}\;.
\end{equation*}
For $X$ to be a Segalic pseudo-functor we require $\ovl{p}X$ to be (the multinerve of) a weakly globular \nfol category (more precisely, in the image of $\Nb{n}:\catwg{n}\rw\funcat{n-1}{\Set}$). That is we require to have a functor
\begin{equation*}
    \p{n+1}:\segpsc{n}{\Cat}\rw\catwg{n}\;.
\end{equation*}

\subsection{The formal definition of Segalic pseudo-functor}
\begin{definition}\label{def-seg-ps-fun}
\index{Segalic pseudo-functors}

    We define the subcategory $\segpsc{n}{\Cat}$ of $\psc{n}{\Cat}$ as follows:

    For $n=1$, $H\in \segpsc{}{\Cat}$ if $H_0$ is discrete and the Segal maps are isomorphisms: that is, for all $k\geq 2$
    \begin{equation*}
        H_k\cong\pro{H_1}{H_0}{k}
    \end{equation*}
    Note that, since $p$ commutes with pullbacks over discrete objects, there is a functor
   \begin{equation*}
   \begin{split}
       & \p{2}:\segpsc{}{\Cat} \rw \Cat\;, \\
       & (\p{2}X)_{k}=p X_k\;.
   \end{split}
   \end{equation*}
   That is the following diagram commutes:
   \begin{equation*}
    \xymatrix{
    \segpsc{}{\Cat} \ar@{^{(}->}[rr]\ar_{\p{2}}[d] && \psc{}{\Cat} \ar^{\ovl{p}}[d]\\
    \Cat \ar[rr] && \funcat{}{\Set}
    }
   \end{equation*}
   When $n>1$, $\segpsc{n}{\Cat}$ is the full subcategory of $\psc{n}{\Cat}$ whose objects $H$ satisfy the following: \mk
   \begin{itemize}
     \item [a)] \emph{Discreteness condition}:  $H_{\uk(0,i)}$ is discrete for all $\uk\in\Dmenop$ and $1 \leq i \leq n$. \mk
     \item [b)] \emph{Segal maps condition}: As explained above, by a) there are Segal maps for each $k_i\geq 2$
\index{Segal maps condition} \index{Discreteness condition}
\begin{equation*}
    H_{\uk}\rw \pro{H_{\uk(1,i)}}{H_{\uk(0,i)}}{k_i}
\end{equation*}
 We require that all these Segal maps are isomorphisms
      \begin{equation*}
      H_{\uk}\cong\pro{H_{\uk(1,i)}}{H_{\uk(0,i)}}{k_i}
      \end{equation*}
      for all $\uk\in\Dmenop$, $1 \leq i \leq n$ and $k_i\geq2$.\mk
     \item [c)] \emph{Truncation functor}: There is a functor \index{Truncation functor}
     \begin{equation*}
        \p{n+1}:\segpsc{n}{\Cat}\rw \catwg{n}
     \end{equation*}
   \end{itemize}
    making the following diagram commute:
    \begin{equation*}
        \xymatrix{
        \segpsc{n}{\Cat} \ar@{^(->}[r] \ar_{\p{n+1}}[d] & \psc{n}{\Cat} \ar^{\ovl{p}}[d]\\
        \catwg{n} \ar_{\Nb{n}}[r] & \funcat{n}{\Set}
        }
    \end{equation*}
\end{definition}
\begin{lemma}\label{lem-seg-pse-fun}
    Let $X\in\segpsc{n}{\Cat}$\; $n\geq 2$. Then for each $j\geq 0$
    \begin{equation*}
      X_{j*}\in\segpsc{n-1}{\Cat}.
    \end{equation*}
\end{lemma}
\begin{proof}
By induction on $n$. Let $X\in\segpsc{2}{\Cat}$. Since $X\in\psc{2}{\Cat}$, for each $j\geq 0$ $X_{j*}\in\psc{}{\Cat}$. By definition of Segalic pseudo-functor, $X_{j0}$ is discrete and for each $r\geq 2$
\begin{equation*}
    X_{jr}\cong \pro{X_{j1}}{X_{j0}}{r}\;.
\end{equation*}
By definition this means that $X_{j*}\in\segpsc{}{\Cat}$. Suppose, inductively, that the lemma holds for $(n-1)$ and let $X\in\segpsc{n}{\Cat}$. For each $j\geq 0$, $X_{j*}\in\psc{n-1}{\Cat}$.

Given $\ur\in\dop{n-1}$ denote $\uk=(j,\ur)\in\dop{n}$. Then for each $1 \leq i\leq n-2$,
\begin{equation*}
    X_{\uk(i+1,0)}=(X_j)_{\ur(i,0)}
\end{equation*}
is discrete since $X\in\segpsc{n}{\Cat}$; further, by hypothesis there are isomorphisms:
\begin{align*}
    (X_j)_{\ur} = X_{\uk}\cong & \pro{X_{\uk(i+1,1)}}{X_{\uk(i+1,0)}}{k_{i+1}}\cong\\
    \cong & \pro{(X_j)_{\ur(i,1)}}{(X_j)_{\ur(i,0)}}{r_i}\;.
\end{align*}
To show that $X_j\in\segpsc{n-1}{\Cat}$ it remains to show that $\p{n}X_j\in\catwg{n-1}$ where
\begin{equation*}
    (\p{n}X_j)_{\ur}=p X_{j\ur}
\end{equation*}
for each $\ur\in\dop{n-1}$.

Since $X\in\segpsc{n}{\Cat}$, by definition $\p{n+1}X\in\catwg{n}$ where $(\p{n+1}X)_{\uk}=p X_{\uk}$ for all $\uk\in\dop{n}$. We also observe that, for each $j\geq 0$
\begin{equation}\label{eq-lem-seg-pse-fun}
    (\p{n+1}X)_j=\p{n}X_j
\end{equation}
since, for each $\ur\in\dop{n-1}$,
\begin{equation*}
    (\p{n+1}X)_{j\ur}=p X_{j\ur}=(\p{n}X_j)_{\ur}\;.
\end{equation*}
Since $\p{n+1}X\in\catwg{n}$, $(\p{n+1}X)_j\in\catwg{n-1}$ so by \eqref{eq-lem-seg-pse-fun} we conclude that $\p{n}X_j\in\catwg{n-1}$ as required.
\end{proof}
\begin{example}\label{ex-pseudo-funct}
Let $X\in\segpsc{2}{\Cat}$. Then $X\in\Ps\funcat{2}{\Cat}$ with $X_{k0}, X_{0s}$ discrete categories for each $k,s\in\Dop$ and $\p{2}X\in\Nb{2}\catwg{2}$ where for all $i,j\geq 0$
\begin{equation*}
  (\p{3}X)_{ij}=p X_{ij}\;.
\end{equation*}
Below are two illustrative figures: in Figure \ref{cornerXsegps} on page \pageref{fig-seg-ps} we have depicted the corner of a pseudo-functor $X\in\segpsc{2}{\Cat}$ with pseudo-commuting squares containing the symbol $\cong$ since the simplicial relations hold only up to isomorphisms; the structures in red color  are discrete categories; Figure \ref{cornerXp3} on page \pageref{fig-seg-ps} depicts the corner of $\p{3}X$, which is a bisimplicial set, double nerve of a weakly globular double category. The structure in green color  is homotopically discrete.

\end{example}

\section{Strictification of Segalic pseudo-functors}\label{sec-stric-pse-fun}
In this section we prove the main result of this chapter, Theorem \ref{the-strict-funct}, that the strictification functor applied to the category of Segalic pseudo-functors gives a weakly globular \nfol category; that is, there is a functor
\begin{equation*}
        \St: \segpsc{n-1}{\Cat}\rw J_n\catwg{n}\cong \catwg{n}.
    \end{equation*}
     The strategy to prove this result is based of the following main steps:
\begin{itemize}
  \item [a)]In Proposition \ref{pro-transf-wg-struc} we show that if an \nfol category is levelwise equivalent (as a diagram in $\funcat{n-1}{\Cat}$) to a Segalic pseudo-functor, then it is a weakly globular \nfol category. The proof of this result uses crucially the criterion for an \nfol category to be weakly globular given in Proposition \ref{pro-crit-ncat-be-wg} b).
      \medskip
  \item [b)] We show in the proof of Theorem \ref{the-strict-funct} that the strictification of a Segalic pseudo-functor is an \nfol category and that it satisfies the hypotheses of Proposition \ref{pro-transf-wg-struc}. The proof depends on some properties of the monad corresponding to Segalic pseudo-functors which we establish in Lemma \ref{lem-strict-mon-wg-seg-ps-fun}.
      \medskip
   \item [c)]  We immediately deduce from a) and b) that the strictification of a Segalic pseudo-functor is a weakly globular \nfol category.
\end{itemize}

\bigskip
\begin{proposition}\label{pro-transf-wg-struc}
    Let $H\in \segpsc{n-1}{\Cat}$ and let $L\in\Catn$ be such that there is an equivalence of categories $(J_n L)_{\uk}\simeq H_{\uk}$ for all $\uk\in\Dmenop$, then
    \begin{itemize}
      \item [a)] $L\in\catwg{n}$.
      \item [b)] If, further, $H_{\uk}\in \cathd{}$ for all $\uk\in\Dmenop$ and $\pn H\in \cathd{n-1}$, then $L\in\cathd{n}$.
    \end{itemize}
\end{proposition}
\begin{proof}
The proof of a) is based on showing that $X\in \Catn$ satisfies the hypotheses of Proposition \ref{pro-crit-ncat-be-wg} b), which then implies that $X\in\catwg{n}$. The proof of b) is based on showing that $X\in \catwg{n}$ satisfies the hypotheses of Corollary \ref{cor-crit-nequiv-rel}, which then implies $X\in\cathd{n}$.

We proceed by induction on $n$. For $n=2$, if $H\in\segpsc{}{\Cat}$, then by definition $H_0$ is discrete; thus, since by hypothesis there is an equivalence of categories $L_0\simeq H_0$, $L_0\in\cathd{}$. By hypothesis $L_k\simeq H_k$ for all $k\in\Dop$, so $p L_k\cong p H_k$, and therefore $\bar p L\cong \bar p H =\p{2} H$ is the nerve of a category. So by Proposition \ref{pro-crit-ncat-be-wg} b) $L\in\catwg{2}$.

If, further, $H_k\in\cathd{}$ for all $k$ and $\p{2}H\in\cathd{}$, then $L_1 \in \cathd{}$ (since $L_1 \sim H_1$) and $\p{2}L \cong\p{2}H \in \cathd{}$. Therefore, by Corollary \ref{cor-crit-nequiv-rel}, $L\in\cathd{2}$.

Suppose, inductively, that the lemma holds for $(n-1)$ and let $L$ and $H$ be as in the hypothesis a).

We are going to show that $L\in\cat{n}$ satisfies the hypotheses of Proposition \ref{pro-crit-ncat-be-wg} b) which then implies that $L\in\catwg{n}$.

Let $\ur\in\dop{n-2}$ and denote $\uk=(i,\ur)\in\dop{n-1}$. By hypothesis, there are equivalences of categories
\begin{equation*}
    (J_{n-1}L_i)_{\ur}=(J_n L)_{\uk}\simeq H_{\uk}=(H_i)_{\ur}\;.
\end{equation*}
Since $L\in\cat{n}$, $L_i\in\cat{n-1}$ and since $H\in\segpsc{n}{\Cat}$, by Lemma \ref{lem-seg-pse-fun},
\begin{equation*}
H_i\in\segpsc{n-1}{\Cat}.
\end{equation*}

%Then, by definition, for all $i\geq 0$
%%
%\begin{equation*}
%    H_{i\bl}\in \segpsc{n-2}{\Cat}\;.
%\end{equation*}
%%
%Also, for all $\uk\in \dop{n-2}$, we have equivalences of categories
%%
%\begin{equation*}
%    (J_{n-1}L_{i\bl})_{\uk}=L_{\uk(i,0)}\simeq H_{\uk(i,0)}=(H_{i\bl})_{\uk}\;.
%\end{equation*}
%
Thus $H_{i}$ and $L_{i}$ satisfy the inductive hypothesis a) and we conclude that $L_{i}\in \catwg{n-1}$. In particular, this implies that $L_{i0}\in \cathd{n-2}$. Thus, by Proposition \ref{pro-crit-ncat-be-wg} b), to show that $L\in \catwg{n}$ it is enough to prove that $L_{0}\in\cathd{n-1}$ and that $\bar p J_n L\in \Nb{n-1}\catwg{n-1}$.

We have $(H_{0})_{\uk}=H_{\uk(0,0)}$ discrete and
\begin{equation*}
    \p{n-1}H_{0}=(\p{n}H)_0 \in\cathd{n-2}
\end{equation*}
since by hypothesis $\p{n}H\in\catwg{n}$ as $H\in\segpsc{n-1}{\Cat}$. Thus $H_{0}$ and $L_{0}$ satisfy the inductive hypothesis b) and we conclude that $L_{0}\in\cathd{n-1}$. For each $\uk \in\dop{n-2}$,
\begin{equation*}
    (\bar p J_n L)_{\uk} = p L_{\uk} \cong p H_{\uk} =(\bar p J_n H)_{\uk}=(\p{n}H)_{\uk}\;.
\end{equation*}
Since $\p{n}H\in\catwg{n}$, $(\p{n}H)_{\uk}\in\catwg{n-1}$ and thus we conclude that
\begin{equation*}
 \bar p J_n L\in \Nb{n-1}\catwg{n-1}.
\end{equation*}
  By Proposition \ref{pro-crit-ncat-be-wg} b), we conclude that $L\in\catwg{n}$, proving a) at step $n$.

Suppose that $H$ is as in b). By Corollary \ref{cor-crit-nequiv-rel}, to show that $L\in\cathd{n}$, it is enough to show that $L_1 \in\cathd{n-1}$ and $\p{n}L \in \cathd{n-1}$.
For all $\uk\in\dop{n-2}$ there is an equivalence of categories
\begin{equation}\label{eq-eqcat}
    (L_{1})_{\uk}\simeq (H_{1})_{\uk}
\end{equation}
Since by hypothesis $(H_{1})_{\uk}\in\cathd{}$ we conclude from \eqref{eq-eqcat} that $(L_{1})_{\uk}\in\cathd{}$. Further, since $\p{n}H\in \cathd{n-1}$, then
\begin{equation*}
    \p{n-1}H_{1}=(\p{n}H)_{1}\in \cathd{n-2}\;.
\end{equation*}
Thus $L_{1}$ and $H_{1}$ satisfy induction hypothesis and we conclude that $L_{1}\in\cathd{n-1}$. Finally,
\begin{equation*}
    \p{n}L\cong\p{n}H\in\cathd{n-1}\;.
\end{equation*}
Thus by Corollary \ref{cor-crit-nequiv-rel} we conclude that $L\in\cathd{n}$.
\end{proof}
\bigskip
In the next Lemma we show some properties of the monad corresponding to Segalic pseudo-functors. This will be used in the proof of the main result of this chapter, Theorem \ref{the-strict-funct} on the strictification of Segalic pseudo-functors. We refer to Section \ref{sbs-pseudo-functors} in chapter \ref{chap-ps} for background regarding the monad for pseudo-functors and about strictification of pseudo-functors.
\begin{lemma}\label{lem-strict-mon-wg-seg-ps-fun}
    Let $T$ be the monad corresponding to the adjunction given by the forgetful functor
    \begin{equation*}
        U:\funcat{n}{\Cat}\rw [ob(\Dnop),\Cat ]
    \end{equation*}
    and its left adjoint. Let $H\in \segpsc{n}{\Cat}$, then
    \begin{itemize}
      \item [a)] The pseudo $T$-algebra corresponding to $H$ has structure map $h: TUH\rw H$ as follows:
          \begin{equation*}
            (TUH)_{\uk}=\underset{\ur\in\zD^n}{\coprod}\zD^n(\uk,\ur)\times H_{\ur}=\underset{\ur\in\zD^n}{\coprod}\underset{\zD^n(\uk,\ur)}{\coprod} H_{\ur}\;.
          \end{equation*}
          If $f\in\zD^n(\uk,\ur)$, let
          \begin{equation*}
            i_{\ur}=\underset{\zD^n(\uk,\ur)}{\coprod}H_{\ur}\rw \underset{\ur\in\zD^n}{\coprod}\underset{\zD^n(\uk,\ur)}{\coprod} H_{\ur}=(TUH)_{\ur}
          \end{equation*}
          \begin{equation*}
            j_f :H_{\ur}\rw \underset{\zD^n(\uk,\ur)}{\coprod}H_{\ur}\;,
          \end{equation*}
          then
          \begin{equation*}
            h_{\uk}\,i_{\ur}\,j_f= H(f)\;.
          \end{equation*}
      \item [b)] There are functors for each $\uk \in \zD^n$, $0\leq i \leq n$,
      \begin{equation*}
        \pt_{i1}, \pt_{i0}:(TUH)_{\uk(1,i)}\rightrightarrows (TUH)_{\uk(0,i)}
      \end{equation*}
      such that the following diagram commutes
      \begin{equation}\label{fig.pro.strict-mon}
        \xymatrix@R=35pt @C=50pt{
        (TUH)_{\uk(1,i)} \ar^{h_{\uk(1,i)}}[r] \ar^{\pt_{i1}}[d]<1ex> \ar_{\pt_{i0}}[d]<-1ex> & H_{\uk(1,i)} \ar^{d_{i1}}[d]<1ex> \ar_{d_{i0}}[d]<-1ex>  \\
        (TUH)_{\uk(0,i)} \ar_{h_{\uk(0,i)}}[r] & H_{\uk(0,i)}
        }
      \end{equation}
      \item [c)] For all $\uk=(k_1,\ldots, k_n)\in\zD^n$ there are  isomorphisms
      \begin{equation*}
        (TUH)_{\uk}=\pro{(TUH)_{\uk(1,i)}}{(TUH)_{\uk(0,i)}}{k_i}\;.
      \end{equation*}
      \item [d)] The morphism $h_{\uk}:(TUH)_{\uk}\rw H_{\uk}$ is given by
      \begin{equation*}
        h_{\uk}=(h_{\uk(1,i)},\ldots,h_{\uk(1,i)})
      \end{equation*}
    \end{itemize}
\end{lemma}
\begin{proof}\

a) From the general correspondence between pseudo $T$-algebras and pseudo-functors, the pseudo $T$-algebra corresponding to $H$ has a structure map $h:TUH\rw H$ as stated. The rest follows from the fact that, if $X$ is a set and $\clC$ is a category, $X\times \clC \cong \uset{X}{\coprod}\clC$.\bk

b) Let $\nu_j:[0]\rw [1]$, $\nu_j(0)=0$. $\nu_j(1)=i$ for $j=0,1$ and let $\zd_{ij}:\uk(0,i)\rw\uk(1,i)$ be given by
\begin{equation*}
    \zd_{ij}(k_s)=\left\{
    \begin{array}{ll}
     k_s, & \hbox{$s\neq i$;} \\
     \nu_{j}(k_i), & \hbox{$s=i$.}
     \end{array}
     \right.
\end{equation*}
Given $f\in \zD^n(\uk(1,i), \ur)$ let $j_f$ and $i_{\ur}$ be the corresponding coproduct injections as in a). Let
\begin{equation*}
    \pt_{ij}:(TUH)_{\uk(1,i)}\rw (TUH)_{\uk(0,i)}
\end{equation*}
be the functors determined by
\begin{equation*}
    \pt_{ij}i_{\ur}j_f = i_{\ur} j_{f\zd_{ij}}\;.
\end{equation*}
From a), we have
\begin{align*}
   & h_{\uk(0,i)}\,\pt_{ij}\, i_{\ur}\,j_f = h_{\uk(0,i)}\, i_{\ur} \, j_{f\zd_{ij}}= H(f\zd_{ij})\\
& d_{ij}\,h_{\uk(1,i)}\, i_{\ur}\, j_{f} = H(\zd_{ij})\, H(f)\;.
\end{align*}
Since $H\in\Ps\funcat{n}{\Cat}$ and $H_{\uk(0,i)}$ is discrete, it is
 \begin{equation*}
 H(f\zd_{ij})=H(\zd_{ij})H(f)
 \end{equation*}
  so that, from above,
\begin{equation*}
    h_{\uk(0,i)}\pt_{ij}i_{\ur}j_f = d_{ij} h_{\uk(1,i)}i_{\ur} j_f
\end{equation*}
for each $\ur,f$. We conclude that
\begin{equation*}
     h_{\uk(0,i)}\pt_{ij}=d_{ij} h_{\uk(1,i)}\;.
\end{equation*}
That is, diagram \eqref{fig.pro.strict-mon} commutes.\bk

c) Since, for each $k_i\geq 2$
\begin{equation*}
    [k_i] = [1]\uset{[0]}{\coprod}\; \oset{k_i}{\ldots}\;
    \uset{[0]}{\coprod}[1]
\end{equation*}
we have, for each $\uk = (k_1,...,k_n)\in\zD^n$ with $k_i\geq 2$ for all $1\leq i\leq n$
\begin{equation*}
    \uk = \uk(1,i)\uset{\uk(0,i)}{\coprod}\; \oset{k_i}{\ldots}\;
    \uset{\uk(0,i)}{\coprod}\uk(1,i)\;.
\end{equation*}
Therefore there is a bijection
\begin{equation*}
    \zD^n(\uk,\ur)=\pro{\zD^n\big(\uk(1,i),\ur\big)}{\zD^n(\uk(0,i),\ur)}{k_i}\;.
\end{equation*}
From the proof of b), the functors
 \begin{equation*}
 \pt_{ij}:(TUH)_{\uk(1,i)}\rw (TUH)_{\uk(0,i)}
 \end{equation*}
  for $j=0,1$ are determined by the functors
\begin{equation*}
    (\overline{\zd_{ij}}, id):\zD^n\big(\uk(1,i),\ur\big)\times H_{\ur} \rw \zD^n\big(\uk(0,i),\ur\big)\times H_{\ur}
\end{equation*}
where $\overline{\zd_{ij}}(g)=g\zd_{ij}$ for $g\in \zD^n\big(\uk(1,i),\ur\big)$ and
\begin{align*}
    & (TUH)_{\uk(1,i)}=\uset{\ur}{\coprod}\zD^n\big(\uk(1,i),\ur\big)\times H_{\ur}\\
    & (TUH)_{\uk(0,i)}=\uset{\ur}{\coprod}\zD^n\big(\uk(0,i),\ur\big)\times H_{\ur}\;.
\end{align*}
It follows that
\begin{align*}
    & \pro{(TUH)_{\uk(1,i)}}{(TUH)_{\uk(0,i)}}{k_i}=\\
    & \uset{\ur}{\coprod}\{ \pro{\zD^n\big(\uk(1,i),\ur\big)}{\zD^n(\uk(0,i),\ur)}{k_i} \}\times H_{\ur}=\\
    & \uset{\ur}{\coprod}\zD^n(\uk,\ur)\times H_{\ur}=(TUH)_{\uk}\;.
\end{align*}
This proves c). \bk

d) From a),
 \begin{equation*}
 {h_k\,i_{r}\,j_f}=H(f)
 \end{equation*}
  for $f \in \zD^n(\uk,\ur)$. Let $f$ correspond to $(\zd_1,\ldots,\zd_{ki})$ in the isomorphism
\begin{equation*}
    \zD^n(\uk,\ur)=\pro{\zD^n\big(\uk(1,i),\ur\big)}{\zD^n(\uk(0,i),\ur)}{k_i}\;.
\end{equation*}
Then $j_f=(j_{\zd_1},\ldots,j_{\zd_{k_i}})$. Since
\begin{equation*}
    H_{\uk} \cong \pro{H_{\uk(1,i)}}{H_{\uk(0,i)}}{k_i}
\end{equation*}
then $H(f)$ corresponds to $(H(\zd_1),\ldots,H(\zd_{k_i}))$ with $p_i H(f)=H(\zd_i)$. Then for all $f$ we have
\begin{align*}
    h_{\uk}\,i_{\ur}\,j_f &=(H(\zd_1),\ldots,H(\zd_{k_i}))=(h_{\uk(1,i)}\,i_{\ur}\,j_{\zd_1},\ldots, h_{\uk(1,i)}\,i_{\ur}\,j_{\zd_{k_i}})=\\
    & =(h_{\uk(1,i)},\ldots,h_{\uk(1,i)})\,i_{\ur}(j_{\zd_1},\ldots, j_{\zd_{k_i}})=(h_{\uk(1,i)},\ldots,h_{\uk(1,i)})\,i_{\ur}\,j_f\;.\\
\end{align*}
It follows that $h_{\uk}=(h_{\uk(1,i)},\ldots,h_{\uk(1,i)})$.
\end{proof}
We now show the main result of this Chapter that the strictification of a Segalic pseudo-functor is a weakly globular \nfol category. This result will be used in Chapter \ref{chap7} in the construction of the rigidification functor from $\tawg{n}$ to $\catwg{n}$.

\begin{theorem}\label{the-strict-funct}
    The strictification functor
    \begin{equation*}
        \St  : \psc{n-1}{\Cat}\rw \funcat{n-1}{\Cat}
    \end{equation*}
    restricts to a functor
    \begin{equation*}
        \St: \segpsc{n-1}{\Cat}\rw \catwg{n}
    \end{equation*}
    (where we identify \catwg{n} with the image $J_{n}\catwg{n}$ of the multinerve functor $J_n:\catwg{n}\rw\funcat{n-1}{\Cat}$\;).

    Further, for each $H\in \segpsc{n-1}{\Cat}$ and $\uk\in\dop{n-1}$, the map
     \begin{equation*}
       (\St H)_{\uk}\rw H_{\uk}
     \end{equation*}
      is an equivalence of categories.
\end{theorem}
\begin{proof}
Let $h:TUH\rw UH$ be as in Section \ref{sbs-pseudo-functors}. As recalled there, to construct the strictification $L=\St H$ of a pseudo-functor $H$ we need to factorize $h=gv$ in such a way that for each $\uk \in \dop{n-1}$, $h_{\uk}$ factorizes as
\begin{equation*}
    (TUH)_{\uk}\xrw{v_{\uk}} L_{\uk} \xrw{g_{\uk}} H_{\uk}
\end{equation*}
with $v_{\uk}$ bijective on objects and $g_{\uk}$ fully faithful. As explained in \cite{PW}, $g_{\uk}$ is in fact an equivalence of categories.

Since the bijective on objects and fully faithful functors form a factorization system in $\Cat$, the commutativity of \eqref{fig.pro.strict-mon} implies that there are functors
\begin{equation*}
    \tilde d_{ij}:L_{\uk(1,i)}\rightrightarrows L_{\uk(0,i)}\qquad j=0,1
\end{equation*}
such that the following diagram commutes:
\begin{equation*}
\xymatrix{
(TUH)_{\uk(1,i)} \ar^{v_{\uk(1,i)}} [rr] \ar_{\pt_{i0}}[d]<-2ex>\ar^{\pt_{i1}}[d]&& L_{\uk(1,i)} \ar^{g_{\uk(1,i)}} [rr] \ar_{\tilde d_{i0}}[d]<-2ex>\ar^{\tilde d_{i1}}[d] && H_{\uk(1,i)}\ar_{d_{i0}}[d]<-2ex>\ar^{d_{i1}}[d]\\
(TUH)_{\uk(0,i)} \ar_{v_{\uk(0,i)}} [rr] && L_{\uk(0,i)} \ar_{g_{\uk(0,i)}} [rr] && H_{\uk(0,i)}\;.
}
\end{equation*}
By Lemma \ref{lem-strict-mon-wg-seg-ps-fun}, $h_{\uk}$ factorizes as
\begin{align*}
    (TUH)_{\uk}=\; & \pro{(TUH)_{\uk(1,i)}}{(TUH)_{\uk(0,i)}}{k_i}\rw\\
    &\xrw{(v_{\uk(1,i)},\ldots,v_{\uk(1,i)})}\pro{L_{\uk(1,i)}}{L_{\uk(0,i)}}{k_i}\rw\\
    & \xrw{(g_{\uk(1,i)},\ldots,g_{\uk(1,i)})}\pro{H_{\uk(1,i)}}{H_{\uk(0,i)}}{k_i}\cong H_{\uk}\;.
\end{align*}
Since $v_{\uk(1,i)}$ and $v_{\uk(0,i)}$ are bijective on objects, so is $(v_{\uk(1,i)},\ldots,v_{\uk(1,i)})$. Since  $g_{\uk(1,i)}$, $g_{\uk(0,i)}$ are fully faithful, so is $(g_{\uk(1,i)},\ldots,g_{\uk(1,i)})$. Therefore the above is the required factorization of $h_{\uk}$ and we conclude that
\begin{equation*}
    L_{\uk}\cong\pro{L_{\uk(1,i)}}{L_{\uk(0,i)}}{k_i}\;.
\end{equation*}
Since $L\in\funcat{n-1}{\Cat}$ by Lemma \ref{lem-multin-iff} this implies that $L\in\cat{n}$. By \cite{PW}, there is an equivalence of categories $L_{\uk}\simeq H_{\uk}$ for all $\uk\in\dop{n-1}$. Therefore, by Proposition \ref{pro-transf-wg-struc}, $L\in\catwg{n}$.
\end{proof}
%%
%%%%%%%%%%%%%%%%%%%%%%%%%%%%%%%%%%%%%%%%%%%%%%%%%%%%%%%%%%%%%%%%%%%%%
\clearpage

\vspace{10mm}

\begin{figure}[ht]
  \centering
\begin{equation*}
\def\objectstyle{\scriptstyle}
\def\labelstyle{\scriptstyle}
%\entrymodifiers={++[o]}
\ssr \xymatrix@R=35pt@C=20pt{
   \bm{\tens{(\tens{X_{11}}{X_{01}})}{(\tens{X_{10}}{X_{00}})}} \ar@<1ex>[r] \ar[r] \ar@<-1ex>[r] \ar@<1ex>[d] \ar[d] \ar@<-1ex>[d] \ar@{} [dr] |{\cong}
   &  \bm{\tens{X_{11}}{X_{10}}} \ar@<-0.5ex>[r] \ar@<0.5ex>[r] \ar@<1ex>[d] \ar[d] \ar@<-1ex>[d] \ar@{} [dr] |{\cong}
   & \textcolor[rgb]{1.00,0.00,0.00}{ \bm{(\tens{X_{01}}{X_{00}})}}  \ar@<1ex>@[red][d] \ar@[red][d] \ar@<-1ex>@[red][d] \\
%%%%%%
     \bm{\cdots\quad\tens{X_{11}}{X_{01}}\quad}\ar[r] \ar@<1ex>[r]  \ar@<-1ex>[r] \ar@<0.5ex>[d] \ar@<-0.5ex>[d] \ar@{} [dr] |{\cong}
    &  \bm{\quad  X_{11}\quad }\ar@<-0.5ex>[r] \ar@<0.5ex>[r]   \ar@<0.5ex>[d] \ar@<-0.5ex>[d] \ar@{} [dr] |{\cong}
    &  \textcolor[rgb]{1.00,0.00,0.00}{\bm{\quad X_{01}\quad}}  \ar@<0.5ex>@[red][d] \ar@<-0.5ex>@[red][d]  \\
%%%%%%
      \bm{\cdots\quad\tens{X_{10}}{X_{00}}\quad} \ar[r] \ar@<1ex>[r]  \ar@<-1ex>[r]
    &  \bm{\quad X_{10}\quad} \ar@<-0.5ex>[r] \ar@<0.5ex>[r]
    &  \textcolor[rgb]{1.00,0.00,0.00}{\bm{\quad X_{00}\quad}} \\
}
\end{equation*}
  \caption{Picture of the corner of $X\in\segpsc{2}{\Cat}$}  \label{cornerXsegps}
\end{figure}

\vspace{10mm}

\definecolor{greenbl}{rgb}{0.0, 0.59, 0.0}

\begin{figure}[ht]
  \centering
\begin{equation*}
\def\objectstyle{\scriptstyle}
\def\labelstyle{\scriptstyle}
%\entrymodifiers={+++[o]}
 \xymatrix@R=20pt@C=20pt{
  \bm{\quad\cdots\quad} \ar@<1ex>[r] \ar[r] \ar@<-1ex>[r] \ar@<1ex>[d] \ar[d] \ar@<-1ex>[d]
   &   \bm{\tens{pX_{11}}{pX_{10}}} \ar@<-0.5ex>[r] \ar@<0.5ex>[r] \ar@<1ex>[d] \ar[d] \ar@<-1ex>[d]
   & \textcolor[rgb]{0.00,0.59,0.00}{\bm{\tens{X_{01}}{X_{00}}}}  \ar@<1ex>@[greenbl][d] \ar@[greenbl][d] \ar@<-1ex>@[greenbl][d] \\
%%%%%%%%%%%%%%%%%%%%
    \bm{\cdots\quad \tens{pX_{11}}{pX_{01}} \quad} \ar[r] \ar@<1ex>[r]  \ar@<-1ex>[r] \ar@<0.5ex>[d] \ar@<-0.5ex>[d]
    & \bm{p X_{11}} \ar@<-0.5ex>[r] \ar@<0.5ex>[r]   \ar@<0.5ex>[d] \ar@<-0.5ex>[d]
    & \textcolor[rgb]{0.00,0.59,0.00}{\bm{ X_{01}}}  \ar@<0.5ex>@[greenbl][d] \ar@<-0.5ex>@[greenbl][d]  \\
%%%%%%%%%%%%%%%%%%%%%%%
     \bm{\cdots\quad \tens{pX_{10}}{pX_{00}} \quad} \ar[r] \ar@<1ex>[r]  \ar@<-1ex>[r]
    & \bm{p X_{10}} \ar@<-0.5ex>[r] \ar@<0.5ex>[r]
    & \textcolor[rgb]{0.00,0.59,0.00}{\bm{ X_{00}}} \\
}
\end{equation*}
  \caption{Picture of the corner of  $\p{3}X$  \label{cornerXp3}, for $X\in\segpsc{2}{\Cat}$}.
\end{figure}

\label{fig-seg-ps}

\clearpage
%%%%%%%%%%%%%%%%%%%%%%%%%%%%%%%%%%%%%%%%%%%%%%%%%%%%%%%%%%%%%%%%%%%%%%%%%%%

\part{Weakly globular Tamsamani $\pmb{n}$-categories and their rigidification}\label{part-4}

\begin{figure}[ht]
  \centering
\vspace{10mm}
\begin{tikzpicture}[node distance = 15mm and 20mm]% First vertical distance and second horizontal from border
\node [block] (A_1_1) {Definition \ref{def-wg-ps-cat}\\Category $\tawg{n}$ of weakly globular Tamsamani $n$-categories};
\node [block, right= of A_1_1] (A_1_2) {Sections \ref{sec-func-qn}, \ref{sec-pull-qn}\\Functors\\$\q{n}:\tawg{n}\rw\tawg{n-1}$\\$\q{n}:\catwg{n}\rw\catwg{n-1}$\\ and their properties};
\node [block3, below= of A_1_1] (A_2_1) {Definition \ref{def-ind-sub-ltawg}\\The category $\lta{n}$};
\node [block2, below= of A_1_2] (A_2_2) {Theorem \ref{the-repl-obj-1}\\Approximating $\tawg{n}$ with $\lta{n}$};
\node [nullblock, below of= A_2_2, node distance=3.7cm] (null) {};

\node [block2, left of=null, node distance=6.8cm ] (A_3_1) {Theorem \ref{the-XXXX}\\Functor\\$\tr{n}:\lta{n}\rw\segpsc{n-1}{\Cat}$};

\node [block2, right  of= null, node distance=8mm] (A_3_2) {Theorem \ref{the-strict-funct}\\$\St:\segpsc{n-1}{\Cat}\rw\catwg{n}$};

\node [nullblock, right of= A_3_1, node distance=3.8cm] (null1) {};

\node [block5, below  of= null1, node distance=38mm] (A_4_1) {Theorem \ref{the-funct-Qn} Rigidification functor\\ \  \\ $Q_2:\tawg{2}\xrw{\tr{2}}\segpsc{}{\Cat}\xrw{\St}\catwg{2}$ \\ \  \\ \smallskip  \mbox{For $n>2$} \hfill\rule{0ex}{0ex} \\ $\Qn:\tawg{n}\xrw{P_n}\lta{n}\xrw{\tr{n}}\segpsc{n-1}{\Cat}\xrw{\St}\catwg{n}$};

% Draw edges
\path [line] (A_1_1) -- (A_1_2);
\path [line] (A_1_2) -- (A_2_2);
\path [line] (A_2_1) -- (A_2_2);
\path [line] (A_2_1) -- (A_3_1.north -| A_2_1);
\path [line] (A_3_1) -- (A_4_1.north -| A_3_1);
\path [line] (A_3_2) -- (A_4_1.north -| A_3_2);

\draw [line] (A_2_2.south)++(-3.0,0) -- (A_4_1);

\end{tikzpicture}
\bk

  \caption{The construction of the rigidification functor $\Qn$.}
  \label{FigIntro-7}
\end{figure}
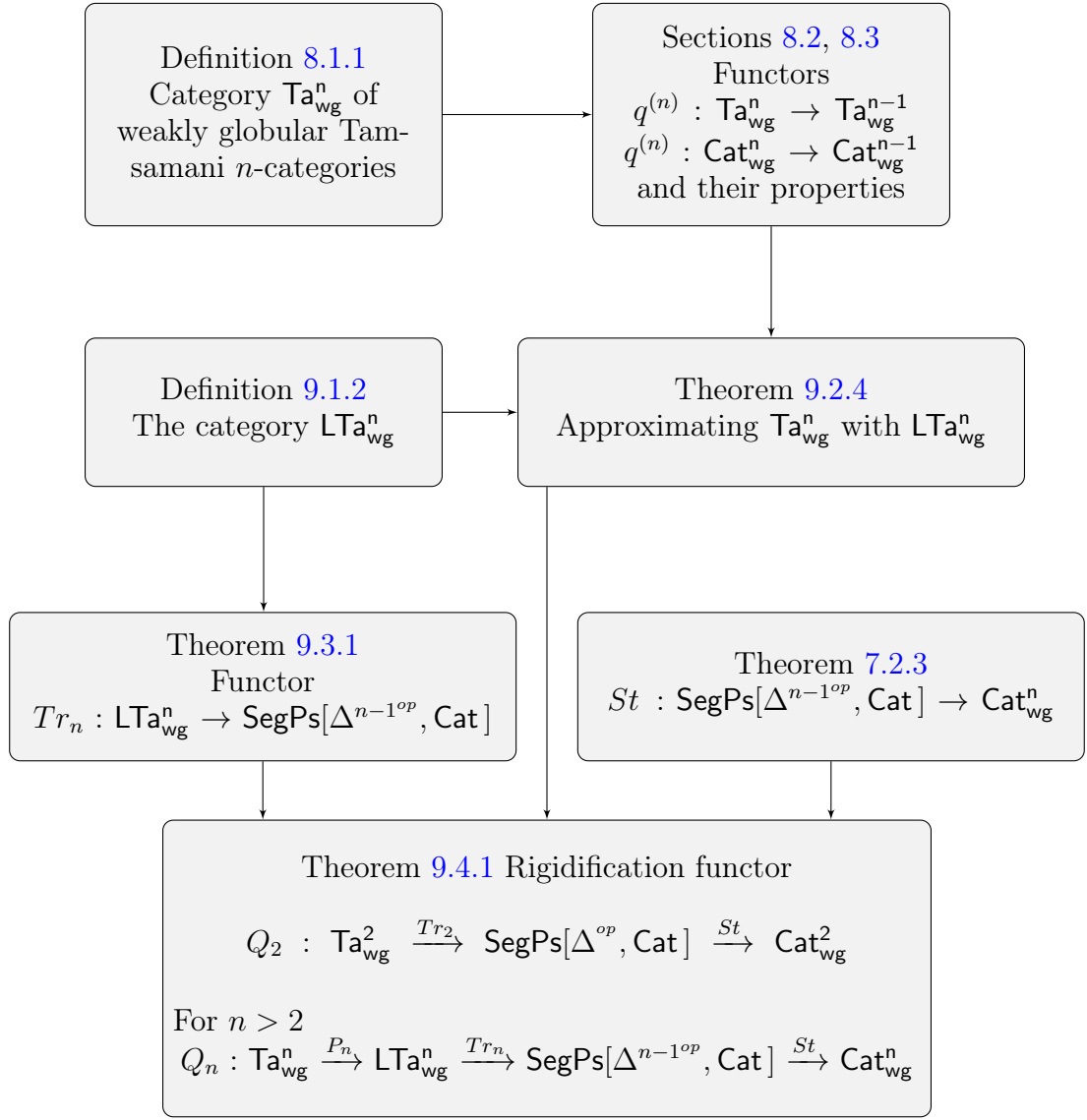

\clearpage
  Part \ref{part-4} is devoted to the category $\tawg{n}$ of weakly globular Tamsamani $n$-categories. The main result of this part is Theorem \ref{the-funct-Qn} which constructs the rigidification functor
    \begin{equation*}
        Q_n:\tawg{n} \rw \catwg{n}
    \end{equation*}
    A schematic summary of the main results of this part is given in Figure \ref{FigIntro-7}.

    In Chapter \ref{chap6} we define the category $\tawg{n}$ of weakly globular Tamsamani $n$-categories and we study its properties. The idea of the category $\tawg{n}$ is explained in Section \ref{sub-idea-tawgn}, before the final definition. Some of the properties of the category $\tawg{n}$ are used throughout the rest of the work: Lemma \ref{lem-x-in-tawg-x-in-catwg} gives a sufficient criterion for an object of $\tawg{n}$ to be in $\catwg{n}$, Proposition \ref{pro-n-equiv} establishes important properties of $n$-equivalences, Proposition \ref{pro-crit-lev-nequiv} gives a sufficient criterion for a $n$-equivalence to be a levelwise equivalence of categories. Proposition \ref{pro-post-trunc-fun} establishes the existence of the functor
     \begin{equation*}
        \qn:\tawg{n}\rw \tawg{n-1}.
     \end{equation*}
      The properties of the latter in relation to certain pullbacks are further studied in Section \ref{sec-pull-qn}, and they play a key role in the proof of Theorem \ref{the-repl-obj-1}, leading to the rigidification functor $Q_n$.

    The construction of the rigidification functor when $n=2$ is quite straightforward and was already done in \cite{PP}. The construction of $Q_n$ when $n>2$ is much more complex and is new to this work: it needs in particular the subcategory $\lta{n}$ of $\tawg{n}$. We introduce the idea of this subcategory in Section \ref{subs-idea-ltan}, before the formal definition.

     In Chapter \ref{chap7} we prove two important results involving this subcategory which are used in the construction of the rigidification functor: Theorem \ref{the-repl-obj-1} and Theorem \ref{the-XXXX}.

      Theorem \ref{the-repl-obj-1} establishes a procedure to approximate, up to $n$-equivalence, an object of $\tawg{n}$ with an object of $\lta{n}$: its proof is based on the properties of the pullback constructions of Section \ref{sec-pull-qn} as well as on the criterion given in Proposition \ref{pro-crit-lev-nequiv} (used in the proof of Lemma \ref{lem-jn-alpha}, leading to Theorem \ref{the-repl-obj-1}). The main steps needed in these constructions are explained in Section \ref{subs-approx-idea}, before the formal proofs are given.

   In  Theorem \ref{the-XXXX} we construct the functor
     \begin{equation*}
        Tr_{n}: \lta{n} \rw \segpsc{n-1}{\Cat}
    \end{equation*}
    The idea if the functor $Tr_n$ is explained in Section \ref{subs-idea-trn}, before the formal proofs.
     The proof of  Theorem \ref{the-XXXX}  relies on a technique to produce pseudo-functors that is an instance of  'transport of structure along an adjunction' recalled in Section \ref{transport-structure}, as well as on the properties of the category $\lta{n}$, and especially Proposition \ref{pro-maps-nu-eqcat}, giving a sufficient criterion for an object of $\tawg{n}$ to be in $\lta{n}$.

We finally construct the rigidification functor $Q_n: \tawg{n} \rw \catwg{n}$. The idea of the construction of $Q_n$ is given in Section \ref{sub-idea-qn}, before the formal proof of Theorem \ref{the-funct-Qn}.
In the case $n=2$, the rigidification functor $Q_2$ is the composite
\begin{equation*}
    Q_2: \tawg{2}\xrw{Tr_{2}} \segpsc{}{\Cat} \xrw{\St} \catwg{2}
\end{equation*}
where $Tr_{2}$ is as in Theorem \ref{the-XXXX}.
    When $n>2$ the functor $Q_n$ is given as a composite
    \begin{equation*}
    \tawg{n}\xrw{P_n} \lta{n} \xrw{\Tr_{n}} \segpsc{n-1}{\Cat}\xrw{\St} \catwg{n}
\end{equation*}
This relies on Theorem \ref{the-repl-obj-1} (for the construction of the functor $P_n$), Theorem \ref{the-XXXX} (for the functor $Tr_n$) and Theorem \ref{the-strict-funct} (for the functor $\St$).\bk

\clearpage

%%%%%%%%%%%%%%%%%%%%%%%%%%%%%%%%%%%%%%%%%%%%%%%%%%%%%%%%%%%%%%%%%%%%%%%%%%%%
\chapter{Weakly globular Tamsamani $\pmb{n}$-categories}\label{chap6}

In this chapter we introduce the most general of the three Segal-type models of this work, the category $\tawg{n}$ of weakly globular Tamsamani $n$-categories, and we study its main properties.

Weakly globular Tamsamani $n$-categories, like weakly globular \nfol categories, satisfy the weak globularity condition, which is formulated using the notion of  homotopically discrete \nfol category introduced in Chapter \ref{chap3}. However, unlike in weakly globular $n$-fold categories, the Segal maps are no longer required to be isomorphisms. The resulting structure is therefore no longer an $n$-fold category.

The behaviour of the compositions is controlled by the induced Segal maps, which for $X\in\tawg{n}$ have the form for each $k\geq 2$\;
\begin{equation*}
  X_k\rw\pro{X_1}{X_0^d}{k}.
\end{equation*}
where $X_0^d$ is the discretization of the homotopically discrete $(n-1)$-fold category $X_0$. In a weakly globular Tamsamani $n$-category the  induced Segal maps are required to be suitable higher categorical equivalences. The formal definition of $\tawg{n}$ is given by induction on dimension.

The category of weakly globular Tamsamani $n$-categories satisfying the globularity condition is the same as the category $\ta{n}$ of Tamsamani $n$-categories, though the original definition \cite{Ta} did not use the bigger category $\tawg{n}$, which is new to this work.

Weakly globular Tamsamani $n$-categories are a generalization both of the Tamsamani model and of weakly globular $n$-fold categories, and there are embeddings
\begin{equation*}
  \catwg{n}\subset\tawg{n} \qquad \text{and}\qquad \ta{n}\subset\tawg{n}\;.
\end{equation*}

This chapter is organized as follows. In Section \ref{sec-WG-Tam-cat} we introduce the category $\tawg{n}$ of weakly globular Tamsamani $n$-categories and of $n$-equivalences and we prove their main properties. In particular we show in Proposition \ref{pro-crit-lev-nequiv} a criterion for an $n$-equivalence in $\tawg{n}$ to be a levelwise equivalence of categories. This is used repeatedly in later chapters.

In Section \ref{sec-func-qn} we continue the study of the category $\tawg{n}$ and we show in Proposition \ref{pro-post-trunc-fun} the existence of the functor
\begin{equation*}
    \q{n}:\tawg{n}\rw \tawg{n-1}\;.
\end{equation*}
We also show in Corollary \ref{pro-post-wg-ncat} that this restricts to a functor
\begin{equation*}
    \q{n}:\catwg{n}\rw \catwg{n-1}\;.
\end{equation*}
The functor $\q{n}$ generalizes to higher dimensions the connected components functor $q:\Cat\rw\Set$. For each $X\in\tawg{n}$ there is a map, natural in $X$
 \begin{equation*}
 \zg\up{n}: X\rw \di{n}\q{n}X
 \end{equation*}
  where $\di{n}$ is the inclusion of $\tawg{n-1}$ in $\tawg{n}$ as a discrete structure in the top dimension. In Section \ref{sec-pull-qn} we study the properties of pullbacks along this map. This will be needed for one of the main constructions in Chapter \ref{chap7} given in the proof of Theorem \ref{the-repl-obj-1}, which will lead to the rigidification functor $Q_n$.

%%%%%%%%%%%%%%%%%%%%%%%%%%%%%%%%%%%%%%%%%%%%%%%%%%%%%%%%%%%%%%%%%%%%%%%%%%%%%%%%%%%%%%%%%%%%%%%%%%%%%%%%%

\section{Weakly globular Tamsamani $\pmb{n}$-categories}\label{sec-WG-Tam-cat}
In this  section we introduce the category $\tawg{n}$ of weakly globular Tamsamani $n$-categories and discuss their properties.

\subsection{The idea of weakly globular Tamsamani $\pmb{n}$-categories }\label{sub-idea-tawgn}
The definition of the category $\tawg{n}$ is inductive on dimension starting with $\Cat$ when $n=1$. In dimension $n>1$, a weakly globular Tamsamani $n$-category $X$ is a simplicial object in the category $\tawg{n-1}$ satisfying additional conditions. These conditions encode the weakness in the structure in two ways: one is the weak globularity condition, requiring $X_{0}$ to be a homotopically discrete $\nm$-fold category. The second is the induced Segal maps condition, requiring that the induced Segal maps for all $k\geq 2$
\begin{equation}\label{eq1-sec-prelim}
  \hmu{k}:X_k \rw \pro{X_1}{X_0^d}{k}
\end{equation}
are $\nm$-equivalences.

By unravelling the inductive Definition \ref{def-wg-ps-cat} we obtain the embedding
\begin{equation*}
  J_n : \tawg{n}\rw \funcat{n-1}{\Cat}\;.
\end{equation*}
The inductive Definition \ref{def-wg-ps-cat} also requires the existence of a truncation functor
\begin{equation*}
  \p{n}:\tawg{n}\rw \tawg{n-1}
\end{equation*}
obtained by applying levelwise to $J_n X$ the isomorphism classes of object functor $p:\Cat\rw\Set$.

 This truncation functor is used  to define $\nequ$s in $\tawg{n}$. This notion is a higher dimensional generalization of a functor which is fully faithful and essentially surjective on objects.
\subsection{The formal definition of the category $\pmb{\tawg{n}}$}
\begin{definition}\rm\label{def-wg-ps-cat}\index{Functor!- $\p{n}$}
    We define the category $\tawg{n}$ by induction on $n$. For $n=1$, $\tawg{1}=\Cat$ and 1-equivalences are equivalences of categories. We denote by $\p{1}=p:\Cat\rw \Set$ the isomorphism classes of object functor.

    \index{Weakly globular!- Tamsamani $n$-category} \index{n-equivalences}

    Suppose, inductively, that we defined for each $1 < k\leq n-1$
    \begin{equation*} \index{Functor!- $J_n$}
        \xymatrix{J_k: \tawg{k}\;\ar@{^{(}->}^(0.4){}[r] & \;\funcat{k-1}{\Cat}}
    \end{equation*}
    and $k$-equivalences in $\tawg{k}$ as well as a functor
    \begin{equation*}
        \p{k}:\tawg{k}\rw \tawg{k-1}
    \end{equation*}
    sending $k$-equivalences to $(k-1)$-equivalences and making the following diagram commute:
    \begin{equation}\label{eq-wg-ps-cat}
    \xymatrix@C=30pt{
    \tawg{k} \ar^{J_{k}}[rr]\ar_{\p{k}}[d] && \funcat{k-1}{\Cat} \ar^{\ovl{p}}[d]\\
    \tawg{k-1}  \ar^{\Nb{k-1}}[rr] && \funcat{k-1}{\Set}
    }
    \end{equation}
    An object $X\in \tawg{k}$ is called discrete if $J_k X$ is constant taking value in a discrete category.

    Define $\tawg{n}$ to be the full subcategory of $\funcat{}{\tawg{n-1}}$ whose objects $X$ are such that
        \begin{itemize}
      \item [a)] \textsl{Weak globularity condition} \;$X_0\in\cathd{n-1}$; Here $X_0 \in \funcat{n-2}{\Cat}$ so we identify $\catwg{n-1}$ with $J_{n-1} \catwg{n-1}$.
      \bk \index{Weak globularity condition}

      \item [b)] \textsl{Induced Segal maps condition}. For all $s\geq 2$ the induced Segal maps \index{Induced Segal maps condition}
      \begin{equation*}
       X_s  \rw \pro{X_1}{X^d_0}{s}
      \end{equation*}
      (induced by the map $\zg:X_0\rw X_0^d$) are $(n-1)$-equivalences.
    \end{itemize}
\medskip
To complete the inductive step, we need to define the \emph{truncation functor} $\p{n}$ and $n$-equivalences. \index{Truncation functor}
Note that the functor
\begin{equation*}
    \ovl{\p{n-1}}:\funcat{}{\tawg{n-1}}\rw \funcat{}{\tawg{n-2}}
\end{equation*}
restricts to a functor
\begin{equation*}
    \p{n}:\tawg{n}\rw \tawg{n-1}\;.
\end{equation*}
In fact, since $X_0\in \catwg{n-1}$, $(\p{n}X)_0=\p{n-1}X_0 \in \catwg{n-2}$. Further, by \eqref{eq-wg-ps-cat} $\p{n-1}$ preserves pullbacks over discrete objects (as the same is true for $p$, see Lemma \ref{lem-p-pres-fib-pro}) so that
\begin{equation*}
  \p{n-1}(\pro{X_1}{X_0^d}{s}) \cong\pro{\p{n-1}X_1}{(\p{n-1}X_0^d)}{s}\; .
\end{equation*}
 Further,
 \begin{equation*}
   \p{n-1}X_0^d = (\p{n-1}X_0)^d
 \end{equation*}
   and $\p{n-1}$ sends $(n-1)$-equivalences to $(n-2)$-equivalences.

 Therefore, the induced Segal maps for $s\geq 2$
\begin{equation*}
    X_s \rw\pro{X_1}{X_0^d}{s}
\end{equation*}
being $(n-1)$-equivalences, give rise to $(n-2)$-equivalences
\begin{equation*}
     \p{n-1}X_s \rw \pro{\p{n-1}X_1}{(\p{n-1}X_0)^d}{s}\; .
\end{equation*}
This shows that $\p{n}X \in \tawg{n-1}$. It is immediate that \eqref{eq-wg-ps-cat} holds at step $n$.

Given $a,b\in X_0^d$, denote by $X(a,b)$ the fiber at $(a,b)$ of the map
    \begin{equation*}
         X_1\xrw{(\pt_0,\pt_1)} X_0\times X_0 \xrw{\zg\times \zg}  X^d_0\times X^d_0\;.
    \end{equation*}
    The object $X(a,b)\in\tawg{n-1}$ should be thought of as hom-$(n-1)$-category. \index{Hom-$(n-1)$-category}
    We define a map $f:X\rw Y$ in $\tawg{n}$ to be an $n$-equivalence if \index{n-equivalences}
    \begin{itemize}
      \item [i)] For all $a,b\in X_0^d$
      \begin{equation*}
        f(a,b): X(a,b) \rw Y(fa,fb)
      \end{equation*}
      is an $(n-1)$-equivalence.\mk

      \item [ii)] $\p{n}f$ is an $(n-1)$-equivalence.
    \end{itemize}
    This completes the inductive step in the definition of $\tawg{n}$.

\end{definition}
\begin{remark}\label{rem-wg-ps-cat}
    It follows by Definition \ref{def-n-equiv} that $\catwg{n}\subset\tawg{n}$.
\end{remark}
\begin{example}\label{ex-tam}
Tamsamani $n$-categories. \index{Tamsamani $n$-category}
\mk

A special case of a weakly globular Tamsamani $n$-category occurs when $X\in\tawg{n}$ is such that $X_0$ and $X_{\oset{r}{1...1}0}$ are discrete for all $1 \leq r \leq n-2$. The resulting category is the category $\Tan$ of Tamsamani's $n$-categories. Note that, if $X\in\Tan$ then $X_s \in \ta{n-1}$ for all $n$, the induced Segal maps $\hmu{s}$ coincide with the Segal maps
\begin{equation*}
  \nu_s:X_s\rw \pro{X_1}{X_0}{s}
\end{equation*}
and $\p{n}X\in\ta{n-1}$. That is, we have a \emph{truncation functor} \index{Truncation functor}
\begin{equation*}
  \p{n}: \ta{n}\rw \ta{n-1}.
\end{equation*}
 A map in $\ta{n}$ is an $n$-equivalence it is so in $\tawg{n}$. \index{n-equivalences}
 Hence this recovers the original definition of Tamsamani's weak $n$-category \cite{Ta}. \index{Weakly globular!- Tamsamani $n$-category}

 In the case $n=2$, the relation between Tamsamani 2-categories \index{Tamsamani 2-categories} and bicategories was shown by Lack and the author \cite{LackPaoli2008} who introduced a 2-nerve functor from the 2-category of bicategories, normal homomorphism and icons to $\ta{2}$.

\end{example}
\begin{example}\label{ex-wg-2-3}
    Weakly globular Tamsamani 2-categories. \index{Weakly globular!- Tamsamani $2$-category}
    \mk

    \nid From the definition, $X\in \tawg{2}$ consists of a simplicial object $X\in\funcat{}{\Cat}$ such that $X_0\in\cathd{}$ and the induced Segal maps
    \begin{equation*}
        \hmuk:X_k \xrw{\muk}\pro{X_1}{X_0}{k}\xrw{\nu_k}\pro{X_1}{X^d_0}{k}
    \end{equation*}
    are equivalences of categories.

    The existence of the functor $\p{2}:\tawg{2}\rw \Cat$ making diagram \eqref{eq-wg-ps-cat} commute is in this case automatic so this condition does need to be imposed as part of the definition. In fact, since $p:\Cat\rw \Set$ sends equivalences to isomorphisms and preserves pullbacks over discrete objects, we obtain isomorphisms for each $k\geq 2$
    \begin{flalign*}
        p(\hmuk): p X_k & \cong p(\pro{X_1}{X_0^d}{k})\cong &\\
        & \cong \pro{p X_1}{p X_{0}^{d}}{k}\cong \pro{p X_1}{(p X_{0})^{d}}{k}\;.&
    \end{flalign*}
    Since $p X_k =(\bar p X)_k$ for all $k$,  the Segal maps for the simplicial set $\bar p X$ are isomorphisms. Therefore $\bar p X$ is the nerve of a category. We therefore have a functor $\p{2}$ given as composite
    \begin{equation*}
        \p{2}:\tawg{2}\xrw{\bar p} \mbox{\sf{Ner}} \Cat \subset \funcat{}{\Set}\xrw{P}\Cat
    \end{equation*}
    where $\mbox{\sf{Ner}} \Cat$ is the full subcategory of $\funcat{}{\Set}$ consisting of nerves of categories and $P$ is the left adjoint to the nerve functor $N:\Cat\rw\funcat{}{\Set}$.
    \end{example}

\subsection{Properties of weakly globular Tamsamani $\pmb{n}$-categories}
\vspace{10mm}

In this section we establish some important properties of the category $\tawg{n}$. Lemma \ref{lem-x-in-tawg-x-in-catwg} and Corollary \ref{cor-crit-ncat-be-wg} give criteria for objects of $\tawg{n}$ to be in $\catwg{n}$. This will be used in Section \ref{sec-func-qn} to establish the existence for the functor $\q{n}$, as well as in Chapters \ref{chap7} and \ref{chap8}, in particular in the proof of Theorem \ref{the-repl-obj-1} which will lead to approximating objects of $\tawg{n}$ with objects of $\lta{n}$ and thus to the rigidification functor.

In Proposition \ref{pro-n-equiv} we prove useful properties of $n$-equivalences in $\tawg{n}$, which will be used throughout the rest of this work. Later on, we will refine part of Proposition \ref{pro-n-equiv} by proving in Corollary \ref{cor-2-out-of} that $n$-equivalences in $\tawg{n}$ have the 2-out-of-3 property.

In Proposition \ref{pro-crit-lev-nequiv} we give a sufficient condition for an $\nequ$ $f$ in $\tawg{n}$ to be such that $J_n f$ is a levelwise equivalence of categories. This will be used in Section \ref{sec-func-qn} in the proof of Lemma \ref{lem-jn-alpha} leading to Theorem \ref{the-repl-obj-1}.

\begin{lemma}\label{lem-wg-ps-cat-b}
  For each $n>1$ there is a functor
  \begin{equation*}
    \di{n}:\tawg{n-1}\rw\tawg{n}
  \end{equation*}
  making the following diagram commute
  \begin{equation}\label{eq-lem-wg-ps-cat-b}
  \xymatrix@C=30pt{
    \tawg{n-1} \ar^{\Nb{n-1}}[rr]\ar_{\di{n}}[d] && \funcat{n-1}{\Set} \ar^{\ovl{d}}[d]\\
    \tawg{n}  \ar^{J_n}[rr] && \funcat{n-1}{\Cat}
    }
  \end{equation}
  where $d: \Set \rw \Cat$ is the discrete category functor. Further, $\di{n}$ sends $(n-1)$-equivalences to $n$-equivalences.
\end{lemma}
\begin{proof}
By induction on $n$. When $n=1$ let
\begin{equation*}
  \di{2}:\Cat\rw\tawg{2}
\end{equation*}
be given, for each $X\in\Cat$ and $k\geq 0$ by
\begin{equation*}
  (\di{2} X)_k=d X_k\;.
\end{equation*}
Clearly $\di{2}X\in\tawg{2}$ since $d X_0$ is discrete, so in particular homotopically discrete and the Segal maps (which coincide in this case with the induced Segal maps) are isomorphisms. Also,
\begin{equation*}
  (J_2 \di{2}X)_k=d X_k=(\ovl{d}NX)_k
\end{equation*}
so that $J_2 \di{2}=\ovl{d}N$.

 Let $f:X\rw Y$ be an equivalence of categories. Then $\di{2}f$ is a 2-equivalence since for each $x,y\in X_0$
\begin{equation*}
  (\di{2}X)(a,b)=d X(a,b)\cong dY(fa,fb)\cong (\di{2}Y)(fa,fb)
\end{equation*}
and $\p{2}\di{2}f=f$ is an equivalence of categories.

Suppose, inductively, that the lemma holds for $(n-1)$ and let $X\in\tawg{n-1}$. Then for each $\uk\in\dop{n-1}$,
\begin{equation*}
  (\ovl{d}\Nb{n-1}X)_{\uk}=d X_{\uk}\;.
\end{equation*}
Therefore, by induction hypothesis, for each $s>0$
\begin{equation*}
  (\ovl{d}\Nb{n-1}X)_{s*}=\di{n-1}X_s\in\tawg{n-1}
\end{equation*}
and $\di{n-1}X_0\in\cathd{n-1}$ since $X_0\in\cathd{n-2}$. The induced Segal maps of $\ovl{d}\Nb{n-1}X$, for each $s\geq 2$ are given by
\begin{equation*}
\begin{split}
    & (\ovl{d}\Nb{n-1}X)_{s}=\di{n-1}X_s \rw\pro{\di{n-1}X_1}{(\di{n-1}X_0)^d}{s}= \\
    & =\di{n-1}(\pro{X_1}{X_0^d}{s})\;.
\end{split}
\end{equation*}
This is a $(n-1)$-equivalence by induction hypothesis and the fact that
\begin{equation*}
  X_s\rw\pro{X_1}{X_0^d}{s}
\end{equation*}
is a $(n-2)$-equivalence. Thus $\ovl{d}\Nb{n-1}X\in J_n\tawg{n}$. If we define
 \begin{equation*}
   (\di{n}X)_{\uk}=dX_{\uk}
 \end{equation*}
 then \eqref{eq-lem-wg-ps-cat-b} commutes.

 Suppose $f:X \rw Y$ is a $(n-1)$-equivalence in $\tawg{n-1}$. Then for each $a,b\in X_0^d$, $X(a,b)\rw Y(fa, fb)$ is a $(n-1)$-equivalence, so by inductive hypothesis
  \begin{equation*}
    (\di{n}X)(a,b)=\di{n-1}X(a,b)\rw \di{n-1}Y(fa, fb)=(\di{n}Y)(fa,fb)
  \end{equation*}
  is a $(n-1)$-equivalence. Also, $\p{n}\di{n}f=f$ is a $(n-1)$-equivalence. We conclude that $\di{n}f$ is a $n$-equivalence.

\end{proof}
\begin{remark}\label{rem-wg-ps-cat-b}
It is immediate that $\di{n}:\tawg{n-1}\rw\tawg{n}$ restricts to functors
\begin{equation*}
  \di{n}: \catwg{n-1}\rw\catwg{n},\qquad \di{n}:\cathd{n-1}\rw\cathd{n},\qquad \di{n}:\ta{n-1}\rw\ta{n}\;.
\end{equation*}
\end{remark}
\begin{notation}\label{not-ex-tam}\index{Functor!- $\p{j,n}$}
    For each $1 \leq j\leq n-1$ denote
    \begin{equation*}
      \p{j,n}=\p{j}\p{j+1}...\p{n}:\tawg{n}\rw \tawg{j-1},\qquad \p{n,n}=\p{n}.
    \end{equation*}
\medskip
  Note that this restricts to functors
  \begin{equation*}
\begin{split}
    &  \p{j,n}=\p{j}\p{j+1}...\p{n}:\ta{n}\rw \ta{j-1}\\
    &  \p{j,n}=\p{j}\p{j+1}...\p{n}:\catwg{n}\rw \catwg{j-1}\\
    & \p{j,n}=\p{j}\p{j+1}...\p{n}:\cathd{n}\rw \cathd{j-1}\;.
\end{split}
  \end{equation*}

\end{notation}

\begin{notation}\label{dnot-ex-tam}\index{Functor!- $\di{n,j}$}
For each $1 \leq j\leq n-1$ denote\mk
\begin{equation*}
  \di{n,j}=\di{n}...\di{j+1}\di{j}: \tawg{j-1}\rw \tawg{n}
\end{equation*}
Note that $\di{n,j}$ restricts to functors
\begin{equation*}
  \di{n,j}: \catwg{j-1}\rw\catwg{n},\qquad \di{n,j}:\cathd{j-1}\rw\cathd{n},\qquad \di{n,j}:\ta{j-1}\rw\ta{n}\;.
\end{equation*}
\end{notation}

\bk
In the following lemma we give a criterion for a weakly globular Tamsamani $n$-category to be in $\catwg{n}$.
\begin{lemma}\label{lem-x-in-tawg-x-in-catwg}
  Let $X\in\tawg{n}$ be such that
      \begin{itemize}
      \item [a)] $X_s\in \catwg{n-1}$ for all $s$.\mk

      \item [b)] $X_s \cong \pro{X_1}{X_0}{s}$ for all $s\geq 2$.\mk

      \item [c)] For all $s\geq 2$ and $1\leq j\leq n-1$
          \begin{equation}\label{eq-lem-eq-def-wg}
          \begin{split}
              & \p{j,n-1}X_s \cong  \p{j,n-1}(\pro{X_1}{X_0}{s})= \\
             = \ &\pro{\p{j,n-1} X_1}{\p{j,n-1} X_0}{s}
          \end{split}
          \end{equation}
    \end{itemize}
    then $X\in\catwg{n}$.
\end{lemma}
\begin{proof}
By induction on $n$. When $n=2$, let $X\in\tawg{2}$ satisfy a), b), c). Then $X_0\in\cathd{}$ and by b), $X\in\cat{2}$. Since $X\in\tawg{2}$ the induced Segal maps
\begin{equation*}
  X_s \rw \pro{X_1}{X_0^d}{s}
\end{equation*}
are equivalences of categories for all $s\geq 2$. Thus, by definition, $X\in\catwg{2}$.

Suppose, inductively, that the lemma holds for $(n-1)$ and let $X\in\tawg{n}$ satisfy a), b), c). Then $X_0\in\cathd{n-1}$, $X_k\in\catwg{n-1}$ for all $k\geq 0$ and, since $X\in\tawg{n}$, the induced Segal maps
\begin{equation*}
  X_s \rw \pro{X_1}{X_0^d}{s}
\end{equation*}
are $(n-1)$-equivalences for all $s\geq 2$.

By hypothesis b) and the definition of $\catwg{n}$, to show that $X\in\catwg{n}$ it is enough to prove that $\p{n}X\in\catwg{n-1}$. We do so by proving that $\p{n-1}X$ satisfies the inductive hypothesis.

For all $s\geq 0$,
 \begin{equation*}
   (\p{n}X)_s=\p{n-1}X_s\in\catwg{n-2}
 \end{equation*}
 since, by hypothesis a), $X_s\in\catwg{n-1}$. Thus $\p{n}X$ satisfies the inductive hypothesis a). Also, from hypothesis c) for $X$ in the case $j=n-1$  for all $s\geq 2$,
\begin{equation*}
\begin{split}
    & \p{n-1}X_s \cong\p{n-1}(\pro{X_1}{X_0}{s})\cong \\
   \cong \ & \pro{\p{n-1}X_1}{\p{n-1}X_0}{s}
\end{split}
\end{equation*}
that is, $\p{n}X$ satisfies inductive hypothesis b). From this and from hypothesis c) for $X$, using the fact that $\p{j,n-1}=\p{j,n-2}\p{n-1}$ we deduce, for all $s\geq 2$ and $1\leq j\leq n-2$
\begin{align*}
    & \p{j,n-2}(\pro{(\p{n}X)_1}{(\p{n}X)_0}{s})=\\
    = \ & \p{j,n-2}(\pro{\p{n-1}X_1}{\p{n-1}X_0}{s})=\\
    = \ & \p{j,n-2}\p{n-1}(\pro{X_1}{X_0}{s})=\\
    = \ & \pro{\p{j,n-1}X_1}{\p{j,n-1}X_0}{s}=\\
    = \ & \pro{\p{j,n-2}(\p{n}X)_1}{\p{j,n-2}(\p{n}X)_0}{s}\;.
\end{align*}
This shows that $\p{n}X$ satisfies inductive hypothesis c). Hence we conclude that $\p{n}X\in\catwg{n-1}$, as required.

\end{proof}

We now deduce another useful criterion for a weakly globular Tamsamani $n$-category to be in $\catwg{n}$.

\begin{corollary}\label{cor-crit-ncat-be-wg}
    Let $X\in \tawg{n}$. Then $X\in\catwg{n}$ if and only if
    \begin{itemize}
      \item [a)]
      $X_s\in\catwg{n-1}$ for all $s\geq 0$.

      \item [b)] $(\Nu{2}X)_k = X_k\up{2}\in\catwg{n-1}$ for all $k\geq 1$.
    \end{itemize}
\end{corollary}
\begin{proof}
If $X\in\catwg{n}$ then a) and b) hold by Proposition \ref{pro-crit-ncat-be-wg}. Suppose conversely that $X\in\tawg{n}$ satisfies a) and b). We show that $X$ satisfies the hypotheses of Lemma \ref{lem-x-in-tawg-x-in-catwg} and therefore deduce that $X\in\catwg{n}$. By hypothesis b), for all $k\geq 0$ and $s\geq 2$
\begin{equation*}
    X_{sk}\cong \pro{X_{1k}}{X_{0k}}{s}\;.
\end{equation*}
Therefore $ X_{s}\cong \pro{X_{1}}{X_{0}}{s}$ for $s\geq 2$ so hypothesis b) in Lemma \ref{lem-x-in-tawg-x-in-catwg} holds. Also by hypothesis b) and by Lemma \ref{lem-prop-pn} we have, for all $k\geq 0$ and $s\geq 2$
\begin{equation*}
\begin{split}
   & \p{j,n-2}(\pro{X_{1k}}{X_{0k}}{s})\cong \\
    & \cong \pro{\p{j,n-2}X_{1k}}{\p{j,n-2}X_{0k}}{s}
\end{split}
\end{equation*}
and therefore
\begin{equation*}
    \p{j,n-1}(\pro{X_{1}}{X_{0}}{s})\cong \pro{\p{j,n-1}X_{1}}{\p{j,n-1}X_{0}}{s}
\end{equation*}
which is hypothesis c) in Lemma \ref{lem-x-in-tawg-x-in-catwg}. We conclude that $X\in\catwg{n}$, as required.
\end{proof}

\begin{lemma}\label{lem-flevel-fneq}
    Let $f:X\rw Y$ in $\tawg{n}$ be a levelwise $(n-1)$-equivalence in $\tawg{n-1}$. Then $f$ is an $n$-equivalence.
\end{lemma}
\begin{proof}
By induction on $n$. Let $n=2$. If $f_0$ is an equivalence of categories, $X_0^d\cong Y_0^d$. Hence
\begin{align}\label{eq1-lem-flevel-fneq}
    & Y_1=\uset{a',b'\in Y_0^d}{\coprod} Y(a',b') \cong \uset{fa,fb\in Y_0^d}{\coprod} Y(fa,fb)\;.
\end{align}
\begin{align}\label{eq2-lem-flevel-fneq}
    & X_1=\uset{a,b\in X_0^d}{\coprod}X(a,b)\;.
\end{align}

Since $f_1$ is an equivalence of categories it follows from \eqref{eq1-lem-flevel-fneq} and \eqref{eq2-lem-flevel-fneq} that $f(a,b)$ is an equivalence of categories. Further, $f_k$ is an equivalence of categories for all $k\geq 0$ so that $p f_k=(\p{2}f)_k$ is an isomorphism, hence $\p{2}f$ is an isomorphism; we conclude that $f$ is a 2-equivalence.

Suppose the lemma holds for $(n-1)$ and let $f$ be as in the hypothesis. Since $f_0$ is a $(n-1)$-equivalence in $\cathd{n-1}$, $X_0^d \cong Y_0^d$, so that \eqref{eq1-lem-flevel-fneq} holds. Since $f_1$ is a $(n-1)$-equivalence  it follows from \eqref{eq1-lem-flevel-fneq} and \eqref{eq2-lem-flevel-fneq} that $f(a,b)$ is a $(n-1)$-equivalence for all $a,b\in X_0^d$.

Since $f_k$ is an $(n-1)$-equivalence for all $k\geq 0$, $\p{n-1}f_k=(\p{n}f)_k$ is a $(n-2)$-equivalence. So $\p{n}f$ satisfies the induction hypothesis and is therefore a $(n-1)$-equivalence. In conclusion, $f$ is an $n$-equivalence.
\end{proof}
\begin{remark}\label{rem-local-equiv}
    Applying inductively Lemma \ref{lem-flevel-fneq} it follows immediately that if a morphism $f$ in $\tawg{n}$ is such that $J_n f$ is a levelwise equivalence of categories, then $f$ is an $n$-equivalence.
\end{remark}
\begin{definition}\label{def-local-equiv}
\index{Local $(n-1)$-equivalence}

    A morphism $f:X\rw Y$ in $\tawg{n}$ is said to be a local $\nm$-equivalence if
  for all $a,b \in X_0^d$,
  \begin{equation*}
    f(a,b):X(a,b)\rw Y(fa,fb)
  \end{equation*}
   is a $\nm$-equivalence in $\tawg{n-1}$.

\end{definition}
\nid In the following proposition, we describe some useful properties of $n$-equivalences in $\tawg{n}$. Using later results we will prove in Corollary \ref{cor-2-out-of} that $n$-equivalences in $\tawg{n}$ have the $2$-out-of-$3$ property.
\begin{proposition}\label{pro-n-equiv}\

\sk
\begin{itemize}
  \item [a)] Let $f$ be a morphism in $\tawg{n}$ which is an $\nequ$. Then $f$ is a local $(n-1)$-equivalence and $\p{1,n} f$ is an isomorphism.\smallskip
  \item [b)]  Let $f$ be a morphism  in $\tawg{n}$ which is a local $(n-1)$-equivalence and is such that $\p{1,n} f$ is surjective. Then $f$ is an $\nequ$.\smallskip
  \item [c)] Let $X\xrw{g} Z \xrw{h} Y$ be morphisms in $\tawg{n}$, $f= hg$ and suppose that $f$ and $h$ are $n$-equivalences. Then $g$ is an $n$-equivalence.\smallskip
  \item [d)]  Let $X\xrw{g} Z \xrw{h} Y$ be morphisms in $\tawg{n}$, $f= hg$ and suppose that $g$ and $h$ are $n$-equivalences. Then $f$ is an $n$-equivalence.\smallskip
  \item [e)]  Let $X\xrw{g} Z \xrw{h} Y$ be morphisms in $\tawg{n}$, $f= hg$ and let $g_0^d: X_0^d \rw Z_0^d$ be surjective; suppose that $f$ and $g$ are $n$-equivalences. Then $h$ is an $n$-equivalence.

\end{itemize}
\end{proposition}
\begin{proof}
By induction on $n$. It is clear for $n=1$. Suppose it is true for $n-1$.

a) Let $f:X\rw Y$ be an $\nequ$ in $\tawg{n}$. Then, by definition, $f$ is a local $(n-1)$-equivalence and $\p{n}f$ is a $\equ{n-1}$. Therefore, by induction hypothesis applied to $\p{n}f$, $\p{1,n} f$ is an isomorphism.

\mk
b) Suppose that $f:X\rw Y$ is a local $(n-1)$-equivalence in $\tawg{n}$ and $\p{1,n}f$ is surjective. To show that $f$ is a $\nequ$ we need to show that $\p{n}f$ is a $\equ{n-1}$. For each $a,b\in X_0^d$
\begin{equation*}
    (\p{n}f)(a,b)=\p{n-1}f(a,b)
\end{equation*}
Since $f(a,b)$ is a $(n-1)$-equivalence, $\p{n-1}f(a,b)$ is a $(n-2)$-equivalence; that is, $\p{n}f$ is a local $(n-2)$-equivalence.

Since $\p{1,n} f = \p{1,n-1}\p{n} f$ is surjective, by inductive hypothesis applied to $\p{n}f$ we conclude that $\p{n}f$ is a $\equ{n-1}$ as required.

\mk
c) For all $a,b\in X_0^d$,
\begin{equation}\label{eq1-pro-n-equiv}
    f(a,b)=h(ga,gb)g(a,b)
\end{equation}
with $f(a,b)$ and $h(ga,gb)$ \equ{n-1}s. By inductive hypothesis, $g(a,b)$ is therefore a $\equ{n-1}$.

By hypothesis and by part a), $\p{1,n} f$ and $\p{1,n} h$ are isomorphisms. Since
\begin{equation}\label{eq2-pro-n-equiv}
    \p{1,n} f=(\p{1,n} h)(\p{1,n} g)
\end{equation}
it follows that $\p{1,n} g$ is an isomorphism, hence in particular it is surjective. By part b), this implies that $g$ is a $\equ{n-1}$.
\mk

d) Suppose that $h$ and $g$ are $\nequ$s. By \eqref{eq1-pro-n-equiv}, $f$ is a local $(n-1)$-equivalence and by \eqref{eq2-pro-n-equiv} $\p{1,n} f$ is an isomorphism. By b), $f$ is thus a $\nequ$.
\mk

e) By hypothesis, for each $a',b'\in Z_0^d$, $a'=ga$, $b'=gb$ for $a,b\in X_0^d$. It follows that $h(a',b')=h(ga,gb)$. Since, by induction hypothesis and by \eqref{eq1-pro-n-equiv}, $h(ga,gb)$ is a $\nm$-equivalence, it follows that such is $h(a',b')$. That is, $h$ is a local equivalence.

By hypothesis and by part a), $\p{1,n} f$ and $\p{1,n} g$ are isomorphisms, so by \eqref{eq2-pro-n-equiv}, such is $\p{1,n} h$.

 We conclude by part b) that $h$ is a $\nequ$.
\end{proof}

\bk
In the following Lemma we consider pullbacks in $\funcat{n-1}{\Cat}$ of morphisms $J_n f$ where $f$ is a morphism in $\tawg{n}$, but we omit writing $J_n$ for ease of notation. This is justified since the functor $J_n:\tawg{n}\rw \funcat{n-1}{\Cat}$ is fully faithful.
\begin{lemma}\label{lem-crit-lev-nequiv}
    Consider the diagram in $\tawg{n}$
    \begin{equation*}
    \xymatrix{
    X \ar^{f}[r]\ar_{\za}[d] & Z \ar@{=}[d] & Y \ar_{g}[l]\ar^{\zb}[d]\\
    X' \ar_{f'}[r] & Z & Y' \ar^{g'}[l]
    }
    \end{equation*}
    with $Z$ discrete. Then:
    \begin{itemize}
      \item [a)] $X\tiund{Z} Y\,,\,X'\tiund{Z'} Y'\in\tawg{n}$, where the pullback are taken in $\funcat{n-1}{\Cat}$.\mk

      \item [b)] $\p{n}(X\tiund{Z}Y)\cong \p{n}X\tiund{\p{n}Z}\p{n}Y$.\mk

      \item [c)] If $\za,\,\zb$ are $n$-equivalences such is
      \begin{equation*}
        (\za,\zb):X\tiund{Z} Y\rw X'\tiund{Z'} Y'\;.
      \end{equation*}
    \end{itemize}
\end{lemma}
\begin{proof}
By induction on $n$. It is clear for $n=1$ since the maps $f,g,f',g'$ are isofibrations as their target is discrete. Suppose, inductively, that the lemma holds for $(n-1)$.

\mk
a) Since pullbacks in $\funcat{n-1}{\Cat}$ are computed pointwise, for each $k\geq 0$
\begin{equation*}
    (X\tiund{Z}Y)_k = X_k \tiund{Z} Y_k \in \funcat{n-2}{\Cat}
\end{equation*}
with $X_k,Y_k\in\tawg{n-1}$. It follows from inductive hypothesis a) that $(X\tiund{Z}Y)_k\in\tawg{n-1}$. Also,
 \begin{equation*}
   (X\tiund{Z}Y)_0 = X_0 \tiund{Z} Y_0 \in\cathd{n-1}
 \end{equation*}
  since $X_0,Y_0\in\cathd{n-1}$ and $Z$ is discrete, (see Lemma \ref{lem-copr-hom-disc}).

To show that $X\tiund{Z}Y\in\tawg{n}$ it remains to prove that the induced Segal maps $\hmu{k}$ for $X\tiund{Z}Y$ are $\nm$-equivalences. We prove this for $k=2$, the case $k>2$ being similar. Note that, by Lemma \ref{lem-copr-hom-disc}.
\begin{equation}\label{eq1-crit-lev-nequiv}
    \tens{(X\tiund{Z}Y)_1}{(X\tiund{Z}Y)_0^d} \cong (\tens{X_1}{X_0^d})\tiund{Z}(\tens{Y_1}{Y_0^d})\;.
\end{equation}
Consider the commutative diagram in $\tawg{n-1}$
\begin{equation}\label{eq2-crit-lev-nequiv}
\xymatrix{
X_2 \ar^{}[r]\ar_{\hmu{2}(X)}[d] & Z \ar@{=}[d] & Y_2 \ar^{}[l] \ar^{\hmu{2}(Y)}[d]\\
\tens{X_1}{X_0^d} \ar^{}[r] & Z & \tens{Y_1}{Y_0^d} \ar^{}[l]
}
\end{equation}
The vertical maps are the induced Segal maps for $X$ and $Y$, hence they are $\nm$-equivalences. By inductive hypothesis b) applied to \eqref{eq2-crit-lev-nequiv} we conclude that the induced map of pullbacks is a $\nm$-equivalence. By \eqref{eq1-crit-lev-nequiv} the latter is the induced Segal map $\hmu{2}$ for $X\tiund{Y}Z$. The proof for $k>2$ is similar and we conclude that $X\tiund{Y}Z\in\tawg{n}$.
\mk

b) By \eqref{eq-wg-ps-cat}, for all $\uk\in\dop{n-2}$
\begin{equation}\label{eq1-lem-crit-lev-nequiv}
    (\p{n}(X\tiund{Z}Y))_{\uk} = p (X\tiund{Z}Y)_{\uk}\;.
\end{equation}
Since pullbacks in $\funcat{n-1}{\Cat}$ are computed pointwise and $p$ commutes with pullbacks over discrete objects, we have
\begin{equation}\label{eq2-lem-crit-lev-nequiv}
\begin{split}
    & p (X\tiund{Z}Y)_{\uk} = p (X_{\uk}\tiund{Z_{\uk}}Y_{\uk})_{\uk} = pX_{\uk}\tiund{pZ_{\uk}}pY_{\uk}=\\
   = & (\p{n} X)_{\uk}\tiund{(\p{n}Z)_{\uk}} (\p{n} Y)_{\uk}\;.
\end{split}
\end{equation}
Since this holds for all $\uk$, \eqref{eq1-lem-crit-lev-nequiv} and \eqref{eq2-lem-crit-lev-nequiv} imply b).

\mk
c) For each $(a,b),\, (c,d)\in (X\tiund{Z}Y)^d_0 = X_0^d\tiund{Z}Y_0^d$ we have %
\begin{equation*}
\begin{split}
    & (X\tiund{Z}Y)((a,b),(c,d))=X(a,c)\times Y(b,d) \\
    & (X'\tiund{Z}Y')((fa,fb),(gc,gd))=X'(fa,gc)\times Y'(fb,gd)\;.
\end{split}
\end{equation*}
Since $\za,\zb$ are $n$-equivalences, $\za(a,c)$ and $\zb(b,d)$ are $\nm$-equivalences, hence such is
\begin{equation*}
    (\za,\zb)((a,b),(c,d))=\za(a,c)\times \zb(b,d)\;.
\end{equation*}
Since $\p{n}$ commutes with pullbacks over discrete objects for each $n$, so does $\p{1,n}$, hence
\begin{equation*}
    \p{1,n}(X\tiund{Z}Y) = \p{1,n}X\tiund{Z}\p{1,n}Y\;.
\end{equation*}
From the hypothesis and from Proposition \ref{pro-n-equiv} a), $\p{1,n}\za$ and $\p{1,n}\zb$ are isomorphisms, thus so is $\p{1,n}(\za,\zb)$. Since, from above, $(\za,\zb)$ is a local $\nm$-equivalence, we conclude by Proposition \ref{pro-n-equiv} b) that $(\za,\zb)$ is a $n$-equivalence.
\end{proof}

The following is a useful criterion for an $n$-equivalence in $\tawg{n}$ to be a levelwise equivalence of categories. It will be used in the rest of this chapter and in Chapter \ref{chap7}.

\begin{proposition}\label{pro-crit-lev-nequiv}
    Let $f:X\rw Y$ be a morphism in $\tawg{n}$ with $n\geq 2$, such that
    \begin{itemize}
      \item [a)] $f$ is a $\nequ$.
      \item [b)] $\p{n-1}X_0 \cong \p{n-1}Y_0$,
      \item [c)] For each $1\leq r < n-1$ and all $k_1, \ldots, k_r \geq 0$,
      \begin{equation*}
        \p{n-r-1}X_{k_1, \ldots, k_r,\, 0}\cong \p{n-r-1}Y_{k_1, \ldots, k_r,\, 0}\;.
      \end{equation*}
      Then $J_nf$ is a levelwise equivalence of categories.
    \end{itemize}
\end{proposition}
\begin{proof}
By induction on $n$. Let $f:X\rw Y$ be a 2-equivalence in $\tawg{2}$ such that
 \begin{equation*}
   X_0^d=p X_0\cong p Y_0=Y_0^d.
 \end{equation*}
 Since $f(a,b)$ is an equivalence of categories for all $a,b\in X_0^d$ we deduce that there is an equivalence of categories
\begin{equation*}
    f_1:X_1=\underset{a,b\in X_0^d}{\cop}X(a,b)\rw Y_1= \underset{a',b'\in Y_0^d}{\cop}Y(a',b')=\underset{fa,fb\in Y_0^d}{\cop}Y(fa,fb).
\end{equation*}

Hence there are equivalences of categories for $k\geq 2$
\begin{equation*}
    X_k \sim \pro{X_1}{X_0^d}{k}\sim \pro{Y_1}{Y_0^d}{k}\sim Y_k\;.
\end{equation*}
In conclusion $X_k\sim Y_k$ for all $k\geq 0$.

Suppose, inductively, that the statement holds for $(n-1)$ and let $f:X\rw Y$ be as in the hypothesis. We show that $J_{n-1}f_k$ is a levelwise equivalence of categories for each $k\geq 0$ by showing that $f_k$ satisfies the inductive hypothesis. It then follows that $J_n f$ is a levelwise equivalence of categories since
\begin{equation*}
    (J_n f)_{k_1...k_{n-1}}=(J_{n-1}f_{k_1})_{k_2...k_{n-1}}\;.
\end{equation*}

Since  $X_0\in\cathd{n-1}$, from b) and Lemma \ref{lem-neq-hom-disc} we obtain
\begin{equation}\label{eq1-pro-crit-lev-nequiv}
    X_0^d=\p{1,n}X_0 \cong \p{1,n}Y_0 \cong Y_0^d\;.
\end{equation}
Thus, by Lemma \ref{lem-neq-hom-disc} again, $f_0:X_0\rw Y_0$ is a $\equ{n-1}$. Further, by hypothesis c),
\begin{align*}
    & \p{n-2}X_{00} \cong \p{n-2}Y_{00}\\
    & \p{n-r-2}X_{0\,k_1...k_r\,0} \cong \p{n-r-2}Y_{0\,k_1...k_r\,0}
\end{align*}
for each $1\leq r < n-2$ and all $k_1...k_r$. Thus $f_0:X_0\rw Y_0$ satisfies the inductive hypothesis and we conclude that $f_0$ is a levelwise equivalence of categories. By \eqref{eq1-pro-crit-lev-nequiv} we also have
\begin{equation*}
    f_1:\underset{a,b\in X_0^d}{\cop}X(a,b)\rw Y_1 = \underset{a',b'\in Y_0^d}{\cop}Y(a',b')=\underset{fa,fb\in Y_0^d}{\cop}Y(fa,fb)\;.
\end{equation*}
Since $f$ is a local $\equ{n-1}$, it follows that $f_1:X_1 \rw Y_1$ is a $\equ{n-1}$. Further, by hypothesis c)
\begin{align*}
    & \p{n-2}X_{10} \cong \p{n-2}Y_{10}\\
    & \p{n-r-2}X_{1\,k_1...k_r\,0} \cong \p{n-r-2}Y_{1\,k_1...k_r\,0}
\end{align*}
for all $1 \leq r < n-2$. Thus $f_1$ satisfies the inductive hypothesis, and is therefore a levelwise equivalence of categories.

For each $k\geq 2$ consider the map
\begin{equation*}
    (f_1,...,f_1):\pro{X_1}{X_0^d}{k}\rw \pro{Y_1}{Y_0^d}{k}\;.
\end{equation*}
Since $X_0^d \cong Y_0^d$ and, from above, $f_1$ is a $\equ{n-1}$, then $(f_1,...,f_1)$ is also a $\equ{n-1}$.

There is a commutative diagram in $\tawg{n-1}$
\begin{equation*}
\xymatrix{
X_k \ar^(0.3){\hmuk}[rr] \ar_{f_k}[d] && \pro{X_1}{X_0^d}{k} \ar^{(f_1,...,f_1)}[d]\\
Y_k \ar_(0.3){\hmuk}[rr] && \pro{Y_1}{Y_0^d}{k}
}
\end{equation*}
where the horizontal induced Segal maps are $(n-1)$-equivalences since $X,Y\in\tawg{n-1}$ and the right vertical map is a $(n-1)$-equivalence from above. It follows from Proposition \ref{pro-n-equiv} c) and d) that $f_k$ is a $(n-1)$-equivalence. Further, from hypothesis c),
\begin{align*}
    & \p{n-1}X_{k0} \cong \p{n-2}Y_{k0}\\
    & \p{n-r-2}X_{k\,k_1...k_r\,0} \cong \p{n-r-2}Y_{k\,k_1...k_r\,0}\;.
\end{align*}
Thus $f_k$ satisfies the induction hypothesis and we conclude that $f_k$ is a levelwise equivalence of categories.

In conclusion, $f_k$ is a levelwise equivalence of categories for all $k\geq 0$. Since this holds for each $k\geq0$ this implies that $f$ is a levelwise equivalence of categories.
\end{proof}
%%
%% TO BE ACTIVATED WHEN SECTION 5 IS COMPLETED
%\begin{corollary}\label{cor-pro-crit-lev-nequiv}
%    Let $f:X\rw Y$ be a $\nequ$ in $\tawg{n}$ such that $X_0 \cong Y_0$, $X_{\seq{k}{1}{r}0} \cong Y_{\seq{k}{1}{r}0}$ for each $1\leq r \leq n-1$ and all $\seqc{k}{1}{r}$. Then $X \cong Y$.
%\end{corollary}
%%
%\begin{proof}
%By hypothesis, in particular the hypotheses of Proposition \ref{pro-crit-lev-nequiv} are satisfied and $f$ is a levelwise equivalence of categories. That is, for each $\seqc{k}{1}{n-1}$ $f_{\seq{k}{1}{n-1}}$ is an equivalence of categories. On the other hand, by hypothesis, $f_{\seq{k}{1}{n-1}}$ is a bijection on objects, thus is an isomorphism. Thus $f$ is an isomorphism.
%\end{proof}

%%%%%%%%%%%%%%%%%%%%%%%%%%%%%%%%%%%%%%%%%%%%%%%%%%%%%%%%%%%%%%%%%%%%%%%%%
\section{The functor $\pmb{\q{n}}$.}\label{sec-func-qn} \index{Functor! - $\q{n}$}
This section introduces the functor
\begin{equation*}
  \q{n}:\tawg{n}\rw \tawg{n-1}
\end{equation*}

This functor is a higher dimensional generalization of the connected component functor $q:\Cat\rw\Set$ and comes equipped with a morphism
 \begin{equation*}
   \zg\up{n}:X\rw \di{n}\q{n}X
 \end{equation*}
   natural in $X\in \tawg{n}$, where $\di{n}$ is as in Lemma \ref{lem-wg-ps-cat-b} . It will be used crucially in Section \ref{sec-cat-lta} to replace a weakly globular \nfol category $X$ with a simpler one (Theorem \ref{the-repl-obj-1}).
\begin{proposition}\label{pro-post-trunc-fun}
    There is a functor
     \begin{equation*}
     \qn:\tawg{n}\rw \tawg{n-1}
     \end{equation*}
      making the following diagram commute
    \begin{equation}\label{eq-lem-post-trunc}
    \xymatrix@C=30pt{
    \tawg{n} \ar^{J_{n}}[rr]\ar_{\qn}[d] && \funcat{n}{\Cat} \ar^{\ovl{q}}[d]\\
    \tawg{n-1} \ar_{\Nb{n-1}}[rr] & & \funcat{n-1}{\Set}
    }
    \end{equation}
    where $q:\Cat\rw \Set$ is the connected component functor. The functor $\qn$ sends $n$-equivalences to $\nm$-equivalences and preserves pullbacks over discrete objects. If $X\in\cathd{n}$, then $\q{n}X=\p{n}X$; further, for each $X\in \tawg{n}$, there is a map
     \begin{equation*}
     \zgu{n}:X\rw \dn\qn X
     \end{equation*}
      natural in $X$.
\end{proposition}
\begin{proof}
By induction on $n$; for $n=1$, $\q{1}=q:\Cat\rw\Set$ is the connected components functor which, by Lemma \ref{lem-q-pres-fib-pro}, has the desired properties. The map $\zgu{1}$ is the unit of the adjunction $q\dashv \di{1}$. If $X\in\cathd{}$, in particular $X$ is a groupoid, so $pX=qX$.

Suppose we defined $\q{n-1}$ with the desired properties and let $X\in\tawg{n}$. We claim that
 \begin{equation*}
   \ovl{\q{n-1}}X\in\tawg{n-1}.
 \end{equation*}
 In fact, for each $s\geq 0$, by induction hypothesis
\begin{equation*}
    ( \ovl{\q{n-1}}X)_s=\q{n-1}X_s\in\tawg{n-2}\;.
\end{equation*}
Also, by induction hypothesis and by Definition \ref{def-hom-dis-ncat},
\begin{equation*}
    ( \ovl{\q{n-1}}X)_0 =\q{n-1}X_0 = \p{n-1}X_0 \in \cathd{n-1}
\end{equation*}
as $X_0 \in \cathd{n-1}$. Further, since
\begin{equation*}
    \hmuk : X_k\rw \pro{X_1}{X_0^d}{k}
\end{equation*}
is a $\equ{n-1}$, by induction hypothesis the map
\begin{equation*}
    \q{n-1}X_k\rw\q{n-1}(\pro{X_1}{X_0^d}{k})\cong \pro{\q{n-1}X_1}{X_0^d}{k}
\end{equation*}
is a $\equ{n-2}$, where we used the fact that $\q{n-1}X_0^d\cong X_0^d$, which follows from diagram \eqref{eq-lem-post-trunc} at step $(n-1)$.

This shows that $ \ovl{\q{n-1}}X\in\tawg{n-1}$. We therefore define
\begin{equation*}
    \q{n}X = \ovl{\q{n-1}}X\;.
\end{equation*}
The fact that $\q{n}$ satisfies diagram \eqref{eq-lem-post-trunc} is immediate from the definitions and the induction hypothesis.

If $X\in\cathd{n}$, by definition $X_k\in\cathd{n-1}$ for each $k$, so by induction hypothesis $\p{n-1}X_k=\q{n-1}X_k$. It follows that $(\p{n}X)_k=(\q{n}X)_k$ for all $k$. That is $\p{n}X=\q{n}X$.

Let $f:X\rw Y$ be a $\nequ$ in $\tawg{n}$ and let $a,b\in X_0^d$. Then from the definitions
\begin{equation*}
    (\q{n}f)(a,b)=\q{n-1}f(a,b)\;.
\end{equation*}
Since $f(a,b)$ is a $\equ{n-1}$, by induction hypothesis $\q{n-1} f(a,b)$ is a $\equ{n-2}$. By Proposition \ref{pro-n-equiv}, to prove that $\q{n}f$ is a $\equ{n-1}$, it is enough to show that $\p{1,n-1}\q{n}f$ is surjective.

Recall that for any category $\clC$ there is a surjective map $p\,\clC\rw q\,\clC$ natural in $\clC$. Applying this levelwise to $J_nX$ we obtain a map
\begin{equation*}
    \za_X\up{n}:\p{n}X \rw \q{n}X
\end{equation*}
natural in $X$. The map $\za\up{n}$ induces  a functor
\begin{equation}\label{eq-post-trunc-fun}
 \p{2,n-1}\za_X\up{n}: \p{2,n-1}\p{n}X \rw \p{2,n-1}\q{n}X\;,
\end{equation}
which is identity on objects. In fact, on objects this map is given by
\begin{equation*}
   \p{1,n-2}\p{n-1}X_0 \rw \p{1,n-2}\q{n-1}X_0
\end{equation*}
and since $X_0\in\cathd{n-1}$, $\p{n-1}X_0=\q{n-1}X_0$ so this map is the identity. It follows that the map in $\Set$
\begin{equation*}
    \p{1,n-1}\za_X\up{n-1}: \p{1,n-1}\p{n}X \rw \p{1,n-1}\q{n}X
\end{equation*}
is surjective. We thus have a commuting diagram
\begin{equation*}
\xymatrix{
\p{1,n-1}\p{n}X \ar^{\p{1,n-1}\p{n}f}[rr] \ar_{\p{1,n-1}\za_X\up{n}}[d] && \p{1,n-1}\p{n}Y \ar^{\p{1,n-1}\za_Y\up{n}}[d]\\
\p{1,n-1}\q{n}X \ar^{\p{1,n-1}\q{n}f}[rr] && \p{1,n-1}\q{n}Y
}
\end{equation*}
in which the top arrow is an isomorphism (by Proposition \ref{pro-n-equiv}) and from above the vertical arrows are surjective. It follows that the bottom map is also surjective. By Proposition \ref{pro-n-equiv} b) we conclude that $\q{n}f$ is a $\equ{n-1}$.

Finally, the map $\zgu{n}:X\rw \di{n}\q{n}X$ is given levelwise by the maps $X_s\rw \di{n-1}\q{n-1}X_s$, which exist by induction hypothesis.
\end{proof}
\begin{remark}\label{rem-spec-isofib}
    For each $X\in \tawg{n}$, from the proof of Proposition \ref{pro-post-trunc-fun} the functor $\p{2,n-1}\za_X\up{n}$ is identity on objects. It is also surjective on morphisms since, by the proof of Proposition \ref{pro-post-trunc-fun} (taking $X_1\in\tawg{n-1}$ instead of $X\in \tawg{n}$), the map
    \begin{equation*}
        (\p{2,n-1}\za_X\up{n})_1:\p{1,n-2}\p{n-1} X_1 \rw \p{1,n-2}\q{n-1} X_1 \;,
    \end{equation*}
    is surjective. It follows that $\p{2,n-1}\za_X\up{n}$ is an isofibration.
\end{remark}

\begin{corollary}\label{pro-post-wg-ncat} \index{Functor! - $\q{n}$}
    The functor
    \begin{equation*}
         \qn:\tawg{n}\rw \tawg{n-1}\;,
    \end{equation*}
    restricts to functors
    \begin{equation*}
        \qn:\catwg{n}\rw \catwg{n-1}\qquad \qn:\cathd{n}\rw \cathd{n-1}\qquad \qn:\ta{n}\rw \ta{n-1}
    \end{equation*}
\end{corollary}
\begin{proof}

By induction on $n$. When $n=2$, $\q{2}X\in\Cat$. Suppose the statement holds for $(n-1)$ and let $X\in\catwg{n}$. To prove that $\q{n}X\in\catwg{n-1}$ we show that it satisfies the hypotheses of Corollary \ref{cor-crit-ncat-be-wg}.

For each $s\geq 0$, $(\q{n}X)_s=\q{n-1}X_s\in \catwg{n-2}$, by inductive hypothesis applied to $X_s$. Thus condition a) in Corollary \ref{cor-crit-ncat-be-wg} holds. Also
\begin{equation*}
    (\q{n}X)_k\up{2}=\q{n-2}X_k\up{2}\in\catwg{n-2}
\end{equation*}
since $X_k\up{2}\in\catwg{n-1}$ (as $X\in \catwg{n}$ using Proposition \ref{pro-crit-ncat-be-wg}). Thus condition b) in Corollary \ref{cor-crit-ncat-be-wg} is satisfied and we conclude that $\q{n}X\in\catwg{n-1}$. In conclusion
\begin{equation*}
  \qn:\catwg{n}\rw \catwg{n-1}.
\end{equation*}
If $X\in \cathd{n}$, $\q{n}X=\p{n}X \in \cathd{n-1}$ by definition of $\cathd{n}$.

Let $X\in\ta{n}$. We show by induction that $\q{n}X \in \ta{n-1}$ and that $\q{n}X$ is discrete if $X$ is discrete. This is clear for $n=2$, since $\q{2}X \in \Cat$.

 Inductively, if $X\in\ta{n}$, $X_{k-1}\in \ta{n-1}$ so by inductive hypothesis
  \begin{equation*}
    (\q{n}X)_k=\q{n-1}X_k\in \ta{n-1}
  \end{equation*}
   with $(\q{n}X)_0=\q{n-1}X_{0}$ discrete since $X_0$ is discrete (using the inductive hypothesis). It follows that $(\q{n}X)_0$ and $(\q{n}X)_{\oset{r}{1...1}0}$ are discrete for all $1 \leq r \leq n-2$, and thus by definition (see Example \ref{ex-tam}) $\q{n}X\in\ta{n-1}$.
It is straightforward that if $X$ is discrete such is $\q{n}X$.

\end{proof}

\begin{notation}\label{not-ex-tam-q}\index{Functor!- $\q{j,n}$}
    For each $1 \leq j\leq n-1$ denote
    \begin{equation*}
      \q{j,n}=\q{j}\q{j+1}...\q{n}:\tawg{n}\rw \tawg{j-1}, \qquad \q{n,n}=\q{n}:\tawg{n}\rw \tawg{n-1}.
    \end{equation*}
    Note that this restricts to functors
  \begin{equation*}
\begin{split}
    &  \q{j,n}=\q{j}\q{j+1}...\q{n}:\ta{n}\rw \ta{j-1}\\
    &  \q{j,n}=\q{j}\q{j+1}...\q{n}:\catwg{n}\rw \catwg{j-1}\\
    & \q{j,n}=\q{j}\q{j+1}...\q{n}:\cathd{n}\rw \cathd{j-1}\;.
\end{split}
  \end{equation*}
\end{notation}

\begin{remark}\label{rem-qrn}\index{Functor!- $\q{j,n}$}
  It is immediate from Corollary \ref{pro-post-wg-ncat} that the functor $\q{j,n}$ restricts to functors
  \begin{equation*}
        \q{j,n}:\catwg{n}\rw \catwg{j-1}\qquad \q{j,n}:\cathd{n}\rw \cathd{j-1} \qquad \q{j,n}:\ta{n}\rw \ta{j-1}
    \end{equation*}
\end{remark}

\begin{lemma}\label{lem-prop-qn}
    For each $X\in\catwg{n}$, $1\leq j < n$ and $s\geq 2$ it is
    \begin{equation}\label{eq-lem-prop-qn}
    \begin{split}
        &\q{j,n-1}X_s \cong  \q{j,n-1}(\pro{X_1}{X_0}{s})=\\
         & =\pro{\q{j,n-1} X_1}{\q{j,n-1} X_0}{s}\;.
    \end{split}
    \end{equation}
\end{lemma}
\begin{proof}
Since $X\in\catwg{n}$ by Corollary \ref{pro-post-wg-ncat} $\q{n}X\in\catwg{n-1}$, hence
\begin{equation*}
    \q{n-1}(\pro{X_1}{X_0}{s})= \pro{\q{n-1} X_1}{\q{n-1} X_0}{s}
\end{equation*}
which is \eqref{eq-lem-prop-qn} for $j=n-1$. Since $\q{j+1,n}X \in\catwg{j}$ for $1\leq j\leq (n-1)$,  its Segal maps are isomorphisms. Further for all $s\geq 0$
\begin{equation*}
    (\q{j+1,n}X)_s=(\q{j+1}...\q{n}X)_s = \q{j}...\q{n-1}X_s = \q{j,n-1} X_s
\end{equation*}
with $X_s=\pro{X_1}{X_0}{s}$ for $s\geq 2$. Hence \eqref{eq-lem-prop-qn} follows for $1\leq j < n$.
\end{proof}

\section{Pullback constructions using $\pmb{\qn}$}\label{sec-pull-qn}
In Section \ref{sec-cat-lta}, Theorem \ref{the-repl-obj-1}, we will replace a weakly globular \nfol category $X$ with a simpler one, and this will play a crucial role in the construction of the rigidification functor $Q_n$. The proof of Theorem \ref{the-repl-obj-1} will involve taking pullbacks along the map $\zg\up{n}$ of Proposition \ref{pro-post-trunc-fun}. In this section we establish several properties of these pullbacks which will be needed later on.

We will always consider $\tawg{n}$ (as well as $\catwg{n}$) as embedded in $\funcat{n-1}{\Cat}$ via the functor $J_n$ and our pullbacks will be taken in $\funcat{n-1}{\Cat}$ (and thus they are levelwise pullback in $\Cat$). To ease the notation, we omit writing explicitly $J_n$ in these pullbacks, which is justified since $J_n$ is fully faithful.
\begin{lemma}\label{lem-spec-pulbk-eqr}
    Let $X\in\cathd{n}$, $Z\in\cathd{n-1}$, $r: Z\rw \qn X$. Consider the pullback in $\funcat{n-1}{\Cat}$
    \begin{equation*}
        \xymatrix@R=35pt @C=40pt{
        P \ar^{}[r] \ar^{}[d] & X \ar^{\zgu{n}_X}[d] \\
        \dn Z \ar_{\dn r}[r] & \dn \qn X
        }
    \end{equation*}
    then $P\in\cathd{n}$ and $\p{n}P=Z$.
\end{lemma}
\begin{proof}
By induction on $n$. For $n=1$, since $\di{1}\q{1}X$ is discrete, the map
 \begin{equation*}
   \zgu{1}:X\rw \di{1}\q{1}X=X^d
 \end{equation*}
  is an isofibration. Therefore, since $\zgu{1}$ is an equivalence of categories (as $X\in\cathd{}$) we have an equivalence of categories
\begin{equation*}
    P=\di{1}Z\tiund{\di{1}\q{1}X} X\simeq \di{1}Z\tiund{\di{1}\q{1}X} \di{1}\q{1}X = \di{1}Z
\end{equation*}
Thus $P\in\cathd{}$ and $pP=Z$.

Suppose, inductively, that the lemma holds for $n-1$ and let $P$ be as in the hypothesis. Since pullbacks in $\funcat{n-1}{\Cat}$ are computed pointwise, for each $k\geq 0$ we have a pullback in $\funcat{n-2}{\Cat}$
\begin{equation*}
\xymatrix{
P_k \ar^{}[rr] \ar^{}[d] && X_k \ar^{}[d]\\
\di{n-1}Z_k \ar^{}[rr] && \di{n-1}\q{n-1}X_k
}
\end{equation*}
where $X_k\in\cathd{n-1}$ (since $X\in \cathd{n}$) and $Z_k \in \cathd{n-2}$ (since $Z\in\cathd{n-1}$). By induction hypothesis, we conclude that $P_k\in\cathd{n-1}$.

We now show that, for each $k \geq 2$
\begin{equation}\label{eq1-lem-spec-pulbk-eqr}
    P_k\cong \pro{P_1}{P_0}{k}\;.
\end{equation}
We illustrate this for $k=2$, the case $k>2$ being similar. Since $X\in\cathd{n}$, $\q{n}X = \p{n}X\in\cathd{n-1}$, so
\begin{equation*}
\begin{split}
    & \q{n-1}X_2=\p{n-1}X_2 =\p{n-1}(\tens{X_1}{X_0})= \\
   =\ & \tens{\p{n-1} X_1}{\p{n-1} X_0}=\tens{\q{n-1} X_1}{\q{n-1} X_0}\;.
\end{split}
\end{equation*}
Since $X_2 \cong \tens{X_1}{X_0}$ and $Z_2 \cong \tens{Z_1}{Z_0}$, it follows from Lemma \ref{lem-char-obj-II} that
 \begin{equation*}
   P_2 \cong \tens{P_1}{P_0}.
 \end{equation*}

To prove that $P\in\cathd{n}$ it remains to show that $\p{n}P\in\cathd{n-1}$. Since $p$ commutes with fiber products over discrete objects, for each $\us\in\dop{n-1}$ we have
\begin{align*}
    &(\p{n}P)_{\us}=p\,P_{\us}=
    p(d\,Z_{\us}\tiund{d q X_{\us}} X_{\us})
     = Z_{\us}\tiund{qX_{\us}}p X_{\us}=Z_{\us}
\end{align*}
where we used the fact that, since $X_{\us}$ is a groupoid, $pX_{\us}=qX_{\us}$. Since this holds for each ${\us}$ we conclude that $\p{n}P=Z\in\cathd{n-1}$ as required.
\end{proof}
\begin{lemma}\label{lem-p2-n-1}
    Let $Y\in\tawg{n}$ and let
    \begin{equation}\label{eq1-lem-p2-n-1}
        X\rw \q{n}Y \lw \p{n}Y
    \end{equation}
    be a diagram in $\tawg{n-1}$ such that $X\tiund{\q{n}Y}\p{n}Y\in \tawg{n-1}$. Then for all $1\leq j \leq n-1$
      \begin{equation*}
        \p{j,n-1}(X\tiund{\q{n}Y}\p{n}Y)=\p{j,n-1}X\tiund{\p{j,n-1}\q{n}Y}\p{j,n-1}\p{n}Y\;.
      \end{equation*}
\end{lemma}
\begin{proof}\

By induction on $n$. For $n=2$, the functor $\p{2}Y\rw \q{2}Y$ is the identity on objects, therefore by Lemma \ref{lem-iso-cla-obj-fib-pro}
\begin{equation*}
    p(X\tiund{\q{2}Y}\p{2}Y)=pX\tiund{p\q{2}Y}p\p{2}Y\;.
\end{equation*}
Suppose, inductively, that the lemma holds for $n-1$. Then for each $k\geq 0$
\begin{equation}\label{eq2-lem-p2-n-1}
\begin{split}
   & (\p{j,n-1}(X\tiund{\q{n}Y}\p{n}Y))_k= \\
   =\; &\p{j-1,n-2}(X_k\tiund{\q{n-1}Y_k}\p{n-1}Y_k)=\\
   =\; &\p{j-1,n-2}X_k \tiund{\p{j-1,n-2}\q{n-1}Y_k} \p{j-1,n-2}\p{n-1}Y_k=\\
   =\; &(\p{j,n-2}X)_k \tiund{(\p{j,n-2}\q{n-1}Y)_k} (\p{j,n-2}\p{n-1}Y)_k\;.
\end{split}
\end{equation}
Since this holds for each $k\geq 0$, the lemma  follows.
\end{proof}
\begin{proposition}\label{pro-spec-plbk-pscatwg}
    Let
     \begin{equation*}
       \di{n}A\xrw{\ \di{n}f\ }\di{n}\q{n}C \xlw{\zgu{n}} C
     \end{equation*}
      be a diagram in $\tawg{n}$ where $f: A\rw\q{n}C$ is a morphism in $\tawg{n-1}$ and consider the pullback in $\funcat{n-1}{\Cat}$

      \begin{equation*}
        \xymatrix@R=35pt @C=40pt{
        P \ar^{w}[r] \ar^{}[d] & C \ar^{\zgu{n}_C}[d] \\
        \dn A \ar_{\dn f}[r] & \dn \qn C
        }
    \end{equation*}
    \begin{itemize}
    \item [a)] Then $P\in \tawg{n}$ and
    \begin{equation}\label{eq-pro-spec-plbk-pscatwg0}
        \p{1,n}P = \p{1,n-1} A \tiund{\p{1,n-1}\q{n}C} \p{1,n} C\;.
    \end{equation}
    \item [b)]

    Consider the commutative diagram in $\tawg{n}$
    \begin{equation}\label{eq-pro-spec-plbk-pscatwg}
        \xymatrix@R=35pt @C=60pt{
        \di{n}A \ar^{\di{n}f}[r] \ar_{a}[d] & \di{n}\q{n} C \ar^{b}[d] & C \ar^{c}[d] \ar_{\zgu{n}_C}[l] \\
        \di{n}D \ar_{\di{n}h}[r] & \di{n}\q{n} F  & F \ar^{\zgu{n}_F}[l]
        }
    \end{equation}
    where $a,b,c$ are $\nequ$s. Then the induced maps of pullbacks
    \begin{equation*}
    (a,c): \di{n}A\tiund{\di{n}\q{n} C} C \rw \di{n}D \tiund{\di{n}\q{n} F} F
    \end{equation*}
    is a $n$-equivalence in $\tawg{n}$.

    \medskip
    \item [c)] If $f$ is an $(n-1)$-equivalence, $P\xrw{w}C$ is an $n$-equivalence.

    \end{itemize}
\end{proposition}
\begin{proof}
By induction on $n$. When $n=1$, the maps $f,\zgu{n}_C,h,\zgu{n}_F$ are isofibrations since their targets is a discrete category. Therefore, since $a,b,c$ are equivalences of categories, the induced map of pullbacks
\begin{equation*}
    (a,c):d A\tiund{d q C} C \rw d D \tiund{d q F} F
\end{equation*}
is an equivalence of categories and
\begin{equation*}
  p (d A\tiund{dqC}C) = A \tiund{qC} pC.
\end{equation*}
 Suppose, inductively, that the proposition holds for $n-1$.

\bk
a) We have
\begin{equation*}
   P_0=(\di{n}A\tiund{\di{n}\q{n} C} C)_0=\di{n-1}A_0\tiund{\di{n-1}\q{n-1} C_0} C_0\;.
\end{equation*}
Since, by hypothesis, $A \in\tawg{n-1}$ and $C\in\tawg{n}$, by definition $A_0 \in\cathd{n-2}$ and $C_0 \in \cathd{n-1}$. Therefore by Lemma \ref{lem-spec-pulbk-eqr}
\begin{equation*}
    \di{n-1}A_0\tiund{\di{n-1}\q{n-1} C_0} C_0\in \cathd{n-1}\;.
\end{equation*}
Further, for each $k\geq 1$,
\begin{equation*}
    P_k=(\di{n}A\tiund{\di{n}\q{n} C} C)_k = \di{n-1}A_k\tiund{\di{n-1}\q{n-1} C_k} C_k
\end{equation*}
where, by hypothesis, $A_k \in\tawg{n-2}$ and $C_k \in \tawg{n-1}$. It follows by inductive hypothesis that
\begin{equation*}
    \di{n-1}A_k\tiund{\di{n-1}\q{n-1} C_k} C_k \in \tawg{n-1}\;.
\end{equation*}
To prove that $\di{n}A\tiund{\di{n}\q{n} C} C \in \tawg{n}$, it remains to show that its induced Segal maps $\hmu{s}$ are $\equ{n-1}$s for all $s\geq 2$. We show this for $s=2$, the case $s>2$ being similar. We have
\begin{equation*}
\begin{split}
    & \hmu{2}:(\di{n}A\tiund{\di{n}\q{n} C} C)_2 \rw \\
    & \rw \tens{(\di{n}A\tiund{\di{n}\q{n} C} C)_1}{(\di{n}A\tiund{\di{n}\q{n} C} C)_0^d}\;.
\end{split}
\end{equation*}
By Lemma \ref{lem-spec-pulbk-eqr},
\begin{equation*}
\begin{split}
    & (\di{n}A\tiund{\di{n}\q{n} C} C)_0^d= (\p{n-1}(\di{n-1}A_0\tiund{\di{n-1}\q{n-1} C_0} C_0))^d=A_0^d = \\
    & = A_0^d\tiund{\di{n-1}\q{n-1} C_0^d} C_0^d
\end{split}
\end{equation*}
where we used the fact that, since $\q{n-1}\di{n-1}=\Id$,
\begin{equation*}
\begin{split}
    & \di{n-1}\q{n-1}C_0^d= \di{n-1}\q{n-1}\di{n-1}...\di{1}p...\p{n-1}C_0= \\
   =\; &\di{n-1}\di{n-2}...\di{1}p...\p{n-1}C_0=C_0^d\;.
\end{split}
\end{equation*}
Recalling that $\q{n-1}$ commutes with pullbacks over discrete objects, we obtain
\begin{equation*}
\begin{split}
    & \tens{(\di{n}A\tiund{\di{n}\q{n}C}C)_1}{(\di{n}A\tiund{\di{n}\q{n}C}C)_0^d}= \\
    & \resizebox{1.0\hsize}{!}{$ =\tens{(\di{n-1}A_1\tiund{\di{n-1}\q{n-1}C_1}C_1)}{(A_0^d\tiund{\di{n-1}\q{n-1}C_0^d}C_0^d)}=$} \\
   & = \di{n-1}(\tens{A_1}{A_0^d})\tiund{\di{n-1}\q{n-1}(\tens{C_1}{C_0^d})}(\tens{C_1}{C_0^d})\;.
\end{split}
\end{equation*}
On the other hand
\begin{equation*}
    (\di{n}A\tiund{\di{n}\q{n}C} C)_2 = \di{n-1}A_2 \tiund{\di{n-1}\q{n-1}C_2} C_2 \;.
\end{equation*}
Hence we see that the map $\hmu{2}$ is the induced map on pullbacks from the diagram in $\tawg{n-1}$
\begin{equation}\label{eq1-pro-spec-plbk-pscatwg}
\xymatrix@C=15pt{
\di{n-1}A_2 \ar^{}[rr] \ar^{}[d] && \di{n-1}\q{n-1}C_2 \ar^{}[d] && C_2 \ar^{}[ll] \ar^{}[d]\\
\di{n-1}(\tens{A_1}{A_0^d}) \ar^{}[rr] && \di{n-1}\q{n-1}(\tens{C_1}{C_0^d}) && \tens{C_1}{C_0^d} \ar^{}[ll]
}
\end{equation}
In this diagram, the left and right vertical maps are $\nm$-equivalences since they are induced Segal maps of $\di{n}A$ and $C$ respectively. The map
\begin{equation*}
    \q{n-1}C_2 \rw \q{n-1}(\tens{C_1}{C_0^d}) = \tens{\q{n-1}C_1}{(\q{n-1}C_0)^d}
\end{equation*}
is the induced Segal map for $\q{n}C\in\tawg{n-1}$, and is therefore a $(n-2)$-equivalence. The central vertical map in \eqref{eq1-pro-spec-plbk-pscatwg} is therefore in particular a $(n-1)$-equivalence.  We can therefore apply the induction hypothesis b) to the diagram \eqref{eq1-pro-spec-plbk-pscatwg} and conclude that the induced map on pullbacks is a $\equ{n-1}$. That is, $\hmu{2}$ is a $\nm$-equivalence.

Similarly one shows that $\hmu{s}$ is a $\nm$-equivalence for all $s \geq 2$. This concludes the proof that
 \begin{equation*}
   \di{n}A \tiund{\di{n}\q{n}C} C \in \tawg{n}
 \end{equation*}
 and \eqref{eq-pro-spec-plbk-pscatwg0} follows from Lemma \ref{lem-p2-n-1}.

\bk
b) We now show that $(a,c)$ is a $\nequ$. By Proposition \ref{pro-n-equiv} b), it is enough to show that it is a local $\equ{n-1}$ and that $\p{1,n}(a,c)$ is an isomorphism.

Let
\begin{equation*}
    (x_1,y_1),(x_2,y_2) \in (\di{n}A \tiund{\di{n}\q{n}C} C)_0^d = A_0^d \tiund{\di{n-1}\q{n-1}C_0^d} C_0^d.
\end{equation*}
Then the map
\begin{equation*}
    (a,c)((x_1,y_1),(x_2,y_2))
\end{equation*}
is the induced map on pullbacks in the diagram
\begin{equation}\label{eq2-pro-spec-plbk-pscatwg}
\xymatrix@C=14pt{
\di{n}A(x_1,x_2) \ar^(0.38){f(x_1,x_2)}[rr] \ar_{a(x_1,x_2)}[d] && \di{n-1}\q{n-1}C(fx_1,fx_2) \ar^{b(fx_1,fx_2)}[d] && C(y_1,y_2) \ar_(0.38){\zgu{n}_C(y_1,y_2)}[ll] \ar_{c(y_1,y_2)}[d]\\
\di{n-1}D(a x_1,a x_2) \ar_(0.45){h(a x_1,a x_2)}[rr] && \di{n-1}\q{n-1}F(h a x_1,h a x_2) && F(cy_1,cy_2) \ar^(0.35){\zgu{n}_F(cy_1,cy_2)}[ll]
}
\end{equation}
The vertical maps in \eqref{eq2-pro-spec-plbk-pscatwg} are $\equ{n-1}$s. By the inductive hypothesis b), since $a,b,c$ are $n$-equivalences we conclude that the induced map of pullbacks is a $\equ{n-1}$. This shows that $(a,c)$ is a local $\equ{n-1}$. By Lemma \ref{lem-p2-n-1},
\begin{equation}\label{pac}
  \p{1,n}(a,c)=(\p{1,n}a,\p{1,n}c).
\end{equation}

Applying the functor $\p{1,n}$ to the diagram  \eqref{eq-pro-spec-plbk-pscatwg} we obtain a commutative diagram in $\Set$
\begin{equation*}
\xymatrix{
\p{1,n-1}A \ar^{}[rr] \ar_{\p{1,n}a}[d] && \p{1,n-1}\q{n}C \ar^{\p{1,n}b}[d] && \p{1,n-1}\p{n}C \ar^{}[ll] \ar^{\p{1,n}c}[d]\\
\p{1,n-1}D \ar^{}[rr] && \p{1,n-1}\q{n}F && \p{1,n-1}\p{n}F \ar^{}[ll]
}
\end{equation*}
Since, by hypothesis, $a,b$ and $c$ are $\nequ$s, by Proposition \ref{pro-n-equiv} a) the vertical maps are isomorphisms, hence by \eqref{pac} such is $\p{1,n}(a,c)$. We conclude from Proposition \ref{pro-n-equiv} b) that $(a,c)$ is an $n$-equivalence.

\bk
c) By a), we have a commutative diagram in $\tawg{n}$
\begin{equation*}
    \xymatrix{
    \di{n}A \ar^{f}[rr] \ar_{\di{n}}[d] && \di{n}\q{n}C \ar@{=}[d] && C \ar_{\zgu{n}_C}[ll] \ar@{=}[d] \\
    \di{n}\q{n}C \ar@{=}[rr] && \di{n}\q{n}C  && C \ar^{\zgu{n}_C}[ll]
    }
\end{equation*}
in which the vertical maps are $n$-equivalences. It follows by b) that the induced map of pullbacks
\begin{equation*}
    P= \di{n}A\tiund{\di{n}\q{n}}C \xrw{w} \di{n}\q{n} C \tiund{\di{n}\q{n}} C = C
\end{equation*}
is an $n$-equivalence.

\end{proof}
%%
%%

%%%%%%%%%%%%%%%%%%%%%%%%%%%%%%%%%%%%%%%%%%%%%%%%%%%%%%%%%%%%%%%%%%%%%%%%%%%%%
\chapter{Rigidifying weakly globular Tamsamani $\pmb{n}$-categories}\label{chap7}
In this chapter we continue the study of the category $\tawg{n}$ of weakly globular Tamsamani $n$-categories introduced in Chapter \ref{chap6}. The main result of this chapter, Theorem \ref{the-funct-Qn}, asserts the existence of a \emph{rigidification functor}
\begin{equation*}
  Q_n:\tawg{n}\rw\catwg{n}
\end{equation*}
such that for each $X\in\tawg{n}$ there is an $n$-equivalence
\begin{equation*}
 Q_n X\rw X.
\end{equation*}
  This result means that $X$ can be approximated up to $n$-equivalence with the more rigid and therefore simpler structure $Q_n X$. In particular, this implies (see Corollary \ref{cor-the-funct-Qn}) that the two categories $\catwg{n}$ and $\tawg{n}$ are equivalent after localization with respect to the $n$-equivalences.

The functor $Q_n$ restricts in particular to a functor
\begin{equation*}
  Q_n:\ta{n}\rw\catwg{n}
\end{equation*}
from Tamsamani $n$-categories to weakly globular \nfol categories.
In Chapter \ref{chap9} we will show that this leads to an equivalence after localization between $\catwg{n}$ and $\ta{n}$, exhibiting $\catwg{n}$ as a new model of weak $n$-categories satisfying, in particular, the homotopy hypothesis.

The rigidification functor factors through the category
 \begin{equation*}
   \segpsc{n-1}{\Cat}
 \end{equation*}
 of Segalic pseudo-functors introduced in Chapter \ref{chap5}. More precisely, $Q_n$ is the composite
\begin{equation*}
    Q_n:\tawg{n}\xrw{\ \;\;} \segpsc{n-1}{\Cat}\xrw{\St} \catwg{n}\subset \funcat{n-1}{\Cat}.
\end{equation*}
In the case $n=2$, it is easy to build pseudo-functors from $\tawg{2}$, and was already done by Pronk and the author in \cite{PP}. More precisely, given $X\in\tawg{2}$, define $\tr{2}X\in[ob(\Dop),\Cat]$ by
\begin{equation}\label{eq1-the-XXXX}
(\tr{2}X)_k=
\left\{
  \begin{array}{ll}
    X_0^d & k=0 \\
    X_1 & k=1 \\
    \pro{X_1}{X_0}{k} & k>1\;.
  \end{array}
\right.
\end{equation}
Since $X\in\tawg{2}$, $X_0\in\cathd{}$ so there are equivalences of categories
\begin{equation*}
\begin{split}
    & X_0\simeq X_0^d \\
    & X_k\simeq \pro{X_1}{X_0^d}{k}\quad \text{for }\; k>1.
\end{split}
\end{equation*}
Thus, for all $k\geq 0$ there is an equivalence of categories
\begin{equation*}
 (\tr{2}X)_k\simeq X_k\;.
\end{equation*}
By using transport of structure (more precisely Lemma \ref{lem-PP} with $\clC=\Dop$) we can lift $\tr{2}X$ to a pseudo-functor
\begin{equation*}
    \tr{2}X \in \psc{}{\Cat}
\end{equation*}
and by construction $\tr{2}X \in\segpsc{}{\Cat}$.

Building pseudo-functors from $\tawg{n}$ when $n>2$ is much more complex, and it is new to this work. The above approach cannot be applied directly because the induced Segal maps, when $n>2$ are $(n-1)$-equivalences but not in general levelwise equivalence of categories. For this reason we introduce an intermediate category $\lta{n}$, from which it is possible to build pseudo-functors using transport of structure. The functor from $\tawg{n}$ to Segalic pseudo-functors factorizes as
\begin{equation*}
    \tawg{n}\xrw{P_n} \lta{n} \xrw{\Tr_{n}} \segpsc{n-1}{\Cat}
\end{equation*}
The functor $P_n$ produces a functorial approximation (up to $n$-equivalence) of an object of $\tawg{n}$ with an object of $\lta{n}$, while $\Tr_{n}$ is built using transport of structure.

This chapter is organized as follows: In Section \ref{sec-cat-lta} we introduce the subcategory $\lta{n}\subset \tawg{n}$ and establish its properties. In particular we show in Section \ref{catnwg-ltawg} that $\catwg{n}\subset\lta{n}$ and we discuss in Section \ref{subs-geom-int} a corresponding geometric interpretation implied by this fact.

 In Section \ref{sec-approx} we show how to approximate up to $n$-equivalence objects of $\tawg{n}$ with objects of \lta{n}. More precisely, using a pullback construction and the properties established in Section \ref{sec-func-qn} we show in Theorem \ref{the-repl-obj-1} that if  $X\in\tawg{n}$ is such that $\qn X$ can be approximated up to $(n-1)$-equivalence with an object of $\catwg{(n-1)}$, then $X$ can be approximated up to an $n$-equivalence with an object of $\lta{n}$. In Section \ref{sec-wg-tam-to-psefun} this is used in the proof of Theorem \ref{the-funct-Qn} to construct the functor
 \begin{equation*}
    P_n:\tawg{n}\rw\lta{n}.
\end{equation*}

In Section \ref{sec-from-pseu} we show how to construct pseudo-functors from the category $\lta{n}$, proving in Theorem \ref{the-XXXX} the existence of the functor

\begin{equation*}
    Tr_n: \lta{n}\rw \segpsc{n-1}{\Cat}
\end{equation*}

In Theorem \ref{the-funct-Qn} we use the functors $P_n$ and $Tr_{n}$ to build the rigidification functor $Q_n$ as the composite
\begin{equation*}
    \tawg{n}\xrw{P_n} \lta{n} \xrw{\Tr_{n}} \segpsc{n-1}{\Cat}\xrw{\St} \catwg{n}.
\end{equation*}

\section{The category $\pmb{\lta{n}}$.}\label{sec-cat-lta}
In this section we introduce the subcategory $\lta{n}$ of $\tawg{n}$. We will show in Section \ref{sec-wg-tam-to-psefun} how to build a functor from this category to the category of Segalic pseudo-functors, which in turn will lead to the construction of the rigidification functor $Q_n$.

\subsection{The idea of the category $\pmb{\lta{n}}$}\label{subs-idea-ltan}
The idea of the category $\lta{n}$ is best conveyed by its characterization given in Proposition \ref{pro-maps-nu-eqcat}, while Definition \ref{def-ind-sub-ltawg} captures its inductive nature which is useful in the proof of Theorem \ref{the-repl-obj-1}.

 Given $X\in \tawg{n} \subset \funcat{n-1}{\Cat}$, if we fix all but one of the $(n-1)$ simplicial directions we obtain a simplicial object in $\Cat$ which at level zero is homotopically discrete, and we can therefore consider the corresponding induced Segal maps (see notation \ref{not-ind-seg-map})
  \begin{equation*}
    \nu(\uk,i):X_{\uk}\rw \pro{X_{\uk(1,i)}}{X^d_{\uk(0,i)}}{k_i}\;.
  \end{equation*}

  Objects $X \in \lta{n}$ are such that all these induced Segal maps  are equivalences of categories and, further, $\p{n}X\in\catwg{n-1}$. We show in Proposition \ref{pro-maps-nu-eqcat} that these two properties characterize objects of $\tawg{n}$ which are in $\lta{n}$.

This characterization is useful to gain an intuition about the category $\lta{n}$; further, we will see in Section \ref{sec-from-pseu} that it allows to apply transport of structure to objects of $\lta{n}$ and thus build Segalic pseudo-functors from them.

\subsection{The formal definition of the category $\pmb{\lta{n}}$}
\begin{notation}\label{not-for-lta}

Let $X\in\tawg{n}$, then $J_n X\in\funcat{n-1}{\Cat}$ and
\begin{equation*}
  Z=\xi_1 J_n X\in\funcat{}{\funcat{n-2}{\Cat}}
\end{equation*}
where $\xi_1$ is as in Lemma \ref{lem-multi-simpl-as}. Since $X_0\in\cathd{n-1}$ there is a map
\begin{equation*}
  X_0\rw\di{n-1}\p{n-1}X_0
\end{equation*}
and therefore a corresponding map in $\funcat{n-2}\Cat$
\begin{equation*}
  J_{n-1}X_0\rw J_{n-1}\di{n-1}\p{n-1}X_0=Y\;.
\end{equation*}
We also have for each $k\in\Dop$
\begin{equation*}
  Z_k=J_{n-1}X_k\;.
\end{equation*}
For each $k\geq 2$ we therefore obtain induced Segal maps (in the sense of Definition \ref{def-ind-seg-map}) for $Z\in\funcat{}{\funcat{n-2}{\Cat}}$
\begin{equation*}
  v_k:Z_k\rw\pro{Z_1}{Y}{k}
\end{equation*}
that is, from above
\begin{equation*}
  v_k: J_{n-1}X_k \rw \pro{J_{n-1}X_1}{J_{n-1}\di{n}\p{n}X_0}{k}\;.
\end{equation*}
\end{notation}
\begin{definition}\label{def-ind-sub-ltawg}
    Define inductively the subcategory
     \begin{equation*}
       \lta{n}\subset\tawg{n}.
     \end{equation*}
      For $n=2$, $\lta{2}=\tawg{2}$. Suppose, inductively, that we defined $\lta{n-1}\subset \tawg{n-1}$. Let $\lta{n}$ be the full subcategory of $\tawg{n}$ whose objects $X$ are such that
    \begin{itemize}
      \item [i)] $X_k\in\lta{n-1}$ for all $k\geq 0$.\bk

      \item [ii)] The maps in $\funcat{n-2}{\Cat}$ (see Notation \ref{not-for-lta})
    \begin{equation*}
        v_k: J_{n-1} X_k \rw \pro{J_{n-1}X_1}{J_{n-1}\di{n-1}\p{n-1}X_0}{k}
    \end{equation*}
    are levelwise equivalences of categories for all $k\geq 2$\bk

     \item [iii)] $\p{n}X\in\catwg{n-1}$.
    \end{itemize}
\end{definition}

\begin{remark}\label{rem-def-ltawg}
  Let $\us\in \dop{n-2} $. Since $(J_{n-1} X_k)_{\us}=X_{k \us}$ and
   \begin{equation*}
     ({J_{n-1}\di{n-1}\p{n-1}X_0})_{\us}=dpX_{0 \us}=X_{0 \us}^d
   \end{equation*}
   condition ii) in Definition \ref{def-ind-sub-ltawg} is equivalent to
\begin{equation*}
X_{k\us}\simeq \pro{X_{1\us}}{X_{0\us}^d}{k}\;.
\end{equation*}
\end{remark}

\subsection{Properties of the category $\pmb{\lta{n}}$}\label{proper-ltawg}
\begin{definition}\label{def-ltawg-equiv}
    Let $\lnta{n}{n}$ be the full subcategory of $\tawg{n}$ whose objects $X$ are such that
    \begin{itemize}
      \item [i)] The maps in $\funcat{n-2}{\Cat}$ (see Notation \ref{not-for-lta})
      \begin{equation*}
        v_k: J_{n-1}X_k \rw \pro{J_{n-1}X_1}{J_{n-1}\di{n-1}\p{n-1}X_0}{k}
      \end{equation*}
      are levelwise equivalences of categories for all $k\geq 2$\bk

      \item [ii)] $\p{n}X\in\catwg{n-1}$.
    \end{itemize}
\end{definition}

\bk

\nid  Note that, by definition of $\tawg{2}$ and $\lta{2}$
      \begin{equation*}
        \lnta{2}{2}=\lta{2}=\tawg{2}.
      \end{equation*}

\begin{lemma}\label{lem-lev-wg-pscat} %5.12
    Let $X\in\tawg{n}$. Then $X\in \lta{n}$ if and only if
    \begin{itemize}
      \item [a)] $X\in \lnta{n}{n}$.\bk

      \item [b)] For each $1 < r \leq n-1$ and each $k_1,...,k_{n-r}\in\Dop$,
      \begin{equation*}
        X_{k_1,...,k_{n-r}}\in \lnta{r}{r}
      \end{equation*}
    \end{itemize}
\end{lemma}
\begin{proof}
By induction on $n$. It holds for $n=2$ since $\lta{2}=\Lb{2}\tawg{2}=\tawg{2}$. Suppose it holds for $(n-1)$ and let $X\in\lta{n}$. Then by definition $X\in\Lb{n}\tawg{n}$. Also by definition $X_{k_1\cdots k_{n-r}}\in\lta{r}$ and therefore
 \begin{equation*}
   X_{k_1\cdots k_{n-r}}\in\Lb{r}\tawg{r}.
 \end{equation*}

Conversely, suppose that $X\in\tawg{n}$ satisfies a) and b). By a), $v_k$ is a levelwise equivalence of categories and $\p{n}X\in\catwg{n-1}$; by b), $X_k\in\Lb{n-1}\tawg{n-1}$ and, further, $X_k$ itself satisfies b). Thus by induction hypothesis applied to $X_k$ we conclude that $X_k\in\lta{n-1}$. By definition, this shows that $X\in\lta{n}$.
\end{proof}

Next we establish a property of the category $\lta{n}$ which will be needed in Section \ref{sec-from-pseu} to build from it  Segalic pseudo-functors. We first fix a notation for the induced Segal maps of simplicial objects in $\Cat$ obtained from objects of $\tawg{n}$ in which all except for one simplicial directions in $J_nX$ are fixed.
\begin{notation}\label{not-ind-seg-map}
    Let $X\in\tawg{n}$, $\uk=(k_1,\ldots,k_{n-1})\in\dop{n-1}$, $1\leq i\leq n-1$. Then there is $X^i_{\uk}\in\funcat{}{\Cat}$ with
    \begin{equation*}
        (X^i_{\uk})_r=X_{\uk(r,i)}=X_{k_1...k_{i-1}r k_{i+1}...k_{n-1}}
    \end{equation*}
    so that $(X^i_{\uk})_{k_i}=X_{\uk}$. Since $X_{k_1,\ldots,k_{i-1}}\in\tawg{n-i+1}$, $X_{k_1,\ldots,k_{i-1}0}\in\cathd{n-i}$ and thus by Lemma \ref{lem-pos-grou-hom-disc}
    \begin{equation*}
        X_{k_1...k_{i-1} 0 k_{i+1}...k_{n-1}}= X_{\uk(0,i)}\in\cathd\;.
    \end{equation*}
    We therefore obtain induced Segal maps in $\Cat$ for all $k_i\geq 2$.
    \begin{equation}\label{eq1-not-ind-seg-map}
        \nu(\uk,i):X_{\uk}\rw \pro{X_{\uk(1,i)}}{X^d_{\uk(0,i)}}{k_i}\;.
    \end{equation}
\end{notation}
\begin{proposition}\label{pro-maps-nu-eqcat}
    Let $X\in\tawg{n}$. Then $X\in\lta{n}$ if and only if the following conditions hold:

      \begin{itemize}
        \item[a)] For each $\uk\in\dop{n-1}$, $1\leq i\leq n-1$ and $k_i\geq 2$ the maps $ \nu(\uk,i)$ in \eqref{eq1-not-ind-seg-map} are equivalences of categories.\bk

        \item[b)] $\p{n}X\in\catwg{n-1}$.
      \end{itemize}
\end{proposition}
\begin{proof}
Let $X\in\lta{n}$. Then b) holds by Definition \ref{def-ltawg-equiv}. We prove that a) holds by induction on $n$.

Consider first the case $1<i\leq n-1$. Since $X\in\lta{n}$, by definition $X_{k_1}\in\lta{n-1}$. Denoting $\ur=(k_2...k_{n-1})$ we have
\begin{equation}\label{eq2-not-ind-seg-map}
\begin{split}
    & X_{\uk}=(X_{k_1})_{\ur} \\
    & X_{\uk(1,i)}=(X_{k_1})_{\ur(1,i-1)}\qquad X_{\uk(0,i)}=(X_{k_1})_{\ur(0,i-1)}\;.
\end{split}
\end{equation}
The induction hypothesis applied to $Y=X_{k_1}$ implies the equivalence of categories for each $r_{i-1}\geq 2$
\begin{equation*}
    Y_{\ur}\simeq \pro{Y_{\ur(1,i-1)}}{Y^d_{\ur(0,i-1)}}{r_{i-1}}\;.
\end{equation*}
By \eqref{eq2-not-ind-seg-map} this means that $\nu(\uk,i)$ is an equivalence of categories for all $1<i\leq n-1$.

\nid Consider the case $i=1$. By definition of $\lta{n}$ the map in $\funcat{n-2}{\Cat}$
\begin{equation*}
    J_{n-1}X_{k_1}\rw\pro{J_{n-1}X_1}{\bd\bp J_{n-1}X_0}{k_1}
\end{equation*}
is a levelwise equivalence of categories for each $k_1\geq 2$. Therefore, for each $\ur=(k_2...k_{n-1})$ there is an equivalence of categories.
\begin{equation*}
    (X_{k_1})_{\ur}\simeq \pro{(X_1)_{\ur}}{d p (X_0)_{\ur}}{k_1}\;.
\end{equation*}
Since $(X_1)_{\ur}=X_{\uk(1,1)}$ and $dp(X_0)_{\ur}=X^d_{\uk(0,1)}$ we obtain the equivalence of categories
\begin{equation*}
    X_{\uk}\simeq \pro{X_{\uk(1,1)}}{X^d_{\uk(0,1)}}{k_1}\;.
\end{equation*}
In conclusion $\nu(\uk,i)$ is an equivalence of categories for all $1\leq i\leq n-1$.

Conversely, let $X\in\tawg{n}$ satisfy a) and b). We show that $X\in\lta{n}$ by induction on $n$.

 When $n=2$, $\tawg{2}=\lta{2}$ so the lemma holds. Suppose it is true for $(n-1)$ and let $X\in\tawg{n}$ satisfy the hypothesis. We show that conditions i), ii), iii) in Definition \ref{def-ind-sub-ltawg} are satisfied and therefore $X\in\lta{n}$. Condition iii) holds by assumption. For each $k\geq 0$, $X_k$ satisfies hypotheses a) and b), since
\begin{equation*}
  \p{n-1}X_k=(\p{n}X)_k \in\catwg{n-2}
\end{equation*}
as $\p{n}X\in\catwg{n-1}$. Thus by induction hypothesis $X_k\in\lta{n-1}$, which is condition i).
As for conditions ii) note that, if $\us=(\seqc{k}{2}{n-1})\in\dop{n-2}$
\begin{equation*}
\begin{split}
   & X_{\uk}=(J_{n-1}X_{k_1})_{\us} \\
   & X_{\uk(1,1)}= (J_{n-1}X_1)_{\us}\\
   & X_{\uk(0,1)}= (J_{n-1}X_0)_{\us}\;.
\end{split}
\end{equation*}
Therefore the maps $\nu(\uk,1)$ are given by
\begin{equation*}
  (J_{n-1}X_{k_1})_{\us}\rw \pro{(J_{n-1}X_1)_{\us}}{(J_{n-1}X_0)_{\us}}{k_1}
\end{equation*}
and these are equivalence of categories by hypothesis. Since this holds for each $\us\in\dop{n-2}$, this means  that
\begin{equation*}
  J_{n-1}X_{k_1}\rw \pro{J_{n-1}X_1}{J_{n-1}X_0}{k_1}
\end{equation*}
is a levelwise equivalence of categories for all $k_1\geq 2$, which is condition ii) in Definition \ref{def-ind-sub-ltawg}. We conclude that $X\in\lta{n}$.

\end{proof}
\subsection{$\pmb{\catwg{n}}$ and the category $\pmb{\lta{n}}$}\label{catnwg-ltawg}

In this section we show that $\catwg{n}$ is a full subcategory of $\lta{n}$.

\begin{proposition}\label{cor-lev-wg-ncat}
    Let $X\in\catwg{n}$. Then $X\in \lta{n}$.
\end{proposition}
\begin{proof}
We claim that $X\in\lnta{n}{n}$. We prove the claim by induction on $n$. The proposition holds for $n=2$ because $\catwg{2}\subset \tawg{2}=\lnta{2}{2}$. Suppose, inductively, that the statement holds for $n-1$ and let $X\in\catwg{n}$. Then, by Proposition \ref{pro-crit-ncat-be-wg}, $ X_k\up{2}\in \catwg{n-1}$ for all $k\geq 0$, where we denote $X_k\up{2}= (N_{(2)}X)_k$. Hence by inductive hypothesis
 \begin{equation*}
   X_k\up{2}\in\lnta{n-1}{n-1}.
 \end{equation*}
 From the definition of $X_k\up{2}$ this means that, for all $s\geq 2$ and $k\geq 0$, the map
\begin{equation*}
   J_{n-2} X_{sk}\rw \pro{J_{n-2}X_{1k}}{J_{n-2}\di{n-2}\p{n-2}X_{0k}}{s}
\end{equation*}
is a levelwise equivalence of categories. As this holds for every $k\geq 0$, we conclude that the map
\begin{equation*}
   J_{n-1} X_s \rw \pro{J_{n-1}X_1}{J_{n-1}\di{n-1}\p{n-1}X_0}{s}
\end{equation*}
is a levelwise equivalence of categories. Further, by hypothesis $\p{n}X\in\catwg{n-1}$ and thus we conclude that $X\in\lnta{n}{n}$, as claimed.

Since $X\in\catwg{n}$, by definition $X_{\seq{k}{1}{n-r}}\in \catwg{r}$ for each $1 < r \leq n-1$ and each $k_1,...,k_{n-r}\in\Dop$, thus from above
  \begin{equation*}
    X_{\seq{k}{1}{n-r}}\in \lnta{r}{r}.
  \end{equation*}
   Hence, by Lemma \ref{lem-lev-wg-pscat}, $X\in\lta{n}$.
\end{proof}
\begin{corollary}\label{rem-lev-wg-ncat}
Let $X\in\catwg{n}$, $1\leq j\leq n-1$, $k\geq 2$. Denote
\begin{equation*}
  Y=\pro{X_1}{\di{n-1,j}\p{j,n-1}X_0}{k}\in \cat{n-1}.
\end{equation*}
Then
\begin{itemize}
  \item [a)] $Y\in\catwg{n-1}$\mk

  \item [b)] If $X_t\in\cathd{n-1}$ for all $t\geq 0$, then $Y\in\cathd{n-1}$.
\end{itemize}
\end{corollary}
\begin{proof}
We prove this for $k=2$, the case $k>2$ being similar. When $j=1$, $\di{n-1,j}\p{j,n-1}X_0=X^d_0$ so a) holds since $X\in\catwg{n}$. As for b), since $X_1\in\cathd{n-1}$ then $Y\in\cathd{n-1}$ by Lemma \ref{lem-copr-hom-disc}.

When $j=n-1$,
\begin{equation*}
  Y=X_1\tiund{\di{n-1,j}\p{j,n-1}X_0}X_1\;.
\end{equation*}
By Proposition \ref{cor-lev-wg-ncat} the map $f: X_2\rw Y$ in $\catwg{n-1}$ is such that $J_{n-1}f$ is a levelwise equivalence of categories and $X_2\in\catwg{n-1}$. Thus by Corollary \ref{rem-last-ins}, $Y\in\catwg{n-1}$, proving a). As for b), if $X_2\in\cathd{n-1}$ by Corollary \ref{rem-last-ins} also $Y\in\cathd{n-1}$.

We now proceed to prove the corollary by induction on $n$. For $n=2$ it holds since $Y=X_1\tiund{X_0^d}X_1$. Suppose, inductively, that it holds for each $r\leq (n-1)$. From above, we can assume that $2\leq j\leq n-2$.
\medskip

a) We verify that $Y\in\cat{n-1}$ satisfies the hypotheses of Proposition \ref{pro-crit-ncat-be-wg} b). In fact, by \ref{pro-crit-ncat-be-wg} a), $(\Nu{2}X)_0\in\catwg{n-1}$, with
\begin{equation*}
  (\Nu{2}X)_{0t}=X_{t0}\in\cathd{n-2}
\end{equation*}
for all $t\geq 0$. Thus $(\Nu{2}X)_0$ satisfies the inductive hypothesis b) and we conclude by induction that
\begin{equation*}
  Y_0= \tens{X_{10}}{\di{n-2,j-1}\p{j-1,n-2}X_{00}}\in\cathd{n-2}\;.
\end{equation*}
Similarly, for all $s\geq 0$, the inductive hypothesis a) applied to $(\Nu{2}X)_{s}\in\catwg{n-1}$ yields
\begin{equation}\label{eq1-rem-lev-wg-ncat}
    Y_s= \tens{X_{1s}}{\di{n-2,j-1}\p{j-1,n-2}X_{0s}}\in\catwg{n-2}
\end{equation}
which implies $Y_{s0}\in\cathd{n-3}$. Thus $Y$ satisfies the hypothesis i) in Proposition \ref{pro-crit-ncat-be-wg} b).

By inductive hypothesis a) applied to $\p{n}X\in\catwg{n-1}$, using the fact that $j\leq n-2$, we obtain
\begin{equation}\label{eq2-rem-lev-wg-ncat}
  \bp Y=\tens{(\p{n}X)_1}{\di{n-2,j}\p{j,n-2}(\p{n}X)_0}\in\Nb{n-2}\catwg{n-2}\;.
\end{equation}
This shows that $Y$ satisfies hypothesis ii) of Proposition \ref{pro-crit-ncat-be-wg} b). We conclude from Proposition \ref{pro-crit-ncat-be-wg} that $Y\in \catwg{n-1}$.\bk

b) We verify that $Y$ satisfies the hypotheses of Corollary \ref{cor-crit-nequiv-rel}. From a), $Y\in\catwg{n-1}$. Further, for all $t\geq 0$
\begin{equation*}
  (\Nu{2}X)_{1t}=X_{t1}\in \cathd{n-2}
\end{equation*}
since, by hypothesis, $X_t\in\cathd{n-1}$. Thus $(\Nu{2}X)_{1}\in\catwg{n-1}$ satisfies the inductive hypothesis b) and we conclude from \eqref{eq1-rem-lev-wg-ncat} and the induction that $Y_1\in\cathd{n-2}$. Also, for each $t\geq 0$
\begin{equation*}
  (\p{n}X)_t=\p{n-1}X_t\in\cathd{n-2}
\end{equation*}
since $X_t\in\cathd{n-1}$. Thus $\p{n}X\in\catwg{n-1}$ satisfies the inductive hypothesis b), so that
$\p{n}Y\in\cathd{n-2}$.

Hence $Y\in\catwg{n-1}$ satisfies all the hypotheses of Corollary \ref{cor-crit-nequiv-rel} and we conclude that $Y\in\cathd{n-1}$.
\end{proof}
%%%%%%%%%%%%%%%%%%%%%%%%%%%%%%%%%%%%%%%%%%%%%%%%%%%%%%%%%%%%%%%%%%%%%%

\subsection{A geometric interpretation}\label{subs-geom-int}
The fact, proved in the previous section, that $\catwg{n}\subset \lta{n}$ implies the equivalence of categories $\nu(\uk,i)$ of Proposition \ref{pro-maps-nu-eqcat} for objects of $\catwg{n}$. We discuss in this section a geometric interpretation implied by this fact. Although this geometric interpretation is interesting in its own right, it is not strictly needed for the rest of this work, so this section can be skipped at first reading.

 We first give some preliminary definitions, which hold in any $n$-fold category.
\begin{definition}\label{def-1-proper-ltawg}
  Let $\us \in\dop{n}$. We call $\us$ an orientation \index{Orientation}if $s_j\in\{0,1\}$ for all $j=1,\ldots,n$ and we denote $t(\us)=\sum^n_{j=1}s_j$.
\end{definition}
\begin{definition}\label{def-2-proper-ltawg}
  Let $X\in\cat{n}$ and let $\us\in\dop{n}$ be an orientation. A $(n,t(\us))$-hypercube \index{Hypercube}in $X$ with orientation $\us$ is an element of the set
  \begin{equation*}
    (\Nb{n}X)_{\us}=X_{\us}
  \end{equation*}
  where $\Nb{n}$ is the multinerve $\Nb{n}:\cat{n}\rw \funcat{n}{\Set}$, and we use the shorthand notation $X_{\us}$ as adopted elsewhere.
\end{definition}
\begin{definition}\label{def-3-proper-ltawg}
  Let $0\leq t\leq n$. We denote by $\Cube{n,t}$ the set of all $(n,t)$-hypercubes \index{Set of $(n,t)$-hypercubes}in $X$ with orientations $\us$ such that $t(\us)=t$, that is
  \begin{equation*}
  \Cube{n,t}=\bigcup_{\substack{\us\in \dop{n}\text{ orientation} \\ \sum^n_{j=1}s_j=t}} X_{\us}
  \end{equation*}
\end{definition}
\begin{definition}\label{def-4-proper-ltawg}
  Let $X\in\cat{n}$, $1\leq i\leq n$ and let $\us$ be an orientation with $s_i=1$ and $t(\us)=t$. Denote
  \begin{equation*}
    \us(0,i)=(s_1,\ldots,s_{i-1},0,s_{i+1},\ldots,s_n)\;.
  \end{equation*}
  Note that $\us(0,i)$ is also an orientation, with $t(\us(0,i))=t-1$. Let $r\geq 2$. An $(r,i)$-string of $(n,t)$-hypercubes \index{String of $(n,t)$-hypercubes}with orientations $\us$ is an element of
  \begin{equation}\label{eq0-proper-ltawg}
    \pro{X_{\us}}{X_{\us(0,i)}}{r}\cong X_{(s_1,\ldots,s_{i-1},r,s_{i+1},\ldots,s_n)}
  \end{equation}
  where the isomorphism in \eqref{eq0-proper-ltawg} follows by the Segal condition characterizing multinerves of $n$-fold categories (see Lemma \ref{lem-multin-iff}).
\end{definition}
\begin{example}\label{ex1-proper-ltawg}
Let $X\in\cat{2}$. Then
\begin{equation*}
  \Cube{2,0}=X_{00}
\end{equation*}
consists of the set of objects in the double category $X$. Also,
\begin{equation*}
  \Cube{2,1}=X_{10}\cup X_{01}
\end{equation*}
consists of the set of arrows with orientations $(1,0)$ and $(0,1)$ which are respectively direction 1 and 2 in Figures \ref{corner-double-nerve} and \ref{corner2} on page \pageref{corners2}. Also
%ù
\begin{equation*}
  \Cube{2,2}=X_{11}
\end{equation*}
is the set of squares.

If $r\geq 2$, a $(r,1)$-string of $(2,1)$-hypercubes is an element of
\begin{equation*}
  \pro{X_{10}}{X_{00}}{r}\cong X_{r0}
\end{equation*}
that is, a sequence of $r$ composable arrows in direction 1. A $(r,1)$-string of $(2,2)$-hypercubes is an element of
\begin{equation*}
  \pro{X_{11}}{X_{01}}{r}\cong X_{r1}
\end{equation*}
that is, a sequence of $r$ composable squares in direction 1.

 Similarly, a $(r,2)$-string of $(2,1)$-hypercubes is an element of
\begin{equation*}
  \pro{X_{01}}{X_{00}}{r}\cong X_{0r}
\end{equation*}
and a $(r,2)$-string of $(2,2)$-hypercubes is an element of
\begin{equation*}
  \pro{X_{11}}{X_{10}}{r}\cong X_{1r}
\end{equation*}
with geometric interpretations similar to the above, but in direction 2 rather than direction 1.
\end{example}
\begin{example}\label{ex2-proper-ltawg}
Let $X\in\cat{3}$. Then
\begin{equation*}
  \Cube{3,0}=X_{000}
\end{equation*}
is the set of objects if the 3-fold category $X$. Also,
\begin{equation*}
  \Cube{3,1}=X_{100}\cup X_{010}\cup X_{001}
\end{equation*}
is the set of edges in the three orientations (see Figures \ref{corner3X} and \ref{corner3} on page \pageref{corners3}). Similarly,
\begin{equation*}
  \Cube{3,2}=X_{110}\cup X_{101}\cup X_{011}
\end{equation*}
is the set of squares, and
\begin{equation*}
  \Cube{3,3}=X_{111}
\end{equation*}
is the set of cubes.

Given $r\geq 2$, a $(r,1)$-string of $(3,1)$-hypercubes is an element of
\begin{equation*}
  \pro{X_{100}}{X_{000}}{r} \cong X_{r00}
\end{equation*}
that is a sequence of $r$-composable edges in direction 1. Similarly for $(r,2)$-strings and $(r,3)$-strings of $(3,1)$-hypercubes.

A $(r,1)$-string of $(3,2)$-hypercubes is an element of
\begin{equation*}
  \pro{X_{110}}{X_{010}}{r}\cong X_{r10}
\end{equation*}
that is a sequence of $r$ composable squares in direction 1. Similarly for $(r,2)$-strings and $(r,3)$-strings of $(3,2)$-hypercubes.

A $(r,1)$-string of $(3,3)$-hypercubes is an element of
\begin{equation*}
  \pro{X_{111}}{X_{011}}{r}\cong X_{r11}
\end{equation*}
that is a sequence of $r$ composable cubes in direction 1. Similarly for $(r,2)$-strings and $(r,3)$-strings of $(3,3)$-hypercubes.
\end{example}

\bk

Let $X\in\cat{n}$, $\us\in\dop{n}$ be an orientation and $1\leq i\leq n-1$. Let
\begin{equation*}
  s\up{i,n}(t,q)=(s_1,\ldots,s_{i-1},t,s_{i+1},\ldots,s_{n-1},q)\in \dop{n}
\end{equation*}
clearly $s\up{i,n}(0,0)$, $s\up{i,n}(0,1)$, $s\up{i,n}(1,0)$, $s\up{i,n}(1,1)$ are orientations and if we denote
\begin{equation*}
  t=1+\sum^n_{\substack{{j=1}\\j\neq i,n}} s_j
\end{equation*}
we have
\begin{equation*}
\begin{split}
    & t(s\up{i,n}(0,0))=t-1 \\
    & t(s\up{i,n}(1,0))=t(s\up{i,n}(0,1))=t\\
    & t(s\up{i,n}(1,1))=t+1\;.
\end{split}
\end{equation*}
Let $X\in\cat{n}$ and let $X(\us,i,n)\in\funcat{2}{\Set}$ with
\begin{equation*}
  X(\us,i,n)(v,w)=X_{(s_1,\ldots,s_{i-1},v,s_{i+1},\ldots,s_{n-1},w)}
\end{equation*}
then $X(\us,i,n)$ is the double nerve of a double category, which we still denote by $X(\us,i,n)$. The objects of this double category are
\begin{equation*}
  \begin{split}
   \{X(\us,i,n)\}_{00}=X_{\us\up{i,n}(0,0)}= & \{(n,t-1)\text{-hypercubes in }X  \\
      & \text{ with orientation }\us\up{i,n}(0,0)\}.
  \end{split}
\end{equation*}
The arrows in the double category are
\begin{equation*}
  \begin{split}
   \{X(\us,i,n)\}_{10}=X_{\us\up{i,n}(1,0)}= & \{(n,t)\text{-hypercubes in }X  \\
      & \text{ with orientation }\us\up{i,n}(1,0)\},
  \end{split}
\end{equation*}
and
\begin{equation*}
  \begin{split}
   \{X(\us,i,n)\}_{01}=X_{\us\up{i,n}(0,1)}= & \{(n,t)\text{-hypercubes in }X  \\
      & \text{ with orientation }\us\up{i,n}(0,1)\}.
  \end{split}
\end{equation*}
The squares in the double category are
\begin{equation*}
  \begin{split}
   \{X(\us,i,n)\}_{11}=X_{\us\up{i,n}(1,1)}= & \{(n,t+1)\text{-hypercubes in }X  \\
      & \text{ with orientation }\us\up{i,n}(1,1)\}\;.
  \end{split}
\end{equation*}
If $r\geq 2$, a sequence of $r$ composable arrows in the double category in direction 1 is an element of
\begin{equation*}
  \pro{X_{\us\up{i,n}(1,0)}}{X_{\us\up{i,n}(0,0)}}{r}\;.
\end{equation*}
By the above geometric interpretation of arrows of $X(\us,i,n)$ it follows that this is a $(r,i)$-string of $(n,t)$-hypercubes in $X$ with orientation $\us\up{i,n}(1,0)$ in the sense of Definition \ref{def-4-proper-ltawg}.

\bk

Suppose now that $X\in\catwg{n}$. If we let
\begin{equation*}
  \uk=(\seqc{s}{1}{n-1})\in\dop{n-1}
\end{equation*}
then
\begin{equation*}
\begin{split}
    & (\uk(0,i),0)=\us\up{i,n}(0,0) \\
    & (\uk(1,i),0)=\us\up{i,n}(1,0) \\
    & (\uk(0,i),1)=\us\up{i,n}(0,1) \\
    & (\uk(1,i),1)=\us\up{i,n}(1,1)\;.
\end{split}
\end{equation*}
Thus
\begin{equation}\label{eq1-proper-ltawg}
\begin{split}
    & X_{\us\up{i,n}(0,0)}=obj(J_n X)_{\uk(0,i)} \\
    & X_{\us\up{i,n}(1,0)}=obj(J_n X)_{\uk(1,i)} \\
    & X_{\us\up{i,n}(0,1)}=mor(J_n X)_{\uk(0,i)} \\
    & X_{\us\up{i,n}(1,1)}=mor(J_n X)_{\uk(1,i)}
\end{split}
\end{equation}
By Proposition \ref{cor-lev-wg-ncat}, $X\in\lta{n}$ so by Proposition \ref{pro-maps-nu-eqcat} there is an equivalence of categories for each $r\geq 2$
\begin{equation*}
\begin{split}
    & \pro{(J_n X)_{\uk(1,i)}}{(J_n X)_{\uk(0,i)}}{r} \rw\\
    & \rw \pro{(J_n X)_{\uk(1,i)}}{(J_n X)^d_{\uk(0,i)}}{r}\;.
\end{split}
\end{equation*}
By \eqref{eq1-proper-ltawg} and the above this means that the double category $X(\us,i,n)$ is weakly globular. As discussed in Example \ref{ex1-n-equival} this implies that every staircase of length $r$ of horizontal arrows in $X(\us,i,n)$ can be lifted to a string of $r$ composable horizontal arrows by strings of vertically invertible squares. To give a geometric interpretation of this lifting condition for the weakly globular double category $X(\us,i,n)$ we introduce the following definition:

\begin{definition}\label{def-5-proper-ltawg}
A staircase of length $r$ of horizontal morphisms in $X(\us,i,n)$ is called a $(r,i)$-staircase of $(n,t)$-hypercubes \index{Staircase of $(n,t)$-hypercubes}in $X$ with orientations $\us\up{i,n}(1,0)$.
\end{definition}
Definition \ref{def-5-proper-ltawg} means that the $(n,t)$-hypercubes are not composable in direction $i$, but are joined by strings of $(n,t)$-hypercubes with orientations $\us\up{i,n}(0,1)$.

  From the geometric interpretation of arrows and squares in the double category $X(\us,i,n)$ given above, we deduce that the lifting condition in the weakly globular double category $X(\us,i,n)$ implies the following:

  \begin{proposition}\label{prop-geom-int}
   Given $X\in \catwg{n}$, every  $(r,i)$-staircase of $(n,t)$-hypercubes in $X$ with orientations $\us\up{i,n}(1,0)$ can be lifted to a $(r,i)$-string of $(n,t)$-hypercubes in $X$ via strings of $(n,t+1)$-hypercubes with orientations $\us\up{i,n}(1,1)$ which are invertible in direction $n$.
  \end{proposition}

    We illustrate below the case $n=3$.
\begin{example}\label{ex3-proper-ltawg}
Let $X\in\catwg{3}$, $r\geq 2$. A $(r,1)$-staircase of $(3,1)$-hypercubes with orientation $(100)$ is an element of
\begin{equation*}
  \pro{X_{100}}{X^d_{000}}{r}.
\end{equation*}
This can be represented as a staircase
\bk

%%
%% Create 1-edge
\tikzset{edge/.pic={
\filldraw (0,0)  circle[radius=0.035cm] -- (1,0)  circle[radius=0.035cm]; %% centered
}}
\begin{figure}[h]
\centering
\begin{tikzpicture}[thick,scale=1]]
\pic [black][scale=1] at (0,0) {edge};
\pic [black][scale=1] at (1,1) {edge};
\pic [black][scale=1] at (2,0.5) {edge};
\pic [black][scale=1] at (3,1.5) {edge};
\pic [black][scale=1] at (5,0) {edge};
\draw[thick, densely dashed] (1,0) -- (1,1);
\draw[thick, densely dashed] (2,1) -- (2,0.5);
\draw[thick, densely dashed] (3,1.5) -- (3,0.5);
\draw[thick, densely dashed] (4,1.5) -- (4,0.5);
\draw[thick, densely dashed] (5,0.5) -- (5,0);
\draw[ultra thick, loosely  dotted] (4.3,0.5) -- (4.65,0.5);
\end{tikzpicture}
\caption{$(r,1)$-staircase of $(3,1)$-hypercubes in $X$}
\label{stair1}
\end{figure}
\nid with edges in direction 1 and vertical sides in direction 3. The equivalence of categories
\begin{equation*}
  \pro{X_{100}}{X_{000}}{r} \simeq \pro{X_{100}}{X^d_{000}}{r}
\end{equation*}
means that this staircase can be lifted to a string of $r$ composable arrows in direction 1 through strings of vertically invertible squares.
\bk

\tikzset{edge/.pic={
\filldraw (0,0)  circle[radius=0.04cm] -- (1,0)  circle[radius=0.04cm]; %% centered
}}
\begin{figure}[h]
\centering
\begin{tikzpicture}[thick,scale=1]]
\pic [black][scale=1] at (0,0) {edge};\pic [black][scale=1] at (0,2.5) {edge};
\node[rotate=-90,scale=1.2] at (0.5,1.5) () {$\bm\Rightarrow$};
\pic [black][scale=1] at (1,1) {edge};\pic [black][scale=1] at (1,2.5) {edge};
\node[rotate=-90,scale=1.2] at (1.5,1.8) () {$\bm\Rightarrow$};
\pic [black][scale=1] at (2,0.5) {edge};\pic [black][scale=1] at (2,2.5) {edge};
\node[rotate=-90,scale=1.2] at (2.5,1.5) () {$\bm\Rightarrow$};
\pic [black][scale=1] at (3,1.5) {edge};\pic [black][scale=1] at (3,2.5) {edge};
\node[rotate=-90,scale=1.2] at (3.5,2.0) () {$\bm\Rightarrow$};
\pic [black][scale=1] at (5,0) {edge};\pic [black][scale=1] at (5,2.5) {edge};
\node[rotate=-90,scale=1.2] at (5.5,1.5) () {$\bm\Rightarrow$};
\draw[ultra thick, loosely dotted] (4.3,2.5) -- (4.6,2.5);
\draw[thick, densely dashed] (0,0) -- (0,2.5);
\draw[thick, densely dashed] (1,0) -- (1,2.5);
\draw[thick, densely dashed] (2,0.5) -- (2,2.5);
\draw[thick, densely dashed] (3,0.5) -- (3,2.5);
\draw[thick, densely dashed] (4,0.5) -- (4,2.5);
\draw[thick, densely dashed] (5,0) -- (5,2.5);
\draw[thick, densely dashed] (6,0) -- (6,2.5);
\draw[ultra thick, loosely dotted] (4.3,0.5) -- (4.6,0.5);
\end{tikzpicture}
\caption{Lifting condition for staircase in figure \ref{stair1}}
\label{stair2}
\end{figure}

A $(r,2)$-staircase of $(3,1)$-hypercubes with orientation $(010)$ is an element of
\begin{equation*}
  \pro{X_{010}}{X^d_{000}}{r}\;.
\end{equation*}
This can be pictured as a staircase like the one in Figure \ref{stair1} but with horizontal edges in direction 2 and vertical edges in direction 3. The lifting condition is similar to Figure \ref{stair2}, but now the vertical squares have orientation in directions 2 and 3 rather than 1 and 3.

A $(r,1)$-staircase of $(3,2)$-hypercubes with orientation $(110)$ is an element of
\begin{equation*}
  \pro{X_{110}}{X^d_{010}}{r}\;.
\end{equation*}
This can be pictured as in Figure \ref{stair3D1} on page \pageref{3Dstair} (where for simplicity we choose $r=3$)

The equivalence of categories
\begin{equation*}
\pro{X_{110}}{X_{010}}{r} \simeq \pro{X_{110}}{X^d_{010}}{r}
\end{equation*}
means that this staircase can be lifted to a sequence of $r$ composable squares in direction 1 via strings of cubes which are vertically invertible. See Figure \ref{stair3D2} on page \pageref{3Dstair} (when $r=3$).

A $(r,2)$-staircase of $(3,2)$-hypercubes with orientation $(110)$ is an element of
\begin{equation*}
  \pro{X_{110}}{X^d_{100}}{r}\;.
\end{equation*}
The geometric picture is obtained by rotating Figure \ref{stair3D1} by $90^\circ$ around the vertical axis, and the lifting condition is similar to Figure \ref{stair3D2} in the rotated picture.

\end{example}
%%%%%%%%%%%%%%%%%%%%%%%%%%%%%%%%%%%%%%%%%%%%%%%%%%%%%%%%%%%%%%%%%%%%%%%%%%%%%%%%%%%

\section{Approximating $\pmb{\tawg{n}}$ with $\pmb{\lta{n}}$.}\label{sec-approx}

The main result of this section, Theorem \ref{the-repl-obj-1}, states that if  $X\in\tawg{n}$ is such that $\qn X$ can be approximated up to $(n-1)$-equivalence with an object of $\catwg{(n-1)}$, then $X$ can be approximated up to an $n$-equivalence with an object of $\lta{n}$. In the next Section this is used in the proof of Theorem \ref{the-funct-Qn} to construct the functor
 \begin{equation*}
    P_n:\tawg{n}\rw\lta{n}
\end{equation*}
from which the rigidification functor $Q_n$ will be built.

\subsection{Main steps in approximating $\pmb{\tawg{n}}$ with $\pmb{\lta{n}}$.}\label{subs-approx-idea}
The basic construction is the pullback in $\funcat{n-1}{\Cat}$
\begin{equation}\label{eq1-subs-approx-idea}
\xymatrix@R=30pt@C=45pt{
P \ar^{w}[r] \ar[d] & X \ar^{\zg\up{n}}[d]\\
\di{n} Z \ar_{\di{n}r}[r] & \di{n}\q{n}X
}
\end{equation}
with $X\in\tawg{n}$, $Z\in\catwg{n-1}$ and $r: Z\rw \q{n} X$ a $(n-1)$-equivalence in $\tawg{n-1}$. As in Section \ref{sec-pull-qn}, in the pullback \eqref{eq1-subs-approx-idea} we have omitted writing $J_n$ for ease of notation. We show in Theorem \ref{the-repl-obj-1} that $P\in\lta{n}$ and $w$ is an $n$-equivalence.

This construction is crucial to proving the existence of the rigidification functor $\Qn$: in the
proof of Theorem \ref{the-funct-Qn} we will use the above construction with $Z=\Qnm\q{n}X$ for
\begin{equation*}
  \Qnm : \tawg{n-1}\rw \catwg{n-1}
\end{equation*}
inductively defined and the $(n-1)$-equivalence $r$ also constructed from the inductive hypothesis.

By Proposition \ref{pro-spec-plbk-pscatwg}, $P\in\tawg{n}$. Thus, by Lemma \ref{lem-lev-wg-pscat}, to show that $P\in\lta{n}$ it is enough to show that $P\in\lnta{n}{n}$  and $P_{\seq{k}{1}{n-r}}\in\lnta{r}{r}$ for each $1 < r\leq n-1$, $\seq{k}{1}{n-r}\in\Dop$. By definition, to prove that $P\in\lnta{n}{n}$ we need to show that

\mk
\nid (i) For each $k\geq 2$ the maps
\begin{equation}\label{eq2-subs-approx-idea}
    v_k : J_{n-1}X_k \rw \pro{J_{n-1}X_1}{J_{n-1}\di{n-1}\p{n-1}X_0}{k}
\end{equation}
  are levelwise equivalence of categories.

\bk
\nid (ii)  $\p{n}P\in\catwg{n-1}$\;.

\bk
The idea to verify condition (i) is to apply the criterion of Proposition \ref{pro-crit-lev-nequiv} (for a $n$-equivalence in $\tawg{n}$ to be a levelwise equivalence of categories) to the map
\begin{equation*}
  P_k \rw \pro{P_1}{\di{n-1}\p{n-1}P_0}{k}\;.
\end{equation*}
That is, we want to show this is a $(n-1)$-equivalence in $\tawg{n-1}$ satisfying the additional conditions in the hypotheses of Proposition \ref{pro-crit-lev-nequiv}: the idea is that these additional conditions are forced by the fact that $Z\in\catwg{n-1}$.

The intuition behind this idea is as follows. Since \eqref{eq1-subs-approx-idea} is a pullback in $\funcat{n-1}{\Cat}$, it is computed levelwise, so for each $k\geq 0$ there is a pullback in $\funcat{n-3}{\Cat}$
\begin{equation}\label{eq3-subs-approx-idea}
\xymatrix@R=30pt@C=45pt{
P_{k0} \ar^{}[r] \ar[d] & X_{k0} \ar^{{(\zg\up{n})}_{k0}}[d]\\
\di{n-2} Z_{k0}=(\di{n}Z)_{k0} \ar_{}[r] & (\di{n}\q{n}X)_{k0}
}
\end{equation}
As $X_{k0}$ is homotopically discrete,
\begin{equation*}
  (\di{n}\q{n}X)_{k0}=\di{n-2}\q{n-2}X_{k0}= \di{n-2}\p{n-2}X_{k0}\;.
\end{equation*}
Since $\p{n-2}$ commutes with pullbacks over discrete objects, from \eqref{eq3-subs-approx-idea} we deduce
\begin{equation*}
  \p{n-2}P_{k0}= Z_{k0}\;.
\end{equation*}
Thus, for instance
\begin{equation*}
\begin{split}
    & \p{n-2}(\tens{P_{1}}{\di{n-1}\p{n-1}P_0})_0 = \p{n-2}(\tens{P_{10}}{\di{n-2}\p{n-2}P_{00}})= \\
   = & \tens{\p{n-2}P_{10}}{\p{n-2}P_{00}}=\tens{Z_{10}}{Z_{00}}\;.\\
\end{split}
\end{equation*}
Since $Z\in\catwg{n}$, $Z_{20}\cong \tens{Z_{10}}{Z_{00}}$ so from above
\begin{equation*}
  \p{n-2}P_{20}\cong \p{n-2}(\tens{P_{1}}{\di{n-1}\p{n-1}P_0})_0 \;.
\end{equation*}
Thus condition b) in the hypotheses of Proposition \ref{pro-crit-lev-nequiv} holds for the map
\begin{equation*}
  P_2\rw \tens{P_{1}}{\di{n-1}\p{n-1}P_0}\;.
\end{equation*}
The main steps in the formal proofs are as follows:
\begin{itemize}
  \item [a)] We apply Proposition \ref{pro-spec-plbk-pscatwg} to deduce that $P\in\tawg{n}$.
  \bk

  \item [b)] We show in Proposition \ref{cor-rep-ob} that the maps \eqref{eq2-subs-approx-idea} are levelwise equivalence of categories. The proof of this fact is based on Lemma \ref{lem-jn-alpha} applied to the induced Segal maps of $P$ which are shown to factor as a composite of maps in $\tawg{n-1}$
      \begin{equation*}
        \hmu{k}: P_k \xrw{\za} \pro{P_1}{\di{n-1}\p{n-1}P_0}{k} \xrw{\zb} \pro{P_1}{P_0^d}{k}\;.
      \end{equation*}
      The fact that $Z\in\catwg{n-1}$ forces additional conditions (hypotheses i) and ii) of Lemma \ref{lem-jn-alpha}) which allow to apply the criterion of Proposition \ref{pro-crit-lev-nequiv} and show that $\za$ is a levelwise equivalence of categories.
      \bk

  \item [c)] In Theorem \ref{the-repl-obj-1} we use again Proposition \ref{cor-rep-ob} and Lemma \ref{lem-x-in-tawg-x-in-catwg} to show that $\p{n}P$ satisfies the hypotheses of Lemma  \ref{lem-x-in-tawg-x-in-catwg} and therefore $\p{n}P \in \catwg{n}$. Together with Proposition \ref{cor-rep-ob} this implies by definition that $P\in \lnta{n}{n}$. Working inductively we then easily establish that $P_{\seq{k}{1}{n-r}}\in \lnta{r}{r}$ for all $1 < r \leq n-1$ and thus we conclude (by Lemma \ref{lem-lev-wg-pscat}) that $P\in \lta{n}$.
\end{itemize}
%%%%%%%%%%%%%%%%%%%%%%%%%%%%%%%%%%%%%%%%%%%%%%%%%%%%%%%%%%%%%%%%%%%%%%%%%%%%%%%

\subsection{Approximating $\pmb{\tawg{n}}$ with $\pmb{\lta{n}}$: the formal proofs}

The following lemma and its corollary are used in the proof of Theorem \ref{the-repl-obj-1}. Their proof use the criterion given in Proposition \ref{pro-crit-lev-nequiv} for an $n$-equivalence in $\tawg{n}$ to be a levelwise equivalence of categories and the properties of pullbacks along the map $\zg\up{n}$ established in Section \ref{sec-pull-qn}.
\begin{lemma}\label{lem-jn-alpha}
    Let
    \begin{equation*}
      B \xrw{\pt_0} X \xlw{\pt_1}B
    \end{equation*}
      be a diagram in $\tawg{n}$ with $X\in\cathd{n}$, $\tens{B}{X}\in \tawg{n}$ and let
    \begin{equation*}
        A\xrw{\za}\tens{B}{X}\xrw{\zb}\tens{B}{X^d}
    \end{equation*}
    be maps $\tawg{n}$ (where $\zb$ is induced by the map $\zg: X\rw X^d$)  such that\mk
    \begin{itemize}
      \item [i)] $\p{n-1}\za_0$, and $\p{n-r-1}\za_{k_1...k_r\,0}$  are isomorphisms for all $1\leq r< n-1$.\mk
      \item [ii)] $(\tens{B}{X})^d_0\cong\tens{B_0^d}{X_0^d}$ \mk

     \nid             $(\tens{B}{X})^d_{k_1...k_r\,0}\cong\tens{B^d_{k_1...k_r\,0}}{X^d_{k_1...k_r\,0}} $ for all $1\leq r< n-1$\mk
      \item [iii)] $\zb\za$ is a $n$-equivalence.\mk
    \end{itemize}
    Then $J_{n}\za$ is a levelwise equivalence of categories.
\end{lemma}
\begin{proof}
Let $x,x'\in A_0^d$. By hypothesis i) and ii),
\begin{equation*}
    A_0^d\cong(\tens{B}{X})^d_0=\tens{B_0^d}{X_0^d}\subset \tens{B_0^d}{X^d}
\end{equation*}
where the last inclusion holds because the map $\zg_0:X_0\rw (X^d)_0$ factors through $X_0^d$. Let
 \begin{equation*}
 \za x=(a,b)\in \tens{B_0^d}{X_0^d},\qquad \za x'=(a',b')\in \tens{B_0^d}{X_0^d}.
 \end{equation*}
 We claim that the composite map
\begin{equation}\label{eq1-lem-jn-alpha}
    A(x,x')\xrw{\za(x,x')} B(a,a')\tiund{X(\pt a,\pt a')}B(b,b')\xrw{s} B(a,a')\tiund{X(\pt a,\pt a')^d}B(b,b')
\end{equation}
is a $\equ{n-1}$. We have
\begin{flalign*}
&(\tens{B}{X})(\za x,\zb x')=B(a,a')\tiund{X(\pt_0 a,\pt_0 a')}B(b,b')&\\
&(\tens{B}{X^d})(\zb \za x,\zb \za x')=B(a,a')\tiund{X^d(\zg a,\zg a')}B(b,b')=B(a,a')\times B(b,b')&
\end{flalign*}
where in the last equality we used the fact that $\zg a=\zg a'$ so that $X^d(\zg a,\zg a')=\{\cdot\}$ since $X^d$ is discrete.

The map $\zg:X\rw X^d$ factors as
\begin{equation*}
    X\rw \di{n}...\di{2}\p{2}...\p{n}X\rw \di{n}...\di{2}d p...\p{n}X=X^d
\end{equation*}
and we have
\begin{equation*}
    (\di{n}...\di{2}\p{2}...\p{n}X)(p \pt_0 a,p \pt_0 a')=X(\pt_0 a,\pt_0 a')^d\;.
\end{equation*}
Thus the map $\zb(x,x')$ factors as
\begin{equation*}
\resizebox{1.0\hsize}{!}{$
    B(a,a')\tiund{X(\pt_0 a,\pt_0 a')}B(b,b')\xrw{s} B(a,a')\tiund{X(\pt_0 a,\pt_0 a')^d}B(b,b')\xrw{t}B(a,a')\times B(b,b')\;.$}
\end{equation*}
On the other hand, since $\p{2,n}X\in\cathd{}$, the set
 \begin{equation*}
 \p{2,n}X(p\pt_0 a,p\pt_0 a')
 \end{equation*}
  contains only one element. Thus $X(\pt_0 a,\pt_0 a')^d$ is the terminal object, and $t=\Id$, so that
  \begin{equation}\label{eq-compzazb}
    (\zb\za)(x,x')=s\za(x,x').
  \end{equation}
Since, by hypothesis, $\zb\za$ is a $n$-equivalence $\zb\za(x,x')$ is a $\nm$-equivalence, so by \eqref{eq-compzazb} the composite \eqref{eq1-lem-jn-alpha} is a $\equ{n-1}$. This proves the claim.

\medskip
\nid We now proceed to the rest of the proof by induction on $n$. The strategy is to show that $\za$ satisfies the hypotheses of Proposition \ref{pro-crit-lev-nequiv}.
\bk

 When $n=2$, since $X(\pt_0 a,\pt_0 a')\in\cathd{}$, the map $s$ is fully faithful. Since  $s\za(x,y)$
 is an equivalence of categories, it is essentially surjective on objects, and therefore $s$ is essentially surjective on objects. It follows that $s$ is an equivalence of categories, and therefore such is $\za(x,y)$. Since by hypothesis $p\za_0$ is a bijection, the map $\p{2}\za$ is bijective on objects, thus $p\p{2}\za$ is surjective. From Proposition \ref{pro-n-equiv} we deduce that $\za$ is a 2-equivalence. Thus $\za$ satisfies the hypotheses of Proposition \ref{pro-crit-lev-nequiv} and we conclude that it is a levelwise equivalence of categories.
 \medskip

\nid Suppose, inductively, that the lemma holds for $n-1$.

 We show that the maps \eqref{eq1-lem-jn-alpha} satisfy the inductive hypothesis.
Since, as proved above, $s\za(x,y)$ is a $\equ{n-1}$, inductive hypothesis iii) holds. Since by hypothesis $\p{n-1}\za_0$ is an isomorphism, so is
\begin{equation*}
    \p{n-2}\za_0(x,x')=(\p{n-1}\za_0)(x,x')
\end{equation*}
as well as
\begin{equation*}
    \p{n-r-2}\{\za_{k_1...k_r\,0}(x,x')\}=(\p{n-r-1}\za_{k_1...k_r\,0})(x,x')\;.
\end{equation*}
Thus inductive hypothesis i) holds for the maps \eqref{eq1-lem-jn-alpha}. Further, using hypothesis ii) we compute
\begin{equation*}
\begin{split}
    & ((\tens{B}{X})(\za x,\za x'))^d_0=(\tens{B_1}{X_1})^d_0(\za x,\za x')= \\
    =\ & (\tens{B^d_{10}}{X^d_{10}})(\za x,\za x')=\tens{B(\za x,\za x')^d_0}{X(\za x,\za x')^d_0}
\end{split}
\end{equation*}
and similarly
\begin{equation*}
\begin{split}
    & ((\tens{B}{X})(\za x,\za x'))^d_{k_1...k_r\,0}=(\tens{B}{X})^d_{1 k_1...k_r\,0}(\za x,\za x')= \\
    \cong\ & \tens{B^d_{1  {k_1...k_r\,0}}(\za x,\za x')}{X^d_{1 {k_1...k_r\,0}}(\za x,\za x')}=\\
    =\ &\tens{B_{k_1...k_r\,0}(\za x,\za x')^d}{X_{k_1...k_r\,0}(\za x,\za x')^d}\;.
\end{split}
\end{equation*}
Thus the inductive hypothesis ii) holds for the maps \eqref{eq1-lem-jn-alpha}.

We conclude by induction that $J_{n-1}\za(x,x')$ is a levelwise equivalence of categories. It follows by Remark \ref{rem-local-equiv}, $\za(x,x')$ is a $\equ{n-1}$. That is, $\za$ is a local $\equ{n-1}$.

 Since by hypothesis $\p{n-1}\za_0$ is an isomorphism, so is
\begin{equation*}
    (\p{2,n}\za)_0= \p{1,n-1}\za_0
\end{equation*}
so that $\p{1,n}\za$ is surjective. Since, from above, $\za$ is a local $\nm$-equivalence, from Proposition \ref{pro-n-equiv} we conclude that $\za$ is a $\nequ$. Together with hypothesis i) this shows that $\za$ satisfies the hypotheses of Proposition \ref{pro-crit-lev-nequiv} and we conclude that $J_{n}\za$ is a levelwise equivalence of categories.
\end{proof}
\begin{remark}\label{rem-approx-tawg}
Lemma \ref{lem-jn-alpha} also generalizes as follows, with completely analogous proof. If $B$ and $X$ are as in Lemma  \ref{lem-jn-alpha}, $k\geq 2$ and
\begin{equation*}
  A\xrw{\za}\pro{B}{X}{k} \xrw{\zb} \pro{B}{X^d}{k}
\end{equation*}
are maps in $\tawg{n}$ such that
\begin{itemize}
  \item [i)] $\p{n-1}\za_0$ \,and\, $\p{n-r-1}\za_{\seq{k}{1}{r} 0}$ are isomorphisms for all $1\leq r \leq n-1$.\bk
  \item [ii)] $(\pro{B}{{X}}{k})_0^d \cong \pro{B_0^d}{X^d}{k}$.
  \nid $(\pro{B}{X}{k})_{k_1...k_r\,0}^d \cong \pro{B_{k_1...k_r\,0}^d}{X_{k_1...k_r\,0}^d}{k}$ for all $1\leq r \leq n-1$.\bk
  \item [iii)] $\zb\za$ is a $n$-equivalence.\bk
\end{itemize}
  Then $J_n\za$ is a levelwise equivalence of categories.
\end{remark}

Using the lemma and remark above, we now prove the following Proposition, which will be used in the proof of Theorem \ref{the-repl-obj-1}. In the proof of this result we use the properties of pullbacks along the map $\zg\up{n}$ established in Section \ref{sec-pull-qn}. As in Section \ref{sec-pull-qn}, we will always consider $\tawg{n}$ (as well as $\catwg{n}$) as embedded in $\funcat{n-1}{\Cat}$ via the functor $J_n$ and our pullbacks will be taken in $\funcat{n-1}{\Cat}$. To ease the notation, we omit writing explicitly $J_n$ in these pullbacks.
\begin{proposition}\label{cor-rep-ob}
  Let $X\in \tawg{n}$, and let
     \begin{equation*}
       r:Z\rw \di{n} \qn X
     \end{equation*}
      be a map in $\tawg{n-1}$ with $Z\in\catwg{n-1}$ and consider the pullback in $\funcat{n-1}{\Cat}$
    \begin{equation*}
        \xymatrix@R=35pt @C=40pt{
        P \ar^{w}[r] \ar^{}[d] & X \ar^{\zgu{n}}[d] \\
        \dn Z \ar_{\dn r}[r] & \dn \qn X
        }
    \end{equation*}
  Then $P\in \tawg{n}$ and for all $ k\geq 2$ and
      \begin{equation*}
        v_k: J_{n-1}P_k\rw \pro{J_{n-1}P_1}{J_{n-1}\di{n-1}\p{n-1}P_0}{k}
      \end{equation*}
      is a levelwise equivalence of categories.
\end{proposition}
\begin{proof}
Throughout this proof we will, for simplicity, denote the map $v_k$ by $\za$.
 By Proposition \ref{pro-spec-plbk-pscatwg}, $P\in\tawg{n}$, therefore its induced Segal maps $\hmu{k}$ are $(n-1)$-equivalences. The strategy of the proof is to show that for each $k\geq 2$ $\hmu{k}$ factors in $\tawg{n-1}$ as
\begin{equation*}
  \hmu{k}: P_k\xrw{\za}\pro{P_1}{\di{n-1}\p{n-1}P_0}{k} \xrw{\zb}\pro{P_1}{P_0^d}{k}
\end{equation*}
and this satisfies the hypotheses of Lemma \ref{lem-jn-alpha} (see also Remark \ref{rem-approx-tawg}): the hypothesis iii) of Lemma \ref{lem-jn-alpha} holds since $P\in\tawg{n}$ while hypotheses i) and ii) will be a direct consequence of the fact that $Z\in\catwg{n-1}$, as illustrated below.

We first show that for each $k\geq 2$
\begin{equation*}
    \pro{P_1}{\di{n-1}\p{n-1}P_0}{k}\in \tawg{n-1}\;.
\end{equation*}
We illustrate this for $k=2$, the case $k>2$ being similar. Since $p$ commutes with pullbacks over discrete objects we have
\begin{equation*}
    \p{n}P=Z\tiund{\q{n}X}\p{n}X
\end{equation*}
and, (since $\q{n-1}X_0=\p{n-1}X_0$ as $X_0\in\cathd{n-1}$)
\begin{equation*}
    \p{n-1}P_0=Z_0\tiund{\q{n-1}X_0}\p{n-1}X_0=Z_0\;.
\end{equation*}
Also, $P_1=\di{n-1}Z_1\tiund{\di{n-1}\q{n-1}X_1}X_1$. Therefore
\begin{equation*}
\begin{split}
    & \tens{P_1}{\di{n-1}\p{n-1}P_0}= \\
    =\ & \tens{(\di{n-1}Z_1\tiund{\di{n-1}\q{n-1}X_1}X_1)}{\di{n-1}Z_0}=\\
    =\ & \di{n-1}(\tens{Z_1}{Z_0})\tiund{\di{n-1}\q{n-1}({X_1}\times{X_1})} ({X_1}\times{X_1})\;.
\end{split}
\end{equation*}
 By Proposition \ref{pro-spec-plbk-pscatwg}, this is an object of $\tawg{n-1}$.

 The induced Segal map $\hmu{2}$ for $P$ can therefore be written as composite of maps in $\tawg{n-1}$
\begin{equation}\label{eq-the-repl-obj-1}
    P_2 \xrw{\za} \tens{P_1}{\di{n-1}\p{n-1}P_0} \xrw{\zb} \tens{P_1}{P_0^d}\;.
\end{equation}
We show that the maps \eqref{eq-the-repl-obj-1}  satisfies the hypotheses of Lemma \ref{lem-jn-alpha}.

Since $P\in\tawg{n}$, $\hmu{2}=\zb\za$ is a $\equ{n-1}$, so hypothesis iii) of Lemma \ref{lem-jn-alpha} holds for the maps \eqref{eq-the-repl-obj-1} .

To check hypothesis i), note that
\begin{equation*}
\begin{split}
    & \p{n-2}P_{20}=\p{n-2}(\di{n-2}Z_{20}\tiund{\di{n-2}\q{n-2}X_{10}} \di{n-2}\p{n-2}X_{10})=Z_{20} \\
    & \p{n-2}(\tens{P_1}{\di{n-1}\p{n-1}P_0})_0=\\
    =\ &\tens{\p{n-2}P_{10}}{\p{n-2}P_{00}}=\tens{Z_{10}}{Z_{00}}\cong Z_{20}
\end{split}
\end{equation*}
where the last isomorphism holds since $Z\in\catwg{n-1}$.
Hence $\p{n-2}\za_0$ is an isomorphism. Similarly
\begin{equation*}
    P_{s_1...s_r}=Z_{s_1...s_r}\tiund{\q{n-r-1}X_{s_1...s_r}}X_{s_1...s_r}
\end{equation*}
\begin{equation*}
    \p{n-r-1}P_{s_1...s_r}=Z_{s_1...s_r}\tiund{\q{n-r-1}X_{s_1...s_r}}X_{s_1...s_r}\;.
\end{equation*}
Thus
\begin{equation*}
\begin{split}
    & \p{n-r-2}P_{2 k_1...k_r 0}=Z_{2 k_1...k_r 0} \\
    & \p{n-r-2}(\tens{P_1}{\di{n-1}\p{n-1}P_0})_{k_1...k_r 0}=\\
    =\ & \p{n-r-2}(\tens{P_{1 k_1...k_r 0}}{\di{n-r-1}\p{n-r-1}P_{0 k_1...k_r 0}})=\\
    =\ & \tens{\p{n-r-2}P_{1 k_1...k_r 0}}{\p{n-r-1}P_{0 k_1...k_r 0}}=\\
    =\ & \tens{Z_{1 k_1...k_r 0}}{Z_{0 k_1...k_r 0}}\cong Z_{2 k_1...k_r 0}\;.
\end{split}
\end{equation*}
where the last isomorphism holds since $Z\in\catwg{n-1}$. This shows that $\p{n-r-2}\za_{k_1...k_r 0}$ is an isomorphism, proving hypothesis i) of Lemma \ref{lem-jn-alpha} for the maps \eqref{eq-the-repl-obj-1}.

To check hypothesis ii) of Lemma \ref{lem-jn-alpha} note that $P^d_{20}=Z^d_{20}$ while
\begin{equation*}
\begin{split}
    & (\tens{P_1}{\di{n-1}\p{n-1}P_0})^d_0=(\tens{P_{10}}{\di{n-2}\p{n-2}P_{00}})^d_0= \\
    =\ & (\tens{Z_{10}}{\di{n-2}\p{n-2}Z_{00}})^d=\tens{Z_{10}^d}{Z_{00}^d}\cong Z_{20}^d
\end{split}
\end{equation*}
where the last isomorphism holds because $Z\in\catwg{n}$ (apply Remark \ref{rem-eq-def-wg-ncat} to $Z_0\up{2}$ which is an object of $\catwg{n-2}$ by Proposition \ref{pro-crit-ncat-be-wg}). Similarly
\begin{equation*}
\begin{split}
    & (P_2)^d_{k_1...k_r 0}=Z^d_{2 k_1...k_r 0} \\
    & (\tens{P_1}{\di{n-1}\p{n-1}P_0})^d_{k_1...k_r 0}=\\
    =\ & (\tens{Z_{1 k_1...k_r 0}}{\di{n-1}\p{n-1}Z_{0 k_1...k_r 0}})^d=\\
    =\ & \tens{Z^d_{1 k_1...k_r 0}}{Z^d_{0 k_1...k_r 0}}=Z^d_{2 k_1...k_r 0}\;.
\end{split}
\end{equation*}
where the last equality holds because $Z\in\catwg{n-1}$ and therefore $Z_{k_0...k_r}\in\catwg{n-r-1}$ by applying Remark \ref{rem-eq-def-wg-ncat} to $(Z_{k_0...k_r})_0\up{2}$ which is an object of $\catwg{n-r-1}$ by Proposition \ref{pro-crit-ncat-be-wg}. This proves that hypothesis ii) of Lemma \ref{lem-jn-alpha} holds for \eqref{eq-the-repl-obj-1}.

So all hypotheses of Lemma \ref{lem-jn-alpha} holds for the maps \eqref{eq-the-repl-obj-1}  and we conclude that $J_{n-1}\za$ is a levelwise equivalence of categories.

\end{proof}

\begin{theorem}\label{the-repl-obj-1}
    Let $X\in \tawg{n}$, and let
     \begin{equation*}
       r:Z\rw \di{n} \qn X
     \end{equation*}
      be a map in $\tawg{n-1}$ with $Z\in\catwg{n-1}$ and consider the pullback in $\funcat{n-1}{\Cat}$
    \begin{equation*}
        \xymatrix@R=35pt @C=40pt{
        P \ar^{w}[r] \ar^{}[d] & X \ar^{\zgu{n}}[d] \\
        \dn Z \ar_{\dn r}[r] & \dn \qn X
        }
    \end{equation*}
    Then
    \begin{itemize}
      \item [a)] $\q{n}P$ and $\p{n}P$ are in $\catwg{n-1}$.\mk

      \item [b)] $P\in\lnta{n}{n}$.\mk

      \item [c)] If $r$ is a $\nm$-equivalence then $w$ is a $\nequ$.\mk

      \item [d)] $P\in\lta{n}$.
    \end{itemize}
\end{theorem}
\begin{proof}
By induction in $n$. When $n=2$, we know by  Proposition \ref{pro-spec-plbk-pscatwg} that $P\in\tawg{2}=\lnta{2}{2}$, and that c) holds. Part a) is trivial since $\p{2}P$ and $\q{2}P$ are in $\Cat$. Part d) holds since $\tawg{2}=\lta{2}$

\medskip
Suppose, inductively that the theorem holds for $n-1$.

\bk

a) We have $\q{n}P=Z\in\catwg{n}$, $\p{n}P\in\tawg{n-1}$ since $P\in \tawg{n}$ by Proposition \ref{pro-spec-plbk-pscatwg}. We now show that $\p{n}P$ satisfies the hypotheses of Lemma \ref{lem-x-in-tawg-x-in-catwg}, which then shows that $\p{n}P\in\catwg{n-1}$.

We have the pullback in $\funcat{n-2}{\Cat}$ for each $k\geq 0$,
\begin{equation*}
\xymatrix{
P_k \ar[rr] \ar[d] && X_k \ar[d]\\
Z_k \ar[rr] && \di{n-1}\p{n-1}X_k
}
\end{equation*}
which\, satisfies\, the\, induction\, hypothesis. \,Therefore, \,by \,induction,\, $\p{n-1}P_k\in\catwg{n-2}$ which is hypothesis a) of Lemma \ref{lem-x-in-tawg-x-in-catwg}. Since, by Proposition \ref{cor-rep-ob}, the map
\begin{equation*}
    v_2:P_2\rw \tens{P_1}{\di{n-1}\p{n-1}P_0}
\end{equation*}
is a levelwise equivalence of categories, it induces an isomorphism
\begin{equation*}
    \p{n-1}P_2\cong\p{n-1}(\tens{P_1}{\di{n-1}\p{n-1}P_0})\cong\tens{\p{n-1}P_1}{\p{n-1}P_0}
\end{equation*}
and similarly all the other Segal maps of $\p{n}P$ are isomorphisms. This proves hypothesis b) of Lemma \ref{lem-x-in-tawg-x-in-catwg} for $\p{n}X$.

To prove hypothesis c) of Lemma \ref{lem-x-in-tawg-x-in-catwg}, we first note that
\begin{equation}\label{eq0-the-repl-obj-1}
    \tens{P_1}{\di{n-1,j}\p{j,n-1}P_0}\in \tawg{n-1}
\end{equation}
where $\di{n-1,j}$ is and in Notation \ref{dnot-ex-tam}.

In fact, since $\p{n-1}P_0=Z_0$ we have
\begin{equation}\label{eq1-the-repl-obj-1}
\begin{split}
    & \tens{P_1}{\di{n-1,j}\p{j,n-1}P_0}= \\
    & \resizebox{1.0\hsize}{!}{$ \;= \tens{(\di{n-1}Z_1 \tiund{\di{n-1}\q{n-1}X_1} X_1)}{\di{n-1,j}\p{j,n-1}Z_0}=$}\\
   & = \{\tens{\di{n-1}Z_1}{\di{n-1,j}\p{j,n-1}Z_0}\}\tiund{\di{n-1}\q{n-1}(X_1\times X_1)} (X_1\times X_1)=\\
   & = \di{n-1}(\tens{Z_1}{\di{n-1,j}\p{j,n-1}Z_0})\tiund{\di{n-1}\q{n-1}(X_1\times X_1)} (X_1\times X_1)\;.
\end{split}
\end{equation}
By Corollary \ref{rem-lev-wg-ncat}, $\tens{Z_1}{\di{n-2,j}\p{j,n-2}Z_0}\in\catwg{n-2}$; by Proposition \ref{pro-spec-plbk-pscatwg} and  \eqref{eq1-the-repl-obj-1} we conclude that \eqref{eq0-the-repl-obj-1} holds. It follows that
\begin{equation*}
\begin{split}
    & \p{j+1,n-1}(\tens{P_1}{\di{n-1,j}\p{j,n-1}P_0})= \\
    =\ & \tens{\p{j+1,n-1}P_1}{\di{j}\p{j}\p{j+1,n-1}P_0}\in\tawg{j}\;.
\end{split}
\end{equation*}
Since $P\in\tawg{n}$, $\p{j+2,n}P\in\tawg{j+1}$ so that the induced Segal map
\begin{equation*}
    \p{j+1,n-1}P_2\rw \tens{\p{j+1,n-1}P_1}{(\p{j+1,n-1}P_0)^d}
\end{equation*}
are $j$-equivalences in $\tawg{j}$. From above, this map factorizes as composite of maps in $\tawg{j}$
\begin{equation}\label{eq-ab-comp}
\begin{split}
    & \p{j+1,n-1}P_2 \xrw{\;\za\;} \tens{\p{j+1,n-1}P_1}{\di{j}\p{j}\p{j+1,n-1}P_0} \xrw{\;\zb\;}\\
   \xrw{\;\zb\;}\ &\tens{\p{j+1,n-1}P_1}{(\p{j+1,n-1}P_0)^d}\;.
\end{split}
\end{equation}
We check that the maps \eqref{eq-ab-comp} satisfy the hypotheses of Lemma  \ref{lem-jn-alpha}.
In fact since (as shown in the proof of Proposition \ref{cor-rep-ob}), $\p{n-2}P_{k0}=Z_{k0}$, from Lemma \ref{lem-prop-pn} applied to $Z\up{2}_{0}\in\catwg{n-2}$ we obtain
\begin{align*}
    & \p{j-1}(\p{j+1,n-1}P_2)_0 = \p{j-1,n-2}P_{20}=\p{j-1,n-3}Z_{20}\cong \\
    & \cong \tens{\p{j-1,n-3}Z_{10}}{\p{j-1,n-3}Z_{00}}=\\
    & = \p{j-1}(\tens{\p{j,n-3}Z_{10}}{\p{j,n-3}Z_{00}})=\\
    & = \p{j-1}(\tens{\p{j+1,n-1}P_1}{\p{j+1,n-1}P_0})_0\;.
\end{align*}
That is, $\p{j-1}\za_0$ is an isomorphism. Similarly one shows that $\p{n-r-1}$ $\za_{\seq{k}{1}{r}0}$ is an isomorphism; thus hypothesis i) of Lemma \ref{lem-jn-alpha} holds. As for hypothesis ii) since (as shown in the proof of Proposition \ref{cor-rep-ob}), $P^d_{k0}\cong Z^d_{k0}$, we have
\begin{align*}
  & (\tens{\p{j+1,n-1}P_1}{\p{j+1,n-1}P_0})^d_0\cong \\
  & \cong (\tens{\p{j,n-3}Z_{10}}{\p{j,n-3}Z_{00}})^d \cong \\
  & \cong (\p{j,n-3}Z_{20})^d \cong Z^d_{20} \cong \tens{Z^d_{10}}{Z^d_{00}}\cong \tens{P^d_{10}}{P^d_{00}}\cong\\
  & \cong \tens{(\p{j+1,n-1}P_1)^d_0}{(\p{j+1,n-1}P_0)^d_0}\;.
\end{align*}
The rest of hypothesis ii) of Lemma \ref{lem-jn-alpha} is checked similarly, while hypothesis iii) of Lemma \ref{lem-jn-alpha} holds from above.

We can therefore apply Lemma  \ref{lem-jn-alpha} to \eqref{eq-ab-comp} and conclude that $J_j\za$ is a levelwise equivalence of categories. Therefore $\p{j}\za$ is an isomorphism, that is
\begin{equation*}
\begin{split}
    & \p{j,n-1}P_2 \cong \p{j}(\tens{\p{j+1,n-1}P_1}{\di{j}\p{j}\p{j+1,n-1}P_0})= \\
    =\ & \tens{\p{j,n-1}P_1}{\p{j,n-1}P_0}\;.
\end{split}
\end{equation*}
Similarly one shows that all the other Segal maps for $\p{j+1,n}P$ are isomorphisms, which proves condition c) in Lemma \ref{lem-x-in-tawg-x-in-catwg} for $\p{n}P$. Thus by Lemma \ref{lem-x-in-tawg-x-in-catwg} we conclude that $\p{n}P\in\catwg{n}$, proving b).

\bk
b) By definition of $\lnta{n}{n}$, this follows from Proposition \ref{cor-rep-ob} and a).

\bk

c) Consider the commuting diagram in $\tawg{n}$
\begin{equation*}
\xymatrix@R=35pt @C=40pt{
\di{n}Z \ar^{\di{n}r}[r] \ar_{\di{n}r}[d] & \di{n}\q{n}X \ar@{=}^{}[d] & X \ar^{}[l] \ar@{=}^{}[d] \\
\di{n}\q{n}X \ar@{=}^{}[r] & \di{n}\q{n}X & X  \ar^{}[l]
}
\end{equation*}
By hypothesis, $\di{n}r$ is a $\nequ$. Thus by Proposition \ref{pro-spec-plbk-pscatwg} the induced map of pullbacks
\begin{equation*}
    P=\di{n}Z \tiund{\di{n}\q{n}X} X \xrw{w} \di{n}\q{n}X \tiund{\di{n}\q{n}X} X =X
\end{equation*}
is a $\nequ$, as required.

\bk

d) By b) and by Lemma \ref{lem-lev-wg-pscat}, to show that $P\in\lta{n}$ it is enough to show that, for each $\seqc{k}{1}{n-s}\in\Dop$, $P_{\seqc{k}{1}{n-s}}\in\lnta{s}{s}$, $1< s \leq n-1$. We have a pullback in $\funcat{s} {\Cat}$
\begin{equation}\label{eq2-the-repl-obj-1}
\xymatrix@R=35pt @C=45pt{
P_{\seq{k}{1}{n-s}} \ar^{w_{\seq{k}{1}{n-s}}}[r] \ar^{}[d] & X_{\seq{k}{1}{n-s}} \ar^{}[d] \\
 \di{s} Z_{\seq{k}{1}{n-s}} \ar_{\di{s}r_{\seq{k}{1}{n-s}}}[r] & \di{s}\q{s} X_{\seq{k}{1}{n-s}}
 }
\end{equation}
where $X_{\seq{k}{1}{n-s}}\in \tawg{s}$ (since $X\in\tawg{n}$) and $Z_{\seq{k}{1}{n-s}}\in \catwg{s-1}$ (since $Z\in\catwg{n-1}$).

Thus \eqref{eq2-the-repl-obj-1} satisfies the hypotheses of the theorem and we conclude from b) that
\begin{equation*}
  P_{\seq{k}{1}{n-s}}\in\lnta{s}{s}
\end{equation*}
 as required.

\end{proof}
\begin{corollary}\label{cor-Pcatwg}
Let $X,Z$ be as in the hypothesis of Theorem \ref{the-repl-obj-1} and assume, further, that $X\in\catwg{n}$. Then $P\in\catwg{n}$.
\end{corollary}
\begin{proof}
By induction on $n$. When $n=2$, $P\in\cat{2}$ with $P_0\in\cathd{}$ and $\p{2}P\in\Cat$ (by Theorem \ref{the-repl-obj-1}). Therefore, by Lemma \ref{lem-x-in-tawg-x-in-catwg}, $P\in\catwg{2}$. Suppose, inductively, that the statement holds for $n-1$. We show that $P$ satisfies the hypotheses of Lemma \ref{lem-x-in-tawg-x-in-catwg}. For each $s\geq 0$ we have a pullback in $\funcat{n-2}{\Cat}$
\begin{equation*}
\xymatrix@R=30pt@C=35pt{
P_s \ar[r] \ar[d] & X_s \ar[d]\\
\di{n-1}Z_s \ar[r] & \di{n-1}\q{n-1} X_s
}
\end{equation*}
with $X_s\in\catwg{n-1}$ (since $X\in\catwg{n}$). So by induction hypothesis, $P_s\in\catwg{n-1}$. Thus hypothesis a) in Lemma \ref{lem-x-in-tawg-x-in-catwg} holds for $P$.

Hypothesis b) also holds since $X\in\cat{n}$, $\di{n}\q{n} X\in\cat{n}$, $\di{n}Z\in\cat{n}$ so $P\in \cat{n}$. Finally, hypothesis c) is satisfied because, by  Theorem \ref{the-repl-obj-1}, $\p{n}P\in\catwg{n-1}$ . We conclude by Lemma \ref{lem-x-in-tawg-x-in-catwg} that $P\in\catwg{n}$.

\end{proof}

%%%%%%%%%%%%%%%%%%%%%%%%%%%%%%%%%%%%%%%%%%%%%%%%%%%%%%%%%%%%%%%%%%%%%%%%%%%
\section{From $\pmb{\lta{n}}$ to pseudo-functors}\label{sec-from-pseu}

In this section we show that we can associate functorially to each object of $\lta{n}$ a pseudo-functor which is Segalic. We build in Theorem \ref{the-XXXX} a functor
\begin{equation*}
        \tr{n}: \lta{n} \rw \segpsc{n-1}{\Cat}.
    \end{equation*}
  together with a pseudo-natural transformation
\begin{equation*}
  t_n(X): \tr{n}X\rw X
\end{equation*}
for each $X\in\lta{n}$ which is a levelwise equivalence of categories. The functor $\tr{n}$ will be used in Section \ref{sec-wg-tam-to-psefun} to build the rigidification functor $\Qn$.

\subsection{The idea of the functor $\pmb{\tr{n}}$ }\label{subs-idea-trn}
The idea of the construction of the functor $\tr{n}$ in the proof of Theorem \ref{the-XXXX} is to use the property of objects $X\in\lta{n}$ proved in Proposition \ref{pro-maps-nu-eqcat} that for each $\uk\in\dop{n-1}$, $1\leq i\leq n-1$ and $k_i\geq 2$ there is an equivalence of categories
\begin{equation*}
 \nu(\uk,i):X_{\uk}\rw \pro{X_{\uk(1,i)}}{X^d_{\uk(0,i)}}{k_i}\;.
\end{equation*}
Using this property and working by induction we build a diagram
\begin{equation*}
  \tr{n}X\in [ob(\dop{n-1}),\Cat]
\end{equation*}
in which
 \begin{itemize}
 \item [i)] For all $\uk\in\dop{n-1}$, $1\leq i\leq n-1$
  \begin{equation*}
    (\tr{n}X)_{\uk (0,i)}=X^d_{\uk(0,i)}
  \end{equation*}
   is discrete .\medskip

  \item [ii)] $(\tr{n}X)_{\oset{n-1}{1\,...\,1}}= X_{\oset{n-1}{1\,...\,1}}$.\medskip

  \item [iii)] For $k_1\geq 2$, $\us=(\seqc{k}{2}{n-1})$, $\uk=(k_1,\us)$,
  \begin{equation*}
    (\tr{n}X)_{\uk}= \pro{(\tr{n-1}X_1)_{\us}}{X_{\uk(0,1)}^d}{k_1}
  \end{equation*}

\end{itemize}
For instance, when $n=2$, we set
\begin{equation*}
  (\tr{2}X)_k=
\left\{
  \begin{array}{ll}
    X_0^d, & k=0 \\
    X_1, & k=1 \\
    \pro{X_1}{X_0^d}{k}, & k>1\;.
  \end{array}
\right.
\end{equation*}
When $n=3$ we set
\begin{align*}
    & (\tr{3}X)_{k_1 k_2}= X_{k_1 k_2}^d \quad \mbox{if $k_1=0$ or $k_2=0$} \\
    & (\tr{3}X)_{11}= X_{11}\\
    & (\tr{3}X)_{k_1 1}= \pro{X_{11}}{X_{01}^d}{k_1} \quad \mbox{if $k_1\geq 2$}\\
    & (\tr{3}X)_{1 k_2}= \pro{X_{11}}{X_{10}^d}{k_2} \quad \mbox{if $k_2\geq 2$}\;.
\end{align*}
If both $k_1\geq 2$ and $k_2\geq 2$, we set
\begin{align*}
    & (\tr{3}X)_{k_1 k_2}=\pro{(\tr{2}X_1)_{1 k_2}}{X^d_{0 k_2}}{k_1} =\\
    & = {(\pro{X_{11}}{X_{10}^d}{k_2})}\tiund{(\pro{X_{01}^d}{X_{00}^d}{k_2})}\oset{k_1}{\cdots}\\
   & \oset{k_1}{\cdots} \tiund{(\pro{X_{01}^d}{X_{00}^d}{k_2})} (\pro{X_{11}}{X_{10}^d}{k_2})
\end{align*}
where we used the fact that, since $X_0\in\cathd{2}$,
\begin{equation*}
  X^d_{0 k_2}\cong \pro{X^d_{01}}{X^d_{00}}{k_2}\;.
\end{equation*}
Note also that, since
\begin{equation*}
  X^d_{k_1 0}\cong \pro{X^d_{10}}{X^d_{00}}{k_1}\;.
\end{equation*}
by the commutation of pullbacks we obtain from above
\begin{align*}
    & (\tr{3}X)_{k_1,k_2}=\\
   & =  {(\pro{X_{11}}{X_{01}^d}{k_1})}\tiund{(\pro{X_{10}^d}{X_{00}^d}{k_1})}\oset{k_2}{\cdots}\\
   & \oset{k_2}{\cdots} \tiund{(\pro{X_{10}^d}{X_{00}^d}{k_1})} (\pro{X_{11}}{X_{01}^d}{k_1})= \\
   & = \pro{(\tr{3}X)_{k_1 1}}{X^d_{k_1 0}}{k_2}\;.
\end{align*}

\bk

After defining $\tr{n}X\in [ob(\dop{n-1}),\Cat]$ we show that, for each $\uk\in\dop{n-1}$, there is an equivalence of categories
\begin{equation}\label{eq01-the-XXXX}
  (\tr{n}X)_{\uk}\simeq X_{\uk}\;.
\end{equation}
Using the 'transport of structure' technique of Lemma \ref{lem-PP}, we then lift $\tr{n}X$ to a pseudo-functor
\begin{equation*}
  \tr{n}X\in\psc{n-1}{\Cat}
\end{equation*}
and we show that this is in fact a Segalic pseudo-functor. Conditions a) and b) in the definition of Segalic pseudo-functor depend on the conditions i), ii), iii) in the definition of $(\tr{n}X)_{\uk}$, while condition c) is a straightforward consequence of the equivalence of categories \eqref{eq01-the-XXXX} and the fact that, since $X\in\lta{n}$, $\p{n}X\in\catwg{n-1}$.
\subsection{The formal construction of the functor $\pmb{\tr{n}}$ }
\begin{theorem}\label{the-XXXX}
    There is a functor
    \begin{equation*}
        Tr_{n}: \lta{n} \rw \segpsc{n-1}{\Cat}
    \end{equation*}
    together with a pseudo-natural transformation
     \begin{equation*}
       t_n (X):Tr_{n}X\rw X
     \end{equation*}
      for each $X\in\lta{n}$ which is a levelwise equivalence of categories.
\end{theorem}
\begin{proof}
By induction on $n$. For $n=2$, let $X\in\lta{2}=\tawg{2}$. Define $\tr{2}X\in[ob(\Dop),\Cat]$
\begin{equation}\label{eq1-the-XXXX}
(\tr{2}X)_k=
\left\{
  \begin{array}{ll}
    X_0^d & k=0 \\
    X_1 & k=1 \\
    \pro{X_1}{X_0}{k} & k>1\;.
  \end{array}
\right.
\end{equation}
Since $X\in\tawg{2}$, $X_0\in\cathd{}$ so there are equivalences of categories
\begin{equation*}
\begin{split}
    & X_0\simeq X_0^d \\
    & X_k\simeq \pro{X_1}{X_0^d}{k}\quad \text{for }\; k>1.
\end{split}
\end{equation*}
Thus, for all $k\geq 0$ there is an equivalence of categories
\begin{equation*}
 (\tr{2}X)_k\simeq X_k\;.
\end{equation*}
We can therefore apply Lemma \ref{lem-PP} with $\clC=\Dop$ and conclude that $\tr{2}X$ lifts to a pseudo-functor
\begin{equation*}
    \tr{2}X \in \psc{}{\Cat}
\end{equation*}
and there is a pseudo-natural transformation
\begin{equation*}
  t_2(X):\tr{2}X\rw X
\end{equation*}
which is a levelwise equivalence of categories. By \eqref{eq1-the-XXXX}, $(\tr{2}X)_0$ is discrete and the Segal maps are isomorphisms. Therefore, by Definition \ref{def-seg-ps-fun},
\begin{equation*}
  \tr{2}X\in\segpsc{}{\Cat}\;.
\end{equation*}
Suppose, inductively, that the theorem holds for $(n-1)$ and let $X\in\lta{n}$. By definition of $\lta{n}$ (see Remark \ref{rem-def-ltawg}), for each $\us\in\dop{n-2}$, $j\geq 2$ there is an equivalence of categories
\begin{equation}\label{eq2-the-XXXX}
X_{j\us}\simeq \pro{X_{1\us}}{X_{0\us}^d}{j}\;.
\end{equation}
Also, by inductive hypothesis applied to $X_j\in\lta{n-1}$ there is an equivalence of categories for all $j\geq 0$ and $\us\in\dop{n-1}$
\begin{equation}\label{eq3-the-XXXX}
  X_{j\us}\simeq (\tr{n-1}X_j)_{\us}\;.
\end{equation}
It follows from \eqref{eq2-the-XXXX} that for each $j\geq 2$ there is an equivalence of categories
\begin{equation}\label{eq4-the-XXXX}
  \pro{X_{1\us}}{X_{0\us}^d}{j}\simeq \pro{(\tr{n-1}X_1)_{\us}}{X_{0\us}^d}{j}\;.
\end{equation}
Thus \eqref{eq2-the-XXXX}, \eqref{eq3-the-XXXX}, \eqref{eq4-the-XXXX} imply that for each $j\geq 2$, $\us\in\dop{2}$ there is an equivalence of categories
\begin{equation}\label{eq5-the-XXXX}
  (\tr{n-1}X_j)_{\us}\simeq \pro{(\tr{n-1}X_1)_{\us}}{X_{0\us}^d}{j}\;.
\end{equation}
Define $\tr{n}X\in [ob(\dop{n-1}),\Cat]$ as follows: for each $\uk=(k_1,\us)\in\dop{n-1}$ (with $k_1\in\Dop$, $\us\in\dop{n-2}$)
\begin{equation}\label{eq6-the-XXXX}
\resizebox{.87\hsize}{!}{$   %%%%%%%%%%%% To reduce the size to 90% of normal
  (\tr{n}X)_{\uk}=
\left\{
\begin{array}{ll}
    X_{\uk}^d & \mbox{for $\uk=(0,\us)$} \\
    (\tr{n-1}X_1)_{\us} & \mbox{for $\uk=(1,\us)$} \\
    \pro{(\tr{n-1}X_1)_{\us}}{X_{0\us}^d}{k_1} & \mbox{for $\uk=(k_1,\us)$},\\
   &  k_1\geq 2.
  \end{array}
\right.$}
\end{equation}
We claim that there is an equivalence of categories for all $\uk\in\dop{n-1}$
\begin{equation}\label{eq7-the-XXXX}
  (\tr{n}X)_{\uk}\simeq X_{\uk}\;.
\end{equation}
In fact, since $X_{0\us}\in\cathd{}$,
\begin{equation*}
  (\tr{n}X)_{0\us}=X_{0\us}^d\simeq X_{0\us}\;.
\end{equation*}
By inductive hypothesis applied to $\tr{n-1}X_1$,
\begin{equation*}
 (\tr{n}X)_{1\us}=(\tr{n-1}X_1)_{\us}\simeq X_{1\us}\;.
\end{equation*}
This implies, when $k_1\geq 2$
\begin{equation*}
  (\tr{n}X)_{\uk}=\pro{(\tr{n-1}X_1)_{\us}}{X_{0\us}^d}{k_1}\simeq \pro{X_{1\us}}{X_{0\us}^d}{k_1}
\end{equation*}
and together with \eqref{eq2-the-XXXX} it follows that
\begin{equation*}
  (\tr{n}X)_{\uk} \simeq X_{\uk}
\end{equation*}
when $k_1\geq 2$. This concludes the proof of \eqref{eq7-the-XXXX}.

 We can therefore apply Lemma \ref{lem-PP} with $\clC=\dop{n-1}$ and conclude that $\tr{n}X$ lifts to a pseudo-functor
\begin{equation*}
  \tr{n}X\in\psc{n-1}{\Cat}
\end{equation*}
with $(\tr{n}X)_{\uk}$ as in \eqref{eq6-the-XXXX}.

 We now show that $\tr{n}X$ is a Segalic pseudo-functor, by checking the conditions in Definition \ref{def-seg-ps-fun}.

We first check condition a) that $(\tr{n}X)_{\uk(0,i)}$ is discrete for all $1\leq i\leq n-1$. By construction \eqref{eq6-the-XXXX}

\begin{equation*}
  (\tr{n}X)_{\uk(0,i)}=X^d_{\uk(0,i)}
\end{equation*}
is discrete. Also by construction and by inductive hypothesis, when $k_1=1$ and $i>1$,
\begin{equation*}
  (\tr{n}X)_{\uk(0,i)}=(\tr{n-1}X_1)_{\us(0,i-1)}
\end{equation*}
is discrete.
Finally, if $k_1>1$ and $i>1$, by construction and by the inductive hypothesis
\begin{equation*}
   (\tr{n}X)_{\uk(0,i)}=\pro{(\tr{n-1}X_1)_{\us(0,i-1)}}{X_{0\us}^d}{k_1}
\end{equation*}
is discrete. This shows that condition a) in Definition \ref{def-seg-ps-fun} is satisfied. We now show condition b) that for each $\uk\in\dop{n-1}$, $1\leq i\leq n-1$ and $k_i\geq 2$,
\begin{equation}\label{eq8-the-XXXX}
  (\tr{n}X)_{\uk} \cong \pro{(\tr{n}X)_{\uk(1,i)}}{(\tr{n}X)_{\uk(0,i)}}{k_i}\;.
\end{equation}
We distinguish various cases:\bk

 i) Let $\uk\in\dop{n-1}$ be such that $k_j=0$ for some $1\leq j\leq n-1$ and let $k_i\geq 2$ for $1\leq i\leq n$. Since, by Proposition \ref{pro-maps-nu-eqcat} there is an equivalence of categories
\begin{equation*}
  X_{\uk}\simeq \pro{X_{\uk(1,i)}}{X^d_{\uk(0,i)}}{k_i}
\end{equation*}
and since $X_{\uk}\in\cathd{}$, $X_{\uk(1,i)}\in \cathd{}$ (as $k_j=0$ and $k_i\geq 2$ so $i\neq j$), there is an isomorphism
\begin{equation}\label{eq9-the-XXXX}
\begin{split}
    & X_{\uk}^d \cong dp(\pro{X_{\uk(1,i)}}{X^d_{\uk(0,i)}}{k_i}) \cong \\
    & \cong \pro{X^d_{\uk(1,i)}}{X^d_{\uk(0,i)}}{k_i}
\end{split}
\end{equation}
which, by \eqref{eq6-the-XXXX}, is the same as \eqref{eq8-the-XXXX} in this case.\bk

 ii) Suppose $k_j\neq 0$ for all $1\leq j\leq n-1$ and let $i=1$. Then by \eqref{eq6-the-XXXX}, if $\us=(\seqc{k}{2}{n-1})$
\begin{equation*}
 \begin{split}
     & (\tr{n}X)_{\uk(1,1)} = (\tr{n-1}X_1)_{\us}\\
     & (\tr{n}X)_{\uk(0,1)} = X^d_{\uk(0,1)}=X_{0\us}^d\;.
 \end{split}
\end{equation*}
 Therefore by \eqref{eq6-the-XXXX}, if $k_1\geq 2$,

\begin{align*}
    & (\tr{n}X)_{\uk}=\pro{(\tr{n-1}X_1)_{\us}}{X_{0\us}^d}{k_1}= \\
  = & \pro{(\tr{n}X)_{\uk(1,1)}}{(\tr{n}X)_{\uk(0,1)}}{k_1}
\end{align*}
which is \eqref{eq8-the-XXXX} in this case.\bk

 iii) Suppose $k_j\neq 0$ for all $1\leq j\leq n-1$, $i>1$ and $k_1=1$. Then by  \eqref{eq6-the-XXXX}, if $\us=(\seqc{k}{2}{n-1})$
\begin{equation*}
  (\tr{n}X)_{\uk}=(\tr{n-1}X_1)_{\us}
\end{equation*}
so in particular
\begin{equation*}
\begin{split}
    & (\tr{n}X)_{\uk(1,i)}=(\tr{n-1}X_1)_{\us(1,i-1)} \\
    & (\tr{n}X)_{\uk(0,i)}=(\tr{n-1}X_1)_{\us(0,i-1)}\;.
\end{split}
\end{equation*}
By induction hypothesis applied to $X_1$, it follows that, since $k_i=s_{i-1}\geq 2$,
\begin{align*}
    & (\tr{n}X)_{\uk}=(\tr{n-1}X_1)_{\us}= \\
  = & \pro{(\tr{n-1}X_1)_{\us(1,i-1)}}{(\tr{n-1}X_1)_{\us(0,i-1)}}{s_{i-1}}=\\
  = & \pro{(\tr{n}X)_{\uk(1,i)}}{(\tr{n}X)_{\uk(0,i)}}{k_i}
\end{align*}
which is \eqref{eq8-the-XXXX} in this case.\bk

 iv) Suppose $k_j\neq 0$ for all $1\leq j\leq n-1$, $i>1$ and $k_1=2$. By  \eqref{eq6-the-XXXX}, if $\us=(\seqc{k}{2}{n-1})$ so that $k_i=s_{i-1}$
\begin{equation}\label{eq10-the-XXXX}
  (\tr{n}X)_{\uk}=\tens{(\tr{n-1}X_1)_{\us}}{X_{0\us}^d}\;.
\end{equation}
By induction hypothesis applied to $X_1$,
\begin{equation}\label{eq11-the-XXXX}
\begin{split}
    & (\tr{n-1}X_1)_{\us}= \\
   = & \pro{(\tr{n-1}X_1)_{\us (1,i-1)}}{(\tr{n-1}X_1)_{\us (0,i-1)}}{s_{i-1}}
\end{split}
\end{equation}
while by \eqref{eq9-the-XXXX}
\begin{equation}\label{eq12-the-XXXX}
  X_{0\us}^d=\pro{X_{0\us(1,i-1)}^d}{X_{0\us(0,i-1)}^d}{s_{i-1}}\;.
\end{equation}
Replacing \eqref{eq11-the-XXXX} and \eqref{eq12-the-XXXX} in \eqref{eq10-the-XXXX}, using the commutation of limits and the fact that
\begin{align*}
& (\tr{n}X)_{\uk(1,i)}=\tens{(\tr{n-1}X_1)_{\us (1,i-1)}}{X_{0\us(1,i-1)}^d}\\
& (\tr{n}X)_{\uk(0,i)}=\tens{(\tr{n-1}X_1)_{\us (0,i-1)}}{X_{0\us(0,i-1)}^d}
\end{align*}
we obtain
\begin{equation*}
 (\tr{n}X)_{\uk}=\pro{(\tr{n}X)_{\uk(1,i)}}{(\tr{n}X)_{\uk(0,i)}}{k_i}
\end{equation*}
which is \eqref{eq8-the-XXXX} in this case.\bk

v) Suppose $k_j\neq 0$ for all $1\leq j\leq n-1$, $i>1$ and $k_1>2$. The proof of \eqref{eq8-the-XXXX} is completely analogous to the one of case iv).\bk

This concludes the proof that condition b) in Definition \ref{def-seg-ps-fun} is satisfied for $\tr{n}$.

 To show that condition c) in Definition \ref{def-seg-ps-fun} holds for $\tr{n}$ we note that the equivalence of categories \eqref{eq7-the-XXXX} implies the isomorphism
\begin{equation*}
  p(\tr{n}X)_{\uk}\cong p X_{\uk}= (\p{n}X)_{\uk}\;.
\end{equation*}
Since $X\in\lta{n}$, $\p{n}X\in\catwg{n-1}$, hence
\begin{equation*}
  \p{n}\tr{n}X \in \catwg{n-1}
\end{equation*}
which is condition c) in Definition \ref{def-seg-ps-fun}. We conclude that
\begin{equation*}
  \tr{n}X \in \segpsc{n-1}{\Cat}\;.
\end{equation*}
By Lemma \ref{lem-PP} there is a morphism in $\psc{n-1}{\Cat}$
\begin{equation*}
  t_n(X): \tr{n}X \rw X
\end{equation*}
which is levelwise the equivalence of categories
\begin{equation*}
  (\tr{n}X)_{\uk} \simeq X_{\uk}\quad \text{for}\;\; \uk\in\dop{n-1}\;.
\end{equation*}

\end{proof}
\begin{corollary}\label{cor-trn}
Let $X \in\lta{n}$ and $\uk\in\dop{n-1}$ be such that $k_j=0$ for some $1\leq j\leq n-1$. Then
\begin{equation*}
  (\tr{n} X)_{\uk}= X_{\uk}^d\;.
\end{equation*}
\end{corollary}
\begin{proof}
By induction on $n$. When $n=2$, by definition of $\tr{2}X$,
\begin{equation*}
  (\tr{2}X)_0=X_0^d\;.
\end{equation*}
Suppose the statement hold for $(n-1)$ and denote $\us=(\seqc{k}{2}{n-1})\in \dop{n-2}$ so that
\begin{equation*}
  \uk=(k_1,\us)\in\dop{n-1}\;.
\end{equation*}
We distinguish three cases:
\begin{itemize}
  \item [i)] When $k_1=0$, by definition of $\tr{n}$
  \begin{equation*}
    (\tr{n}X)_{\uk}=(\tr{n}X)_{(0,\us)}=X_{0\us}^d=X_{\uk}^d\;.
  \end{equation*}

  \item [ii)] Let $k_1=1$ and suppose $s_j=0$ for some $2\leq j\leq n-1$. Then by definition of $\tr{n}X$ and by inductive hypothesis applied to $X_1\in\lta{n-1}$,
      \begin{equation*}
        (\tr{n}X)_{\uk}=(\tr{n-1}X_1)_{\us}=(X_1)_{\us}^d=X_{\uk}^d\;.
      \end{equation*}

  \item [iii)] Let $k_1>1$ and suppose $s_j=0$ for some  $2\leq j\leq n-1$. Then by definition of $\tr{n}X$ and by inductive hypothesis applied to $X_1\in\lta{n-1}$, we have
      \begin{equation}\label{eq1-cor-trn}
      \begin{split}
         & (\tr{n}X)_{\uk}=\pro{(\tr{n-1}X_1)_{\us}}{X_{0\us}^d}{k_1}= \\
          & =\pro{X_{1\us}}{X_{0\us}^d}{k_1}\;.
      \end{split}
      \end{equation}
      Since $X\in\lta{n}$, by Remark \ref{rem-def-ltawg} there is an equivalence of categories
      \begin{equation*}
        X_{\uk}=X_{k_1\us}\simeq \pro{X_{1\us}}{X_{0\us}^d}{k_1}
      \end{equation*}
      and therefore, since $X_{\uk}, X_{1\us}\in\cathd{}$,
      \begin{equation}\label{eq2-cor-trn}
        X_{\uk}^d\cong \pro{X^d_{1\us}}{X^d_{0\us}}{k_1}\;.
      \end{equation}
      We deduce from \eqref{eq1-cor-trn} and \eqref{eq2-cor-trn} that
      \begin{equation*}
        (\tr{n}X)_{\uk}=X_{\uk}^d\;.
      \end{equation*}

\end{itemize}

\end{proof}

\begin{example}\label{ex-from-lta-pseudo}
The functor
\begin{equation*}
  \tr{3}:\lta{3}\rw\segpsc{2}{\Cat}
\end{equation*}
is given as follows: $X\in\lta{3}$ consists of $X\in\funcat{2}{\Cat}$ such that
\begin{itemize}
  \item [i)] $X_k\in\tawg{2}$ for all $k\geq 0$ with $X_0\in\cathd{2}$. Thus in particular $X_{k0}\in \cathd{}$ and the induced Segal maps
\begin{equation*}
  X_{ks}\rw \pro{X_{k 1}}{X^d_{k 0}}{s}
\end{equation*}
are equivalences of categories for all $s\geq 2$.\bk

  \item [ii)] $\p{3}X\in\catwg{2}$.\bk

  \item [iii)] For each $s\geq 0$ and $k\geq 2$ the maps
\begin{equation*}
  X_{ks}\rw \pro{X_{1s}}{X^d_{0s}}{k}
\end{equation*}
are equivalences of categories.
\end{itemize}
Below is a picture of the corner of $\tr{3}X$, where the symbol $\cong$ indicates that the squares pseudo-commute, that is $\tr{3}X$ is not a bisimplicial object in $\Cat$ but a pseudo-functor from $\dop{2}$ to $\Cat$.

\begin{equation*}
\def\objectstyle{\scriptstyle}
\def\labelstyle{\scriptstyle}
%\entrymodifiers={++[o]}
\ssr \xymatrix@R=35pt@C=20pt{
   \bm{\tens{(\tens{X_{11}}{X^d_{01}})}{(\tens{X^d_{10}}{X^d_{00}})}} \ar@<1ex>[r] \ar[r] \ar@<-1ex>[r] \ar@<1ex>[d] \ar[d] \ar@<-1ex>[d] \ar@{} [dr] |{\cong}
   &  \bm{(\tens{X_{11}}{X^d_{10}})} \ar@<-0.5ex>[r] \ar@<0.5ex>[r] \ar@<1ex>[d] \ar[d] \ar@<-1ex>[d] \ar@{} [dr] |{\cong}
   &  \bm{(\tens{X^d_{01}}{X^d_{00}})}  \ar@<1ex>[d] \ar[d] \ar@<-1ex>[d] \\
     \bm{\cdots\quad(\tens{X_{11}}{X^d_{01}})\quad}\ar[r] \ar@<1ex>[r]  \ar@<-1ex>[r] \ar@<0.5ex>[d] \ar@<-0.5ex>[d] \ar@{} [dr] |{\cong}
    &  \bm{\quad  X_{11}\quad }\ar@<-0.5ex>[r] \ar@<0.5ex>[r]   \ar@<0.5ex>[d] \ar@<-0.5ex>[d] \ar@{} [dr] |{\cong}
    &  \bm{\quad X^d_{01}\quad}  \ar@<0.5ex>[d] \ar@<-0.5ex>[d]  \\
      \bm{\cdots\quad(\tens{X^d_{10}}{X^d_{00}})\quad} \ar[r] \ar@<1ex>[r]  \ar@<-1ex>[r]
    &  \bm{\quad X^d_{10}\quad} \ar@<-0.5ex>[r] \ar@<0.5ex>[r]
    &  \bm{\quad X^d_{00}\quad} \\
}
\end{equation*}

\bk
Note that
\begin{align*}
    & \tens{(\tens{X_{11}}{X^d_{01}})}{(\tens{X^d_{10}}{X^d_{00}})}\cong \\
    & \cong\tens{(\tens{X_{11}}{X^d_{10}})}{(\tens{X^d_{10}}{X^d_{00}})}
\end{align*}
so that the Segal maps of $\tr{3}X$ in both horizontal and vertical directions are isomorphisms.

\end{example}

\bk
The following Lemma will be used in the proof of Proposition \ref{pro-fta-tam-2}. The latter will be crucial in proving the properties of the discretization functor in the proof of Theorem \ref{the-disc-func}.

\begin{lemma}\label{lem-from-lta-to-pseu}
Let $X\in\catwg{n}$, $Y\in\lta{n}$ be such that $Y_{\uk}$ is discrete for all $\uk\in\dop{n-1}$ such that $k_j=0$ for some $1\leq j\leq n-1$. Let $\uk,\;\us \in\dop{n-1}$ and let $\uk \rw\us$ be a morphism in $\dop{n-1}$. Suppose that the following conditions hold:
\begin{itemize}
  \item [i)] If $k_j,s_j \neq 0$ for all $1\leq j\leq n-1$, then $ X_{\uk}=Y_{\uk},\quad X_{\us}=Y_{\us}$ and the maps
   \begin{equation*}
   X_{\uk}\rw X_{\us},\quad Y_{\uk}\rw Y_{\us}\;
   \end{equation*}
    coincide.

  \bk
  \item [ii)] If $k_j=0$ for some $1\leq j\leq n-1$ and $s_t=0$ for some $1\leq t\leq n-1$, then
  $X^d_{\uk}=Y_{\uk},\quad X^d_{\us}=Y_{\us}$ and the two maps
    \begin{equation*}
    X^d_{\uk}\rw X^d_{\us},\quad Y_{\uk}\rw Y_{\us}
    \end{equation*}
coincide, where $f^d:X^d_{\uk}\rw X^d_{\us}$ is induced by $f:X_{\uk}\rw X_{\us}$ and thus also coincides with the composite
  \begin{equation*}
    X^d_{\uk}\xrw{\zg'_{X_k}} X_{\uk} \xrw{f} X_{\us}\xrw{\zg_{X_s}}X^d_{\us}
  \end{equation*}
  (where $\zg$ is the discretization map and $\zg'$ a section), since $f^d=f^d \zg_{X_k}\zg'_{X_k}=\zg_{X_s}f\zg'_{X_k}$.
  \bk

  \item [iii)] If $k_j\neq 0$ for all $1\leq j\leq n-1$ and $s_t=0$ for some $1\leq t\leq n-1$, the following diagram commutes
  \begin{equation*}
  \qquad\qquad\qquad\xymatrix{
  X_{\uk} \ar[r] \ar@{=}[d] & X_{\us} \ar^(0.3){\zg_{X_{\us}}}[rr] && X^d_{\us}=Y_{\us}\mbox{\qquad\qquad\qquad\qquad\qquad} \\
  Y_{\uk}\ar[rrru]
  }
  \end{equation*}
  where $\zg_{X_s}$ is the discretization map.

  \bk
  \item [iv)] If $k_j=0$ for some $1\leq j\leq n-1$ and $s_t\neq 0$ for all $1\leq t\leq n-1$ then the following diagram commutes
  \begin{equation*}
  \qquad\qquad\qquad\xymatrix{
  X^d_{\uk} \ar^{\zg'_{X_{\uk}}}[rr] \ar@{=}[d] && X_{\uk} \ar[r] & X_{\us}=Y_{\us}\mbox{\qquad\qquad\qquad\qquad} \\
  Y_{\uk}\ar[rrru]
  }
  \end{equation*}
  Then
  \begin{itemize}
    \item [a)] For all $\uk\in\dop{n-1}$, $(\tr{n}X)_{\uk}=(\tr{n}Y)_{\uk}$\bk

    \item [b)] For all $\uk\in\dop{n-1}$ such that $k_j\neq 0$ for all $1\leq j\leq n-1$, the maps
    \begin{equation*}
      (\tr{n}X)_{\uk}\rightleftarrows X_{\uk}, \qquad (\tr{n}Y)_{\uk}\rightleftarrows Y_{\uk}
    \end{equation*}
    coincide. \bk

    \item [c)] $\tr{n}X=\tr{n}Y$\;.
  \end{itemize}

\end{itemize}

\end{lemma}
\begin{proof}
See Appendix \ref{app-a} on page \pageref{page-app-a}.
\end{proof}
%%

%%%%%%%%%%%%%%%%%%%%%%%%%%%%%%%%%%%%%%%%%%%%%%%%%%%%%%%%%%%%%%%%%%%%%%%%%%%

%%
\section{Rigidifying weakly globular Tamsamani ${\pmb{n}}$-categories}\label{sec-wg-tam-to-psefun}
In this section we prove the main result of the chapter, Theorem \ref{the-funct-Qn}, establishing the existence of a rigidification functor
\begin{equation*}
  Q_n:\tawg{n}\rw \catwg{n}
\end{equation*}
replacing $X\in\tawg{n}$ with an $n$-equivalent object $Q_n X$.

 \subsection{The rigidification functor $\pmb{Q_n}$: main steps}\label{sub-idea-qn}
 The construction of the functor $Q_n$ is inductive and uses three main ingredients:

  \begin{itemize}
 \item [a)] The approximation up to $n$-equivalence of an object of $\tawg{n}$ with an object of $\lta{n}$ using the pullback construction of  Theorem \ref{the-repl-obj-1}. We showed in Theorem \ref{the-repl-obj-1} that if  $X\in\tawg{n}$ is such that $\qn X$ can be approximated up to $(n-1)$-equivalence with an object of $\catwg{(n-1)}$, then $X$ can be approximated up to an $n$-equivalence with an object of $\lta{n}$. Thus, given inductively the functor $Q_{(n-1)}$, for each $X\in\tawg{n}$ we can approximate $\qn X$ with $Q_{(n-1)}\qn X \in \catwg{(n-1)}$; thus by the above we can approximate $X$ with an object $P_n X\in\lta{n}$ and we obtain a functor
        \begin{equation*}
          P_n: \tawg{n}\rw \lta{n}.
        \end{equation*}
      \medskip
    \item [b)] The functor $Tr_n$  from the category $\lta{n}$ to the category of Segalic pseudo-functors, which we established in Theorem \ref{the-XXXX}.
      \medskip
    \item[c)] The functor $\St$ from Segalic pseudo-functors to weakly globular $n$-fold categories from Theorem \ref{the-strict-funct}.
  \end{itemize}
We then define the rigidification functor $Q_2$ to be
\begin{equation*}
    Q_2: \tawg{2}\xrw{Tr_{2}} \segpsc{}{\Cat} \xrw{\St} \catwg{2}
\end{equation*}
 The rigidification functor $Q_n$, when $n>2$ is defined as the composite
\begin{equation*}
    \tawg{n}\xrw{P_n} \lta{n} \xrw{\Tr_{n}} \segpsc{n-1}{\Cat}\xrw{\St} \catwg{n}.
\end{equation*}

\bigskip

\subsection{The rigidification functor: the formal proof}
\begin{theorem}\label{the-funct-Qn}
\index{Rigidification functor}

    There is a functor, called rigidification,
    \begin{equation*}
        Q_n:\tawg{n} \rw \catwg{n}
    \end{equation*}
    and for each $X\in\tawg{n}$ a morphism in $\tawg{n}$
     \begin{equation*}
       s_n(X):Q_n X\rw X
     \end{equation*}
      natural in $X$, such that $(s_n(X))_{k}$ is a $(n-1)$-equivalence for all $k\geq0$. In particular, $s_{n}(X)$ is an $n$-equivalence.
\end{theorem}
\begin{proof}
By induction on $n$. When $n=2$, let $Q_2$ be the composite
\begin{equation*}
    Q_2: \tawg{2}\xrw{Tr_{2}} \segpsc{}{\Cat} \xrw{\St} \catwg{2}
\end{equation*}
where $Tr_{2}$ is as in Theorem \ref{the-XXXX}. By Theorem \ref{the-strict-funct}, $Q_2X \in \catwg{2}$. Recall \cite{Lack} that the strictification functor
\begin{equation*}
    \St:\psc{}{\Cat} \rw \funcat{}{\Cat}
\end{equation*}
is left adjoint to the inclusion
\begin{equation*}
    J: \funcat{}{\Cat} \rw \psc{}{\Cat}
\end{equation*}
and that the components of the unit are equivalences in $\psc{}{\Cat}$. By Theorem \ref{the-XXXX} there is a morphism in $\psc{}{\Cat}$
\begin{equation*}
    t_2(X):Tr_{2}X \rw JX \;.
\end{equation*}
By adjunction this corresponds to a morphism in $\funcat{}{\Cat}$
\begin{equation*}
    Q_2 X=\St Tr_{2} X\xrw{s_2(X)} X
\end{equation*}
making the following diagram commute
\begin{equation*}
    \xymatrix@C=50pt@R=50pt{
    Tr_{2}X \ar^(0.4){\eta}[r] \ar_{t_2(X)}[dr] & J\St Tr_{2}X
    \ar^{J s_2(X)}[d] \\
    &  JX
    }
\end{equation*}
Since $\eta$ and $t_2(X)$ are levelwise equivalences of categories, such is $J s_2(X)$.

Suppose, inductively, that we defined $Q_{n-1}$. Define the functor
\begin{equation*}
    P_n:\tawg{n}\rw\lta{n}
\end{equation*}
as follows. Given $X\in\tawg{n}$, consider  the pullback in $\funcat{n-1}{\Cat}$
\begin{equation*}
    \xymatrix@C=70pt@R=30pt{
    P_n X \ar^{w_{n}(X)}[r] \ar[d] & X \ar^{\zg_{n}}[d] \\
    \di{n}Q_{n-1}\q{n}X \ar_{\di{n}s_{n-1}(\q{n}X)}[r] & \di{n}\q{n}X
    }
\end{equation*}
By Theorem \ref{the-repl-obj-1}, $P_n X\in\lta{n}$. Define
\begin{equation*}
Q_n X=\St Tr_{n} P_n X.
\end{equation*}
By Theorem \ref{the-strict-funct}, $Q_nX \in \catwg{n}$. Let $s_n(X):Q_n X\rw X$ be the composite
\begin{equation*}
s_n(X): Q_n X \xrw{h_n(P_n X)} P_n X \xrw{w_{n}(X)} X\;.
\end{equation*}
 where  the morphism in $\funcat{n-1}{\Cat}$
\begin{equation*}
    Q_n X= \St Tr_{n}P_n X \xrw{h_n(P_n X)} P_n X
\end{equation*}
 corresponds by adjunction to the morphism in $\psc{n-1}{\Cat}$
\begin{equation*}
    Tr_{n}P_n X \xrw{t_n(P_n X)} J P_n X
\end{equation*}
(where $t_n(P_n X)$ is as in Theorem \ref{the-XXXX}) such that the following diagram commutes
\begin{equation*}
    \xymatrix@C=50pt@R=50pt{
    Tr_{n}P_n X \ar^(0.35){\eta}[r] \ar_{t_n(P_n X)}[dr] & J \St Tr_{n}P_n
    X= J Q_n X \ar^{J h_n(P_n X)}[d]\\
    & J P_n X
    }
\end{equation*}
We need to show that $(s_n(X))_{k}$ is an $(n-1)$-equivalence. Since $\eta$ and $t_n(P_n X)$ are levelwise equivalences of categories, such is $h_n(P_n X)$ so in particular $(h_n(P_n X))_{k}$ is a levelwise equivalence of categories, and thus is a $(n-1)$-equivalence (see Remark \ref{rem-local-equiv}).

Since pullbacks in $\funcat{n-1}{\Cat}$ are computed pointwise, there is a pullback in $\funcat{n-2}{\Cat}$
\begin{equation*}
\xymatrix@C=90pt@R=40pt{
(P_n X)_k \ar^{(w_{n}(X))_k}[r] \ar[d] & X_k \ar[d]\\
\di{n-1}(Q_{n-1}\q{n}X)_k \ar_(0.55){\di{n-1}(s_{n-1}(\q{n}X))_k}[r] & \di{n-1}\q{n-1}X_k
}
\end{equation*}
where $X_k\in\tawg{n-1}$ (since $X\in\tawg{n}$) and
 \begin{equation*}
 (Q_{n-1}\q{n}X)_k\in\catwg{n-2}
 \end{equation*}
  (since $Q_{n-1}\q{n}X\in\catwg{n-1}$) and, by induction hypothesis, $(s_{n-1}(\q{n}X))_k$ is a $(n-2)$-equivalence. It follows by Theorem \ref{the-repl-obj-1} that $(w_{n}(X))_k$ is a $(n-1)$-equivalence.

In conclusion both $(h_n(P_n X))_k$ and $(w_{n}(X))_k$ are $(n-1)$-equivalences so by Proposition \ref{pro-n-equiv} such is their composite
\begin{equation*}
    (s_n(X))_k : (Q_n X)_k \xrw{(h_n(P_n X))_k} (P_n X)_k \xrw{(w_{n}(X))_k} X_k
\end{equation*}
as required. By Lemma \ref{lem-flevel-fneq}, it follows that $s_n(X)$ is an $n$-equivalence.

\end{proof}
\begin{corollary}\label{cor-the-funct-Qn}
    The functors $Q_n:\tawg{n}\rw \catwg{n}$ and the embedding $i:\catwg{n}\hookrightarrow \tawg{n}$ induce an equivalence of categories
    \begin{equation}\label{eq-cor-the-funct-Qn}
        \tawg{n}\bsim^n\; \simeq\; \catwg{n}\bsim^n
    \end{equation}
    after localization with respect to the $n$-equivalences.
\end{corollary}
\begin{proof}
Given $X\in\tawg{n}$, by Theorem \ref{the-funct-Qn} there is an $n$-equivalence in $\tawg{n}$ $i Q_n X\rw X$, therefore
 \begin{equation*}
 i Q_n X\cong X
 \end{equation*}
  in $\tawg{n}\bsim^n$.

Let $Y\in\catwg{n}$; then $i Y\in\tawg{n}$ so by Theorem \ref{the-funct-Qn} there is an $n$-equivalence in $\tawg{n}$ $i Q_n i Y \rw i Y$. Since $i$ is fully faithful, $ Q_n i Y \rw Y$ is an $n$-equivalence in $\catwg{n}$. It follows that
 \begin{equation*}
 Q_n i Y \cong Y
 \end{equation*}
  in $\catwg{n}\bsim^n$. In conclusion $Q_n$ and $i$ induce the equivalence of categories \eqref{eq-cor-the-funct-Qn}.
\end{proof}

\begin{remark}\label{rem-pn-catwg}
  It follows from Corollary \ref{cor-Pcatwg} that given $X\in\catwg{n}$, $n>2$, $P_n X \in \catwg{n}$.
\end{remark}
%%%%%%%%%%%%%%%%%%%
\clearpage

\begin{figure}
\vspace{10mm}
  \centering
  \begin{tikzpicture}[thick, mark options={color=black}] %% fill opacity=0.6
\filldraw[fill=black!6!white, draw=white][thick] (3,3) -- (1,1) -- (1,-2) -- (3,0) -- (3,3);
\filldraw[fill=black!30!white, draw=black][thick] (0,3) -- (3,3) -- (1,1) -- (-2,1) -- (0,3);
\filldraw[fill=black!30!white, draw=black][thick] (3,0) -- (6,0) -- (4,-2) -- (1,-2) -- (3,0);
\filldraw[fill=black!6!white, draw=white][thick] (4,-2) -- (4,5) -- (6,7) -- (6,0) -- (4,-2);
\filldraw[fill=black!30!white, draw=black][thick] (6,7) -- (9,7) -- (7,5) -- (4,5) -- (6,7);
\draw[black, dotted] (1,1) -- (1,-2);
\draw[black, dotted] (3,3) -- (3,0);
\draw[black, dotted] (4,5) -- (4,-2);
\draw[black, dotted] (6,7) -- (6,0);
\draw[black] (3,0) -- (6,0);
\draw[black] (4,-2) -- (6,0);

\end{tikzpicture}
  \caption{A $(3,1)$-staircase of $(3,2)$-hypercubes in $X\in\catwg{3}$ with orientation $(1,1,0)$}\label{stair3D1}
\end{figure}

\bk

\begin{figure}
\vspace{15mm}
  \centering
  \begin{tikzpicture}[thick, mark options={color=black}]  fill opacity=0.6
\filldraw[fill=black!6!white, draw=white][thick] (3,3) -- (1,1) -- (1,-2) -- (3,0) -- (3,3);
\filldraw[fill=black!30!white, draw=black][thick] (0,3) -- (3,3) -- (1,1) -- (-2,1) -- (0,3);
\filldraw[fill=black!30!white, draw=black][thick] (3,0) -- (6,0) -- (4,-2) -- (1,-2) -- (3,0);
\filldraw[fill=black!6!white, draw=white][thick] (4,-2) -- (4,5) -- (6,7) -- (6,0) -- (4,-2);
\filldraw[fill=black!30!white, draw=black][thick] (6,7) -- (9,7) -- (7,5) -- (4,5) -- (6,7);

\filldraw[fill=black!30!white, draw=black][thick] (-2,9) -- (1,9) -- (3,11) -- (0,11) -- (-2,9);
\filldraw[fill=black!30!white, draw=black][thick] (1,9) -- (4,9) -- (6,11) -- (3,11) -- (1,9);
\filldraw[fill=black!30!white, draw=black][thick] (4,9) -- (7,9) -- (9,11) -- (6,11) -- (4,9);

\draw[black] (3,0) -- (6,0);
\draw[black] (4,-2) -- (6,0);

\draw[black, dotted] (1,1) -- (1,-2);
\draw[black, dotted] (3,3) -- (3,0);
\draw[black, dotted] (4,5) -- (4,-2);
\draw[black, dotted] (6,7) -- (6,0);
\draw[black, dotted] (-2,9) -- (-2,1);
\draw[black, dotted] (1,9) -- (1,1);
\draw[black, dotted] (4,9) -- (4,5);
\draw[black, dotted] (7,9) -- (7,5);
\draw[black, dotted] (0,11) -- (0,3);
\draw[black, dotted] (3,11) -- (3,3);
\draw[black, dotted] (6,11) -- (6,7);
\draw[black, dotted] (9,11) -- (9,7);

\end{tikzpicture}
  \caption{Lifting condition for staircase in Figure \ref{stair3D1}}\label{stair3D2}
\end{figure}

\label{3Dstair}

\clearpage
%%%%%%%%%%%%%%%%%%%%%%%%%%%%%%%%%%%%%%%%%%%%%%%%%%%%%%%%%%%%%%%%%%%%%%%%%%%%%%%

%%%%%%%%%%%%%%%%%%%%%%%%%%%%%%%%%%%%%%%%%%%%%%%%%%%%%%%%%%%%%%%%%%%%%%%%%%%%%%%%

\part{Weakly globular n-fold categories as a model of weak n-categories}\label{part-5}

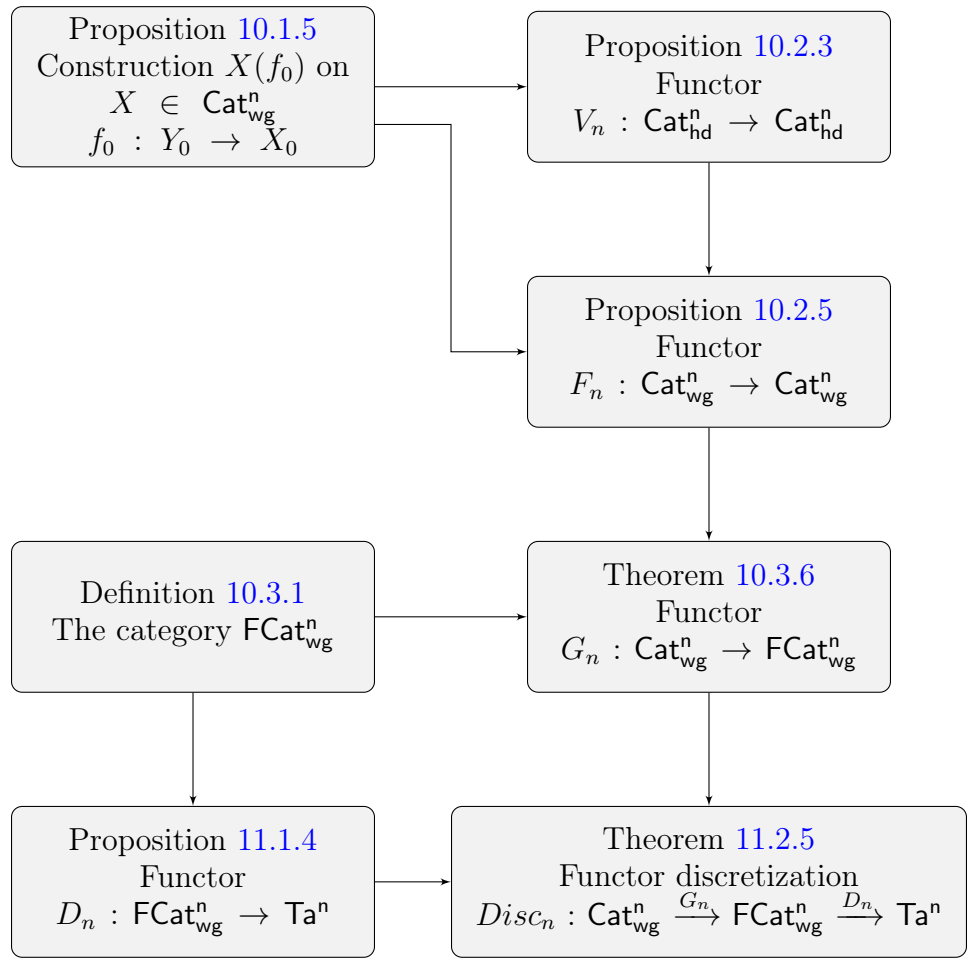
\begin{figure}[ht]
  \centering
\vspace{10mm}
\begin{tikzpicture}[node distance = 15mm and 20mm]% First vertical distance and second horizontal from border
\node [block3] (pro8_1_5) {Proposition \ref{pro-gen-const-1-new}\\Construction $X(f_0)$ on\\$X\in\catwg{n}$\\$f_0:Y_0\rw X_0$};
\node [block3, right= of pro8_1_5] (pro8_2_3) {Proposition \ref{pro-gen-const-2}\\Functor\\$V_n:\cathd{n}\rw\cathd{n}$}; %
\node [nullblock, below = of pro8_1_5] (null) {};
\node [block3, below= of pro8_2_3] (pro8_2_5) {Proposition \ref{pro-gen-const-1}\\Functor\\$F_n:\catwg{n}\rw \catwg{n}$};
\node [block3, below= of pro8_2_5] (the8_3_6) {Theorem \ref{pro-fta-1}\\Functor\\$G_n:\catwg{n}\rw\ftawg{n}$};
\node [block3, left= of the8_3_6] (def8_3_1) {Definition \ref{def-fta-1}\\The category $\ftawg{n}$};
\node [block3, below= of def8_3_1] (pro9_1_4) {Proposition \ref{pro-fta-tam-1}\\Functor\\$\Dn:\ftawg{n}\rw\ta{n}$};
\node [block2, below= of the8_3_6] (the9_1_9) {Theorem \ref{the-disc-func}\\Functor discretization\\$ \Discn:\catwg{n}\xrw{G_n}\ftawg{n}\xrw{\Dn}\ta{n}$};

% Draw edges
\draw [line] (pro8_1_5.east)++(0,-0.5) -- ++(1,0) --++(0,-3) |- (pro8_2_5.west);
\draw [line] (pro8_1_5.east) -- (pro8_2_3.west);
\draw [line] (pro8_2_3) -- (pro8_2_5);
\draw [line] (pro8_2_5) -- (the8_3_6);
\draw [line] (the8_3_6) -- (the9_1_9);
\draw [line] (def8_3_1) -- (the8_3_6);
\draw [line] (pro9_1_4) -- (the9_1_9);
\draw [line] (def8_3_1) -- (pro9_1_4);

\end{tikzpicture}
  \caption{Construction of the discretization functor.}
  \label{FigIntro-9}
\end{figure}

In Part \ref{part-5} we construct the discretization functor
 \begin{equation*}
 \Discn : \catwg{n}\rw\ta{n}
 \end{equation*}
 from weakly globular $n$-fold categories to Tamsamani $n$-categories, and we prove the final results: the comparison between $\catwg{n}$ and $\ta{n}$, exhibiting $\catwg{n}$ as a model of weak $n$-categories (Theorem \ref{cor-the-disc-func}), and the homotopy hypothesis for groupoidal weakly globular \nfol categories (Theorem \ref{cor-gta-2} and Corollary \ref{cor1-mod-fund-wg-group}). A schematic summary of the main results of this part is contained in Figures \ref{FigIntro-9} and \ref{FigIntro-5}.

 In Chapter \ref{chap8} we introduce the category $\ftawg{n}$, whose objects have homotopically discrete substructures with functorial sections to the discretization maps. The idea of the category $\ftawg{n}$ is introduced in Section \ref{sub-idea-ftawgn}, before the formal definitions. The main result of Chapter \ref{chap8} is Theorem \ref{pro-fta-1} on the existence of the functor
  \begin{equation*}
    G_n:\catwg{n}\rw \ftawg{n}.
  \end{equation*}
   This functor is constructed inductively using the functor
    \begin{equation*}
      F_n:\catwg{n}\rw \catwg{n}
    \end{equation*}
    of Proposition \ref{pro-gen-const-1}. The latter approximates up to $n$-equivalence a weakly globular \nfol category $X$ with a better behaved one $F_n X$ in which the homotopically discrete object at level 0 admits a functorial choice of section to the discretization map. This is based on a general construction on $X\in\catwg{n}$ and $f_0:Y_0\rw X_0$ given in Proposition \ref{pro-gen-const-1-new}, for an appropriate choice of the map $f_0:Y_0\rw X_0$, given in Proposition \ref{pro-gen-const-2}.

    The ideas of these constructions are explained in Section \ref{subs-idea-constrxf} (for the construction $X(f_0)$), in Section \ref{subs-idea-fn} (for the functors $V_n$ and $F_n$) and in Section \ref{subs-idea-gn} (for the functor $G_n$).

 In Chapter \ref{chap9} we construct the discretization functor and we obtain the main results of this work. In Proposition \ref{pro-fta-tam-1} we build a functor
  \begin{equation*}
    \Dn:\ftawg{n}\rw \ta{n}
  \end{equation*}
  which discretizes the homotopically discrete substructures of the objects of $\ftawg{n}$: because of the properties of the category $\ftawg{n}$, this can be done in a functorial way. The idea of the functor $D_n$ is explained in Section \ref{subs-fn-dn-idea}, before the formal definition.

   In Definition \ref{def-disc-func} we define the discretization functor as the composite
\begin{equation*}
     \Discn:\catwg{n}\xrw{G_n}\ftawg{n}\xrw{\Dn}\ta{n}\;
\end{equation*}
and we establish its properties in Theorem \ref{the-disc-func}. The idea of the functor $\Discn$ is explained in Section \ref{sec-disc-final}, before the formal definition.

 The properties of the discretization functor also depends on the properties of the functor rigidification $Q_n$ (see Proposition \ref{pro-fta-tam-2}). Finally, the functors $Q_n$ and $\Discn$ lead to the main comparison result (Theorem \ref{cor-the-disc-func}) on the equivalence of categories
\begin{equation*}
  \ta{n}\bsim^n\;\simeq \; \catwg{n}\bsim^n.
\end{equation*}

In the last part of Chapter \ref{chap9} we define the groupoidal version of the three Segal-type models, that is the categories $\gcatwg{n}$, $\gtawg{n}$ and $\gta{n}$. We show in Theorem \ref{cor-gta-2} that the category $\gcatwg{n}$ of groupoidal weakly globular \nfol categories gives a model of $n$-types.

 In Corollary \ref{cor1-mod-fund-wg-group} we exhibit an alternative and more convenient functor from spaces to $\gcatwg{n}$ using the results of Blanc and the author \cite{BP}.

%%%%%%%%%%%%%%%%%%%%%%%%%%%%%%%%%%%%%%%%%%%%%%%%%%%%%%%%%%%%%%%%%%%%%%%%%%%%%

% Figure

\clearpage

\begin{figure}[ht]
  \centering
\vspace{10mm}
\begin{tikzpicture}[node distance = 15mm and 20mm]% First vertical distance and second horizontal from border
\node [block3] (the7_4_1) {Theorem \ref{the-funct-Qn}\\$\Qn:\tawg{n}\rw\catwg{n}$};
\node [nullblock, below of= the7_4_1, node distance=15mm] (null) {};
\node [block3, below of = null, node distance=15mm] (the9_1_9) {Theorem \ref{the-disc-func}\\$\Discn:\catwg{n}\rw\ta{n}$};
\node [block3, right of = null, node distance=68mm] (the9_1_10) {Theorem \ref{cor-the-disc-func}\\$\ta{n}\bsim^{n}\;\simeq \catwg{n}\bsim^{n}$};
\node[block3, below = of the9_1_9] (def9_2_1) {Definitions \ref{def-gta-1}, \ref{def-gta-2}\\Categories\\$\gcatwg{n},\;\;\gta{n},\;\;\gtawg{n}$};
\node[block2, right of = def9_2_1, node distance=68mm] (cor9_2_8) {Theorem \ref{cor-gta-2}\\$\gcatwg{n}\bsim^n\;\simeq \mbox{Ho($n$-types)}$};
\node [rectangle, draw, fill=gray!10,
    text width=75mm, text centered, rounded corners, minimum height=36mm] [below = of cor9_2_8, yshift=0mm] (cor11_4_6){Corollary \ref{cor1-mod-fund-wg-group}\\ \mk $\mbox{$n$-types} \xrw{\clH_n}\gpdwg{n}\xrw{j}\gcatwg{n}$\\ \mk $\widetilde{B}:\gcatwg{n}\rw\mbox{$n$-types}$\\ \mk induce \hfill\rule{0ex}{0ex} \\ $\gcatwg{n}\bsim^n\; \simeq \; \mbox{Ho($n$-types)}$};

% Draw edges
%\draw [line] (pro8_1_5.east)++(0,-0.5) -- ++(1,0) --++(0,-3) |- (pro8_2_5.west);
%
\draw [line] (the9_1_10) -- (cor9_2_8);
\draw [line] (the7_4_1.east) -- (the9_1_10.west);
\draw [line] (the9_1_9.east) -- (the9_1_10.west);
\draw [line] (def9_2_1.east) -- (cor9_2_8.west);
\draw [line] (cor9_2_8) -- (cor11_4_6);

\end{tikzpicture}
  \caption{\catwg{n} as a model of weak $n$-categories}
  \label{FigIntro-5}
\end{figure}

\chapter{Functorial choices of homotopically discrete objects}\label{chap8}

In Chapter \ref{chap7} we built a rigidification functor from weakly globular Tamsamani $n$-categories to weakly globular $n$-fold categories, which in particular affords a functor
\begin{equation*}
  Q_n:\ta{n}\rw \catwg{n}
\end{equation*}
producing $n$-equivalent objects in $\tawg{n}$.

To reach the full comparison between $\ta{n}$ and $\catwg{n}$ we need a functor in the other direction, namely a \emph{discretization functor}

\begin{equation*}
 \Discn : \catwg{n}\rw\ta{n}
 \end{equation*}

The idea of the functor $\Discn$ is to replace the homotopically discrete sub-structures in $X\in\catwg{n}$ by their discretizations in order to recover the globularity condition. This affects the Segal maps, which from being isomorphisms in $X$ become $(n-1)$-equivalences in $\Discn X$.

We illustrate this idea in the case $n=2$. Given $X\in\catwg{2}$, by definition $X_0\in\cathd{}$, so there is a discretization map $\zg:X_0\rw X_0^d$ which is an equivalence of categories. Given a choice $\zg'$ of pseudo-inverse, we have $\zg\zg'=\Id$ since $X_0^d$ is discrete.

We can therefore construct $D_0 X\in\funcat{}{\Cat}$ as follows
\begin{equation*}
    (D_0 X)_k = \left\{
                  \begin{array}{ll}
                    X_0^d, & k=0 \\
                    X_k, & k>0\;.
                  \end{array}
                \right.
\end{equation*}
The face maps
 \begin{equation*}
 (D_0 X)_1\rightrightarrows (D_0 X)_0
 \end{equation*}
  are given by $\zg \pt_i$\; $i=0,1$ (where $\pt_i:X_1 \rightrightarrows X_0$ are face maps of $X$) while the degeneracy map
   \begin{equation*}
   (D_0 X)_0 \rw (D_0 X)_1
   \end{equation*}
    is $\zs_0\zg'$ (where $\zs_0:X_0 \rw X_1$ if the degeneracy map of $X$). All other face and degeneracy maps in $D_0 X$ are as in $X$. Since $\zg\zg'=\Id$, all simplicial identities are satisfied for $D_0X$. By construction, $(D_0 X)_0$ is discrete while the Segal maps are given, for each $k\geq 2$, by
\begin{equation*}
    \pro{X_1}{X_0}{k}\rw \pro{X_1}{X_0^d}{k}
\end{equation*}
and these are equivalences of categories since $X\in\catwg{2}$. Thus, by definition, $D_0 X\in\ta{2}$. This construction however does not afford a functor
 \begin{equation*}
 D_0:\catwg{2}\rw\ta{2}
 \end{equation*}
  but only a functor
\begin{equation*}
    D_0:\catwg{2}\rw(\ta{2})_{\ps}
\end{equation*}
where $(\ta{2})_{\ps}$ is the full subcategory of $\psc{}{\Cat}$ whose objects are in $\ta{2}$. This is because, for any morphism $F:X\rw Y$ in $\ta{2}$, the diagram in $\Cat$
\begin{equation*}
\xymatrix@C=35pt{
X_0^d \ar^{f^d}[r] \ar_{\zg'(X_0)}[d] & Y_0^d \ar^{\zg'(Y_0)}[d]\\
X_0 \ar_{f}[r] & Y_0
}
\end{equation*}
in general only pseudo-commutes. Hence $D_0$ cannot be used as a definition of the discretization functor $Disc_2$.

To overcome this problem we introduce the category $\ftawg{n}$ whose objects are weakly globular \nfol categories in which there are functorial choices of sections to the discretization maps of the homotopically discrete sub-structures. We then show that we can approximate any object of $\catwg{n}$ with an $n$-equivalent object of $\ftawg{n}$. Namely we prove in Theorem \ref{pro-fta-1} that there is a functor
\begin{equation*}
  G_n:\catwg{n}\rw \ftawg{n}
\end{equation*}
and an $n$-equivalence $G_n X\rw X$.
 In the next chapter we  build a functor
 \begin{equation*}
 D_n : \ftawg{n}\rw\ta{n}
 \end{equation*}
  and construct the discretization functor
 \begin{equation*}
 \Discn : \catwg{n}\rw\ta{n}
 \end{equation*}
  as the composite
\begin{equation*}
    \catwg{n}\xrw{G_n}\ftawg{n}\xrw{D_n}\ta{n}\;.
\end{equation*}
This chapter is organized as follows.

In Section \ref{sec-canonical} we develop a general construction on the category $\catwg{n}$ that allows to replace $X\in \catwg{n}$ with a $n$-equivalent one $X(f_0)$ by modifying $X_0\in \catwg{n-1}$ via a map $f_0:Y_0\rw X_0$ satisfying certain properties (see Proposition \ref{pro-gen-const-1-new}). In Section \ref{wg-canon-hom}  we make an appropriate choice of the map $f_0$ (see Proposition \ref{pro-gen-const-2}) to construct in Proposition \ref{pro-gen-const-1} a functor
 \begin{equation*}
   F_n:\catwg{n}\rw \catwg{n}
 \end{equation*}
  such that, for each $X\in\catwg{n}$, $(F_n X)_0\in\cathd{n-1}$ admits a functorial (that is, natural in $X$) choice of section to the discretization map $(F_n X)_0\rw{(F_n X)}^d_0$.

In Section \ref{sec-fta}, Theorem \ref{pro-fta-1}, we define the category $\ftawg{n}$ and we use the functor $F_n$ to build inductively the functor
 \begin{equation*}
   G_n:\catwg{n}\rw\ftawg{n}.
 \end{equation*}
Namely, we define $G_2=F_2$ and given $G_{n-1}$,
     \begin{equation*}
       G_n=\ovl{G}_{n-1}\circ F_n
     \end{equation*}
    see Definition \ref{def-Gn-new} and Theorem \ref{pro-fta-1}.
 In Section \ref{sec-fta-to-tam} the functor $G_n$ will be used to define the discretization functor from $\catwg{n}$ to $\ta{n}$.

 \vspace{5mm}

 \textbf{Important note:} Throughout this chapter, we will consider pullbacks in $\funcat{n-1}{\Cat}$ of maps in $\tawg{n}$ or in $\ftawg{n}$ after applying the functor $J_n$. For ease of notation, we omit writing explicitly $J_n$ in these diagrams. This is justified since $J_n$ is fully faithful. Similarly, when referring to maps in $\tawg{n}$ or in $\ftawg{n}$ to be 'levelwise equivalence of categories' or 'levelwise isofibrations' we always mean after applying the functor $J_n$.

%%
%%%%%%%%%%%%%%%%%%%%%%%%%%%%%%%%%%%%%%%%%%%%%%%%%%%%%%%%%%%%%
\section{A construction on $\pmb{\catwg{n}}$}\label{sec-canonical}
In this section we develop a general construction on the category $\catwg{n}$ that allows to replace $X\in \catwg{n}$ with a $n$-equivalent one $X(f_0)$ by modifying $X_0\in \cathd{n-1}$ in an appropriate way, via a map $f_0:Y_0\rw X_0$ in $\cathd{n-1}$ satisfying additional conditions (see Proposition \ref{pro-gen-const-1}).
 Proposition \ref{pro-gen-const-1} will be used in the next section to functorially approximate up to $n$-equivalence a weakly globular \nfol category $X$ with a better behaved one in which the homotopically discrete $(n-1)$-fold category $X_0$ admits a functorial choice of section to the discretization map.

 \subsection{The idea of the construction $\pmb{X(f_0)}$ }\label{subs-idea-constrxf}
 The construction $X(f_0)$ of Proposition \ref{pro-gen-const-1} is based on an application of a well known construction on internal categories (Lemma \ref{lem-gen-constr}), and on conditions on the map $f_0$ to ensure that $X(f_0)$ and $X$ are $n$-equivalent.

 For any internal category $X\in\Cat\clC$ and morphism $f_0:X'_0\rw X_0$, pulling back $f_0$ along the map $X_1\xrw{(\pt_0,\pt_1)} X_0\times X_0$ gives rise to an internal category $X(f_0)$ with $(X(f_0))_k$ given by the pullbacks \eqref{genconstr_1} and \eqref{genconstr_2}. This construction is also well behaved with respect to pullbacks, as spelled out in Lemma \ref{lem-gen-constr}.

 In Proposition \ref{pro-gen-const-1-new} we apply this general construction to $X\in\catwg{n}$, viewed as an internal category in $\catwg{n-1}$ in direction 1. Upon application of the functor $J_{n-1}:\catwg{n-1}\rw\funcat{n-2}{\Cat}$ (which, being fully faithful, preserves pullbacks) the pullbacks
\eqref{genconstr_1} and \eqref{genconstr_2} give rise to pullbacks in $\funcat{n-2}{\Cat}$ (we omit writing $J_{n-1}$ explicitly for ease of notation). These pullbacks are computed levelwise, that is for each $\uk\in\dop{n-2}$ they give rise to a pullback in $\Cat$.

The additional conditions imposed in the hypotheses of Proposition \ref{pro-gen-const-1-new} are such that the above levelwise pullbacks in $\Cat$ are pullbacks along isofibrations which are surjective on objects and the same is true after application of the functor $\p{r,n-1}$ for each $1<r<n$.

Two properties of pullbacks in $\Cat$ along isofibrations are particularly relevant for us:
\begin{itemize}
  \item [a)] They are preserved by $p$ (Lemma \ref{lem-gen-const-1})
  \item [b)] They preserve objects of $\cathd{}$ (Lemma \ref{lem-gen-const-3})
\end{itemize}

In the proof of Proposition \ref{pro-gen-const-1-new} we show that property a) implies that the \nfol category $X(f_0)$ satisfies the hypotheses of Proposition \ref{pro-crit-ncat-be-wg} b), and thus $X(f_0)\in \catwg{n}$, while property b) implies (via Lemma \ref{lem-gen-const-4}) that $X(f_0)\rw X$ is an $n$-equivalence.

From the fact that isofibrations are stable under pullbacks we also deduce in the proof of Proposition \ref{pro-gen-const-1-new} that the map $V(X): X(f_0)\rw X$ is levelwise an isofibration in $\Cat$ which is surjective on objects, and the same holds for $\p{r,n}V(X)$ and, under additional conditions on $f_0$, for $\q{r,n}V(X)$. This will be used in Section \ref{subs-vnfn} in the definition of the functors $V_n$ and $F_n$, where the construction $X(g_0)$ will be used for a map $g_0$ of the form $V(X)$.

In Corollary \ref{cor-gen-const-1-pb} we show that the construction $X(f_0)$ is well behaved with respect to pullbacks. The proof relies on the corresponding property of the construction of Lemma \ref{lem-gen-constr}. These properties will be used to study the behaviour with respect to pullbacks of $V_n$ and of $F_n$ in Section \ref{subs-vnfn}, in Corollary \ref{cor-gen-const-4} and
Corollary \ref{cor-fn-pull}. In turn, this will play an important role in the construction of the functor $G_n$ in Section \ref{gn-formal}.

 \begin{lemma}\label{lem-gen-constr}
 Let $\clC$ be a category with finite limits; let $X\in\Cat\clC$ and $f_0:X'_0\rw X_0$ be a morphism in $\clC$. There is $X(f_0)\in\Cat\clC$ with $X(f_0)_0=X'_0$ and $X(f_0)_1$ given by the pullback in $\clC$
\begin{equation}\label{genconstr_1}
    \xymatrix{
    X(f_0)_1 \ar^{}[rr] \ar^{}[d] && X'_0\times X_0 \ar^{f_0\times f_0}[d]\\
    X_1 \ar^{(\pt_0, \pt_1)}[rr] && X_0\times X_0
    }
\end{equation}
and a morphism in $\Cat\clC$
\begin{equation*}
  V(X):X(f_0)\rw X.
\end{equation*}

Further, given a diagram $X\rw Z \lw Y$ in $\Cat\clC$ and morphisms in $\clC$ $f_0:X'_0 \rw X_0$, $g_0:Y'_0 \rw Y_0$, $h_0:Z'_0\rw Z_0$, $X'_0\rw Z'_0$, $Y'_0\rw Z'_0$ marking the following diagram commute
\begin{equation*}
\xymatrix{
X'_0 \ar[r]\ar^{f_0}[d] & Z'_0 \ar^{h_0}[d] & Y'_0 \ar[l]\ar^{g_0}[d]\\
X_0 \ar[r] & Z_0 & Y_0 \ar[l]
}
\end{equation*}
Then
\begin{equation}\label{eq1-lem-gen-constr}
  (X\tiund{Z}Y)(f_0\tiund{h_0}g_0)=X(f_0)\tiund{Z(h_0)} Y(g_0)
\end{equation}

\begin{equation}\label{eq1-lem-gen-constr-map}
  V(X\tiund{Z}Y)=V(X)\tiund{V(Z)} V(Y)
\end{equation}
 \end{lemma}
\begin{proof}
Let $\pt'_i=\pr_i v:X(f_0)_1\rw X'_0$, $i=0,1$ where $\pr_0, \pr_1$ are the two projections, so that $v=(\pt'_0,\pt'_1)$ and $\pt_i f_1=f_0\pt'_i$, $i=0,1$. We have
\begin{equation*}
\begin{split}
    & (\pt_0,\pt_1)c (f_1,f_1)=(\pt_0\pr_0(f_1,f_1),\pt_1\pr_1(f_1,f_1))= \\
    & = (f_0\pt'_1\pr_0,f_0\pt'_2\pr_1)=(f_0\times f_0)(\pt'_0\pr_0,\pt'_1\pr_1)\;.
\end{split}
\end{equation*}
So there is a map
\begin{equation*}
  c':X(f_0)_1\tiund{Y_0} X(f_0)_1\rw X(f_0)_1
\end{equation*}
making the following diagram commute:
\begin{equation*}
\xymatrix@R=40pt@C=50pt{
X(f_0)_1\tiund{Y_0} X(f_0)_1 \ar_(0.6){c'}[r] \ar_{(f_1,f_1)}[d] \ar@/^4ex/^{(\pt'_0\pr_0,\pt'_1\pr_1)}[rr] & X(f_0)_1 \ar_{v=(\pt'_0,\pt'_1)}[r] \ar_{f_1}[d] & X'_0\times X'_0 \ar^{f_0\times f_0}[d]\\
\tens{X_1}{X_0}\ar_{c}[r] & X_1 \ar_{(\pt_0,\pt_1)}[r] & X_0\times X_0
}
\end{equation*}
In particular
\begin{equation*}
  \pt'_i c'=\pt'_i\pr_i\qquad i=0,1\;.
\end{equation*}
The other axioms of internal category for $X(f_0)$ follow immediately from the axioms for $X$ and the universal property of pullbacks. The morphism $V(X):X(f_0)\rw X$ is given by $(f_0, f_1)$.

Given $X\rw Z \lw Y$  as in the hypothesis, we have
\begin{equation*}
\{(X\tiund{Z}Y)(f_0\tiund{h_0}g_0)\}_0=X'_0\tiund{Z'_0}Y'_0=\{X(f_0)\}_0\tiund{\{Z(h_0)\}_0} \{Y(g_0)\}_0
\end{equation*}
and the pullback in $\clC$
\begin{equation*}
\xymatrix{
\{(X\tiund{Z}Y)(f_0\tiund{h_0}g_0\}_1 \ar[r] \ar[d] & \tens{(X'_0\tiund{Z'_0}Y'_0)}{}\ar[d]\\
X_1\tiund{Z_1}Y_1 \ar[r] & \tens{(X_0\tiund{Z_0}Y_0)}{}
}
\end{equation*}
Since
\begin{equation*}
 \tens{(X_0\tiund{Z_0}Y_0)}{}=(X_0\times X_0)\tiund{(Z_0\times Z_0)}(Y_0\times Y_0)
\end{equation*}
and similarly for $X'_0, Z'_0, Y'_0$ we conclude that
\begin{equation*}
  \{(X\tiund{Z}Y)(f_0\tiund{h_0}g_0)\}_1=\{X(f_0)\}_1\tiund{\{Z(h_0)\}_1} \{Y(g_0)\}_1\;.
\end{equation*}
Thus \eqref{eq1-lem-gen-constr} and \ref{eq1-lem-gen-constr-map} follow.

\end{proof}

\begin{remark}\label{rem-gen-contr}
Note that for each $k\geq 2$, there is a pullback in $\clC$
\begin{equation}\label{genconstr_2}
    \xymatrix@R=40pt{
    X(f_0)_k=\pro{X(f_0)_1}{Y_0}{k} \ar^{}[rr] \ar^{}[d] && \pro{Y_0}{}{k+1} \ar^{\pro{f_0}{}{k+1}}[d]\\
    X_k=\pro{X_1}{X_0}{k} \ar^{}[rr] && \pro{X_0}{}{k+1}
    }
\end{equation}
\end{remark}
\begin{lemma}\label{lem-gen-const-3}
    Let
\begin{equation*}
    \xymatrix{
    P \ar[r] \ar[d] & C \ar^{f}[d]\\
    A \ar_{s}[r] & B
    }
\end{equation*}
be a pullback in $\Cat$ with $f$ an isofibration and with $A,B,C \in \cathd{}$. Then $P\in\cathd{}$.
\end{lemma}
\begin{proof}
Since $f$ is an isofibration and $A\simeq A^d$, $B\simeq B^d$, $C\simeq C^d$, we have
\begin{equation*}
    P\simeq A \tms{ps}{B}C \simeq A^d\tiund{B^d} C^d
\end{equation*}
and therefore $P\in\cathd{}$.
\end{proof}
\begin{lemma}\label{lem-gen-const-4}
        Let
\begin{equation*}
    \xymatrix{
    P \ar[r] \ar_{h}[d] & C \ar^{f}[d]\\
    A \ar_{g}[r] & B
    }
\end{equation*}
be a pullback in $\funcat{n-1}{\Cat}$ with $A,B,C \in \catwg{n}$ and $f$ an $n$-equivalence which is a levelwise isofibration in $\Cat$ and the same holds for $p\up{r,n}f$ for all $1< r\leq n$. Then
\begin{itemize}
  \item [a)] For each $1< r\leq n$ there is a pullback in $\funcat{r-2}{\Cat}$

  \begin{equation*}
    \xymatrix{
    p\up{r,n}P \ar[rr] \ar[d] && p\up{r,n}C \ar[d]\\
    p\up{r,n}A \ar[rr] && p\up{r,n}B
    }
\end{equation*}

  \item [b)] $h$ is an $n$-equivalence.\mk
  \item [c)] Suppose, further, that $q\up{r,n}f$  is a levelwise isofibration in $\Cat$ for all $1< r\leq n$. Then for each $1< r\leq n$ there is a pullback in $\funcat{r-2}{\Cat}$

      \begin{equation*}
    \xymatrix{
    q\up{r,n}P \ar[rr] \ar[d] && q\up{r,n}C \ar[d]\\
    q\up{r,n}A \ar[rr] && q\up{r,n}B
    }
\end{equation*}

\end{itemize}

\end{lemma}
\begin{proof}
By induction on $n$. When $n=2$, for each $k\in\dop{}$ there is a pullback in $\Cat$
\begin{equation*}
\xymatrix@R=30pt@C=40pt{
P_{k} \ar[r] \ar_{h_{k}}[d] & C_{k} \ar^{f_{k}}[d]\\
A_{k} \ar_{g_{k}}[r] & B_{k}
}
\end{equation*}
where $f_{k}$ is an isofibration. Thus, by Lemma \ref{lem-gen-const-1}, there is a pullback in $\Set$
\begin{equation*}
\xymatrix@R=30pt@C=40pt{
(\p{2}P)_{k}= p P_{k} \ar[r] \ar[d] & p C_{k}=(\p{2}P)_{k} \ar[d]\\
(\p{2}A)_{k}= p A_{k} \ar[r] & p B_{k}=(\p{2}B)_{k}
}
\end{equation*}
Since this holds for each $k$, there is a pullback in $\Cat$
\begin{equation*}
\xymatrix@R=30pt@C=40pt{
\p{2}P \ar[r] \ar[d] & \p{2}C \ar[d] \\
\p{2}A \ar[r] & \p{2}B
}
\end{equation*}
which proves a). The proof of c) is similar. As for b), since $f$ is an isofibration, $P$ is equivalent to the pseudo-pullback $A\tms{ps}{B}C$; since $f$ is an equivalence of categories, the latter is equivalent to $A$.

\medskip
Suppose, inductively, that the lemma holds for $(n-1)$.

a) For each $\uk\in\dop{n-1}$ there is a pullback in $\Cat$
\begin{equation*}
\xymatrix@R=30pt@C=40pt{
P_{\uk} \ar[r] \ar_{h_{\uk}}[d] & C_{\uk} \ar^{f_{\uk}}[d]\\
A_{\uk} \ar_{g_{\uk}}[r] & B_{\uk}
}
\end{equation*}
where $f_{\uk}$ is an isofibration. Thus, by Lemma \ref{lem-gen-const-1}, there is a pullback in $\Set$
\begin{equation*}
\xymatrix@R=30pt@C=40pt{
(\p{n}P)_{\uk}= p P_{\uk} \ar[r] \ar[d] & p C_{\uk}=(\p{n}P)_{\uk} \ar[d]\\
(\p{n}A)_{\uk}= p A_{\uk} \ar[r] & p B_{\uk}=(\p{n}B)_{\uk}
}
\end{equation*}
Since this holds for each $\uk$, there is a pullback in $\funcat{n-2}{\Cat}$
\begin{equation}\label{eq1-lem-gen-const-4}
\xymatrix@R=30pt@C=40pt{
\p{n}P \ar[r] \ar[d] & \p{n}C \ar[d] \\
\p{n}A \ar[r] & \p{n}B
}
\end{equation}
which is a) for $r=n$. Since \eqref{eq1-lem-gen-const-4} satisfies the inductive hypothesis we deduce a pullback in $\funcat{r-2}{\Cat}$ for each $1< r \leq n-1$
\begin{equation}\label{eq1-lem-gen-const-4-new}
\xymatrix@R=30pt@C=40pt{
\p{r,n}P \ar[r] \ar[d] & \p{r,n}C \ar[d] \\
\p{r,n}A \ar[r] & \p{r,n}B
}
\end{equation}
\bk

b) By part a), there is a pullback in $\Cat$
\begin{equation*}
    \xymatrix{
    p\up{2,n}P \ar[rr] \ar[d] && p\up{2,n}C \ar[d]\\
    p\up{2,n}A \ar[rr] && p\up{2,n}B
    }
\end{equation*}
and therefore, at object level, a pullback in $\Set$
\begin{equation*}
    \xymatrix{
    P_0^d \ar[r] \ar_{}[d] & C_0^d \ar^{}[d]\\
    A_0^d \ar_{}[r] & B_0^d
    }
\end{equation*}
Let $(a,c),(a',c')\in P_0^d$. Then there is a pullback in $\catwg{n-1}$
\begin{equation*}
    \xymatrix{
    P((a,c),(a',c')) \ar[rr] \ar_{h((a,c),(a',c'))}[d] && C(c,c') \ar^{}[d]\\
    A(a,a') \ar_{}[rr] && B(ga,ga')
    }
\end{equation*}
By hypothesis, this satisfies the induction hypothesis so $h((a,c),(a',c'))$ is a $(n-1)$-equivalence.

We also have the pullback
\begin{equation*}
    \xymatrix{
    p\up{n}P \ar[rr] \ar_{p\up{n}h}[d] && p\up{n}C \ar^{p\up{n}f}[d]\\
    p\up{n}A \ar[rr] && p\up{n}B
    }
\end{equation*}
with $p\up{n}f$ a $(n-1)$-equivalence and thus by inductive hypothesis $p\up{n}h$ is a $(n-1)$-equivalence. We conclude that $h$ is an $n$-equivalence.

c) The proof is completely analogous to the one of part a).
\end{proof}

In the following Proposition we use the construction of Lemma \ref{lem-gen-constr} for $X\in \catwg{n}$ (viewed as internal category in $\catwg{n-1}$ in direction 1) for a particular choice of a map $f_0:Y_0 \rw X_0$ in $\cathd{n-1}$ such that the map $V(X):X(f_0)\rw X$ is an $n$-equivalence and has other desirable properties. This result will be used in the next section in the proofs of Propositions \ref{pro-gen-const-2} and \ref{pro-gen-const-1}.

\begin{proposition}\label{pro-gen-const-1-new}
    Let $X\in\catwg{n}$ and $f_0:Y_0\rw X_0$ be a morphism in $\cathd{n-1}$ such that, for each $1< r < n$, $f_0$ and $p\up{r,n-1}f_0$ are levelwise isofibrations in $\Cat$ which are surjective on objects. Then
    \begin{itemize}
      \item [a)] $X(f_0)\in\catwg{n}$  and
      \begin{equation*}
        p\up{r,n}(X (f_0))=(p\up{r,n}X)(p\up{r-1,n-1}f_0);
      \end{equation*}

      \item [b)] The map $V(X):X(f_0)\rw X$ is an $n$-equivalence;\mk

      \item [c)]$V(X)$ is levelwise an isofibration in $\Cat$  surjective on objects, and the same is true for $p\up{r,n}V(X)$ for all $1< r \leq n $.\mk

      \item [d)] If $X\in\cathd{n}$, $X(f_0)\in\cathd{n}$.\mk
      \item [e)] Suppose, further, that for each $1< r < n$, $q\up{r,n-1}f_0$ is a levelwise isofibration in $\Cat$ which is surjective on objects. Then
          \begin{equation*}
        q\up{r,n}(X(f_0))=(q\up{r,n}X)(q\up{r-1,n-1}f_0);
      \end{equation*}
      and $q\up{r,n}V(X)$ is a levelwise isofibration surjective on objects, for all $1< r\leq n$.

    \end{itemize}
\end{proposition}
\begin{proof}
By induction on $n$. Let $n=2$. Since $X(f_0)\in \cat{2}$ and  $Y_0\in\cathd{}$, to show that $X(f_0)\in\catwg{2}$ it is enough to show (by Lemma \ref{lem-crit-doucat-wg}) that $\ovl{p}X(f_0)\in N\Cat$. That is, for each $k\geq 2$
\begin{equation*}
    p(X(f_0))_k=\pro{p(X(f_0))_1}{p(X(f_0))_0}{k}\;.
\end{equation*}
We show this for $k=2$, the case $k>2$ being similar.
From Remark \ref{rem-gen-contr},
\begin{equation*}
    X(f_0)_2=(\tens{X_1}{X_0})\tiund{X_0\times X_0\times X_0} (Y_0\times Y_0\times Y_0)\;.
\end{equation*}
Since $f_0$ is an isofibration, using Lemma \ref{lem-gen-const-1}, the fact that $p\up{2}X\in\Cat$ and the fact that $p$ preserves products, we obtain
\begin{align}\label{part1-n2}
    &pX(f_0)_2=p(\tens{X_1}{X_0})\tiund{p(X_0\times X_0\times X_0)} p(Y_0\times Y_0\times Y_0)=\\
    =&(\tens{pX_1}{pX_0})\tiund{pX_0\times pX_0\times pX_0} pY_0\times pY_0\times pY_0\;.
\end{align}
On the other hand,
\begin{align*}
    & \ \tens{pX(f_0)_1}{pX(f_0)_0}=\\
    = & \ \tens{(pX_1\tiund{\tens{pX_0}{}}(\tens{pY_0}{}))}{pX_0\tiund{pX_0}pY_0}=\\
    = & \ (\tens{pX_1}{pX_0})\tiund{\tens{(\tens{pX_0}{}){}} {pX_0}}\tens{(\tens{pY_0}{}){}}{pY_0}=\\
    = & \ p(\tens{X_1}{X_0})\tiund{pX_0\times pX_0\times pX_0}(pY_0\times pY_0\times pY_0)\;.
\end{align*}
Therefore
\begin{equation*}
    pX(f_0)_2=\tens{pX(f_0)_1}{pX(f_0)_0}\;.
\end{equation*}
The case $k>2$ is similar. This shows $X(f_0)\in\catwg{2}$. From \eqref{part1-n2} we see that
 \begin{equation*}
   p\up{2}(X (f_0))=(p\up{2}X)(pf_0).
 \end{equation*}
This concludes the proof of a) when $n=2$.

We now show that $X(f_0)\rw X$ is a 2-equivalence. Let $a,b \in Y_0^d$. We have a pullback in $\Cat$
\begin{equation}\label{eq1-gen-const}
    \xymatrix{
    X(f_0)(a,b) \ar[rr] \ar_{V(X)(a,b)}[d] && Y_0(a)\times Y_0(b) \ar[d]\\
    X_1(f_0 a,f_0 b) \ar[rr] && X_0(f_0 a)\times X_0(f_0 b)
    }
\end{equation}
Since $X_0,Y_0\in\cathd{}$, $Y_0(a)\rw X_0(f_0(a))$ is an equivalence of categories hence it is in particular fully faithful. Thus by Lemma \ref{lem-gen-const-2}, $V(X)(a,b)$ is also fully faithful. Applying Lemma \ref{lem-gen-const-1} to \eqref{eq1-gen-const} we obtain a pullback in $\Set$
\begin{equation}\label{eq2-gen-const}
    \xymatrix{
   pX(f_0)(a,b) \ar[rr] \ar_{}[d] && pY_0(a)\times pY_0(b) \ar[d]\\
    pX_1(f_0 a,f_0 b) \ar[rr] && pX_0(f_0 a)\times pX_0(f_0 b)
    }
\end{equation}
Since, by hypothesis, $Y_0\rw X_0$ is surjective on objects, the right vertical map in \eqref{eq2-gen-const} is surjective, therefore such is the left vertical map in \eqref{eq2-gen-const}.
Thus
 \begin{equation*}
   X(f_0)(a,b) \rw X_1(f_0 a,f_0 b)
 \end{equation*}
 is essentially surjective on objects and in conclusion it is an equivalence of categories.

To show that $V(X):X(f_0)\rw X$ is a 2-equivalence, it is enough to show (by Proposition \ref{pro-n-equiv}) that $pp\up{2}V(X)$ is surjective. This follows from the fact that $p\up{2}V(X)$ is surjective on objects as $pf_0$ is surjective (since by hypothesis $f_0$ is surjective on objects). This concludes the proof of b) in the case $n=2$. Using Remark \ref{rem-gen-contr} and the fact that isofibrations are stable under pullbacks, c) follows.

Finally, if $X\in\cathd{2}$, since by a)  $X(f_0)\in\catwg{2}$ and $V(X)$ is a 2-equivalence, it follows from Proposition \ref{pro-nequiv-to-obj} that $X(f_0)\in\cathd{2}$, proving d) in the case $n=2$. The proof of part e) for the case $n=2$ is completely similar to the one of part a).

Suppose, inductively, that the proposition holds for $(n-1)$, let $X\in\catwg{n}$ and $f_0$ be as in the hypothesis.\bk

a) We show that $X(f_0)\in\catwg{n}$ by proving that it satisfies the hypothesis of Proposition \ref{pro-crit-ncat-be-wg} b). By the general construction of Lemma \ref{lem-gen-constr}, $X(f_0)\in\cat{n}$ and $(X(f_0))_0=Y_0\in\cathd{n-1}$. Since $X\in\catwg{n}$, $X_{\bl 0}\in\cathd{n-1}$ and we have a pullback in $\funcat{n-2}{\Cat}$
\begin{equation}\label{eq3-gen-const}
    \xymatrix{
    X(f_0)_{10} \ar[rr] \ar[d] && Y_{00}\times Y_{00}\ar^{f_{00}\times f_{00}}[d]\\
    X_{10} \ar[rr] && X_{00}\times X_{00}
    }
\end{equation}
where, by hypothesis, $f_{00}$ and $p\up{r,n-2} f_{00}$ is levelwise an isofibration in $\Cat$, which is surjective on objects. Thus \eqref{eq3-gen-const} satisfies the induction hypothesis c) and we conclude that $(X(f_0))_{\bl 0}\in\cathd{n-1}$. In particular, $(X(f_0))_{s0}\in\cathd{n-2}$. Thus hypothesis i) in Proposition \ref{pro-crit-ncat-be-wg} holds.

 To show that hypothesis ii) holds, we need to show that $\ovl{p}J_n X(f_0)\in N\lo{n-1}\catwg{n-1}$. Let $\ur\in \dop{n-2}$; by construction,
\begin{equation*}
    (X(f_0))_{2\ur}=  \{\tens{X_{1\ur}}{X_{0\ur}}\}\tiund{X_{0\ur}\times X_{0\ur}\times X_{0\ur}} \{Y_{0\ur}\times Y_{0\ur}\times Y_{0\ur}\}\;.
\end{equation*}
Since, by hypothesis, $(f_0)_{\ur}$ is an isofibration, by Lemma \ref{lem-gen-const-1} we obtain
\begin{equation*}
    p(X(f_0))_{2\ur}= p \{\tens{X_{1\ur}}{X_{0\ur}}\}\tiund{p \{X_{0\ur}\times X_{0\ur}\times X_{0\ur}\}} p\{Y_{0\ur}\times Y_{0\ur}\times Y_{0\ur}\}\;.
\end{equation*}
Since $X\in \catwg{n}$, $\p{n}X \in \catwg{n-1}$ and since $(\p{n-1}X_j)_{\ur}=pX_{j\ur}$ we have
\begin{equation*}
    p \{\tens{X_{1\ur}}{X_{0\ur}}\}=\tens{p(X_{1\ur})}{p(X_{0\ur})}
\end{equation*}
As $p$ commutes with products, we obtain
\begin{equation*}
    p(X(f_0))_{2\ur}=\{\tens{pX_{1\ur}}{pX_{0\ur}}\}\tiund{pX_{0\ur}\times pX_{0\ur}\times pX_{0\ur}}\{pY_{0\ur}\times pY_{0\ur}\times pY_{0\ur}\}\;.
\end{equation*}
Since this holds for all $\ur$, it follows that
\begin{equation*}
    \ovl{p}X(f_0)_2=\tens{\ovl{p}X(f_0)_1}{\ovl{p}X(f_0)_0}\;.
\end{equation*}
Similarly for each $k>2$,
\begin{equation*}
    \ovl{p}X(f_0)_k=\pro{\ovl{p}X(f_0)_1}{\ovl{p}X(f_0)_0}{k}\;.
\end{equation*}
We conclude that
\begin{equation}\label{pin}
    \ovl{p}J_n X(f_0)=(\p{n}X)(\p{n-1}f_0)\;.
\end{equation}

The map $\p{n}f_0$ satisfies the inductive hypothesis a) and we therefore conclude that $\ovl{p}J_n X(f_0)\in\catwg{n-1}$. Thus $X(f_0)$ satisfies the hypotheses of Proposition \ref{pro-crit-ncat-be-wg}b) and we conclude that $X(f_0)\in\catwg{n}$.

It also follows from \eqref{pin} that

\begin{equation*}
        \p{n}(X (f_0))=(\p{n}X)(\p{n-1}f_0).
      \end{equation*}
Using the induction hypothesis we further obtain

\begin{align*}
p\up{r,n}(X(f_0))& =p\up{r,n-1}\p{n}(X (f_0))=\\
&=p\up{r,n-1}(\p{n}X)(\p{n-1}f_0)=\\
&=(p\up{r,n-1}\p{n}X)(p\up{r-1,n-2}p\up{n-1}f_0)=\\
&=(p\up{r,n}X)(p\up{r-1,n-1}f_0).
\end{align*}
This completes the proof of a).
\bk

b) Let $a,b\in Y_0^d$. There is a pullback in $\catwg{n-1}$
\begin{equation*}
    \xymatrix{
    X(f_0)(a,b) \ar[rr]\ar[d] && Y_0(a)\times Y_0(b) \ar[d]\\
    X(f_0 a,f_0 b)\ar[rr] && X_0(f_0 a)\times X_0(f_0 b)
    }
\end{equation*}
Since $X_0(f_0 a),\, Y_0(a)\in\cathd{n-1}$ and $Y_0(a)^d=\{a\}\cong \{f_0 a\}=X_0(f_0 a)^d$, the map $Y_0(a)\rw X_0(f_0 a)$ is an $n$-equivalence (by Lemma \ref{lem-neq-hom-disc}). Also, this map is a levelwise isofibration in $\Cat$ (as such is $Y_0\rw X_0$). It follows from Lemma \ref{lem-gen-const-4} that
\begin{equation*}
  X(f_0)(a,b)\rw X(f_0 a,f_0 b)
\end{equation*}
 is a $(n-1)$-equivalence.

Finally, since by hypothesis $\p{2,n-1}X(f_0)$ is surjective on objects, $\p{1,n-1}X(f_0)$  is surjective. By a), this is the object part of the map  $\p{2,n}X(f_0)\rw \p{2,n}X$ which implies that $\p{1,n}X(f_0)\rw \p{1,n}X$ is surjective. By Proposition \ref{pro-n-equiv}, we conclude that $X(f_0)\rw X$ is an $n$-equivalence.\bk

c) By construction there are pullbacks in $\funcat{n-2}{\Cat}$, for each $k\geq 2$
\begin{equation*}
\xymatrix@R=33pt@C=40pt{
X(f_0)_1 \ar[r] \ar_{V(X)_1}[d] & Y_0\times Y_0 \ar[d]\\
X_1 \ar[r]  & X_0\times X_0
}
\qquad
\xymatrix@R=30pt@C=40pt{
X(f_0)_k \ar[r] \ar_{V(X)_k}[d] & \pro{Y_0}{}{k+1} \ar[d]\\
X_k \ar[r]  &  \pro{X_0}{}{k+1}
}
\end{equation*}
Since pullbacks in $\funcat{n-2}{\Cat}$ are computed pointwise, for each $\ur\in\dop{n-2}$ there are pullbacks in $\Cat$
\begin{equation*}
\xymatrix@R=30pt@C=30pt{
X(f_0)_{1\ur} \ar[r] \ar_{V(X)_{1\ur}}[d] & Y_{0\ur}\times Y_{0\ur} \ar^{f_{0\ur}\times f_{0\ur}}[d]\\
X_{1\ur}\ar[r] & Y_{0\ur}\times Y_{0\ur}
}
\qquad
\xymatrix@R=30pt@C=30pt{
X(f_0)_{k\ur} \ar[r] \ar_{V(X)_{k\ur}}[d] & \pro{Y_{0\ur}}{}{k+1} \ar^{\pro{f_{0\ur}}{}{k+1}}[d]\\
X_{k\ur}\ar[r] & \pro{X_{0\ur}}{}{k+1}
}
\end{equation*}
Since, by hypothesis, $f_{0\ur}$ is an isofibration and isofibrations are stable under pullbacks, $V(X)_{1\ur}$ and $V(X)_{k\ur}$ are isofibrations. Since the object functor $ob:\Cat\rw\Set$ preserves pullbacks, we also have pullbacks in $\Set$
\begin{equation*}
\begin{split}
   & \xymatrix@R=30pt@C=40pt{
\ob X(f_0)_{1\ur} \ar[r] \ar_{\ob V(X)_{1\ur}}[d] & \ob Y_{0\ur}\times \ob Y_{0\ur} \ar^{\ob f_{0\ur}\times \ob f_{0\ur}}[d]\\
\ob X_{1\ur}\ar[r] & \ob Y_{0\ur}\times \ob Y_{0\ur}
} \\
    & \xymatrix@R=30pt@C=40pt{
\ob X(f_0)_{k\ur} \ar[r] \ar_{\ob V(X)_{k\ur}}[d] & \pro{\ob Y_{0\ur}}{}{k+1} \ar^{\pro{\ob f_{0\ur}}{}{k+1}}[d]\\
\ob X_{k\ur}\ar[r] & \pro{\ob X_{0\ur}}{}{k+1}
}
\end{split}
\end{equation*}
Since, by hypothesis, $\ob f_{0\ur}$ is surjective, such are $\ob(V(X)_{1\ur})$ and $\ob(V(X)_{k\ur})$. In conclusion, $V(X)$ is levelwise an isofibration surjective on objects.

By part a)
\begin{equation*}
  \p{r,n}(X(f_0))=(\p{r,n}X)(\p{r-1,n-1}f_0)
\end{equation*}
and by hypothesis $\p{r-1,n-1}f_0$ is levelwise an isofibration surjective on objects. By induction hypothesis applied to $\p{r,n}X$ we deduce that $\p{r,n}V(X)$ is also a levelwise isofibration surjective on objects.

\bk

d) This follows from a) and b) using Proposition \ref{pro-nequiv-to-obj}.\bk

e) Reasoning as in part a) it is easy to see that
  $$\ovl{q}J_n X(f_0)=(\q{n}X)(\q{n-1}f_0)$$
   from which we deduce that
\begin{equation*}
        \q{n}(X (f_0))=(\q{n}X)(\q{n-1}f_0).
      \end{equation*}
The rest follows by a similar argument as in parts a) and c).

\end{proof}
We next study the behaviour of the construction of Proposition \ref{pro-gen-const-1-new} with respect to pullbacks.

\begin{corollary}\label{cor-gen-const-1-pb}
\

\begin{itemize}
  \item [a)] Suppose that $X\rw Z\lw Y$ is a diagram in $\catwg{n}$ such that $X\tiund{Z}Y\in\catwg{n}$ and suppose we have a commuting diagram in $\catwg{n-1}$
\begin{equation*}
\xymatrix{
X'_0 \ar[r]\ar^{f_0}[d] & Z'_0 \ar^{h_{0}}[d] & Y'_0 \ar[l]\ar^{g_0}[d]\\
X_0 \ar[r] & Z_0 & Y_0 \ar[l]
}
\end{equation*}
such that $X'_0\tiund{Z'_0}Y'_0\in\catwg{n-1}$, then
      \begin{equation*}
        (X\tiund{Z} Y)(f_0\tiund{h_0}g_0)\cong X(f_0)\tiund{Z(h_0)}Y(g_0)\;,
      \end{equation*}
      \begin{equation*}
  V(X\tiund{Z}Y)=V(X)\tiund{V(Z)} V(Y)
\end{equation*}

  \item [b)] Suppose, further, that $f_0$, $g_0$, $h_0$, $\p{r,n-1}f_0$, $\p{r,n-1}g_0$, $\p{r,n-1}h_0$ are levelwise isofibrations in $\Cat$ surjective on objects for all $1< r < n$ and
      \begin{equation*}
        \p{r,n}(X\tiund{Z} Y)=\p{r,n}X\tiund{\p{r,n}Z}\p{r,n}Y
      \end{equation*}
      for all $1\leq r \leq n$. Then
      \begin{equation}\label{eq1-cor-gen-const-1-pb}
        \begin{split}
            & \p{r,n}\{(X\tiund{Z} Y)(f_0\tiund{h_0}g_0)\}\cong \\
            & \cong \p{r,n}\{X(f_0)\}\tiund{\p{r,n}\{Z(h_0)\}} \p{r,n}\{Y(g_0)\}
        \end{split}
      \end{equation}
  \item [c)] If, in addition, $\q{r,n-1}f_0$, $\q{r,n-1}h_0$, $\q{r,n-1}g_0$, are levelwise isofibrations in $\Cat$ surjective on objects for all $1< r < n$ and
      \begin{equation*}
        \q{r,n}(X\tiund{Z}Y)=\q{r,n}X\tiund{\q{r,n}Z}\q{r,n}Y
      \end{equation*}
      for all $1\leq r \leq n$, then
      \begin{equation*}
      \begin{split}
          & \q{r,n}\{(X\tiund{Z}Y)(f_0\tiund{h_0}g_0)\}\cong \\
          & \cong\q{r,n}\{X(f_0)\}\tiund{\q{r,n}\{Z(h_0)\}}\q{r,n}\{Y(g_0)\}
      \end{split}
      \end{equation*}
\end{itemize}
\end{corollary}
\begin{proof}
\

 a) This follows from Lemma \ref{lem-gen-constr}, taking $\clC=\catwg{n-1}$.

\bk
 b) By hypothesis and by part a)
\begin{equation*}
\begin{split}
    & (\p{r,n}(X\tiund{Z}Y))(\p{r-1,n-1}f_0\tiund{\p{r-1,n-1}h_0}\p{r-1,n-1}g_0)\cong \\
    & \cong (\p{r,n}X\tiund{\p{r,n}Z}\p{r,n}Y) (\p{r-1,n-1}f_0\tiund{\p{r-1,n-1}h_0}\p{r-1,n-1}g_0)\cong \\
    & \cong (\p{r,n}X)(\p{r-1,n-1}f_0)\tiund{(\p{r,n}Z)(\p{r-1,n-1}h_0)}(\p{r,n}Y)(\p{r-1,n-1}g_0)\;.
\end{split}
\end{equation*}
By Proposition \ref{pro-gen-const-1-new}
\begin{equation*}
  (\p{r,n}X)(\p{r-1,n-1}f_0)=\p{r,n}(X(f_0))
\end{equation*}
and similarly for $Y,Z$ and $X\tiund{Z}Y$. Therefore from above we obtain \ref{eq1-cor-gen-const-1-pb}.

\bk
 c) The proof is completely analogous to the one of b), with $\q{r,n}$ in place of $\p{r,n}$, and using Proposition \ref{pro-gen-const-1-new} e).

\end{proof}
%%%
%%%%

\section{Weakly globular $\pmb{n}$-fold categories and functorial choices of homotopically discrete objects}\label{wg-canon-hom}

\medskip

In this section we show (Proposition \ref{pro-gen-const-1}) that we can functorially approximate up to $n$-equivalence a weakly globular \nfol category $X$ with a better behaved one $F_n X$ in which the homotopically discrete $(n-1)$-fold category at level $0$ admits a functorial choice of section to the discretization map. The functor $F_n$ will be used in Section \ref{sec-fta} to construct the functor $G_n:\catwg{n}\rw\ftawg{n}$ which will then lead in Section \ref{sec-fta-to-tam} to the discretization functor from $\catwg{n}$ to $\ta{n}$.

\subsection{The idea of the functors $\pmb{V_n}$ and $\pmb{F_n}$}\label{subs-idea-fn}
 The functor $F_n$ is based on the construction $X(f_0)$ of the previous section given in Proposition \ref{pro-gen-const-1-new} for an appropriate choice of the map $f_0:Y_0\rw X_0$.

 When $n=2$ we define
  \begin{equation*}
    F_2 X=X(u_{X_0})
  \end{equation*}
   where $u_{X_0}:\Dec X_0 \rw X_0$ is as in Section \ref{decalage}. Since $u_{X_0}$ is an isofibration surjective on objects (see Lemma \ref{lem-gen-const-5}), by Proposition \ref{pro-gen-const-1-new} this replaces $X$ with a $2$-equivalent $F_2 X$ in which $(F_2 X)_0=\Dec X_0$. As observed in Section \ref{decalage},  $\Dec X_0$ is homotopically discrete and has a functorial choice of section to the discretization map.

When $n>1$ the construction of $F_n X$ is again based on Proposition \ref{pro-gen-const-1-new} but for a more complex choice of the map $f_0:Y_0\rw X_0$. We build the appropriate map in Proposition \ref{pro-gen-const-2} where we construct a functor
 \begin{equation*}
 V_n:\cathd{n}\rw \cathd{n}
 \end{equation*}
 together with a map
  \begin{equation*}
    v_{n}(X):V_n X \rw X
  \end{equation*}
   for all $X\in \cathd{n}$ satisfying the hypotheses needed to apply Proposition \ref{pro-gen-const-1-new} and such that if $h:X\rw Y$ is a morphism in $\cathd{n}$, the following diagram  commutes for appropriate choices of sections to the discretization maps $V_n X \rw (V_n X)^d$ and $V_n Y \rw (V_n Y)^d$.
          \begin{equation}\label{eq-v-functo}
            \xymatrix{
            V_n X \ar[rr] && V_n Y\\
            (V_n X)^d \ar[rr]\ar[u] && (V_n Y)^d\ar[u]
            }
          \end{equation}
 The construction of $V_n$ is inductive on dimension, starting with  $V_1 X=\Dec X$ and $v_{1}(X)=u_X:\Dec X\rw X$ as in Section \ref{decalage}; By Lemma \ref{lem-gen-const-5}, $v_{1}(X)$ is an isofibration and is surjective on objects.

 Having defined inductively $V_{n-1}$ and $v_{n-1}{(X_0)}:V_{n-1}(X_0) \rw X_0$, we define
  \begin{equation*}
  F_n X=X(v_{n-1}(X_0))
  \end{equation*}
   and we define $V_{n}(X)$ via the pullback in $\funcat{n-1}{\Cat}$
\begin{equation}\label{eq-Vn_bis}
    \xymatrix{
    V_{n}(X) \ar^{h}[rr]\ar^{l}[d] && F_n X\ar^{r}[d]\\
    \di{n,2}\Dec\q{2,n}F_n X\ar^{u'}[rr] && \di{n,2}\q{2,n}F_n X
    }
\end{equation}
where
 \begin{equation*}
   u'=u_{(\q{2,n}F_n X)}.
 \end{equation*}
 The map $v_{n}(X):V_n X \rw X$ is defined by the composite
 \begin{equation*}
    v_{n}(X):V_n X \xrw{\;\;h\;\;}F_n X\xrw{f_{n}(X)} X\;.
\end{equation*}
where the map $f_{n}(X)$ is as in Proposition \ref{pro-gen-const-1-new} (we adopted the notation $f_{n}(X)$ instead of $V(X))$.
By the definition of $V_n(X)$ in \eqref{eq-Vn_bis} in order to construct the map
\begin{equation*}
  \zg'_{V_n X}: (V_n X)^d\rw V_n X
\end{equation*}
we need to construct maps
\begin{equation*}
  t: (V_n X)^d \rw \di{n-2}\Dec \q{2,n}F_n X,\qquad s: (V_n X)^d \rw F_n X
\end{equation*}
such that $u't=rs$, where $u', r$ are as in \eqref{eq-Vn_bis}. It is not difficult to show that
\begin{equation*}
  (V_n X)^d=\di{n,2}(\Dec\q{2,n}F_n X)^d=\di{1}(V_{n-1} X_0)^d.
\end{equation*}
This gives rise to natural maps
\begin{equation*}
  t: (V_n X)^d=\di{n,2}(\Dec\q{2,n}F_n X)^d\rw\di{n,2}\Dec\q{2,n}F_n X\;.
\end{equation*}
and
\begin{equation*}
  s: (V_{n}X)^d = \di{1}(V_{n-1}X_0)^d \xrw{\di{1} v} \di{1}(F_n X)_0 \xrw{\za_{F_n X}}F_n X\;.
\end{equation*}
where
\begin{equation*}
  (F_n X)_0^d = (V_{n-1}X_0)^d \xrw{\;v\;} V_{n-1}X_0 =(F_n X)_0\;.
\end{equation*}
is given by the inductive hypothesis and $\za$ is the counit of the adjunction $d \dashv \ob$. In the proof of Proposition \ref{pro-gen-const-2} we show that $t$ and $s$ define $\zg'_{V_n X}: (V_n X)^d\rw V_n X$ with the required properties.

In the course of this proof we thus define a functor
\begin{equation*}
        F_n:\cathd{n}\rw \cathd{n}
    \end{equation*}
    by
\begin{equation*}
  F_n X=X(v_{n-1}(X_0))
  \end{equation*}
for  $v_{n-1}(X_0):V_{n-1}X_0 \rw X_0$ in $\cathd{n-1}$. In Proposition \ref{pro-gen-const-1} we extend this to a functor
\begin{equation*}
        F_n:\catwg{n}\rw \catwg{n}
    \end{equation*}
 again defined by
\begin{equation*}
  F_n X=X(v_{n-1}(X_0))
  \end{equation*}
Since $(F_n X)_0=V_{n-1}X_0$ the functoriality of the sections to the discretization map of $(F_n X)_0$ follows immediately from \eqref{eq-v-functo}.

The functors $V_n$ and $F_n$ are well behaved with respect to pullbacks, and this plays an important role in the construction of the functor $G_n$ in Section \ref{gn-formal} (see also the informal discussion in Section \ref{subs-idea-gn} about this point). The proofs of these facts are based on the properties with respect to pullbacks of the general constructions of Lemma \ref{lem-gen-constr} and Corollary \ref{cor-gen-const-1-pb}. In Corollary \ref{cor-gen-const-4}, we show that given a diagram in $\cathd{n}$
\begin{equation*}
X\rw Z \lw Y
\end{equation*}
 such that $X\tiund{Z}Y \in \cathd{n}$ and such that this pullback is preserved by $\p{r,n}$ and $\q{r,n}$ (for each $1\leq r\leq n$), then this pullback is also preserved by $V_n$, $\p{r,n}V_n$ and $\q{r,n}V_n$.

We then deduce in Corollary \ref{cor-fn-pull} a similar property for the functor $F_n$. Namely, given a diagram in $\catwg{n}$
\begin{equation*}
X\rw Z \lw Y
\end{equation*}
 such that $X\tiund{Z}Y \in \catwg{n}$ and such that this pullback is preserved by $\p{r,n}$ and $\q{r,n}$ (for each $1\leq r\leq n$), then this pullback is also preserved by $F_n$, $\p{r,n}F_n$ and $\q{r,n}F_n$.

\subsection{The functors $\pmb{V_n}$ and $\pmb{F_n}$}\label{subs-vnfn} In this section we construct the functors $V_n$ and $F_n$ and we study their properties.

\begin{lemma}\label{lem-gen-const-5}
    If $X\in\cathd{}$, $\Dec X\in\cathd{}$ and the map $u_X:\Dec X\rw X$ as in Section \ref{decalage} is an isofibration surjective on objects.
\end{lemma}
\begin{proof}
Since $X\in\cathd{}$, $X=A[f]$ for a surjective map of sets $f:A\rw B$, where $A[f]$ is as in Definition \ref{def-int-eq-rel}. Thus $\Dec X=(\tens{A}{B})[d_0]$ where $d_0:\tens{A}{B}\rw A$, $d_0(x,y)=x$. The source and target maps
\begin{equation*}
    {\tld}_0, {\tld}_1:(\Dec X)_1=A\tiund{B}A\tiund{B}A\rw(Dec X)_0=\tens{A}{B}
\end{equation*}
are ${\tld}_0(x,y,z)=(x,y)$, ${\tld}_1(x,y,z)=(x,z)$.

Given $(x,y)\in(\Dec X)_0$ and an isomorphism $(u_X(x,y)=y,z)\in X_1$, we have $(x,y,z)\in(\Dec X)_1$ with $\tld_1(x,y,z)=(x,z)$. In picture:
\begin{equation*}
    \xymatrix{
    (x,y)=\tld_0(x,y,z) \ar^{(x,y,z)}[rr] \ar[d] && (x,z)=\tld_1(x,y,z)\ar[d]\\
    u_X(x,y)=y \ar_{(y,z)=u_X(x,y,z)}[rr] && u_X(x,z)=z
    }
\end{equation*}
By definition, this shows that $u_X:\Dec X\rw X$ is an isofibration. It is also surjective on objects since $(u_X)_1=pr_2:\tens{A}{B}\rw A$ is surjective.
\end{proof}
In the following the notation $\di{n,j}$ is as in \ref{dnot-ex-tam}.

\begin{lemma}\label{lem-gen-const-6}\
    \begin{itemize}
      \item [a)] Let $A\in\cathd{n}$, $B,C\in\Set$ and consider the pullback in $\funcat{n-1}{\Cat}$
          \begin{equation*}
            \xymatrix{
            Q \ar[r] \ar[d] & A\ar[d]\\
            \di{n,1}C \ar[r] & \di{n,1}B
            }
          \end{equation*}
          then $Q\in\cathd{n}$.
          \bk

      \item [b)] Let $X\in\cathd{n}$, $Z\in\cathd{}$ and consider the pullback in $\funcat{n-1}{\Cat}$
          \begin{equation*}
            \xymatrix{
            P \ar[r] \ar[d] & X\ar[d]\\
            \di{n,2}Z \ar[r] & \di{n,2}\q{2,n}X
            }
          \end{equation*}
          Then $P\in\cathd{n}$.
    \end{itemize}
\end{lemma}
\begin{proof}
\nid By induction on $n$.

\nid In the case $n=1$ for a), since $dB$ is discrete, $A\rw dB$ is an isofibration, therefore
\begin{equation*}
   Q=A\tiund{dB}dC \simeq    A\tms{ps}{dB}dC \simeq A^d\tiund{dB}dC\;.
\end{equation*}
Hence $Q\in\cathd{}$. The case $n=2$ for b) is Lemma \ref{lem-spec-pulbk-eqr}.

\bk
Suppose, inductively, that the lemma holds for $(n-1)$.

\nid a) for each $k\geq 0$ there is a pullback in $\funcat{n-2}{\Cat}$
\begin{equation*}
    \xymatrix{
    Q_k \ar[r] \ar[d] & A_k\ar[d]\\
            \di{n-1,1}C_k \ar[r] & \di{n-1,1}B_k
    }
\end{equation*}
Therefore, by inductive hypothesis a), $Q_k\in\cathd{n-1}$. For each $\ur=(r_1,...,r_{n-1})\in\dop{n-1}$, we have a pullback in $\Cat$
\begin{equation*}
    \xymatrix{
    Q_{\ur} \ar[r] \ar[d] & A_{\ur}\ar[d]\\
            d C \ar[r] & d B
    }
\end{equation*}
Since $p$ commutes with fiber products over discrete objects, we have a pullback in $\Set$
\begin{equation*}
    \xymatrix{
    p Q_{\ur} \ar[r] \ar[d] & p A_{\ur}\ar[d]\\
             C \ar[r] &  B
    }
\end{equation*}
It follows that there is a pullback in $\funcat{n-2}{\Cat}$
\begin{equation*}
    \xymatrix{
    \p{n} Q \ar[r] \ar[d] & \p{n} A\ar[d]\\
            \di{n-1,1} C \ar[r] & \di{n-1,1} B
    }
\end{equation*}
By inductive hypothesis a) we conclude that $\p{n}Q\in\cathd{n}$. By definition, this means that $Q\in\cathd{n}$.
\bk

\nid b) For each $k\geq 0$, there is a pullback in $\funcat{n-2}{\Cat}$
\begin{equation*}
    \xymatrix{
    P_k \ar[r] \ar[d] & X_k \ar[d]\\
            \di{n-1,1} Z_k \ar[r] & \di{n-1,1}\q{1,n-1}X_k
    }
\end{equation*}
Since $X_k\in\cathd{n-1}$ (as $X\in\cathd{n}$), by part a) this implies that $P_k\in\cathd{n-1}$. Since $\p{n}$ commutes with fiber products over discrete objects, we also have a pullback in $\funcat{n-2}{\Cat}$
\begin{equation*}
    \xymatrix{
    \p{n} P \ar[r] \ar[d] & \p{n} X \ar[d]\\
            \di{n-1,2}\q{2,n} Z \ar[r] & \di{n-1,2}\q{2,n}X
    }
\end{equation*}
where $\p{n}X\in\cathd{n-1}$. By inductive hypothesis b), $\p{n}P\in\cathd{n-1}$. Hence by definition, $P\in\cathd{n}$.
\end{proof}
\begin{proposition}\label{pro-gen-const-2}
    For each $n\geq 1$, there is a functor
     \begin{equation*}
       V_n:\cathd{n}\rw\cathd{n}
     \end{equation*}
      with a map
      \begin{equation*}
        v_{n}(X):V_n X\rw X
      \end{equation*}
       natural in $X\in\cathd{n}$ such that
    \begin{itemize}
      \item [a)] $v_{n}(X)$ is a levelwise isofibration in $X\in\cathd{n}$ which is surjective on objects and, for all $1<r\leq n$, the same holds for $\p{r,n} v_{n}(X)$ and $\q{r,n} v_{n}(X)$ for all $1<r\leq n$.\bk

      \item [b)] $V_n$ is identity on discrete objects and preserves pullbacks over discrete objects.\bk

      \item [c)] If $h:X\rw Y$ is a morphism in $\cathd{n}$, the following diagram  commutes for appropriate choices of sections to the discretization maps $V_n X \rw (V_n X)^d$ and $V_n Y \rw (V_n Y)^d$.
          \begin{equation*}
            \xymatrix{
            V_n X \ar[rr] && V_n Y\\
            (V_n X)^d \ar[rr]\ar[u] && (V_n Y)^d\ar[u]
            }
          \end{equation*}
    \end{itemize}
\end{proposition}
\begin{proof}
By induction on $n$. For $n=1$, let $V_1 X=\Dec X$ and $v_{1}(X)=u_X:\Dec X\rw X$ as in Section \ref{decalage}. By Lemma \ref{lem-gen-const-5}, $v_{1}(X)$ is an isofibration and is surjective on objects. Also $\Dec$ preserves pullbacks. Given a morphism $h:X\rw Y$ in $\cathd{}$, we have a commutative diagram
\begin{equation*}
    \xymatrix{
    \Dec X \ar^{\Dec  h}[rr]\ && \Dec Y \\
    (\Dec X)^d=d X_0 \ar^{dh_0}[rr]\ar[u] && (\Dec Y)^d=d Y_0\ar[u]
    }
\end{equation*}
where $dX_0 \rw \Dec X$ and $dY_0 \rw \Dec Y$ are the functorial sections to the discretization maps as in Section \ref{decalage}.This proves the lemma in the case $n=1$. Suppose, inductively, that it holds for $(n-1)$ and let $X\in\cathd{n}$.\bk

\nid a) Let
 \begin{equation*}
  F_n X=X(v_{n-1}(X_0))
  \end{equation*}
  where $v_{n-1}(X_0):V_{n-1}X_0 \rw X_0$ is given by the inductive hypothesis applied to $X_0\in\cathd{n-1}$ and $X(v_{n-1}(X_0))$ is as in Lemma \ref{lem-gen-constr}. By inductive hypothesis a), $v_{n-1}(X_0)$ satisfies the hypothesis of Proposition \ref{pro-gen-const-1-new} and thus $F_n X\in\cathd{n}$.

We define $V_n X$ via the pullback in $\funcat{n-1}{\Cat}$
\begin{equation}\label{eq-Vn}
    \xymatrix{
    V_n X \ar^{h_n(X)}[rr]\ar^{l}[d] && F_n X\ar^{r}[d]\\
    \di{n,2}\Dec\q{2,n}F_n X\ar^{u'}[rr] && \di{n,2}\q{2,n}F_n X
    }
\end{equation}
where
 \begin{equation*}
   u'=u_{(\q{2,n}F_n X)}
 \end{equation*}
  Since $F_n X\in\cathd{n}$, $\q{2,n}F_n X\in\cathd{}$, hence $\Dec \q{2,n}F_n X\in\cathd{}$. Thus by Lemma \ref{lem-gen-const-6}, $V_n X\in\cathd{n}$.

For each $\uk\in\dop{n-1}$, there is a pullback in $\Cat$
\begin{equation}\label{eq4-gen-const}
    \xymatrix{
    (V_n X)_{\uk} \ar^{(h_n (X))_{\uk}}[rr]\ar[d] && (F_n X)_{\uk}\ar^{r_{\uk}}[d]\\
    d(\di{n,2}\Dec\q{2,n}F_n X)_{\uk} \ar^{u'_{\uk}}[rr] && d(\di{n,2}\q{2,n}F_n X)_{\uk}\;.
    }
\end{equation}
The bottom horizontal map is an isofibration since the target is discrete; hence, since isofibrations are stable under pullback, $(h_n (X))_{\uk}$ is also an isofibration. The bottom horizontal map in \eqref{eq4-gen-const} is also surjective on objects since
\begin{equation*}
   ob(u'_{\uk}): (N\Dec \q{2,n}F_n X)_{k_1} \rw (N\q{2,n}F_n X)_{k_1}
\end{equation*}
is surjective for all $k_1\geq 0$, where $N:\Cat\rw\funcat{}{\Set}$ is the nerve functor. Since $ob((h_n (X))_{\uk})$ is the pullback of $ob(u'_{\uk})$, it follows that $(h_n (X))_{\uk}$ is also surjective on objects.

 By inductive hypothesis a) the map $v_{n-1}(X_0):V_{n-1}X_0 \rw X_0$ satisfies the hypothesis of Proposition \ref{pro-gen-const-1-new}, and thus the map (which we now denote by $f_n(X)$ instead of $V(X)$),
  \begin{equation*}
  f_{n}(X):F_n X=X(v_{n-1}(X_0)) \rw X
  \end{equation*}
   is a levelwise isofibration surjective on objects. We conclude from above that the same holds for the composite map
\begin{equation*}
    v_{n}(X):V_n X \xrw{\;\;h_n (X)\;\;}F_n X\xrw{f_{n}(X)} X\;.
\end{equation*}
We now show that $\p{n,r}v_{n}(X)$ is a levelwise isofibration surjective on objects for all $r$. Since, by Proposition \ref{pro-gen-const-1-new} c), this holds for $\p{r,n}V(X)$, it is sufficient to show this for $\p{r,n}h_n (X)$.

Since $\p{r,n}$ commutes with pullbacks over discrete objects, we have a pullback in $\funcat{r-1}{\Cat}$
\begin{equation*}
    \xymatrix{
    p\up{r,n}V_n X \ar^{p\up{r,n}h_n (X)}[rr]\ar[d] && \p{r,n}F_n X\ar[d]\\
    \di{r,2}\Dec q\up{2,n}F_n X\ar[rr] && \di{r,2} q\up{2,n}F_n X\;.
    }
\end{equation*}
Using a similar argument as above we conclude that $\p{r,n}h_n (X)$ is a levelwise isofibration surjective on objects.

Finally, by inductive hypothesis a) the map $v_{n-1}(X_0):V_{n-1}X_0 \rw X_0$ satisfies the hypothesis of Proposition \ref{pro-gen-const-1-new} e), and thus the map $q\up{r,n}f_n(X)$ is a levelwise isofibration surjective on objects. To show that the same holds for $\q{r,n}v_{n}(X)$ it is therefore enough to show it for $\q{r,n}h_n (X)$.

Since $\q{r,n}$ commutes with pullbacks over discrete objects, we have a pullback in $\funcat{r-1}{\Cat}$
\begin{equation*}
    \xymatrix{
    q\up{r,n}V_n X \ar^{\p{r,n}h_n (X)}[rr]\ar[d] && q\up{r,n}F_n X\ar[d]\\
    \di{r,2}\Dec\q{2,n}F_n X\ar[rr] && \di{r,2}\q{2,n}F_n X\;.
    }
\end{equation*}
Using a similar argument as above we conclude that $\q{r,n}h_n (X)$ is a levelwise isofibration surjective on objects. This proves a).\bk

\nid b) If $X$ is discrete, $F_n X=X=\di{n,2}\q{2,n}X$, thus $V_n X=X$. Since, by inductive hypothesis, $V_{n-1}$ commutes with pullbacks over discrete objets, so does $F_n$ as easily seen. Since $\q{2,n}$ commutes with pullbacks over discrete objects and $\Dec$ commutes with pullbacks, it follows by construction that $V_n$ commutes with pullbacks over discrete objects.\bk

\nid c) By the definition of $V_n(X)$ in \eqref{eq-Vn} in order to construct the map
\begin{equation*}
  \zg'_{V_n X}: (V_n X)^d\rw V_n X
\end{equation*}
we need to construct maps
\begin{equation*}
  t: (V_n X)^d \rw \di{n-2}\Dec \q{2,n}F_n X,\qquad s: (V_n X)^d \rw F_n X
\end{equation*}
such that $u't=rs$, where $u', r$ are as in \eqref{eq-Vn}. By \eqref{eq-Vn} , since $\q{2,n}$ commutes with pullbacks over discrete objects,
 \begin{equation*}
 \q{2,n}V_n X=\Dec \q{2,n}F_n X,
 \end{equation*}
 so that
\begin{equation*}
\begin{split}
    & \q{1,n}V_n X=q\Dec \q{2,n}F_n X=(\q{2,n}F_n X)_0=\\
    & =\q{1,n-1}(F_n X)_0=\q{1,n-1}V_{n-1}X_0\;.
\end{split}
\end{equation*}
It follows that
\begin{equation}\label{eq5-gen-const}
  (V_n X)^d=\di{n,2}(\q{2,n}V_n X)^d=\di{n,2}(\Dec\q{2,n}F_n X)^d
\end{equation}
as well as
\begin{equation}\label{eq6-gen-const}
  (V_n X)^d=\di{1}(F_n X)_0^d=\di{1}(V_{n-1} X_0)^d
\end{equation}
where we denote $\di{1}:\cat{n-1}\xrw{\;d\;}\Cat(\cat{n-1})\oset{\xi_1^{-1}}{\cong}\cat{n}$.

There is a natural map
\begin{equation*}
  (\Dec\q{2,n}F_n X)^d\rw\Dec\q{2,n}F_n X
\end{equation*}
and therefore, using \eqref{eq5-gen-const} a corresponding map
\begin{equation*}
  t: (V_n X)^d=\di{n,2}(\Dec\q{2,n}F_n X)^d\rw\di{n,2}\Dec\q{2,n}F_n X\;.
\end{equation*}
Note that $t$ is also natural in $X$. The composite map
\begin{equation*}
  (\Dec \q{2,n}F_n X)^d = d(\q{2,n}F_n X)_0 \rw \Dec\q{2,n}F_n X\rw \q{2,n}F_n X
\end{equation*}
is the counit $\za$ of the adjunction $d \dashv \ob$ (see Remark \ref{rem-obj-fun}) at $(\q{2,n}F_n X)_0$, so $u't$ is the corresponding map
\begin{equation}\label{eq7-gen-const}
  (V_n X)^d = \di{n,2}d(\q{2,n}F_n X)_0 \xrw{\di{n,2}\za_{(\q{2,n}F_n X)}} \di{n,2}(\q{n,2}F_n X)\;.
\end{equation}
Using \eqref{eq6-gen-const} and the inductive hypothesis on $X_0$, we obtain a natural map
\begin{equation*}
  (F_n X)_0^d = (V_{n-1}X_0)^d \xrw{\;v\;} V_{n-1}X_0 =(F_n X)_0\;.
\end{equation*}
We define $s$ to be the composite
\begin{equation*}
  s: (V_{n}X)^d = \di{1}(V_{n-1}X_0)^d \xrw{\di{1} v} \di{1}(F_n X)_0 \xrw{\za_{F_n X}}F_n X\;.
\end{equation*}
Since $\za$ is natural in $X$ (as it is the counit of the adjunction $d \dashv \ob$) and $v$ is natural by induction hypothesis, $s$ is also natural in $X$.

The map $r$ in \eqref{eq-Vn} is natural (since $J_n r$ is levelwise unit of the adjunction $q\dashv d$), therefore we have a commuting diagram
\begin{equation*}
\xymatrix@R=40pt@C=50pt{
\di{1}(F_n X)_0 \ar^{\za_{F_n X}}[r] \ar_{r'}[d] & F_n X \ar^{r}[d]\\
(V_n X)^d = \di{n,2}\q{2,n} \di{1}(F_n X)_0 \ar[r] & \di{n,2}\q{2,n}F_n X
}
\end{equation*}
where the equality on the bottom left corner follows from \eqref{eq6-gen-const}.

Thus $r'=\di{1}\zg_{(F_n X)_0}$ where
 \begin{equation*}
   \zg_{(F_n X)_0}: (F_n X)_0 \rw (F_n X)^d_0
 \end{equation*}
 is the discretization map; $\di{1}v$ is a section for $r'$ since $v$ is a section for $\zg_{(F_n X)_0}$. Hence
\begin{equation*}
  rs=r\za_{F_n X}\di{1}v=(\di{n,2}\za_{\q{2,n}F_n X})r'(\di{1}v)=\di{n,2}\za_{\q{2,n}F_n X}\;.
\end{equation*}
By \eqref{eq7-gen-const} we conclude that $u's=rt$.

Since $(V_n X)^d=\di{n,2}(\Dec \q{2,n} F_n X)^d$, the discretization map
\begin{equation*}
\zg_{V_n X}: V_n X\rw (V_n X)^d
\end{equation*}
is the composite
\begin{equation*}
  V_n X\xrw{\;\ell\;} \di{n,2}\Dec\q{2,n}F_n X\xrw{\di{n,2}\zg}\di{n,2}(\Dec \q{2,n}F_n X)^d
\end{equation*}
where $\zg$ is the discretization map for $\Dec \q{2,n}F_n X$ and $l$ is as in \eqref{eq-Vn}. Thus
\begin{equation*}
  \zg_{Vn X}\zg'_{Vn X}=\di{n,2}\zg \ell\zg'_{Vn X}=(\di{n,2}\zg)t=\Id\;.
\end{equation*}
As observed above, by construction both $t$ and $s$ are natural in $X$, so by elementary properties of pullbacks, such is $\zg'_{Vn X}$. Thus, given a morphism $h:X\rw Y$ in $\cathd{n}$, there is a commuting diagram
\begin{equation*}
\xymatrix@R=40pt@C=50pt{
V_n X \ar^{V_n h}[r] & V_n Y\\
(V_n X)^d \ar_{(V_n h)^d}[r]\ar^{\zg'_{Vn X}}[u] & (V_n Y)^d\ar_{\zg'_{Vn Y}}[u]
}
\end{equation*}
as required.

\end{proof}
In the following Corollary we establish the properties of $V_n$ with respect to pullbacks.
\begin{corollary}\label{cor-gen-const-4}
  Let $X\rw Z\lw Y$ be a diagram in $\cathd{n}$ such that $X\tiund{Z}Y\in\cathd{n}$ and such that, for all $1\leq r\leq n$
  \begin{equation}\label{eq-cor-gen-const-4-1}
    \p{r,n}(X\tiund{Z}Y)=\p{r,n}X\tiund{\p{r,n}Z}\p{r,n}Y\;,
  \end{equation}
  \begin{equation}\label{eq-cor-gen-const-4-2}
    \q{r,n}(X\tiund{Z}Y)=\q{r,n}X\tiund{\q{r,n}Z}\q{r,n}Y\;.
  \end{equation}
  Then
  \begin{itemize}

    \item [a)] For all $n$,
    \begin{equation*}
     \begin{split}
     & V_n(X\tiund{Z}Y)\cong V_n X\tiund{V_n Z}V_n Y, \\
     & v_n(X\tiund{Z}Y)\cong v_n (X)\tiund{v_n (Z)}v_n (Y).
   \end{split}
  \end{equation*}

    \item [b)] For  all $1\leq r\leq n$,
    \begin{equation*}
    \begin{split}
        & \p{r,n}\{V_n(X\tiund{Z}Y)\} = \p{r,n}(V_n X)\tiund{\p{r,n}(V_n Z)} \p{r,n}(V_n Y)\\
        & \q{r,n}\{V_n(X\tiund{Z}Y)\} = \q{r,n}(V_n X)\tiund{\q{r,n}(V_n Z)} \q{r,n}(V_n Y)
    \end{split}
    \end{equation*}
  \end{itemize}
\end{corollary}
\begin{proof}
\

a) Since $V_1=\Dec$ commutes with pullbacks, a) holds for $n=1$. Suppose, inductively, that it holds for $(n-1)$. We claim that
\begin{equation}\label{eq-cor-gen-const-4-3}
\begin{split}
  & F_n(X\tiund{Z}Y)=F_n X\tiund{F_n Z}F_n Y\;,\\
  & f_n(X\tiund{Z}Y)=f_n (X)\tiund{f_n (Z)}f_n (Y)\;.
  \end{split}
\end{equation}
In fact, by construction
\begin{equation*}
   F_n(X\tiund{Z}Y)=(X\tiund{Z}Y)(v_{n-1}(X_0\tiund{Z_0}Y_0))
\end{equation*}
where
\begin{equation*}
  v_{n-1}(X_0\tiund{Z_0}Y_0): V_{n-1}(X_0\tiund{Z_0}Y_0)\rw X_0\tiund{Z_0}Y_0\;.
\end{equation*}
Observe that the diagram in $\cathd{n-1}$ $X_0\rw Z_0\lw Y_0$ satisfies the inductive hypothesis. In fact, since $X\tiund{Z}Y\in\cathd{n}$,
 \begin{equation*}
 (X\tiund{Z}Y)_0=X_0\tiund{Z_0}Y_0\in\cathd{n-1}
 \end{equation*}
  while hypotheses \eqref{eq-cor-gen-const-4-1}, \eqref{eq-cor-gen-const-4-2} imply
\begin{equation*}
\begin{split}
    & \p{r-1,n-1}(X_0\tiund{Z_0}Y_0) = \{\p{r,n}(X\tiund{Z}Y)\}_0= \\
    & =\{\p{r,n} X\tiund{\p{r,n} Z} \p{r,n} Y)\}_0=\{\p{r-1,n-1} X_0\tiund{\p{r-1,n-1} Z_0} \p{r-1,n-1} Y_0)\}
\end{split}
\end{equation*}
and similarly for $\q{r-1,n-1}$. By induction hypothesis we deduce
\begin{equation*}
\begin{split}
    & V_{n-1}(X_0\tiund{Z_0}Y_0)=V_{n-1}X_0\tiund{V_{n-1}Z_0}V_{n-1}Y_0 \\
    & v_{n-1}(X_0\tiund{Z_0}Y_0)=v_{n-1}(X_0)\tiund{v_{n-1}(Z_0)}v_{n-1}(Y_0)\;.
\end{split}
\end{equation*}
By Corollary \eqref{cor-gen-const-1-pb} a) it follows that
\begin{equation*}
\begin{split}
    & F_n(X\tiund{Z}Y)=(X\tiund{Z}Y)(v_{n-1}(X_0\tiund{Z_0}Y_0))= \\
    & =X(v_{n-1}(X_0))\tiund{Z(v_{n-1}(Z_0))} Y(v_{n-1}(Y_0))= F_n X\tiund{F_n Z}F_n Y\;
\end{split}
\end{equation*}
and
\begin{equation*}
  f_n(X\tiund{Z}Y)=f_n (X)\tiund{f_n (Z)}f_n (Y)\;.
\end{equation*}
which is \eqref{eq-cor-gen-const-4-3}.

By Proposition \ref{pro-gen-const-2} a) the following maps
\begin{equation*}
\begin{split}
    & \q{2,n}f_0,\quad \q{2,n}g_0,\quad \q{2,n}h_0\\
    & \p{2,n}f_0,\quad \p{2,n}g_0,\quad \p{2,n}h_0
\end{split}
\end{equation*}
are isofibrations surjective on objects. Together with hypotheses \eqref{eq-cor-gen-const-4-1} and \eqref{eq-cor-gen-const-4-2} this means that the hypotheses of Corollary \ref{cor-gen-const-1-pb} c) are satisfied and therefore we conclude that
\begin{equation}\label{eq-cor-gen-const-4-4}
\begin{split}
    & \q{2,n}F_n(X\tiund{Z}Y)=\q{2,n}\{(X\tiund{Z}Y)(f_{X_0}\tiund{f_{Z_0}}f_{Y_0})\}\cong \\
    & \cong \q{2,n}\{X(f_{X_0})\}\tiund{\q{2,n}\{Z(f_{Z_0})\}}\q{2,n}\{Y(f_{Y_0})\}\cong\\
    & \cong \q{2,n}F_n X\tiund{\q{2,n}F_n Z} \q{2,n}F_n Y\;.
\end{split}
\end{equation}
Since $\Dec$ commutes with pullbacks, it also follows that
\begin{equation}\label{eq-cor-gen-const-4-4a}
\begin{split}
    & \Dec \q{2,n}F_n(X\tiund{Z}Y)=\Dec(\q{2,n}F_n X\tiund{\q{2,n}F_n Z} \q{2,n}F_n Y)=\\
    & =\Dec\q{2,n}F_n X\tiund{\Dec\q{2,n}F_n Z} \Dec\q{2,n}F_n Y\;.
\end{split}
\end{equation}
Since, by construction,
\begin{equation*}
  V_n(X\tiund{Z}Y)=\di{n,2}\Dec\q{2,n}F_n(X\tiund{Z}Y) \tiund{\di{n,2}\q{2,n}F_n(X\tiund{Z}Y)} F_n(X\tiund{Z}Y)
\end{equation*}
using \eqref{eq-cor-gen-const-4-3}, \eqref{eq-cor-gen-const-4-4}, \eqref{eq-cor-gen-const-4-4a} and the commutation of pullbacks we conclude that
\begin{equation*}
  V_n(X\tiund{Z}Y)= V_n X\tiund{V_n Z}V_n Y\;.
\end{equation*}
as well as
\begin{equation*}
  h_n(X\tiund{Z}Y)= h_n (X)\tiund{h_n (Z)}h_n (Y)\;.
\end{equation*}
Since, by definition, $v_n (X)= h_n (X) f_n (X)$ by \eqref{eq-cor-gen-const-4-3} we deduce that
\begin{equation*}
  v_n(X\tiund{Z}Y)= v_n (X)\tiund{v_n (Z)}v_n (Y)\;.
\end{equation*}
This concludes the proof of a).

\bk
b) By \eqref{eq-cor-gen-const-4-3} and Corollary \ref{cor-gen-const-1-pb} we have
\begin{equation}\label{eq-cor-gen-const-4-5}
\begin{split}
    & \p{r,n}\{F_n (X\tiund{Z}Y)\}=\p{r,n}\{(X\tiund{Z}Y)(f_{X_0}\tiund{f_{Z_0}}f_{Y_0})\}= \\
    & =\p{r,n}\{X(f_{X_0})\} \tiund{\p{r,n}\{Z(f_{Z_0})\}} \p{r,n}\{Y(f_{Y_0})\}= \\
    & =\p{r,n}F_n X \tiund{\p{r,n}F_n Z} \p{r,n}F_n Y\;.
\end{split}
\end{equation}
Similarly, using Corollary \ref{cor-gen-const-1-pb} c) one shows
\begin{equation}\label{eq-cor-gen-const-4-6}
 \q{r,n}\{F_n (X\tiund{Z}Y)\}= \q{r,n}F_n X \tiund{\q{r,n}F_n Z} \q{r,n}F_n Y\;.
\end{equation}
In particular this holds for $r=2$ from which we deduce
\begin{equation}\label{eq-cor-gen-const-4-7}
  \Dec\q{2,n}\{F_n (X\tiund{Z}Y)\}= \Dec\q{2,n}F_n X \tiund{\Dec\q{2,n}F_n Z} \Dec\q{2,n}F_n Y\;.
\end{equation}
From the definition of $V_n X$, since $\p{r,n}$ commutes with pullbacks over discrete objects, we have a pullback in $\funcat{r-2}{\Cat}$
\begin{equation*}
\xymatrix@R=40pt@C=50pt{
\p{r,n}V_n X \ar[r] \ar[d] & \p{r,n}F_n X \ar[d]\\
\di{r,2}\Dec \q{2,n}F_n X \ar[r] & \di{r,2} \q{2,n}F_n X
}
\end{equation*}
So in particular
\begin{equation*}
\begin{split}
    & \p{r,n}V_n (X\tiund{Z}Y)= \\
    & = \di{r,2}\Dec \q{2,n} \{ F_n(X\tiund{Z}Y) \} \tiund{\di{r,2}\Dec \q{2,n} F_n(X\tiund{Z}Y)} \p{r,n}\{ F_n(X\tiund{Z}Y) \}\;.
\end{split}
\end{equation*}
Using \eqref{eq-cor-gen-const-4-5}, \eqref{eq-cor-gen-const-4-6}. \eqref{eq-cor-gen-const-4-7} and the commutation of pullbacks we deduce
\begin{equation*}
  \p{r,n}V_n (X\tiund{Z}Y)= \p{r,n}V_n X \tiund{\p{r,n}V_n Z} \p{r,n}V_n Y\;.
\end{equation*}
Similarly, since $\q{r,n}$ commutes with pullbacks over discrete objects, we have a pullback in $\funcat{r-2}{\Cat}$
\begin{equation*}
\xymatrix@R=40pt@C=50pt{
\q{r,n}V_n X \ar[r] \ar[d] & \q{r,n}F_n X \ar[d]\\
\di{r,2}\Dec \q{2,n}F_n X \ar[r] & \di{r,2} \q{2,n}F_n X
}
\end{equation*}
Using \eqref{eq-cor-gen-const-4-6}. \eqref{eq-cor-gen-const-4-7} and the commutation of pullbacks we deduce
\begin{equation*}
  \q{r,n}V_n (X\tiund{Z}Y)= \q{r,n}V_n X \tiund{\q{r,n}V_n Z} \q{r,n}V_n Y\;.
\end{equation*}
\end{proof}
In the proof of Proposition \ref{pro-gen-const-2} we defined a functor
\begin{equation*}
        F_n:\cathd{n}\rw \cathd{n}
    \end{equation*}
    by
\begin{equation*}
  F_n X=X(v_{n-1}(X_0)).
\end{equation*}
where the map  $v_{n-1}(X_0):V_{n-1}X_0 \rw X_0$ in $\cathd{n-1}$ is as in Proposition \ref{pro-gen-const-2}. We now extend this to a functor
\begin{equation*}
        F_n:\catwg{n}\rw \catwg{n}
    \end{equation*}
 defined by
\begin{equation*}
  F_n X=X(v_{n-1}(X_0)).
\end{equation*}
where $v_{n-1}(X_0):V_{n-1}X_0 \rw X_0$ is as in Proposition \ref{pro-gen-const-2}.

\begin{proposition}\label{pro-gen-const-1}
    For each $n\geq 2$, there is a functor
    \begin{equation*}
        F_n:\catwg{n}\rw \catwg{n}
    \end{equation*}
    and a map
    \begin{equation*}
      f_n(X):F_n X\rw X
    \end{equation*}
     natural in $X\in\catwg{n}$, such that
    \begin{itemize}
      \item [i)] $f_n(X)$ is an $n$-equivalence.\bk

      \item [ii)] $F_n$ is identity on discrete objects and preserves pullbacks over discrete objects.\bk

      \item [iii)] If $f:X\rw Y$ is a morphism in $\catwg{n}$ the following diagram commutes for appropriate choices of sections to the discretization maps $(F_n X)_0 \rw (F_n X )^d_0$ and $(F_n Y)_0 \rw (F_n Y )^d_0$.
          \begin{equation*}
            \xymatrix{
            (F_n X)_0 \ar[rr] && (F_n Y)_0\\
            (F_n X )^d_0 \ar[rr]\ar[u]  && (F_n Y)^d_0\ar[u]\;.
            }
          \end{equation*}\bk

          \item [iv)] If $X\in \cathd{n}$, $F_n X\in \cathd{n}$.
    \end{itemize}
\end{proposition}
\begin{proof}
Given $X\in\catwg{n}$, since $X_0\in\cathd{n-1}$ by Proposition \ref{pro-gen-const-2} there  is a map
\begin{equation*}
    v_{n-1}(X_0):V_{n-1}X_0 \rw X_0
\end{equation*}
such that $v_{n-1}(X_0)$ and $\p{r,n-1}v_{n-1}(X_0)$ are levelwise isofibrations in $\Cat$ surjective on objects for all $1<r\leq (n-1)$. Let
 \begin{equation*}
  F_n X=X(v_{n-1}(X_0)).
\end{equation*}
 By Proposition \ref{pro-gen-const-1-new}, $F_nX\in\catwg{n}$ and there is an $n$-equivalence which we now denote by $f_n(X)$,
  \begin{equation*}
f_n(X)=V(X):F_n X\rw X
\end{equation*}
proving i).

  If $X$ is discrete, so is $X_0$, thus by Proposition \ref{pro-gen-const-2} $ v_{n-1}(X_0)=\Id$ and therefore $F_n X=X$.

Let $X\rw Z \lw Y$ be a pullback in $\catwg{n}$ with $Z$ discrete. By Proposition \ref{pro-gen-const-2},
\begin{equation*}
  V_{n-1}(X_0\tiund{Z}Y_0)= V_{n-1}X_0\tiund{Z}V_{n-1}Y_0
\end{equation*}
and therefore, as easily checked,
\begin{equation*}
    (F_n(X\tiund{Z}Y))_1=(F_n X)_1\tiund{Z}(F_n Y)_1\;.
\end{equation*}
It follows that
\begin{equation*}
    F_n(X\tiund{Z}Y)=F_n X\tiund{Z}F_n Y
\end{equation*}
which is ii).

 Since $(F_n X)_0=V_{n-1}X_0$, iii) follows from Proposition \ref{pro-gen-const-2}. Since, by i), $f_{n}(X)$ is an $n$-equivalence, iv) follows from Proposition \ref{pro-nequiv-to-obj}.
\end{proof}

We next study the behaviour of $F_n$ with respect to pullbacks.

\begin{corollary}\label{cor-fn-pull}
Let $X\rw Z \lw Y$ be a diagram in $\catwg{n}$ such that $X\tiund{Z}Y \in \catwg{n}$ and such that, for each $1\leq r\leq n$
\begin{equation}\label{eq-cor-fn-pull-1}
  \p{r,n} (X\tiund{Z}Y)= \p{r,n} X \tiund{\p{r,n} Z} \p{r,n} Y\;.
\end{equation}
\begin{equation}\label{eq-cor-fn-pull-2}
  \q{r,n} (X\tiund{Z}Y)= \q{r,n} X \tiund{\q{r,n} Z} \q{r,n} Y\;.
\end{equation}
Then
\begin{itemize}
  \item [a)] For all $n$,
  \begin{equation*}
  \begin{split}
    & F_n(X\tiund{Z}Y)\cong F_n X \tiund{F_n Z} F_n Y\\
    & f_n(X\tiund{Z}Y)\cong f_n (X) \tiund{f_n (Z)} f_n (Y).
    \end{split}
  \end{equation*}
  \item [b)] For all $1\leq r \leq n$,
  \begin{equation*}
  \p{r,n} \{F_n(X\tiund{Z}Y)\}\cong \p{r,n} F_n X \tiund{\p{r,n}F_n Z} \p{r,n}F_n Y\;,
  \end{equation*}
  \begin{equation*}
  \q{r,n} \{F_n(X\tiund{Z}Y)\}\cong \q{r,n} F_n X \tiund{\q{r,n}F_n Z} \q{r,n}F_n Y\;.
  \end{equation*}
\end{itemize}
\end{corollary}
\begin{proof}
\

a) By construction $F_n X=X(f_{X_0})$ where $ v_{n-1}(X_0):V_{n-1}X_0\rw X_{0}$. Since $X\tiund{Z}Y\in\catwg{n}$, $X_0\tiund{Z_0}Y_0\in\cathd{n-1}$ and it satisfies the hypotheses of Corollary \ref{cor-gen-const-4}. Therefore
\begin{equation*}
  V_{n-1}(X_0\tiund{Z_0}Y_0)\cong V_{n-1}X_0\tiund{V_{n-1}Z_0}V_{n-1}Y_0
\end{equation*}
so that
\begin{equation*}
  v_{n-1}(X_0\tiund{Z_0}Y_0)=v_{n-1}(X_0)\tiund{v_{n-1}(Z_0)}v_{n-1}(Y_0)\;.
\end{equation*}
By Corollary \ref{cor-gen-const-1-pb} a) it follows that
\begin{equation*}
\begin{split}
    & F_n(X\tiund{Z}Y)=(X\tiund{Z}Y)(v_{n-1}(X_0\tiund{Z_0}Y_0))= \\
    & =X(v_{n-1}(X_0))\tiund{Z(v_{n-1}(Z_0))} Y(v_{n-1}(Y_0))= F_n X\tiund{F_n Z}F_n Y\;
\end{split}
\end{equation*}
and
\begin{equation*}
f_n(X\tiund{Z}Y)\cong f_n (X) \tiund{f_n (Z)} f_n (Y).
\end{equation*}

b) By Proposition \ref{pro-gen-const-2} a) the maps
\begin{equation*}
v_{n-1}(X_0),\quad v_{n-1}(Y_0),\quad v_{n-1}(Z_0),\quad v_{n-1}(X_0\tiund{Z_0}Y_0)
\end{equation*}
 satisfy the hypotheses of Corollary \ref{cor-gen-const-1-pb} b). Further, by hypothesis \eqref{eq-cor-fn-pull-1}
\begin{equation*}
  \p{r-1,n-1}(X_0\tiund{Z_0}Y_0)= \p{r-1,n-1} X_0\tiund{\p{r-1,n-1} Z_0} \p{r-1,n-1} Y_0\;.
\end{equation*}
Thus all the hypotheses of Corollary \ref{cor-gen-const-1-pb} b) are satisfied and we conclude that
\begin{equation*}
\begin{split}
    & \p{r,n} \{F_n(X\tiund{Z}Y)\} = \p{r,n} \{(X\tiund{Z}Y)(v_{n-1}(X_0)\tiund{v_{n-1}(Z_0)}v_{n-1}(Y_0))\}= \\
    & = \p{r,n}\{X(v_{n-1}(X_0))\}\tiund{\p{r,n}\{Z(v_{n-1}(Z_0))\}} \p{r,n}\{Y(v_{n-1}(Y_0))\}=\\
    & = \p{r,n} F_n X \tiund{\p{r,n}F_n Z} \p{r,n}F_n Y\;.
\end{split}
\end{equation*}
The proof for $\q{r,n}$ is similar.

\end{proof}
%%
%%
%%%%%%%%%%%%%%%%%%%%%%%%%%%%%%%%%%%%%%%%%%%%%%%%%%%%%%%%%%%%%%%%%%%%%%%%%%%%%%%%%%%
\section{The category $\pmb{\ftawg{n}}$}\label{sec-fta}

In this section we introduce the category $\ftawg{n}$. We then show in Theorem \ref{pro-fta-1} that there is a functor
\begin{equation*}
        G_n:\catwg{n}\rw\ftawg{n}
    \end{equation*}
    and an $n$-equivalence, $G_nX \rw X$ for each $X\in \catwg{n}$. The construction of $G_n$ is inductive and uses the functor $F_n$ built in Section \ref{sec-canonical}. Namely, we define $G_2=F_2$ and given $G_{n-1}$,
     \begin{equation*}
       G_n=\ovl{G}_{n-1}\circ F_n
     \end{equation*}
    see Definition \ref{def-Gn-new} and Theorem \ref{pro-fta-1}.
    In the next chapter, $G_n$ will be used to define the discretization functor $\Discn:\catwg{n}\rw \ta{n}$.

\subsection{The idea of the category $\pmb{\ftawg{n}}$}\label{sub-idea-ftawgn}
The idea of the construction of the discretization functor from $\catwg{n}$ to $\ta{n}$ is to replace the homotopically discrete substructures in $X\in\catwg{n}$ by their discretization. As outlined in the introduction to this chapter, this cannot be done in a functorial way unless there are functorial sections to the discretization maps of the homotopically discrete substructures.

For this reason, we introduce in this section a new category $\ftawg{n}$ in which we have functorial sections to the discretization maps for the homotopically discrete substructures.

The idea of the category $\ftawg{n}$ is to modify objects and morphisms in the category $\catwg{n}$ by imposing extra structure giving functorial sections to the discretization maps of the homotopically discrete substructures in the multinerve of objects of $\catwg{n}$. So a map $f: A\rw B$ of homotopically discrete $k$-fold categories (for the appropriate dimension $k$) gives rise to a corresponding commuting diagram
\begin{equation*}
\xymatrix@R=25pt@C=45pt{
A \ar^{f}[r] & B\\
A^d \ar^{\zg'(A)}[u] \ar_{f^d}[r] & B^d \ar_{\zg'(B)}[u]
}
\end{equation*}
where $\zg(A)\zg'(A)=\Id$, $\zg(B)\zg'(B)=\Id$, $\zg(A): A\rw A^d$ and $\zg(B): B\rw B^d$ being the discretization maps.

The choice of the maps $f: A\rw B$ with respect to which we require this functorial behavior is as follows. Given $X\in \catwg{n}$,\: $\uk,\ur \in\Delta^{s^{op}}$, and a morphism $\uk \rw \ur$ in $\Delta^{s^{op}}$ we have a corresponding morphism in $\cathd{n-s-1}$
\begin{equation}\label{eq1-idea-ftwg}
  f: X_{\uk 0}\rw X_{\ur 0}\;.
\end{equation}
In the definition of the category $\ftawg{n}$ (see Definition \ref{def-fta-1}), we impose the commutativity of diagram \eqref{eq-def-fta-1} to require the existence of sections to the discretization maps of $X_{\uk 0}$ and $X_{\ur 0}$ that behave functorially with respect to the map \eqref{eq1-idea-ftwg}.
We next consider the maps
\begin{equation*}
  \zg'(X_0): X_0^d \rw X_0,\qquad \zg'(X_{\uk 0}):X_{\uk 0}^d\rw X_{\uk 0}
\end{equation*}
and the corresponding maps of homotopically discrete objects
\begin{equation}\label{eq2-idea-ftwg}
  (\zg'(X_0))_0: X_0^d \rw X_{00},\qquad (\zg'(X_{\uk 0}))_{\uh 0}:X_{\uk 0}^d\rw X_{\uk 0\uh 0}
\end{equation}
for each $\uh\in\Delta^{t^{op}}$, \; $1\leq t\leq n-3$. We then impose the functoriality of the sections to the discretization maps with respect to the maps \eqref{eq2-idea-ftwg}. That is, we require the commutativity of diagrams \eqref{eq-def-fta-1a} and \eqref{eq-def-fta-1b} in Definition \ref{def-fta-1}. This defines objects of $\ftawg{n}$.

As for morphisms of $\ftawg{n}$, given a morphism $F:X\rw Y$ in $\catwg{n}$ and $\uk\in\Delta^{s^{op}}$ we obtain a map of homotopically discrete structures
\begin{equation}\label{eq3-idea-ftwg}
  F_0:X_0\rw Y_0,\qquad F_{\uk 0}:X_{\uk 0}\rw Y_{\uk 0}
\end{equation}
for all $\uk\in\Delta^{s^{op}}$. We impose functoriality of the sections to the discretization maps with respect to the maps \eqref{eq3-idea-ftwg}. This translates into the commutativity of diagrams \eqref{eq-def-fta-2b} and \eqref{eq-def-fta-2} in Definition \ref{def-fta-1}.

As we will see in the next chapter, the definition of $\ftawg{n}$ is exactly what is needed to functorially discretize the homotopically discrete substructures and thus build a functor $\Dn:\ftawg{n}\rw \ta{n}$.

%%%%%%%%%%%%%%%%%%%%%%%%%%%%%%%%%%%%%%%%%%%%%%%%%%%%%%%%%%%%%%%%%%%%%%%%%%%%%
\subsection{The formal definition of the category $\pmb{\ftawg{n}}$} We now give the formal definition of the category $\ftawg{n}$ and we establish its properties.

\begin{definition}\label{def-fta-1}
    Define the category $\ftawg{n}$ as follows. Let $\ftawg{1}=\Cat$.
    Let $\ftawg{2}$ have the following objects and morphisms:
     \begin{itemize}

      \item [i)] Objects of $\ftawg{2}$ are objects of $\catwg{2}$.
     \item [ii)]A morphism $F:X\rw Y$ in $\ftawg{2}$ is a morphism in $\catwg{2}$ such that the following diagram commutes
          \begin{equation}\label{eq-def-fta-2a}
          \xymatrix@C=50pt@R=30pt{
          X_{0} \ar^{F_{0}}[r] & Y_{0}\\
          X^d_{0} \ar_{F^d_{0}}[r] \ar^{\zg'(X_{0})}[u] & Y^d_{0} \ar_{\zg'(Y_{0})}[u]
          }
          \end{equation}
where $\zg'(X_{0})$ and $\zg'(Y_{0})$ are sections to the discretization maps.

     \end{itemize}
     \bk

     For each $n> 2$, let $\ftawg{n}$ have the following objects and morphisms:

    An object $X$ of $\ftawg{n}$ consists of $X\in\catwg{n}$ such that for all $\uk=(\seqc{k}{1}{s})\in\dop{s}$, $\ur=(\seqc{r}{1}{s})\in\dop{s}$, $(1\leq s \leq n-2)$ and morphisms $\uk\rw\ur$ in $\dop{s}$, the corresponding morphisms
    \begin{equation*}
        f:X_{\uk 0}\rw X_{\ur 0}
    \end{equation*}
    in $\cathd{n-s-1}$ are such that there are choices of sections to the discretization maps
    \begin{align*}
        &\zg(X_{\uk}):X_{\uk 0} \rw X_{\uk 0}^d\\
        &\zg(X_{\ur}):X_{\ur 0} \rw X_{\ur 0}^d\\
    \end{align*}
    making the following diagrams commute

    \begin{itemize}

      \item [i)]
    \begin{equation}\label{eq-def-fta-1}
    \xymatrix@C=50pt@R=30pt{
    X_{\uk 0} \ar^{f}[r] & X_{\ur 0}\\
    X^d_{\uk 0} \ar_{f^d}[r] \ar^{\zg'(X_{\uk 0})}[u] & X^d_{\ur 0} \ar_{\zg'(X_{\ur 0})}[u]
    }
    \end{equation}

   \item [ii)]   For all $\uh=(h_1,\ldots,h_t)\in\Delta^{t^{op}}$, $1\leq t\leq n-3$,
   \begin{equation}\label{eq-def-fta-1a}
   \xymatrix@R=25pt@C=55pt{
   X_0^d \ar^{(\zg'(X_0))_{0}}[r] \ar@{=}[d] & X_{00}\\
   X_0^d \ar_{(\zg'(X_0))^d_0}[r] & X_{00}^d \ar_{\zg'(X_{00})}[u]
   }
   \qquad
   \xymatrix@R=25pt@C=50pt{
   X_0^d \ar^{(\zg'(X_0))_{\uh 0}}[r] \ar@{=}[d] & X_{0\uh 0}\\
   X_0^d \ar_{(\zg'(X_0))^d_{\uh 0}}[r] & X_{0\uh 0}^d \ar_{\zg'(X_{0\uh 0})}[u]
   }
   \end{equation}
   \smallskip
   \begin{equation}\label{eq-def-fta-1b}
   \xymatrix@R=25pt@C=55pt{
   X_{\uk 0}^d \ar^{(\zg'(X_{\uk 0}))_{0}}[r] \ar@{=}[d] & X_{\uk 00}\\
   X_{\uk 0}^d \ar_{(\zg'(X_{\uk 0}))^d_0}[r] & X_{\uk 00}^d \ar_{\zg'(X_{\uk 00})}[u]
   }
   \qquad
   \xymatrix@R=25pt@C=50pt{
   X_{\uk 0}^d \ar^{(\zg'(X_{\uk 0}))_{\uh 0}}[r] \ar@{=}[d] & X_{\uk 0\uh 0}\\
   X_{\uk 0}^d \ar_{(\zg'(X_{\uk 0}))^d_{\uh 0}}[r] & X_{\uk 0\uh 0}^d \ar_{\zg'(X_{\uk 0\uh 0})}[u]
   }
   \end{equation}
    \end{itemize}

       A morphism $F:X\rw Y$ in $\ftawg{n}$ is a morphism in $\catwg{n}$ such that the following diagram commutes
          \begin{equation}\label{eq-def-fta-2b}
          \xymatrix@C=50pt@R=30pt{
          X_{0} \ar^{F_{0}}[r] & Y_{0}\\
          X^d_{0} \ar_{F^d_{0}}[r] \ar^{\zg'(X_{0})}[u] & Y^d_{0} \ar_{\zg'(Y_{0})}[u]
          }
          \end{equation}
        and such that for all $\uk=(\seqc{k}{1}{s})\in\dop{s}$, $1\leq s \leq n-2$, the following diagram commutes
          \begin{equation}\label{eq-def-fta-2}
          \xymatrix@C=50pt@R=30pt{
          X_{\uk 0} \ar^{F_{\uk 0}}[r] & Y_{\uk 0}\\
          X^d_{\uk 0} \ar_{F^d_{\uk 0}}[r] \ar^{\zg'(X_{\uk 0})}[u] & Y^d_{\uk 0} \ar_{\zg'(Y_{\uk 0})}[u]
          }
          \end{equation}
where $\zg'(X_{0})$, $\zg'(Y_{0})$, $\zg'(X_{\uk 0})$, $\zg'(Y_{\uk 0})$ are sections to the corresponding discretization maps.
\end{definition}
\begin{remark}\label{rem-fta-2}
It is immediate from Definition \ref{def-fta-1} that
    \begin{equation*}
      \ftawg{n} \subset \funcat{}{\ftawg{n-1}}
    \end{equation*}
Also, from from Definition \ref{def-fta-1} we see that $\zg'(X_{0}):X_{0}^d \rw X_{0}$ is a map in $\ftawg{n-1}$ and $\zg'(X_{\uk 0}):X_{\uk 0}^d \rw X_{\uk 0}$ is a map in $\ftawg{n-s-1}$.

Further, we observe that if $X\in \ftawg{n}$ and $X\in\cathd{n}$, the discretization map $\zg: X\rw X^d$  is a map in $\ftawg{n}$. In fact, by naturality of discretization map, for each $\uk\in\dop{s}\quad (1\leq s \leq n-2)$ the following diagrams commute
\begin{equation*}
\xymatrix@R=30pt@C=50pt{
X_0 \ar^{\zg_0}[r] \ar_{\zd_0}[d] & (X^d)_0 \ar@{=}[d] \\
X_0^d \ar[r] & (X^d)_0^d
}
\qquad
\xymatrix@R=30pt@C=50pt{
X_{\uk 0} \ar^{\zg_{\uk 0}}[r] \ar_{\zd_{\uk 0}}[d] & (X^d)_{\uk 0} \ar@{=}[d] \\
X_{\uk 0}^d \ar_{\zd_{\uk 0}^d}[r] & (X^d)_{\uk 0}^d
}
\end{equation*}
where $\zd_0$ (resp- $\zd_{\uk 0}$) is the discretization map for $X_0$ (resp. $X_{\uk 0}$). Thus we have the following commuting diagrams
\begin{equation*}
\xymatrix@R=20pt@C=30pt{
X_0 \ar^{\zg_0}[rr] \ar^{\zd_0}[dr] & & (X^d)_0 \ar@{=}[dd] \\
 & X_0^d \ar[ur]&\\
X_0^d \ar^{\zd'_0}[uu]\ar_{\zg_0^d}[rr] \ar@{=}[ur] & & (X^d)_0^d
}
\qquad
\xymatrix@R=20pt@C=30pt{
X_{\uk 0} \ar^{\zg_{\uk 0}}[rr] \ar^{\zd^{\uk 0}}[dr] & & (X^d)_{\uk 0} \ar@{=}[dd] \\
 & (X_{\uk 0})^d \ar[ur] & \\
X_{\uk 0}^d \ar^{\zd'_{\uk 0}}[uu]\ar_{\zg_{\uk 0}^d}[rr]\ar@{=}[ur] & & (X^d)_{\uk 0}^d
}
\end{equation*}
By Definition \ref{def-fta-1}, the commutativity of the outer part of these diagrams mean that $\zg$ is a morphism in $\ftawg{n}$.
\end{remark}
\begin{lemma}\label{lem-fta-1}
    The functors $\p{n},\q{n}:\catwg{n}\rw \catwg{n-1}$ induce functors
    \begin{equation*}
        \p{n},\q{n}:\ftawg{n}\rw \ftawg{n-1}\;.
    \end{equation*}
\end{lemma}
\begin{proof}
Let $X\in\ftawg{n}$ and $\uk\rw\ur$ be a morphism in $\dop{s}$. By applying the functor $\p{n-s-1}$ to the commuting diagram \eqref{eq-def-fta-1} and using the fact that
\begin{equation*}
    \p{n-s-1}X_{\uk 0}=(\p{n}X)_{\uk 0},\qquad \p{n-s-1}X^d_{\uk 0}= X^d_{\uk 0}=(\p{n}X)^d_{\uk 0}
\end{equation*}
we obtain the commuting diagram
\begin{equation}\label{eq-fp}
\xymatrix@C=60pt@R=40pt{
(\p{n}X)_{\uk 0} \ar^{\p{n-s-1}f}[r] & (\p{n}X)_{\ur 0} \\
(\p{n}X)^d_{\uk 0} \ar_{\p{n-s-1}f^d}[r] \ar[u] & (\p{n}X)^d_{\ur 0} \ar[u]
}
\end{equation}
By applying $\p{n-2}$ to the diagram on the left of \eqref{eq-def-fta-1a} and $\p{n-t-2}$ to the diagram on the right of \eqref{eq-def-fta-1a} and using the fact that
\begin{equation*}
  \p{n-2}X_{00}=(\p{n}X)_{00} \qquad \p{n-t-2}X_{0\uh 0}=(\p{n}X)_{0\uh 0}
\end{equation*}
we obtain commuting diagrams
\begin{equation}\label{eq-def-fta-4}
\xymatrix@R=25pt@C=35pt{
(\p{n}X)^d_0 \ar[r] \ar@{=}[d] & (\p{n}X)_{00} \\
(\p{n}X)^d_0  \ar[r] & (\p{n}X)^d_{00} \ar[u]
}
\qquad
\xymatrix@R=25pt@C=35pt{
(\p{n}X)^d_0 \ar[r] \ar@{=}[d] & (\p{n}X)_{0\uh 0} \\
(\p{n}X)^d_0  \ar[r] & (\p{n}X)^d_{0\uh 0} \ar[u]
}
\end{equation}
Similarly, applying $\p{n-s-2}$ to the diagram on the left of \eqref{eq-def-fta-1b} and $\p{n-s-t-2}$ to the diagram on the right of \eqref{eq-def-fta-1b} and using the fact that
\begin{equation*}
  \p{n-s-2}X_{\uk 00}=(\p{n}X)_{\uk 00} \qquad \p{n-s-t-2}X_{\uk 0\uh 0}=(\p{n}X)_{\uk 0\uh 0}
\end{equation*}
we obtain commuting diagrams
\begin{equation}\label{eq-def-fta-5}
\xymatrix@R=25pt@C=35pt{
(\p{n}X)^d_{\uk 0} \ar[r] \ar@{=}[d] & (\p{n}X)_{\uk 00} \\
(\p{n}X)^d_{\uk 0}  \ar[r] & (\p{n}X)^d_{\uk 00} \ar[u]
}
\qquad
\xymatrix@R=25pt@C=35pt{
(\p{n}X)^d_{\uk 0} \ar[r] \ar@{=}[d] & (\p{n}X)_{\uk 0\uh 0} \\
(\p{n}X)^d_{\uk 0}  \ar[r] & (\p{n}X)^d_{\uk 0\uh 0} \ar[u]
}
\end{equation}
Together with \eqref{eq-fp}, \eqref{eq-def-fta-4} and \eqref{eq-def-fta-5} mean by definition that $\p{n}X\in\ftawg{n}$.

Given $F:X\rw Y$ in $\ftawg{n}$,
by applying $\p{n-1}$ to the commuting diagram \eqref{eq-def-fta-2b} we obtain the commuting diagram
\begin{equation*}
\xymatrix@C=60pt@R=40pt{
(\p{n}X)_{0} \ar^{(\p{n}F)_{0}}[r] & (\p{n}Y)_{0} \\
(\p{n}X)^d_{0} \ar_{(\p{n}F)^d_{0}}[r] \ar[u] & (\p{n}Y)^d_{0} \ar[u]
}.
\end{equation*}
By applying $\p{n-s-1}$ to the commuting diagram \eqref{eq-def-fta-2} we obtain the commuting diagram
\begin{equation*}
\xymatrix@C=60pt@R=40pt{
(\p{n}X)_{\uk 0} \ar^{(\p{n}F)_{\uk 0}}[r] & (\p{n}Y)_{\uk 0} \\
(\p{n}X)^d_{\uk 0} \ar_{(\p{n}F)^d_{\uk 0}}[r] \ar[u] & (\p{n}Y)^d_{\uk 0} \ar[u]
}
\end{equation*}
By definition this means that $\p{n}F\in\ftawg{n}$.

 The proof for $\q{n}$ is analogous.
\end{proof}

\subsection{The idea of the functor $\pmb{G_n}$}\label{subs-idea-gn} As shown in Proposition \ref{pro-gen-const-1}, the functor $F_n:\catwg{n}\rw \catwg{n}$ replaces $X\in\catwg{n}$ with an $n$-equivalent $F_n X$ in which $(F_n X)_0$ admits a functorial choice of section to the discretization map.
By definition of $\ftawg{2}$, the functor $F_2$ is in fact a functor $F_2: \catwg{2}\rw \ftawg{2}$ and thus we define $G_2=F_2$.

When $n>2$, recall that by definition objects of $\ftawg{n}$ are such that the homotopically discrete substructures in the multinerve have functorial sections to the discretization maps, as well as other functoriality properties. The idea of the functor $G_n$ is to inductively apply $F_s$ to every sub-simplicial dimension $s$. That is, we define inductively
  \begin{equation*}
    G_n:\catwg{n}\rw \funcat{}{\catwg{n-1}}
  \end{equation*}
  by $G_2=F_2$; when $n>2$, given $G_{n-1}$ we let
     \begin{equation*}
       G_n=\ovl{G}_{n-1}\circ F_n.
     \end{equation*}
Showing that in fact $G_n$ lands in $\ftawg{n}$ involves several steps, as developed in the proof of Theorem \ref{pro-fta-1}.
First we show that $G_n X \in \catwg{n}$ by proving that it satisfies the hypotheses of Lemma \ref{lem-x-in-tawg-x-in-catwg}.

For this we first check that $G_n X \in \tawg{n}$. The main ingredients here are the inductive hypotheses on $G_{n-1}$ that it preserves $(n-1)$-equivalences and pullbacks over discrete objects and it is identity on discrete objects.
This easily implies that there is an $(n-1)$-equivalence
\begin{equation*}
    (G_n X)_0=G_{n-1}(F_n X)_0 \rw G_{n-1}(F_n X)^d_0 = (F_n X)^d_0\;.
\end{equation*}
so that $(G_n X)_0$ is homotopically discrete. It also implies that the induced Segal maps of $G_n X$ are $(n-1)$-equivalences, thus in conclusion $X\in \tawg{n}$.

Hypothesis a) of Lemma \ref{lem-x-in-tawg-x-in-catwg} is immediate by induction since $(G_n X)_k=G_{n-1}(F_n X)_k \in \ftawg{n-1}$ so in particular $(G_n X)_k\in\catwg{n-1}$.

As for hypotheses b) and c) in Lemma \ref{lem-x-in-tawg-x-in-catwg}, these are shown to hold for $G_n X$ by assuming, inductively, that $G_{n-1}$ satisfies the same properties of $F_n$ with respect to pullbacks established in Corollary \ref{cor-fn-pull}. Recall that this states that given a diagram in $\catwg{n}$
\begin{equation*}
X\rw Z \lw Y
\end{equation*}
 such that $X\tiund{Z}Y \in \catwg{n}$ and such that this pullback is preserved by $\p{r,n}$ and $\q{r,n}$ (for each $1\leq r\leq n$), then this pullback is also preserved by $F_n$, $\p{r,n}F_n$ and $\q{r,n}F_n$.

 We require the same property to hold, inductively, for $G_{n-1}$. We then apply $G_{n-1}$ to the diagram in $\catwg{n-1}$
 \begin{equation}\label{eq-pro-fta-1-1}
 (F_n X)_1 \xrw{\;\pt_0 \;} (F_n X)_0 \xlw{\;\pt_{0}\;}(F_n X)_1
\end{equation}
Since $F_n X\in \catwg{n}$, this diagram is such that
\begin{equation*}
  (F_n X)_2 \cong (F_n X)_1 \tiund{(F_n X)_0} (F_n X)_1 \in\catwg{n-1}\;.
\end{equation*}
and it commutes with $\p{r,n}$ and $\q{r,n}$. Thus from inductive hypothesis,
 \begin{equation*}
\begin{split}
    & (G_n X)_2 = G_{n-1}(F_n X)_2= G_{n-1}((F_n X)_1 \tiund{(F_n X)_0} (F_n X)_1)\cong\\
    & \cong G_{n-1}(F_n X)_1 \tiund{G_{n-1}(F_n X)_0} G_{n-1}(F_n X)_1\cong (G_n X)_1 \tiund{(G_n X)_0} (G_n X)_1\;.
\end{split}
\end{equation*}
Similarly one shows that, for each $s\geq 2$
\begin{equation*}
  (G_n X)_s \cong \pro{(G_n X)_1}{(G_n X)_0}{s} \;.
\end{equation*}
This proves that $G_n X$ satisfies hypothesis b) of Lemma \ref{lem-x-in-tawg-x-in-catwg}. Hypothesis c) of Lemma \ref{lem-x-in-tawg-x-in-catwg} is checked similarly and we conclude that $G_n X \in \catwg{n}$.

For the inductive step we then need to show that $G_n$ satisfies the commutation properties with respect to pullbacks at step $n$. This is done by using Corollary \ref{cor-fn-pull}, as detailed in the proof of Theorem \ref{pro-fta-1} e).

In the proof of Theorem \ref{pro-fta-1} we also check the remaining functoriality conditions and show that $G_n: \catwg{n} \rw \ftawg{n}$.

\subsection{The functor $\pmb{G_n}$: the formal proof}\label{gn-formal}
\begin{definition}\label{def-Gn-new}
  Define inductively
  \begin{equation*}
    G_n:\catwg{n}\rw \funcat{}{\catwg{n-1}}
  \end{equation*}
  by $G_2=F_2$; given $G_{n-1}$ let
     \begin{equation*}
       G_n=\ovl{G}_{n-1}\circ F_n
     \end{equation*}
     where $F_n$ is an in Proposition \ref{pro-gen-const-1}.
\end{definition}

\begin{lemma}\label{lem-zn}
  For $(\seqc{k}{1}{s})\in\dop{s}$, $1\leq s\leq n-2$, $Y\in\catwg{n}$ define inductively $Z_{k_1 k_2\ldots k_s}Y\in\catwg{n-s}$ by
\begin{equation*}
  Z_{k_1}Y=(F_n Y)_{k_1}
\end{equation*}
and for $1< i\leq s-1$, given $Z_{k_1 k_2\ldots k_{i-1}}Y$, let
\begin{equation}\label{eq-pro-fta-8}
  Z_{k_1 k_2\ldots k_i}Y=(F_{n-i+1}Z_{k_1 k_2\ldots k_{i-1}}Y)_{k_i}\;.
\end{equation}
Then for each $2\leq s\leq n-2$
\begin{equation*}
  (G_n Y)_{\seq{k}{1}{s}}=G_{n-s}(F_{n-s+1}Z_{\seq{k}{1}{s-1}}Y)_{k_s}\;.
\end{equation*}
where $G_n$ is as in Definition \ref{def-Gn-new}.
\end{lemma}

\begin{proof}
We show this by induction on $s$. By Definition of $G_n$,
\begin{equation*}
\begin{split}
    & (G_n Y)_{k_1}=G_{n-1}(F_n Y)_{k_1} \\
    & (G_n Y)_{k_1 k_2}=G_{n-2}(F_{n-1}(F_n Y)_{k_1})_{k_2}=G_{n-2}(F_{n-1}Z_{k_1}Y)_{k_2}\;.
\end{split}
\end{equation*}
Suppose, inductively, that the lemma holds for $s-1$. Then, by inductive hypothesis and by \eqref{eq-pro-fta-8}
\begin{equation*}
\begin{split}
    & (G_n Y)_{\seq{k}{1}{s}}= \{(G_n Y)_{\seq{k}{1}{s-1}} \}_{k_s}=\\
    & \{ G_{n-s+1}(F_{n-s+2}Z_{\seq{k}{1}{s-2}}Y)_{k_{s-1}} \}_{k_s}= \\
    & = (G_{n-s+1} Z_{\seq{k}{1}{s-1}}Y)_{k_s}=\\
    &  =G_{n-s}(F_{n-s+1}Z_{\seq{k}{1}{s-1}}Y)_{k_s}
\end{split}
\end{equation*}
proving the lemma.

\end{proof}

\begin{theorem}\label{pro-fta-1}
Let $G_n$ be as in Definition \ref{def-Gn-new}. Then

    \begin{itemize}
      \item [a)] $G_n:\catwg{n}\rw \ftawg{n}$.\mk

      \item [b)] There is an $n$-equivalence
       \begin{equation*}
       g_n(X):G_n X\rw X
       \end{equation*}
        natural in $X\in\catwg{n}$ and such that $(g_n(X))_0^d:(G_n X)_0^d\rw X_0^d$ is surjective.\mk

      \item [c)] $G_n$ preserves $n$-equivalences.\mk

      \item [d)] $G_n$ is identity on discrete objects and preserves pullbacks over discrete objects.\mk

      \item [e)] Let $X\rw Z \lw Y$ be a diagram in $\catwg{n}$ such that $X\tiund{Z}Y\in\catwg{n}$ and such that, for each $1\leq r\leq n$
    \begin{equation*}
    \p{r,n} (X\tiund{Z}Y)= \p{r,n} X \tiund{\p{r,n} Z} \p{r,n} Y\;,
    \end{equation*}
    \begin{equation*}
      \q{r,n} (X\tiund{Z}Y)= \q{r,n} X \tiund{\q{r,n} Z} \q{r,n} Y\;.
    \end{equation*}
Then
  \begin{equation*}
  G_n(X\tiund{Z}Y)=G_n X\tiund{G_n Z}G_n Y\;,
\end{equation*}
\begin{equation*}
  g_n(X\tiund{Z}Y)=g_n (X)\tiund{g_n (Z)}g_n (Y)\;,
\end{equation*}
  \begin{equation*}
  \p{r,n} G_n(X\tiund{Z}Y)\cong \p{r,n} G_n X \tiund{\p{r,n}G_n Z} \p{r,n}G_n Y\;,
  \end{equation*}
  \begin{equation*}
  \q{r,n} G_n(X\tiund{Z}Y)\cong \q{r,n} G_n X \tiund{\q{r,n}G_n Z} \q{r,n}G_n Y\;.
  \end{equation*}
  \mk

  \item [f)] If $Y\in\cathd{n}$ and $\gamma': Y^d\rw Y$ is a section to the discretization map, then $G_n\zg'$ is a morphism in $\ftawg{n}$.

    \end{itemize}
\end{theorem}
\begin{proof}
By induction on $n$. For $n=2$, by Proposition \ref{pro-gen-const-1} and Corollary \ref{cor-fn-pull} the functor $G_2=F_2:\catwg{2}\rw \catwg{2}$ is in fact a functor $\catwg{2}\rw \ftawg{2}$ satisfying a) - f), where to show b) we use the fact that $(G_2 X)_0^d=X_{00} \rw X_0^d$ is surjective.

Suppose we defined $G_{n-1}$ satisfying the above properties and let $X\in\catwg{n}$.

 a) First we show that $G_n X \in \catwg{n}$  using the criterion given in Lemma \ref{lem-x-in-tawg-x-in-catwg}.

 We first check that $G_n X\in\tawg{n}$. By construction we have
\begin{equation*}
    (G_n X)_0=G_{n-1}(F_n X)_0\;.
\end{equation*}
Since $(F_n X)_0\in\cathd{n-1}$ (as $F_n X\in\catwg{n}$), there is a $(n-1)$-equivalence $(F_n X)\rw(F_n X)^d_0$. Thus by inductive hypothesis c) and d) this induces an $(n-1)$-equivalence
\begin{equation*}
    (G_n X)_0=G_{n-1}(F_n X)_0 \rw G_{n-1}(F_n X)^d_0 = (F_n X)^d_0\;.
\end{equation*}
Therefore $(G_n X)_0\in\cathd{n-1}$ and
\begin{equation*}
    (G_n X)^d_0=(F_n X)^d_0\;.
\end{equation*}
For each $k>0$ by inductive hypothesis we also have
\begin{equation*}
    (G_n X)_k = G_{n-1}(F_n X)_k\in\catwg{n-1}\;.
\end{equation*}
To show that $G_n X\in\tawg{n}$ it remains to prove that the induced Segal maps are $(n-1)$-equivalences. Since $F_n X\in\catwg{n}$ there are $(n-1)$-equivalences
\begin{equation*}
    (F_n X)_2 \rw \tens{(F_n X)_1}{(F_n X)^d_0}\;.
\end{equation*}
Using the induction hypotheses c) and d) this induces an $(n-1)$-equivalence
\begin{equation*}
\begin{split}
  (G_n X)_2 & =G_{n-1}(F_n X)_2 \rw G_{n-1} \{\tens{(F_n X)_1}{(F_n X)^d_0}\}\cong \\
    & \cong \tens{(G_n X)_1}{(G_n X)^d_0}\;.
\end{split}
\end{equation*}
Similarly one shows that all other induced Segal maps for $G_n X$ are $(n-1)$-equivalences. We conclude that $G_n X\in\tawg{n}$.

We now check the rest of the hypotheses in Lemma \ref{lem-x-in-tawg-x-in-catwg}. Hypothesis a) holds since, as seen above, for each $k\geq 0$, $(G_n X)_k\in\catwg{n-1}$.

Note that the diagram in $\catwg{n-1}$
\begin{equation}\label{eq-pro-fta-1-1}
 (F_n X)_1 \xrw{\;\pt_0 \;} (F_n X)_0 \xlw{\;\pt_{0}\;}(F_n X)_1
\end{equation}
satisfies the inductive hypothesis e). In fact, since $F_n X\in\catwg{n}$
\begin{equation*}
  (F_n X)_2 \cong (F_n X)_1 \tiund{(F_n X)_0} (F_n X)_1 \in\catwg{n-1}\;.
\end{equation*}
The rest of the inductive hypothesis e) for the diagram \eqref{eq-pro-fta-1-1} follows from the fact that $F_n X\in\catwg{n}$ using Lemma \ref{lem-prop-pn}  and Lemma \ref{lem-prop-qn}. Thus by inductive hypothesis e) we obtain
\begin{equation*}
\begin{split}
    & (G_n X)_2 = G_{n-1}(F_n X)_2= G_{n-1}((F_n X)_1 \tiund{(F_n X)_0} (F_n X)_1)\cong\\
    & \cong G_{n-1}(F_n X)_1 \tiund{G_{n-1}(F_n X)_0} G_{n-1}(F_n X)_1\cong (G_n X)_1 \tiund{(G_n X)_0} (G_n X)_1\;.
\end{split}
\end{equation*}
Similarly one shows that, for each $s\geq 2$
\begin{equation*}
  (G_n X)_s \cong \pro{(G_n X)_1}{(G_n X)_0}{s} \;.
\end{equation*}
This proves that $G_n X$ satisfies hypothesis b) of Lemma \ref{lem-x-in-tawg-x-in-catwg}.
 By inductive hypothesis e) applied to the diagram \eqref{eq-pro-fta-1-1} we obtain, for each $1\leq r< n$
\begin{equation*}
\begin{split}
    & \p{r,n-1}(G_n X)_2= \p{r,n-1}G_{n-1}(F_n X)_2 \cong\\
    & \cong \p{r,n-1}G_{n-1} \{ (F_n X)_1\tiund{(F_n X)_0}(F_n X)_1 \} \cong\\
    & \cong\p{r,n-1}G_{n-1}(F_n X)_1\tiund{\p{r,n-1}G_{n-1}(F_n X)_0}\p{r,n-1}G_{n-1}(F_n X)_1=\\
    & =\p{r,n-1}(G_n X)_1\tiund{\p{r,n-1}(G_n X)_0}\p{r,n-1}(G_n X)_1\;.
\end{split}
\end{equation*}
Similarly one shows that, for each $s\geq 2$
\begin{equation*}
\begin{split}
    & \p{r,n-1}(G_n X)_s \cong \\
    & \cong \pro{\p{r,n-1}(G_n X)_1}{\p{r,n-1}(G_n X)_0}{s}\;,
\end{split}
\end{equation*}
which is hypothesis c) of Lemma \ref{lem-x-in-tawg-x-in-catwg} for $G_n X$. Thus all hypotheses of Lemma \ref{lem-x-in-tawg-x-in-catwg} hold and we conclude that $G_n X\in\catwg{n}$.

 To show that $G_n X\in\ftawg{n}$ we need to prove the commutativity of the diagrams in i) and ii) of Definition \ref{def-fta-1}. Let $\uk=(\seqc{k}{1}{s})$, $\ur=(\seqc{r}{1}{s})$ in $\dop{s}$, denote $\uk'=(\seqc{k}{2}{s})$, $\ur'=(\seqc{r}{2}{s})$ and suppose we have a morphism $\uk\rw\ur$ in $\dop{s}$. By factoring this as
\begin{equation*}
    \xymatrix@R=25pt{
    \uk=(k_1,\uk') \ar[rr] \ar[dr] && \ur=(r_1,\ur')\\
    & (r_1,\uk')\ar[ur]
    }
\end{equation*}
we obtain a factorization
\scriptsize
\begin{equation}\label{eq-pro-fta-1}
\xymatrix@C=5pt{
(G_n X)_{\uk 0}=\{G_{n-1}(F_n X)_{k_1}\}_{\uk' 0} \ar[rr]\ar[dr] && \{G_{n-1}(F_n X)_{r_1}\}_{\ur' 0}=(G_n X)_{\ur 0}\\
& \{G_{n-1}(F_n X)_{r_1}\}_{\uk' 0}\ar[ur]
}
\end{equation}
\normalsize
Consider the morphism $(F_n X)_{k_1}\rw (F_n X)_{r_1}$ in $\catwg{n-1}$. Since, by induction hypothesis a), $G_{n-1}:\catwg{n-1}\rw \ftawg{n-1}$ there is a commuting diagram
\begin{equation}\label{eq-pro-fta-2}
\xymatrix@C=50pt@R=40pt{
(G_n X)_{\uk 0}=\{G_{n-1}(F_n X)_{k_1}\}_{\uk' 0} \ar[r] & \{G_{n-1}(F_n X)_{r_1}\}_{\uk' 0}\\
{(G_n X)}^d_{\uk 0}=\{G_{n-1}(F_n X)_{k_1}\}^d_{\uk' 0} \ar[r]\ar[u] & \{G_{n-1}(F_n X)_{r_1}\}^d_{\uk' 0}\ar[u]
}
\end{equation}
Since, by induction hypothesis, $G_{n-1}(F_n X)_{r_1}\in\ftawg{n-1}$ we also have a commuting diagram
\begin{equation}\label{eq-pro-fta-3}
\xymatrix@C=50pt@R=40pt{
\{G_{n-1}(F_n X)_{r_1}\}_{\uk' 0} \ar[r] & \{G_{n-1}(F_n X)_{r_1}\}_{\ur' 0}=(G_n X)_{\ur 0}\\
\{G_{n-1}(F_n X)_{r_1}\}^d_{\uk' 0} \ar[r]\ar[u] & \{G_{n-1}(F_n X)_{r_1}\}^d_{\ur' 0}={(G_n X)}^d_{\ur 0}\ar[u]
}
\end{equation}
Combining \eqref{eq-pro-fta-1}, \eqref{eq-pro-fta-2}, \eqref{eq-pro-fta-3} we obtain a commuting diagram
\begin{equation}\label{eq-new}
\xymatrix@C=40pt@R=35pt{
(G_n X)_{\uk 0} \ar[r] & (G_n X)_{\ur 0}\\
(G_n X)^d_{\uk 0} \ar[r]\ar[u] & (G_n X)^d_{\ur 0}\ar[u]
}
\end{equation}
Thus $G_n X$ satisfies condition \eqref{eq-def-fta-1} in the definition of $\ftawg{n}$.

Applying inductive hypothesis f) to the map in $\cathd{n-1}$ $(F_n X)_{0}^d\rw(F_n X)_{0} $, we see that the map
\begin{equation*}
  (G_n X)_{0}^d=(F_n X)_{0}^d=G_{n-1}(F_n X)_{0}^d\rw G_{n-1}(F_n X)_{0}=(G_n X)_0.
\end{equation*}
is a morphism in $\ftawg{n-1}$.
This means that the two diagrams \eqref{eq-def-fta-1a} in ii) of Definition \ref{def-fta-1} commute for $G_n X$.

Applying the inductive hypothesis f) to the map in $\cathd{n-1+1}$
\begin{equation*}
  (F_{n-s+1}Z_{\uk}(X))_0^d\rw(F_{n-s+1}Z_{\uk}(X))_0
\end{equation*}
and using the fact that, by Lemma \ref{lem-zn} it is
\begin{equation*}
  (G_n X)_{\uk 0}=G_{n-s}(F_{n-s+1}Z_{\uk}(X))_0.
\end{equation*}
we deduce that the map
\begin{equation*}
  (G_n X)_{\uk 0}^d\rw (G_n X)_{\uk 0}.
\end{equation*}
is a morphism in $\ftawg{n-s-1}$. This means that the two diagrams \ref{eq-def-fta-1b} in ii) of Definition \ref{def-fta-1} commute for $G_n X$.
Together with \eqref{eq-new} we conclude that $G_n X\in \ftawg{n}$.

Finally, let $F:X\rw Y$ be a morphism in $\catwg{n}$. Then
 \begin{equation*}
  (F_n F)_{0}:(F_n X)_{0}\rw (F_n Y)_{0}
 \end{equation*}
 is a morphism in $\catwg{n-1}$. Thus by induction hypothesis it induces a morphism in $\ftawg{n-1}$

 \begin{equation*}
    G_{n-1}(F_n X)_{0}\rw G_{n-1}(F_n Y)_{0}
\end{equation*}
such that the following diagram commutes:
\begin{equation*}
\xymatrix@R=35pt{
(G_n X)_{0}=G_{n-1}(F_n X)_{0} \ar[r] & G_{n-1}(F_n Y)_{0}=(G_n Y)_{0}\\
(G_n X)^d_{0}\ar[r]\ar[u] & (G_n Y)^d_{0}\ar[u]
}
\end{equation*}

 It also induces a morphism
\begin{equation*}
    G_{n-1}(F_n X)_{k_1}\rw G_{n-1}(F_n Y)_{k_1}
\end{equation*}
such that the following diagram commutes:
\begin{equation*}
\xymatrix@R=35pt{
(G_n X)_{\uk 0}=\{G_{n-1}(F_n X)_{k_1}\}_{\uk' 0} \ar[r] & \{G_{n-1}(F_n Y)_{k_1}\}_{\uk' 0}=(G_n Y)_{\uk 0}\\
(G_n X)^d_{\uk 0}\ar[r]\ar[u] & (G_n Y)^d_{\uk 0}\ar[u]
}
\end{equation*}
This shows that $G_n F$ is a morphism in $\ftawg{n}$.

 In conclusion
\begin{equation*}
    G_n:\catwg{n}\rw \ftawg{n}\;.
\end{equation*}
The fact that $(G_n X)_k=G_{n-1}(F_n X)_k\in\catwg{n-1}$ follows by induction.
\bk

b) The morphism $g_n(X): G_n X\rw X$ is given levelwise  by
 \begin{equation*}
 (g_n(X))_k=g_{n-1}(X_k):(G_n X)_k= G_{n-1} X_k\rw X_k.
 \end{equation*}
  By inductive hypothesis this is an $(n-1)$-equivalence for each $k$, hence $g_n(X)$ is a $n$-equivalence by Lemma \ref{lem-flevel-fneq}.

We now show that $(g_n(X))_0^d$ is surjective. As in the proof of a), $(G_n X)^d_0=(F_n X)^d_0$ and by Proposition \ref{pro-gen-const-1}, $(F_n X)^d_0=(V_{n-1}X_0)^d$. So we need to show that
\begin{equation*}
  (V_{n-1}X_0)^d\rw X_0^d
\end{equation*}
is surjective. By Proposition \ref{pro-gen-const-2} applied to $X_0\in\cathd{n-1}$, the functor
\begin{equation*}
  \p{2,n-1}V_{n-1} X_0\rw\p{2,n-1} X_0
\end{equation*}
is surjective on objects. Hence
\begin{equation*}
  (V_{n-1}X_0)^d=p\p{2,n-1}V_{n-1}X_0 \rw p\p{2,n-1}X_0=X_0^d
\end{equation*}
is surjective, as required.

\bk

c) Let $F:X\rw Y$ be an $n$-equivalence in $\catwg{n}$. There is a commuting diagram
\begin{equation*}
\xymatrix@R=30pt@C=40pt{
G_n X \ar^{G_n F}[r] \ar_{g_n(X)}[d] & G_n Y \ar^{g_n(Y)}[d]\\
X \ar_{F}[r] & Y
}
\end{equation*}
in which, by b), the vertical maps and the bottom horizontal map are $n$-equivalences. By Proposition \ref{pro-n-equiv} c) it follows that $G_n F$ is also an $n$-equivalence.
\bk

d) This follows immediately by the analogous properties of $F_n$ and by the inductive hypothesis.\bk

e) By hypothesis the diagram in $\catwg{n}\;\; X\rw Z\lw Y$ satisfies the hypotheses of Corollary \ref{cor-fn-pull} so that
\begin{equation}\label{eq-pro-fta-4}
F_n(X\tiund{Z}Y)\cong F_n X \tiund{F_n Z} F_n Y\;,
\end{equation}
\begin{equation}\label{eq-pro-fta-5}
\p{r,n} F_n(X\tiund{Z}Y)\cong \p{r,n} F_n X \tiund{\p{r,n}F_n Z} \p{r,n}F_n Y\;,
\end{equation}
\begin{equation}\label{eq-pro-fta-6}
\q{r,n} F_n(X\tiund{Z}Y)\cong \q{r,n} F_n X \tiund{\q{r,n}F_n Z} \q{r,n}F_n Y\;.
\end{equation}
We claim that the diagram in $\catwg{n-1}$
\begin{equation}\label{eq-pro-fta-7}
  (F_n X)_k \rw (F_n Z)_k \lw (F_n Y)_k
\end{equation}
satisfies the inductive hypothesis e), In fact
\begin{equation*}
  (F_n X)_k \tiund{(F_n Z)_k} (F_n Y)_k=(F_n(X\tiund{Z}Y))_k\in\catwg{n-1}
\end{equation*}
since $F_n(X\tiund{Z}Y)\in\catwg{n}$. Also, taking the $k^{th}$-component in \eqref{eq-pro-fta-5} we obtain
\begin{equation*}
\begin{split}
    & \p{r-1,n-1}\{ (F_n X)_k \tiund{(F_n Z)_k} (F_n Y)_k\} = \p{r-1,n-1} (F_n(X\tiund{Z}Y))_k =\\
    & \cong \p{r-1,n-1} (F_n X)_k \tiund{\p{r-1,n-1}(F_n Z)_k} \p{r-1,n-1}(F_n Y)_k
\end{split}
\end{equation*}
and similarly for $\q{r-1,n-1}$.

Thus by inductive hypothesis e) applied to \eqref{eq-pro-fta-7} we obtain
\begin{equation*}
\begin{split}
    & (G_n(X\tiund{Z}Y))_k = G_{n-1}(F_n(X\tiund{Z}Y))_k= \\
    & =G_{n-1} ((F_n X)_k \tiund{(F_n Z)_k} (F_n Y)_k)=\\
    & = G_{n-1} (F_n X)_k \tiund{G_{n-1}(F_n Z)_k} G_{n-1}(F_n Y)_k =\\
    & = (G_n X)_k\tiund{(G_n Z)_k}(G_n Y)_k\;.
\end{split}
\end{equation*}
and

\begin{equation*}
\begin{split}
    & (g_n(X\tiund{Z}Y))_k = g_{n-1}(F_n(X\tiund{Z}Y))_k= \\
    & =g_{n-1} ((F_n X)_k \tiund{(F_n Z)_k} (F_n Y)_k)=\\
    & = g_{n-1} ((F_n X)_k) \tiund{g_{n-1}((F_n Z)_k)} g_{n-1}((F_n Y)_k) =\\
    & = (g_n X)_k\tiund{(g_n Z)_k}(g_n Y)_k\;.
\end{split}
\end{equation*}

Since this holds for each $k\geq 0$ it follows that
\begin{equation*}
  G_n(X\tiund{Z}Y)=G_n X\tiund{G_n Z}G_n Y\;.
\end{equation*}
and
\begin{equation*}
  g_n(X\tiund{Z}Y)=g_n (X)\tiund{g_n (Z)}g_n (Y)\;.
\end{equation*}
\bk
By inductive hypothesis e) applied to \eqref{eq-pro-fta-7} we also obtain
\begin{equation*}
\begin{split}
    & \{\p{r,n} G_n(X\tiund{Z}Y)\}_k = \p{r-1,n-1} G_{n-1}(F_n(X\tiund{Z}Y))_k= \\
    & =\p{r-1,n-1} G_{n-1} ((F_n X)_k \tiund{(F_n Z)_k} (F_n Y)_k)=\\
    & = \p{r-1,n-1} G_{n-1} (F_n X)_k \tiund{ \p{r-1,n-1} G_{n-1}(F_n Z)_k}  \p{r-1,n-1} G_{n-1}(F_n Y)_k =\\
    & = (\p{r,n}G_n X)_k\tiund{(\p{r,n}G_n Z)_k}(\p{r,n}G_n Y)_k\;.
\end{split}
\end{equation*}
Since this holds for each $k\geq 0$ it follows that
\begin{equation*}
  \p{r,n} G_n(X\tiund{Z}Y)=\p{r,n}G_n X\tiund{\p{r,n}G_n Z}\p{r,n}G_n Y\;.
\end{equation*}
Similarly, taking the $k^{th}$ component in \eqref{eq-pro-fta-6} and using inductive hypothesis e) on \eqref{eq-pro-fta-7} one obtain
\begin{equation*}
  \q{r,n} G_n(X\tiund{Z}Y)=\q{r,n}G_n X\tiund{\q{r,n}G_n Z}\q{r,n}G_n Y\;.
\end{equation*}
This concludes the proof of e) at step $n$.

\bk

f) By Proposition \ref{pro-gen-const-1} the morphism $Y^d \rw Y$ induces a commuting diagram
\begin{equation*}
\xymatrix@R=30pt@C=50pt{
(F_n Y^d)_0=Y^d \ar[r] \ar@{=}[d] & (F_n Y)_0 \\
(F_n Y^d)_0=Y^d \ar[r] & (F_n Y)^d_0 \ar[u]
}
\end{equation*}
Thus applying $G_{n-1}$ to this diagram and using the fact (from above) that
\begin{equation*}
  (G_n Y)^d_0 =(F_n Y)^d_0=G_{n-1}(F_n Y)^d_0
\end{equation*}
we obtain
\begin{equation}\label{eq-pro-fta-10}
\xymatrix@R=30pt@C=40pt{
(G_n Y^d)_0 = G_{n-1}(F_n Y^d)_0 \ar[r] \ar@{=}[d] & G_{n-1}(F_n Y)_0=(G_n Y)_0 \\
(G_n Y^d)_0^d = G_{n-1}(F_n Y^d)_0 \ar[r] & G_{n-1}(F_n Y)^d_0=(G_n Y)^d_0 \ar[u]
}
\end{equation}
For each $\uk=(\seqc{k}{1}{s-1})\in\dop{s-1}$, by Lemma \ref{lem-zn}
\begin{equation*}
  (G_n Y)_{\uk 0}=G_{n-s}(F_{n-s+1}Z_{\uk}Y)_0
\end{equation*}
where we abbreviated $Z_{\seq{k}{1}{s-1}}Y=Z_{\uk}Y$. The map $Y^d\rw Y$ induces the map
\begin{equation*}
Z_{\uk}Y^d \rw Z_{\uk}Y
\end{equation*}
and by Proposition \ref{pro-gen-const-1} this induces the commuting diagram
\begin{equation*}
\xymatrix@R=30pt@C=50pt{
(F_{n-s+1}Z_{\uk}Y^d)_0 \ar[r] \ar@{=}[d] & (F_{n-s+1}Z_{\uk}Y)_0\\
(F_{n-s+1}Z_{\uk}Y^d)_0 \ar[r] & (F_{n-s+1}Z_{\uk}Y)_0^d \ar[u]
}
\end{equation*}
Applying $G_{n-s}$ to this diagram and recalling that by Lemma \ref{lem-zn} and the above
\begin{equation*}
\begin{split}
    & (G_n Y)_{\uk 0}=G_{n-s}(F_{n-s+1}Z_{\uk}Y)_0 \\
    & (G_n Y)^d_{\uk 0}=(F_{n-s+1}Z_{\uk}Y)^d_0
\end{split}
\end{equation*}
we obtain the commuting diagram
\begin{equation}\label{eq-pro-fta-11}
\xymatrix@R=30pt@C=50pt{
(G_n Y^d)_{\uk 0} \ar[r] \ar@{=}[d] & (G_n Y)_{\uk 0} \\
(G_n Y^d)^d_{\uk 0} \ar[r] & (G_n Y)^d_{\uk 0}\ar[u]
}
\end{equation}
By definition, \eqref{eq-pro-fta-10} and \eqref{eq-pro-fta-11} mean that $G_n\zg'$ is a morphism in $\ftawg{n}$.

\end{proof}

%%%%%%%%%%%%%%%%%%%%%%%%%%%%%%%%%%%%%%%%%%%%%%%%%%%%%%%%%%%%%%%%%%%%%%%%%%%%%%
\chapter{Weakly globular $\pmb{n}$-fold categories as a model of weak $\pmb{n}$-categories}\label{chap9}

In this chapter we prove that the category $\catwg{n}$ of weakly globular \nfol categories constitutes a model of weak $n$-categories. We show this by proving that there is an equivalence of categories between the localizations of $\ta{n}$ and of $\catwg{n}$ with respect to the $n$-equivalences. This shows a type of equivalence (up to higher categorical equivalence) between $\catwg{n}$ and $\ta{n}$.

 We also show that the category $\catwg{n}$ satisfies the homotopy hypothesis. As explained in Part \ref{part-1}, the latter is one of the main desiderata for a model of weak $n$-categories, while the comparison with the Tamsamani model is a contribution to the still largely open problem of comparing between different models of higher categories.

 The homotopy hypothesis is shown by introducing the full subcategory
\begin{equation*}
    \gcatwg{n}\subset \catwg{n}
\end{equation*}
of groupoidal weakly globular $n$-fold categories and showing (Theorem \ref{cor-gta-2}) that there is an equivalence of categories
\begin{equation}\label{fineq}
    \gcatwg{n}\bsim^n\;\simeq\;\Ho\mbox{\rm\text{($n$-types)}}.
\end{equation}

In Corollary \ref{cor1-mod-fund-wg-group} the equivalence of categories \eqref{fineq} is realized by a different pair of functors, that uses the functor $\Top\rw\gcatwg{n}$ of Blanc and the author in \cite{BP}: this provides a more explicit form for the fundamental groupoidal weakly globular \nfol groupoid of a space, which is independent on \cite{Ta}.

Our main result, Theorem \ref{cor-the-disc-func}, is that there are comparison functors
\begin{equation*}
    \Qn:\Tan \leftrightarrows \catwg{n}:\Discn
\end{equation*}
inducing equivalences of categories
\begin{equation*}
    \Tan\bsim^n\;\simeq\;\catwg{n}\bsim^n
\end{equation*}
after localization with respect to the $n$-equivalences. The rigidification functor  $\Qn$ is from Theorem \ref{the-funct-Qn} while the functor $\Discn$, called \emph{discretization functor}, is built in this chapter.

The idea of the functor $\Discn$ is to replace the homotopically discrete sub-structures in $X\in\catwg{n}$ by their discretizations in order to recover the globularity condition. This affects the Segal maps, which from being isomorphisms in $X$ become $(n-1)$-equivalences in $\Discn X$.

 However, as outlined in the introduction to chapter \ref{chap8}, for this method to work the discretization maps need to have functorial sections.

For this reason, we use the category $\ftawg{n}$ introduced in Definition \ref{def-fta-1}. We build in Proposition \ref{pro-fta-tam-1} a functor
 \begin{equation*}
   \Dn:\ftawg{n}\rw \ta{n}.
 \end{equation*}
We showed in Theorem \ref{pro-fta-1} that there is a functor
\begin{equation*}
        G_n:\catwg{n}\rw\ftawg{n}
    \end{equation*}
    and an $n$-equivalence $G_nX \rw X$ for each $X\in \catwg{n}$. The discretization functor
 \begin{equation*}
 \Discn : \catwg{n}\rw\ta{n}
 \end{equation*}
  is defined to be the composite
\begin{equation*}
    \catwg{n}\xrw{G_n}\ftawg{n}\xrw{D_n}\ta{n}\;.
\end{equation*}
This chapter is organized as follows. In Section \ref{sec-fta-to-tam}, Proposition \ref{pro-fta-tam-1} we define the functor $\Dn:\ftawg{n}\rw \ta{n}$ and establish its properties. As a consequence, and using the previous results of Proposition \ref{pro-n-equiv}, we show in Corollary \ref{cor-2-out-of} that $n$-equivalences in $\tawg{n}$ have the 2-out-of-3 property.

In Section \ref{sec-disc-final} the functor $D_n$ is used in Definition \ref{def-disc-func} to build the discretization functor from $\catwg{n}$ to $\ta{n}$. Together with the results of chapters \ref{chap7} and \ref{chap8} this leads to the main result Theorem \ref{cor-the-disc-func}.

In Section \ref{sec-group-wg-nfol-cat} we define groupoidal weakly globular $n$-fold categories and, using the results of the previous sections we show in Theorem \ref{cor-gta-2} that they are an algebraic model of $n$-types. In Section \ref{mod-fund-wg-group} we realize the equivalence of categories of Theorem \ref{cor-gta-2} through a different pair of functors, using the results of Blanc and the authors in \cite{BP}. This provides a more convenient model for the fundamental groupoidal weakly globular \nfol groupoid of a space, which is very explicit and is independent on \cite{Ta}. We illustrate this with some examples in low dimensions.

\section{From ${\pmb{\ftawg{n}}}$ to Tamsamani $n$-categories}\label{sec-fta-to-tam}
In this section we define a  functor
\begin{equation*}
    \D{n}:\ftawg{n}\rw \ta{n}
\end{equation*}
and we study its properties. As a corollary, using our previous results, we also establish that $n$-equivalences in $\tawg{n}$ have the 2-out-of-3 property.

\subsection{The idea of the functor $\pmb{\D{n}}$}\label{subs-fn-dn-idea}
 The idea of the functor $\D{n}$ is to replace the homotopically discrete sub-structures in $X\in\ftawg{n}$ by their discretization, thus recovering the globularity condition. From the definition of $\ftawg{n}$, this can be done in a functorial way. This discretization process goes at the expenses of the Segal maps, which from being isomorphisms in $\ftawg{n}$ become higher categorical equivalences, so we obtain objects of $\ta{n}$.

The construction of $\D{n}$ is inductive, and we first discretize the structure at level 0 via a functor
\begin{equation*}
  \rz : \ftawg{n} \rw \funcat{}{\ftawg{n-1}}
\end{equation*}
such that $(\rz X)_0$ is discrete for all $X \in \ftawg{n}$.

 The definition of $\rz$ is based on the following general construction. Let $Y\in\funcat{}{\clC}$, $Y_0^d\in\clC$ and suppose there are maps in $\clC$
\begin{equation*}
  \zg(Y_0): Y_0 \rw Y_0^d \qquad  \zg'(Y_0): Y^d_0 \rw Y_0
\end{equation*}
natural in $Y$, such that $\zg(Y_0)\zg'(Y_0)=\Id$ and such that a morphism $F: Y\rw Z$ in $\funcat{}{\clC}$ induces commuting diagrams
\begin{equation}\label{eq1-subs-fn-dn-idea}
\xymatrix@C=50pt@R=30pt{
Y_{0} \ar^{F_{0}}[r] \ar_{\zg(Y_{0})}[d]& Z_{0}\ar^{\zg(Z_{0})}[d]\\
Y^d_{0} \ar_{F^d_{0}}[r] & Z^d_{0}
}
\qquad
\xymatrix@C=50pt@R=30pt{
Y^d_{0} \ar^{F^d_{0}}[r] \ar_{\zg'(Y_{0})}[d]& Z^d_{0}\ar^{\zg'(Z_{0})}[d]\\
Y_{0} \ar_{F_{0}}[r] & Z_{0}
}
\end{equation}
Define $\rz Y$ as follows:
\begin{equation*}
  (\rz Y)_k =
  \left\{
    \begin{array}{ll}
      Y_0^d, & k=0 \\
      Y_k, & k>0
    \end{array}
  \right.
\end{equation*}
The face operators $\pt'_0,\pt'_1:Y_1\rightrightarrows Y_0^d$ are given by $\pt'_i=\zg(Y_0)\pt_i$\; $i=0,1$ and the degeneracy operator $\zs':Y^d_0\rw Y_1$ by $\zs'=\zs\zg'(Y_0)$ where $\pt_0,\pt_1:Y_1\rightrightarrows Y_0$ and $\zs:Y_0\rw Y_1$ are the face and degeneracy operators of $Y$. All other face and degeneracy operators of $\rz Y$ are as in $Y$. Since $\zg(Y_0)\zg'(Y_0)=\Id$, all simplicial identities for $\rz Y$ hold so that $\rz Y\in\funcat{}{\clC}$.

Let $F:Y\rw Z$ be a map in $\funcat{}{\clC}$. From the commutativity of \eqref{eq1-subs-fn-dn-idea} this induces a map in $\funcat{}{\clC}$, $\rz F: \rz Y\rw \rz Z$, so $\rz$ is a functor
\begin{equation*}
  \rz:\funcat{}{\clC}\rw \funcat{}{\clC}\;.
\end{equation*}
We apply this construction to the case where $Y=\Nu{1}X$ with $X\in \ftawg{n}$, $\zg: X_0\rw X_0^d$ is the discretization map and $\zg'=X^d_0 \rw X_0$ a functorial section. As observed in Remark \ref{rem-fta-2}, $\Nu{1}X\in\funcat{}{\ftawg{n-1}}$ and $\zg,\zg'$ are maps in $\ftawg{n-1}$; thus by definition of morphism in $\ftawg{n}$ (see Definition \ref{def-fta-1}) a morphism $F: X\rw Y$ in $\ftawg{n}$ induces commuting diagrams as in \eqref{eq1-subs-fn-dn-idea}.

So all the conditions to apply the previous construction are met and we define the functor
\begin{equation*}
  \rz : \ftawg{n}\rw \funcat{}{\ftawg{n-1}}
\end{equation*}
(see Definition \ref{def-fta-tam-1}). The effect of $\rz$ is to discretize the object $X_0$ to $X_0^d$. This, however, does not yet produce an object of $\ta{n}$ since, for $k\geq 0$, $(\rz X)_k=X_k$ is in $\ftawg{n-1}$, not in $\ta{n-1}$.

We perform the rest of the discretization of $X$ inductively. Namely, we define inductively
\begin{equation*}
  \Dn:\ftawg{n}\rw \ta{n}
\end{equation*}
by
\begin{equation*}
  D_2=\rz, \qquad \Dn=\overline{D}_{n-1}\rz
\end{equation*}
The effect of $\Dn$ is to discretize all the homotopically discrete substructures of $X\in\ftawg{n}$, thus recovering the globularity condition. The proof that $\Dn X \in \ta{n}$ needs checking the Segal maps condition, and this is done inductively in the proof of Proposition \ref{pro-fta-tam-1}.

%%%%%%%%%%%%%%%%%%%%%%%%%%%%%%%%%%%%%%%%%%%%%%%%%%%%%%%%%%%%%%%%%%%%%%%%%%%%%%%%%%%%%%%%%%%%%%%%
\subsection{The functor $\pmb{\D{n}}$: definition and properties}
\begin{definition}\label{def-fta-tam-1}
    Let $\rz:\ftawg{n}\rw\funcat{}{\ftawg{n-1}}$ be given by
    \begin{equation*}
        (\rz X)_k=
        \left\{
          \begin{array}{ll}
            X_0^d, & k=0 \\
            X_k, & k>0\;.
          \end{array}
        \right.
    \end{equation*}
    The face operators $\pt'_0,\pt'_1:X_1\rightrightarrows X_0^d$ are given by $\pt'_i=\zg\pt_i$, $i=0,1$ and the degeneracy $\zs':X_0^d\rw X_1$ by $\zs'=\zs\zg'$ where $\pt_0,\pt_1: X_1 \rightrightarrows X_0$ and $\zs:X_0 \rw X_1$ are the corresponding face and degeneracy operators of $X$, $\zg:X_0\rw X_0^d$ is the discretization map and $\zg':X^d_0\rw X_0$ is a functorial section. All other face and degeneracy maps of $\rz X$ are as in $X$.

    Note that by definition of $\ftawg{n}$ and by Remark \ref{rem-fta-2} the maps $\zg$ and $\zg'$ are morphisms in $\ftawg{n-1}$, therefore such are $\pt'_i$ and $\zs'$. Since $\zg \zg'=id$, all simplicial identities are satisfied, thus $R_0 X\in\funcat{}{\ftawg{n-1}}$.
\end{definition}
\begin{remark}\label{rem-fta-tam-1}
    By definition of $\ftawg{n}$, given $f:X\rw Y$ in $\ftawg{n}$ there is a commuting diagram
    \begin{equation}\label{eq-rem-fta-tam-1}
    \xymatrix@C=40pt{
    X_0 \ar^{f_0}[r] & Y_0 \\
    X_0^d \ar^{\zg'(X_0)}[u] \ar[r] & Y_0^d \ar_{\zg'(Y_0)}[u]
    }
    \end{equation}
    and this induces a morphism in $\funcat{}{\ftawg{n-1}}$
    \begin{equation*}
      \rz f:\rz X\rw \rz Y.
    \end{equation*}
    Thus $\rz$ is a functor. Note that while $\rz X$ could be defined for any $X\in\tawg{n}$, given a morphism $f$ in $\tawg{n}$ since in general \eqref{eq-rem-fta-tam-1} does not commute, one cannot define $\rz f$ as above.
\end{remark}
\begin{lemma}\label{lem-fta-tam-1}
Let $\rz$ be as in Definition \ref{def-fta-tam-1}. Then:\mk
\begin{itemize}
  \item [a)] $\rz$ is identity on objects and commutes with pullbacks over discrete objects.\bk

  \item [b)] $\p{2}=\ovl{p}\rz$, \; $\q{2}=\ovl{q}\rz$ \; while for $n>2$\mk

             $\ovl{\p{n-1}}\rz=\rz\p{n}$, \; $\ovl{\q{n-1}}\rz=\rz\q{n}$.\bk

  \item [c)] For each $X\in\ftawg{n}$ the Segal maps of $\rz X$
  \begin{equation*}
    (\rz X)_{k} \rw \pro{(\rz X)_{1}}{(\rz X)_{0}}{k}
  \end{equation*}
  are $(n-1)$-equivalences for all $k\geq 2$.
\end{itemize}
\end{lemma}
\begin{proof}
\

a) This is immediate by the definition of $\rz$ since, if $X\rw Z\lw Y$ is a pullback in $\ftawg{n}$ with $Z$ discrete, $(X\tiund{Z}Y)_0^d=X_0^d\tiund{Z}Y_0^d$ by Lemma \ref{lem-copr-hom-disc}.\bk

b) If $X\in\ftawg{2}$.
\begin{equation*}
  (\ovl{p}\rz X)_0= p X_0^d = X_0^d =(\p{2}X)_0
\end{equation*}
while for $k>1$
\begin{equation*}
  (\ovl{p}\rz X)_k= p X_k = (\p{2}X)_k
\end{equation*}
so that $\ovl{p}\rz X=\p{2}X$. Similarly one shows that $\ovl{q}\rz X=\q{2}X$.

If $X\in\ftawg{n}$ for $n>2$, we have
\begin{equation*}
  (\ovl{\p{n-1}}\rz X)_0=X_0^d=(\p{n}X)_0^d=(\rz\p{n}X)_0
\end{equation*}
while for $k>0$,
\begin{equation*}
  (\ovl{\p{n-1}}\rz X)_k=\p{n-1}(\rz X)_k=\p{n-1}X_k=\rz(\p{n}X)_k\;.
\end{equation*}
In conclusion $\ovl{\p{n-1}}\rz X=\rz\p{n}X$. Similarly one shows that
 \begin{equation*}
 \ovl{\q{n-1}}\rz X=\rz\q{n}X.
 \end{equation*}

\bk
c) For each $k\geq 2$ the Segal maps for $\rz X$ are
\begin{equation*}
  (\rz X)_k =X_k\rw \pro{X_1}{X_0^d}{k}=\pro{(\rz X)_1}{(\rz X)_0}{k}
\end{equation*}
and these are $(n-1)$-equivalences since $X\in\ftawg{n}$.

\end{proof}
\begin{proposition}\label{pro-fta-tam-1}
There is a functor
\begin{equation*}
    \Dn:\ftawg{n}\rw \ta{n}
\end{equation*}
defined inductively by
\begin{equation*}
    D_2=\rz, \qquad \Dn=\ovl{D}_{n-1}\circ \rz \qquad \text{for }\; n>2
\end{equation*}
where $\rz$ is as in Definition \ref{def-fta-tam-1}, such that
\begin{itemize}
  \item [a)] $\Dn$ is identity on discrete objects and commutes with pullbacks over discrete objects.\mk

  \item [b)] $\p{2}D_2=\p{2}$, $\q{2}D_2=\q{2}$, while for $n>2$
  \begin{equation*}
    \p{n}\Dn X=\Dnm \p{n} X,\quad \q{n}\Dn X=\Dnm \q{n} X.
  \end{equation*}

  \item [c)] For each $X\in\ftawg{n}$ and $a,b\in X_0^d$,
  \begin{equation*}
    (\Dn X)(a,b)=\Dnm X(a,b)\;.
  \end{equation*}

  \item [d)] $\Dn$ preserves and reflects $n$-equivalences.
\end{itemize}
\end{proposition}
\begin{proof}
By induction on $n$. When $n=2$, $D_2 X=\rz X\in\funcat{}{\Cat}$ is such that $(D_2 X)_0=X_0^d$ is discrete and, by Lemma \ref{lem-fta-tam-1} c), the Segal maps are equivalences of categories. Thus $D_2 X\in\ta{2}$. Note also that by Lemma \ref{lem-fta-tam-1} b)
\begin{equation*}
  \p{2}D_2 X=\ovl{p}\rz X=\ovl{p}X=\p{2}X\;.
\end{equation*}
Similarly, $\q{2}D_2=\q{2}$. By Lemma \ref{lem-fta-tam-1} a), $D_2$ satisfies a).

Let $f:X\rw Y$ be a 2-equivalence. For all $a,b\in (D_2 X)^d_0=X_0^d$, $(D_2 f)(a,b)=f(a,b)$ is an equivalence of categories. Also, $\p{2}D_2 f=\p{2}f$ is an equivalence of categories. Thus by definition $D_2 f$ is a 2-equivalence.

Suppose that $f:X\rw Y$ is such that $D_2 f$ a 2-equivalence. Then for all $a,b\in (D_2 X)^d_0=X_0^d$, $(D_2 f)(a,b)=f(a,b)$ is an equivalence of categories. Also, $\p{2}D_2 f=\p{2}f$ is an equivalence of categories. Thus by definition $f$ is a 2-equivalence. This completed the proof of d) when $n=2$.

Suppose, inductively, that the proposition holds for $(n-1)$ and let $X\in\ftawg{n}$. Then by induction hypothesis a)
\begin{equation*}
    (\Dn X)_k =
    \left\{
      \begin{array}{ll}
        \Dnm X_0^d=X_0^d & k=0 \\
        \Dnm X_k, & k>0\;.
      \end{array}
    \right.
\end{equation*}
Thus by induction hypothesis $(\Dn X)_k\in\ta{n-1}$ for all $k\geq 0$ with $(\Dn X)_0$ discrete.

To show that $\Dn X\in\ta{n}$ it remains to show that the Segal maps are $(n-1)$-equivalences. Since $X\in\ftawg{n}$, for each $k\geq 2$ the map
\begin{equation*}
    \mu_k: X_k\rw \pro{X_1}{X_0^d}{k}
\end{equation*}
is a $(n-1)$-equivalence. By inductive hypotheses a) and d) this induces a $(n-1)$-equivalence
\begin{align*}
&\Dnm \mu_k:\Dnm X_k=(\Dn X)_k \rw \Dnm(\pro{X_1}{X_0^d}{k})\cong\\
& \cong \pro{(\Dn X)_1}{(\Dn X)_0}{k}\;.
\end{align*}
This shows that the Segal maps of $\Dn X$ are $(n-1)$-equivalences. We conclude that $\Dn X\in\ta{n}$.
\bk

a) This follows from Lemma \ref{lem-fta-tam-1} and the inductive hypothesis.
\bk

b) Recalling that $\p{n}=\ovl{\p{n-1}}:\ta{n}\rw\ta{n-1}$, using the inductive hypothesis and Lemma \ref{lem-fta-tam-1} b) we obtain
\begin{align*}
&\p{n}\Dn=\ovl{\p{n-1}}\ovl{\Dnm}\rz \\
& = \ovl{D_{n-2}}\ovl{\p{n-1}}\rz = \ovl{D_{n-2}}\rz\p{n}=\Dnm\p{n}\;.
\end{align*}
The proof for $\q{n}\Dn X$ is similar.
\bk

c) By definition of $X(a,b)$, we have a pullback in $\ftawg{n-1}$
\begin{equation*}
\xymatrix@R=35pt@C=40pt{
X(a,b) \ar[r]\ar@{_{(}->}[d] & \di{n-1}\{a\}\times \di{n-1}\{b\} \ar@{_{(}->}[d]\\
X_1 \ar[r] & X_0^d \times X_0^d
}
\end{equation*}
Therefore, by a), we also have a pullback
\begin{equation*}
\xymatrix@R=38pt@C=40pt{
\Dnm X(a,b) \ar[r]\ar[d] & \di{n-1}\{a\}\times \di{n-1}\{b\} \ar@{_{(}->}[d]\\
(\Dn X)_1=\Dnm X_1 \ar[r] & X_0^d \times X_0^d=(\Dn X)_0\times (\Dn X)_0
}
\end{equation*}
so that
\begin{equation*}
  (\Dn X)(a,b)=\Dnm X(a,b)\;.
\end{equation*}
\bk

d) Let $f:X\rw Y$ be an $n$-equivalence in $\ftawg{n}$. By c), for each $a,b\in(\Dn X)_0=X_0^d$ it is
\begin{equation*}
    (\Dn f)(a,b)=\Dnm f(a,b)
\end{equation*}
and this is a $(n-1)$-equivalence by the inductive hypothesis applied to the $(n-1)$-equivalence $f(a,b)$. Further, by b)
\begin{equation*}
  \p{n}\Dn f=\Dnm \p{n} f
\end{equation*}
 is also a $(n-1)$-equivalence by inductive hypothesis applied to the $(n-1)$-equivalence $\p{n}f$. This shows that $\Dn f$ is a $n$-equivalence.

 Let $f:X\rw Y$ be such that $D_n f$ is an $n$-equivalence. By c), for each $a,b\in(\Dn X)_0=X_0^d$ it is
\begin{equation*}
    (\Dn f)(a,b)=\Dnm f(a,b)
\end{equation*}
and this is a $(n-1)$-equivalence. By inductive hypothesis we conclude that $f(a,b)$ is an $(n-1)$-equivalence. Further, by b)
\begin{equation*}
  \p{n}\Dn f=\Dnm \p{n} f
\end{equation*}
 is also a $(n-1)$-equivalence. So by inductive hypothesis $ \p{n}f$ is a $(n-1)$-equivalence. In conclusion $f$ is a $n$-equivalence.

\end{proof}
\begin{example}\label{example-fta-tam-1}
  Let $X\in\ftawg{3}$, so that
  \begin{equation*}
    J_3 X\in\funcat{2}{\Cat}\;.
  \end{equation*}
A picture of the corner of $J_3 X$ is found in Figure \ref{Corner-J3-1} on page \pageref{Corner-J3} where the structures in red are homotopically discrete and they are equipped with functorial sections to the discretization maps. In Figure \ref{Corner-J3-2} on page \pageref{Corner-J3} we depict the corner of
\begin{equation*}
  J_3 D_3 X\in\funcat{2}{\Cat}\;.
\end{equation*}
We see that the homotopically discrete substructures in $X$ have been replaced by discrete ones (also in red).

\end{example}
\begin{corollary}\label{cor-2-out-of}
  $n$-Equivalences in $\tawg{n}$ have the 2-out-of-3 property.
\end{corollary}
\index{Two-out-of three property}
\begin{proof}
By Proposition \ref{pro-n-equiv} the only case that remains to be checked is when we have morphisms
\begin{equation*}
  X\xrw{g} Z \xrw{h} Y
\end{equation*}
in $\tawg{n}$ such that $hg$ and $g$ are $n$-equivalences, but (unlike in Proposition \ref{pro-n-equiv} c)), no further assumptions are required on $g$. We need to show that $h$ is an $n$-equivalence. Since, by Theorem \ref{pro-fta-1}, $G_n$ preserves $n$-equivalences and, by Proposition \ref{pro-fta-tam-1}, $\Dn$ preserves $n$-equivalences, we have morphisms in $\ta{n}$
\begin{equation*}
  \Dn G_n X \xrw{\Dn G_n g} \Dn G_n Z \xrw{\Dn G_n h} \Dn G_n Y
\end{equation*}
in which $\Dn G_n h$ and the composite $(\Dn G_n h)(\Dn G_n g)$ are $n$-equivalences. Since by \cite{S2} $n$-equivalences in $\ta{n}$ have the 2-out-of-3 property, this implies that $\Dn G_n h$ is an $n$-equivalence. Since, by Proposition \ref{pro-fta-tam-1} d), $\Dn$ reflects $n$-equivalences, we conclude that $G_n h$ is an $n$-equivalence.

On the other hand, we have a commuting diagram
\begin{equation*}
\xymatrix@R=30pt@C=35pt{
G_n Z \ar^{G_n h}[r] \ar^{g_{n}(Z)}[d] & G_n Y \ar^{g_{n}(Y)}[d]\\
Z \ar_{h}[r] & Y
}
\end{equation*}
in which the vertical maps and the top horizontal maps are $n$-equivalences. Thus, by Proposition \ref{pro-n-equiv} d), we conclude that the composite
\begin{equation}\label{eq-cor-2-out-of}
  G_n Z \xrw{g_{n}(Z)} Z \xrw{h} Y
\end{equation}
is an $n$-equivalence. By Theorem \ref{pro-fta-1} b),the map $g_{n}(Z)^d:(G_n Z)^d_0\rw Z_0^d$ is surjective. Hence the morphisms \eqref{eq-cor-2-out-of} satisfies the hypotheses of Proposition \ref{pro-n-equiv} e) and we conclude that $h$ is an $n$-equivalence.
\end{proof}

\section{The discretization functor and the comparison result}\label{sec-disc-final}
In this section we define the discretization functor
\begin{equation*}
  \Discn:\catwg{n}\rw \ta{n}
\end{equation*}
and we establish the main result of this work, theorem \ref{cor-the-disc-func}, asserting that the functors $\Discn$ and $Q_n$ induce an equivalence of categories
 \begin{equation*}
        \ta{n}\bsim^n\;\simeq \; \catwg{n}\bsim^n
    \end{equation*}

\subsection{The idea of the functor $\pmb{\Discn}$}\label{subs-idea-discn} The idea of the discretization functor
\begin{equation*}
  \Discn: \catwg{n} \rw \ta{n}
\end{equation*}
is to discretize the homotopically discrete sub-structures in the multinerve of objects of $\catwg{n}$ to recover the globularity condition. As discussed in Chapter \ref{chap8}, this needs functorial sections to the discretization maps of the homotopically discrete sub-structures.  For this reason we introduced the category $\ftawg{n}$ and the discretization process from this category is the functor
\begin{equation*}
  \Dn:\ftawg{n}\rw \ta{n}
\end{equation*}
of Section \ref{def-fta-tam-1}.

We define the discretization functor to be the composite
\begin{equation*}
   \Discn: \catwg{n}\xrw{G_n}\ftawg{n}\xrw{\Dn}\ta{n}\;.
\end{equation*}
where $G_n$ is an in Theorem \ref{pro-fta-1} and $D_n$ as in Proposition  \ref{pro-fta-tam-1}.

This realizes the idea of discretizing the homotopically discrete sub-structures in each $X\in\catwg{n}$, but after replacing $X$ with the $n$-equivalent $G_n X \in \ftawg{n}$.

The main property of $\Discn$ is that, for each $X\in\catwg{n}$, $\Discn X$ and $X$ are suitably equivalent in $\tawg{n}$. We show this fact in Theorem \ref{the-disc-func}, where we prove that there is a zig-zag of $n$-equivalences in $\tawg{n}$, of the form
\begin{equation}\label{zizag}
  \Discn \lw Q_n \Discn X \rw X.
\end{equation}
This relies on Proposition \ref{pro-fta-tam-2}, establishing that, for each $X\in \ftawg{n}$, $Q_n \Dn X_n=Q_n X$. When applied to $G_n X$ (for $X\in \catwg{n}$) this fact implies
\begin{equation*}
  Q_n \Discn X= Q_n\Dn G_n X=Q_n G_n X.
\end{equation*}
The zig-zag \eqref{zizag} is then obtained using the maps $s_n(\Discn X)$, $s_n (G_n X)$, $g_n(X)$ as follows
\begin{equation*}
  \Discn \xlw{s_n(\Discn X)} Q_n\Discn X=Q_n G_n X \xrw{s_n (G_n X)} G_n X \xrw{g_n(X)} X
\end{equation*}

The proof of Proposition \ref{pro-fta-tam-2} relies on the definition of $Q_n$ as well as on the property of the functor $Tr_n$ established in Lemma \ref{lem-from-lta-to-pseu}. The latter gives conditions on $X\in\catwg{n}$ and $Y \in \lta{n}$ under which $Tr_n X =Tr_n Y$. Using the definition of $P_n: \tawg{n} \rw \lta{n}$ in Theorem \ref{the-funct-Qn}, we show that, given $X\in\catwg{n},$ $P_n \Dn X\in \lta{n}$ and $P_n X\in\cat{n}$ satisfy these conditions, and therefore
\begin{equation*}
    Tr_n P_n \Dn X= Tr_n P_n X.
\end{equation*}
In turn, this implies
\begin{equation*}
    \Qn \Dn X = \St Tr_n P_n \Dn X= \St Tr_n P_n X= \Qn X\;.
\end{equation*}
which is Proposition \ref{pro-fta-tam-2}.

\subsection{The comparison result} In this section we prove our main comparison result between Tamsamani $n$-categories and weakly globular \nfol categories. We first need to establish a number of properties about the functor $D_n$ of the previous section, see Proposition \ref{pro-fta-tam-2} and Lemma \ref{lem-fta-tam-3} below.

\begin{definition}\label{def-disc-func}
Define the discretization functor \index{Discretization!- functor}
\begin{equation*}
  \Discn:\catwg{n}\rw \ta{n}
\end{equation*}
  to be the composite
\begin{equation*}
    \catwg{n}\xrw{G_n}\ftawg{n}\xrw{\Dn}\ta{n}\;.
\end{equation*}
where $G_n$ is an in Theorem \ref{pro-fta-1} and $D_n$ as in Proposition  \ref{pro-fta-tam-1}.
\end{definition}

The following Lemmas and Proposition establish some facts about the functor $D_n$ of Proposition \ref{pro-fta-tam-1} which will be needed to study the properties of the discretization functor $\Discn$.
\begin{lemma}\label{lem-fta-tam-2}
Let $X\in\ftawg{n}$
\begin{itemize}
  \item [a)]If $\uk\in\dop{n-1}$ is such that $k_j\neq 0$ for all $1\leq j\leq n-1$. Then $(\Dn X)_{\uk}=X_{\uk}$.
  \item [b)]If $\uk \rw \us$ is a morphism in $\dop{n-1}$ with $k_j= 0$ for some $1\leq j\leq n-1$ and $s_i\neq 0$ for all $1\leq i\leq n-1$, then the map $(\Dn X)_{\uk} \rw (\Dn X)_{\us}=X_{\us}$ factors as
      \begin{equation*}
        (\Dn X)_{\uk} \rw X_{\uk} \rw X_{\us}.
      \end{equation*}
\end{itemize}

\end{lemma}
\begin{proof}
By induction on $n$. When $n=2$, $D_2 X=\rz X$ has $(D_2 X)_k=X_k$ for all $k\neq 0$, proving a). By construction of $R_0 X$, if $s>0$, the map
\begin{equation*}
  (D_2 X)_0=X_0^d \rw (D_2 X)_s=X_s
\end{equation*}
factors as
\begin{equation*}
  X_0^d \rw X_0 \rw X_s
\end{equation*}
proving b).

Suppose, inductively, that the lemma holds for $(n-1)$.

a) Denote $\ur=(\seqc{k}{2}{n-1})$. Then by inductive hypothesis applied to $X_{k_1}$ we have
\begin{equation*}
    (\Dn X)_{\uk}=(\Dnm X_{k_1})_{\ur}=(X_{k_1})_{\ur}=X_{\uk}\;.
\end{equation*}

b) By the construction of $R_0 X$, for each $s_1 >0$ the map
\begin{equation*}
  (R_0 X)_0=X_0^d \rw (R_0 X)_{s_1}=X_{s_1}
\end{equation*}
factors as
\begin{equation*}
  X_0^d \rw X_0 \rw X_{s_1}.
\end{equation*}
Thus applying $D_{n-1}$ we see that the map
\begin{equation*}
  (D_n X)_0=X_0^d \rw (D_n X)_{s_1}=D_{n-1}X_{s_1}
\end{equation*}
factors as
\begin{equation}\label{eq-prop-dn}
  X_0^d \rw D_{n-1}X_0 \rw D_{n-1}X_{s_1}
\end{equation}
Let $\uk \rw \us$ be a morphism in $\dop{n-1}$ satisfying hypotheses b), and denote $\uk= (k_1, \ur)$, $\us= (s_1, \uv)$ so we have a corresponding morphism in $\dop{n-2}$ $\ur \rw \uv$. We distinguish two cases:

i) Consider first the case $k_1=0$. By naturality, \eqref{eq-prop-dn} gives a commuting diagram in $\Cat$

\begin{equation}\label{diag1-prop-dn}
(\Dn X)_{\uk}=
\xymatrix{
X_0^d \ar[r] \ar[dr] & (\Dnm X_0)_{\underline{v}} \ar[r]& (\Dnm X_{s_1})_{\underline{v}}\\
& (\Dnm X_0)_{\ur} \ar[u] \ar[r] &  (\Dnm X_{s_1})_{\ur} \ar[u]
}
=(\Dn X)_{\us}
\end{equation}
where we used the fact that, since $X_0^d$ is discrete $(X_0^d)_{\uv}=(X_0^d)_{\us}=X_0^d$ and $(\Dn X)_{\uk}=X_0^d$ since $k_1=0$.

Suppose that $r_i \neq 0$ for all $1\leq i\leq (n-2)$. Then by part a) $(D_{n-1}X_0)_{\ur}=X_{0\ur}=X_{\uk}$. So by diagram \ref{diag1-prop-dn}  we see that the map $ (\Dn X)_{\uk} \rw (\Dn X)_{\us}$ factors through $X_{\uk}$.

Suppose $r_i=0$ for some $1\leq i\leq (n-2)$. Then, since by hypothesis $v_j \neq 0$ for all $1\leq j\leq (n-2)$ we can apply the inductive hypothesis to $X_0$ and deduce that the map
\begin{equation*}
  (D_{n-1}X_0)_{\ur} \rw (D_{n-1}X_0)_{\uv}
\end{equation*}
factors through $X_{0\ur}=X_{\uk}$. From the commuting diagram \ref{diag1-prop-dn} we deduce that the map $ (\Dn X)_{\uk} \rw (\Dn X)_{\us}$ factors through $X_{\uk}$.
\medskip

ii) Consider now the case $k_1>0$. By hypothesis $r_j=0$ for some $1\leq j\leq (n-2)$ and $(\Dn X)_{\uk}=D_{n-1}X_{k_1}$. The morphism
 $ (\Dn X)_{\uk} \rw (\Dn X)_{\us}$ factors as
\begin{equation}\label{eq-prop-dn-2}
  (\Dn X)_{\uk}=(\Dn X)_{(k_1,\ur)} \rw (\Dn X)_{(k_1,\uv)} \rw (\Dn X)_{(s_1,\uv)}=(\Dn X)_{\us}.
\end{equation}
But $(\Dn X)_{(k_1,\ur)}=(D_{n-1}X_{k_1})_{\ur} $ and $(\Dn X)_{(k_1,\uv)}=(D_{n-1}X_{k_1})_{\uv} $. By induction hypothesis applied to $X_{k_1}$ the map
\begin{equation*}
  (D_{n-1}X_{k_1})_{\ur} \rw (D_{n-1}X_{k_1})_{\uv}
\end{equation*}
factors through $X_{k_1 \ur}=X_{\uk}$. Thus by \eqref{eq-prop-dn-2} we see that the map $ (\Dn X)_{\uk} \rw (\Dn X)_{\us}$ factors through $X_{\uk}$.

\end{proof}
\begin{proposition}\label{pro-fta-tam-2}
Let $X\in\ftawg{n}$, then $\Qn\Dn X=\Qn X$.
\end{proposition}
\begin{proof}
By induction on $n$. Let $X\in\ftawg{2}$. It is immediate that $\rz X$ and $X$ satisfy the hypotheses of Lemma \ref{lem-from-lta-to-pseu} so that
\begin{equation*}
    Tr_2 \rz X= Tr_2 X\;.
\end{equation*}
Hence
\begin{equation*}
    Q_2D_2 X=\St Tr_2 \rz X= \St Tr_2 X=Q_2 X\;.
\end{equation*}
Suppose, inductively, the statement holds for $(n-1)$ and let $X\in\ftawg{n}$.

We claim that $P_n\Dn X$ and $P_n X$ satisfy the hypotheses of Lemma \ref{lem-from-lta-to-pseu}, where $P_n:\tawg{n}\rw \lta{n}$ is as in Theorem \ref{the-funct-Qn}.
First note that since $X\in\ftawg{n}$, in particular $X\in\catwg{n}$, so by remark \ref{rem-pn-catwg} $P_n X \in \catwg{n}$. We now show that $(P_n D_n X)_{\uk}$ is discrete for all $\uk \in \dop{n-1}$ such that $k_j=0$ for some $1\leq j \leq (n-1)$.

 By definition of $P_n$ there is a pullback in $\funcat{n-1}{\Cat}$
\begin{equation}\label{eq-lem-fta-tam-5}
\xymatrix@R=35pt{
P_n\Dn X \ar[r]\ar[d]_{z} & \Dn X\ar[d]\\
\di{n}\Qnm \q{n}\Dn X \ar[r] & \di{n}\q{n}\Dn X
}
\end{equation}
On the other hand, by Proposition \ref{pro-fta-tam-1} and the inductive hypothesis
\begin{equation*}
    \Qnm \q{n}\Dn X=\Qnm \Dnm \q{n}X=\Qnm \q{n}X
\end{equation*}
so that \eqref{eq-lem-fta-tam-5} coincides with
\begin{equation}\label{eq-lem-fta-tam-6}
\xymatrix@R=35pt{
P_n\Dn X \ar[r]\ar[d]_{z} & \Dn X\ar[d]\\
\di{n}\Qnm \q{n} X \ar[r] & \di{n}\q{n}\Dn X
}
\end{equation}

Since pullbacks in $\funcat{n-1}{\Cat}$ are computed pointwise, for each $\uk\in\dop{n-1}$ the diagram \eqref{eq-lem-fta-tam-6}  gives rise to a pullback in $\Cat$
\begin{equation}\label{eq-lem-fta-tam-8}
\xymatrix@R=35pt{
(P_n \Dn X)_{\uk} \ar[r]\ar[d]_{z_{\uk}} & (\Dn X)_{\uk}\ar[d]\\
d(\Qnm \q{n} X)_{\uk} \ar[r] & dq(\Dn X)_{\uk}
}
\end{equation}

If $k_j=0$ for some $1\leq j \leq (n-1)$, then $(\Dn X)_{\uk}$ is discrete (since $\Dn X\in\ta{n}$) hence the right vertical map in \eqref{eq-lem-fta-tam-8} is the identity, and thus so is the left vertical map in \eqref{eq-lem-fta-tam-8}. That is
\begin{equation*}
    (P_n \Dn X)_{\uk}= d(\Qnm \q{n} X)_{\uk}
\end{equation*}
and
\begin{equation}\label{eq-lem-fta-tam-10}
    p(P_n \Dn X)_{\uk}= (\Qnm \q{n} X)_{\uk}\;.
\end{equation}
We also have a pullback in $\funcat{n-1}{\Cat}$
\begin{equation}\label{eq-lem-fta-tam-7}
\xymatrix@R=35pt{
P_n X \ar[r]\ar[d]_{t} & X\ar[d]\\
\di{n}\Qnm \q{n} X \ar[r] & \di{n}\q{n} X
}
\end{equation}
and for each $\uk \in \dop{n-1}$ a pullback in $\Cat$
\begin{equation}\label{eq-lem-fta-tam-9}
\xymatrix@R=35pt{
(P_n X)_{\uk} \ar[r]\ar[d]_{t_{\uk}} & X_{\uk}\ar[d]\\
d(\Qnm \q{n} X)_{\uk} \ar[r] & dq X_{\uk}
}
\end{equation}
We now check hypotheses i) through iv) of Lemma \ref{lem-from-lta-to-pseu} for $P_n X\in\catwg{n}$ and $P_n \Dn X\in\lta{n}$.

i) Let $\uk\in\dop{n-1}$ and $\us\in\dop{n-1}$ be such that $k_j\neq 0$  and $s_j\neq 0$ for all $1 \leq j \leq (n-1)$. Then by Lemma \ref{lem-fta-tam-2} $(\Dn X)_{\uk}=X_{\uk}$. Hence the right vertical maps in \eqref{eq-lem-fta-tam-8} and \eqref{eq-lem-fta-tam-9} coincide. It follows that $z_{\uk}=t_{\uk}$ and
\begin{equation*}
    (P_n\Dn X)_{\uk} = (P_n X)_{\uk}\;.
\end{equation*}
Similarly $(P_n\Dn X)_{\us} = (P_n X)_{\us}$.
Given a morphism $\uk \rw \us$ in $\dop{n-1}$, clearly the maps
 \begin{equation*}
 (P_n X)_{\uk} \rw (P_n X)_{\us}, \qquad(P_n \Dn X)_{\uk} \rw (P_n \Dn X)_{\us}
 \end{equation*}
  coincide.

ii) Let $\uk\in\dop{n-1}$ and $\us\in\dop{n-1}$ be such that $k_j=0$ for some $1 \leq j \leq (n-1)$ and $s_i=0$ for some $1 \leq i\leq (n-1)$. Then $X_{\uk}\in\cathd{}$ so $q X_{\uk}=p X_{\uk}$. Thus from \eqref{eq-lem-fta-tam-9}, using the fact that $p$ commutes with pullbacks over discrete objects, we obtain
\begin{equation}\label{eq-lem-fta-tam-11}
    p(P_n X)_{\uk}= (\Qnm \q{n} X)_{\uk}\;.
\end{equation}
It follows from \eqref{eq-lem-fta-tam-10} and \eqref{eq-lem-fta-tam-11} that
\begin{equation*}
    (P_n \Dn X)^d_{\uk}=dp(P_n \Dn X)_{\uk}= dp(P_n X)_{\uk}=(P_n X)^d_{\uk}\;.
\end{equation*}
Similarly $(P_n \Dn X)^d_{\us}=(P_n X)^d_{\us}$.
Given a morphism $\uk \rw \us$ in $\dop{n-1}$, the maps
 \begin{equation*}
 (P_n X)_{\uk}^d \rw (P_n X)_{\us}^d, \qquad(P_n \Dn X)_{\uk} \rw (P_n \Dn X)_{\us}
 \end{equation*}
  coincide, and they are equal to the maps
   \begin{equation*}
   d(Q_{n-1} \q{n}X)_{\uk} \rw d(Q_{n-1} \q{n}X)_{\us}.
   \end{equation*}

iii) Let $\uk \rw\us$ be a morphism in $\dop{n-1}$ and suppose that $k_j \neq 0$ for all $1 \leq j \leq n-1$ and $s_i =0$ for some $1 \leq i \leq n-1$. From above, $z_{\us}=\Id$ while $t_{\us}: (P_n X)_{\us}\rw (P_n X)^d_{\us}$ is the discretization map. By construction we have commuting diagrams
\begin{equation}\label{eq-lem-fta-tam-11a}
\xymatrix@C=12pt{
(P_n \Dn X)_{\uk} \ar[r] \ar_{z_{\uk}}[d] & (P_n \Dn X)_{\us} \ar^{z_{\us}=\Id}[d]\\
d(Q_{n-1}\q{n}X)_{\uk} \ar[r] & d(Q_{n-1}\q{n}X)_{\us}
}
\quad\;
\xymatrix@C=12pt{
(P_n  X)_{\uk} \ar[r] \ar_{t_{\uk}}[d] & (P_n  X)_{\us} \ar^{t_{\us}}[d]\\
d(Q_{n-1}\q{n}X)_{\uk} \ar[r] & d(Q_{n-1}\q{n}X)_{\us}
}
\end{equation}
\bk
and, from above, $(P_n \Dn X)_{\uk}=(P_n  X)_{\uk}$, $z_{\uk}=t_{\uk}$ while
 \begin{equation*}
   (P_n \Dn X)_{\us}=d(Q_{n-1}\q{n}X)_{\us}=(P_n  X)_{\us}^d.
 \end{equation*}

We therefore see from \eqref{eq-lem-fta-tam-11a} that the map $(P_n \Dn X)_{\uk}\rw (P_n \Dn X)_{\us}$ factors as
\begin{equation*}
  (P_n \Dn X)_{\uk}=(P_n X)_{\uk} \xrw{z_{\uk}=t_{\uk}} d(Q_{n-1}\q{n}X)_{\uk} \rw d(Q_{n-1}\q{n}X)_{\us}
\end{equation*}
which is the same as
\begin{equation*}
   (P_n \Dn X)_{\uk}=(P_n X)_{\uk} \rw (P_n X)_{\us} \xrw{t_{\us}} (P_n X)^d_{\us}= d(Q_{n-1}\q{n}X)_{\us}\;.
\end{equation*}
This proves hypothesis iii) in Lemma \ref{lem-from-lta-to-pseu} for $P_n X$ and $P_n \Dn X$.
\mk

iv) Let $\uk \rw \us$ be a morphism in $\dop{n-1}$ and suppose that $k_j=0$ for some $1\leq j\leq n-1$ and $s_i\neq 0$ for all $1 \leq i \leq n-1$. Since, by Lemma \ref{lem-fta-tam-2} b) the map $(\Dn X)_{\uk}\rw (\Dn X)_{\us}=X_{\us}$ factors as
\begin{equation*}
  (\Dn X)_{\uk}\rw X_{\uk}\rw X_{\us}
\end{equation*}
by the definitions it follows that the map of pullbacks
\begin{equation*}
\begin{split}
    & (P_n X)^d_{\uk} = (P_n \Dn X)_{\uk} = d(Q_{n-1}\q{n}X)_{\uk} \tiund{dq(\Dn X)_{\uk}} (\Dn X)_{\uk} \\
    & \hspace{6.9cm}\downarrow \\
    & (P_n \Dn X)_{\us}= (P_n X)_{\us} = d(Q_{n-1}\q{n}X)_{\us} \tiund{dq X_{\us}} X_{\us}
\end{split}
\end{equation*}
factors through
\begin{equation*}
  (P_n X)_{\uk}= d(Q_{n-1}\q{n}X)_{\uk} \tiund{dq X_{\uk}} X_{\uk}
\end{equation*}

 Thus all the hypotheses of Lemma \ref{lem-from-lta-to-pseu} are satisfied and we conclude that
\begin{equation*}
    Tr_n P_n \Dn X= Tr_n P_n X
\end{equation*}
which implies
\begin{equation*}
    \Qn \Dn X = \St Tr_n P_n \Dn X= \St Tr_n P_n X= \Qn X\;.
\end{equation*}
\end{proof}
\begin{lemma}\label{lem-fta-tam-3}
  Let $f:Z\rw X$ be a map in $\tawg{n}$ with $Z\in\ftawg{n}$ and $X\in\ta{n}$. Then
  \begin{itemize}
    \item [a)] There is a map in $\ta{n}$  $g:\Dn Z\rw X$, natural in $Z\rw X$.\bk

    \item [b)] If $f$ is a $n$-equivalence, such is $\Dn Z\rw X$.\bk
  \end{itemize}
\end{lemma}
\begin{proof}
Denote by $\pt_i,\zs_i$ the face and degeneracy operators of $Z$, $\pt'_i,\zs'_i$ those of $X$. Let $\zg': Z^d_0\rw Z_0$ the functorial section to the discretization map $\zg: Z_0\rw Z_0^d$. Denote $f_0^d:Z_0^d\rw X_0^d=X_0$. Then
\begin{equation}\label{eq-lem-fta-tam-3-1}
\begin{split}
    & f_0^d\zg=\Id f_0=f_0 \\
    & f_0^d=f_0^d \zg\zg'=f_0\zg'\;.
\end{split}
\end{equation}
This implies
\begin{equation}\label{eq-lem-fta-tam-3-2}
\begin{split}
    & f_0^d(\zg\pt_i)= f_0\pt_i=\pt'_i f_1 \\
    & f_1(\zs_0\zg')=\zs'_0 f_0 \zg'=\zs'_0 f_0^d\;.
\end{split}
\end{equation}
We prove the lemma by induction on $n$. When $n=2$, $\D{2}=\rz$ and we define
\begin{equation*}
  g_k:(\rz Z)_k\rw X_k
\end{equation*}
to be $f_0^d$ when $k=0$ and $f_k$ when $k>0$. From \eqref{eq-lem-fta-tam-3-2}, this is a simplicial map $g:\rz Z\rw X$.

Suppose $f:Z\rw X$ is a 2-equivalence. Then, for each $a,b\in Z_0^d$, there are equivalences of categories
\begin{equation*}
\begin{split}
    & (\rz Z)(a,b)=Z(a,b) \simeq X(fa,fb) \\
    & \p{2}\rz Z =\p{2}Z \simeq \p{2}X
\end{split}
\end{equation*}
so that $g$ is also a 2-equivalence. Suppose, inductively, that the lemma holds for $n-1$.

\bk
a) By \eqref{eq-lem-fta-tam-3-2} there is a map
\begin{equation*}
  h:\rz Z\rw X
\end{equation*}
given by $h_0=f_0^d$, $h_k=f_k$ when $k>0$. By induction hypothesis, we have maps for each $k>0$
\begin{equation*}
\begin{split}
    & (\Dn Z)_k=(\ovl{\Dnm}\rz Z)_k=\Dnm Z_k \xrw{v_k} X_k \\
    & (\Dn Z)_0=Z_0^k=(\ovl{\Dnm}\rz Z)_0 \xrw{f^d_0} X_0\;.
\end{split}
\end{equation*}
Therefore, since $h:\rz Z \rw X$ is a simplicial map and $v_k$ are natural in $Z_k \rw X_k$, we obtain a map
\begin{equation*}
  g:\Dn Z \rw X
\end{equation*}
given by $g_k=v_k$ for $k>0$, $g_0=f^d_0$.\bk

b) If $f$ is an $n$-equivalence, for all $a,b \in Z_0^d$, $Z(a,b)\rw X(fa,fb)$ is a $(n-1)$-equivalence, thus by inductive hypothesis and by Proposition \ref{pro-fta-tam-1} c), such is
\begin{equation*}
  (\Dn Z)(a,b) =\Dnm Z(a,b) \rw X(fa,fb)\;.
\end{equation*}
Since $f$ is an $n$-equivalence, $\p{n}Z\rw \p{n}X$ is a $(n-1)$-equivalence, so using Proposition \ref{pro-fta-tam-1} b) and the inductive hypothesis we obtain a $(n-1)$-equivalence
\begin{equation*}
  \p{n}\Dn Z=\Dnm \p{n}Z \rw \p{n}X\;.
\end{equation*}
By definition we conclude that $\Dn Z\rw X$ is an $n$-equivalence.

\end{proof}
We now establish the main properties of the discretization functor $\Discn:\catwg{n}\rw \ta{n}$ of Definition \ref{def-disc-func}. The proof of this result relies on the properties of $D_n$ established in this chapter, as well as on the properties of the functor $G_n$ studied in Chapter \ref{chap8} and of the functor $Q_n$ studied in Chapter \ref{chap7}.
\begin{theorem}\label{the-disc-func}
Let $\Discn:\catwg{n}\rw \ta{n}$ be as in Definition \ref{def-disc-func}. Then

  \begin{itemize}
    \item [a)] $\Discn$ is identity on discrete objects and commutes with pullbacks over discrete objects.\mk

    \item [b)] For each $X\in\catwg{n}$ there is a zig-zag of $n$-equivalences in $\tawg{n}$ between $X$ and $\Discn X$.\mk

    \item [c)] $\Discn$ preserves $n$-equivalences.
  \end{itemize}
\end{theorem}
\begin{proof}
\

a) This follows from the fact that the same is true for $G_n$ and $\Dn$ (see Theorem \ref{pro-fta-1} and Proposition \ref{pro-fta-tam-1}).\bk

b) Let $X\in\catwg{n}$, then by Proposition \ref{pro-fta-tam-2}
\begin{equation*}
    \Qn\Discn X=\Qn \Dn G_n X = \Qn G_n X\;.
\end{equation*}
Hence by Theorem \ref{the-funct-Qn} there are $n$-equivalences in $\tawg{n}$
\begin{equation*}
  \Discn \xlw{s_n(\Discn X)} Q_n\Discn X=Q_n G_n X \xrw{s_n (G_n X)} G_n X.
\end{equation*}
On the other hand by Theorem \ref{pro-fta-1} there is an $n$-equivalence
 \begin{equation*}
 G_n X\xrw{g_n(X)}X.
 \end{equation*}
  So by composition we obtain a zig-zag of $n$-equivalences
\begin{equation*}
    \Discn X \lw \Qn \Discn X \rw X
\end{equation*}
as required.
\bk

c) This follows from the fact that, by Theorem \ref{pro-fta-1} and Proposition \ref{pro-fta-tam-1}, the same is true for $G_n$ and $\Dn$.
\end{proof}
We now prove our main comparison result between weakly globular \nfol categories and Tamsamani $n$-categories.

\begin{theorem}\label{cor-the-disc-func}
    The functors
    \begin{equation*}
        \Qn:\ta{n}\leftrightarrows \catwg{n}:\Discn
    \end{equation*}
    induce an equivalence of categories after localization with respect to the $n$-equivalences
    \begin{equation*}
        \ta{n}\bsim^n\;\simeq \; \catwg{n}\bsim^n
    \end{equation*}
\end{theorem}
\begin{proof}
Let $X\in\catwg{n}$. By Theorem \ref{the-funct-Qn} and Theorem \ref{pro-fta-1} there are $n$-equivalences
\begin{equation*}
    \Qn\Discn X=\Qn \Dn G_n X\cong\Qn G_n X\rw G_n X\rw X\;.
\end{equation*}
So there is an $n$-equivalence in $\catwg{n}$
\begin{equation*}
    \beta_X:\Qn \Discn X\rw X\;.
\end{equation*}
It follows that $\Qn \Discn X\cong X$ in $\catwg{n}\bsim^n$.

Let $Y\in\ta{n}$. By Theorem \ref{the-funct-Qn} and the above there are $n$-equivalences in $\tawg{n}$
\begin{equation*}
    \Discn \Qn Y\xlw{s_n(\Discn Q_n Y)}\Qn \Discn \Qn Y \xrw{\beta_{Q_n Y}} \Qn Y \xrw{s_{n}(Y)} Y\;.
\end{equation*}
Composing this with the $n$-equivalence
\begin{equation*}
    Z=G_n \Qn \Discn \Qn Y\xrw{g_n(Q_n \Discn Q_n Y)} \Qn \Discn \Qn Y
\end{equation*}
we obtain $n$-equivalences in $\tawg{n}$
\begin{equation}\label{eq-zigcomp-theor}
\Discn \Qn Y\xlw{a} Z \xrw{b} Y\;.
\end{equation}
where
\begin{equation*}
  a= s_n(\Discn Q_n Y) g_n(Q_n \Discn Q_n Y), \quad b= s_{n}(Y)\beta_{Q_n Y} g_n(Q_n \Discn Q_n Y).
\end{equation*}
Since $Z\in\ftawg{n}$, $\Discn \Qn Y\in\ta{n}$ and $Y\in\ta{n}$, applying Lemma \ref{lem-fta-tam-3} to $a$ and $b$ in \eqref{eq-zigcomp-theor} we obtain a zig-zag of $n$-equivalences in $\ta{n}$
\begin{equation*}
    \Discn \Qn Y \lw \Dn Z \rw Y\;.
\end{equation*}
It follows that $ \Discn \Qn Y \cong Y$ in $\ta{n}\bsim^n$.

\end{proof}

\begin{remark}\label{rem-final}
  From Corollary \ref{cor-the-funct-Qn} and Theorem \ref{cor-the-disc-func} we have equivalences of categories
  \begin{equation}\label{eq-final}
        \tawg{n}\bsim^n\; \simeq\; \catwg{n}\bsim^n\; \simeq\;  \ta{n}\bsim^n\;
    \end{equation}
  This means the three Segal-type models are equivalent after localization by the $n$-equivalences. Since both $\ta{n}$ and $\catwg{n}$ are embedded in $\tawg{n}$ this can be interpreted as a kind of partial strictification result for the larger model $\tawg{n}$. Namely, in $\tawg{n}$ the weakening occurs in two ways: with the weakening of the Segal maps and with the weak globularity condition. The equivalences of categories \eqref{eq-final} shows that only one of these two is necessary to obtain a model of weak $n$-categories: the weak globularity condition only, giving rise to the model $\catwg{n}$ or the Segal maps condition only, giving rise to the model $\tawg{n}$.

\end{remark}

%%
%%%%%%%%%%%%%%%%%%%%%%%%%%%%%%%%%%%%%%%%%%%%%%%%%%%%%%%%%%%%%%%%%%%%%%%
\section{Groupoidal weakly globular $\pmb{n}$-fold categories}\label{sec-group-wg-nfol-cat}

In this section we introduce the subcategory $\gcatwg{n}\subset\catwg{n}$ of groupoidal weakly globular \nfol categories and we show that it is an algebraic model of $n$-types. \index{Algebraic model of $n$-types}This means that weakly globular $n$-fold categories satisfy the homotopy hypothesis.
\begin{definition}\label{def-gta-1}
    The full subcategory $\gtawg{n}\subset\tawg{n}$ of groupoidal weakly globular Tamsamani $n$-categories is defined inductively as follows. \index{Groupoidal!- weakly globular Tamsamani $n$-categories}

    For $n=1$, $\gtawg{1}=\Gpd$. Note that $\cathd{}\subset \gtawg{1}$. Suppose inductively we defined $\gtawg{n-1}\subset\tawg{n-1}$. We define $X\in\gtawg{n}\subset\tawg{n}$ such that
    \begin{itemize}
      \item [i)] $X_k\in \gtawg{n-1}$ for all $k\geq 0$.\mk

      \item [ii)] $\p{n}X\in \gtawg{n-1}$.
    \end{itemize}
\end{definition}
\begin{remark}\label{rem-def-gwgn}
  If $X\in\gtawg{n}$, then by definition $X_1\in \gtawg{n-1}$, thus for each $a,b \in X_0^d$, $X(a,b)\in \gtawg{n-1}$. Also note that $\gtawg{n}$ is closed under products, as easily seen by induction.
\end{remark}
\begin{lemma}\label{lem-gta-1}
    Let $f:X\rw Y$ be an $n$-equivalence in $\tawg{n}$
    \begin{itemize}
      \item [i)] If $Y\in\gtawg{n}$ then $X\in\gtawg{n}$.\mk

      \item [ii)] If $X\in\gtawg{n}$ then $Y\in\gtawg{n}$.
    \end{itemize}
\end{lemma}
\begin{proof}
By induction on $n$. The case $n=1$ holds since a category equivalent to a groupoid is itself a groupoid. Suppose, inductively, that the lemma holds for $n-1$ and let $f:X\rw Y$ be an $n$-equivalence.

i) For each $a,b\in X_0^d$ the map
\begin{equation*}
    f(a,b):X(a,b)\rw Y(fa,fb)
\end{equation*}
is a $(n-1)$-equivalence in $\tawg{n-1}$ with $Y(fa,fb)\in\gtawg{n-1}$ (see Remark \ref{rem-def-gwgn}). So by induction hypothesis $X(a,b)\in \gtawg{n-1}$. Since
\begin{equation*}
    X_1=\uset{a,b\in X_0^d}{\coprod}X(a,b)
\end{equation*}
it follows that $X_1\in\gtawg{n-1}$. We also have
\begin{equation*}
    \tens{X_1}{X_0^d}=\uset{a,b,c\in X_0^d}{\coprod}X(a,b)\times X(b,c)
\end{equation*}
and, by Remark \ref{rem-def-gwgn}, $X(a,b)\times X(b,c)\in \gtawg{n-1}$. Thus $\tens{X_1}{X_0^d}\in\gtawg{n-1}$. Similarly one can show that
 \begin{equation*}
   \pro{X_1}{X_0^d}{k}\in\gtawg{n-1}
 \end{equation*}
 for all $k\geq 2$. On the other hand, there is an $(n-1)$-equivalence in $\tawg{n-1}$
  \begin{equation*}
    X_k \rw \pro{X_1}{X_0^d}{k}
  \end{equation*}
  Thus by inductive hypothesis we conclude that $X_k\in\gtawg{n-1}$ for all $k\geq 0$.

By definition there is a $(n-1)$-equivalence
\begin{equation*}
    \p{n}f:\p{n}X\rw\p{n}Y
\end{equation*}
with $\p{n}Y\in\gtawg{n-1}$ since by hypothesis $Y\in\gtawg{n}$. Hence by inductive hypothesis $\p{n}X\in\gtawg{n-1}$. We conclude that $X\in\gtawg{n-1}$.
\medskip

ii) The proof is completely similar to one of i).

\end{proof}
\begin{remark}\label{rem-gta-1}
    It follows immediately from the definition of $\gtawg{n}$ that the embedding
    \begin{equation*}
        J_n:\tawg{n}\hookrightarrow \funcat{n-1}{\Cat}
    \end{equation*}
    restricts to the embedding
    \begin{equation*}
        J_n:\gtawg{n}\hookrightarrow\funcat{n-1}{\Gpd}\;.
    \end{equation*}
    Since $p=q:\Gpd\rw\Set$ it follows that for each $X\in\gtawg{n}$ there is a morphism, natural in $X$,
    \begin{equation*}
        X\rw\di{n}\p{n}X\;.
    \end{equation*}
\end{remark}
\begin{definition}\label{def-gta-2}
    The category $\gcatwg{n}\subset\catwg{n}$ of groupoidal weakly globular \nfol categories is the full subcategory of $\catwg{n}$ whose objects $X$ are in $\gtawg{n}$.
    \index{Groupoidal!- weakly globular $n$-fold categories} \index{Groupoidal!- Tamsamani $n$-categories}

    The category $\gta{n}\subset \ta{n}$ of groupoidal Tamsamani $n$-categories is the full subcategory of $\ta{n}$ whose objects $X$ are in $\gtawg{n}$.
\end{definition}
\begin{remark}\label{rem-gta-2}
    The following facts are immediate from the definitions:
    \begin{itemize}
      \item [a)] $X\in\gcatwg{n}$ (resp. $X\in\gta{n}$) if and only if for each $k\geq 0$ $X_k\in\gcatwg{n-1}$ (resp. $X_k\in\gta{n-1}$) and $\p{n}X\in\gcatwg{n}$ (resp. $\p{n}X\in\gta{n-1}$).\mk

      \item [b)] Let $f:X\rw Y$ be an $n$-equivalence in $\tawg{n}$ and suppose that $Y\in\gtawg{n}$. Then if $X\in\catwg{n}$ it is $X\in\gcatwg{n}$  and  if $X\in\ta{n}$ then $X\in\gta{n}$.
    Similarly if $f$ is an $n$-equivalence in $\tawg{n}$ and $X\in\gtawg{n}$.

    \end{itemize}
\end{remark}
\begin{corollary}\label{cor-gta-1}
The following facts hold:
    \begin{itemize}
      \item [a)] The functor
      \begin{equation*}
        \Qn:\tawg{n}\rw\catwg{n}
      \end{equation*}
      restricts to a functor
      \begin{equation*}
        \Qn:\gtawg{n}\rw\gcatwg{n}
      \end{equation*}
      such that for each $X\in\gtawg{n}$ there is a $n$-equivalence in $\gtawg{n}$ $s_n(X):\Qn X\rw X$.\mk

      \item [b)] The functor
      \begin{equation*}
        \Discn:\catwg{n}\rw\ta{n}
      \end{equation*}
      restricts to a functor
      \begin{equation*}
        \Discn:\gcatwg{n}\rw\gta{n}
      \end{equation*}
      such that for each $X\in\gcatwg{n}$ there is a zig-zag of $n$-equivalences in $\gtawg{n}$ between $X$ and $\Discn X$.
    \end{itemize}
\end{corollary}

\begin{proof}
By Theorem \ref{the-funct-Qn} there is an $n$-equivalence in $\tawg{n}$
 \begin{equation*}
   s_n(X):\Qn X\rw X.
 \end{equation*}
 Since $X\in\gtawg{n}$ and $Q_n X\in \catwg{n}$, by Lemma \ref{lem-gta-1} and Remark \ref{rem-gta-2} $Q_n X\in \gcatwg{n}$. By Theorem \ref{the-disc-func} there is a zig-zag of $n$-equivalences between $X$ and $\Discn X$. Since $X\in\gtawg{n}$, by Lemma \ref{lem-gta-1} this is a zig-zag of $n$-equivalences in $\gtawg{n}$, and since $X\in \ta{n}$, by Remark \ref{rem-gta-2} $\Discn X\in \gta{n}$.
\end{proof}

In the next Proposition we specialize the comparison result of Theorem \ref{cor-the-disc-func} to the higher groupoidal setting.

\begin{proposition}\label{pro-gta-1}
    The functors
    \begin{equation*}
        \Qn:\gta{n}\leftrightarrows \gcatwg{n}:\Discn
    \end{equation*}
    induce an equivalence of categories after localization with respect to the $n$-equivalences
    \begin{equation*}
        \gta{n}\bsim^n \;\simeq\; \gcatwg{n}\bsim^n\;.
    \end{equation*}
\end{proposition}
\begin{proof}
Let $X\in\gcatwg{n}$. As in the proof of Theorem \ref{cor-the-disc-func} there is an $n$-equivalence in $\catwg{n}$
\begin{equation*}
    \Qn\Discn X \rw X\;.
\end{equation*}
Since $X\in\catwg{n}$, by Remark \ref{rem-gta-2}, $\Qn \Discn X\in\gcatwg{n}$, so this is an $n$-equivalence in $\gcatwg{n}$. It follows that there is an isomorphism in $\gcatwg{n}\bsim^n$
\begin{equation*}
    \Qn\Discn X\cong X.
\end{equation*}

Let $Y\in\gta{n}$. By the proof of Theorem \ref{cor-the-disc-func} there is a zig-zag of $n$-equivalences in $\ta{n}$
\begin{equation*}
    \Discn \Qn Y\lw  \Dn Z \rw Y\;.
\end{equation*}
Since $Y\in\gta{n}$, by Remark \ref{rem-gta-2} this is a zig-zag of $n$-equivalences in $\gta{n}$. It follows that there is an isomorphism in $\gta{n}\bsim^n$
\begin{equation*}
    \Discn \Qn Y\cong Y
\end{equation*}
\end{proof}
As a consequence of the previous proposition and of the results of Tamsamani \cite{Ta}, we obtain that groupoidal weakly globular \nfol categories are an algebraic model of $n$-types. That is our model $\catwg{n}$ of weak $n$-categories satisfies the homotopy hypothesis.\index{Homotopy hypothesis}

In what follows let
\begin{equation*}
\xymatrix{
\gta{n} \ar@/_/[rr]_{B} && \nty \ar@/_/[ll]_{\clT_{n}}
}
\end{equation*}
be the fundamental Tamsamani $n$-groupoid functor $\clT_n$ \index{Fundamental Tamsamani $n$-groupoid} and the classifying space functor $B$ as in \cite{S2}.\index{Classifying space}
\begin{theorem}\label{cor-gta-2}
The functors
\begin{equation*}
\xymatrix{
\gcatwg{n} \ar@/_/[rr]_{B {\cirsm }\Discn} && \nty \ar@/_/[ll]_{\Qn {\cirsm} \clT_{n}}
}
\end{equation*}
induce an equivalence of categories
\begin{equation*}
    \gcatwg{n}\bsim^n\;\simeq\;\Ho(\nty)\;.
\end{equation*}
\end{theorem}
\begin{proof}
By \cite{S2} the functors $\clT_n$ and $B$ induce an equivalence of categories
\begin{equation}\label{eq1-cor-gta-2}
    \gta{n}\bsim^n\;\simeq\;\Ho(\nty)\;.
\end{equation}
while by Proposition \ref{pro-gta-1} the functors $\Qn$ and $\Discn$ induce an equivalence of categories
\begin{equation}\label{eq2-cor-gta-2}
  \gcatwg{n}\bsim^n\;\simeq\;\gta{n}\bsim^n\;.
\end{equation}
By \eqref{eq1-cor-gta-2} and \eqref{eq2-cor-gta-2} the result follows.
\end{proof}
\begin{remark}\label{rem-cor-gta-2}\index{Geometric weak equivalences}
  We call a map $f$ in $\gcatwg{n}$ a geometric weak equivalence if $(B {\cirsm} \Discn)(f)$ is a weak homotopy equivalence of spaces. We note that a map $f$ in $\gcatwg{n}$ is an $n$-equivalence if and only if it is a geometric weak equivalence. In fact, if $f$ is an $n$-equivalence, it is an isomorphism in $\gcatwg{n}\bsim^n$ so by Theorem \ref{cor-gta-2}, $(B {\cirsm} \Discn)(f)$ is an isomorphism in $\Ho(\nty)$, thus it is a weak homotopy equivalence in $n$-types.

  Conversely, if $f$ is a geometric weak equivalence, $(B {\cirsm} \Discn)(f)$ is an isomorphism in $\Ho(\nty)$, so by Theorem  \ref{cor-gta-2} (since equivalence of categories reflect isomorphisms), $f$ is an isomorphism in $\gcatwg{n}\bsim^n$, so $f$ is an $n$-equivalence in $\gcatwg{n}$.
\end{remark}

We finally observe that, as an immediate consequence of our results, all the three Segal-type models of this work are a model of weak $n$-category satisfying the homotopy hypothesis. In what follows $\gseg{n}$ denotes any of the three groupoidal Segal-type models $\gcatwg{n}$, $\gta{n}$, $\gtawg{n}$.

\begin{corollary}\label{cor-three-models}\index{Homotopy hypothesis}
  Each of the three Segal-type models $\seg{n}$ is a model of weak $n$-categories satisfying the homotopy hypothesis, that is there is an equivalence of categories
 \begin{equation*}
  \gseg{n}\bsim^n\;\simeq\;\Ho(\nty)\;.
\end{equation*}
\end{corollary}

\begin{proof}
In the case $\seg{n}=\ta{n}$ this is the result of \cite{Ta}. When $\seg{n}=\catwg{n}$ this is the content of Theorems \ref{cor-the-disc-func} and \ref{cor-gta-2}. In the case $\seg{n}=\tawg{n}$ by Corollary \ref{cor-the-funct-Qn} there is an equivalence of categories
\begin{equation}\label{eq-cor-the-funct-Qn}
        \tawg{n}\bsim^n\; \simeq\; \catwg{n}\bsim^n
    \end{equation}
    Thus by Remark \ref{rem-gta-2} this restricts to an equivalence of categories
 \begin{equation*}
        \gtawg{n}\bsim^n\; \simeq\; \gcatwg{n}\bsim^n
    \end{equation*}
    so by Theorem \ref{cor-gta-2} we conclude that there is an equivalence of categories
    \begin{equation*}
    \gtawg{n}\bsim^n\;\simeq\;\Ho(\nty)\;.
\end{equation*}
\end{proof}

%%%%%%%%%%%%%%%%%%%%%%%%%%%%%%%%%%%%%%%%%%%%%%%%%%%%

\section{A convenient model for the fundamental groupoidal weakly globular $\pmb{n}$-fold category functor}\label{mod-fund-wg-group}
Theorem \ref{cor-gta-2} exhibits the fundamental groupoidal weakly globular $n$-fold category functor \index{Fundamental groupoidal weakly globular $n$-fold category functor}
 \begin{equation}\label{fund-wg}
  \clG_n:\nty\rw \gcatwg{n}
 \end{equation}
 as the composite $\Qn {\cirsm} \clT_n$ where $\clT_n$ is the Tamsamani $n$-groupoid functor from \cite{S2} and $\Qn$ is the rigidification functor. Using the results of Blanc and the author \cite{BP} we exhibit an alternative functor
\begin{equation*}
  j \clH_n:\nty\rw \gcatwg{n}
\end{equation*}
which is simpler than $\clG_n$ and whose definition is independent on \cite{S2}.

Using our previous results, we show in Corollary \ref{cor1-mod-fund-wg-group} that $j \clH_n$ and the classifying space functor \index{Classifying space} $\btil:\gcatwg{n}\rw \nty$ induce equivalence of categories
\begin{equation*}
  \gcatwg{n}\bsim^n\;\simeq\;\Ho(\nty)\;.
\end{equation*}
Thus $j \clH_n$ can be used as a fundamental functor instead of $\clG_n$.

\subsection{The functor \pmb{$\clH_n$}}
In \cite[Definition 3.19]{BP} Blanc and the author introduced the category $\gpdwg{n}$ of weakly globular $n$-fold groupoids,\index{Weakly globular $n$-fold groupoids} which is a full subcategory of the category $\gpd{n}$ of $n$-fold groupoids. It is immediate from the definitions that $\gpdwg{n}$ is a full subcategory of $\gcatwg{n}$ and that a map in $\gpdwg{n}$ is a $n$-equivalence if and only if it is so in $\gcatwg{n}$.

As for $\gcatwg{n}$ (see Remark \ref{rem-cor-gta-2}), it was shown in \cite{BP} that $n$-equivalences in $\gpdwg{n}$ are the same as geometric weak equivalences. Let
\begin{equation*}
  \clS:\Top \rw\funcat{}{\Set}
\end{equation*}
be the singular functor,\index{Singular functor} whose image consists of fibrant simplicial sets. Let
\begin{equation*}
  \orn{n}:\funcat{}{\Set}\rw\funcat{n}{\Set}
\end{equation*}
be the functor induced by the ordinal sum $or_{n}:\Delta^n\rw\Delta$ \index{Ordinal sum}. Thus
\begin{equation*}
  (\orn{n}X)_{p_1\ldots p_n}=X_{n-1+p_1+\cdots +p_n}\;.
\end{equation*}
 \begin{remark}\label{rem-or}
 The functor $\orn{n}$ produces an $n$-fold simplicial resolution of a simplicial set, since it can be shown \cite[Lemma 2.13]{BP} that for any simplicial set $X$, there is a natural weak equivalence
\begin{equation*}
\zve\lo{n}:Diag_n \orn{n} X \rw X
\end{equation*}
where $Diag_n:\funcat{n}{\Set}\rw\funcat{}{\Set}$ is the multi-diagonal functor,\index{Multi-diagonal} given by
\begin{equation*}
  (Diag_n X)_m = X_{m,m,\ldots,m}\;.
\end{equation*}
We also showed in \cite[Section 2.9]{BP} that
\begin{equation}\label{eq6-mod-fund-wg-group}
  \orn{n}Y = \overline{\ord}\lo{n-1}\up{2}\orn{2} Y
\end{equation}
and we proved in \cite[Lemma 2.28]{BP} that if $Y$ is a Kan complex, for each $n\geq 2$
\begin{equation}\label{eq7-mod-fund-wg-group}
  \overline{\ord}\lo{n-1}\up{2} \Nu{2} \hat\pi_{1}\up{2} \orn{2} X \cong \Nu{n} \hat\pi_{1}\up{n} \orn{n} X\;.
\end{equation}
 \end{remark}
Let
\begin{equation*}
  \clP_n:\funcat{n}{\Set}\rw\gpd{n}
\end{equation*}
be the left adjoint to the $n$-fold nerve \index{Left adjoint to the $n$-fold nerve}
\begin{equation*}
  \N{n}:\gpd{n}\rw\funcat{n}{\Set}\;.
\end{equation*}
\begin{definition}{\rm \cite[Definition 2.30]{BP}}\label{def-hn}
The fundamental weakly globular $n$-fold groupoid \index{fundamental weakly globular $n$-fold groupoid} functor is given by the composite
  \begin{equation}
 \clH_n: \nty\xrw{\clS} \funcat{}{\Set} \xrw{\orn{n}} \funcat{n}{\Set} \xrw{\clP_n} \gpd{n}
\end{equation}
\end{definition}

For a general $n$-fold simplicial set $Y$, $\clP_n Y$ does not have a simple and explicitly computable expression. However, we showed that, given a space $X$, the fibrancy of $\clS X$ induces a property of $\orn{n}\clS X$ which we called in \cite{BP} $(n,2)$-fibrancy (see \cite[Definition 2.3.1 and Proposition 2.3.9]{BP}). We then showed that to apply $\clP_n$ to a $(n,2)$-fibrant $n$-fold simplicial set we need only apply the usual fundamental groupoid \index{Fundamental groupoid}in each of the $(n-1)$-simplicial directions. Thus we have
\begin{theorem}{\rm \cite[Theorem 2.40]{BP}}\label{the-2.40-BP}
Let $\clH_n$ be as in Definition \ref{def-hn} and $X$ be a space. Then
\begin{equation}\label{eq1-mod-fund-wg-group}
  \clH_n X=\hat{\pi}\up{1} \hat{\pi}\up{2} \cdots \hat{\pi}\up{n} \orn{n} \clS X\;.
\end{equation}

\end{theorem}

Using this explicit description of $\clH_n$ we showed in \cite{BP} that the functor $\clH_n$ in fact lands in $\gpdwg{n}$. Further, we proved
\begin{theorem}{\rm \cite[Theorem 4.32]{BP}}\label{the-4.42-BP}
  The functors
\begin{equation*}
  \clH_n:\nty\rw \gpdwg{n}, \qquad B:\gpdwg{n}\rw \nty
\end{equation*}
induce functors
  \begin{equation*}
    \Ho(\nty) \leftrightarrows \gpdwg{n}\bsim^n
  \end{equation*}
  with $B\clH_n\cong \Id$.
\end{theorem}
\begin{remark}\label{rem1-mod-fund-wg-group}
  Let
  \begin{equation*}
    \btil = B{\cirsm}\Discn:\gcatwg{n} \rw \nty.
  \end{equation*}
   From \cite{BP}, the following diagram commutes
\begin{equation*}
\xymatrix{
\gcatwg{n}\bsim^n \ar^{\btil}[rr] &&  \Ho(\nty)\\
\gpdwg{n}\bsim^n\ar@{^{(}->}_{j}[u] \ar_{B}[urr]
}
\end{equation*}
\end{remark}

Using our previous results we deduce the following corollary, which shows that $j \clH_n$ can be used as an alternative fundamental functor from $n$-types to $\gcatwg{n}$.
\begin{corollary}\label{cor1-mod-fund-wg-group}
  Let $j\clH_n$ be the composite
\begin{equation*}
 j\clH_n: \nty\xrw{\clS} \funcat{}{\Set} \xrw{\orn{n}} \funcat{n}{\Set} \xrw{\clP_n}\gpdwg{n} \overset{j}{\hookrightarrow} \gcatwg{n}
\end{equation*}
and let $\btil: \gcatwg{n} \rw \nty$ be as in Remak \ref{rem1-mod-fund-wg-group}. Then $j\clH_n$ and $\btil$ induce an equivalence of categories
\begin{equation}\label{eq2-mod-fund-wg-group}
  \gcatwg{n}\bsim^n\;\simeq\;\Ho(\nty)\;.
\end{equation}
\begin{proof}
Let $X\in\gcatwg{n}$; by Theorem \ref{cor-gta-2}, $\btil X$ is an $n$-type and by Theorem \ref{the-4.42-BP} and Remark \ref{rem1-mod-fund-wg-group}
\begin{equation}\label{eq3-mod-fund-wg-group}
 \btil j \clH_n \btil X \cong B\clH_n \btil X \cong \btil X
\end{equation}
in $\Ho(\nty)$. Since both $X$ and $j \clH_n \btil X$ are in $\gcatwg{n}$, the equivalence of categories
\begin{equation*}
  \gcatwg{n}\bsim^n\;\cong\; \Ho(\nty)
\end{equation*}
of Theorem  \ref{cor-gta-2} induced by $\Qn\clT_n$ and $\btil$, together with \eqref{eq3-mod-fund-wg-group} imply that
\begin{equation}\label{eq4-mod-fund-wg-group}
  j\clH_n \btil X \cong X
\end{equation}
in $\gcatwg{n}\bsim^n$.

 Let $Y\in\Ho(\nty)$. By Theorem \ref{the-4.42-BP} and Remark \ref{rem1-mod-fund-wg-group}
\begin{equation}\label{eq5-mod-fund-wg-group}
  \btil j \clH_n Y \cong B \clH_n Y \cong Y
\end{equation}
in $\Ho(\nty)$. By \eqref{eq4-mod-fund-wg-group} and \eqref{eq5-mod-fund-wg-group} we conclude that $j \clH_n$ and $\btil$ induce the equivalence of categories \eqref{eq2-mod-fund-wg-group}.
\end{proof}

\end{corollary}
The following corollary shows that the functor
\begin{equation*}
  \p{n}:\gcatwg{n}\rw\gcatwg{n-1}
\end{equation*}
is the algebraic version of the Postnikov truncation functor \index{Postnikov truncation}
\begin{equation*}
  \mbox{$n$-types$\rw (n-1)$-types}\;.
\end{equation*}
%%%
\begin{corollary}\label{cor2-mod-fund-wg-group}
Let $X\in\gcatwg{n}$. The map
\begin{equation*}
    X\rw\di{n}\q{n}X=\di{n}\p{n}X
 \end{equation*}
induces a map of spaces
\begin{equation*}
  \btil X\rw\btil \p{n}X
\end{equation*}
such that, for each $0\leq i \leq n-1$, $x\in\btil X$
\begin{equation*}
  \pi_i(\btil X,x)\cong\pi_i(\btil \p{n}X,x)\;.
\end{equation*}
\end{corollary}
\begin{proof}
It is shown in \cite{BP} that the functor
\begin{equation*}
  \p{n}:\gpdwg{n}\rw\gpdwg{n-1}
\end{equation*}
(denoted $\Pi_0\up{n}$ in \cite{BP}) is such that, for each $Y\in\gpdwg{n}$, the map $Y\rw\di{n}\p{n}Y$ induces a map of spaces
\begin{equation*}
  BY\rw B\p{n}Y
\end{equation*}
such that, for each $0\leq i \leq n-1$, $y\in BY$
\begin{equation}\label{eq1-cor2-mod-fund-wg-group}
  \pi_i(BY,y)\cong \pi_i(B\p{n}Y,y)
\end{equation}
Let $X\in\gcatwg{n}$. By Corollary \ref{cor2-mod-fund-wg-group} there is a zigzag of $n$-equivalences between $X$ and $j \clH_n\btil X$, and thus also a zig-zag of $(n-1)$-equivalences between $\p{n}X$ and $\p{n}j \clH_n\btil X=j\p{n}\clH_n\btil X$. This implies that there are zig-zag of weak homotopy equivalences between $\btil X$ and $\btil j\clH_n \btil X=B\clH_n \btil X$ as well as between $\btil \p{n} X$ and $\btil j\p{n}\clH_n\btil X=B\p{n}\clH_n\btil X$. Therefore, for all $x\in\btil X$
\begin{equation}\label{eq2-cor2-mod-fund-wg-group}
\begin{split}
    & \pi_i(\btil X,x)\cong \pi_i(B\clH_n\btil X,x') \\
    & \pi_i(\btil\p{n} X,x)\cong \pi_i(B\p{n}\clH_n\btil X,x')\;.
\end{split}
\end{equation}
By \eqref{eq1-cor2-mod-fund-wg-group}, taking $Y=\clH_n\btil X\in\gpdwg{n}$ we obtain for each $0\leq i \leq n-1$
\begin{equation*}
  \pi_i(B\clH_n\btil X,x')=\pi_i(B\p{n}\clH_n\btil X,x')\;.
\end{equation*}
By \eqref{eq2-cor2-mod-fund-wg-group} this implies, for all $0\leq i \leq n-1$
\begin{equation*}
  \pi_i(\btil X,x)\cong\pi_i(\btil\p{n} X,x)\;.
\end{equation*}
\end{proof}

\bk

\subsection{Some examples}\label{sub-last-examples}

 We now illustrate the fundamental weakly globular $n$-fold groupoid of a space in some low dimensional cases. This shows how explicit and convenient is the use of the functor $j \clH_n $ to produce the fundamental groupoidal weakly globular \nfol category of a space.

For each $n>1$ denote by $\Rbt{n}$ the composite
\begin{equation}\label{eq-rn}
  \funcat{}{\Set} \xrw{\orn{n}} \funcat{n}{\Set} \xrw{\clP_n} \gpdwg{n}
\end{equation}
and let $\Rbt{1}=\hat\pi_1:\funcat{}{\Set}\rw\Gpd$ be the fundamental groupoid, so that, using our previous notation
\begin{equation*}
  j\clH_n = j \Rbt{n}\clS : \nty \rw \gcatwg{n}\;.
\end{equation*}
In \cite[Lemma 4.14]{BP} we give an iterative description of $\Rbt{n}Y$ for a Kan complex $Y$ which is more transparent that the formula \ref{eq1-mod-fund-wg-group}. More precisely, let $u_Y:\Dec Y\rw Y$ be as in Section \ref{decalage} and consider the corresponding internal equivalence relation
\begin{equation*}
  (\Dec Y)[u] \in\Gpd(\funcat{}{\Set})
\end{equation*}
as in Definition \ref{def-int-eq-rel}. Denote by
\begin{equation*}
  \Lbt{\bl} Y\in \funcat{}{\funcat{}{\Set}}
\end{equation*}
the nerve of $(\Dec Y)[u]$, so that
\begin{equation}\label{eq-l}
\Lbt{k} Y=
\left\{
  \begin{array}{ll}
    \Dec Y, & \mbox{if $k=0$} \\
    \pro{\Dec Y}{Y}{k+1}, & \mbox{if $k\geq 1$}
  \end{array}
\right.
\end{equation}
A picture of the corner of $\Lbt{\bl}Y$ is given below.

\begin{figure}[ht]
\begin{center}
$$
\entrymodifiers={++++[]}
\xymatrix@R=15pt@C=10pt{
\cdots~~~
Y_{3}\oset{d_3}\times_{Y_{2}}~Y_{3}\oset{d_3}\times_{Y_{2}}~Y_{3}\;
\ar@<2ex>[rr] \ar[rr]\ar@<-2ex>[rr]
     \ar@<2.5ex>[d] \ar[d]\ar@<-2.5ex>[d] &&
Y_{3}\oset{d_3}\times_{Y_{2}}~Y_{3} \ar@<1ex>[rr]^(0.6){p\sb{2}} \ar@<-1ex>[rr]^(0.6){p\sb{1}}
     \ar@<2.5ex>[d] \ar[d]\ar@<-2.5ex>[d] &&
Y_{3}  \ar@<2.5ex>[d]^{d_{2}} \ar[d]^{d_{1}}\ar@<-2.5ex>[d]^{d_{0}} \\
\cdots~~~
Y_{2}\oset{d_2}\times_{Y_{1}}~Y_{2}\oset{d_2}\times_{Y_{1}}~Y_{2}\;
\ar@<2ex>[rr] \ar[rr] \ar@<-2ex>[rr]
     \ar@<0.5ex>[d] \ar@<-0.5ex>[d] &&
Y_{2}\oset{d_2}\times_{Y_{1}}~Y_{2} \ar@<1ex>[rr]^(0.6){p\sb{2}} \ar@<-1ex>[rr]^(0.6){p\sb{1}}
     \ar@<0.5ex>[d] \ar@<-0.5ex>[d] &&
Y\sb{2} \ar@<0.5ex>[d]^{d_{1}} \ar@<-0.5ex>[d]_{d_{0}} \\
\cdots~~~
Y_{1}\oset{d_1}\times_{Y_{0}}~Y_{1}\oset{d_1}\times_{Y_{0}}~Y_{1}\;
\ar@<2ex>[rr] \ar[rr] \ar@<-2ex>[rr]
    \ar@{}[u] &&
Y_{1}\oset{d_1}\times_{Y_{0}}~Y_{1} \ar@<1ex>[rr]^(0.6){p\sb{2}} \ar@<-1ex>[rr]^(0.6){p\sb{1}} && Y_{1}
}
$$
\end{center}
\caption{Corner of $\Lbt{\bl} Y$}
\label{cornerldot}
\end{figure}

\begin{proposition}{\rm \cite[Lemma 4.14]{BP}}\label{pro-4.14-BP}
\begin{itemize}
  \item [a)] For each $k\geq 0$
\begin{equation}\label{eq7bis-mod-fund-wg-group}
  (\Nu{n}\Rbt{n}Y)_k\cong\Rbt{n-1}\Lbt{n} Y\;.
\end{equation}
Thus, for each $k\geq 1$
\begin{equation}\label{eq8-mod-fund-wg-group}
  \Rbt{n}\Lbt{n}Y\cong\pro{\Rbt{n-1}\Lbt{1}Y}{\Rbt{n-1}\Dec Y}{k}\;.
\end{equation}

  \item [b)] If $Y$ is homotopically trivial for $k\geq 1$
\begin{equation*}
  \Rbt{n}\Lbt{n}Y\cong\pro{\Rbt{n}\Dec Y}{\Rbt{n} Y}{k+1}\;.
\end{equation*}

\end{itemize}

\end{proposition}

We also proved in  \cite[Proposition 4.28]{BP} that for every Kan complex $Y$,
\begin{equation}\label{eq9-mod-fund-wg-group}
  \p{n}\Rbt{n}Y\cong \Rbt{n-1}Y
\end{equation}
and thus, if $X$ is a $n$-type and $Y=\clS X$,
\begin{equation}\label{eq10-mod-fund-wg-group}
  \p{n} j \clH_n X = \p{n} j \Rbt{n} \clS X = j \clH_{n-1}X\;.
\end{equation}

\begin{example}\label{ex1-mod-fund-wg-group}
  \emph{The fundamental weakly globular double groupoid of a space.}\index{Fundamental weakly globular double groupoid of a space}

Let $X$ be a space and $Y=\clS X$ its singular simplicial set. The bisimplicial set $\orn{2}Y$ can be described as follows.

Let $\Dec$ and $\Dec'$ be the two d\'{e}calage comonads as in Section \ref{decalage}. The comonad $\Dec$ yields a simplicial resolution $Z\in\funcat{}{\funcat{}{\Set}}$ for any $Y\in\funcat{}{\Set}$ with
\begin{equation*}
  Z_{k-1} =\Dec^k Y = \underset{k}{\underbrace{\Dec(\Dec\ldots\;\Dec\ldots)}}\in\funcat{}{\Set}
\end{equation*}
It can be shown (see for instance \cite{Illusie1972}) that
\begin{equation*}
  \orn{2}Y=Z\;.
\end{equation*}
The bisimplicial set $\orn{2}Y$ is depicted in Figure \ref{CornerOr2} on page \pageref{corner2page}, viewed as a horizontal simplicial object in $\funcat{}{\Set}$ (the degeneracy maps are not shown). The corresponding resolution using $\Dec'$ is also depicted in Figure \ref{CornerOr2}, viewed as a vertical simplicial object in $\funcat{}{\Set}$.

 From Theorem \ref{the-2.40-BP} it is
\begin{equation}\label{eq11-mod-fund-wg-group}
  j \clH_n X= \clP_2 \orn{2}\clS X =\hat\pi_1\up{1} \hat\pi_1\up{2} \orn{2} Y
\end{equation}
where $\hat\pi_1\up{1}$ and $\hat\pi_1\up{2}$ are the fundamental groupoids in the two simplicial directions. Since $Y=\clS X$ is a Kan complex, such are $\Dec Y$ and $\Dec' Y$, so $(\orn{2}Y)_{k*}$ and $(\orn{2}Y)_{*k}$ are Kan complexes for all $k$, and taking their fundamental groupoids amounts to dividing out the 1-simplices by the relations given by the 2-simplicies. Using the formula \eqref{eq7bis-mod-fund-wg-group} we obtain
\begin{equation*}
  (\Nu{2} j \clH_{2} X)_0 =\hat\pi_1 \Dec Y
\end{equation*}
which is the homotopically discrete groupoid corresponding to the surjective map of sets $d_1: X_1\rw X_0$. By \eqref{eq7bis-mod-fund-wg-group} we also have
\begin{equation*}
  (\Nu{2} j \clH_{2} X)_k =\hat\pi_1 \Lbt{k} Y
\end{equation*}
where $\Lbt{k} Y$ is as in Figure \ref{cornerldot}.

In Figure \ref{CornerDblNerve} on page \pageref{corner2page} we display the corner of the double nerve of $j \clH_2 X$ where $(\tens{Y_2}{Y_1})\bsim$ denotes the result modding out by the relations of the 2-simplices.

Note that $(j \clH_2 X)_0$ is homotopically discrete while $\bar p j \clH_2 X$ is the nerve of the groupoid $\hat\pi_1 Y$. Therefore, by Lemma \ref{lem-crit-doucat-wg}, $j \clH_2 X\in\catwg{2}$ and, given the groupoidal structures, $j \clH_2 X\in\gcatwg{2}$. In fact this is also a double groupoid, so that $\clH_2 X\in\gpdwg{2}$.

\end{example}
%%%
%%
\begin{example}\label{ex2-mod-fund-wg-group}
 \emph{ The fundamental weakly globular 3-fold groupoid of a space.}\index{Fundamental weakly globular 3-fold groupoid of a space}

Let $X$ be a space and $Y=\clS X$ its singular simplicial set. By \eqref{eq6-mod-fund-wg-group}
\begin{equation*}
  \orn{3}Y=\ovl\ord\lo{2}\up{2}\orn{2}Y\;.
\end{equation*}
See Figure \ref{CornerOr3} on page \pageref{corner3page} for a picture of the corner of $\orn{3}Y$. In Figure \ref{CornerOr3nerve} on page \pageref{corner3page} we have a picture of $j \clH_3 X$, where $\clS X=Y$.

The isomorphisms describing $Z_4\up{0,1}$, $Z_4\up{0,2}$, $Z_4\up{1,2}$ in Figure \ref{CornerOr3nerve} are derived from the simplicial identities. Namely, the simplicial identity $d_0 d_0=d_0 d_1$ implies that the limit of the following diagrams are isomorphic:
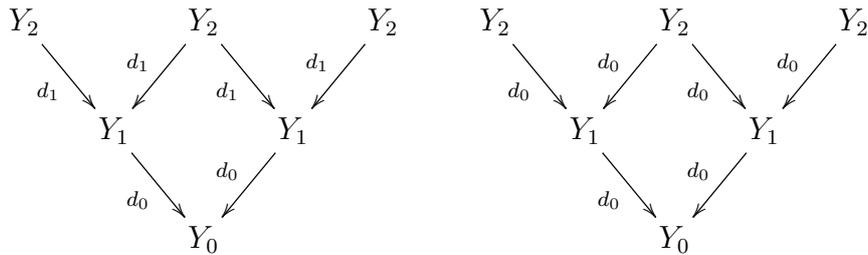
\begin{equation*}
\xymatrix@C=15pt{
Y_2 \ar_{d_1}[dr] && Y_2 \ar_{d_1}[dl] \ar_{d_1}[dr] && Y_2 \ar_{d_1}[dl]\\
& Y_1 \ar_{d_0}[dr] && Y_1 \ar_{d_0}[dl] &\\
&& Y_0
}
\qquad
\xymatrix@C=15pt{
Y_2 \ar_{d_0}[dr] && Y_2 \ar_{d_0}[dl] \ar_{d_0}[dr] && Y_2 \ar_{d_0}[dl]\\
& Y_1 \ar_{d_0}[dr] && Y_1 \ar_{d_0}[dl] &\\
&& Y_0
}
\end{equation*}
That is,
\begin{equation*}
  \tens{(Y_2 \oset{d_1}{\times}_{Y_1}Y_2)}{(Y_1 \oset{d_0}{\times}_{Y_0}Y_1)} \cong
  \tens{(Y_2 \oset{d_0}{\times}_{Y_1}Y_2)}{(Y_1 \oset{d_0}{\times}_{Y_0}Y_1)}\;.
\end{equation*}
Similarly, the simplicial identity $d_0d_2=d_1d_0$ implies the isomorphism
\begin{equation*}
  \tens{(Y_2 \oset{d_2}{\times}_{Y_1}Y_2)}{(Y_1 \oset{d_0}{\times}_{Y_0}Y_1)} \cong
  \tens{(Y_2 \oset{d_0}{\times}_{Y_1}Y_2)}{(Y_1 \oset{d_1}{\times}_{Y_0}Y_1)}
\end{equation*}
and the simplicial identity $d_1d_2\cong d_1d_1$ implies the isomorphism
\begin{equation*}
  \tens{(Y_2 \oset{d_2}{\times}_{Y_1}Y_2)}{(Y_1 \oset{d_1}{\times}_{Y_0}Y_1)} \cong
  \tens{(Y_2 \oset{d_1}{\times}_{Y_1}Y_2)}{(Y_1 \oset{d_1}{\times}_{Y_0}Y_1)}
\end{equation*}
\end{example}
The face operators indicated in the picture are the respective projections, while we omitted drawing the degeneracies.

\clearpage

%%%%%%%%%%%%%%%%%%%%%%%%%%%%%%%%%%%%%%%%%%%%%%%%%%%%%%
% Figure (Corner Or2)

\begin{figure}[ht]
\bk
\begin{center}
$$
\entrymodifiers={+++[o]}
\xymatrix@R=30pt@C=10pt{
\cdots \ssr &
Y_{5} \ar@<2ex>[rrr]^{d_{5}} \ar[rrr]^{d_{4}}\ar@<-2ex>[rrr]^{d_{3}}
     \ar@<2.5ex>[d]^{d_{2}} \ar[d]^{d_{1}}\ar@<-2.5ex>[d]^{d_{0}} &&&
Y_{4} \ar@<1ex>[rrr]^{d_{4}} \ar@<-1ex>[rrr]^{d_{3}}
     \ar@<2.5ex>[d]^{d_{2}} \ar[d]^{d_{1}}\ar@<-2.5ex>[d]^{d_{0}} &&&
Y_{3}  \ar@<2.5ex>[d]^{d_{2}} \ar[d]^{d_{1}}\ar@<-2.5ex>[d]^{d_{0}} \\
\cdots \ssr  &
Y_{4} \ar@<2ex>[rrr]^{d_{4}} \ar[rrr]^{d_{3}}\ar@<-2ex>[rrr]^{d_{2}}
     \ar@<0.5ex>[d]^{d_{1}} \ar@<-0.5ex>[d]_{d_{0}} &&&
Y_{3} \ar@<1ex>[rrr]^{d_{3}} \ar@<-1ex>[rrr]^{d_{2}}
     \ar@<0.5ex>[d]^{d_{1}} \ar@<-0.5ex>[d]_{d_{0}} &&&
Y_{2}       \ar@<0.5ex>[d]^{d_{1}} \ar@<-0.5ex>[d]_{d_{0}} \\
\cdots \ssr &
Y_{3} \ar@<2ex>[rrr]^{d_{3}} \ar[rrr]^{d_{2}}\ar@<-2ex>[rrr]^{d_{1}} &&&
Y_{2} \ar@<1ex>[rrr]^{d_{2}} \ar@<-1ex>[rrr]^{d_{1}} &&& Y_{1}
}
$$
\end{center}
\caption{Corner of $\orn{2} Y$}
\label{CornerOr2}
\end{figure}

\bk

%%%%%%%%%%%%%%%%%%%%%%%%%%%%%%%%%%%%%%%
% Figure Corner double nerve

\begin{figure}[ht]
\bk
\begin{center}
$$
\entrymodifiers={++[o]}
\xymatrix@R=25pt@C=5pt{
 &&
Z \ar@<1ex>[rr]^{} \ar@<-1ex>[rr]^{}
 \ar@<2.5ex>[d]^{} \ar[d]^{}\ar@<-2.5ex>[d]^{} &&
Y_1 \oset{d_0}{\times}_{Y_0}Y_1 \oset{d_0}{\times}_{Y_0}Y_1 \ar@<2.5ex>[d]^{} \ar[d]^{}\ar@<-2.5ex>[d]^{} \\
%%%
W\ar@<2ex>[rr]^{} \ar[rr]^{}\ar@<-2ex>[rr]^{} \ar@<1ex>[d]^{} \ar@<-1ex>[d]_{} &&
(Y_2 \oset{d_2}{\times}_{Y_1}Y_2)\bsim \ar@<1ex>[rr]^{} \ar@<-1ex>[rr]^{}
\ar@<1ex>[d]^{} \ar@<-1ex>[d]_{} &&
Y_1 \oset{d_0}{\times}_{Y_0}Y_1 \ar@<1ex>[d]^{p_2} \ar@<-1ex>[d]_{p_1} \\
%%%
Y_1 \oset{d_1}{\times}_{Y_0}Y_1 \oset{d_1}{\times}_{Y_0}Y_1 \ar@<2ex>[rr]^{} \ar[rr]^{}\ar@<-2ex>[rr]^{} &&
Y_1 \oset{d_1}{\times}_{Y_0}Y_1  \ar@<1ex>[rr]^{p_2} \ar@<-1ex>[rr]_{p_1} && Y_{1}
}
$$
\end{center}
\caption{Corner of the double nerve  of $\clH_2 X$ for $Y=\clS X$}
\label{CornerDblNerve}
\end{figure}
where:
$$Z=\tens{(Y_2 \oset{d_2}{\times}_{Y_1}Y_2)\bsim}{(Y_1 \oset{d_1}{\times}_{Y_0}Y_1)}$$
$$W=\tens{(Y_2 \oset{d_2}{\times}_{Y_1}Y_2)\bsim}{(Y_1 \oset{d_0}{\times}_{Y_0}Y_1)}$$

\label{corner2page}

%%%%%%%%%%%%%%%%%%%%%%%%%%%%%%%%%%%%%%%%%%%%%%%%%%%%%%%
\clearpage

% Figure (Corner Or3)

\begin{figure}[ht]
\vspace{-5mm}
\begin{center}
$$
\entrymodifiers={+++[o]}
\xymatrix @R=20pt @C=35pt{
Y_{5}\ar@<0.7ex>[rr]^{d_{5}} \ar@<-0.7ex>[rr]_{d_{4}} \ar@<0.7ex>[dd]^{d_{1}}
    \ar@<-0.7ex>[dd]_{d_{0}} \ar@<0.7ex>[rd]^{d_{3}} \ar@<-0.7ex>[rd]_{d_{2}}
&&
Y_{4} \ar@{-}[d]<0.7ex> \ar@{-}[d]<-0.7ex> \ar@<0.7ex>[rd]^{d_{3}}
     \ar@<-0.7ex>[rd]_{d_{2}}\\
&
Y_{4} \ar@<0.7ex>[rr]^(0.3){d_{4}} \ar@<-0.7ex>[rr]_(0.3){d_{3}}
     \ar@<0.7ex>[dd]^(0.3){d_{1}} \ar@<-0.7ex>[dd]_(0.3){d_{0}}
&
    \text{  }\ar@<0.7ex>[d]^{d_{1}} \ar@<-0.7ex>[d]_{d_{0}}
&
Y_{3} \ar@<0.7ex>[dd]^{d_{1}} \ar@<-0.7ex>[dd]_{d_{0}}\\
Y_{4} \ar@<0.7ex>[rd]^{d_{2}} \ar@<-0.7ex>[rd]_{d_{1}} \ar@{-}[r]<0.7ex>
   \ar@{-}[r]<-0.7ex>  &\text{  }\ar@<0.7ex>[r]^{d_{4}} \ar@<-0.7ex>[r]_{d_{3}}
&
Y_{3} \ar@<0.7ex>[rd]^{d_{2}} \ar@<-0.7ex>[rd]_{d_{1}} \\
&
Y_{3} \ar@<0.7ex>[rr]^{d_{3}} \ar@<-0.7ex>[rr]_{d_{2}}
&&
Y_{2} }
$$
\end{center}
\vspace{-5mm}
\caption{Corner of $\orn{3} Y$}
\label{CornerOr3}
\end{figure}
%%%%%%%%%%%%%%%%%%%%%%%%%%%%%%%%%%%%%%%%%%%%%%%%%%%%%%%%%%%

% Figure (Corner Or3 nerve)

\begin{figure}[ht]
\vspace{-4mm}
\begin{center}
$$
\entrymodifiers={+++[o]}
\xymatrix @R=20pt @C=25pt{
Z_{5}\ar@<0.7ex>[rr]^{} \ar@<-0.7ex>[rr]_{} \ar@<0.7ex>[dd]^{}
    \ar@<-0.7ex>[dd]_{} \ar@<0.7ex>[rd]^{} \ar@<-0.7ex>[rd]_{}
&&
Z_{4}\up{0,1} \ar@{-}[d]<0.7ex> \ar@{-}[d]<-0.7ex> \ar@<0.7ex>[rd]^{}
     \ar@<-0.7ex>[rd]_{}\\
&
Z_{4}\up{0,2} \ar@<0.7ex>[rr]^(0.3){} \ar@<-0.7ex>[rr]_(0.3){}
     \ar@<0.7ex>[dd]^(0.3){} \ar@<-0.7ex>[dd]_(0.3){}
&
    \text{  }\ar@<0.7ex>[d]^{} \ar@<-0.7ex>[d]_{}
&
Z_{3}\up{0} \ar@<0.7ex>[dd]^{} \ar@<-0.7ex>[dd]_{}\\
%%%%%
Z_{4}\up{1,2} \ar@<0.7ex>[rd]^{} \ar@<-0.7ex>[rd]_{} \ar@{-}[r]<0.7ex>
   \ar@{-}[r]<-0.7ex>
&
\text{  }\ar@<0.7ex>[r]^{} \ar@<-0.7ex>[r]_{}
&
Z_{3}\up{1} \ar@<0.7ex>[rd]^{} \ar@<-0.7ex>[rd]_{} \\
&
Z_{3}\up{2} \ar@<0.7ex>[rr]^{} \ar@<-0.7ex>[rr]_{}
&&
Y_{2} }
$$
\end{center}

\caption{Corner of the 3-fold nerve of $\clH_3 X$, with $Y=\clS X$}
\label{CornerOr3nerve}
\end{figure}
%%Entries definitions
%
\vspace{-3mm}
\begin{equation*}
  Z_3\up{0}=Y_2 \oset{d_0}{\times}_{Y_1}Y_2, \quad Z_3\up{1}=Y_2 \oset{d_1}{\times}_{Y_1}Y_2, \quad Z_3\up{2}=Y_2 \oset{d_2}{\times}_{Y_1}Y_2\;.
\end{equation*}
\begin{equation*}
  Z_4\up{0,1}=\tens{(Y_2 \oset{d_1}{\times}_{Y_1}Y_2)}{(Y_1 \oset{d_0}{\times}_{Y_0}Y_1)}\cong
              \tens{(Y_2 \oset{d_0}{\times}_{Y_1}Y_2)}{(Y_1 \oset{d_0}{\times}_{Y_0}Y_1)}
\end{equation*}
\begin{equation*}
  Z_4\up{0,2}=\tens{(Y_2 \oset{d_2}{\times}_{Y_1}Y_2)}{(Y_1 \oset{d_0}{\times}_{Y_0}Y_1)}\cong
              \tens{(Y_2 \oset{d_0}{\times}_{Y_1}Y_2)}{(Y_1 \oset{d_1}{\times}_{Y_0}Y_1)}
\end{equation*}
\begin{equation*}
  Z_4\up{1,2}=\tens{(Y_2 \oset{d_2}{\times}_{Y_1}Y_2)}{(Y_1 \oset{d_1}{\times}_{Y_0}Y_1)}\cong
              \tens{(Y_2 \oset{d_1}{\times}_{Y_1}Y_2)}{(Y_1 \oset{d_1}{\times}_{Y_0}Y_1)}
\end{equation*}

\begin{equation*}
Z_5=Z_4\up{0,1}\bsim\;\cong\;Z_4\up{0,2}\bsim\;\cong\;Z_4\up{1,2}\bsim\;.
\end{equation*}

%
%%%%%%%%%%%%%%%%%%%%%%%%%%%%%%%%%%%%%%%%%%%%%%%%%%%%%%
% End figures
\label{corner3page}
\clearpage

%%%%%%%%%%%%%%%%%%%%%%%%%%%%%%%%%%%%%%%%%%%%%%%%%%%%%%%%%%%%%%%%%%%%%%%%%%%%%%
% Figure Corner J3 and J3D3

\begin{figure}
\vspace{20mm}
\begin{center}
\begin{equation*}
\def\labelstyle{\scriptstyle}
%\entrymodifiers={++[o]}
\ssr \xymatrix@R=35pt@C=20pt{
    \bm{\cdots}\ar@<1ex>[r] \ar[r] \ar@<-1ex>[r] \ar@<1ex>[d] \ar[d] \ar@<-1ex>[d] \ar@{}[dr] |{}
   &  \bm{\tens{X_{11}}{X_{10}}} \ar@<-0.5ex>[r] \ar@<0.5ex>[r] \ar@<1ex>[d] \ar[d] \ar@<-1ex>[d] \ar@{} [dr] |{}
   &  \textcolor[rgb]{1.00,0.00,0.00}{\bm{\tens{X_{01}}{X_{00}}}}  {\ar@<1ex>@[red][d] \ar@<0ex>@[red][d] \ar@<-1ex>@[red][d]} \\
     \bm{\cdots\quad \tens{X_{11}}{X_{01}} \quad}\ar[r] \ar@<1ex>[r]  \ar@<-1ex>[r] \ar@<0.5ex>[d] \ar@<-0.5ex>[d] \ar@{} [dr] |{}
    &  \bm{\quad  X_{11}\quad }\ar@<-0.5ex>[r] \ar@<0.5ex>[r]   \ar@<0.5ex>[d] \ar@<-0.5ex>[d] \ar@{} [dr] |{}
    &  \textcolor[rgb]{1.00,0.00,0.00}{\bm{\quad X_{01}}\quad}  \ar@<0.5ex>@[red][d] \ar@<-0.5ex>@[red][d]  \\
      \bm{\cdots\quad \textcolor[rgb]{1.00,0.00,0.00}{\tens{X_{10}}{X_{00}}} \quad} \ar@<0ex>@[red][r] \ar@<1ex>@[red][r]  \ar@<-1ex>@[red][r]
    &  \bm{\quad \textcolor[rgb]{1.00,0.00,0.00}{X_{10}}\quad} \ar@<-0.5ex>@[red][r] \ar@<0.5ex>@[red][r]
    &  \bm{\quad \textcolor[rgb]{1.00,0.00,0.00}{X_{00}}\quad} \\
}
\end{equation*}
\end{center}
\caption{Corner of $J_3 X\in\funcat{2}{\Cat}$ for $X\in\ftawg{3}$}
\label{Corner-J3-1}
\vspace{20mm}
\end{figure}

\bk

\begin{figure}
\begin{center}
\begin{equation*}
\def\labelstyle{\scriptstyle}
%\entrymodifiers={++[o]}
\ssr \xymatrix@R=35pt@C=20pt{
    \bm{\cdots}\ar@<1ex>[r] \ar[r] \ar@<-1ex>[r] \ar@<1ex>[d] \ar[d] \ar@<-1ex>[d] \ar@{} [dr] |{}
   &  \bm{\tens{X_{11}}{X_{10}}} \ar@<-0.5ex>[r] \ar@<0.5ex>[r] \ar@<1ex>[d] \ar[d] \ar@<-1ex>[d] \ar@{} [dr] |{}
   &  \textcolor[rgb]{1.00,0.00,0.00}{\bm{X^d_{0*}}}  {\ar@<1ex>@[red][d] \ar@<0ex>@[red][d] \ar@<-1ex>@[red][d]} \\
     \bm{\cdots\quad \tens{X_{11}}{X_{01}} \quad}\ar[r] \ar@<1ex>[r]  \ar@<-1ex>[r] \ar@<0.5ex>[d] \ar@<-0.5ex>[d] \ar@{} [dr] |{}
    &  \bm{\quad  X_{11}\quad }\ar@<-0.5ex>[r] \ar@<0.5ex>[r]   \ar@<0.5ex>[d] \ar@<-0.5ex>[d] \ar@{} [dr] |{}
    &  \textcolor[rgb]{1.00,0.00,0.00}{\bm{\quad X^d_{0*}}\quad}   \ar@<0.5ex>@[red][d] \ar@<-0.5ex>@[red][d]  \\
      \bm{\cdots\quad \textcolor[rgb]{1.00,0.00,0.00}{(\tens{X_{10}}{X_{00}})^d} \quad} \ar@<0ex>@[red][r] \ar@<1ex>@[red][r]  \ar@<-1ex>@[red][r]
    &  \bm{\quad \textcolor[rgb]{1.00,0.00,0.00}{X^d_{10}}\quad} \ar@<-0.5ex>@[red][r] \ar@<0.5ex>@[red][r]
    &  \bm{\quad \textcolor[rgb]{1.00,0.00,0.00}{X^d_{0*}}\quad} \\
}
\end{equation*}
\end{center}
\caption{Corner of $J_3 D_3 X\in\funcat{2}{\Cat}$ for $X\in\ftawg{3}$}
\label{Corner-J3-2}
\end{figure}

\label{Corner-J3}
\clearpage
%%%%%%%%%%%%%%%%%%%%%%%%%%%%%%%%%%%%%%%%%%%%%%%%%%%%%%%%%%%%%%%%%%%%%%%%%%%%%%

%%%%%%%%%%%%%%%%%%%%%%%%%%%%%%%%%%%%%%%%%%%%%%%%%%%%%%%%%%%%%%%%%%%%%%%%%

\appendix
  \chapter{Proof of Lemma \ref{lem-from-lta-to-pseu}}\label{app-a}

\pagelabel{page-app-a}

\textbf{Lemma \ref{lem-from-lta-to-pseu}}\it
\
 Let $X\in\catwg{n}$, $Y\in\lta{n}$ be such that $Y_{\uk}$ is discrete for all $\uk\in\dop{n-1}$ such that $k_j=0$ for some $1\leq j\leq n-1$. Let $\uk,\;\us \in\dop{n-1}$ and let $\uk \rw\us$ be a morphism in $\dop{n-1}$. Suppose that the following conditions hold:
\begin{itemize}
  \item [i)] If $k_j,s_j \neq 0$ for all $1\leq j\leq n-1$, then $ X_{\uk}=Y_{\uk},\quad X_{\us}=Y_{\us}$ and the maps
   \begin{equation*}
   X_{\uk}\rw X_{\us},\quad Y_{\uk}\rw Y_{\us}\;
   \end{equation*}
    coincide.

  \bk
  \item [ii)] If $k_j=0$ for some $1\leq j\leq n-1$ and $s_t=0$ for some $1\leq t\leq n-1$, then
  $X^d_{\uk}=Y_{\uk},\quad X^d_{\us}=Y_{\us}$ and the two maps
    \begin{equation*}
    X^d_{\uk}\rw X^d_{\us},\quad Y_{\uk}\rw Y_{\us}
    \end{equation*}
coincide, where $f^d:X^d_{\uk}\rw X^d_{\us}$ is induced by $f:X_{\uk}\rw X_{\us}$ and thus also coincides with the composite
  \begin{equation*}
    X^d_{\uk}\xrw{\zg'_{X_k}} X_{\uk} \xrw{f} X_{\us}\xrw{\zg_{X_s}}X^d_{\us}
  \end{equation*}
  (where $\zg$ is the discretization map and $\zg'$ a section), since $f^d=f^d \zg_{X_k}\zg'_{X_k}=\zg_{X_s}f\zg'_{X_k}$.
  \bk

  \item [iii)] If $k_j\neq 0$ for all $1\leq j\leq n-1$ and $s_t=0$ for some $1\leq t\leq n-1$, the following diagram commutes
  \begin{equation*}
  \qquad\qquad\qquad\xymatrix{
  X_{\uk} \ar[r] \ar@{=}[d] & X_{\us} \ar^(0.3){\zg_{X_{\us}}}[rr] && X^d_{\us}=Y_{\us}\mbox{\qquad\qquad\qquad\qquad\qquad} \\
  Y_{\uk}\ar[rrru]
  }
  \end{equation*}
  where $\zg_{X_s}$ is the discretization map.

  \bk
  \item [iv)] If $k_j=0$ for some $1\leq j\leq n-1$ and $s_t\neq 0$ for all $1\leq t\leq n-1$ then the following diagram commutes
  \begin{equation*}
  \qquad\qquad\qquad\xymatrix{
  X^d_{\uk} \ar^{\zg'_{X_{\uk}}}[rr] \ar@{=}[d] && X_{\uk} \ar[r] & X_{\us}=Y_{\us}\mbox{\qquad\qquad\qquad\qquad} \\
  Y_{\uk}\ar[rrru]
  }
  \end{equation*}
  Then
  \begin{itemize}
    \item [a)] For all $\uk\in\dop{n-1}$, $(\tr{n}X)_{\uk}=(\tr{n}Y)_{\uk}$\bk

    \item [b)] For all $\uk\in\dop{n-1}$ such that $k_j\neq 0$ for all $1\leq j\leq n-1$, the maps
    \begin{equation*}
      (\tr{n}X)_{\uk}\rightleftarrows X_{\uk}, \qquad (\tr{n}Y)_{\uk}\rightleftarrows Y_{\uk}
    \end{equation*}
    coincide. \bk

    \item [c)] $\tr{n}X=\tr{n}Y$\;.
  \end{itemize}

\end{itemize}
\rm

\begin{proof}
By induction on $n$. Let $n=2$. By definition of $\tr{2}$ and by conditions i) and ii) in the hypothesis,
\begin{equation*}
\begin{split}
    & (\tr{2}X)_0=X_0^d =Y_0 =(\tr{2}Y)_0 \\
    & (\tr{2}X)_1=X_1 =Y_1 =(\tr{2}Y)_1\;.
\end{split}
\end{equation*}
Further, by hypothesis iii) the maps $\pt'_0,\pt'_1:Y_1\rw Y_0$ are the composites
\begin{equation*}
\xymatrix{
Y_1=X_1 \ar@<1ex>^{\pt_0}[rr] \ar@<-1ex>_{\pt_0}[rr] && X_0 \ar^(0.4){\zg_{X_0}}[rr] && X_0^d=Y_0\\
}
\end{equation*}
that is $\pt'_i=\zg_{X_0}\pt_i$, $i=0,1$. This implies that for each $k\geq 2$
\begin{equation*}
  (\tr{2}X)_k = \pro{X_1}{X_0^d}{k} = \pro{Y_1}{Y_0}{k} = (\tr{2}Y)_k \;.
\end{equation*}
This proves a) when $n=2$.

We now prove b) when $n=2$. The map
\begin{equation*}
  (\tr{2}X)_1 = X_1\rw X_1
\end{equation*}
is the identity and by the hypothesis i) coincides with
\begin{equation*}
  (\tr{2}Y)_1 = Y_1\rw Y_1 = X_1\;.
\end{equation*}
When $k>1$, the maps
\begin{equation}\label{eq1-lem-from-lta-00}
  (\tr{2}X)_k = \pro{X_1}{X_0^d}{k}\rightleftarrows \pro{X_1}{X_0}{k}
\end{equation}
are the induced Segal maps for $X$ and their pseudo-inverses. The induced Segal maps of $X$ arise from the commuting diagram (see also Definition \ref{def-ind-seg-map}).
\begin{equation}\label{eq2-lem-from-lta-0}
\xymatrix@C=10pt{
&&& X_k \ar_{\nu_1}[dll] \ar^{\nu_k}[drr]&&&\\
& X_1 \ar_{\zg\pt_0}[dl] \ar^{\zg\pt_1}[dr] && \cdots && X_1 \ar_{\zg\pt_0}[dl] \ar^{\zg\pt_1}[dr] &\\
X_0^d && X_0^d  & \cdots & X_0^d && X_0^d \\
}
\end{equation}
By hypothesis i) the maps
\begin{equation*}
  \nu_i:X_k\rw X_1 \qquad \nu_i:Y_k\rw Y_1
\end{equation*}
coincide; by hypothesis iii) the maps
\begin{equation*}
  X_1 \xrw{\zg\pt_i} X_0^d \qquad Y_1 \xrw{\pt'_i} Y_0
\end{equation*}
coincide. Thus \eqref{eq2-lem-from-lta-0} coincides with
\begin{equation*}
\xymatrix@C=10pt{
&&& Y_k \ar_{\nu_1}[dll] \ar^{\nu_k}[drr]&&&\\
& Y_1 \ar_{\pt'_{0}}[dl] \ar^{\pt'_{1}}[dr] && \cdots && Y_1 \ar_{\pt'_{0}}[dl] \ar^{\pt'_{1}}[dr] &\\
Y_0 && Y_0  & \cdots & Y_0 && Y_0 \\
}
\end{equation*}
so the induced Segal maps of $X$ and $Y$ coincide. So by \eqref{eq1-lem-from-lta-00} the maps
\begin{equation*}
  (\tr{2}X)_k \rightleftarrows X_k, \qquad (\tr{2}Y)_k \rightleftarrows Y_k
\end{equation*}
coincide, proving b) when $n=2$.

To show c) when $n=2$ we first show that, for each morphism $k \rw s$ in $\dop{}$, the maps
\begin{equation}\label{eqre-maps-same}
   (\tr{2}X)_k  \rw  (\tr{2}X)_s,\qquad (\tr{2}Y)_k  \rw  (\tr{2}Y)_s
\end{equation}
coincide. By the proof of Lemma \ref{lem-PP} these maps are the composites
\begin{equation}\label{eq3-lem-from-lta}
  (\tr{2}X)_k \rw X_k \rw X_s \rw (\tr{2}X)_s
\end{equation}
\begin{equation}\label{eq4-lem-from-lta}
  (\tr{2}Y)_k \rw Y_k \rw Y_s \rw (\tr{2}Y)_s
\end{equation}
Let $k>0$ and $s>0$. Then by b) and by hypothesis i), the maps \eqref{eq3-lem-from-lta} and \eqref{eq4-lem-from-lta} coincide.

Let $k=0$ and $s>0$. Then \eqref{eq3-lem-from-lta} and \eqref{eq4-lem-from-lta} are given by
\begin{equation*}
  (\tr{2}X)_0 = X^d_0 \rw X_0 \rw X_s=(\tr{2}X)_s
\end{equation*}
\begin{equation*}
  (\tr{2}Y)_0 = Y_0 \rw Y_s
\end{equation*}
and these coincide by hypothesis iv).

Suppose $k>0$ and $s=0$. Then  \eqref{eq3-lem-from-lta} and \eqref{eq4-lem-from-lta}  are given by
\begin{equation}\label{eq5-lem-from-lta}
  (\tr{2}X)_k \rw X_k \rw X_0 \rw X_0^d=(\tr{2}X)_0
\end{equation}
\begin{equation}\label{eq6-lem-from-lta}
  (\tr{2}Y)_k \rw Y_k \rw Y_0 = (\tr{2}Y)_0\;.
\end{equation}
By b), the maps
\begin{equation*}
  (\tr{2}X)_k \rw X_k,\qquad (\tr{2}Y)_k \rw Y_k
\end{equation*}
coincide while by hypothesis iii) the maps
\begin{equation*}
  X_k\rw X_0 \rw X_0^d,\qquad Y_k\rw Y_0
\end{equation*}
coincide. Therefore \eqref{eq5-lem-from-lta} and \eqref{eq6-lem-from-lta} coincide.

If $k=s=0$, the composite
\begin{equation*}
  (\tr{2}X)_0= X^d_0 \rw X_0 \xrw{f_{X_0}} X_0 \rw X_0^d
\end{equation*}
is equal to $f^d_{X_0}$, which coincides with $f^d_{Y_0}$ by hypothesis ii).

 We conclude that \eqref{eq3-lem-from-lta} and \eqref{eq4-lem-from-lta} always coincide.

By the definition of pseudo-functor (see Definition \ref{def-pseudo-fun}) in order to prove that $Tr_2 X=Tr_2 Y$ it remains to show that, given morphisms $k \rw s \rw r$ in $\dop{}$ the 2-dimensional pasting diagrams

\begin{equation*}
\begin{tikzcd}
(\tr{2}X)_{k} \arrow[r] \arrow[rr, bend right=25, "\Downarrow"]
&
(\tr{2}X)_{s} \arrow[r]
&
(\tr{2}X)_{r}
\end{tikzcd}
\end{equation*}

\begin{equation*}
\begin{tikzcd}
(\tr{2}Y)_{k} \arrow[r] \arrow[rr, bend right=25, "\Downarrow"]
&
(\tr{2}Y)_{s} \arrow[r]
&
(\tr{2}Y)_{r}
\end{tikzcd}
\end{equation*}
coincide and, given $id_{k}: k \rw k$ in $\dop{}$, the 2-dimensional pasting diagrams

\begin{equation*}
\begin{tikzcd}
(\tr{2}X)_{k} \arrow[rr, "(\tr{2}X)(\Id_{k})"] \arrow[rr, bend right=25, "\Downarrow", "\Id"']
&&
(\tr{2}X)_{k}
\end{tikzcd}
\end{equation*}

\begin{equation*}
\begin{tikzcd}
(\tr{2}Y)_{k} \arrow[rr, "(\tr{2}Y)(\Id_{k})"] \arrow[rr, bend right=25, "\Downarrow", "\Id"']
&&
(\tr{2}Y)_{k}
\end{tikzcd}
\end{equation*}
coincide.

The proof of this is as in the inductive step on pages \pageref{pr-lem-part2} to \pageref{pr-lem-part3}. In fact, the proof of these parts of the inductive step only uses a) and b) and the equality of the maps \eqref{eqre-maps-same}, all of which have been proved for the case $n=2$, but it does not use the equality of the 2-dimensional pasting diagrams at step $(n-1)$. We therefore refer the reader to the later part of this proof for this step. This concludes the proof of the lemma in the case $n=2$.

\bk

Suppose, inductively, that the lemma holds for $(n-1)$ and let $X,Y$ be as in the hypothesis.

a) By definition of $\tr{n}$ and by hypothesis ii), for all $\us\in\dop{n-2}$,
\begin{equation*}
 (\tr{n}X)_{(0,\us)}=X^d_{(0,\us)}=Y^d_{(0,\us)}=(\tr{n}Y)_{(0,\us)}\;.
\end{equation*}
Clearly $X_1\in\catwg{n-1}$ and $Y_1\in\lta{n-1}$ satisfy the inductive hypothesis. Thus, using the definition of $\tr{n}$ and the inductive hypothesis a) on $X_1,Y_1$ we obtain
\begin{equation*}
  (\tr{n}X)_{(1,\us)}=(\tr{n-1}X)_{\us}=(\tr{n-1}Y)_{\us}=(\tr{n}Y)_{(1,\us)}\;.
\end{equation*}
We claim that, for all $\us\in\dop{n-2}$, the maps
\begin{equation}\label{eq7-lem-from-lta}
  (\tr{n-1}X_1)_{\us} \rw X_{1\us}\rw X_{0\us}\rw X^d_{0\us}
\end{equation}
\begin{equation}\label{eq8-lem-from-lta}
  (\tr{n-1}Y_1)_{\us} \rw Y_{1\us}\xrw{\hspace{18mm}} Y_{0\us}
\end{equation}
coincide. In fact, suppose $s_j\neq 0$ for all $1\leq j\leq n-2$. Then by inductive hypothesis b) applied to $X_1,Y_1$ the maps
\begin{equation*}
 (\tr{n-1}X_1)_{\us} \rw X_{1\us}, \qquad (\tr{n-1}Y_1)_{\us} \rw Y_{1\us}
\end{equation*}
coincide, while by hypothesis iii) the maps
\begin{equation*}
  X_{1\us}\rw X_{0\us}\rw X^d_{0\us}, \qquad Y_{1\us}\rw Y_{0\us}
\end{equation*}
coincide. Thus the composites \eqref{eq7-lem-from-lta} and \eqref{eq8-lem-from-lta} coincide.

Suppose $s_j=0$ for some $1\leq j\leq n-2$. Then by Corollary \ref{cor-trn}
\begin{equation*}
   (\tr{n-1}X_1)_{\us}=X^d_{1\us}, \qquad (\tr{n-1}Y_1)_{\us}=Y^d_{1\us}\;.
\end{equation*}
Thus the maps \eqref{eq7-lem-from-lta} and \eqref{eq8-lem-from-lta} are given by
\begin{equation*}
   X^d_{1\us}\rw X_{1\us}\rw X_{0\us}\rw X^d_{0\us}, \qquad  Y_{1\us}\rw Y_{0\us}
\end{equation*}
and these coincide by hypothesis ii). This proves the claim. From this claim, the inductive hypothesis a) on $X_1,Y_1$ and the definition of $\tr{n}$ it follows that, for each $k_1>1$, $\uk=(k_1,\us)$
\begin{equation*}
\begin{split}
    & (\tr{n}X)_{\uk}=\pro{(\tr{n-1}X_1)_{\us}}{X^d_{0\us}}{k_1}= \\
    & =\pro{(\tr{n-1}Y_1)_{\us}}{Y_{0}}{k_1}=(\tr{n}Y)_{\uk}
\end{split}
\end{equation*}
This proves a).
\bk

b) Let $\uk=(k_1, \us)$ with $k_1 \neq 0$ and $s_j\neq 0$ for all $1\leq j\leq n-2$. By induction hypothesis b)
applied to $X_1$ and $Y_1$, the maps
\begin{equation*}
\begin{split}
    & (\tr{n}X)_{(1,\us)}=(\tr{n-1}X_1)_{\us}\rightleftarrows X_{1\us} \\
    & (\tr{n}Y)_{(1,\us)}=(\tr{n-1}Y_1)_{\us}\rightleftarrows Y_{1\us}
\end{split}
\end{equation*}
coincide. This implies that, for each $k_1>1$, the maps
\begin{equation}\label{eq9-lem-from-lta}
  (\tr{n}X)_{\uk}=\pro{(\tr{n-1}X_1)_{\us}}{X^d_{0\us}}{k_1}\rw \pro{X_{1\us}}{X^d_{0\us}}{k_1}
\end{equation}
\begin{equation}\label{eq10-lem-from-lta}
  (\tr{n}Y)_{\uk}=\pro{(\tr{n-1}Y_1)_{\us}}{Y_{0\us}}{k_1}\rw \pro{Y_{1\us}}{Y_{0\us}}{k_1}
\end{equation}
coincide. On the other hand, reasoning as in the case $n=2$, we see that the maps
\begin{equation}\label{eq11-lem-from-lta}
  \pro{X_{1\us}}{X^d_{0\us}}{k_1}\rightleftarrows \pro{X_{1\us}}{X_{0\us}}{k_1}=X_{k_1\us}
\end{equation}
\begin{equation}\label{eq12-lem-from-lta}
  \pro{Y_{1\us}}{Y_{0\us}}{k_1}\rightleftarrows Y_{k_1\us}
\end{equation}
coincide. Composing \eqref{eq9-lem-from-lta} with \eqref{eq11-lem-from-lta} and \eqref{eq10-lem-from-lta} with \eqref{eq12-lem-from-lta} we therefore conclude that the maps
\begin{equation*}
  (\tr{n}X)_{\uk} \rightleftarrows X_{\uk}, \qquad (\tr{n}Y)_{\uk} \rightleftarrows Y_{\uk}
\end{equation*}
coincide for each $\uk$ such that $k_j\neq 0$ for all $1\leq j\leq n-1$. This proves b).
\bk

c) By a) and by the definition of pseudo-functor (see Definition \ref{def-pseudo-fun}) in order to prove that $Tr_n X=Tr_n Y$ it remains to show that:
\begin{itemize}
  \item [i)] For each morphism $\uk\rw\us$ in $\dop{n-1}$, the maps
\begin{equation*}
  (\tr{n}X)_{\uk}\rw (\tr{n}X)_{\us}, \qquad (\tr{n}Y)_{\uk}\rw (\tr{n}Y)_{\us}
\end{equation*}
coincide.

  \item [ii)] Given morphisms $\uk \rw \us \rw \ur$ in $\dop{n-1}$ the 2-dimensional pasting diagrams

\begin{equation}\label{diagram-2a}
\begin{tikzcd}
(\tr{n}X)_{\uk} \arrow[r] \arrow[rr, bend right=25, "\Downarrow"]
&
(\tr{n}X)_{\us} \arrow[r]
&
(\tr{n}X)_{\ur}
\end{tikzcd}
\end{equation}

\begin{equation}\label{diagram-2b}
\begin{tikzcd}
(\tr{n}Y)_{\uk} \arrow[r] \arrow[rr, bend right=25, "\Downarrow"]
&
(\tr{n}Y)_{\us} \arrow[r]
&
(\tr{n}Y)_{\ur}
\end{tikzcd}
\end{equation}

  coincide.

  \item [iii)] Given $id_{\uk}: \uk \rw \uk$ in $\dop{n-1}$, the 2-dimensional pasting diagrams

\begin{equation}\label{diagram-3a}
\begin{tikzcd}
(\tr{n}X)_{\uk} \arrow[rr, "(\tr{n}X)(\Id_{\uk})"] \arrow[rr, bend right=25, "\Downarrow", "\Id"']
&&
(\tr{n}X)_{\uk}
\end{tikzcd}
\end{equation}

\begin{equation}\label{diagram-3b}
\begin{tikzcd}
(\tr{n}Y)_{\uk} \arrow[rr, "(\tr{n}Y)(\Id_{\uk})"] \arrow[rr, bend right=25, "\Downarrow", "\Id"']
&&
(\tr{n}Y)_{\uk}
\end{tikzcd}
\end{equation}

  coincide.

\end{itemize}

i) By the proof of Lemma \ref{lem-PP}, these maps are given as composites
\begin{equation}\label{eq13-lem-from-lta}
  (\tr{n}X)_{\uk}\rw X_{\uk}\rw X_{\us}\rw (\tr{n}X)_{\us}
\end{equation}
\begin{equation}\label{eq14-lem-from-lta}
  (\tr{n}Y)_{\uk}\rw Y_{\uk}\rw Y_{\us}\rw (\tr{n}Y)_{\us}\;.
\end{equation}
Suppose that $k_j\neq 0$, $s_j\neq 0$ for all $1\leq j\leq n-1$. Then by b) and by hypothesis i), \eqref{eq13-lem-from-lta} and \eqref{eq14-lem-from-lta} coincide.

Suppose that $k_j=0$ for some $1\leq j\leq n-1$ and $s_t=0$ for some $1\leq t\leq n-1$. Then by Corollary \ref{cor-trn}, \eqref{eq13-lem-from-lta} and \eqref{eq14-lem-from-lta} are given by the composites
\begin{equation*}
\begin{split}
    & (\tr{n}X)_{\uk}=X^d_{\uk}\rw X_{\uk}\rw X_{\us}\rw (\tr{n}X)_{\us}=X^d_{\us} \\
    & (\tr{n}Y)_{\uk}=Y_{\uk}\xrw{\hspace{49mm}} Y_{\us}
\end{split}
\end{equation*}
and these coincide by hypothesis ii).
Suppose that $k_j=0$ for some $1\leq j\leq n-1$ and $s_i\neq 0$ for all $1\leq i\leq n-1$. Then \eqref{eq13-lem-from-lta} and \eqref{eq14-lem-from-lta} are given by the composites
\begin{equation}\label{eq15-lem-from-lta}
  (\tr{n}X)_{\uk}=X^d_{\uk}\rw X_{\uk}\rw X_{\us}\rw (\tr{n}X)_{\us}
\end{equation}
\begin{equation}\label{eq16-lem-from-lta}
  (\tr{n}Y)_{\uk}=Y_{\uk}\xrw{\hspace{16mm}} Y_{\us}\rw(\tr{n}Y)_{\us} \;.
\end{equation}
By a) the maps
\begin{equation*}
  X_{\us}\rw (\tr{n}X)_{\us}, \qquad Y_{\us}\rw(\tr{n}Y)_{\us}
\end{equation*}
coincide, with the hypothesis iv) the maps
\begin{equation*}
  X^d_{\uk}\rw X_{\uk}\rw X_{\us}, \qquad Y_{\uk}\rw Y_{\us}
\end{equation*}
coincide. Hence by composing, we deduce that \eqref{eq15-lem-from-lta} and \eqref{eq16-lem-from-lta} coincide.
In conclusion the maps \eqref{eq13-lem-from-lta} and \eqref{eq14-lem-from-lta} always coincide, proving i).

ii)\pagelabel{pr-lem-part2}  We distinguish the following eight cases. For each we refer to diagrams on pages \pageref{cases1-4} and \pageref{cases5-8}. Using a), b) and the hypotheses we see that in each case the left and right pasting diagrams coincide.\bk

\begin{itemize}
  \item [\emph{Case 1}:]$k_i\neq 0$ for all $1 \leq i \leq (n-1)$; $s_j\neq 0$ for all $1 \leq j \leq (n-1)$; $r_t\neq 0$ for all $1 \leq t \leq (n-1)$.\bk

  \item [\emph{Case 2}:]  $k_i= 0$ for some $1 \leq i \leq (n-1)$; $s_j\neq 0$ for all $1 \leq j \leq (n-1)$; $r_t\neq 0$ for all $1 \leq t \leq (n-1)$.

      \nid Note that by hypothesis iv) the map $Y_{\uk}\rw Y_{\us}$ factors as
      \begin{equation*}
        Y_{\uk}=X_{\uk}^d \rw X_{\uk} \rw X_{\us}=Y_{\us}.
      \end{equation*}
 \bk

  \item [\emph{Case 3}:]  $k_i\neq 0$ for all $1 \leq i \leq (n-1)$; $s_j= 0$ for some $1 \leq j \leq (n-1)$;  $r_t\neq 0$ for all $1 \leq t \leq (n-1)$.

      \nid Note that by hypotheses iii) and iv) the maps
      \begin{equation*}
        Y_{\uk}\rw Y_{\us}, \qquad Y_{\us}\rw Y_{\ur}
      \end{equation*}
      factor as
      \begin{equation*}
        Y_{\uk}=X_{\uk}\rw X_{\us} \rw X_{\us}^d=Y_{\us}, \qquad Y_{\us}=X_{\us}^d\rw X_{\us} \rw X_{\ur}=Y_{\ur}.
      \end{equation*}
      \bk
 \item [\emph{Case 4}:]  $k_i\neq 0$ for all $1 \leq i \leq (n-1)$; $s_j\neq 0$ for all $1 \leq j \leq (n-1)$; $r_t= 0$ for some $1 \leq t \leq (n-1)$.

     \nid Note that by hypothesis iv) the map $Y_{\us}\rw Y_{\ur}$ factors as
      \begin{equation*}
        Y_{\us}=X_{\us} \rw X_{\ur} \rw X_{\ur}^d=Y_{\ur}.
      \end{equation*}
\bk
  \item [\emph{Case 5}:] $k_i\neq 0$ for all $1 \leq i \leq (n-1)$; $s_j= 0$ for some $1 \leq j \leq (n-1)$; $r_t= 0$ for some $1 \leq t \leq (n-1)$.

      \nid Note that by hypotheses ii) and iv) the maps
      \begin{equation*}
        Y_{\uk}\rw Y_{\us}, \qquad Y_{\us}\rw Y_{\ur}
      \end{equation*}
      factor as
      \begin{equation*}
        Y_{\uk}=X_{\uk}\rw X_{\us} \rw X_{\us}^d=Y_{\us}, \qquad Y_{\us}=X_{\us}^d\rw X_{\us} \rw X_{\ur}\rw \rw X_{\ur}^d=Y_{\ur}.
      \end{equation*}
\bk
  \item [\emph{Case 6}:]  $k_i= 0$ for some $1 \leq i \leq (n-1)$; $s_j\neq 0$ for all $1 \leq j \leq (n-1)$; $r_t= 0$ for some $1 \leq t \leq (n-1)$.

  \nid Note that by hypotheses iii) and iv) the maps
      \begin{equation*}
        Y_{\uk}\rw Y_{\us}, \qquad Y_{\us}\rw Y_{\ur}
      \end{equation*}
      factor as
      \begin{equation*}
        Y_{\uk}=X_{\uk}^d\rw X_{\uk} \rw X_{\us}=Y_{\us}, \qquad Y_{\us}=X_{\us}\rw X_{\ur} \rw X_{\ur}^d=Y_{\ur}.
      \end{equation*}
  \bk
  \item [\emph{Case 7}:]$k_i= 0$ for some $1 \leq i \leq (n-1)$; $s_j= 0$ for some $1 \leq j \leq (n-1)$; $r_t\neq 0$ for all $1 \leq t \leq (n-1)$.

      \nid Note that by hypotheses ii) and iv) the maps
      \begin{equation*}
        Y_{\uk}\rw Y_{\us}, \qquad Y_{\us}\rw Y_{\ur}
      \end{equation*}
      factor as
      \begin{equation*}
        Y_{\uk}=X_{\uk}^d\rw X_{\uk} \rw X_{\us}\rw X_{\us}^d=Y_{\us}, \qquad Y_{\us}=X_{\us}^d\rw X_{\us}\rw X_{\ur}=Y_{\ur}.
      \end{equation*}
\bk
   \item [\emph{Case 8}:] $k_i= 0$ for some $1 \leq i \leq (n-1)$; $s_j= 0$ for some $1 \leq j \leq (n-1)$; $r_t=0$ for some $1 \leq t \leq (n-1)$.

   \nid Note that by hypothesis ii) the maps
      \begin{equation*}
        Y_{\uk}\rw Y_{\us}, \qquad Y_{\us}\rw Y_{\ur}
      \end{equation*}
      factor as
      \begin{equation*}
        Y_{\uk}=X_{\uk}^d\rw X_{\uk} \rw X_{\us}\rw X_{\us}^d=Y_{\us}, \qquad Y_{\us}=X_{\us}^d\rw X_{\us}\rw X_{\ur}\rw X_{\ur}^d=Y_{\ur}.
      \end{equation*}

\end{itemize}

\bk

iii)\pagelabel{pr-lem-part3}  Suppose that $k_i \neq 0$ for all $1 \leq i \leq (n-1)$. Then \ref{diagram-3a} and \ref{diagram-3b} are given by
\begin{equation*}
\begin{tikzcd}
(\tr{n}X)_{\uk} \arrow[r] \arrow[d]\arrow[dd, bend right=40, "\Rightarrow","\Id"']
&
(\tr{n}X)_{\uk} \arrow[d, shift left=0.5ex] \arrow[dd, bend left=40, "\Rightarrow"',"\Id"]\\
%%%%%%%%%%%%%%%%%%%%%%%%%%%%%%%%%%%%%%%%%%%%%%%%%%%%%%%%%%%%%%%%%%%%%%%%%%%%%%%
X_{\uk} \arrow[d] \arrow[r, "\Id"]
&
X_{\uk} \arrow[d] \arrow[u, shift left=0.5ex]\\
%%%%%%%%%%%%%%%%%%%%%%%%%%%%%%%%%%%%%%%%%%%%%%%%%%%%%%%%%%%%%%%%%%%%%%%%%%%%%%%
(\tr{n}X)_{\uk} \arrow[r, equal]
&
(\tr{n}X)_{\uk}
\end{tikzcd}
\qquad
\begin{tikzcd}
(\tr{n}Y)_{\uk} \arrow[r] \arrow[d]\arrow[dd, bend right=40, "\Rightarrow","\Id"']
&
(\tr{n}Y)_{\uk} \arrow[d, shift left=0.5ex] \arrow[dd, bend left=40, "\Rightarrow"',"\Id"]\\
%%%%%%%%%%%%%%%%%%%%%%%%%%%%%%%%%%%%%%%%%%%%%%%%%%%%%%%%%%%%%%%%%%%%%%%%%%%%%%%
Y_{\uk} \arrow[d] \arrow[r, "\Id"]
&
Y_{\uk} \arrow[d] \arrow[u, shift left=0.5ex]\\
%%%%%%%%%%%%%%%%%%%%%%%%%%%%%%%%%%%%%%%%%%%%%%%%%%%%%%%%%%%%%%%%%%%%%%%%%%%%%%%
(\tr{n}Y)_{\uk} \arrow[r, equal]
&
(\tr{n}Y)_{\uk}
\end{tikzcd}
\end{equation*}
and these coincide by a), b) and hypothesis i).

Suppose that $k_i=0$ for some $1 \leq i \leq (n-1)$. Then \ref{diagram-3a} and \ref{diagram-3b} are given by
\begin{equation*}
\begin{tikzcd}
X^d_{\uk} \arrow[r, "\Id"] \arrow[d]\arrow[dd, bend right=40, "\Id"']
&
X^d_{\uk} \arrow[d, shift left=0.5ex] \arrow[dd, bend left=40, "\Id"]\\
%%%%%%%%%%%%%%%%%%%%%%%%%%%%%%%%%%%%%%%%%%%%%%%%%%%%%%%%%%%%%%%%%%%%%%%%%%%%%%%
X_{\uk} \arrow[d] \arrow[r, "\Id"]
&
X_{\uk} \arrow[d] \arrow[u, shift left=0.5ex]\\
%%%%%%%%%%%%%%%%%%%%%%%%%%%%%%%%%%%%%%%%%%%%%%%%%%%%%%%%%%%%%%%%%%%%%%%%%%%%%%%
X_{\uk} \arrow[r, "\Id"']
&
X_{\uk}
\end{tikzcd}
\qquad
\begin{tikzcd}
Y_{\uk} \arrow[rr, "\Id"] \arrow[d, equal]
& &
Y_{\uk} \arrow[d, equal] \\
%%%%%%%%%%%%%%%%%%%%%%%%%%%%%%%%%%%%%%%%%%%%%%%%%%%%%%%%%%%%%%%%%%%%%%%%%%%%%%%
Y_{\uk} \arrow[r] \arrow[d, equal] \arrow[rr, bend right=30, "\Id"']
&
X_{\uk} \arrow[r]
&
Y_{\uk} \arrow[d, equal]\\
%%%%%%%%%%%%%%%%%%%%%%%%%%%%%%%%%%%%%%%%%%%%%%%%%%%%%%%%%%%%%%%%%%%%%%%%%%%%%%%
Y_{\uk} \arrow[rr, "\Id"']
& &
Y_{\uk}
\end{tikzcd}
\end{equation*}

and they coincide as $Y_{\uk}= X_{\uk}^d$.

\clearpage

\tikzcdset{column sep/tiny=0.9ex}
\tikzcdset{column sep/small=1.5ex}
CASE 1
\begin{equation*}
\begin{tikzcd}[column sep=tiny]
\Sc{(\tr{n}X)_{\uk}}\arrow[rr]\arrow[d] \arrow[dd, bend right=45, "\Id"', "{\Rightarrow}"]
& &
|[alias=DZ]|\Sc{(\tr{n}X)_{\us}} \arrow[rr] \arrow[dr]
&&
\Sc{(\tr{n}X)_{\ur}} \arrow[d, shift left=0.5ex] \arrow[dd, bend left=45, "\Id", "{\Rightarrow}"' ]\\
%%%%
|[alias=BZ]|\Sc{X_{\uk}} \arrow[r]\arrow[d, shift right=0.5ex]\arrow[rrrr, bend right]
&
\Sc{X_{\us}} \arrow[ru]\arrow[rr, equal,""{name=EQ}] \arrow[phantom, "{\Downarrow}", from=DZ, to=EQ]
& &
\Sc{X_{\us}} \arrow[r] & |[alias=HZ]|\Sc{X_{\ur}} \arrow[d] \arrow[u, shift left=0.5ex]\\ %\arrow[Rightarrow, from=HZ, to=R]
%%%%
\Sc{(\tr{n}X)_{\uk}} \arrow[rrrr]\arrow[u, shift right=0.5ex] &&&& \Sc{(\tr{n}X)_{\ur}}
\end{tikzcd}
\quad
\begin{tikzcd}[column sep=tiny]
\Sc{(\tr{n}Y)_{\uk}}\arrow[rr]\arrow[d] \arrow[dd, bend right=45, "\Id"' , "{\Rightarrow}"]
& &
|[alias=DZ]|\Sc{(\tr{n}Y)_{\us}} \arrow[rr] \arrow[dr]
&&
\Sc{(\tr{n}Y)_{\ur}} \arrow[d, shift left=0.5ex] \arrow[dd, bend left=45, "\Id", "{\Rightarrow}"' ]\\
%%%%
|[alias=BZ]|\Sc{Y_{\uk}} \arrow[r]\arrow[d, shift right=0.5ex]\arrow[rrrr, bend right]
&
\Sc{Y_{\us}} \arrow[ru]\arrow[rr, equal,""{name=EQ}] \arrow[phantom, "{\Downarrow}", from=DZ, to=EQ]
& &
\Sc{Y_{\us}} \arrow[r] & |[alias=HZ]|\Sc{Y_{\ur}}  \arrow[d] \arrow[u, shift left=0.5ex]\\
%%%%
\Sc{(\tr{n}Y)_{\uk}} \arrow[rrrr]\arrow[u, shift right=0.5ex] &&&& \Sc{(\tr{n}Y)_{\ur}}
\end{tikzcd}
\end{equation*}

\vspace{10mm}
CASE 2
\begin{equation*}
\begin{tikzcd}[column sep=small]
\Sc{X^d_{\uk}}\arrow[rr]\arrow[d] \arrow[dd, bend right=45, "\Id"' ]
& &
|[alias=DZ]|\Sc{(\tr{n}X)_{\us}} \arrow[rr] \arrow[dr]
&&
\Sc{(\tr{n}X)_{\ur}} \arrow[d, shift left=0.5ex] \arrow[dd, bend left=45, "\Id", "{\Rightarrow}"' ]\\
%%%%
\Sc{\!\!\!\!\!\!X_{\uk}} \arrow[r]\arrow[d]\arrow[rrrr, bend right]
&
\Sc{X_{\us}} \arrow[ru]\arrow[rr, equal,""{name=EQ}] \arrow[phantom, "{\Downarrow}", from=DZ, to=EQ]
& &
\Sc{X_{\us}} \arrow[r] & |[alias=HZ]|\Sc{X_{\ur}}  \arrow[d] \arrow[u, shift left=0.5ex]\\
%%%%
\Sc{X^d_{\uk}} \arrow[rrrr] &&&& \Sc{(\tr{n}X)_{\ur}}
\end{tikzcd}
\quad
\begin{tikzcd}[column sep=tiny]
\Sc{Y_{\uk}}\arrow[rrr]\arrow[d, equal]
& & &
|[alias=DZ]|\Sc{(\tr{n}Y)_{\us}} \arrow[rr] \arrow[dr]
&&
\Sc{(\tr{n}Y)_{\ur}} \arrow[d, shift left=0.5ex] \arrow[dd, bend left=45, "\Id", "{\Rightarrow}"' ]\\
%%%%
|[alias=BZ]|\Sc{Y_{\uk}} \arrow[r]\arrow[d, equal]\arrow[rrrrr, bend right]
&
\Sc{X_{\uk}}  \arrow[r]
&
\Sc{Y_{\us}} \arrow[ru]\arrow[rr, equal,""{name=EQ}] \arrow[phantom, "{\Downarrow}", from=DZ, to=EQ]
& &
\Sc{Y_{\us}} \arrow[r] & |[alias=HZ]|\Sc{Y_{\ur}}  \arrow[d] \arrow[u, shift left=0.5ex]\\
%%%%
\Sc{Y_{\uk}} \arrow[rrrrr]
&&&&&
\Sc{(\tr{n}Y)_{\ur}}
\end{tikzcd}
\end{equation*}

\vspace{10mm}

CASE 3
\begin{equation*}
\begin{tikzcd}[column sep=tiny]
\Sc{(\tr{n}X)_{\uk}}\arrow[rr]\arrow[d] \arrow[dd, bend right=45, "\Id"' , "{\Rightarrow}"]
& &
|[alias=DZ]|\Sc{X^d_{\us}} \arrow[rr] \arrow[dr]
&&
\Sc{(\tr{n}X)_{\ur}} \arrow[d, shift left=0.5ex] \arrow[dd, bend left=45, "\Id", "{\Rightarrow}"' ]\\
%%%%
|[alias=BZ]|\Sc{X_{\uk}} \arrow[r]\arrow[d, shift right=0.5ex]\arrow[rrrr, bend right]
&
\Sc{X_{\us}} \arrow[ru]\arrow[rr, equal,""{name=EQ}] \arrow[phantom, "{\Downarrow}", from=DZ, to=EQ]
& &
\Sc{X_{\us}} \arrow[r] & |[alias=HZ]|\Sc{X_{\ur}}  \arrow[d] \arrow[u, shift left=0.5ex]\\
%%%%
\Sc{(\tr{n}X)_{\uk}} \arrow[rrrr]\arrow[u, shift right=0.5ex] &&&& \Sc{(\tr{n}X)_{\ur}}
\end{tikzcd}
\quad
\begin{tikzcd}[column sep=small]
\Sc{(\tr{n}Y)_{\uk}}\arrow[rr]\arrow[d] \arrow[dd, bend right=45, "\Id"' ,, "{\Rightarrow}"]
& &
|[alias=DZ]|\Sc{Y_{\us}} \arrow[rr] \arrow[d, equal]
&&
\Sc{(\tr{n}Y)_{\ur}} \arrow[d, shift left=0.5ex] \arrow[dd, bend left=45, "\Id", "{\Rightarrow}"' ]\\
%%%%
|[alias=BZ]|\Sc{Y_{\uk}} \arrow[r]\arrow[d, shift right=0.5ex]\arrow[rrrr, bend right=45]
&
\Sc{X_{\us}} \arrow[r] \arrow[rr, bend right=50, "\Id"', "{\Downarrow}"]
&
\Sc{Y_{\us}} \arrow[r]
&
\Sc{X_{\us}} \arrow[r] & |[alias=HZ]|\Sc{Y_{\ur}}  \arrow[d] \arrow[u, shift left=0.5ex]\\
%%%%
\Sc{(\tr{n}Y)_{\uk}} \arrow[rrrr]\arrow[u, shift right=0.5ex]
&&&&
\Sc{(\tr{n}Y)_{\ur}}
\end{tikzcd}
\end{equation*}

\vspace{10mm}

CASE 4
\begin{equation*}
\begin{tikzcd}[column sep=small]
\Sc{(\tr{n}X)_{\uk}}\arrow[rr]\arrow[d] \arrow[dd, bend right=40, "\Id"' , "{\Rightarrow}"]
& &
|[alias=DZ]|\Sc{(\tr{n}X)_{\us}} \arrow[rr] \arrow[dr]
&&
\Sc{X^d_{\ur}} \arrow[d, shift left=0.5ex] \arrow[dd, bend left=40, "\Id"]\\
%%%%
|[alias=BZ]|\Sc{X_{\uk}} \arrow[r]\arrow[d, shift right=0.5ex]\arrow[rrrr, bend right]
&
\Sc{X_{\us}} \arrow[ru]\arrow[rr, equal,""{name=EQ}] \arrow[phantom, "{\Downarrow}", from=DZ, to=EQ]
& &
\Sc{X_{\us}} \arrow[r] & |[alias=HZ]|\Sc{X_{\ur}}  \arrow[d] \arrow[u, shift left=0.5ex]\\
%%%%
\Sc{(\tr{n}X)_{\uk}} \arrow[rrrr]\arrow[u, shift right=0.5ex] &&&& \Sc{X^d_{\ur}}
\end{tikzcd}
\quad
\begin{tikzcd}[column sep=small]
\Sc{(\tr{n}Y)_{\uk}}\arrow[rr]\arrow[d] \arrow[dd, bend right=40, "\Id"' , "{\Rightarrow}"]
& &
|[alias=DZ]|\Sc{(\tr{n}Y)_{\us}} \arrow[rrr] \arrow[dr]
&&&
\Sc{Y_{\ur}} \arrow[d, equal] \\
%%%%
|[alias=BZ]|\Sc{Y_{\uk}} \arrow[r]\arrow[d, shift right=0.5ex]\arrow[rrrrr, bend right]
&
\Sc{Y_{\us}} \arrow[ru]\arrow[rr, equal,""{name=EQ}] \arrow[phantom, "\Downarrow", from=DZ, to=EQ]
&&
\Sc{Y_{\us}} \arrow[r]
&
\Sc{X_{\ur}} \arrow[r]
&
\Sc{Y_{\ur}}  \arrow[d, equal] \\
%%%%
\Sc{(\tr{n}Y)_{\uk}} \arrow[rrrrr]\arrow[u, shift right=0.5ex] &&&&& \Sc{Y_{\ur}}
\end{tikzcd}
\end{equation*}

\pagelabel{cases1-4}

\vspace{10mm}

CASE 5
\begin{equation*}
\begin{tikzcd}[column sep=small]
\Sc{(\tr{n}X)_{\uk}}\arrow[rr]\arrow[d] \arrow[dd, bend right=40, "\Id"' , "{\Rightarrow}"]
& &
|[alias=DZ]|\Sc{X^d_{\us}} \arrow[rr] \arrow[dr]
&&
\Sc{X^d_{\ur}} \arrow[d] \arrow[dd, bend left=40, "\Id"]\\
%%%%
|[alias=BZ]|\Sc{X_{\uk}} \arrow[r]\arrow[d, shift right=0.5ex]\arrow[rrrr, bend right]
&
\Sc{X_{\us}} \arrow[ru]\arrow[rr, equal,""{name=EQ}] \arrow[phantom, "{\Downarrow}", from=DZ, to=EQ]
& &
\Sc{X_{\us}} \arrow[r] & |[alias=HZ]|\Sc{X_{\ur}}  \arrow[d] \\
%%%%
\Sc{(\tr{n}X)_{\uk}} \arrow[rrrr]\arrow[u, shift right=0.5ex] &&&& \Sc{X^d_{\ur}}
\end{tikzcd}
\quad
\begin{tikzcd}[column sep=small]
\Sc{(\tr{n}Y)_{\uk}}\arrow[rr]\arrow[d] \arrow[dd, bend right=40, "\Id"' , "{\Rightarrow}"]
& &
|[alias=DZ]|\Sc{Y_{\us}} \arrow[rrr] \arrow[d, equal]
&&&
\Sc{Y_{\ur}} \arrow[d, equal] \\
%%%%
\Sc{Y_{\uk}} \arrow[r]\arrow[d, shift right=0.5ex]\arrow[rrrrr, bend right=35]
&
\Sc{X_{\us}} \arrow[r] \arrow[rr, bend right=40, "\Scc{\Downarrow}", "{\Scc{\Id}}"']
&
\Sc{Y_{\us}} \arrow[r]
&
\Sc{X_{\us}} \arrow[r]
&
\Sc{X_{\ur}} \arrow[r]
&
\Sc{Y_{\ur}}  \arrow[d, equal] \\
%%%%
\Sc{(\tr{n}Y)_{\uk}} \arrow[rrrrr]\arrow[u, shift right=0.5ex] &&&&& \Sc{Y_{\ur}}
\end{tikzcd}
\end{equation*}

\vspace{10mm}

CASE 6
\begin{equation*}
\begin{tikzcd}[column sep=small]
\Sc{X^d_{\uk}}\arrow[rr]\arrow[d] \arrow[dd, bend right=30, "\Id"']
& &
|[alias=DZ]|\Sc{(\tr{n}X)_{\us}} \arrow[rr] \arrow[dr]
&&
\Sc{X^d_{\ur}} \arrow[d] \arrow[dd, bend left=30, "\Id"]\\
%%%%
|[alias=BZ]|\Sc{X_{\uk}} \arrow[r]\arrow[d, shift right=0.5ex]\arrow[rrrr, bend right]
&
\Sc{X_{\us}} \arrow[ru]\arrow[rr, equal,""{name=EQ}] \arrow[phantom, "{\Downarrow}", from=DZ, to=EQ]
& &
\Sc{X_{\us}} \arrow[r] & |[alias=HZ]|\Sc{X_{\ur}}  \arrow[d] \\
%%%%
\Sc{X^d_{\uk}} \arrow[rrrr]\arrow[u, shift right=0.5ex] &&&& \Sc{X^d_{\ur}}
\end{tikzcd}
\quad
\begin{tikzcd}[column sep=tiny]
\Sc{Y_{\uk}}\arrow[rrr]\arrow[d, equal]
& & &
|[alias=DZ]|\Sc{(\tr{n}Y)_{\us}} \arrow[rrr] \arrow[dr]
& & &
\Sc{Y_{\ur}} \arrow[d, equal] \\
%%%%
\Sc{Y_{\uk}} \arrow[r]\arrow[d, equal]\arrow[rrrrrr, bend right]
&
\Sc{X_{\uk}} \arrow[r]
&
\Sc{Y_{\us}} \arrow[ru]\arrow[rr, equal,""{name=EQ}] \arrow[phantom, "\Downarrow", from=DZ, to=EQ]
&&
\Sc{Y_{\us}} \arrow[r]
&
\Sc{X_{\ur}} \arrow[r]
&
\Sc{Y_{\ur}}  \arrow[d, equal] \\
%%%%
\Sc{Y_{\uk}} \arrow[rrrrrr] &&&&&& \Sc{Y_{\ur}}
\end{tikzcd}
\end{equation*}

\vspace{10mm}

CASE 7
\begin{equation*}
\begin{tikzcd}[column sep=tiny]
\Sc{X^d_{\uk}}\arrow[rr]\arrow[d] \arrow[dd, bend right=30, "\Id"']
& &
|[alias=DZ]|\Sc{X^d_{\us}} \arrow[rr] \arrow[dr]
&&
\Sc{(\tr{n}X)_{\ur}} \arrow[d, shift left=0.5ex] \arrow[dd, bend left=40, "\Id","\Rightarrow"']\\
%%%%
|[alias=BZ]|\Sc{X_{\uk}} \arrow[r]\arrow[d, shift right=0.5ex]\arrow[rrrr, bend right]
&
\Sc{X_{\us}} \arrow[ru]\arrow[rr, equal,""{name=EQ}] \arrow[phantom, "{\Downarrow}", from=DZ, to=EQ]
& &
\Sc{X_{\us}} \arrow[r]
&
|[alias=HZ]|\Sc{X_{\ur}}  \arrow[d] \arrow[u, shift left=0.5ex] \\
%%%%
\Sc{X^d_{\uk}} \arrow[rrrr]\arrow[u, shift right=0.5ex]
&&&&
\Sc{(\tr{n}X){\ur}}
\end{tikzcd}
\quad
\begin{tikzcd}[column sep=small]
\Sc{Y_{\uk}}\arrow[rrr]\arrow[d, equal]
& & &
\Sc{Y_{\us}} \arrow[rr] \arrow[d, equal]
& &
\Sc{(\tr{n}Y)_{\ur}} \arrow[d, shift left=0.5ex] \arrow[dd, bend left=40, "\Id","\Rightarrow"']\\
%%%%%%%%%%%%%%%%%%%%%%%%%%%%%%%%%%%%%%%%%
\Sc{Y_{\uk}} \arrow[r]\arrow[d, equal]\arrow[rrrrr, bend right=40]
&
\Sc{X_{\uk}} \arrow[r]
&
\Sc{X_{\us}} \arrow[r] \arrow[rr, bend right=40, "\Scc{\Downarrow}", "{\Scc{\Id}}"']
&
\Sc{Y_{\us}} \arrow[r]
&
\Sc{X_{\us}} \arrow[r]
&
\Sc{Y_{\ur}} \arrow[u, shift left=0.5ex] \arrow[d] \\
%%%%%%%%%%%%%%%%%%%%%%%%%%%%%%%%%%%%%%%%%%%%%%%%%%%%%
\Sc{Y_{\uk}} \arrow[rrrrr]
&&&&&
\Sc{(\tr{n}Y)_{\ur}}
\end{tikzcd}
\end{equation*}

\vspace{10mm}

CASE 8
\begin{equation*}
\begin{tikzcd}[column sep=small]
\Sc{X^d_{\uk}}\arrow[rr]\arrow[d] \arrow[dd, bend right=30, "\Id"']
& &
|[alias=DZ]|\Sc{X^d_{\us}} \arrow[rr] \arrow[dr]
&&
\Sc{X^d_{\ur}} \arrow[d] \arrow[dd, bend left=30, "\Id"]\\
%%%%
|[alias=BZ]|\Sc{X_{\uk}} \arrow[r]\arrow[d]\arrow[rrrr, bend right]
&
\Sc{X_{\us}} \arrow[ru]\arrow[rr, equal,""{name=EQ}] \arrow[phantom, "{\Downarrow}", from=DZ, to=EQ]
& &
\Sc{X_{\us}} \arrow[r] & |[alias=HZ]|\Sc{X_{\ur}}  \arrow[d] \\
%%%%
\Sc{X^d_{\uk}} \arrow[rrrr] &&&& \Sc{X^d_{\ur}}
\end{tikzcd}
\quad
\begin{tikzcd}[column sep=small]
\Sc{Y_{\uk}} \arrow[rrr] \arrow[d, equal]
&&&
\Sc{Y_{\us}} \arrow[rrr] \arrow[d, equal]
&&&
\Sc{Y_{\ur}} \arrow[d, equal] \\
%%%%%%%%%%%%%%%%%%%%%%%%%%%%%%%%%%%%%%%%%%%%%%%%%%%%%
\Sc{Y_{\uk}} \arrow[r] \arrow[d, equal] \arrow[rrrrrr, bend right=35]
&
\Sc{X_{\uk}} \arrow[r]
&
\Sc{X_{\us}} \arrow[r]  \arrow[rr, bend right=40, "\Scc{\Downarrow}", "{\Scc{\Id}}"']
&
\Sc{Y_{\us}} \arrow[r]
&
\Sc{X_{\us}} \arrow[r]
&
\Sc{X_{\ur}} \arrow[r]
&
\Sc{Y_{\ur}} \arrow[d, equal]\\
%%%%%%%%%%%%%%%%%%%%%%%%%%%%%%%%%%%%%%%%%%%%%%%%%%%%
\Sc{Y_{\uk}} \arrow[rrrrrr]
&&&&&&
\Sc{Y_{\ur}}
\end{tikzcd}
\end{equation*}

\pagelabel{cases5-8}

\end{proof}
\clearpage

%%%%%%%%%%%%%%%%%%%%%%%%%%%%%%%%%%%%%%%%%%%%%%%%%%%%%%%%%%%%%%%%%%%%%%%%%
\backmatter
% Define Page style for backmatter
\pagestyle{fancy}
% Delete the current section for header and footer
\fancyhf{}
% Set custom header
\lhead[]{\thepage}
\rhead[\thepage]{}

%    Bibliography styles amsplain or harvard are also acceptable.
%\bibliographystyle{amsplain}

%%%%%%%%%%%%%%%%%%%%%%%%%%%%%%%%%%%%%%%%%%%%%%%%%%%%%%%%%%%%%%%%%%%%%%%

\listoffigures
%\phantomsection
%\addcontentsline{toc}{section}{List of Figures}

%%%%%%%%%%%%%%%%%%%%%%%%%%%%%%%%%%%%%%%%%%%%%%%%%%%%%%%%%%%%%%%%%%%%%%%

\printindex

%%%%%%%%%%%%%%%%%%%%%%%%%%%%%%%%%%%%%%%%%%%%%%%%%%%%%%%%%%%%%%%%%%%%%

\chapter*{Index of terminology and notation}\vspace{-5mm}
\begin{center}
\renewcommand{\arraystretch}{1.5}
\small
\begin{longtable}{l p{6.7cm} l}

%%%%%%%%%%%%%%%%%%%%%%%%%%%%%%%%%%%%%%%%%%%%%%%%%%%%%%%%%%%%%%%%%%%%
 $\Af$ & Internal equivalence relation corresponding to $f:A\rw B$ & Definition \ref{def-int-eq-rel} \\
 $B$ & Classifying space & Definition \ref{def-class-sp-funct}\\
 &  &  Sections \ref{sec-group-wg-nfol-cat}, \ref{mod-fund-wg-group}\\
 $\cathd{n}$ & Category of homotopically discrete \mbox{$n$-fold} categories & Definition \ref{def-hom-dis-ncat} \\
 $\Cat \clC$ & Category of internal categories in $\clC$ & Definition \ref{def-intercat}\\
 $\cat{n}(\clC)$ & Category of \nfol internal categories in $\clC$ & Section \ref{sbs-nint-cat}\\
 $\catwg{n}$ & Category of weakly globular $n$-fold categories & Definition \ref{def-n-equiv} \\
 $\Cube{n,t}$ & Set of $(n,t)$-hypercubes & Definition \ref{def-3-proper-ltawg} \\
 $Dec$ & Functor d\'{e}calage & Section \ref{decalage} \\
 $\Delta$ & Simplicial category & Section \ref{subs-multisimplob}\\
  $\Dnop$ & Product of $n$ copies of $\dop{}$ & Section \ref{subs-multisimplob}\\
 $\funcat{}{\clC}$ & Category of simplicial objects in $\clC$ & Section \ref{subs-simplob} \\
 $\funcat{n}{\clC}$ & Category of \nfol simplicial objects in $\clC$ & Section \ref{not-simp}\\
  $\zgu{n}_X$ &  Map $X\rw \di{n}\p{n}X$ for $X\in\cathd{n}$ & Definition \ref{def-hom-dis-ncat-2}\\
 $Diag_n$ & Multi-diagonal & Definition \ref{def-class-sp-funct}\\
 $\Discn$  & Discretization functor  &  Definition \ref{def-disc-func} \\
 $\di{n}$ & Functor $\seg{n-1}\rw\seg{n}$   &  Definition \ref{dn}, \\
 & & Lemma \ref{lem-wg-ps-cat-b},\\
  & & Remark \ref{rem-wg-ps-cat-b} \\
  $\di{n,j}$ & Functor $\seg{j-1}\rw\seg{n}$ & Notation \ref{dnot-ex-tam}\\
 $\Dn$   & Functor $\ftawg{n}\rw \ta{n}$  & Proposition  \ref{pro-fta-tam-1} \\
 $\eqr{n}$ & Category of $n$-equivalence relations & Definition \ref{def1-int-eq-rel} \\
 $F_n$ &  Functor $\catwg{n}\rw \catwg{n}$  &  Proposition \ref{pro-gen-const-1} \\
 $f_n(X)$ & Map $F_n X\rw X$ & Proposition \ref{pro-gen-const-1} \\
 $\ftawg{n}$  &   &  Definition \ref{def-fta-1} \\
 $\gcatwg{n}$  & Category of groupoidal weakly globular $n$-fold categories  & Definition \ref{def-gta-2} \\
 $g_n(X)$ & Map $G_n X\rw X$ & Theorem \ref{pro-fta-1}\\
 $G_n$  & Functor $\catwg{n}\rw \ftawg{n}$  &  Theorem \ref{pro-fta-1}  \\
 $\mathcal{G}_n$ & Fundamental groupoidal weakly globular $n$-fold category functor & Section \ref{mod-fund-wg-group}\\
 $\Gpd\,\clC$ & Category of internal groupoids in $\clC$ & Definition \ref{def-intergroup}\\
 $\gpdwg{n}$ & Category of weakly globular \nfol groupoids &  Section \ref{mod-fund-wg-group}\\
 $\gseg{n}$ & Groupoidal Segal-type model & Section \ref{subs-main-res} \\
 $\gta{n}$  & Category of groupoidal Tamsamani \mbox{$n$-categories}  & Definition \ref{def-gta-2} \\
 $\gtawg{n}$  &  Category of groupoidal weakly globular Tamsamani $n$-categories & Definition \ref{def-gta-1} \\
 $\mathcal{H}_n$& Fundamental weakly globular $n$-fold groupoid functor & Definition \ref{def-hn}\\
 $h_n (X)$ & Map $V_n (X)\rw F_n (X)$, $X\in\cathd{n}$ & Proof of Proposition \ref{pro-gen-const-2}\\
 $\Ho\text{($n$-types)}$ & Homotopy category of $n$-types & Section \ref{homotyp}\\
 $J_n$ & Functor $\seg{n}\rw\funcat{n-1}{\Cat}$  &  Notation \ref{not-ner-func-dirk} \\
 &  &  Definition \ref{def-wg-ps-cat}\\
 $\lnta{n}{n}$ &   &  Definition \ref{def-ltawg-equiv} \\
 $\sf{L}\tawg{n}$ &  &  Definition \ref{def-ind-sub-ltawg} \\
 $n\mi\Cat$ & Category of strict $n$-categories & Definition \ref{def-strict}\\
 $n\mi\Gpd$ & Category of strict $n$-groupoids & Section \ref{sbs-multi-strict}\\
 $N$ & Nerve functor $N:\Cat\clC \rw \funcat{}{\clC}$ & Section \ref{sus-ner-funct}\\
 $\Nu{k}$ & Nerve functor in the $k^{th}$ direction & Definition \ref{def-ner-func-dirk} \\
 $\Nb{n}$ & Multinerve functor & Definition \ref{def-multinerve}\\
  $\orn{2}$ & Functor $\funcat{}{\Set}\rw\funcat{2}{\Set}$ & Figure \ref{CornerOr2}\\
 $\orn{3}$ & Functor $\funcat{}{\Set}\rw\funcat{3}{\Set}$& Figure \ref{CornerOr3}\\
  $\orn{n}$ & Functor $\funcat{}{\Set}\rw\funcat{n}{\Set}$ & Section \ref{mod-fund-wg-group}\\
  $or_{n}$  & Ordinal sum $\Delta^n\rw\Delta$ & Section \ref{mod-fund-wg-group}\\
 $p$ & Isomorphism classes of object functor & Section \ref{sbs-funct-cat} \\
 $\p{j,n}$ & Functor $\seg{n}\rw \seg{j-1}$ & Definition \ref{def-pn},\\
 & &  Notation \ref{not-ex-tam}\\
 $\p{n}$ & Functor $\seg{n} \rw \seg{n-1}$ & Definition \ref{def-n-equiv},\\
 & & Definition \ref{def-wg-ps-cat} \\
 &  & Definition \ref{def-hom-dis-ncat}\\
 $P_n$ & Functor $P_n:\tawg{n}\rw\lta{n}$ & Proof of Theorem \ref{the-funct-Qn}\\
 $ \clP_n$ & Left adjoint to \nfol nerve & Section \ref{mod-fund-wg-group}\\
 $\psc{n}{\Cat}$ & Category of pseudo-functors  & Section \ref{sbs-pseudo-functors} \\
 $q$ & Connected components functor & Section \ref{sbs-funct-cat} \\
 $\q{j,n}$ & Functor $\seg{n}\rw \seg{j-1}$ & Proposition \ref{pro-post-trunc-fun},\\
 & &  Remark \ref{rem-qrn}\\
 $\q{n}$ & Functor $\seg{n} \rw \seg{n-1}$ & Proposition \ref{pro-post-trunc-fun},\\
  & & Corollary \ref{pro-post-wg-ncat} \\
 $Q_n$ & Rigidification functor  & Theorem \ref{the-funct-Qn} \\
 $\rz$ & Functor $\ftawg{n}\rw\funcat{}{\ftawg{n-1}}$ & Definition \ref{def-fta-tam-1}\\
 $\Rbt{n}$ & Composite $\mathcal{P}_n \orn{n}$ & Equation Section \ref{sub-last-examples}\\
 $ \clS$ & Singular functor $\Top \rw\funcat{}{\Set}$ & Section \ref{mod-fund-wg-group}\\
 $\seg{n}$ & Segal-type model & Section \ref{sec-comm-fea}\\
 $\segpsc{n}{\Cat}$  & Category of Segalic pseudo-functors & Definition \ref{def-seg-ps-fun} \\
 $s_n(X)$ & Map $Q_n X\rw X$ & Theorem \ref{the-funct-Qn} \\
 $\St$ & Strictification functor & Section \ref{sbs-pseudo-functors} \\
 $\Tan$ & Category of Tamsamani $n$-categories & Example \ref{ex-tam} \\
 $\tawg{n}$ & Category of weakly globular Tamsamani $n$- categories & Definition \ref{def-wg-ps-cat} \\
 $t_n (X)$  & Pseudo-natural transformation  $Tr_{n}X\rw X$ &  Theorem \ref{the-XXXX} \\
 $v_k$ & Map $J_{n-1} X_k \rw $   & Notation \ref{not-for-lta} \\
   & $J_{n-1}(\pro{X_1}{\di{n-1}\p{n-1}X_0}{k})$ &\\
 $\mathcal{T}_n$ & Fundamental Tamsamani $n$-groupoid functor & Theorem \ref{cor-gta-2}\\
 $Tr_{n}$  & Functor $\lta{n} \rw \segpsc{n-1}{\Cat}$ & Theorem \ref{the-XXXX} \\
 $v_{n}(X)$   & Map $V_n X\rw X$ & Proposition \ref{pro-gen-const-2} \\
 $V_n$  & Functor $\cathd{n}\rw\cathd{n}$  & Proposition \ref{pro-gen-const-2} \\
 $V(X)$  & Map $X(f_0)\rw X$ & Lemma \ref{lem-gen-constr}\\
 &  &   Proposition \ref{pro-gen-const-1-new} \\
 $X(a,b)$ & Hom-$(n-1)$-category of $X\in\seg{n}$ & Notation \ref{not-fiber} \\
 &  & Definition \ref{def-n-equiv}\\
 &  &  Definition \ref{def-wg-ps-cat}\\
 $X(f_0)$  & Construction on $X\in\Cat \clC$ and $f_0:X'_0 \rw X_0$ & Lemma \ref{lem-gen-constr} \\
 $X^d$ & Discretization of $X\in \cathd{n}$ &  Definition \ref{def-hom-dis-ncat}  \\
 $X\up{2}_k$  & Same as $(\Nu{2}X)_k$, $X\in \cat{n}$ & Definition \ref{def-kdir-wg-ncat} \\
 $\zg\lo{n}$ & Discretization map $X\rw X^d$ for $X\in\cathd{n}$ & Definition \ref{def-hom-dis-ncat-2}\\
 $\zgu{n}$ &  Map $X\rw \dn\qn X$ for $X\in\tawg{n}$ & Proposition \ref{pro-post-trunc-fun}\\
 $\mu_k$ & $k^{th}$ Segal map & Definition \ref{def-seg-map}\\
 $\hmu{k}$ & $k^{th}$ induced Segal map & Definition \ref{def-ind-seg-map}\\
 $\xi_i$ & Map $\funcat{n}{\clC}\rw\funcat{}{\funcat{n-1}{\clC}}$ & Lemma \ref{lem-multi-simpl-as}\\
 $\tilde{\xi_i}$ & Map $\cat{n}(\clC) \rw \Cat(\cat{n-1}(\clC))$ & Proposition \ref{pro-mult-ner}\\
 $w_{n}(X)$ & Map $P_n X \rw X$ & Proof of Theorem \ref{the-funct-Qn}\\

\end{longtable}
\normalsize
\end{center}
\label{pagina}
%%%%%%%%%%%%%%%%%%%%%%%%%%%%%%%%%%%%%%%%%%%%%%%%%%%%%%%%%%%%%%%%%%%%%%%%%%%%%%%

\printbibliography

\address{{Department of Mathematics, University of Leicester, LE17RH, UK}}

\email{sp424@le.ac.uk}

\end{document}